\documentclass[a4paper, 12pt, twoside]{memo-l}
\usepackage[multiple]{footmisc}
\usepackage{currfile,datetime}
\usepackage{amsthm,amsmath,amssymb,bbm,geometry,epsfig,listings, paralist, enumerate,comment,nicefrac,verbatim,multirow,xcolor}
\usepackage{mathrsfs}  
\usepackage{mathtools}
\usepackage{xparse}

\usepackage[hyphens]{url}
\usepackage[hyphenbreaks]{breakurl}

\usepackage[colorlinks,linkcolor=red,citecolor=green]{hyperref}

\usepackage[maxbibnames=9, backend=biber]{biblatex}
\numberwithin{section}{chapter}
\newcommand{\var}{\vartheta}
\newcommand{\varz}{\vartheta}
\renewcommand{\H}{\mathbb{H}}
\newcommand{\U}{\mathbb{U}}
\newcommand{\Span}{\operatorname{span}}
\newcommand{\Div}{\operatorname{div}}

\newcommand{\id}{\operatorname{Id}}
\newcommand{\citationand}{\&}
\newcommand{\ud}{\,\mathrm{d}}	
\newcommand{\ds}{\ud s}
\newcommand{\dWs}{\ud W_s}
\newcommand{\eps}{{\ensuremath{\varepsilon}}}
\newcommand{\C}{\operatorname{C}}

\providecommand{\dist}{{\operatorname{dist}}}
\allowdisplaybreaks
\def\CO{\mathcal{O}}
\DeclareUnicodeCharacter{0229}{\c{e}}

\usepackage[capitalise,compress]{cleveref} 

\crefname{equation}{}{}

\newtheorem{lemma}{Lemma}[chapter]
\newtheorem{remark}[lemma]{Remark}

\newtheorem{theorem}[lemma]{Theorem}

\newtheorem{definition}[lemma]{Definition}
\newtheorem{prop}[lemma]{Proposition}
\newtheorem{corollary}[lemma]{Corollary}

\newtheorem{sett}[lemma]{Setting}

\makeatletter
\def\blx@maxline{77}
\makeatother

\crefname{subsection}{Subsection}{Subsections}
\crefname{enumi}{item}{items}
\creflabelformat{enumi}{(#2\textup{#1}#3)}

\DeclareFontEncoding{LS1}{}{}
\DeclareFontSubstitution{LS1}{stix}{m}{n}
\DeclareMathAlphabet{\mathscr}{LS1}{stixscr}{m}{n}

\newcommand{\vast}[2]{\left#2 \rule{0pt}{#1}\kern-.25ex\right.}

\newcommand{\R}{\mathbbm{R}}
\newcommand{\N}{\mathbbm{N}}
\newcommand{\Z}{\mathbbm{Z}}
\newcommand{\1}{\mathbbm{1}}
\renewcommand{\P}{\mathbbm{P}}
\newcommand{\E}{\mathbbm{E}}

\newcommand{\funcH}[2]{{\left\vert\kern-0.25ex\left\vert\kern-0.25ex\left\vert #1     \right\vert\kern-0.25ex\right\vert\kern-0.25ex\right\vert}_{2,#2}}

\newcommand{\funcN}[2]{{\left\vert\kern-0.25ex\left\vert\kern-0.25ex\left\vert #1     \right\vert\kern-0.25ex\right\vert\kern-0.25ex\right\vert}_{1,#2}}

\NewDocumentCommand{\fabs}{sO{}m}{%
  {\IfBooleanTF{#1}
    {\fabsaux{\left|}{\right|}{#3}}
    {\fabsaux{#2|}{#2|}{#3}}}
}
\makeatletter
\newcommand{\fabsaux}[3]{\mathpalette\fabsaux@i{{#1}{#2}{#3}}}
\newcommand{\fabsaux@i}[2]{\fabsaux@ii#1#2}
\newcommand{\fabsaux@ii}[4]{%
  \sbox\z@{$\m@th#1#2#4#3$}%
  \sbox\tw@{$\m@th\|$}%
  \mathopen{\hbox to\wd\tw@{\hss\vrule height \ht\z@ depth \dp\z@ width .3\wd\tw@\hss}}%
  \mkern-2mu #4 \mkern-2mu 
  \mathclose{\hbox to\wd\tw@{\hss\vrule height \ht\z@ depth \dp\z@ width .3\wd\tw@\hss}}%
}
\makeatother

\NewDocumentCommand{\ffabs}{som}{%
  {\IfBooleanTF{#1}
    {\fabsaux{\left|}{\right|}{#3}}
    {\IfNoValueTF{#2}
      {\fnabsaux{|}{|}{#3}}
      {\fabsaux{#2|}{#2|}{#3}}
    }
  }
}
\makeatletter
\newcommand{\fnabsaux}[3]{\mathpalette\fnabsaux@i{{#1}{#2}{#3}}}
\newcommand{\fnabsaux@i}[2]{\fnabsaux@ii#1#2}
\newcommand{\fnabsaux@ii}[4]{%
  \sbox\z@{$\m@th#1#2#4#3$}%
  \sbox\tw@{$\m@th\|$}%
  \mathopen{\hbox to\wd\tw@{\hss\vrule height .8\ht\z@ depth .5\dp\z@ width .3\wd\tw@\hss}}%
  \mkern-2mu #4 \mkern-2mu 
  \mathclose{\hbox to\wd\tw@{\hss\vrule height .8\ht\z@ depth .5\dp\z@ width .3\wd\tw@\hss}}%
}
\makeatother

\title
{The Kolmogorov backward equation for stochastic Burgers equations and for
stochastic 2D-Navier-Stokes equations}
\author{Martin Hutzenthaler$^{1}$, Robert Link$^{2}$, \\
\small{$^1$ Faculty of Mathematics, University of Duisburg-Essen,}\\
\small{Essen, Germany; e-mail: \texttt{martin.hutzenthaler}\textcircled{\texttt{a}}\texttt{uni-due.de}}\\
\small{$^2$ Faculty of Mathematics, University of Duisburg-Essen,}\\
\small{Essen, Germany; e-mail: \texttt{robert.link}\textcircled{\texttt{a}}\texttt{uni-due.de}}
}

\numberwithin{equation}{chapter}

\begin{document}

\begin{abstract}
	
In this book we establish under suitable assumptions
the uniqueness and existence of viscosity solutions
of Kolmogorov backward equations for stochastic partial differential equations (SPDEs).
In addition, we show that this solution is the semigroup of the corresponding SPDE.
This generalizes the Feynman-Kac formula to SPDEs and establishes a link
between solutions of Kolmogorov equations and solutions
of the corresponding SPDEs.
In contrast to the literature 
we only assume that the nonlinear part of the drift is Lipschitz continuous
on bounded sets (and not globally Lipschitz continuous) 
and we allow the diffusion coefficient to be degenerate and non-constant.
In the last part of this book we apply our results to stochastic Burgers equations
and to stochastic 2-D Navier-Stokes equations.
\end{abstract}
\maketitle

\tableofcontents

\chapter{Introduction}\label{sec:intro}
Kolmogorov proved in his seminal contribution \cite[Section 13]{kolmogoroff1931analytischen}
that continuous-time real-valued Markov processes solve the Kolmogorov backward equation under suitable assumptions.
It is by now common knowledge that
finite-dimensional linear parabolic partial differential equations (PDEs) without zero order terms
can be represented as semigroups of stochastic differential equations (SDEs) under suitable assumptions.
This connection between macroscopic models (PDEs) and microscopic models (SDEs) is very fruitful.
The case of infinite-dimensional linear PDEs is, however, not well understood.
To the best of our knowledge, there exists no stochastic representation in the literature
for solutions of infinite-dimensional Kolmogorov backward equations with
superlinearly growing drift functions
and unbounded initial functions.
The main contribution of this book is to establish the stochastic representation
\begin{equation}
	\forall(t,x)\in[0,1]\times H\colon u(t,x)=\E[\varphi(X_t^{x})]
\end{equation}
between the solution $u$ of an infinite-dimensional Kolmogorov backward equation
\begin{align} \begin{split}
\label{eq: Kol back eq}
	\forall(t,x)\in(0,1)\times H\colon 
		&\tfrac{ \partial }{ \partial t }
			u(t,x) -
			\big((D u)(t,x) \big)
			 (x''- F(x))
			-\tfrac 12 \,
			\text{tr}\!\left(
				B(x) \, (D^2 \, u)(t,x) \,  [ B(x) ] ^*
			\right) 
			= 0 
		\\ &
			u(0,x) = \varphi(x)
\end{split} \end{align}
and the solution $X$ of the corresponding SPDE in the case of superlinearly growing coefficients
under suitable assumptions.
 In particular, our results
include the benchmark cases of stochastic Burgers equations and stochastic 2D-Navier-Stokes equations.

A major challenge is the notion of solution for the Kolmogorov backward equation \eqref{eq: Kol back eq}.
We cannot expect \eqref{eq: Kol back eq} to have a classical ($C^2$-)solution.
Hairer, Hutzenthaler \& Jentzen \cite{HairerHutzenthalerJentzen2015} proved that even in finite dimensions
and for bounded and smooth functions $F$, $B$, and $\varphi$ the solution of \eqref{eq: Kol back eq}
might not be classical.
For this reason we will consider the notion of viscosity solutions.
A viscosity solution is informally speaking a function which in each point can be approximated from above with classical subsolutions and from below with classical supersolutions;
see Definition \ref{d:viscosity.solution} for a precise definition.
We note that Definition \ref{d:viscosity.solution} is a slight generalization of 
\cite[Definition 3.3]{Ishii1993} in order to allow for unbounded initial function;
see Remark \ref{rem: visc of ishii} below.
To the best of our knowledge there exists no result in the literature establishing existence and uniqueness of solutions of infinite-dimensional Kolmogorov backward equations
if $F$ is unbounded and not globally one-sided Lipschitz continuous 
and $B$ is non-constant.
To illustrate our main results, we consider in the following theorem the stochastic Burgers equation.
\begin{theorem}
\label{t: 1.1}
	Let
	$T \in (0,\infty)$,
	let $H=L^2((0,1);\R)$,
	let
		$ ( \Omega, \mathcal{F}, \P ) $
		be a probability space with a normal filtration\footnote{This means 
$
	\forall t \in [0,T] \colon
	\{
		A \subseteq \Omega \colon 
			\P(A)= 0
	\} 
	\subseteq
	\mathbb{F}_t = \bigcap_{s\in (t, T]} \mathbb{F}_s
$}
		$ ( \mathbb{F}_t )_{ t \in [0,T] } $, 
		let 
		$
			( W_t )_{ t \in [0,T] } 
		$ 
		be an $ \operatorname{Id}_{H} $-cylindrical 
		$ ( \mathbb{F}_t )_{ t \in [0,T] } $-Wiener
		process and for every $x \in H$ let
	$
		X^x \colon [0,T] \times \Omega \rightarrow H
	$
	be a mild solution of the 
	stochastic Burgers equation
	\begin{equation}
			X^x_t 
		= 
			x + \int_0^t \Delta X^x_s - \tfrac 12  ((X_s^x)^2)' \ud s 
			+\int_0^t X^x_s \cdot (-\Delta)^{-1} \dWs,
	\end{equation}
	let $\varphi \in \C(H,\R)$ have at most polynomial growth,
	and let $u \colon [0,T] \times H \to \R$
	be the function satisfying for all $t \in [0,T]$, $x \in H$ that
	$u(t,x)= \E[ \varphi(X^x_t) ]$.
	Then it holds that $u \in \C([0,T] \times H,\R)$
	and that $u$ is the unique viscosity solution (cf.\@ Definition \ref{d:viscosity.solution})
	of
		\begin{equation} 
		\label{eq: Kol PDE}
			\begin{split}
				\tfrac{ \partial }{ \partial t }
				u(t,x) -
				\int_0^1
					(x''(r)- \tfrac 12 (x^2)'(r)) 
					\big((D_{\H} u)(t,x) \big)(r)
				\ud r 
				-\tfrac 12
				\text{tr}\!\left(
					x \cdot (D^2_{\H} \, u)(t,x) \cdot x^* \, (-\Delta)^{-2} 
				\right)
				= 0
			\end{split}
		\end{equation}
		for $t \in (0,T)$, $x \in H$,
		relative to 
		$
			(
				(0,T) \times H \ni (t,u) 
				\to \nicefrac 12 \|u \|^2_{H_{1/2}} \in [0, \infty],
				\R \times \mathbb{H}, 
				\R \times \mathbb{H}_{1/2}
			)
		$
		which grows at most polynomially and satisfies 
		$u(0,\cdot) = \varphi(\cdot)$.
\end{theorem}
Theorem \ref{t: 1.1} follows directly from Corollary \ref{cor: Burgers equation}\footnote{with $B \leftarrow \id_{H_{-1}}$, $U \leftarrow H_{-1}$}
and Remark \ref{rem: sufficient condition for B}\footnote{with $s\leftarrow 1$, $b \leftarrow ([0,1] \times \R \ni (x,y) \to y \in \R)$}.

Stochastic representations of finite-dimensional PDEs are well known.
Theorem 4.3 in Pardoux \& Peng \cite{PardouxPeng1992} 
shows
that the semigroup of an SDE is the unique viscosity solution of the 
corresponding Kolmogorov equation 
if the drift and the diffusion
functions are globally Lipschitz continuous.
Furthermore, Theorem 4.14 and Theorem 4.16
in Hairer, Hutzenthaler \& Jentzen \cite{HairerHutzenthalerJentzen2015}
imply
that the semigroup is the unique viscosity solution 
satisfying a growth condition of the 
corresponding Kolmogorov equation if the drift and the diffusion
functions are locally Lipschitz continuous and there exists a Lyapunov
function.
Moreover, existence and uniqueness of viscosity solutions for 
finite-dimensional Kolmogorov equations have been established
e.g.\@ in
Theorem 8.2 in Crandall, Ishii \& Lions \cite{CrandallIshiiLions1992}
if the drift is globally one-sided Lipschitz continuous
and the diffusion is globally Lipschitz continuous
and for at most polynomially growing solutions in Theorem C.3.4 in Peng \cite{Peng2010}
if the drift is locally Lipschitz continuous and the diffusion is constant.

There are mainly two different approaches to show
existence and uniqueness of solutions for 
infinite-dimensional Kolmogorov equations.
The first one is using the notion of viscosity solution. 
In this approach 
with the exception of Gozzi, Sritharan \& \'{S}wi\c{e}ch \cite{GozziSritharanSwiech2005} 
(which assumes $B$ to be constant)
it is assumed that $F$ and $B$ are globally (one-sided) Lipschitz continuous 
which is a very strong tool.
In particular,
Theorem 2 and Theorem 3 in Lions \cite{Lions1989} assume that the drift and the diffusion function  are globally Lipschitz continuous
and there exists a classical unbounded subsolution,
Theorem 6.3 and Corollary 6.2 in Ishii \cite{Ishii1993}
establish existence and uniqueness of bounded viscosity solutions for 
infinite-dimensional Kolmogorov equations of SPDEs
if the nonlinearity and diffusion functions are globally Lipschitz continuous.
Moreover, Theorem 3.8 in Gozzi, Rouy \& \'{S}wi\c{e}ch \cite{GozziRouySwiech2000}
proves existence and uniqueness of  bounded viscosity solutions of 
elliptic
infinite-dimensional Kolmogorov equations of SPDEs, 
if $F$ is globally one-sided Lipschitz continuous with respect
to suitable norms,
if $B$ is globally Lipschitz continuous with respect
to suitable norms,
 and if
$F$ and $B$ have at most linear growth
and are uniformly continuous.
Theorem 6.3 in Gozzi, Sritharan \& \'{S}wi\c{e}ch \cite{GozziSritharanSwiech2005} 
shows the uniqueness and existence of at most polynomially growing
viscosity solutions of the Hamilton-Jacobi-Bellmann equations associated with the
stochastic 2-D Navier-Stokes equations
with periodic boundary conditions, 
if $B$ is constant, 
if $\varphi$ is
Lipschitz continuous on bounded sets with respect to suitable norms
and at most polynomially growing,
and if some other conditions are fulfilled.
Moreover, Theorem 3.1 in Bang \& Tran \cite{BangTran2006}
shows the uniqueness and existence of at most linear growing viscosity solutions of
infinite-dimensional Kolmogorov equations of SPDEs, 
if $F$ is globally one-sided Lipschitz continuous with respect
to suitable norms, 
if $B$ is globally Lipschitz continuous with respect
to suitable norms,
if
$F$ and $B$ have at most linear growth, and if $\varphi$ is bounded
and uniformly continuous.

The other approach is to study the equation \eqref{eq: Kol PDE} using the 
smoothing property of the semigroup.
These methods need
as a very strong tool
that $B$ is non-degenerate and often assume $B$ to be constant.
For example Theorem 5.3.5 in Cerrai \cite{Cerrai2001} proves that
the semigroup is the
unique classical solution of \eqref{eq: Kol PDE}
if $F$ is globally Lipschitz continuous, $B$ is non-degenerate and constant,
$\varphi$ is bounded,
and if suitable additional assumptions are fulfilled.
In addition, Da Prato \cite{DaPrato2004} studies elliptic Kolmogorov equations
and establishes in
Theorem 5.25 
(resp.\@ Theorem 6.14) 
that the semigroup is the 
unique strong solutions of the 
elliptic Kolmogorov equations of stochastic Burgers equations
(resp.\@ stochastic 2-D Navier-Stokes equations) 
with periodic boundary conditions
if $B$ satisfies $B = (-\Delta)^\gamma$, with $\gamma > \nicefrac 14$
(resp.\@ $B = (-\Delta)^\gamma$, with $\gamma > 1$).
Furthermore, Theorem 2.15 in Flandoli, Luo \& Ricci \cite{FlandoliLuoRicci2021}
shows uniqueness and existence of a suitable mild solution, 
if $F$ is bounded, 
$B$ is non-degenerate and constant, and if $\varphi$ is bounded.
Further result showing uniqueness and existence of solutions 
includes Theorem 9.27 in Da Prato \& Zabczyk \cite{DaPratoZabczyk1992}
for Kolmogorov backward equations of SPDEs if $F$ is bounded and
Lipschitz continuous and $B$ is constant, 
Theorem 3.6 in Gubinelli \& Perkowski
\cite{GubinelliPerkowski2020} for Kolmogorov backward equations of stochastic Burgers equations with additive space-time white noise,
Theorem 2.2 in R\"ockner \& Sobol \cite{RocknerSobol2006}  
for Kolmogorov forward equations of stochastic Burgers with additive trace-class noise
and Theorem 2.2 in R\"ockner \& Sobol \cite{RocknerSobol2007}  
for Kolmogorov forward equations of stochastic 2-D Navier-Stokes equations with multiplicative trace-class noise.

\section{Notation}
\label{sec: Notation}
The following notation 
is used throughout this book.
We denote by $(0, \infty]$, $[0, \infty]$, 
$\N$, and by $\N_0$ the sets satisfying that
$(0, \infty] = (0, \infty) \cup \{ \infty \}$, 
$[0, \infty] = [0, \infty) \cup \{ \infty \}$, 
$\N = \{1, 2, 3, \ldots \}$, 
and that $\N_0 = \N \cup \{ 0 \}$
and we denote by
$\sup (\emptyset)$,
$\inf (\emptyset)$, 
$\tfrac 00$,
$0^0$,
$0 \cdot \infty$,
$1^\infty$,
and by $\infty ^0$
the extended real number satisfying that
$\sup (\emptyset) = - \infty$,
$\inf (\emptyset) = \infty$,
$\tfrac 00 =0$,
$0^0 =1$,
$0 \cdot \infty = 0$, 
$1^\infty = 1$, and that
$\infty^0 =1$.
For all sets $A, B, C, D$ with
$D \subseteq C$ and every function
$ f \colon A \to B$ we denote 
by $\1_{D} \colon C \to \{0,1 \}$
the function satisfying for all $x \in C$ that
\begin{equation}
		\1_{D}(x)
	=
		\begin{cases}
			1 & \textrm{ if } x \in D \\
			0 & \textrm{ if } x \notin D
		\end{cases},
\end{equation}
by $f(A)$ the set satisfying that
$
		f(A)
	=
		\{f(x) \colon x \in A \}
$,
by $\mathcal{P}(A)$ 
the power set of $A$,
by $\mathbb{M}(A,B)$ 
the set of all mappings from $A$ to $B$,
and by
$f|_C \colon C \cap A \to B$ the function satisfying for all $x \in C \cap A$ 
that
$f|_C(x) = f(x)$.
Moreover, 
for real Hilbert spaces 
$\mathbb{H}_1 = (H_1, \langle \cdot, \cdot \rangle_{H_1}, \| \cdot \|_{H_1})$ and
$\mathbb{H}_2 = (H_2, \langle \cdot, \cdot \rangle_{H_2}, \| \cdot \|_{H_2})$ 
we denote by
$\langle \cdot , \cdot \rangle_{H_1 \times H_2} \colon H_1 \times H_2 \to \R$ 
the scalar product satisfying for all $(x_1, x_2)$, $(x_3,x_4) \in H_1 \times H_2$ that
$
		\langle (x_1 , x_2), (x_3 , x_4) \rangle_{H_1 \times H_2} 
	= 
		\langle x_1 , x_3 \rangle_{H_1} + \langle x_2 , x_4 \rangle_{H_2}
$, 
by $\| \cdot \|_{H_1 \times H_2} \colon H_1 \times H_2 \to [0, \infty)$ the norm 
satisfying for all $(x_1, x_2) \in H_1 \times H_2$ that
$
		\| (x_1, x_2) \|^2_{H_1 \times H_2} 
	=
		\|x_1\|^2_{H_1} + \|x_2\|^2_{H_2},
$
by $\mathbb{H}_1 \times \mathbb{H}_2$ the Hilbert space satisfying that
$
		\mathbb{H}_1 \times \mathbb{H}_2
	=
		(
			H_1 \times H_2, 
			\langle \cdot, \cdot \rangle_{H_1 \times H_2},
			\| \cdot \|_{H_1 \times H_2}
		),
$
by 
$
	\| \cdot \|_{L(\mathbb{H}_1, \mathbb{H}_2)} \colon L(\mathbb{H}_1, \mathbb{H}_2) 
		\to [0, \infty)
$ 
the norm satisfying for all 
$A \in L(\mathbb{H}_1, \mathbb{H}_2)$ that 
\begin{equation}
		\| A \|_{L(\mathbb{H}_1, \mathbb{H}_2)}
	= 
		\begin{cases}
			\sup_{x \in H_1 \backslash \{0\}} \frac {\| Ax \|_{H_2}}{\|x\|_{H_1}}, 
			& \textrm{ if } H_1 \neq \{0\}\\
			1, & \textrm{ if }H_1 = \{0\}, 
		\end{cases}
\end{equation}
by
$\mathbb{L}(\mathbb{H}_1, \mathbb{H}_2)$ 
the Banach space satisfying that
$
		\mathbb{L}(\mathbb{H}_1, \mathbb{H}_2)
	=
		(
			L(\mathbb{H}_1, \mathbb{H}_2), 
			\| \cdot \|_{L(\mathbb{H}_1, \mathbb{H}_2)}
		),
$ 
by
$
		\mathbb{H}_1' 
	= 
		\mathbb{L}(\mathbb{H}_1, \R)
$
its dual space,
and by 
$
		\mathbb{HS}(\mathbb{H}_1, \mathbb{H}_2)
	=
		(
			HS(\mathbb{H}_1, \mathbb{H}_2), 
			\langle \cdot, \cdot \rangle_{HS(\mathbb{H}_1, \mathbb{H}_2)},
			\| \cdot \|_{HS(\mathbb{H}_1, \mathbb{H}_2)}
		)
$
the Hilbert space of Hilbert-Schmidt operators from $H_1$ to $H_2$.
Furthermore, for 
Hilbert spaces 
$\mathbb{H}_1=(H_1, \langle \cdot, \cdot \rangle_{H_1}, \| \cdot \|_{H_1})$,
$\mathbb{H}_2=(H_2, \langle \cdot, \cdot \rangle_{H_2}, \| \cdot \|_{H_2})$,
a closed linear subspace $\mathbb{V} \subseteq \mathbb{H}_1$, 
the dual spaces 
$\mathbb{H}_1'=(H_1, \| \cdot \|_{H_1})$,
$
		(\mathbb{H}_1 \times \mathbb{H}_2)'
	= 
		((H_1 \times H_2)', \| \cdot \|_{(H_1 \times H_2)'})
$, and
$\mathbb{V}'= (V', \| \cdot \|_{V'})$,
$r \in \R$,
subsets $A, B \subseteq H_1$,
a dense subset $M \subsetneq H_1$
an open set $O \subseteq H_1$,
a function $f \colon O \to H_2$,
we denote by
$\Span_{\mathbb{H}_1}(A)$ the set satisfying that
$
		\Span_{\mathbb{H}_1}(A) 
	= 
		\{ 
			x \in H \colon 
			(
				\exists n \in \N, \exists \lambda_1, \ldots, \lambda_n \in \R,
				\exists x_1, \ldots, x_n \in A \colon
					x = \sum_{i=1}^n (\lambda_i x_i)
			)
		\}
$,
by
$\overline{A}_{\mathbb{H}_1} \subseteq H_1$ the set satisfying that
$
		\overline{A}_{\mathbb{H}_1}
	=
		\{
			x \in H_1 \colon 
				(\exists (x_n)_{n \in \N} \subseteq A 
					\colon \lim_{n \to \infty} \|x-x_n\|_{H_1} =0)
		\}
$,
by
$A+B$ the set satisfying that
$A+B = \{a+b \colon a \in A, b \in B\}$,
by $A-B$ the set satisfying that
$A-B = \{a-b \colon a \in A, b \in B\}$,
by $r A$ the set satisfying that
$r A = \{r a \colon a \in A\}$,
by
$
	\dist_{\mathbb{H}_1} \colon H_1 \times \mathcal{P}(H_1) \to [0, \infty]
$
the function satisfying for all $x \in H_1$ and all 
$A \in \mathcal{P}(H_1)$ that
$
	\dist_{\mathbb{H}_1}(x,A) 
= 
	\inf \{ \|x-y\|_{H_1} \colon y \in A\}
$,
by $\pi_{M}^{H_1} \colon H_1 \to H_1$ the operator satisfying that
$\pi_{M}^{H_1} = \id_{H_1}$,
by
$\pi_V^{H_1} \colon H_1 \to V$ the projection from $H_1$ to $V$,
by $\pi_{V'}^{H_1'} \colon H_1' \to H_1'$ the operator satisfying for
all $\xi \in H_1'$ that
$\pi_{V'}^{H_1'} \xi = \xi \pi_V^{H_1}$,
by 
$\pi^{(H_1 \times H_2)'}_1 \colon (H_1 \times H_2)' \to H_1'$,
the operator satisfying 
for all $\xi \in (H_1 \times H_2)'$ and all
$x \in H_1$ that
$
		(\pi^{(H_1 \times H_2)'}_1 \xi) (x)
	=
		 \xi (x, 0)
$, by
$\pi^{(H_1 \times H_2)'}_2 \colon (H_1 \times H_2)' \to H_2'$
the operator satisfying 
for all $\xi \in (H_1 \times H_2)'$ and all
$y \in H_2$ that
$
		(\pi^{(H_1 \times H_2)'}_2 \xi) (y)
	=
		\xi (0, y)
$,
by
$\C^{n}_{\mathbb{H}_1, \mathbb{H}_2}(A, H_2)$, $n \in \N_0$,
 the sets satisfying that for all $n \in \N_0$ it holds that
\begin{align}
\nonumber
			\C^{n}_{\mathbb{H}_1, \mathbb{H}_2}(O, H_2)
		=
			\{
				&f \in \mathbb{M}(O, H_2) \colon 
					\textrm{ $f$ is }n 
					\textrm{ times continuously Fr\'echet differentiable with } \\
					&\textrm{ respect to the }
						\| \cdot \|_{H_1} \textrm{-norm and the }
						\| \cdot \|_{H_2} \textrm{-norm}
			\},
\end{align} 
by $\C_{\mathbb{H}_1, \mathbb{H}_2}(O, H_2)$ the set satisfying that
$
		\C_{\mathbb{H}_1, \mathbb{H}_2}(O, H_2)
	=
		\C^{0}_{\mathbb{H}_1, \mathbb{H}_2}(O, H_2),
$ 
by
$\C^{n}_{\mathbb{H}_1}(O, \R)$, $n \in \N_0$, the sets satisfying that
for all $n \in \N_0$ it holds that
$
		\C^{n}_{\mathbb{H}_1}(O, \R)
	=
		\C^{n}_{\mathbb{H}_1, \R}(O, \R),
$
by $\C_{\mathbb{H}_1}(O, \R)$ the set satisfying that
$
		\C_{\mathbb{H}_1}(O, \R)
	=
		\C^{0}_{\mathbb{H}_1}(O, \R),
$ 
by 
$
	D_{\mathbb{H}_1, \mathbb{H}_2} 
$ 
the first Fr\'echet derivative 
with respect to the $\| \cdot \|_{H_1}$-norm
and the $\| \cdot \|_{H_2}$-norm, 
by 
$
	D^2_{\mathbb{H}_1, \mathbb{H}_2}
$
the second  Fr\'echet derivative with respect to the 
$\| \cdot \|_{H_1}$-norm and the $\| \cdot \|_{H_2}$-norm,
by $D_{\mathbb{H}_1}$ the operator satisfying that
$
		D_{\mathbb{H}_1}
	=
		D_{\mathbb{H}_1, \R}
$,
by $D^2_{\mathbb{H}_1}$ the operator satisfying that
$
		D^2_{\mathbb{H}_1}
	=
		D^2_{\mathbb{H}_1, \R}
$,
for an open set $U \subseteq \R$ 
and $u \in \C^1_{\R \times \H_1, \R}(U \times O,H_2)$ 
(resp.\@ $u \in \C^2_{\R \times \H_1, \R}(U \times O,H_2)$)
we denote
by 	
$
	D_{\mathbb{H}_1} u 
		\in \C_{\R \times \H_1, \mathbb{L}(\H_1,\R)}(U \times O, L(\H_1,\R))
$
(resp.\@
	$
		D^2_{\mathbb{H}_1} u 
			\in 
				\C_{\R \times \H_1, \mathbb{L}(\H_1, \mathbb{L}(\H_1,\R))}
					(U \times O,L(\H_1, \mathbb{L}(\H_1,\R)))
	$)
 the partial derivative satisfying for all $(t,x) \in O\times U$ and all $y \in H_1$ that
$
		(D_{\mathbb{H}_1} u) (t,x) (y) 
	= 
		(D_{\R \times \mathbb{H}_1} u)(t,x) (0, y)
$
(resp.\@ the partial derivative satisfying for all $(t,x) \in O\times U$ and all $y,z \in H_1$ that
$
		(D^2_{\mathbb{H}_1} u) (t,x) (y)(z) 
	= 
		(D^2_{\R \times \mathbb{H}_1} u)(t,x) (0, y)(0, z)
$),
and we call a set $B \subseteq H_1$ $\mathbb{H}_1$-bounded if 
$\sup \{ \|x\|_{H_1} \colon x \in B\} < \infty$ 
and $\mathbb{H}_1$-closed if 
$\overline{B}_{\mathbb{H}_1} = B$,
and we call the function $f$ $\mathbb{H}_2$-bounded if
$f(O)$ is $\mathbb{H}_2$-bounded.
Furthermore, we denote by $\Gamma \colon (0, \infty) \to (0, \infty)$
the Gamma function, that is, the function satisfying for all 
$x \in (0, \infty)$ that
$\Gamma(x) = \int _0^\infty t^{x-1} e^{-t} \ud t$.
For 
an $ \R $-vector space $ V $,
and a mapping $ \left\| \cdot \right\| \colon V \to [0,\infty] $
which satisfies for all 
$ 
	v, w \in V 
$,
$ \lambda \in \R  $
that
$( \| v \| = 0 \Leftrightarrow v = 0 )$,
$ \left\| \lambda v \right\|  =  |\lambda | \left\| v \right\| $,
and 
$ \left\| v + w \right\|  \leq \left\| v \right\| + \left\| w \right\| $
we call $ \left\| \cdot \right\| $ an extended norm on $ V $
and we call $ ( V , \left\| \cdot \right\| ) $ an extendedly normed 
$\R$-vector space.
For a metric space
$ ( M, d ) $,
an extendedly normed vector space
$ ( E, \left\| \cdot \right\| ) $,
a real number
$ r \in [0,1] $,
and a set
$ A \subseteq (0,\infty) $ 
we denote by
$
  \left| 
    \cdot  
  \right|_{
    \C^{ r, A }( M, \left\| \cdot \right\| )
  },
  \left| \cdot \right|_{
    \C^{ r }( M, \left\| \cdot \right\| )
  },
  \left\| \cdot \right\|_{
    \C( M, \left\| \cdot \right\| )
  }, $ $
  \left\| \cdot \right\|_{
    \C^{ r }( M, \left\| \cdot \right\| )
  } \colon
  \mathbb{M}( M, E ) \to [0, \infty] $
the mappings satisfying for all $ f \in \mathbb{M}( M, E ) $ that
\begin{align}
  \left| 
    f  
  \right|_{
    \C^{ r, A }( M, \left\| \cdot \right\| )
  }
&  
  =
  \sup\!\left(
    \left\{
      \tfrac{
        \left\| f(e_1) - f(e_2) \right\|
      }{
        \left| d(e_1, e_2) \right|^{ r }
      }
      \colon
      e_1, e_2 \in M,
      d(e_1,e_2) \in A
    \right\} \cup \left\{ 0 \right\}
  \right)
  \in [0,\infty],
\\
  \left| f \right|_{
    \C^{ r }( M, \left\| \cdot \right\| )
  }
&  
  =
  \left| f \right|_{
    \C^{ r, (0,\infty) }( M, \left\| \cdot \right\| )
  }
 \in [0,\infty],
\\
  \left\| f \right\|_{
    \C( M, \left\| \cdot \right\| )
  }
  & 
  =
  \sup\!\left(
    \left\{
      \left\| f(e) \right\|
      \colon 
      e\in M
    \right\}
    \cup
    \{ 0 \}
  \right)
  \in [0,\infty] ,
\\
  \left\| f \right\|_{
    \C^r( M, \left\| \cdot \right\| )
  }
  & 
  =
  \left\| f \right\|_{ \C( M, \left\| \cdot \right\| ) }
  +
  \left| f \right|_{ \C^r( M, \left\| \cdot \right\| ) }
  \in [0,\infty]
\end{align}
and we denote by 
$ 
  \C^r( M, \left\| \cdot \right\| ) 
$ 
the set satisfying that
\begin{equation}
  \C^r( M, \left\| \cdot \right\| ) 
=
  \{ 
    f \in 
    \C(M,E)
    \colon
    \| f \|_{ \C^r( M, \left\| \cdot \right\| ) } < \infty
  \}.
\end{equation}
For a measure space $ ( \Omega, \mathbb{F}, \mu ) $,
an extendedly normed vector space $ ( V , \left\| \cdot \right\| ) $,
and real numbers $ p \in [ 0, \infty ) $, $ q \in ( 0, \infty ) $
we denote by $ L^0( \mu; \left\| \cdot \right\| ) $
the set given by
\begin{equation}
		L^0( \mu; \left\| \cdot \right\| )
	= 
		\bigl\{  
			f \in \mathbb{M}( \Omega, V ) \colon
				f \text{ is } ( \mathbb{F}, \left\| \cdot \right\| ) \text{-strongly measurable} 
		\bigr\},
\end{equation}
we denote by
$ \left\| \cdot \right\|_{ L^q( \mu; \left\| \cdot \right\| ) } \colon
    L^0( \mu; \left\| \cdot \right\| ) \to [0,\infty] $
the mapping which satisfies
for all $ f \in L^0( \mu; \left\| \cdot \right\| ) $ that
\begin{equation}
\left\| f \right\|_{ L^q( \mu; \left\| \cdot \right\| ) }
= \left[  \int_{ \Omega }
                \left\| f( \omega ) \right\|^q   \mu( \ud \omega )
   \right]^{ \nicefrac{ 1 }{ q } }
   \in [0,\infty],
\end{equation}
we denote by
$ L^q( \mu; \left\| \cdot \right\| ) $
the set given by
\begin{equation}
L^q( \mu; \left\| \cdot \right\| ) 
= \bigl\{  f \in L^0( \mu; \left\| \cdot \right\| ) \colon
    \left\| f \right\|_{ L^q( \mu; \left\| \cdot \right\| ) } 
    < \infty \bigr\}.
\end{equation}
In addition, for 
a real Hilbert space $\mathbb{H}=(H, \langle \cdot, \cdot \rangle_H, \| \cdot \|_H),$
a subset $A \subseteq H$,
and a function $f \colon A \to \R \cup \{-\infty, \infty\}$
we call an element $x \in A$ a strict $\mathbb{H}$-maximum of $f$, 
if $x$ is a global maximum of $f$ 
and if for all $(x_n)_{n \in \N} \subseteq A$ with
$\lim_{n \to \infty} f(x_n) = f(x)$ it holds that 
$\lim_{n \to \infty} \| x_n - x\|_H = 0$.
For the rest of this section fix
real Hilbert spaces 
$\mathbb{H}=(H, \langle \cdot, \cdot \rangle_H, \| \cdot \|_H)$
and
$\mathbb{X} = (X, \langle \cdot, \cdot \rangle_X, \| \cdot \|_X)$  
satisfying that $X \subseteq H$, 
$\mathbb{X}$ is embedded continuously in $\mathbb{H}$, and
$X$ is dense in $H$ with respect to the $\| \cdot \|_H$-norm
and fix the dual spaces
$\mathbb{H}'=(H', \| \cdot \|_{H'})$ 
and $\mathbb{X}' = (X', \| \cdot \|_{X'})$.
Then
we denote
by
$\langle \cdot, \cdot \rangle_{H, H'} \colon H \times H' \to \R$ and by
$\langle \cdot, \cdot \rangle_{H', H} \colon H' \times H \to \R$
the functions satisfying
for all $x \in H$ and $y \in H'$ that
$\langle y, x \rangle_{H', H} =\langle x, y \rangle_{H, H'} = y(x)$, by 
$
  \mathbb{S}_{\mathbb{H}, \mathbb{H}'}
  =
  \{ 
    B \in L (\mathbb{H}, \mathbb{H}')
			\colon (
				\forall x, y \in H 
					\colon \langle x, B y \rangle_{H, H'}= \langle y, Bx \rangle_{H,H'}
			)
  \}
$
the set of all
bounded symmetric linear operators from $H$ into $H'$,
and by
$
	\otimes \colon H' \times H' \to L(\mathbb{H}, \mathbb{H}')
$
the operator satisfying for all $x \in H$ 
and all $y_1$, $y_2 \in H'$ that
\begin{equation}
		(y_1 \otimes y_2)(x) 
	= 
		\langle x, y_1 \rangle_{H, H'} \cdot  y_2.
\end{equation}
Moreover, we denote by $I_{\mathbb{H}} \colon H \to H'$
the 
function from $H$ to $H'$ satisfying
for all $x$, $y \in H $ that
\begin{equation}
		\langle I_{\mathbb{H}} x, y \rangle_{H', H}
	=
		\langle x, y \rangle_{H},
\end{equation}
by 
$
	\cdot \colon \mathbb{S}_{\mathbb{H}, \mathbb{H}'} \times  \mathbb{S}_{\mathbb{H}, \mathbb{H}'}
		\to  \mathbb{S}_{\mathbb{H}, \mathbb{H}'}
$
the multiplication operator satisfying for all 
$A$, $B \in \mathbb{S}_{\mathbb{H}, \mathbb{H}'}$ that
$A \cdot B = A I^{-1}_\mathbb{H} B$,
and by 
$E_{\mathbb{X}', \mathbb{H}'} \colon D(E_{\mathbb{X}', \mathbb{H}'}) \subseteq X' \to H'$
the extension operator satisfying that 
\begin{equation}
		D(E_{\mathbb{X}', \mathbb{H}'}) 
	= 
		\{ 
			x \in X' \colon 
				(
					\exists C \in (0, \infty) \colon \forall y \in X \colon 
						\langle x, y \rangle_{X', X} \leq C \|y\|_H 
				) 
		\}
\end{equation}
and satisfying that for all $x \in D(E_{\mathbb{X}', \mathbb{H}'})$ and 
all $y \in X$ it holds that
\begin{equation}
		\langle E_{\mathbb{X}', \mathbb{H}'}(x), y \rangle_{H',H} 
	=  
		\langle x, y \rangle_{X',X}.
\end{equation}
Note that $X \subseteq H$ implies that
for all $y \in H'$ it holds that
$y |_X \in X'$
and thus we have for all $A \in \mathbb{S}_{\mathbb{H}, \mathbb{H}'}$ that
$ (A |_X) |_X \in \mathbb{S}_{\mathbb{X}, \mathbb{X}'}$
and it holds
for all $y \in H'$ and all
$x \in X \subseteq H$ that
$\langle x, y |_X \rangle_{X,X'} = \langle x, y \rangle_{H,H'}$.
Furthermore, for all
$ B, C \in \mathbb{S}_{\mathbb{X}, \mathbb{X}'} $
we write
$ B \leq C $ 
in the following
if for all 
$ x \in X $
it holds that 
$
  \langle x, B x \rangle_{X, X'}
  \leq
  \langle x, C x \rangle_{X, X'}
$.
Moreover, 
for
a set $W \subseteq X$
we call 
a function
$
  F
  \colon   
	W
  \times
  \R
  \times
  H'
  \times
  \mathbb{S}_{\mathbb{X}, \mathbb{X}'}
  \to\R
$
\emph{degenerate elliptic} 
(see, e.g., (F0) in
Ishii~\cite{Ishii1993})
if
for all 
$ x \in  W $,
$ r \in \R $,
$ p \in H' $
and all 
$ 
  A, B \in \mathbb{S}_{\mathbb{X}, \mathbb{X}'}
$
with
$ A \geq B $ we have
$
  F(x, r, p, A) \leq F(x, r, p, B)
$.
In addition, for a set $O \subseteq H$,
a dense subset $W \subseteq O$, 
and a function
$u \colon O \to \R \cup \{ -\infty, \infty \}$ 
we denote by $\overline{u}^W_\mathbb{H} \colon O \to \R \cup \{ -\infty, \infty \}$ 
and by 
$\underline{u}^W_\mathbb{H} \colon O \to \R \cup \{ -\infty, \infty \}$
the functions satisfying for
all $x \in O$ that
\begin{equation} 
		\overline{u}^W_\mathbb{H}(x) 
	= 
		\lim_{\eps \downarrow 0} \sup \{ u(y) \colon y \in W,~ \|x-y\|_H \leq \eps \}
\end{equation}
and that
\begin{equation}
		\underline{u}^W_\mathbb{H}(x) 
	= 
		\lim_{\eps \downarrow 0} \inf \{ u(y) \colon y \in W,~ \|x-y\|_H \leq \eps \}.
\end{equation}

\chapter{Main properties of viscosity solutions}
\label{sec:Kolmogorov.equations}

The classical notion of viscosity solutions introduced in Crandall \& Lions
\cite{CrandallLions1981} 
(see also Crandall, Ishii \& Lions~\cite{CrandallIshiiLions1992})
are not suited to deal with discontinuous differential equation.
For this reason Ishii introduced in \cite{Ishii1993} a more general notion of viscosity 
solutions which overcomes this problem.
In Definition \ref{d:viscosity.solution} below we slightly generalize the notion of 
viscosity solution of Ishii
\cite[Definition 3.3]{Ishii1993} to allow for unbounded initial function;
see Remark \ref{rem: visc of ishii} below.

In the first part of this chapter we
define
viscosity solutions
and we show that many important properties of viscosity solutions
also hold for our more general notion of viscosity solutions.
In the next two parts we prove 
two important lemmas needed in chapter 4.
In particular we prove in the second part
that viscosity solutions are stable under limits
(Lemma \ref{l:limits.of.viscosity.solutions}).
And in the final part of this chapter we prove a chain rule for semijets 
and a lifting result for viscosity solutions 
(Lemma \ref{prop: u proj is viscosity solution}).
This will allow us in chapter 4 to treat finite-dimensional viscosity solutions
as viscosity solutions over an infinite-dimensional Hilbert space.
This chapter is based on Section 4.2 and Section 4.3 in 
Hairer, Hutzenthaler \& Jentzen \cite{HairerHutzenthalerJentzen2015} and
Section 2 in Ishii \cite{Ishii1993}.
Throughout this book we use the notation from Section \ref{ssec: Notation} 
below.
\section{Semijets}
\label{ssec: Notation}
First we recall 
(a minor generalization since we allow 
for all $x \in O$ that $u(x) \in \{-\infty, \infty\}$ of)
the definition of the \emph{semijets}
(see, e.g., Section~2 in 
Ishii \cite{Ishii1993}).
\begin{definition}[Second-order semijets]
\label{d:semijets}
Let $\mathbb{H}=(H, \langle \cdot, \cdot \rangle_H, \| \cdot \|_H)$ 
be a  real Hilbert space
and denote by
$H'$ the set satisfying that
$H' = L(\mathbb{H}, \R)$
and by $\mathbb{H}'$ the Hilbert space satisfying that
$\mathbb{H}'= \mathbb{L}(\mathbb{H}, \R)$,
let $ O \subseteq H $ be an open set
and let $ u \colon O \to \R \cup \{-\infty, \infty\} $ be a function.
Then we denote by
$ 
  J^2_{ \mathbb{H}, + } u 
$,
$
  J^2_{ \mathbb{H}, - } u 
$,
$ 
  \hat{J}^2_{ \mathbb{H}, + } u 
$,
$
  \hat{J}^2_{ \mathbb{H}, - } u 
	\colon O \to 
  \mathcal{P}\!\left( H' \times \mathbb{S}_{\mathbb{H}, \mathbb{H}'} \right) 
$
the functions satisfying
for all $ x \in O $ with $u(x) \notin \R$ that
$ 
		\big( J^2_{ \mathbb{H}, + } u \big)( x ) 
	= \big( \hat{J}^2_{ \mathbb{H}, + } u \big)( x )
	=	\big( J^2_{ \mathbb{H}, - } u \big)( x )
	= \big( \hat{J}^2_{ \mathbb{H}, - } u \big)( x )
	= \emptyset
$
and for all $ x \in O $ with $u(x) \in \R$ that 
\begin{equation}
\begin{split}
  \big(
    J^2_{ \mathbb{H}, + } u
  \big)( x )
& =
  \left\{
    ( p, A )
    \in H' \times \mathbb{S}_{\mathbb{H}, \mathbb{H}'}
    \colon
    \limsup_{
      O
      \ni y \to x
    }
    \left(
    \tfrac{
      u( y ) - u( x )
      -
      \left< p, y - x \right>_{H', H}
      -
      \frac{ 1 }{ 2 }
      \left< y - x, A ( y - x ) \right>_{H, H'}
    }{
      \left\| x - y \right\|_H^2
    }
    \right)
    \leq 0
  \right\}
  ,
\\ 
  \big( 
    \hat{J}^2_{ \mathbb{H}, + } u
  \big)( x )
& =
  \left\{
    (p, A)
    \in H' \times \mathbb{S}_{\mathbb{H}, \mathbb{H}'}
    \colon
    \left(
    \!\!
    \begin{array}{c}
      \exists
      \,
      ( x_n, p_n, A_n )_{ n \in \N } 
      \subseteq
      O \times H' \times \mathbb{S}_{\mathbb{H}, \mathbb{H}'} 
      \colon 
		\\
      \left(
        \forall \, n \in \N \colon
        (p_n, A_n) \in 
        ( J^2_{ \mathbb{H}, + } u)( x_n )
      \right) 
      \text{and } 
		\\
				\lim_{ n \to \infty } (
					\|x-x_n\|_H
					+ |u(x) - u(x_n)|  \quad\\
					\qquad + \|p- p_n\|_{H'} 
					+ \|A-A_n\|_{L(\mathbb{H}, \mathbb{H}')} 
				)
     = 0
    \end{array}
    \!\!
    \right)
  \right\} ,
\\
  \big(
    J^2_{ \mathbb{H}, - } u
  \big)( x )
& =
  \left\{
    ( p, A )
    \in H' \times \mathbb{S}_{ \mathbb{H}, \mathbb{H}'} 
    \colon
    \liminf_{
      O
      \ni y \to x
    }
    \left(
    \tfrac{
      u( y ) - u( x )
      -
      \left< p, y - x \right>_{H', H}
      -
      \frac{ 1 }{ 2 }
      \left< y - x, A ( y - x ) \right>_{H, H'}
    }{
      \left\| y - x \right\|_H^2
    }
    \right)
    \geq 0
  \right\}, 
\\ 
  \big( 
    \hat{J}^2_{ \mathbb{H}, - } u
  \big)( x )
& =
  \left\{
    (p, A)
    \in H' \times \mathbb{S}_{\mathbb{H}, \mathbb{H}'}
    \colon
    \left(
    \!\!
    \begin{array}{c}
      \exists
      \,
      ( x_n, p_n, A_n )_{ n \in \N } 
      \subseteq
      O \times H' \times \mathbb{S}_{\mathbb{H}, \mathbb{H}'}
      \colon
		\\
      \left(
        \forall \, n \in \N \colon
        (p_n, A_n) \in 
        ( J^2_{ \mathbb{H}, - } u)( x_n )
      \right)
      \text{and }
		\\
				\lim_{ n \to \infty } (
					\|x-x_n\|_H
					+ |u(x) - u(x_n)|  \quad\\
					\qquad + \|p- p_n\|_{H'} 
					+ \|A-A_n\|_{L(\mathbb{H}, \mathbb{H}')} 
				)
			= 0
    \end{array}
    \!\!
    \right)
  \right\}.
\end{split}
\end{equation}
\end{definition}
\begin{remark}
	Note that $J^2_{ \mathbb{H}, +} $ and $\hat{J}^2_{ \mathbb{H}, + } $
	are monotone in the second argument and
	that $\hat{J}^2_{ \mathbb{H}, +} $
	is closed under limits. The corresponding statements holds for
	$J^2_{ \mathbb{H}, -} $ resp.\@ $\hat{J}^2_{ \mathbb{H}, -} $.
	Moreover, for every function $u \colon O \to \R \cup \{-\infty, \infty\}$
	it holds that
	$
			J^2_{ \mathbb{H}, - } u
		=
			- J^2_{ \mathbb{H}, + } (-u)
	$
	and that
	$
			\hat{J}^2_{ \mathbb{H}, - } u
		=
			- \hat{J}^2_{ \mathbb{H}, + } (-u).
	$
\end{remark}
\section{Setting}
\label{ssec: Setting H X without t}
Throughout this chapter the following setting is frequently used.
Let
$\mathbb{H} = (H, \langle \cdot, \cdot \rangle_H, \| \cdot \|_H)$ and 
$\mathbb{X} = (X, \langle \cdot, \cdot \rangle_X, \| \cdot \|_X)$  
be real separable Hilbert spaces with the property that 
$X \subseteq H$, that $\mathbb{X}$ is embedded continuously in $\mathbb{H}$, and that
$X$ is dense in $H$ with respect to the $\| \cdot \|_H$-norm,
let $\mathbb{H}'=(H', \| \cdot \|_{H'})$ and
$\mathbb{X}'=(X', \| \cdot \|_{X'})$ be the dual spaces,
let $O \subseteq H$ be an open set,
let $W \subseteq O$ be a subset satisfying that
$W$
is dense in $O$ with respect to the $\| \cdot \|_H$-norm,
let
$h \colon O \to \R \cup \{ \infty \}$
be a
with respect to the $\| \cdot \|_H$-norm lower semicontinuous 
function
with the property that 
$
	h |_{X} \in \C_{\mathbb{X}}^2(X \cap O,\R)
$
and that
$
W=
	\{ 
		y \in X \cap O \colon (D_{\mathbb{X}} (h|_X)) (y) \in D(E_{\mathbb{X}', \mathbb{H}'}), 
		~ (J^{2}_{\mathbb{H}, -}h) (y) \neq \emptyset 
	\}
$,
let $\mathcal U$ be the set satisfying that
$\mathcal U = O \times \R \times H' \times \mathbb{S}_{\mathbb{H}, \mathbb{H}'}$
and let $\mathcal W$ be the set satisfying that
$\mathcal W = W \times \R \times H' \times \mathbb{S}_{\mathbb{H}, \mathbb{H}'}$.
\section{Notation}
Before we now give the definition of a viscosity solution let us introduce
some additional notation.
\begin{definition}
Assume the Setting in Section \ref{ssec: Setting H X without t},
for functions
$F \colon W \times \R \times H' \times \mathbb{S}_{\mathbb{X}, \mathbb{X}'} \to \R$
and
$u \colon O \to \R \cup \{-\infty, \infty\}$, and
for
$\delta \in (0, \infty)$,
we will denote by 
$
	u_{\mathbb{H},\delta,h}^{+, W}, ~u_{\mathbb{H},\delta,h}^{-, W}
		\colon O \to \R \cup \{-\infty, \infty\}
$
the functions satisfying for all $x \in O$ that
\begin{equation}
		u_{\mathbb{H},\delta,h}^{+, W} (x)
	= 
		\underline{(u + \delta h)}_{\mathbb{H}}^{W}(x)
\end{equation}
and that 
\begin{equation}
		u_{\mathbb{H},\delta,h}^{-, W}(x)
	= 
		\overline{(u - \delta h)}_{\mathbb{H}}^{W}(x),
\end{equation}
by 
$
	F_{\mathbb{H},\mathbb{X}, \delta, h}^+, ~F_{\mathbb{H},\mathbb{X}, \delta, h}^- \colon 
		\mathcal{W} \to \R,
$
the functions 
satisfying for all 
$
	(x,r,p,B) \in 
		\mathcal{W}
$ that
\begin{equation}
\label{eq: def F+ delta h}
		F_{\mathbb{H},\mathbb{X}, \delta, h}^+(x,r,p,B)
	= 
		F\big(
			x,
			r+ \delta h(x), 
			p+ \delta E_{\mathbb{X}', \mathbb{H}'}( (D_{\mathbb{X}}(h|_X))(x)), 
			(B|_X) |_X+ \delta (D_{\mathbb{X}}^2(h|_X))(x)
		\big)
\end{equation}
and that
\begin{equation}
		F_{\mathbb{H},\mathbb{X}, \delta, h}^-(x,r,p,B)
	= 
		F\big(
			x,
			r- \delta h(x), 
			p- \delta E_{\mathbb{X}', \mathbb{H}'}( (D_{\mathbb{X}}(h|_X))(x)), 
			(B|_X) |_X- \delta (D_{\mathbb{X}}^2(h|_X))(x)
		\big),
\end{equation}
by 
$
	d^{+, W}_{\mathbb{H}, \delta, h, u},
	~d^{-, W}_{\mathbb{H}, \delta, h, u} \colon 
		\mathcal{U}^2
			\to [0, \infty  ]
$ 
the functions
satisfying for all 
$
		(\xi, \eta)
	= 
		((x,r,p,B),(y,s,q,C)) 
$
$
			\in \mathcal{U}^2
$ 
that
\begin{equation}
\label{eq:def of d delta without t}
			d^{+, W}_{\mathbb{H}, \delta, h, u}(\xi, \eta) 
		=
			\begin{cases}
				&\|x-y\|_H \vee |r-s| \vee \|p-q\|_{H'} 
				\vee \| B-C \|_{L(\mathbb{H},\mathbb{H}')}
				\vee |u^{+, W}_{\mathbb{H}, \delta, h}(x)-u^{+, W}_{\mathbb{H}, \delta, h}(y)| \\
				& \qquad \qquad\qquad \qquad \qquad \qquad\qquad \qquad\qquad \quad
				\text{ if }
					u^{+, W}_{\mathbb{H}, \delta, h}(x), \,
					u^{+, W}_{\mathbb{H}, \delta, h}(y) \in \R \\
				&\infty  \qquad \qquad\qquad \qquad \qquad\qquad \qquad\qquad \quad \quad
				\text{ if } 
					u^{+, W}_{\mathbb{H}, \delta, h}(x) \textrm{ or } 
					u^{+, W}_{\mathbb{H}, \delta, h}(y) \notin \R
			\end{cases} 
\end{equation}
	and that
\begin{equation}
\label{eq:def of d delta - without t}
	\begin{split}
			&d^{-, W}_{\mathbb{H}, \delta, h, u}(\xi, \eta)
		=
			\begin{cases}
				&\|x-y\|_H \vee |r-s| \vee \|p-q\|_{H'} 
				\vee \| B-C \|_{L(\mathbb{H}, \mathbb{H}')} 
				\vee |u^{-, W}_{\mathbb{H}, \delta, h}(x)-u^{-, W}_{\mathbb{H}, \delta, h}(y)|\\
				& \qquad \qquad\qquad \qquad \qquad\qquad \qquad\qquad \qquad \quad
				\text{if } 
					u^{-, W}_{\mathbb{H}, \delta, h}(x), \,
					u^{-, W}_{\mathbb{H}, \delta, h}(y) \in \R \\
				&\infty \qquad \qquad\qquad \qquad \qquad\qquad \qquad\qquad \quad \quad
				\text{if } 
					u^{-, W}_{\mathbb{H}, \delta, h}(x) \textrm{ or }  
					u^{-, W}_{\mathbb{H}, \delta, h}(y) \notin \R
			\end{cases},
	\end{split}
\end{equation}
by 
$
	F_{\mathbb{H}, \mathbb{X}, \delta, h, u}^+,
	~F_{\mathbb{H}, \mathbb{X}, \delta, h, u}^- \colon 
		\mathcal U \to \R \cup \{-\infty, \infty \}
$ 
the functions satisfying for all $\xi \in \mathcal U$ that
\begin{equation}
\label{eq: def F+ delta h u}
		F_{\mathbb{H}, \mathbb{X}, \delta, h, u}^{+}(\xi)
	= 
		\lim_{\eps \downarrow 0} \inf 
			\big \{ 
				F_{\mathbb{H}, \mathbb{X}, \delta, h}^{+}(\eta) \colon 
				\eta \in \mathcal{W},
				~d^{-,W}_{\mathbb{H}, \delta, h, u}(\xi, \eta) \leq \eps 
			\big \}
\end{equation}
and that
\begin{equation}	
\label{eq: def F- delta h u}
		F_{\mathbb{H}, \mathbb{X}, \delta, h, u}^-(\xi)
	= 
		\lim_{\eps \downarrow 0} \sup 
			\big\{ 
				F^-_{\mathbb{H}, \mathbb{X}, \delta, h}(\eta) \colon 
				\eta \in \mathcal{W},
				~d^{+,W}_{\mathbb{H}, \delta, h, u}(\xi, \eta) \leq \eps 
			\big\}.
\end{equation}
\end{definition}
\begin{remark}
\label{rem: F hat is semicontinuous}
	By construction, $F^{+}_{\mathbb{H}, \mathbb{X}, \delta, h, u}$ 
	(resp.\@ $F^{-}_{\mathbb{H}, \mathbb{X}, \delta, h, u}$)
	is lower (resp.\@ upper) semicontinuous with respect to the
	disctance function $d^{-,W}_{\mathbb{H}, \delta, h, u}$ (resp.\@ $d^{+,W}_{\mathbb{H}, \delta, h, u}$).
\end{remark}
\section{Definition and basic properties of viscosity solutions}
\label{ssec:Definition of viscosity solutions}
With this notation we give the definition of viscosity solutions.
\begin{definition}[Viscosity solution]
\label{d:viscosity.solution}
	Assume the setting in Section \ref{ssec: Setting H X without t} and
	let 
	$
		F \colon 
			W  
			\times \R \times H' \times \mathbb{S}_{\mathbb{X}, \mathbb{X}'} \to \R
	$ 
	be a degenerate elliptic function.
  A function $ u \colon O \to \R \cup \{-\infty \}$
  is
  said to be
  a \emph{viscosity subsolution} of $F=0$ relative to $(h, \mathbb{H}, \mathbb{X})$
  (or, equivalently, a viscosity solution 
  of $ F \leq 0 $ relative to $(h, \mathbb{H}, \mathbb{X})$ ) 
	if $ u $ is locally bounded from above and if
  it holds 
  for all $ x \in O $,
	$\delta \in (0, \infty)$,
	and all
  $ 
    \phi \in \C_{\mathbb{H}}^2(O,\R)
  $ 
  with 
  $ \phi \geq u_{\mathbb{H}, \delta, h}^{-, W} $
  and
  $ \phi(x) = u_{\mathbb{H}, \delta, h}^{-, W}(x)$
  that
  \begin{equation}
			F_{\mathbb{H}, \mathbb{X}, \delta, h, u}^+\big(
				x, \phi(x), (D_\mathbb{H} \, \phi)(x), 
				(D_\mathbb{H}^2 \, \phi)(x) 
			\big)
			\leq 0.
  \end{equation}
  Similarly,
  a
  function $ u \colon O \to \R \cup \{\infty \} $
  is said to be
  a \emph{viscosity supersolution} 
  of $ F = 0 $ relative to $(h, \mathbb{H}, \mathbb{X})$ (or, equivalently, 
  a viscosity solution 
  of $F \geq 0$ relative to $(h, \mathbb{H}, \mathbb{X})$ )
  if
  $ u $ is locally bounded from below 
  and if
  it holds 
  for all $ x \in O $
  and all
  $ 
    \phi \in \C_{\mathbb{H}}^2(O,\R)
  $ 
  with 
  $ \phi \leq u_{\mathbb{H}, \delta, h}^{+, W} $
  and
  $ \phi(x) = u_{\mathbb{H}, \delta, h}^{+, W}(x) $
  that
  \begin{equation}  \label{eq:Fleq0}
    F^-_{\mathbb{H}, \mathbb{X}, \delta, h, u}\big(
      x, \phi(x), (D_\mathbb{H} \, \phi)(x), 
      (D_\mathbb{H}^2 \, \phi)(x) 
    \big)
    \geq 0.
  \end{equation}
  Finally, a 
  function $ u \colon O\to\R $ 
  is said to be
  a \emph{viscosity solution} of $F=0$ relative to $(h, \mathbb{H}, \mathbb{X})$
  if $u$ is both a viscosity subsolution
  and a viscosity supersolution of 
  $F=0$ relative to $(h, \mathbb{H}, \mathbb{X})$.
	In addition, we call $u\colon O \to \R$ a classical viscosity solution 
	(resp.\@ subsolution, supersolution)
	if $u \in \C_{\mathbb{H}}(O,\R)$ and if $u$ is a viscosity solution 
	(resp.\@ subsolution, supersolution) relative to 
	$(x \in H \to 0 \in \R, \mathbb{H}, \mathbb{H})$.
\end{definition}

\begin{remark}
\label{rem: visc of ishii}
	With the additional assumption that $h$ is convex and that
	$\{x \in H \colon  ( J^2_{\mathbb{H}, - } )( x ) $ $\neq \emptyset\} \subseteq X$
	we get the Definition 2.3 in Ishii~\cite{Ishii1993}.
\end{remark}
\begin{remark}
	For continuous $F$ and with
	$X \leftarrow H$ and $h \leftarrow 0$ 
	we get the "classical" definition of viscosity solution
	(see, e.g., Section~2 in Crandall, 
	Ishii 
	$\&$  Lions~\cite{CrandallIshiiLions1992}
	and also Definition~1.2
	in Appendix~C
	in Peng~\cite{Peng2010}).
\end{remark}
The following elementary lemma for viscosity solutions slightly
generalizes, e.g., Lemma 4.2 in in Hairer, Hutzenthaler 
	\& Jentzen \cite{HairerHutzenthalerJentzen2015}. The proof is straight-forward
and therefore omitted.

\begin{lemma}[Sign changes of viscosity solutions]
\label{lem:sign_changes}
	Assume the setting in Section \ref{ssec: Setting H X without t},
	let 
	$
		F \colon 
			W \times \R \times H' \times \mathbb{S}_{\mathbb{X}, \mathbb{X}'} \to \R
	$ 
	be a degenerate elliptic function, 
	and let
$ u \colon O \to \R $ be a viscosity solution of
$ F \geq 0 $ relative to $(h, \mathbb{H}, \mathbb{X})$.
Then the function 
$ 
  \tilde{F} \colon 
		W \times \R \times H' \times \mathbb{S}_{\mathbb{X}, \mathbb{X}'} \to \R 
$
satisfying for all 
$ 
	( x, r, p, A ) \in 
		W \times \R \times H' \times \mathbb{S}_{\mathbb{X}, \mathbb{X}'} 
$
that
$
  \tilde{F}( x, r, p, A ) 
  = - F( x, - r , - p, - A )
$
is degenerate elliptic and 
the function 
$
  O \ni x \mapsto - u(x) \in \R
$
is a viscosity solution of $ \tilde{F} \leq 0 $ relative to $(h, \mathbb{H}, \mathbb{X})$.
The corresponding statement holds for viscosity solutions
of $ F \leq 0 $ and $ F = 0 $ relative to $(h, \mathbb{H}, \mathbb{X})$
respectively.
\end{lemma}

In the next lemma we give an alternative characterization of the semijets via test functions
	(see also Remark~2.3 in 
	Crandall, Ishii \& Lions~\cite{CrandallIshiiLions1992}
	and Lemma 4.4 in Hairer, Hutzenthaler 
	\& Jentzen \cite{HairerHutzenthalerJentzen2015}).
More precisely, the upper semijet $(J^2_{\mathbb{H},+}u)(x)$
of a function $u\colon O\to\R\cup\{-\infty,\infty\}$ at a point $x\in O$
is the set of tuples of first and second order derivatives at $x$
of $C^2$-functions lying above $u$ and touching $u$ in $x$.

\begin{lemma}[Properties of semijets]
\label{lem:semijets}
Let $\mathbb{H} = (H, \langle \cdot, \cdot \rangle_H, \| \cdot \|_H)$ be a real
Hilbert space and 
let $\mathbb{H'} = (H', \| \cdot \|_{H'})$ be its dual space,
let $ O \subseteq H $
be an open set, and let 
$ 
  u \colon O \to \R \cup \{ -\infty, \infty \}
$
be locally bounded from above (resp.\@ below).
Then we have for all $ x \in O $ with $u(x) \in \R$ that
\begin{equation}
\begin{split}	
  ( J^2_{\mathbb{H}, + } u)( x )
& =
  \Big\{ 
    \big( 
      (D_\mathbb{H} \, \phi)(x),
      ( D_\mathbb{H}^2 \, \phi)(x)
    \big)
    \in 
    H' \times \mathbb{S}_{\mathbb{H}, \mathbb{H}'} \\
	& \qquad \qquad
    \colon
    \big(
      \phi \in \C_{\mathbb{H}}^2( O, \R )
      \text{ with }
      u(x) = \phi(x)
      \text{ and }
      u \leq \phi
    \big)
  \Big\}
\\ & =
  \Big\{ 
    \big( 
      (D_\mathbb{H} \, \phi)(x),
      ( D_\mathbb{H}^2\, \phi)(x)
    \big)
    \in 
    H' \times \mathbb{S}_{\mathbb{H}, \mathbb{H}'} \\
	& \qquad \qquad
    \colon
    \big(
      \phi \in \C_\mathbb{H}^2( O, \R )
      \text{ and }
      u - \phi 
      \text{ has a local maximum at }
      x
    \big)
  \Big\}
\end{split}
\end{equation}
(resp.\@ that
\begin{equation}
\begin{split}
  ( J^2_{ \mathbb{H}, - } u)( x )
& =
  \Big\{ 
    \big( 
      (D_\mathbb{H} \, \phi)(x),
      ( D_\mathbb{H}^2\, \phi)(x)
    \big)
    \in 
    H' \times \mathbb{S}_{\mathbb{H}, \mathbb{H}'} \\
	& \qquad \qquad
    \colon
    \big(
      \phi \in \C_\mathbb{H}^2( O, \R )
      \text{ with }
      u(x) = \phi(x)
      \text{ and }
      u \geq \phi
    \big)
  \Big\}
\\ & =
  \Big\{ 
    \big( 
      (D_\mathbb{H} \, \phi)(x),
      ( D_\mathbb{H}^2\, \phi)(x)
    \big)
    \in 
    H' \times \mathbb{S}_{\mathbb{H}, \mathbb{H}'} \\
	& \qquad \qquad
    \colon
    \big(
      \phi \in \C_\mathbb{H}^2( O, \R )
      \text{ and }
      u - \phi 
      \text{ has a local minimum at }
      x
    \big)
  \Big\}).
\end{split}
\end{equation}
\end{lemma}

As an immediate consequence we get the following corollary
	(see also Remark~2.3 in 
	Crandall, Ishii \citationand~Lions~\cite{CrandallIshiiLions1992}
	and Corollary 4.5 in Hairer, Hutzenthaler 
	\& Jentzen \cite{HairerHutzenthalerJentzen2015}).

\begin{corollary}[Characterizations of viscosity solutions]
\label{cor:semijets_equivalence}
	Assume the setting in Section \ref{ssec: Setting H X without t},
	let 
	$
		F \colon 
			W \times \R \times H' \times \mathbb{S}_{\mathbb{X}, \mathbb{X}'} \to \R
	$ 
	be a degenerate elliptic function,
  and let $ u \colon O \to \R \cup \{ -\infty \} $ be locally bounded from above.
  Then the following three assertions are equivalent:
  \begin{itemize}
   \item 
  $ u $ is a viscosity subsolution of $ F = 0 $ relative to $(h, \mathbb{H}, \mathbb{X})$
  ($ u $ is a viscosity solution of $ F \leq 0 $ relative to $(h, \mathbb{H}, \mathbb{X})$),
  \item
  for every $ x \in O $, $\delta \in (0,\infty)$ with 
	$u_{\mathbb{H}, \delta, h}^{-, W}(x) \in \R$,
  and every 
  $
    \phi
    \in
    \{ 
      \psi \in \C_\mathbb{H}^2( O, \R) \colon
      x \text{ is a }
	$
	$
				\text{local maximum of }
      ( u_{\mathbb{H}, \delta, h}^{-, W} - \psi ) \colon O \to \R \cup \{ -\infty \}
    \}
  $ 
  it holds that \\
  $
    F^+_{\mathbb{H},\mathbb{X},\delta, h, u}\big(
      x, u_{\mathbb{H}, \delta, h}^{-, W}(x), (D_{\mathbb{H}} \, \phi)(x), 
      (D_{\mathbb{H}}^2\,\phi)(x)
    \big)
    \leq 0
  $,
  \item
  for every $ x \in O $, $\delta \in (0,\infty)$,
  and every $ (p,A) \in ( J^2_{ \mathbb{H}, + } u_{\mathbb{H}, \delta, h}^{-, W})( x ) $
  it holds that \\
  $
			F^+_{\mathbb{H}, \mathbb{X}, \delta, h, u}
				( x, u_{\mathbb{H}, \delta, h}^{-, W}(x), p, A ) 
		\leq 0
  $.
  \end{itemize}
  The corresponding statement holds for 
  viscosity supersolutions relative to $(h, \mathbb{H}, \mathbb{X})$ and
  viscosity solutions relative to $(h, \mathbb{H}, \mathbb{X})$.
\end{corollary}

For convenience of the reader we recall another characterization
of viscosity solutions which corresponds to the
statement given on page 608 and 609 in Ishii \cite{Ishii1993}
(see also Remark~2.4 
in Crandall, Ishii \citationand~Lions~\cite{CrandallIshiiLions1992}
and Corollary 4.6 in Hairer, Hutzenthaler 
\& Jentzen \cite{HairerHutzenthalerJentzen2015}).
It follows immediately from Corollary~\ref{cor:semijets_equivalence},
Remark \ref{rem: F hat is semicontinuous},
and the fact that for all $x \in O$ with 
$u^{-, W}_{\mathbb{H}, \delta, h}(x) \in \{-\infty, \infty\}$ it holds that
$
		(J^2_{\mathbb{H},+} u^{-, W}_{\mathbb{H}, \delta, h})(x)
	=
	 \emptyset.
$

\begin{corollary}[Alternative characterizations of viscosity solutions]
\label{cor:semijets_equivalence2}
	Assume the setting in Section \ref{ssec: Setting H X without t},
	let 
	$F \colon W \times \R \times H' \times \mathbb{S}_{\mathbb{X}, \mathbb{X}'} \to \R$ 
	be a degenerate elliptic function,
  and let $ u \colon O \to \R \cup \{-\infty\} $ be locally bounded from above.
  Then 
  $ u $ is a viscosity solution of $ F \leq 0 $ relative to $(h, \mathbb{H}, \mathbb{X})$
  if and only if it holds
  for every $ x \in O $, $\delta \in (0,\infty)$,
  and every $ (p,A) \in ( \hat{J}^2_{ \mathbb{H}, + } u^{-, W}_{\mathbb{H}, \delta, h})( x ) $
  that
  $
    F^+_{\mathbb{H}, \mathbb{X}, \delta, h, u}
			( x, u^{-, W}_{\mathbb{H}, \delta, h}(x), p, A ) \leq 0
  $.
  The corresponding statement holds for 
  viscosity solutions of $ F \geq 0 $ relative to $(h, \mathbb{H}, \mathbb{X})$
	and of $ F = 0 $ relative to $(h, \mathbb{H}, \mathbb{X})$,
  respectively.
\end{corollary}

For the next proposition we need an additional lemma
from functional analysis.
\begin{lemma}
\label{lem: h weakly lower semicontinuous}
	Let $\mathbb{H} = (H, \langle \cdot, \cdot \rangle_H, \| \cdot \|_H)$
	be a real Hilbert space,
	let $ O \subseteq H $
  be an open and convex set,
	let $h \colon O \to \R \cup \{ -\infty, \infty \}$ be 
	a with respect to the $\| \cdot \|_H$-norm
	lower semicontinuous and convex function, 
	and let
	$x \in O$ and $(x_n)_{n \in \N} \subseteq O$ 
	satisfy that
	$w-\lim_{n \to \infty} x_n = x$.
	Then it holds that
	$\liminf_{n \to \infty} h(x_n) \geq h(x)$.
\end{lemma}
\begin{proof}[Proof
of Lemma~\ref{lem: h weakly lower semicontinuous}]
	Without loss of generality we assume that
	$
			\liminf_{n \to \infty} h(x_n)
		=
			\lim_{n \to \infty} h(x_n)
	$
	(else take a subsequence).
	Next notice that the uniform convergence principle implies that
	$(x_n)_{n \in \N} \subseteq O$ is $\mathbb{H}$-bounded
	and thus
	the Banach-Saks Theorem implies that there exists a 
	subsequence $(n_j)_{j \in \N} \subseteq \N$ satisfying that
	$\lim_{N \to \infty} \frac 1N \sum_{j=1}^{N} x_{n_j} = x$.
	Moreover, the assumption that $O$ is convex ensures that
	for all $N \in \N$ it holds that
	$\frac 1N \sum_{j=1}^{N} x_{n_j} \in O$ and combining this with the assumption that
	$h$ is convex and lower semicontinuous 
	with respect to the $\| \cdot \|_{\mathbb{H}}$-norm shows then that
	\begin{equation}
			h(x) 
		\leq 
			\lim_{N \to \infty} h \left( \frac 1N \sum_{j=1}^{N} x_{n_j} \right)
		\leq
			\lim_{N \to \infty} \frac 1N \sum_{j=1}^{N} h(x_{n_j})
		=
			\lim_{j \to \infty} h(x_{n_j}) 
		=
			\liminf_{n \to \infty} h(x_{n}),
	\end{equation}
	which completes the proof of Lemma \ref{lem: h weakly lower semicontinuous}.
\end{proof}
The next two propositions, which basically are Proposition 2.2 in
Ishii~\cite{Ishii1993}, give sufficient conditions for
the assumptions on $h$ in Definition \ref{d:viscosity.solution}.
\begin{prop}[Points with nonempty semijet are dense for convex, semicontinuous functions]
\label{prop: h assumption}
	Let $\mathbb{H} = (H, \langle \cdot, \cdot \rangle_H, \| \cdot \|_H)$ 
	be a real Hilbert space,
	let $\mathbb{H}'=(H', \| \cdot \|_{H'})$ be its dual space,
	let $ O \subseteq H $
  be an open and convex set,
	and let $h \colon O \to \R \cup \{ \infty \}$ be 
	a with respect to the $\| \cdot \|_H$-norm
	lower semicontinuous and convex function.
	Then $\{ y \in O \colon (J^{2}_{\mathbb{H}, -}h) (y) \neq \emptyset \}$
	is a dense subset of $O$ with respect to the $\| \cdot \|_H$-norm. 
\end{prop}
\begin{proof}[Proof
of Proposition~\ref{prop: h assumption}]
	If $ O = \emptyset $, 
  then the assertion is trivial. 
  So for the rest of 
  the proof, we assume that 
  $ O \neq \emptyset $.
	In addition,
	for the rest of the proof
	fix $\hat{x} \in O$ and 
	denote by
	$K \subseteq H$ the set satisfying that
	$
			K
		= 
			\{  
				y \in O \colon 
					\|y - \hat{x}\|_H \leq 
					(\dist_{\mathbb{H}}(\hat{x},H \backslash O)/2) \wedge 1
			\}.
	$
	Furthermore, 
	for every $\alpha \in (0,\infty)$ let
	$ 
			( x^{(n)}_{\alpha} )_{n \in \N} 
		\subseteq 
			K
	$ 
	be a sequence with the property that 
	\begin{equation}
	\label{eq: minimizing h+norm}
			\lim_{n \to \infty} ( h(x^{(n)}_{\alpha}) 
				+ \alpha \| x^{(n)}_{\alpha} - \hat{x} \|_H^2) 
		= 
			\inf_{x \in K} 
				( h(x) + \alpha \| x - \hat{x} \|_H^2).
	\end{equation}
	Since $K$ is $\mathbb{H}$-closed, $\mathbb{H}$-bounded, and convex
	it follows from the Banach-Saks Theorem and from Lemma 5.1.4 in 
	Kato \cite{Kato1980}
	that it is also weakly compact. Therefore
	for all $\alpha \in (0, \infty)$
	there exists a 
	$
		\hat{y}_{\alpha} 
			\in K
	$ 
	and a sequence $(n^{\alpha}_j)_{j \in \N} \subseteq \N$
	such that
	$w-\lim_{j \to \infty} x^{(n^{\alpha}_j)}_{\alpha} = \hat{y}_{\alpha}$.
	This together with \eqref{eq: minimizing h+norm},
	Lemma \ref{lem: h weakly lower semicontinuous},
	and the fact that the functions 
	$O \ni x \to h(x) + \alpha \| x - \hat{x} \|_H^2 \in \R \cup \{\infty\}$,
	$\alpha \in (0, \infty)$,
	are convex and lower semicontinuous 
	with respect to the $\| \cdot \|_H$-norm shows that
	for all $\alpha \in (0, \infty)$ it holds that
	\begin{equation}
		\begin{split}
				\lim_{j \to \infty} 
					h(x^{(n^{\alpha}_j)}_{\alpha}) 
					+ \alpha \| x^{(n^{\alpha}_j)}_{\alpha} - \hat{x} \|_H^2
			\geq{} 
				&h(\hat{y}_{\alpha}) 
				+ \alpha \| \hat{y}_{\alpha} - \hat{x} \|_H^2
			\geq 
				\inf_{x \in K}
					( h(x) + \alpha \| x - \hat{x} \|_H^2) \\
			={}
				&\lim_{j \to \infty} 
					h(x^{(n^{\alpha}_j)}_{\alpha}) 
					+ \alpha \| x^{(n^{\alpha}_j)}_{\alpha} - \hat{x} \|_H^2.
		\end{split}
	\end{equation}
	Thus, for every $\alpha \in (0,\infty)$,
	$\hat{y}_{\alpha}$ is the minimum of the function 
	$
		K \ni x 
			\to h(x) + \alpha \| x - \hat{x} \|_H^2 \in \R \cup \{ \infty \}
	$
	and we hence have for all 
	$x \in K$ that
	\begin{align}
	\label{eq: max point}
	\nonumber
				&h(x)
			\geq{}
				h(\hat{y}_{\alpha}) 
				+ \alpha \| \hat{y}_{\alpha} - \hat{x} \|_H^2
				- \alpha \| x - \hat{x} \|_H^2 \\
			={}
				&h(\hat{y}_{\alpha})
				+ \alpha \langle 
						(\hat{y}_{\alpha} - \hat{x}) - (x - \hat{x}),
						(\hat{y}_{\alpha} - \hat{x}) + (x - \hat{x})
					\rangle_H \\
		\nonumber
			={}
				&h(\hat{y}_{\alpha})
				- \alpha \langle 
						x- \hat{y}_{\alpha},
						2\hat{y}_{\alpha} - 2\hat{x} + x -\hat{y}_{\alpha}
					\rangle_H
			=
				h(\hat{y}_{\alpha})
			 	- 2\alpha \langle 
						x - \hat{y}_{\alpha},
						\hat{y}_{\alpha}- \hat{x}
					\rangle_H
				- \alpha \| \hat{y}_{\alpha} -x \|^2_H.
	\end{align}
	Moreover, since $h$ is convex and lower semicontinuous 
	with respect to the $\| \cdot \|_H$-norm,
	Lemma \ref{lem: h weakly lower semicontinuous} shows that
	it is also bounded from below on the weakly compact set
	$K$.
	Thus it holds that
	$
		\inf_{x \in K}
			h(x) > -\infty
	$
	and this yields that for
	all $\alpha \in (0, \infty)$ 
	and all $r \in (0, 1 \wedge \dist_{\mathbb{H}}(\hat{x}, H \backslash O)/2)$ 
	it holds that
	\begin{equation}
		\begin{split}
			&(\inf_{x \in K}
				h(x) )
			+\alpha \| \hat{y}_{\alpha}- \hat{x} \|_H^2 
			\leq 
				h(\hat{y}_{\alpha}) 
				+ \alpha \| \hat{y}_{\alpha}- \hat{x} \|_H^2
			=
				\inf_{x \in K} \left (
					h(x) + \alpha \| x- \hat{x} \|_H^2
				\right ) \\
			\leq{}&
				\inf_{\substack{x \in O \\ \| x- \hat{x} \|_H \leq r}} \left (
					h(x) + \alpha \| x- \hat{x} \|_H^2
				\right ) 
			\leq
				\inf_{\substack{x \in O \\ \| x- \hat{x} \|_H \leq r}} \left (
					h(x) + \alpha r^2
				\right ).
		\end{split}
	\end{equation}
	Dividing now by $\alpha$ and taking the limits 
	$\alpha \to \infty$ and $ r \to 0$ we get
	\begin{equation}
		\begin{split}
				&\lim_{\alpha \to \infty} \left (
					\frac{\inf_{x \in K} h(x)}{\alpha}
					+ \| \hat{y}_{\alpha}- \hat{x} \|_H^2
				\right) 
			=
				\lim_{\alpha \to \infty} \left (
					\| \hat{y}_{\alpha}- \hat{x} \|_H^2
				\right) \\
			\leq{}
				&\lim_{r \downarrow 0} \lim_{\alpha \to \infty}
				\inf_{\substack{x \in O \\ \| x- \hat{x} \|_H^2 \leq r}} \left (
						\frac{h(x)}{\alpha} + r^2
					\right )
			=
				0.
		\end{split}
	\end{equation}
	This yields that $\lim_{\alpha \to \infty} \hat{y}_{\alpha} = \hat{x}$.
	Thus choosing $\alpha \in (0, \infty)$ big enough we get that 
	$\hat{y}_{\alpha}$ is an inner point of 
	$K$
	and combining this with \eqref{eq: max point} shows that
	$(-2\alpha I_{\mathbb{H}} (\hat{y}_{\alpha}- \hat{x}) , -2\alpha I_{\mathbb{H}})
		\in \big (J^2_{\mathbb{H}, -}h \big ) (\hat{y}_{\alpha})$.
	Thus,
	since $\hat{x}$ was arbitrary,
	we obtain that $\{ y \in O \colon (J^{2}_{\mathbb{H}, -}h) (y) \neq \emptyset \}$
	is a dense subset of $O$ with respect to the $\| \cdot \|_H$-norm
	and this completes the proof of Proposition \ref{prop: h assumption}.
\end{proof}
\begin{prop}
\label{prop: diff and semijet}
	Assume the setting in Section \ref{ssec: Setting H X without t},
	and let $x \in X \cap O$ 
	and $(p,A) \in (J^{2}_{\mathbb{H}, -}h) (x)$. Then
	it holds that 
	$(D_{\mathbb{X}} (h|_X)) (x) \in D(E_{\mathbb{X}', \mathbb{H}'})$,
	that $E_{\mathbb{X}', \mathbb{H}'} \big( (D_{\mathbb{X}} (h|_X))(x) \big) = p$,
	and that $(A|_{X}) |_X \leq (D^2_{\mathbb{X}} (h|_X))(x)$.
\end{prop}
\begin{proof}[Proof of Proposition~\ref{prop: diff and semijet}]
	Since $\mathbb{X}$ is embedded continuously
	and densely in $\mathbb{H}$ and 
	$(p,A ) \in (J^{2}_{\mathbb{H}, -}h) (x)$ we have that 
	$
			\liminf_{O \cap X \ni y \to x} 
				\frac
					{
						h(y)-h(x) - \langle p, y-x \rangle_{H', H} 
						- \frac 12 \langle A(y-x),(y-x) \rangle_{H' ,H}
					}
					{\|x-y\|^2_X} 
		\geq 
			0
	$.
	This yields that
	\begin{equation}
		\begin{split}
			&\liminf_{O \cap X \ni y \to x} 
				\frac
					{
						\left \langle (D_{\mathbb{X}} (h |_X)) (x) - p |_X, y-x \right \rangle_{X', X}}
					{\|x-y\|^2_X} \\
			&\qquad \qquad
				+\frac
					{\tfrac 12
							\left \langle 
								((D^2_{\mathbb{X}} (h|_X)) (x) - (A|_X) |_X)(y-x), y-x 
							\right \rangle_{X', X}
					}
					{\|x-y\|^2_X} \\
		={}&
			\liminf_{O \cap X \ni y \to x} \bigg (
				\frac
					{
						h(y)-h(x)-\langle p, y-x \rangle_{H', H} 
						- \tfrac 12 \langle A (y-x), y-x \rangle_{H', H}
					}
					{\|x-y\|^2_X} 
				-\frac
					{
						h(y)-h(x)
					}
					{\|x-y\|^2_X}\\
			& \qquad 
				-\frac
					{
						-\langle (D_{\mathbb{X}} (h|_X)) (x), y-x \rangle_{X', X} 
						- \tfrac 12 \langle (D^2_{\mathbb{X}} (h|_X)) (x) (y-x), y-x \rangle_{X', X}
					}
					{\|x-y\|^2_X}
			\bigg )
		\geq
			0.
		\end{split}
	\end{equation}
	Hence by choosing $y= x + t e$ we get for all $e \in X$ with $\|e \|_X =1$
	that
	\begin{equation}
		\begin{split}
			\liminf_{0 \neq t \to 0} 
			\left(
				\frac
					{\left \langle (D_{\mathbb{X}} (h|_X)) (x) - p|_X, e \right \rangle_{X', X} }{t}
				+ \tfrac 12 \left \langle 
						((D^2_{\mathbb{X}} (h|_X)) (x))e - (A e)|_X, e 
					\right \rangle_{X', X}
			\right)
		\geq
			0
		\end{split}
	\end{equation}
	and thus it holds for all $e \in X$ with $\|e \|_X =1$ that
	$\left \langle (D_{\mathbb{X}} (h|_X)) (x) - p |_X, e \right \rangle_{X', X} = 0$ and that \\
	$\left \langle ((D_{\mathbb{X}}^2(h|_X)) (x))e - (A e) |_X, e \right \rangle_{X', X} \geq 0$.
	Combining this with the fact that $X$ is dense in $H$ then shows that 
	$(D_{\mathbb{X}} (h |_X)) (x) = p |_X $ and that 
	$(A|_{X}) |_X \leq (D_{\mathbb{X}}^2 (h |_X)) (x)$.
	In addition, $(D_{\mathbb{X}} (h |_X)) (x) = p |_X $ implies that
	for all $y \in X$ it holds that
	$
			|\langle (D_{\mathbb{X}} (h |_X)) (x), y \rangle_{X',X}|
		=
			|\langle p, y \rangle_{H',H}|
		\leq 
			\|p \|_{H'} \cdot \|y\|_H
	$
	and therefore
	$(D_{\mathbb{X}} (h |_X)) (x) \in D(E_{\mathbb{X}',\mathbb{H'}})$.
	This completes the proof
	of Proposition~\ref{prop: diff and semijet}.
\end{proof}
The next proposition ensures that $h$ still fulfills 
the assumptions in Definition \ref{d:viscosity.solution} if it is multiplied
by a sufficiently nice function.
\begin{prop}[Stability of the set of "h-functions" under multiplication with positive
$\C^2$-functions]
\label{prop: h and Vh}
	Assume the setting in Section \ref{ssec: Setting H X without t}
	and let
  $
    V \in
    \C_{\mathbb{H}}^2( O, (0,\infty))
  $.
	Then $Vh$ is lower semicontinuous 
	with respect to the $\| \cdot \|_{\mathbb{H}}$-norm and
	it holds that
	$(Vh) |_{X} \in \C_{\mathbb{X}}^2(X \cap O, \R)$
	and that
	$
			\{ 
				y \in X \cap O \colon 
					(D_{\mathbb{X}} (Vh)|_X) (y) \in D(E_{\mathbb{X}', \mathbb{H}'}), 
					~ (J^{2}_{\mathbb{H}, -} (Vh)) (y) \neq \emptyset 
			\}
		=
			W.
	$
\end{prop}
\begin{proof} [Proof
of Proposition~\ref{prop: h and Vh}]
	First note that the lower semicontinuity of $h$ 
	with respect to the $\| \cdot \|_{\mathbb{H}}$-norm
	together with 
	$ h |_{X} \in \C_{\mathbb{X}}^2(X \cap O, \R)$ 
	and with
	$ V \in \C_{\mathbb{H}}^2(O, (0,\infty)) $
	ensures that 
	$Vh$ is lower semicontinuous
	with respect to the $\| \cdot \|_{\mathbb{H}}$-norm,
	that
	$
		(Vh) |_{X} \in \C_{\mathbb{X}}^2(X \cap O, \R),
	$
	and that
	$
			\{ 
				z \in O \cap X \colon 
					(D_{\mathbb{X}} ((Vh)|_{X })) (z) \in D(E_{\mathbb{X}', \mathbb{H}'}) 
			\}
		=
			\{ 
				z \in O \cap X \colon 
					(D_\mathbb{X} (h|_{X})) (z) \in D(E_{\mathbb{X}', \mathbb{H}'}) 
			\}.
	$
	Next, Lemma \ref{lem:semijets} implies that 
	for every
	$
		z \in
			\{ y \in O \colon (J^{2}_{\mathbb{H}, -} h) (y) \neq \emptyset \}
	$
	there exists a 
	$ \phi_z \in \C_{\mathbb{H}}^2(O, \R)$
	such that
	$ 
    \phi_z( z ) = 
    h (z)
  $ 
  and that
  $ 
    \phi_z \leq h.
  $
	Thus we have
	for every
	$
		z \in
			\{ y \in O \colon (J^{2}_{\mathbb{H}, -} h) (y) \neq \emptyset \}
	$
	that
	$ 
    V(z) \phi_z( z ) = 
    V(z) h (z)
  $ 
  and that
  $ 
    V \phi_z \leq V h
  $
	and this implies that
	$
			\left( (D_{\mathbb{H}} (\phi_z V))(z),   (D_{\mathbb{H}}^2 (\phi_z V))(z) \right)
		\in (J^{2}_{\mathbb{H}, -} (Vh)) (z).
	$
	Hence we get that
	$
			\{ z \in O \colon (J^{2}_{\mathbb{H}, -} h) (z) \neq \emptyset \}
		\subseteq
			\{ z \in O \colon (J^{2}_{\mathbb{H}, -} (Vh)) (z) \neq \emptyset \}.
	$
	Analogously it follows that
	$
			\{ z \in O \colon (J^{2}_{\mathbb{H}, -} (Vh)) (z) \neq \emptyset \}
		\subseteq
			\{ z \in O \colon (J^{2}_{\mathbb{H}, -} h) (z) \neq \emptyset \}
	$
	and this together with
	$
			\{ 
				z \in X \cap O \colon 
					(D_{\mathbb{X}} ((Vh)|_{X})) (z) \in D(E_{\mathbb{X}', \mathbb{H}'}) 
			\}
		=
			\{ 
				z \in X \cap O \colon 
					(D_\mathbb{X} (h|_{X})) (z) \in D(E_{\mathbb{X}', \mathbb{H}'}) 
			\}
	$
	shows that
	\begin{equation}
		\begin{split}
				&\{ 
					y \in X \cap O \colon 
						(D_{\mathbb{X}} ((Vh)|_X)) (y) \in D(E_{\mathbb{X}', \mathbb{H}'}), 
						~ (J^{2}_{\mathbb{H}, -} (Vh)) (y) \neq \emptyset 
				\} 
			=
				W.
		\end{split}
	\end{equation}
	This completes the proof of Proposition~\ref{prop: h and Vh}.
\end{proof}
Next we show that classical solutions are also viscosity solutions,
cf.\@ Ishii \cite{Ishii1993}.
 \begin{lemma}[Classical solutions are viscosity solutions]
  \label{lem:classicalsolution}
	Assume the setting in Section \ref{ssec: Setting H X without t},
	assume that
	$
		\{ y \in O \colon (J^{2}_{\mathbb{H}, -}h) (y) \neq \emptyset \} \subseteq X
	$,
	let
	$F \colon W  
		\times \R \times H' \times \mathbb{S}_{\mathbb{X}, \mathbb{X}'} \to \R$ 
	be a degenerate elliptic function,
	and
	let
  $ u \in \C_{\mathbb{H}}^2( O, \R ) $ 
  be a classical subsolution 
  of $ F = 0 $, i.e., suppose that for all $ x \in W $ it holds that
  \begin{equation}
  \label{eq:classical}
    F\big(
      x, u(x), (D_{\mathbb{H}} \, u)(x), 
      ((D_{\mathbb{H}}^2 \, u)(x) |_X) |_X
    \big)
    \leq 0.
  \end{equation}
  Then $ u $
  is also
  a viscosity subsolution 
  of $ F  = 0 $ 
	relative to $(h, \mathbb{H}, \mathbb{X})$.
\end{lemma}
\begin{proof}[Proof
of Lemma~\ref{lem:classicalsolution}]
	If $O = \emptyset$, 
	then the assertion is trivial. 
	So for the rest of the proof, we assume that 
	$O \neq \emptyset $.
	In addition, note that the assumption that $h$ is lower semicontinuous 
	with respect to the $\| \cdot \|_{\mathbb{H}}$-norm, and 
	$u$ is continuous with respect to the $\| \cdot \|_{\mathbb{H}}$-norm
	implies that for all $\delta \in (0, \infty)$
	it holds that $u- \delta h$ is upper semicontinuous 
	with respect to the $\| \cdot \|_{\mathbb{H}}$-norm,
	and thus it holds
	for all $\delta \in (0, \infty)$
	and all $x \in W$ that
	$
			u^{-, W}_{\mathbb{H}, \delta, h}(x)
		=
			\overline{(u -\delta h)}_{\mathbb{H}}^W (x)
		= 
			u(x)- \delta h(x)
	$. 
	Moreover, since
	$u \in \C_{\mathbb{H}}^2(O, \R)$ we get for all $x \in H$
	and all $\delta \in (0, \infty)$ that 
	\begin{equation}
	\label{eq: semijet equation}
			\delta (J^2_{\mathbb{H}, -} h) (x) 
		= 
			(J^2_{\mathbb{H}, +} u)(x) - (J^2_{\mathbb{H}, +} (u - \delta h))(x).
	\end{equation}
	This together with the assumption that
	$W$ is dense in $O$ with respect to the $\| \cdot \|_H$-norm
	and with $O \neq \emptyset$ shows
	that for all $\delta \in (0, \infty)$ it holds that
	$
			\{ 
				y \in O \colon 
				~ (J^{2}_{\mathbb{H}, +} (u- \delta h)) (y) \neq \emptyset 
			\}
		\supseteq
			\{ 
				y \in O \colon 
				~ \delta (J^{2}_{\mathbb{H}, -} (h)) (y) \neq \emptyset
			\}
		\supseteq
		 W
		\neq \emptyset.
	$
	In the next step fix $\delta \in (0, \infty)$, 
	$ 
		x \in 
		\{ 
			y \in O \colon 
			~ (J^{2}_{\mathbb{H}, +} (u- \delta h)) (y) \neq \emptyset 
		\},
	$
	and $(p,A) \in (J^2_{\mathbb{H}, +} (u - \delta h))(x)$.
	Then \eqref{eq: semijet equation} implies that 
	$
		( (D_{\mathbb{H}} \, u )(x) - p, (D_{\mathbb{H}}^2 \, u) (x) - A) 
			\in \delta (J^2_{\mathbb{H}, -} h)(x)
	$ 
	and thus the assumption
	$\{ y \in O \colon (J^{2}_{\mathbb{H}, -}h) (y) \neq \emptyset \} \subseteq X$
	implies that $x \in X$.
	Furthermore, 
	it follows from
	Proposition \ref{prop: diff and semijet} that
	\begin{equation}
		\begin{split}
			&(D_{\mathbb{X}} (h|_X))(x) \in D(E_{\mathbb{X}', \mathbb{H}'}), \quad
			(D_{\mathbb{H}} \, u)(x) - p 
		= 
			\delta E_{\mathbb{X}', \mathbb{H}'} \big( (D_{\mathbb{X}} (h|_X))(x) \big),
			\quad \text{and that}\\
			&(((D_{\mathbb{H}}^2 \, u)(x) - A)|_X) |_X \leq \delta (D_{\mathbb{X}}^2 (h|_X))(x).
		\end{split}
	\end{equation}
	This shows that $x \in W$.
	Combining this with
	\eqref{eq: def F+ delta h u}, 
	the fact that
	$u^{-, W}_{\mathbb{H}, \delta, h} (x) = u(x)- \delta h(x)$,
	the assumption that $F$ is degenerate elliptic,
	the assumption that $W \subseteq X$,
	and
	the assumption that $u$ is a classical subsolution of $F = 0$,
	shows that
	\begin{equation}
		\begin{split}
				&F^+_{\mathbb{H}, \mathbb{X}, \delta, h, u}
					(x, u^{-, W}_{\mathbb{H}, \delta, h} (x), p ,A)
			\leq
				F^+_{\mathbb{H}, \mathbb{X}, \delta, h}(x, (u -\delta h)(x), p ,A) \\
			={}&
				F(
					x, 
					(u -\delta h)(x) + \delta h(x), 
					p+ \delta E_{\mathbb{X}', \mathbb{H}'} \big( (D_{\mathbb{X}} (h|_X))(x) \big) ,
					(A |_X) |_X + \delta (D^2_{\mathbb{X}} (h|_X))(x)
				) \\
			={}&
				F(
					x, 
					u, 
					(D_{\mathbb{H}} \, u)(x),
					(A |_X) |_X + \delta (D^2_{\mathbb{X}} (h|_X))(x)
				) \\
			\leq{}&
				F(x, u(x), (D_{\mathbb{H}} \, u)(x) ,((D_{\mathbb{H}}^2 \, u)(x) |_X) |_X)
			\leq 0.
		\end{split}
	\end{equation}
	Thus Corollary \ref{cor:semijets_equivalence} shows that
	u is a viscosity subsolution of 
	$ F = 0 $ 
	relative to $(h, \mathbb{H}, \mathbb{X})$ and this
	completes the proof
	of Lemma \ref{lem:classicalsolution}.
\end{proof}
The next lemma connects the definition of the viscosity solution given above with the
definition given in Crandall, Ishii \citationand\ 
Lions~\cite{CrandallIshiiLions1992}.
\begin{lemma}[Continuous viscosity solutions are
	almost classical viscosity solutions]
\label{lem: ordinary visc. solution}
	Assume the setting in Section \ref{ssec: Setting H X without t},
	assume that for all $x \in X$ it holds that
	\begin{equation}
	\label{eq: h nicely approximated by its derivatives}
		\begin{split}
				\lim_{\eps \to 0}
				\inf \bigg \{
					&\frac
						{
							h(y)-h(x) 
							- \langle (D_{\mathbb{X}}(h|_X))(x), y-x \rangle_{X', X}
							- \frac 12 \langle (D_{\mathbb{X}}^2(h|_X))(x)(y-x), y-x \rangle_{X', X}
						} 
						{\| x- y \|^2_{H}}
					\colon \\ 
					&y \in X, 
					~\| y-x \|_H \leq \eps
				\bigg \}
			\geq 0,
		\end{split}
	\end{equation}
	let
	$
		F \colon W  
			\times \R \times H' \times \mathbb{S}_{\mathbb{X}, \mathbb{X}'} \to \R
	$ 
	be a degenerate elliptic function
	such that 
	for all $\delta \in (0, \infty)$ it holds that
	$F^+_{\mathbb{H}, \mathbb{X}, \delta, h}$ 
	is lower semicontinuous with respect to
	$d^{-, W}_{\mathbb{H}, \delta, h, u}$
	and that
	for all 
	$x \in W$,
	$p \in H'$,
	$r \in \R$,
	$A \in \mathbb{S}_{\mathbb{X}, \mathbb{X}'}$ and all 
	$(A_n)_{n \in \N} \subseteq \mathbb{S}_{\mathbb{H}, \mathbb{H}'}$
	satisfying for all $\xi \in X$ that
	$\lim_{n \to \infty} \| A_n(\xi) |_{X} - A(\xi) \|_{X'} = 0$
	it holds that
	\begin{equation}
	\label{eq: F weak continuous in A}
		\limsup_{n \to \infty}
			F(x,r,p,(A_n |_X) |_X) \geq F(x,r,p,A),
	\end{equation}
	and let
  $ u \in \C_{\mathbb{H}}( O, \R ) $ 
  be a viscosity subsolution 
  of $ F = 0 $ 
	relative to $(h, \mathbb{H}, \mathbb{X})$.
  Then for all $x \in W$ satisfying that
	the set 
	$
		\{
			y \in X \colon 
				((D^2_{\mathbb{X}} (h|_X))(x)) (y) \in D(E_{\mathbb{X}', \mathbb{H}'})
		\}
	$ 
	is dense
	in $X$ with respect to the $\| \cdot \|_X$-norm
	and all $(p, A) \in (J_{\mathbb{H}, +}^2 u)(x)$ it holds that
	$F(x,u(x), p , (A |_X) |_X) \leq 0$.
\end{lemma}
\begin{proof}[Proof
of Lemma~\ref{lem: ordinary visc. solution}]
	As before note hat the assumption that $h$ is lower semicontinuous 
	with respect to the $\| \cdot \|_{\mathbb{H}}$-norm and 
	$u$ is continuous with respect to the $\| \cdot \|_{\mathbb{H}}$-norm
	implies that for all $\delta \in (0, \infty)$ 
	it holds that $u- \delta h$ is upper semicontinuous 
	with respect to the $\| \cdot \|_{\mathbb{H}}$-norm
	and thus it holds
	for all $x \in W$ and for all $\delta \in (0, \infty)$ that
	$
			u^{-, W}_{\mathbb{H}, \delta, h}(x)
		=
			\overline{(u -\delta h)}_{\mathbb{H}}^W (x)
		= 
			u(x)- \delta h(x)$.
	Moreover, without loss of generality we can assume
	that there exist an
	$x \in W$ and $(p, A) \in (J^2_{\mathbb{H}, +} u)(x)$ such that
	the set 
	$\{y \in X \colon ((D^2_{\mathbb{X}} (h|_X))(x)) (y) \in D(E_{\mathbb{X}', \mathbb{H}'})\}$ 
	is dense
	in $X$ with respect to the $\| \cdot \|_X$-norm
	(else the assertion is trivial).
	From the definition
	of $W$ it follows that there exists a $(\tilde{p}, \tilde{A}) \in (J^2_{\mathbb{H}, -} h)(x)$
	and Proposition \ref{prop: diff and semijet} then implies that
	$E_{\mathbb{X}', \mathbb{H}'} \big( (D_{\mathbb{X}} (h|_X))(x) \big) = \tilde{p}$
	and that $(\tilde{A}|_{X}) |_X \leq (D^2_{\mathbb{X}} (h|_X))(x)$.
	Furthermore, it follows that the operator
	$
		\{y \in X \colon ((D^2_{\mathbb{X}} (h|_X))(x)) (y) \in D(E_{\mathbb{X}', \mathbb{H}'})\}
		\ni z \to 
			I_\mathbb{H}^{-1}
				E_{\mathbb{X}', \mathbb{H}'} \Big( \big((D^2_{\mathbb{X}} (h|_X))(x) \big)(z) \Big) 
			- I_\mathbb{H}^{-1} \tilde{A}(z) \in H
	$
	is self-adjoint
	and nonnegative and thus 
	there exists a sequence $(A_n)_{n \in \N} \subseteq \mathbb{S}_{\mathbb{H}, \mathbb{H}}$ 
	satisfying
	for all 
	$
		\xi \in 
			\{y \in X \colon ((D^2_{\mathbb{X}} (h|_X))(x)) (y) \in D(E_{\mathbb{X}', \mathbb{H}'})\}
	$ 
	that 
	\begin{equation}
	\label{eq: dense An limit}
			\lim_{n \to \infty} A_n \xi 
		= 
			I_{\mathbb{H}}^{-1}E_{\mathbb{X}', \mathbb{H}'} 
				\Big( \big((D^2_{\mathbb{X}} (h|_X))(x) \big)(\xi) \Big) 
			- I_{\mathbb{H}}^{-1} \tilde{A}(\xi)
	\end{equation}
	and satisfying for all
	$
		\xi \in 
			\{y \in X \colon ((D^2_{\mathbb{X}} (h|_X))(x)) (y) \in D(E_{\mathbb{X}', \mathbb{H}'})\}
	$ 
	and all $n \in \N$ that 
	\begin{equation}
	\label{eq: An inequality}
		\begin{split}
				\langle A_n \xi, \xi \rangle_H 
			& \leq 
				\langle
					(I_{\mathbb{H}}^{-1}E_{\mathbb{X}', \mathbb{H}'} [ (D^2_{\mathbb{X}} (h|_X))(x) ] 
					- I_{\mathbb{H}}^{-1} \tilde{A}) \xi,
					\xi
				\rangle_H \\
			&=
				\langle
					(E_{\mathbb{X}', \mathbb{H}'} [ (D^2_{\mathbb{X}} (h|_X))(x) ] 
					- \tilde{A}) \xi,
					\xi
				\rangle_{H', H} 
			=
				\langle
					( (D^2_{\mathbb{X}} (h|_X))(x)) \xi
					- (\tilde{A} \xi ) |_X,
					\xi
				\rangle_{X', X}
		\end{split}
	\end{equation}
	(take 
		e.g.\@ $A_n = \int_0^n \lambda \mathrm{d}E_\lambda$
		where $E$ denotes a projection-valued measure associated to the operator
		$	
				\{y \in X \colon ((D^2_{\mathbb{X}} (h|_X))(x)) (y) \in D(E_{\mathbb{X}', \mathbb{H}'})\}
			\ni z \to 
				I_{\mathbb{H}}^{-1}E_{\mathbb{X}', \mathbb{H}'} 
					\Big( \big((D^2_{\mathbb{X}} (h|_X))(x) \big)(z) \Big) 
				- I_{\mathbb{H}}^{-1} \tilde{A}(z) \in H
		$
		(see, e.g., Ch. 6 \S 5 in Kato \cite{Kato1980} )
	).
	In addition, \eqref{eq: An inequality} together with the assumption
	that the set 
	$
		\{y \in X \colon ((D^2_{\mathbb{X}} (h|_X))(x)) (y) \in D(E_{\mathbb{X}', \mathbb{H}'})\}
	$ 
	is dense in $X$ with respect to the $\| \cdot \|_X$-norm implies that
	for all $\xi \in X$ and all $n \in \N$ it holds that
	\begin{equation}
	\label{eq: An inequality in X}
			\langle (I_{\mathbb{H}} A_n \xi)|_X, \xi \rangle_{X', X} 
		=
			\langle A_n \xi, \xi \rangle_{H}
		\leq 
			\langle
					( (D^2_{\mathbb{X}} (h|_X))(x)) \xi
					- (\tilde{A} \xi ) |_X,
					\xi
			\rangle_{X', X}.
	\end{equation}
	This together with \eqref{eq: h nicely approximated by its derivatives} shows
	then that
	for all $n \in \N$ it holds that
	\begin{equation}
		\begin{split}
			&\lim_{\eps \to 0}
				\inf \Bigg \{
					\frac
						{
							h(y)-h(x) 
							- \langle 
									E_{\mathbb{X}', \mathbb{H}'} \big((D_{\mathbb{X}} (h|_X))(x) \big),
									y-x 
								\rangle_{H', H}}
						{\| y-x \|^2_{H}} \\
					& \qquad \qquad
				-\frac
						{
							\frac 12 \langle (\tilde{A}+ I_{\mathbb{H}} A_n)(y-x), y-x \rangle_{H', H}
						}
						{\| y-x \|^2_{H}} \colon 
					y \in X, ~ \| y-x \|_H \leq \eps
				\Bigg \} \\
			\geq{}&
				\lim_{\eps \to 0}
				\inf \Bigg \{
					\frac
						{
							h(y)-h(x) 
							- \langle (D_{\mathbb{X}}(h|_X))(x), y-x \rangle_{X', X}
							- \frac 12 
								\langle \big( (D_{\mathbb{X}}^2(h|_X))(x) \big) (y-x), y-x \rangle_{X', X}
						}
						{\| y-x \|^2_{H}}
					\colon \\
				& \qquad \qquad \qquad
					~ y \in X,
					~ \| y-x \|_H \leq \eps
				\Bigg \} 
			\geq{} 0.
		\end{split}
	\end{equation}
	Combining this with the assumption that 
	$\mathbb{X}$ is continuously embedded in $\mathbb{H}$ 
	shows for all $n \in \N$ that
	$
		\big ( 
			E_{\mathbb{X}', \mathbb{H}'} \big((D_{\mathbb{X}} (h|_X))(x) \big),
			\tilde{A}+ I_{\mathbb{H}} A_n 
		\big )
			\in (J^2_{\mathbb{H}, -} h)(x)
	$.
	Furthermore, \eqref{eq: dense An limit} implies
	that for all 
	$
		\xi \in
			\{
				y \in X \colon
					((D^2_{\mathbb{X}} (h|_X))(x)) (y) \in D(E_{\mathbb{X}', \mathbb{H}'})
			\}
	$ 
	it holds that
	\begin{equation}
	\label{eq: An convergence in X on a dense set}
			\lim_{n \to \infty} \| 
				(I_\mathbb{H} A_n \xi) |_X
				- \big((D^2_{\mathbb{X}} (h|_X))(x) \big)(\xi) 
				+ (\tilde{A}(\xi)) |_X
			\|_{X'}
		= 0.
	\end{equation}
	Moreover,
	Lemma 3.2.4 in Zimmer \cite{Zimmer1990} and
	\eqref{eq: An inequality in X}
	implies that for all $n \in \N$ it holds that
	$
			\| (I_{\mathbb{H}} (A_n|_X))|_X \|_{L(\mathbb{X},\mathbb{X}')} 
		\leq 
			\| (D^2_{\mathbb{X}} (h|_X))(x) \|_{L(\mathbb{X},\mathbb{X}')} 
			+ \| (\tilde{A}|_X) |_X \|_{L(\mathbb{X},\mathbb{X}')}
	$
	and
	this together with \eqref{eq: An convergence in X on a dense set}
	yields that for all 
	$\xi \in X$ and all
	$
		\hat{\xi} \in
			\{
				y \in X \colon
					((D^2_{\mathbb{X}} (h|_X))(x)) (y) \in D(E_{\mathbb{X}', \mathbb{H}'})
			\}
	$ 
	it holds that
	\begin{equation}
		\begin{split}
				&\lim_{n \to \infty} \| 
					(I_\mathbb{H} A_n \xi) |_X
					- \big((D^2_{\mathbb{X}} (h|_X))(x) \big)(\xi) 
					+ (\tilde{A}(\xi)) |_X
				\|_{X'} \\
			\leq{}&
				\lim_{n \to \infty} \Big(
				\| 
					(I_\mathbb{H} A_n \hat{\xi}) |_X
					- \big((D^2_{\mathbb{X}} (h|_X))(x) \big)(\hat{\xi}) 
					+ (\tilde{A}(\hat{\xi})) |_X
				\|_{X'} \\
				& \qquad +
				\| 
					(I_\mathbb{H} (A_n (\xi - \hat{\xi}))) |_X
					- \big((D^2_{\mathbb{X}} (h|_X))(x) \big)(\xi-\hat{\xi}) 
					+ (\tilde{A}(\xi-\hat{\xi})) |_X
				\|_{X'} 
				\Big)\\
			\leq{}&
				\|\xi - \hat{\xi}\|_{X} \cdot
					(2\| (D^2_{\mathbb{X}} (h|_X))(x) \|_{L(\mathbb{X},\mathbb{X}')} 
					+ 2\| (\tilde{A}|_X) |_X \|_{L(\mathbb{X},\mathbb{X}')}) < \infty.
		\end{split}
	\end{equation}
	Combining this with the assumption that 
	$
		\{y \in X \colon ((D^2_{\mathbb{X}} (h|_X))(x)) (y) \in D(E_{\mathbb{X}', \mathbb{H}'})\}
	$ 
	is dense in $X$ with respect to the $\| \cdot \|_X$-norm then shows that for all
	$\xi \in X$ it holds that
	\begin{equation}
	\label{eq: An convergence in X}
			\lim_{n \to \infty} \| 
				(I_\mathbb{H} A_n \xi) |_X
				- \big((D^2_{\mathbb{X}} (h|_X))(x) \big)(\xi) 
				+ (\tilde{A}(\xi)) |_X
			\|_{X'}
		= 0.
	\end{equation}
	In addition, it follows from the fact that for all $\delta \in (0, \infty)$
	it holds that
	$
			(J^2_{\mathbb{H}, +} (u-\delta h))(x) 
		\supseteq 
			(J^2_{\mathbb{H}, +} u)(x) -\delta (J^2_{\mathbb{H}, -} h)(x)$,
	the fact that 
	$
		\big ( 
			E_{\mathbb{X}', \mathbb{H}'} \big((D_{\mathbb{X}} (h|_X))(x) \big),
			\tilde{A}+ I_{\mathbb{H}} A_n 
		\big )
			\in (J^2_{\mathbb{H}, -} h)(x)
	$,
	and the fact that
	$(p, A) \in (J^2_{\mathbb{H},+} u)(x)$
	that for all $\delta \in (0, \infty)$ it holds that
	$
		(
			p-\delta E_{\mathbb{X}', \mathbb{H}'} \big((D_{\mathbb{X}} (h|_X))(x) \big), 
	$
	$
			A-\delta(\tilde{A}+ I_{\mathbb{H}} A_n)
		) 
		\in (J^2_{\mathbb{H},+} (u-\delta h))(x)
	$.
	Thus we get with 
	Corollary \ref{cor:semijets_equivalence}
	and with the assumption that 
	$u$ is a viscosity subsolution 
  of $ F= 0 $ 
	relative to $(h, \mathbb{H}, \mathbb{X})$
	that for all $n \in \N$ and all $\delta \in (0, \infty)$ it holds that
	$
			F^+_{\mathbb{H}, \mathbb{X}, \delta, h, u}(
				x, 
				u^{-,W}_{\mathbb{H},\delta, h}(x), 
				p-\delta E_{\mathbb{X}', \mathbb{H}'} \big((D_{\mathbb{X}} (h|_X))(x) \big),
				A-\delta(\tilde{A} + I_{\mathbb{H}} A_n)) 
		\leq 0.
	$
	Combining this with
	\eqref{eq: def F+ delta h u},
	\eqref{eq: F weak continuous in A},
	\eqref{eq: An inequality in X},
	the assumption that for all
	$\delta \in (0, \infty)$ it holds that
	$F^+_{\mathbb{H}, \mathbb{X}, \delta, h}$ is  lower semicontinuous
	with respect to $d^{-, W}_{\mathbb{H}, \delta, h, u}$,
	and with the fact that for all
	$\delta \in (0, \infty)$ 
	and all $x \in W$
	it holds that
	$u^{-, W}_{\mathbb{H}, \delta, h} (x) = u(x)- \delta h(x)$ shows that
	for all $\delta \in (0, \infty)$ it holds that
	\begin{equation}
		\begin{split}
			0 \geq{}&
				\limsup_{n \to \infty} \left(
					F^+_{\mathbb{H}, \mathbb{X}, \delta, h, u}(
						x, 
						u^{-,W}_{\mathbb{H},\delta, h}(x), 
						p-\delta E_{\mathbb{X}', \mathbb{H}'} \big((D_{\mathbb{X}} (h|_X))(x) \big),
						A-\delta(\tilde{A}+ I_{\mathbb{H}} A_n)
					)
				\right) \\
			={}&
				\limsup_{n \to \infty} \left(
					F^+_{\mathbb{H}, \mathbb{X}, \delta, h}(
						x, 
						(u -\delta h)(x), 
						p-\delta E_{\mathbb{X}', \mathbb{H}'} \big((D_{\mathbb{X}} (h|_X))(x) \big),
						A -\delta(\tilde{A} + I_{\mathbb{H}} A_n )
					) 
				\right) \\
			={}&
				\limsup_{n \to \infty} \left(
					F(
						x, 
						u(x), 
						p,
						(A |_X) |_X
						-\delta( (\tilde{A} |_X) |_X+ (I_{\mathbb{H}}A_n |_X) |_X) 
						+ \delta (D_{\mathbb{X}}^2(h|_X))(x)
					)
				\right) \\
			\geq{}
				&F(x, u(x), p, (A |_X) |_X),
		\end{split}
	\end{equation}
	which completes the proof of lemma \ref{lem: ordinary visc. solution}.
\end{proof}
\section{An approximation result for viscosity solutions}

The next lemma implies that locally uniform limits of viscosity solutions
are again viscosity solutions. Lemma \ref{l:limits.of.viscosity.solutions}
will be used in the proof of
Theorem \ref{thm: existence}.
\begin{lemma}[Stability of viscosity solutions under limits]
  \label{l:limits.of.viscosity.solutions}
	Assume the setting in Section \ref{ssec: Setting H X without t},
	assume that $h$ is bounded from below on $\mathbb{H}$-bounded
	subsets of $O$,
  let 
  $ u_n \in \C_{\mathbb{H}}(O,\R) $, 
  $ n \in \N_0 $, 
	be
	bounded 
	on $\mathbb{H}$-bounded subsets of $O$,
  let
  $ 
    F_n \colon 
    W \times \R \times H' \times
    \mathbb{S}_{\mathbb{X}, \mathbb{X}'} \to \R
  $, 
  $ n \in \N_0 $, be degenerate elliptic
  functions,
  assume that
	for all 
  $ 
		(x_n,r_n,p_n,A_n)_{n \in \N_0} 
			\subseteq \mathcal U
	$
	and all $\delta \in (0, \infty)$ with
	$
		\lim_{n \to \infty} 
			d^{-, W}_{\mathbb{H}, \delta, h, u_0}(
				(x_0,r_0,p_0,A_0), (x_n,r_n,p_n,A_n)
			)
		= 0
	$
	and with 
	$
		\limsup_{n \to \infty}
			(F_n)^+_{\mathbb{H}, \mathbb{X}, \delta, h, u_0}
		(x_n,r_n,p_n,
	$
	$
		A_n) \leq 0
	$
	it holds that
  \begin{equation}  
  \label{eq:convergence F+}
		\begin{split}
					(F_0)^+_{\mathbb{H}, \mathbb{X}, \delta, h, u_0}(x_0,r_0,p_0,A_0)
			\leq
				\limsup_{ n \to \infty } \,
					(F_n)^+_{\mathbb{H}, \mathbb{X}, \delta, h, u_0}(x_n,r_n,p_n,A_n),
		\end{split}
  \end{equation}
	that for all $\tilde{x} \in O$ there exists an $R_{\tilde{x}} \in (0, \infty)$ 
	such that
	\begin{equation}
	\label{eq: uniform convergence of u_n}
		\lim_{R \to \infty}
		\limsup_{ n \to \infty }	
			\sup_{ x \in 
				\{ 
					y \in O \colon 
						\|y-\tilde{x}\|_H \leq R_{\tilde{x}}, 
						~\underline{h}_\mathbb{H}^W (y) \leq R 
				\} 
			}
				\left|
					u_n (x) - u_0(x) 
				\right|
		= 0,
	\end{equation}
	and that
	\begin{equation}
	\label{eq: uniform bound of u_n}
		\sup_{ n \in \N }	
			\sup_{ x \in \{ y \in O \colon \|y-\tilde{x}\|_H \leq R_{\tilde{x}} \} }
				\left(
					u_n (x) - u_n(x_0) 
				\right)
		<\infty,
	\end{equation}
  and assume 
  for every $ n \in \N $ that
  $ u_n $ is a viscosity subsolution
  of $ F_n = 0 $ relative to $(h, \mathbb{H}, \mathbb{X})$.
  Then $ u_0 $ is a viscosity subsolution
  of $ F_0 = 0 $ relative to $(h, \mathbb{H}, \mathbb{X})$.
 \end{lemma}

 \begin{proof}[Proof
  of
  Lemma~\ref{l:limits.of.viscosity.solutions}]
	If $ O = \emptyset $, 
  then the assertion is trivial. 
  So for the rest of 
  the proof, we assume that 
  $ O \neq \emptyset $.
  Then the proof is divided into
  two steps.
  {\it Step~1:}
	For the rest of Step 1 fix
  $ x_0 \in O $, 
	$\delta \in (0, \infty)$,
	and $r \in [0,1]$ 
	such that
  $
    r=
    \min\!\big( 1, R_{x_0}, 
      \frac{ 1 }{ 2 } 
      \dist_{\mathbb{H}}( x_0 ,
      H \backslash O) 
    \big) 
  $.
	Moreover, assume that there exist
	a 
	$
    \phi \in \C_{\mathbb{H}}^2(O,\R)
  $, 
	such that
	$\phi$ is bounded on $\mathbb{H}$-bounded subsets of $O$
  and such that $ x_0 $
  is a strict $\mathbb{H}$-maximum of 
  $ (u_0)^{-,W}_{\mathbb{H}, \delta, h} - \phi $.
	We obtain from the assumption 
	that for all
	$n \in \N_0$ it holds that
	$u_n \in \C_{\mathbb{H}}(O, \R)$ 
	that 
	for all $n \in \N_0$ 
	it holds that
	\begin{equation}
	\label{eq: u continuous}
			(u_n)^{-, W}_{\mathbb{H}, \delta, h} 
		=
			\overline{(u_n - \delta h)}_\mathbb{H}^W
		= 
			u_n - \delta \underline{h}_\mathbb{H}^W.
	\end{equation}
	This together with \eqref{eq:def of d delta without t} and the continuity
	of $u_n$ 
	with respect to the $\| \cdot \|_{\mathbb{H}}$-norm
	ensures that for all $n \in \N_0$ 
	$d^{+, W}_{\mathbb{H}, \delta, h, u_n}$  
	induce the
	same topology as
	$d^{+, W}_{\mathbb{H}, \delta, h, 0}$.
  Since $ O \subseteq H $ 
  is an open set
	and $O \neq \emptyset$,
  we have that $ r \in (0,1] $.
  In addition, the fact that
	$
		\{ 
      y \in H \colon
      \| y - x_0 \|_H \leq r
    \}
	$
	is convex, $\mathbb{H}$-bounded, and
	closed together with
	the Banach-Saks Theorem and with Lemma 5.1.4 in 
	Kato \cite{Kato1980}
	implies that it is weakly compact. Therefore it follows from page $3$
	in Stegall \cite{Stegall1978}  that the set
	$\{ x \in H \colon \| x\|_H \leq r \}$ is also an RNP
	set (for the definition see page $3$ in Stegall \cite{Stegall1978})
	and this together with
	boundedness of the function $ \phi $
  and of the functions 
  $ u_n $, $ n \in \N $, on $\mathbb{H}$-bounded subsets of $O$,
	the assumption that $h$ is bounded from below
	on $\mathbb{H}$-bounded subsets of $O$,
  upper semicontinuity 
	with respect to the $\| \cdot \|_{\mathbb{H}}$-norm
	of the function $ -\phi $
  and of the functions 
  $ (u_n)^{-, W}_{\mathbb{H}, \delta, h} $, $ n \in \N $,
  and with
	the theorem starting on page $4$ in Stegall \cite{Stegall1978}
  proves that for all $n \in \N$ 
  there exists a $p_n \in H$ and a 
  sequence 
  $ 
    x_n \in 
    \{ 
      y \in H \colon
      \| y - x_0 \|_H \leq r
    \}\subseteq O $,
  $ n \in \N $,
  of vectors such that
	$\|p_n\|_H \leq \frac 1n$
	and such that for all $n \in \N$
	it holds that
	$x_n$ is a strict $\mathbb{H}$-maximum of
  $ 
    \{ y \in H \colon \| y - x_0 \|_H \leq r \}
			\ni x \to 
				(u_n)^{-, W}_{\mathbb{H}, \delta, h}(x) - \phi(x) + \langle x, p_n \rangle_H
		\in \R \cup \{- \infty\}.
  $
	Thus \eqref{eq: u continuous} implies that for all $n \in \N$ it holds that
	 \begin{align}
		\nonumber
				&u_n(x_n) - \delta \underline{h}_\mathbb{H}^W (x_n) - \phi(x_n)
			=
				(u_n)^{-, W}_{\mathbb{H}, \delta, h}(x_n) - \phi(x_n)
					+ \langle x_n, p_n \rangle_H - \langle x_n, p_n \rangle_H \\
		\begin{split}
			\geq{}&
				(u_n)^{-, W}_{\mathbb{H}, \delta, h}(x_0) - \phi(x_0)
				+ \langle x_0, p_n \rangle_H - \langle x_n, p_n \rangle_H\\
			\geq{}&
				u_n(x_0) - \delta \underline{h}_\mathbb{H}^W (x_0) - \phi(x_0)
				- (\|x_0\|_H +\|x_n\|_H) \|p_n\|_H
		\end{split} \\ \nonumber
			\geq{}&
				u_n(x_0) - \delta \underline{h}_\mathbb{H}^W (x_0) - \phi(x_0)
				- \tfrac{2\|x_0\|_H +1}{n}
	\end{align}
	and this together with \eqref{eq: uniform bound of u_n}
	and the fact that $\phi$ is bounded on $\mathbb{H}$-bounded subsets of $O$
	shows that
	\begin{equation}
	\label{eq: h xn bound}
		\begin{split}
				&\sup_{n \in \N} \left( \delta \underline{h}_\mathbb{H}^W (x_n) \right)
			\leq
				\sup_{n \in \N} \left(
					u_n(x_n)- u_n(x_0)
					+\delta \underline{h}_\mathbb{H}^W (x_0)
					+\phi(x_0)-\phi(x_n) + \tfrac{2\|x_0\|_H +1}{n}
				\right) \\
			\leq{}&
				\sup_{n \in \N} (
					u_n(x_n)- u_n(x_0)
				)
				+\delta \underline{h}_\mathbb{H}^W (x_0)
				+\phi(x_0)
				-(
					\inf_{ x \in \{ y \in O \colon \|y-x_0\|_H \leq r \}} 
						\phi(x)
				) 
				+ 2\|x_0\|_H +1
				< \infty.\\
		\end{split}
	\end{equation}
  Next we prove that
  the sequence $ (x_n)_{ n \in \N } $
  converges to $ x_0 $. 
	Therefore note, that we get from \eqref{eq: u continuous} and
	from the fact that 
	$x_n$ is a strict $\mathbb{H}$-maximum of
  $ 
    \{ y \in H \colon \| y - x_0 \|_H \leq r \}
			\ni x \to 
				(u_n)^{-, W}_{\mathbb{H}, \delta, h}(x) - \phi(x) + \langle x, p_n \rangle_H
		\in \R \cup \{- \infty\}
  $
	that
  \begin{equation}
	\label{eq: equation u_0 and u_n}
		\begin{split}
				&(u_0)^{-, W}_{\mathbb{H}, \delta, h}(x_n) - \phi(x_n) \\
			={}&
				(u_0)^{-, W}_{\mathbb{H}, \delta, h}(x_n) 
				- (u_n)^{-, W}_{\mathbb{H}, \delta, h}(x_n) 
				+ (u_n)^{-, W}_{\mathbb{H}, \delta, h}(x_n) - \phi(x_n) 
				+ \langle x_n, p_n \rangle_H - \langle x_n, p_n \rangle_H \\
			\geq{}&
				u_0 (x_n) - \delta \underline{h}_\mathbb{H}^W(x_n)
				- u_n(x_n) + \delta \underline{h}_\mathbb{H}^W (x_n)
				+ (u_n)^{-, W}_{\mathbb{H}, \delta, h}(x_0) - \phi(x_0) 
				+ \langle x_0, p_n \rangle_H - \langle x_n, p_n \rangle_H \\
			={}&
				u_0 (x_n) - u_n(x_n) 
				+ (u_n)^{-, W}_{\mathbb{H}, \delta, h}(x_0) - \phi(x_0) 
				+ \langle x_0, p_n \rangle_H - \langle x_n, p_n \rangle_H.
		\end{split}
	\end{equation}
	Moreover, observe that from 
	\eqref{eq: equation u_0 and u_n},
	\eqref{eq: uniform convergence of u_n},
	\eqref{eq: h xn bound},
	\eqref{eq: u continuous},
	and from the facts that
	for all $n \in \N$ it holds that
	$x_n$, $x_0 \in \{ y \in H \colon \| y - x_0 \|_H \leq r \}$, that
	$r \leq 1$, and that $\|p_n\|_H \leq \frac 1n$ it follows that
	\begin{equation}
	\label{eq: xn strict max}
		\begin{split}
				&\liminf_{n \to \infty} (u_0)^{-, W}_{\mathbb{H}, \delta, h}(x_n) - \phi(x_n) \\
			\geq{}&
				 \liminf_{n \to \infty} \left (
					u_0(x_n) - u_n(x_n) 
					+ (u_n)^{-, W}_{\mathbb{H}, \delta, h}(x_0) - \phi(x_0) 
					+ \langle x_0, p_n \rangle_H - \langle x_n, p_n \rangle_H
				\right) \\
			\geq{}&
				 \lim_{n \to \infty} \left (
					(u_0)^{-, W}_{\mathbb{H}, \delta, h}(x_0) - \phi(x_0) 
					- \frac 1n (1 + 2 \|x_0\|_H)
				\right)
			=
				(u_0)^{-, W}_{\mathbb{H}, \delta, h}(x_0) - \phi(x_0).
		\end{split}
	\end{equation}
  Thus  
	the fact that $ x_0 $
  is a strict $\mathbb{H}$-maximum of 
  $ (u_0)^{-,W}_{\mathbb{H}, \delta, h} - \phi $ 
	together with 
	\eqref{eq: xn strict max}
	and with
	the continuity of $\phi$ 
	with respect to the $\| \cdot \|_{\mathbb{H}}$-norm
	implies that
  $ \lim_{ n \to \infty } \|x_n - x_0\|_H =0 $
	and that 
	\begin{equation}
			\lim_{n \to \infty} (u_n)^{-, W}_{\mathbb{H}, \delta, h}(x_n)
		=
			(u_0)^{-, W}_{\mathbb{H}, \delta, h}(x_0).
	\end{equation}
  The continuity of 
	$u_0$ 
	with respect to the $\| \cdot \|_{\mathbb{H}}$-norm,
	of
  $ 
    D_\mathbb{H} \phi \colon O
    \to H'
  $ 
	with respect to the $\| \cdot \|_{\mathbb{H}}$
	and the $\| \cdot \|_{H'}$-norm
  and of 
  $ 
    D^2_{\mathbb{H}} \phi \colon
    O \to \mathbb{S}_{\mathbb{H}, \mathbb{H}'}
  $
	with respect to the $\| \cdot \|_{\mathbb{H}}$
	and the $\| \cdot \|_{L(\mathbb{H},\mathbb{H}')}$-norm
	together with the fact that
	$\forall n \in \N \colon \|p_n\|_H \leq \frac 1n$
  hence implies that
  \begin{equation}
  \label{eq: limit}
		\begin{split}
			&\lim_{ n \to \infty }
				d^{-, W}_{\mathbb{H}, \delta, h, u_0}
					\Big(
						\big( 
							x_n, (u_n)^{-, W}_{\mathbb{H}, \delta, h}(x_n), 
							(D_{\mathbb{H}} \phi)(x_n) +p_n,  
							(D^2_{\mathbb{H}} \phi)(x_n)
						\big), 
			\\ &\qquad\qquad\qquad\qquad\qquad\qquad\qquad\qquad
						\big( 
							x_0, (u_0)^{-, W}_{\mathbb{H}, \delta, h}(x_0), 
							(D_{\mathbb{H}} \phi)(x_0),
							(D^2_{\mathbb{H}} \phi)(x_0)
						\big) 
					\Big)
			= 0.
		\end{split}
  \end{equation}
  In addition, 
  $ \lim_{ n \to \infty } \|x_n - x_0\|_H =0 $
	and the fact that
	for every $n \in \N$ it holds that
	$x_n$ is a strict $\mathbb{H}$-maximum of
  $ 
    \{ y \in H \colon \| y - x_0 \|_H \leq r \}
			\ni x \to 
				(u_n)^{-, W}_{\mathbb{H}, \delta, h}(x) - \phi(x) + \langle x, p_n \rangle_H
					\in \R \cup \{- \infty \}
  $
  show that
  there exists a 
  natural number $ n_0 \in \N $
  such that
  for all 
  $ n \in \{n_0, n_0 + 1, \ldots \} $
  it holds that
  $ \| x_n - x_0 \|_H < r $
  and that
  $ x_n \in O $ is a local maximum 
  of the function 
  $
		O \ni x \to 
			(u_n)^{-, W}_{\mathbb{H}, \delta, h}(x) 
			- \phi(x) + \langle x, p_n \rangle_H 
				\in \R \cup \{- \infty\}
	$.
  Hence, 
	it follows from
  Corollary~\ref{cor:semijets_equivalence},
	the fact that for all $n \in \N_0$ it holds that
	$d^{-, W}_{\mathbb{H}, \delta, h, u_n}$ induce the
	same topology as
	$d^{-, W}_{\mathbb{H}, \delta, h, 0}$,
  and from the assumption that for all $n \in \N$
	it holds that
  $ u_n $ is a
  viscosity subsolution of $ F_n = 0 $  
  that for all $ n \in \{n_0,n_0+1,\ldots\} $ it holds that
  \begin{equation}  
  \label{eq:viscosity0}
		\begin{split}
				&(F_n)^+_{\mathbb{H}, \mathbb{X}, \delta, h, u_0}\big( 
					x_n, 
					(u_n)^{-, W}_{\mathbb{H}, \delta, h}(x_n), 
					(D_{\mathbb{H}} \phi)(x_n) + p_n,
					( D^2_{\mathbb{H}} \phi)(x_n)
				\big) \\
			={}&
				(F_n)^+_{\mathbb{H}, \mathbb{X}, \delta, h, 0}\big( 
					x_n, 
					(u_n)^{-, W}_{\mathbb{H}, \delta, h}(x_n), 
					(D_{\mathbb{H}} \phi)(x_n) + p_n,
					( D^2_{\mathbb{H}} \phi)(x_n)
				\big) \\
			={}&
				(F_n)^+_{\mathbb{H}, \mathbb{X}, \delta, h, u_n}\big( 
					x_n, 
					(u_n)^{-, W}_{\mathbb{H}, \delta, h}(x_n), 
					(D_{\mathbb{H}} \phi)(x_n) + p_n,
					( D^2_{\mathbb{H}} \phi)(x_n)
				\big)
			\leq 0.
		\end{split}
  \end{equation}
	Assumption~\eqref{eq:convergence F+},
	Equation~\eqref{eq: limit}, and
	Inequality~\eqref{eq:viscosity0}
  therefore yield
  that
  \begin{equation}  
  \begin{split}
		&
			 (F_0)^+_{\mathbb{H}, \mathbb{X}, \delta, h, u_0}\big(
				 x_0, (u_0)^{-, W}_{\mathbb{H}, \delta, h}(x_0), (D_{\mathbb{H}} \phi )(x_0),  
				 (D^2_{\mathbb{H}} \phi)(x_0)
			\big) \\
    \leq{}&
			\limsup_{n \to \infty} \left(
				(F_n)^+_{\mathbb{H}, \mathbb{X}, \delta, h, u_0}\big( 
					x_n, 
					(u_n)^{-, W}_{\mathbb{H}, \delta, h}(x_n), 
					(D_{\mathbb{H}} \phi)(x_n) + p_n,
					( D^2_{\mathbb{H}} \phi)(x_n)
				\big)
			\right)
		\leq 
			0.
   \end{split}     
  \end{equation}   
  We thus have proved
  that
	for all $\delta \in (0, \infty)$,
  $ x \in O $, 
	and all
  $ 
    \phi \in 
    \{ 
      \psi \in \C_{\mathbb{H}}^2( O, \R ) 
      \colon
      x 
  $
  is a strict $\mathbb{H}$-maximum of 
  $
      ( (u_0)^{-, W}_{\mathbb{H}, \delta, h} - \psi ) \colon O \to \R 
  $
	and $\psi$ is bounded on $\mathbb{H}$-bounded subsets of $O \}$
	it holds that
	\begin{equation}
     (F_0)^+_{\mathbb{H}, \mathbb{X}, \delta, h, u_0}\big(
				 x_0, (u_0)^{-, W}_{\mathbb{H}, \delta, h}(x_0), (D_{\mathbb{H}} \phi )(x_0),  
				 (D^2_{\mathbb{H}} \phi)(x_0)
			\big)
    \leq 0.
  \end{equation}
  {\it Step~2:}
  For the rest of Step 2 fix 
  $ x_0 \in O $,
	$\delta \in (0, \infty)$,
  and
  $
    \phi \in \C_{\mathbb{H}}^2(O,\R)
  $ 
  such that
  $
    \phi( x_0 ) = (u_0)^{-, W}_{\mathbb{H}, \delta, h}(x_0)
  $
  and 
  $ \phi \geq (u_0)^{-, W}_{\mathbb{H}, \delta, h}(x_0) $.
  Next denote by 
  $ 
    \tilde{\phi} \colon O \to \R
  $
	the function satisfying for all $x \in O$
	that
  $
    \tilde{\phi}(x) = \phi(x) + \| x - x_0 \|_H^4
  $.
  Note that $\tilde{\phi} \in \C_{\mathbb{H}}(O, \R)$
	implies that $\tilde{\phi}$ is locally bounded.
	Therefore there exists an open set $\tilde{O} \subseteq O$
	such that $x_0 \in \tilde{O}$
	and such that $\tilde{\phi}$ is bounded on $\tilde{O}$.
	Moreover, it follows
	from the definition of $\phi$ and $\tilde{\phi}$
	that
  $ x_0 $ is
  a strict 
  $\mathbb{H}$-maximum of the function
  $ ( (u_0)^{-, W}_{\mathbb{H}, \delta, h}- \tilde{\phi} ) \colon \tilde{O} \to \R $.
  Step~1 can thus be applied 
	with $O \leftarrow \tilde{O}$, $\phi \leftarrow \tilde{\phi}|_{\tilde{O}}$ 
	and with $u_0 \leftarrow u_0 |_{\tilde{O}}$ to obtain
  \begin{equation}
  \label{eq:phi.eps.inequality}
		\begin{split}
				&(F_0)^+_{\mathbb{H}, \mathbb{X}, \delta, h, u_0}\big(
					x_0, (u_0)^{-, W}_{\mathbb{H}, \delta, h}(x_0), (D_{\mathbb{H}} \phi )(x_0),  
					(D^2_{\mathbb{H}} \phi)(x_0)
				\big)
			\leq 0.
		\end{split}
  \end{equation} 
  This and Corollary \ref{cor:semijets_equivalence}
	show that
  $ u_0 $ is a viscosity subsolution 
  of $ F_0 = 0 $ relative to $(h, \mathbb{H}, \mathbb{X})$.
  The proof 
  of 
  Lemma~\ref{l:limits.of.viscosity.solutions}
  is thus completed.
 \end{proof}
	\begin{remark}
		Lemma \ref{lem:sign_changes} implies that the corresponding result also holds for
		viscosity supersolutions.
	\end{remark}
\section{Chain rule for semijets}
The following lemma is used in the proof of Lemma \ref{lem: chain rule}
below. 
\begin{lemma}
\label{lem: orthogonal set}
	Let $\lambda \in (0,\infty)$,
	let $\mathbb{H}=(H, \langle \cdot, \cdot \rangle_H, \| \cdot \|_H)$
	be a real Hilbert space, 
	let $H_1 \subseteq H$ be a with respect to the
	$\|\cdot\|_H$-norm closed linear subset of $H$,
	let $g \in \C_{\mathbb{H} \times \mathbb{H}}(H \times H, \R)$ satisfy for all
	$x,y,z \in H$, $x_1\in H_1$, and all $t \in \R$ that
	\begin{equation}
	\label{eq: assumption g}
		\begin{split}
			&tg(x,y)=g(tx,y) = g(y,tx), \quad
			g(x+y,z) = g(x,z) + g(y,z), \quad
			g(x,x) \geq 0, \quad \\
			&g(x_1,x_1) \geq \lambda \|x_1\|^2_H,
		\end{split}
	\end{equation}
	and let $H_2 \subseteq H$ be the set satisfying that
	$H_2 = \{x \in H \colon (\forall x_1 \in H_1 \colon g(x,x_1)= 0)\}$. 
	Then $H_2$ is a with respect to the
	$\|\cdot\|_H$-norm closed linear subset of $H$,
	satisfying that
	$H_1 \cap H_2 = \{0\}$ and that $H_1 + H_2 = H$.
\end{lemma}
\begin{proof}[Proof of Lemma~\ref{lem: orthogonal set}]
	Note that \eqref{eq: assumption g} ensures that $g|_{H_1 \times H_1}$
	is a scalar product on $H_1$ and that $H_1$ is also closed with respect
	to the norm induced by $g|_{H_1 \times H_1}$. Therefore there
	exists an orthonormal basis $\mathfrak{A} \subseteq (H_1, g(\cdot, \cdot))$
	such that for all $\alpha, \beta \in \mathfrak{A}$ it holds that
	$g(\alpha,\beta) = \1_{\{\alpha\}}(\beta)$.
	Moreover we have for every $x \in H$ and every finite set
	$\mathfrak{B} \subseteq \mathfrak{A}$ that
	\begin{equation}
		\begin{split}
				&\sum_{\alpha \in \mathfrak{B}} (g(x,\alpha))^2
			=
				- \sum_{\alpha \in \mathfrak{B}} (g(x,\alpha))^2 
				+ 2\sum_{\alpha, \beta \in \mathfrak{B}} g(x,\alpha)g(x,\beta) g(\alpha,\beta) \\
			={}&
				g(x,x) 
				- g\Big(
						x- \sum_{\alpha \in \mathfrak{B}} g(x,\alpha) \alpha,
						x- \sum_{\alpha \in \mathfrak{B}} g(x,\alpha) \alpha
				\Big)
			\leq g(x,x).
		\end{split}
	\end{equation}
	This shows that for all $x \in H$ it holds that
	$\sum_{\alpha \in \mathfrak{A}} (g(x,\alpha))^2 < \infty$
	and, since $H_1$ is closed with respect to the norm induced by $g|_{H_1 \times H_1}$,
	it follows that
	$\sum_{\alpha \in \mathfrak{A}} g(x,\alpha) \alpha \in H_1$.
	In addition, it holds for all $x \in H$ and all $\beta \in \mathfrak{A}$
	that 
	\begin{equation}
			g \Big(x- \sum_{\alpha \in \mathfrak{A}} g(x,\alpha) \alpha ,\beta \Big)
		=
			g(x, \beta) - g(x,\beta) g(\beta ,\beta)
		= 
			0.
	\end{equation}
	Combining this with the fact that for all $x \in H$ it holds that
	$
			x 
		= 
			(x - \sum_{\alpha \in \mathfrak{A}} g(x,\alpha) \alpha)
			+\sum_{\alpha \in \mathfrak{A}} g(x,\alpha) \alpha
	$
	shows that $H_1 + H_2 =0$.
	Furthermore, $H_1 \cap H_2 =\{0\}$ follows directly from the definition
	of $H_2$ and the fact that for all $x_1 \in H_1$ it holds that
	$g(x_1,x_1) \geq \lambda \|x_1\|^2_H$.
	Finally let $(x^{(n)}_2)_{n \in \N} \subseteq H_2$ 
	be a Cauchy sequence with respect to the $\| \cdot \|_H$-norm.
	Since $\mathbb{H}$ is a Hilbert space we get that there exists an $x^{(0)} \in H$
	such that $\lim_{n \to \infty} \|x_2^{(n)} - x^{(0)}\|_H =0$.
	Moreover, we obtain from $H_1 + H_2 =H$ that there exists
	$x^{(0)}_1 \in H_1$ and $x^{(0)}_2 \in H_2$ such that
	$x^{(0)}_1 +x^{(0)}_2 = x^{(0)}$.
	Combining this with \eqref{eq: assumption g},
	the definition of $H_2$,
	and the continuity of $g$ yields that
	\begin{equation}
		\begin{split}
				0
			={}&
				\lim_{n \to \infty} g \big(x_2^{(n)} -x^{(0)},x_2^{(n)} -x^{(0)} \big)
			=
				\lim_{n \to \infty} 
					g\big(
						(x_2^{(n)} -x_2^{(0)}) -x_1^{(0)},(x_2^{(n)} -x_2^{(0)}) -x_1^{(0)}
					\big) \\
			={}&
				\lim_{n \to \infty} 
				\big(
					g\big(
						(x_2^{(n)} -x_2^{(0)}) ,(x_2^{(n)} -x_2^{(0)})
					\big)
					+ g(x_1^{(0)},x_1^{(0)})
				\big)
			\geq
				\lambda \|x^{(0)}_1\|^2_H.
		\end{split}
	\end{equation}
	Hence we have that $x^{(0)}_1 =0$ and thus that
	$x^{(0)} = x^{(0)}_2 \in H_2$. This proves that
	$H_2$ is closed with respect to the
	$\|\cdot\|_H$-norm and therefore finishes the proof of
	Lemma \ref{lem: orthogonal set}.
\end{proof}
With this lemma we can now generalize the ordinary chain rule to semijets.
To the best of our knowledge the chain rule for semijets is new.
\begin{lemma}[Chain rule for semijets]
\label{lem: chain rule}
	Let $\mathbb{H}=(H, \langle \cdot, \cdot \rangle_H, \| \cdot \|_H)$ and
	$\mathbb{V}=(V, \langle \cdot, \cdot \rangle_V, \| \cdot \|_V)$
	be real Hilbert spaces,
	let
	$\mathbb{H}'=(H', \| \cdot \|_{H'})$ and
	$\mathbb{V}'=(V', \| \cdot \|_{V'})$
	be their dual spaces,
	let
	$O \subseteq H$ and $U \subseteq V$ be open sets,
	and let $u \in \mathbb{M}(U, \R)$ and
	$f \in \C^2_{\mathbb{H}, \mathbb{V}}(O, U)$.
	Then
	\begin{itemize}
	\item[i)]
		for all 
		$x_0 \in O$, 
		$(\tilde{p}, \tilde{A}) \in H' \times \mathbb{S}_{\mathbb{H}, \mathbb{H}'}$
		satisfying that for all
		$\lambda \in (0,\infty)$
		there exist
		$p \in V'$ and
		$
			(A, B) 
				\in \mathbb{S}_{\mathbb{V}, \mathbb{V}'} \times \mathbb{S}_{\mathbb{H}, \mathbb{H}'}
		$
		such that
		$(p, A) \in (J_{\mathbb{V}, +}^2 u)(f(x_0))$,
		$B \geq -\lambda I_\mathbb{H} \pi^H_{\ker((D_{\mathbb{H}, \mathbb{V}} f)(x_0)) }$ 
		and that for all $x_1, x_2 \in H$ 
		it holds that
		\begin{align}
		\nonumber
					&\langle p, (D_{\mathbb{H}, \mathbb{V}} f)(x_0)x_1 \rangle_{V', V} 
				= 
					\langle \tilde{p}, x_1 \rangle_{H', H} \qquad \textrm{ and that} \\
					&\big \langle 
						A \big( (D_{\mathbb{H}, \mathbb{V}} f)(x_0) x_1 \big), 
						(D_{\mathbb{H}, \mathbb{V}} f)(x_0) x_2 
					\big \rangle_{V', V}
					+ \big \langle
							p, \big( \big( (D^2_{\mathbb{H}, \mathbb{V}} f)(x_0) \big) (x_1) \big) (x_2)
						\big \rangle_{V', V} \\
			\nonumber
					& \qquad \qquad\qquad \qquad\qquad \qquad\qquad \qquad \qquad \qquad
					+\langle x_1, B x_2 \rangle_{H, H'}
				=
					\langle x_1, \tilde{A} x_2 \rangle_{H, H'},
		\end{align}
		it holds that
		\begin{equation}
			(\tilde{p}, \tilde{A}) \in (J^2_{\mathbb{H}, +} (u \circ f)) (x_0),
		\end{equation}
		and
	\item[ii)]
		for all 
		$\lambda \in (0, \infty)$,
		$x_0 \in O$, 
		$(\tilde{p}, \tilde{A}) \in H' \times \mathbb{S}_{\mathbb{H}, \mathbb{H}'}$
		satisfying that
		$(\tilde{p}, \tilde{A}) \in (J_{\mathbb{H}, +}^2 (u \circ f))(x_0)$
		and
		that
		$(D_{\mathbb{H}, \mathbb{V}} f)(x_0) \colon H \to V$
		is surjective,
		there exist
		$B \in \mathbb{S}_{\mathbb{H}, \mathbb{H}'}$ and
		$(p, A) \in (J_{\mathbb{V}, +}^2 u) (f(x_0))$
		such that
		$B \geq -\lambda I_\mathbb{H} \pi^H_{\ker((D_{\mathbb{H}, \mathbb{V}} f)(x_0))}$
		and that for all $x_1, x_2 \in H$ it holds that
		\begin{align}
		\nonumber
					&\langle p, (D_{\mathbb{H}, \mathbb{V}} f)(x_0)x_1 \rangle_{V', V} 
				= 
					\langle \tilde{p}, x_1 \rangle_{H', H} \qquad \textrm{ and that} \\
					&\big \langle 
						A \big( (D_{\mathbb{H}, \mathbb{V}} f)(x_0)  x_1 \big), 
						(D_{\mathbb{H}, \mathbb{V}} f)(x_0) x_2 
					\big \rangle_{V', V}
					+ \big \langle
							p, \big( \big( (D^2_{\mathbb{H}, \mathbb{V}} f)(x_0) \big) (x_1) \big) (x_2)
						\big \rangle_{V', V} \\
				\nonumber
				& \qquad \qquad\qquad \qquad\qquad \qquad\qquad \qquad \qquad \qquad
					+ \langle x_1, B x_2 \rangle_{H, H'}
				=
					\langle x_1, \tilde{A} x_2 \rangle_{H, H'}.
		\end{align}
	\end{itemize}
\end{lemma}
\begin{proof}[Proof of Lemma \ref{lem: chain rule}] 
	If $ O = \emptyset $, 
  then the assertion is trivial. 
  So for the rest of 
  the proof, we assume that 
  $ O \neq \emptyset $.
	We start with the first part of the assertion.
	Therefore fix $x_0 \in O$ and 
	$(\tilde{p}, \tilde{A}) \in H' \times \mathbb{S}_{\mathbb{H}, \mathbb{H}'}$
	such that
	for every $\lambda \in (0, \infty)$
	there exist
	$(p_\lambda, A_\lambda) \in V' \times \mathbb{S}_{\mathbb{V}, \mathbb{V}'}$
	and $B_\lambda \in \mathbb{S}_{\mathbb{H}, \mathbb{H}'}$ such that
	$(p_\lambda, A_\lambda) \in (J_{\mathbb{V}, +}^2 u)(f(x_0))$,
	$B_\lambda \geq -\lambda I_\mathbb{H} \pi^H_{\ker((D_{\mathbb{H}, \mathbb{V}} f)(x_0))}$ 
	and such that
	for all $x_1, x_2 \in H$ it holds that
	\begin{align}
	\label{eq: chain rule formula}
	\nonumber
				&\langle p_\lambda, (D_{\mathbb{H}, \mathbb{V}} f)(x_0)x_1 \rangle_{V', V} 
			= 
				\langle \tilde{p}, x_1 \rangle_{H', H} \qquad \textrm{ and that} \\
				&\big \langle 
					A_\lambda \big( (D_{\mathbb{H}, \mathbb{V}} f)(x_0) x_1 \big), 
					(D_{\mathbb{H}, \mathbb{V}} f)(x_0) x_2 
				\big \rangle_{V', V}
				+ \big \langle
						p_\lambda, 
						\big( \big( (D^2_{\mathbb{H}, \mathbb{V}} f)(x_0) \big) (x_1) \big) (x_2)
					\big \rangle_{V', V} \\
		\nonumber
			& \qquad \qquad\qquad \qquad\qquad \qquad\qquad \qquad \qquad \qquad\qquad
				+ \langle x_1, B_\lambda x_2 \rangle_{H, H'}
			=
				\langle x_1, \tilde{A} x_2 \rangle_{H, H'}.
	\end{align}
	Next denote 
	by $H_1$ the set satisfying that
	$H_1 = \ker((D_{\mathbb{H}, \mathbb{V}} f)(x_0))$.
	From $f \in \C^2_{\mathbb{H}, \mathbb{V}}(O, U)$ it follows
	that there exists an $r_2 \in \mathbb{M}(O,U)$ such that
	\begin{equation}
	\label{eq: r2 convergence}
		\lim_{O \ni x \to x_0} 
			\frac{\| r_2(x) \|_V}{\|x-x_0\|_H^2} = 0
	\end{equation}
	and that for all $x \in O$ it holds that
	\begin{equation}
	\label{eq: taylor f}
			f(x) - f(x_0)
		=  
			\big( (D_{\mathbb{H}, \mathbb{V}}f)(x_0) \big)(x- x_0)
			+ \tfrac 12 \big( \big( (D^2_{\mathbb{H}, \mathbb{V}}f)(x_0) \big)(x-x_0) \big)(x-x_0)
			+ r_2(x).
	\end{equation}
	Furthermore, we denote by $r_1 \in \mathbb{M}(O,U)$ the function 
	satisfying for all $x \in O$ that
	\begin{equation}
	\label{eq: def of r_1}
			r_1(x) 
		= 
			r_2(x)
			+ \tfrac 12 \big( \big( (D^2_{\mathbb{H}, \mathbb{V}}f)(x_0) \big)(x-x_0) \big)(x-x_0).
	\end{equation}
	Then we get from \eqref{eq: chain rule formula}, \eqref{eq: taylor f}, 
	from \eqref{eq: def of r_1},
	and from the fact that 
	$
		\forall \lambda \in (0, \infty) 
			\colon B_\lambda \geq -\lambda I_{\mathbb{H}} \pi^H_{H_1}
	$ 
	that 
	for all $\lambda \in (0, \infty)$ it holds that 
	\begin{align}
		\nonumber
							&u(f(x)) - u(f(x_0)) 
							- \langle \tilde{p}, x - x_0 \rangle_{H', H}
							- \tfrac 12 \langle x - x_0, \tilde{A} (x- x_0) \rangle_{H, H'} \\ 
		\nonumber
				={}&
								u(f(x)) - u(f(x_0)) 
								- \langle 
										p_\lambda ,\big( (D_{\mathbb{H}, \mathbb{V}}f)(x_0) \big) (x - x_0) 
									\rangle_{V', V}
								- \tfrac 12 \langle (x-x_0), B_\lambda (x-x_0) \rangle_{H, H'} \\
		\nonumber
								&-\tfrac 12 
									\big \langle 
										\big( (D_{\mathbb{H}, \mathbb{V}}f)(x_0) \big) (x - x_0), 
										A_\lambda \big((D_{\mathbb{H}, \mathbb{V}}f)(x_0) \big) (x- x_0) 
									\big \rangle_{V, V'} \\
		\begin{split}
								&-\tfrac 12
									\big \langle
										p_\lambda, 
										\big( 
											\big( (D^2_{\mathbb{H}, \mathbb{V}} f)(x_0) \big) (x-x_0) 
										\big) (x-x_0)
									\big \rangle_{V', V} \\
				\leq{}&
								u(f(x)) - u(f(x_0)) 
								- \big \langle 
										p_\lambda, f(x) - f(x_0) -r_2(x) 
									\big \rangle_{V', V}
								+ \tfrac \lambda2 \|\pi^H_{H_1} (x-x_0)\|_H^2 
		\end{split}  \\
		\nonumber
								&-\tfrac 12 
									\big \langle 
										f(x)- f(x_0) -r_1(x), 
										A_\lambda (f(x)- f(x_0) -r_1(x)) 
									\big \rangle_{V, V'} \\
				\leq{}&
					\nonumber
								u(f(x)) - u(f(x_0)) 
								- \langle 
										p_\lambda, f(x) - f(x_0)
									\rangle_{V', V}
								-\tfrac 12 
									\big \langle 
										f(x)- f(x_0), 
										A_\lambda (f(x)- f(x_0) ) 
									\big \rangle_{V, V'} \\
					\nonumber
								&+\langle 
										p_\lambda, r_2(x)
								\rangle_{V', V}
								+ \big \langle 
									r_1(x), 
									A_\lambda (f(x)- f(x_0)) 
								\big \rangle_{V, V'}
								-\tfrac 12 
									\langle 
										r_1(x), 
										A_\lambda \, r_1(x) 
									\rangle_{V, V'}
					+ \tfrac {\lambda}{2}  \|x-x_0\|_H^2.
	\end{align}
	This together with 
	\eqref{eq: r2 convergence},
	the fact that 
	$
		\lim_{O \ni x \to x_0} 
			\frac{\| r_1(x) \|_V}{\|x-x_0\|_H} = 0
	$,
	and
	the assumption that for all
	$\lambda \in (0, \infty)$ it holds that
	$(p_\lambda, A_\lambda) \in (J_{\mathbb{V}, +}^2 u)(f(x_0))$
	shows that for all $\lambda \in (0, \infty)$
	it holds that
	\begin{align}
	\nonumber
				&\limsup_{O \ni x \to x_0}
					\frac
						{
							u(f(x)) - u(f(x_0)) 
							- \langle \tilde{p}, x - x_0 \rangle_{H', H}
							- \tfrac 12 \langle x - x_0, \tilde{A} (x- x_0) \rangle_{H, H'}
						}
						{ \| x - x_0 \|^2_H} \\
		\nonumber
				\leq{}&
					\limsup_{O \ni x \to x_0} \Bigg( \Bigg(
						\frac
							{
								u(f(x)) - u(f(x_0)) 
								- \langle 
										p_\lambda, f(x) - f(x_0)
									\rangle_{V', V}
							}
							{ \| f(x) - f(x_0) \|^2_V} 
						\\
		\nonumber
				& \qquad \qquad \qquad
					-\frac
							{
								\tfrac 12 
									\langle 
										f(x)- f(x_0), 
										A_\lambda (f(x)- f(x_0) ) 
									\rangle_{V, V'}
							}
							{ \| f(x) - f(x_0) \|^2_V} 
					\Bigg ) 
					\cdot \frac {\| f(x) - f(x_0) \|^2_V} {\| x - x_0 \|^2_H}
				\Bigg)
				+ \frac {\lambda}{2}\\
		\nonumber
			\leq{}&
				\limsup_{U \ni y \to f(x_0)} \Bigg( \Bigg(
						\frac
							{
								u(y) - u(f(x_0)) 
								- \langle 
										p_\lambda, y - f(x_0)
									\rangle_{V', V}
								-\tfrac 12 
									\langle 
										y- f(x_0), 
										A_\lambda (y- f(x_0) ) 
									\rangle_{V, V'}
							}
							{ \| y - f(x_0) \|^2_V} 
					\Bigg ) \\
			& \qquad \qquad \qquad
					\cdot \| (D_{\mathbb{H}, \mathbb{V}} f) (x_0) \|^2_{L(\mathbb{H}, \mathbb{V})}
				\Bigg) 
				+ \frac {\lambda}{2}
			\leq
				\frac {\lambda}{2}.
	\end{align}
	Letting $\lambda \downarrow 0$ then implies that
	\begin{equation}
		\begin{split}
				\limsup_{ O \ni x \to x_0}
					\frac
						{
							u(f(x)) - u(f(x_0)) 
							- \langle \tilde{p}, x - x_0 \rangle_{H', H}
							- \tfrac 12 \langle x - x_0, \tilde{A} (x- x_0) \rangle_{H, H'}
						}
						{ \| x - x_0 \|^2_H} 
			\leq 
				\lim_{\lambda \downarrow 0}
					\frac {\lambda}{2}
			=
				0.
		\end{split}
	\end{equation}
	Thus we obtain that 
	$(\tilde{p}, \tilde{A}) \in (J^2_{\mathbb{H}, +} (u \circ f)) (x_0)$
	which completes the proof of the first part.
	\newline
	For the second part of the assertion fix
	$\lambda \in (0, \infty)$,
	$x_0 \in O$, and 
	$(\tilde{p}, \tilde{A}) \in H' \times \mathbb{S}_{\mathbb{H}, \mathbb{H}'}$
	such that
	$(\tilde{p}, \tilde{A}) \in (J_{\mathbb{H}, +}^2 (u \circ f))(x_0)$,
	and that
	$(D_{\mathbb{H}, \mathbb{V}} f)(x_0) \colon H \to V$ is surjective.
	Next denote 
	by
	$\mathfrak{A}_1 \subseteq H_1$
	an orthonormal basis of 
	$
		(
			H_1, 
			\langle \cdot, \cdot \rangle_H|_{H_1^2}, 
			\| \cdot \|_H|_{H_1}
		)
	$,
	by $\mathfrak{B}_1 \subseteq H$ an orthonormal basis of $\mathbb{H}$ satisfying that
	$\mathfrak{A}_1 \subseteq \mathfrak{B}_1$,
	by $H_2 \subseteq H$ the set satisfying that
	$
			H_2 
		= 
			\overline{\Span_{\mathbb{H}} (\mathfrak{B}_1 \backslash \mathfrak{A}_1)}_\mathbb{H}
	$
	and by $\mathbb{H}_1, \mathbb{H}_2$ the Hilbert spaces satisfying that
	$
			\mathbb{H}_1 
		= 
			(
				H_1, 
				\langle \cdot, \cdot \rangle_H \, |_{H_1 \times H_1}, 
				\| \cdot \|_H \, |_{H_1}
			)
	$ 
	and that
	$
			\mathbb{H}_2 
		= 
			(
				H_2, 
				\langle \cdot, \cdot \rangle_H \, |_{H_2 \times H_2},
				\| \cdot \|_H \, |_{H_2}
			).
	$
	From this definition,
	the fact that $H_1 = \ker((D_{\mathbb{H}, \mathbb{V}} f)(x_0)$,
	and from the fact that 
	$(D_{\mathbb{H}, \mathbb{V}} f)(x_0) \colon H \to V$ is surjective
	it follows that
	$H_1 \cap H_2 = \{ 0 \}$,
	$H_1 + H_2 = H$,
	$H_1 \perp H_2$,
	and that
	$(D_{\mathbb{H}, \mathbb{V}} f)(x_0)|_{H_2} \colon H_2 \to V$ is bijective.
	Moreover, denote by 
	$F \colon H_1 \times H_2 \to V$ the function satisfying 
	for all $x_1 \in H_1$ and all $x_2 \in H_2$ that
	$F(x_1,x_2) = f(x_0 + x_1 +x_2)$.
	Note that 
	$
		H_2 \ni x_2 \to 
				\big( (D_{\mathbb{H}_1 \times \mathbb{H}_2, \mathbb{V}} F)(0,0) \big)(0,x_2)
			=
				\big( (D_{\mathbb{H}, \mathbb{V}} f)(x_0) \big) (x_2) \in V
	$
	is bijective and thus the implicit function theorem 
	for Hilbert spaces (see eg. page 417 in Edwards \cite{Edwards1994})
	yields that there
	exist a neighborhood $O_1 \subseteq H_1$ of $0$ and a function
	$h \in \C^2_{\mathbb{H}_1, \mathbb{H}_2}(O_1, H_2)$ such that
	$h(0) = 0$,
	that
	for all $x_1 \in O_1$ it holds that
	\begin{equation}
	\label{eq: def of h chain}
		 f(x_0)=F(0,0)=F(x_1, h(x_1)) = f(x_0 + x_1 +h(x_1)),
	\end{equation}
	that for all $x_1 \in H_1$ it holds that
	\begin{equation}
	\label{eq: implicit function first derivative}
			-\big( (D_{\mathbb{H}, \mathbb{V}} f)(x_0) \big) (x_1) 
		= 
			\big( (D_{\mathbb{H}, \mathbb{V}} f)(x_0) \big) 
				\big( \big( (D_{\mathbb{H}_1, \mathbb{H}_2} h)(0) \big) (x_1) \big),
	\end{equation}
	and that for all $x_1, \tilde{x}_1 \in O_1$ it holds that
	\begin{equation}
	\label{eq: implicit function second derivative}
		\begin{split}
				&-\Big( \big( (D^2_{\mathbb{H}, \mathbb{V}} f)(x_0) \big) 
					\big( x_1 + \big( (D_{\mathbb{H}_1, \mathbb{H}_2} h)(0) \big) (x_1) \big) 
				\Big) 
					\big( 
						\tilde{x}_1 + \big( (D_{\mathbb{H}_1, \mathbb{H}_2} h)(0) \big) (\tilde{x}_1) 	
					\big) \\
			={}& 
				\Big( (D_{\mathbb{H}, \mathbb{V}} f)(x_0) \Big) 
					\Big( 
						\big( 
							\big( (D^2_{\mathbb{H}_1, \mathbb{H}_2} h)(0) \big) (x_1) 
						\big) (\tilde{x}_1) 
					\Big).
		\end{split}
	\end{equation}
	In addition, \eqref{eq: implicit function first derivative} together with the
	fact that $H_1 = \ker((D_{\mathbb{H}, \mathbb{V}} f)(x_0)$ and that
	$(D_{\mathbb{H}_1, \mathbb{H}_2} h)(0) \in L(\mathbb{H}_1, \mathbb{H}_2)$ shows that
	for all $x_1 \in H_1$ it holds that
	$
		\big( (D_{\mathbb{H}_1, \mathbb{H}_2} h)(0) \big) (x_1) \in H_1 \cap H_2
	$
	and combining this with the fact that $H_1 \cap H_2 = \{ 0\}$
	yields then
	\begin{equation}
	\label{eq: h first derivative}
		(D_{\mathbb{H}_1, \mathbb{H}_2} h)(0) = 0.
	\end{equation}
	Furthermore, we obtain from \eqref{eq: implicit function second derivative} 
	and from \eqref{eq: h first derivative} that for all
	$x_1, \tilde{x}_1 \in H_1$ it holds that
	\begin{equation}
		\begin{split}
				-\big( \big( (D^2_{\mathbb{H}, \mathbb{V}} f)(x_0) \big) 
					(x_1) 
				\big) (\tilde{x}_1) 
			= 
				\Big( (D_{\mathbb{H}, \mathbb{V}} f)(x_0) \Big) 
					\Big( 
						\big( 
							\big( (D^2_{\mathbb{H}_1, \mathbb{H}_2} h)(0) \big) (x_1) 
						\big) (\tilde{x}_1) 
					\Big)
		\end{split}
	\end{equation}
	and this with the fact that $(D_{\mathbb{H}, \mathbb{V}} f)(x_0)|_{H_2}$ is bijective
	shows that
	for all $x_1, \tilde{x}_1 \in H_1$ it holds that
		\begin{equation}
		\label{eq: h second derivative}
				\big( 
					\big( (D^2_{\mathbb{H}_1, \mathbb{H}_2} h)(0) \big) (x_1) 
				\big) (\tilde{x}_1)
			=
				-((D_{\mathbb{H}, \mathbb{V}} f)(x_0) |_{H_2} )^{-1}
				\Big(
					\big( \big( (D^2_{\mathbb{H}, \mathbb{V}} f)(x_0) \big) 
						(x_1) 
					\big) (\tilde{x}_1)
				\Big).
	\end{equation}
	Moreover, 
	\eqref{eq: h first derivative} together with
	\eqref{eq: h second derivative},
	with
	$h(0) = 0$,
	and with
	$h \in \C^2_{\mathbb{H}_1, \mathbb{H}_2}(O_1, H_2)$
	implies that 
	\begin{equation}
	\label{eq: h convergence in 0}
		\lim_{O_1 \ni x_1 \to 0} 
			\frac{\| h(x_1) \|_H}{\|x_1\|_H} = 0
	\end{equation}
	and that there exists a function 
	$\tilde{r}_2 \in \mathbb{M}(O_1,H_2)$ such that
	$
		\lim_{O_1 \ni x_1 \to 0} 
			\frac{\| \tilde{r}_2(x_1) \|_H}{\|x_1\|_H^2} = 0
	$
	and that for all $x_1 \in O_1$ it holds that
	\begin{equation}
	\label{eq: taylor h}
			h(x_1)
		=  
			-\frac 12
			((D_{\mathbb{H}, \mathbb{V}} f)(x_0) |_{H_2} )^{-1}
				\Big(
					\big( \big( (D^2_{\mathbb{H}, \mathbb{V}} f)(x_0) \big) 
						(x_1) 
					\big) (x_1)
				\Big)
			+ \tilde{r}_2(x_1).
	\end{equation}
	In addition, \eqref{eq: h convergence in 0} yields that
	\begin{equation}
	\label{eq: changing quotient}
		\limsup_{O_1 \ni x_1 \to 0}
			\frac
				{\| x_1 \|_H^2}
				{\|x_1 + h(x_1) \|_H^2}
		= 
			\begin{cases}
				1 &\textrm{ if } H \neq \{0\} \\
				0 &\textrm{ if } H = \{0\}.
			\end{cases}
	\end{equation}
	Combining now \eqref{eq: def of h chain},
	\eqref{eq: h convergence in 0},
	\eqref{eq: taylor h},
	\eqref{eq: changing quotient},
	the fact that
	$
		\lim_{O_1 \ni x_1 \to 0} 
			\frac{\| \tilde{r}_2(x_1) \|_H}{\|x_1\|_H^2} = 0
	$
	and the fact that
	$(\tilde{p}, \tilde{A}) \in (J_{\mathbb{H}, +}^2 (u \circ f))(x_0)$ shows then
	\begin{align}
	 \nonumber
				0
			\geq&{}
				\limsup_{O \ni x \to x_0}
					\frac
						{
							u(f(x)) - u(f(x_0)) 
							- \langle \tilde{p}, x - x_0 \rangle_{H', H}
							- \tfrac 12 \langle x - x_0, \tilde{A} (x- x_0) \rangle_{H, H'}
						}
						{ \| x - x_0 \|^2_H} \\
		\nonumber
			\geq&{}
				\limsup_{O_1 \ni x_1 \to 0} \bigg(
					\frac
						{
							u(f(x_0 + x_1+ h(x_1))) - u(f(x_0)) 
							- \langle \tilde{p}, x_1+ h(x_1) \rangle_{H', H}
						}
						{ \| x_1+ h(x_1)\|^2_H} \\
		\nonumber
			& \qquad \qquad\qquad 
					-\frac
						{
							\tfrac 12 \langle x_1+ h(x_1), \tilde{A} (x_1+ h(x_1)) \rangle_{H, H'}
						}
						{ \| x_1+ h(x_1)\|^2_H} \bigg) \\
			={}&
				\limsup_{O_1 \ni x_1 \to 0}
					\Bigg(
						\Bigg(
							\frac
								{
									u(f(x_0)) - u(f(x_0)) 
									- \tfrac 12 \langle x_1+ h(x_1), \tilde{A} (x_1+ h(x_1)) \rangle_{H, H'}
								}
								{ \| x_1\|^2_H} \\
			\nonumber
						& \qquad \qquad \qquad 
							-\frac
								{
									\Big \langle 
											\tilde{p}, 
											x_1 
											-\tfrac 12
												((D_{\mathbb{H}, \mathbb{V}} f)(x_0) |_{H_2} )^{-1}
													\Big(
														\big( 
															\big( (D^2_{\mathbb{H}, \mathbb{V}} f)(x_0) \big) (x_1) 
														\big) (x_1)
													\Big)
											+ \tilde{r}_2(x_1) 
										\Big \rangle_{H', H}
								}
								{ \| x_1\|^2_H} 
					\Bigg) 
			\\& \nonumber \qquad \qquad \qquad 
					\cdot \frac 
						{ \| x_1\|^2_H}
						{ \| x_1+ h(x_1)\|^2_H}
				\Bigg)\\
		\nonumber
			={}&
				\limsup_{O_1 \ni x_1 \to 0} \bigg(
					\frac
						{ 
							- \langle 
									\tilde{p}, 
									x_1 
								\rangle_{H', H}
							+	\Big \langle 
									\tilde{p},
									\frac 12
										((D_{\mathbb{H}, \mathbb{V}} f)(x_0) |_{H_2} )^{-1}
											\Big(
												\big( 
													\big( (D^2_{\mathbb{H}, \mathbb{V}} f)(x_0) \big) (x_1) 
												\big) (x_1)
											\Big)
								\Big \rangle_{H', H}
						}
						{ \| x_1\|^2_H} \\
		\nonumber
				& \qquad \qquad \qquad
					-\frac
						{
							\tfrac 12 \langle x_1, \tilde{A} (x_1) \rangle_{H, H'}
						}
						{ \| x_1\|^2_H}
			\bigg).
	\end{align}
	Thus we get for every $x_1 \in H_1$ with $\|x_1 \|_H =1$ that
	\begin{equation}
		\begin{split}
				0
			\geq{}&
				\limsup_{t \downarrow 0} \bigg(
					\frac
						{ 
							- \langle 
									\tilde{p}, 
									t x_1 
								\rangle_{H', H}
							+	\Big \langle 
									\tilde{p},
									\tfrac 12
										((D_{\mathbb{H}, \mathbb{V}} f)(x_0) |_{H_2} )^{-1}
											\Big(
												\big( 
													\big( (D^2_{\mathbb{H}, \mathbb{V}} f)(x_0) \big) (t x_1) 
												\big) (t x_1)
											\Big)
								\Big \rangle_{H', H}
						}
						{ \| x_1\|^2_H} \\
				& \qquad \qquad \qquad \qquad
					-\frac
						{
							\tfrac 12 \langle t x_1, \tilde{A} (t x_1) \rangle_{H, H'}
						}
						{ t^2 \| x_1\|^2_H} \bigg)\\
			={}&
				\limsup_{t \downarrow 0}
					\Bigg(
						\frac
							{ 
								- \langle 
										\tilde{p}, 
										x_1 
									\rangle_{H', H}
							}
							{t}
						\Bigg)
				+	\frac 12
					\Big \langle 
						\tilde{p},
						((D_{\mathbb{H}, \mathbb{V}} f)(x_0) |_{H_2} )^{-1}
							\Big(
								\big( 
									\big( (D^2_{\mathbb{H}, \mathbb{V}} f)(x_0) \big) (x_1) 
								\big) (x_1)
							\Big)
					\Big \rangle_{H', H}
					\\& 
				- \tfrac 12 
				\big \langle x_1, \tilde{A} (x_1) \big \rangle_{H, H'}
		\end{split}
	\end{equation}
	and this shows that for every $x_1 \in H_1$ with $\|x_1 \|_H =1$ it holds that
	$
		- \langle \tilde{p}, x_1 \rangle_{H', H} =0
	$ 
	and that
	$
		\langle x_1, \tilde{A} (x_1) \rangle_{H, H'}
		-\Big \langle 
			\tilde{p},
				((D_{\mathbb{H}, \mathbb{V}} f)(x_0) |_{H_2} )^{-1}
					\Big(
						\big( 
							\big( (D^2_{\mathbb{H}, \mathbb{V}} f)(x_0) \big) (x_1) 
						\big) (x_1)
					\Big)
		\Big \rangle_{H', H}
		\geq 0.
	$
	Therefore we obtain that for all $x_1 \in H_1$ it holds that
	\begin{equation}
	\label{eq: tilde p nice}
		\langle \tilde{p}, x_1 \rangle_{H', H} =0
	\end{equation} 
	and that
	\begin{equation}
	\label{eq: tilde A nice}
		\langle x_1, \tilde{A} (x_1) \rangle_{H, H'}
		-\Big \langle 
			\tilde{p},
				((D_{\mathbb{H}, \mathbb{V}} f)(x_0) |_{H_2} )^{-1}
					\Big(
						\big( 
							\big( (D^2_{\mathbb{H}, \mathbb{V}} f)(x_0) \big) (x_1) 
						\big) (x_1)
					\Big)
		\Big \rangle_{H', H}
		\geq 0.
	\end{equation}
	Next denote by $p \in V'$ the function satisfying for
	all $x_1 \in H$ that
	\begin{equation}
	\label{eq: def of p}
				\langle p, (D_{\mathbb{H}, \mathbb{V}} f)(x_0)x_1 \rangle_{V', V} 
			= 
				\langle \tilde{p}, x_1 \rangle_{H', H},
	\end{equation}
	by
	$g \colon \mathbb{M}(\mathbb{V}, \mathbb{V}') \times H \times H \to \R$
	the function satisfying for all 
	$C \in \mathbb{M}(\mathbb{V}, \mathbb{V}')$
	and all $x_1$, $x_2 \in H$ that
	\begin{equation}
	\label{eq: def of g}
		\begin{split}
				g(C,x_1, x_2) 
			={}&
				\langle x_1, \tilde{A} x_2 \rangle_{H, H'}
				- \big \langle 
						C \big( (D_{\mathbb{H}, \mathbb{V}} f)(x_0) x_1 \big), 
						(D_{\mathbb{H}, \mathbb{V}} f)(x_0) x_2 
					\big \rangle_{V', V} \\
				&- \big \langle
						p, \big( \big( (D^2_{\mathbb{H}, \mathbb{V}} f)(x_0) \big) (x_1) \big) (x_2)
					\big \rangle_{V', V}
				+ \lambda \langle \pi^H_{H_1} x_1, x_2 \rangle_{H}
		\end{split}
	\end{equation}
	and 
	by
	$S \subseteq \mathbb{S}_{\mathbb{V}, \mathbb{V}'}$ the set
	satisfying that 
	\begin{equation}
	\label{eq: def of S}
		\begin{split}
				S
			= 
				\{ 
					C \in \mathbb{S}_{\mathbb{V}, \mathbb{V}'} \colon 
						(\forall x \in H \colon
							g(C, x, x)
							\geq 0)
				\}.
		\end{split}
	\end{equation}
	Note that \eqref{eq: tilde p nice} together with the fact that
	$(D_{\mathbb{H}, \mathbb{V}} f)(x_0) \colon H \to V$ 
	is surjective ensures that
	$p$ is well defined and
	note that we have for all $C \in \mathbb{S}_{\mathbb{V}, \mathbb{V}'}$ 
	and all $x_1, x_2, x_3, x_4 \in H$ that
	$
			g(C, x_1+x_2, x_3+ x_4) 
		= 
			g(C,x_1, x_3) + g(C,x_1, x_4) + g(C,x_2, x_3) + g(C,x_2, x_4).
	$
	In addition we get from \eqref{eq: tilde A nice},
	\eqref{eq: def of p},
	\eqref{eq: def of g} and the fact that
	$H_1 = \ker((D_{\mathbb{H}, \mathbb{V}} f)(x_0))$ that for all
	$x_1 \in H_1$ and all $C \in \mathbb{M}(\mathbb{V}, \mathbb{V}')$ 
	it holds that
	\begin{equation}
	\label{eq: g bounded from below on H_1}
			g(C,x_1,x_1) 
		\geq \lambda \|x_1\|_H^2.
	\end{equation}
	Furthermore, it follows from the fact that 
	$(D_{\mathbb{H}, \mathbb{V}} f)(x_0)|_{H_2}$ is bijective
	that for all $x_2 \in H_2$ it holds that
	\begin{equation}
	\label{eq: Df lower bound}
		\begin{split}
				&\|
					((D_{\mathbb{H}, \mathbb{V}} f)(x_0)) (x_2)
				\|_{V} \\
			={}
				&\|
					((D_{\mathbb{H}, \mathbb{V}} f)(x_0)|_{H_2})^{-1} 
				\|^{-1}_{L(\mathbb{V},\mathbb{H}_2)}
				\|
					((D_{\mathbb{H}, \mathbb{V}} f)(x_0)|_{H_2})^{-1} 
				\|_{L(\mathbb{V},\mathbb{H}_2)}
				\|
					((D_{\mathbb{H}, \mathbb{V}} f)(x_0)) (x_2)
				\|_{V} \\
			\geq{}&
				\|
					((D_{\mathbb{H}, \mathbb{V}} f)(x_0)|_{H_2})^{-1} 
				\|^{-1}_{L(\mathbb{V},\mathbb{H}_2)}
				\|
					((D_{\mathbb{H}, \mathbb{V}} f)(x_0)|_{H_2})^{-1} 
						\big(((D_{\mathbb{H}, \mathbb{V}} f)(x_0)|_{H_2}) (x_2) \big)
				\|_{H} \\
			={}&
				\|
					((D_{\mathbb{H}, \mathbb{V}} f)(x_0)|_{H_2})^{-1} 
				\|^{-1}_{L(\mathbb{V},\mathbb{H}_2)}
				\|x_2\|_{H}.
		\end{split}
	\end{equation}
	We now show with Zorn's lemma that 
	$S$ contains a maximal element.
	First observe that we get from 
	\eqref{eq: tilde A nice},
	\eqref{eq: def of p},
	\eqref{eq: def of g},
	\eqref{eq: Df lower bound},
	and from the fact that $H_1 \perp H_2$ that
	for all $x_1 \in H_1$, $x_2 \in H_2$, $R, \Lambda \in (0, \infty)$
	with 
	$
			R
		= 
			\|\tilde{A}\|_{L(\mathbb{H}, \mathbb{H}')}
			+ \| p \|_{V'} 
				\|
					(D^2_{\mathbb{H}, \mathbb{V}} f)(x_0)
				\|_{L(\mathbb{H}, \mathbb{L}(\mathbb{H}, \mathbb{V}))}
	$
	and with
	$
		\Lambda
			\|
				((D_{\mathbb{H}, \mathbb{V}} f)(x_0)|_{H_2})^{-1} 
			\|^{-2}_{L(\mathbb{V},\mathbb{H}_2)}
		\geq
			\frac{R^2}{\lambda} + R
	$
	it holds that
	\begin{align}
	\nonumber
				&g(- \Lambda I_\mathbb{V}, x_1 + x_2, x_1 + x_2)
			=
				g(- \Lambda I_\mathbb{V}, x_1, x_1)
				+2 g(- \Lambda I_\mathbb{V}, x_1 , x_2)
				+g(- \Lambda I_\mathbb{V}, x_2, x_2) \\
	\nonumber
			={}&
				\langle x_1, \tilde{A} x_1 \rangle_{H, H'}
				- \big \langle
						p, \big( \big( (D^2_{\mathbb{H}, \mathbb{V}} f)(x_0) \big) (x_1) \big) (x_1)
					\big \rangle_{V', V}
				+ 
				2 \langle x_1, \tilde{A} x_2 \rangle_{H, H'} \\
	\nonumber
				&- 2\big \langle
						p, \big( \big( (D^2_{\mathbb{H}, \mathbb{V}} f)(x_0) \big) (x_1) \big) (x_2)
					\big \rangle_{V', V} 
				+\langle x_2, \tilde{A} x_2 \rangle_{H, H'}
				+ \Lambda \big \langle 
						(D_{\mathbb{H}, \mathbb{V}} f)(x_0) x_2, 
						(D_{\mathbb{H}, \mathbb{V}} f)(x_0) x_2 
					\big \rangle_{V} \\
	\nonumber
				&- \big \langle
						p, \big( \big( (D^2_{\mathbb{H}, \mathbb{V}} f)(x_0) \big) (x_2) \big) (x_2)
					\big \rangle_{V', V} 
				+ \lambda \|x_1\|_H^2 \\
			\geq{}&
				\langle x_1, \tilde{A} x_1 \rangle_{H, H'}
				- \Big \langle 
						\tilde{p},
						((D_{\mathbb{H}, \mathbb{V}} f)(x_0) |_{H_2} )^{-1}
							\Big(
								\big( 
									\big( (D^2_{\mathbb{H}, \mathbb{V}} f)(x_0) \big) (x_1) 
								\big) (x_1)
							\Big)
					\Big \rangle_{H', H} \\
	\nonumber
				&- 2 \|x_1\|_H \|x_2 \|_H 
						(
							\|\tilde{A}\|_{L(\mathbb{H}, \mathbb{H}')}
							+ \| p \|_{V'} 
								\|
									(D^2_{\mathbb{H}, \mathbb{V}} f)(x_0)
								\|_{L(\mathbb{H}, \mathbb{L}(\mathbb{H}, \mathbb{V}) )}
						) \\
	\nonumber
				&- \|x_2\|^2_H
						\big(
							\|\tilde{A}\|_{L(\mathbb{H}, \mathbb{H}')}
							+ \| p \|_{V'} 
								\|
									(D^2_{\mathbb{H}, \mathbb{V}} f)(x_0)
								\|_{L(\mathbb{H}, \mathbb{L}(\mathbb{H}, \mathbb{V}) )}
						\big) \\
	\nonumber
					&+ \Lambda 
						\|
							((D_{\mathbb{H}, \mathbb{V}} f)(x_0)|_{H_2})^{-1} 
						\|^{-2}_{L(\mathbb{V},\mathbb{H}_2)}
						\|x_2\|^2_{H} + \lambda \|x_1\|^2_{H}\\
	\nonumber
				\geq{}&
					\lambda
					\Big(
						\|x_1\|_H
						-\frac{R}{\lambda} \|x_2\|_H
					\Big)^2
					+\left(
						\Lambda
						\|
							((D_{\mathbb{H}, \mathbb{V}} f)(x_0)|_{H_2})^{-1} 
						\|^{-2}_{L(\mathbb{V},\mathbb{H}_2)}
						-\frac{R^2}{\lambda}
						-R
					\right) \|x_2\|^2_H
				\geq{}
					0.
	\end{align}
	This together with $H_1 + H_2 = H$ shows that
	$S \neq \emptyset$.
	Moreover, let $P \subseteq S$ be a totally ordered subset. We denote by
	$C_P \in \mathbb{M}(\mathbb{V}, \mathbb{V}')$ the function satisfying
	for all $y_1$, $y_2 \in V$ that
	$
			\langle C_P y_1, y_2 \rangle_{V', V}
		=
			\frac 14 (
				\sup_{C \in P} (\langle C(y_1+ y_2), y_1+y_2 \rangle_{V', V}) 
				- \sup_{C \in P} (\langle C(y_1-y_2), y_1-y_2 \rangle_{V', V})
			).
	$
	Note, that we get from \eqref{eq: def of g} and from
	\eqref{eq: def of S} that for all $C \in S$ and all
	$x \in H$ it holds that
	\begin{equation}
		\begin{split}
				&\big \langle 
					C \big( (D_{\mathbb{H}, \mathbb{V}} f)(x_0) x \big), 
					(D_{\mathbb{H}, \mathbb{V}} f)(x_0) x 
				\big \rangle_{V', V} \\
			\leq{}&
				\langle x, \tilde{A} x \rangle_{H, H'}
				- \big \langle
						p, \big( \big( (D^2_{\mathbb{H}, \mathbb{V}} f)(x_0) \big) (x) \big) (x)
					\big \rangle_{V', V}
				+ \lambda \langle \pi^H_{H_1} x, x \rangle_{H}.
		\end{split}
	\end{equation}
	Combining this with the surjectivity of 
	$(D_{\mathbb{H}, \mathbb{V}} f)(x_0) \colon H \to V$
	ensures that $C_P$ is well defined.
	Moreover, 
	\eqref{eq: def of g} together with
	\eqref{eq: def of S} and with
	$P \subseteq S$ shows for all $x \in H$ that
	\begin{equation}
		\begin{split}
				g(C_P,x, x) 
			=
				\inf_{C \in P}
					g(C, x, x)
			\geq 
				0.
		\end{split}
	\end{equation}
	This together with \eqref{eq: def of g},
	and the fact that for all $C \in P$ it holds that $C_P \geq C$
	ensures that for all $x \in H$ and all $C \in P$ it holds that
	\begin{equation}
		\begin{split}
				&\big \langle 
					C \big( (D_{\mathbb{H}, \mathbb{V}} f)(x_0) x \big), 
					(D_{\mathbb{H}, \mathbb{V}} f)(x_0) x 
				\big \rangle_{V', V} 
			\leq
				\big \langle 
					C_P \big( (D_{\mathbb{H}, \mathbb{V}} f)(x_0) x \big), 
					(D_{\mathbb{H}, \mathbb{V}} f)(x_0) x 
				\big \rangle_{V', V} \\
			\leq{}&
				\langle x, \tilde{A} x \rangle_{H, H'}
				- \big \langle
						p, \big( \big( (D^2_{\mathbb{H}, \mathbb{V}} f)(x_0) \big) (x) \big) (x)
					\big \rangle_{V', V}
				+ \lambda \langle \pi^H_{H_1} x, x \rangle_{H}
		\end{split}
	\end{equation}
	and thus the assumption that 
	$(D_{\mathbb{H}, \mathbb{V}} f)(x_0) \colon H \to V$ 
	is surjective implies that
	\begin{equation}
		\sup_{y \in \{ \tilde{y} \in V \colon \|\tilde{y}\|_V \leq 1\}} 
			|\langle C_P y, y \rangle_{V',V}| < \infty.
	\end{equation}
	In addition, it follows directly from the definition that
	$C_P$ is symmetric and therefore
	it remains to prove that $C_P$ is linear to conclude that $C_P \in S$.
	Let $\overline{y}_1$, $\overline{y}_2$, and $\overline{y}_3 \in V$ 
	and denote for every $i \in \{1, 2, 3, 4, 5, 6\}$
	by $(C_{i,n})_{n \in \N} \subseteq P$ sequences satisfying that
	\begin{align}
			\lim_{n \to \infty}
				\langle 
					(\overline{y}_1+\overline{y}_2+\overline{y}_3), 
					(C_{1,n})(\overline{y}_1+\overline{y}_2+\overline{y}_3)
				\rangle_{V, V'} 
		&=
			\sup_{C \in P}
				\langle 
					(\overline{y}_1+\overline{y}_2+\overline{y}_3), 
					C(\overline{y}_1+\overline{y}_2+\overline{y}_3) 
				\rangle_{V, V'}, \\
			\lim_{n \to \infty}
				\langle 
					(\overline{y}_1-\overline{y}_2-\overline{y}_3), 
					(C_{2,n})(\overline{y}_1-\overline{y}_2-\overline{y}_3) 
				\rangle_{V, V'} 
		&=
			\sup_{C \in P}
				\langle 
					(\overline{y}_1-\overline{y}_2-\overline{y}_3), 
					C(\overline{y}_1-\overline{y}_2-\overline{y}_3) 
				\rangle_{V, V'}, \\
			\lim_{n \to \infty}
				\langle 
					(\overline{y}_1+\overline{y}_2), 
					(C_{3,n})(\overline{y}_1+\overline{y}_2) 
				\rangle_{V, V'}
		&=
			\sup_{C \in P}
				\langle 
					(\overline{y}_1+\overline{y}_2),
					C(\overline{y}_1+\overline{y}_2) 
				\rangle_{V, V'}, \\
			\lim_{n \to \infty}
				\langle
					(\overline{y}_1-\overline{y}_2), 
					(C_{4,n})(\overline{y}_1-\overline{y}_2) 
				\rangle_{V, V'}
		&=
			\sup_{C \in P}
				\langle 
					(\overline{y}_1-\overline{y}_2), 
					C(\overline{y}_1-\overline{y}_2) 
				\rangle_{V, V'}, \\
			\lim_{n \to \infty}
				\langle 
					(\overline{y}_1+\overline{y}_3),
					(C_{5,n})(\overline{y}_1+\overline{y}_3) 
				\rangle_{V, V'}
		&=
			\sup_{C \in P}
				\langle 
					(\overline{y}_1+\overline{y}_3), 
					C(\overline{y}_1+\overline{y}_3) 
				\rangle_{V, V'}, \\
	\textrm{and that } \quad
			\lim_{n \to \infty}
				\langle 
					(\overline{y}_1-\overline{y}_3), 
					(C_{6,n})(\overline{y}_1-\overline{y}_3) 
				\rangle_{V, V'}
		&=
			\sup_{C \in P}
				\langle 
					(\overline{y}_1-\overline{y}_3), 
					C(\overline{y}_1-\overline{y}_3) 
				\rangle_{V, V'}.
	\end{align} 
	Furthermore, denote by $(C_{0,n})_{n \in \N} \subseteq P$ the sequence satisfying
	for all $n \in \N$ that
	$
		C_{0,n} = \max\{ C_{i,n} \colon i \in \{1, 2, 3, 4, 5, 6\} \}.
	$
	Note that the fact that $P$ is totally ordered ensures that
	$(C_{0,n})_{n \in \N}$ is well defined.
	In addition, we have
	\begin{equation}
		\begin{split}
			&\liminf_{n \to \infty}
					\langle 
						(\overline{y}_1+\overline{y}_2+\overline{y}_3),
						(C_{0,n})(\overline{y}_1+\overline{y}_2+\overline{y}_3) 
					\rangle_{V, V'} \\
			\geq{}&
				\lim_{n \to \infty}
					\langle 
						(\overline{y}_1+\overline{y}_2+\overline{y}_3), 
						(C_{1,n})(\overline{y}_1+\overline{y}_2+\overline{y}_3) 
					\rangle_{V, V'} \\
			={}&
				\sup_{C \in P}
					\langle 
						(\overline{y}_1+\overline{y}_2+\overline{y}_3), 
						C(\overline{y}_1+\overline{y}_2+\overline{y}_3) 
					\rangle_{V, V'} \\
			\geq{}& 
				\limsup_{n \to \infty}
					\langle 
						(\overline{y}_1+\overline{y}_2+\overline{y}_3), 
						(C_{0,n})(\overline{y}_1+\overline{y}_2+\overline{y}_3) 
					\rangle_{V, V'}
		\end{split}
	\end{equation}
	and thus
	$
			\sup_{C \in P}
				\langle
					(\overline{y}_1+\overline{y}_2+\overline{y}_3),
					C(\overline{y}_1+\overline{y}_2+\overline{y}_3) 
				\rangle_{V, V'}
		= 
			\lim_{n \to \infty}
				\langle 
					(\overline{y}_1+\overline{y}_2+\overline{y}_3),
					(C_{0,n})(\overline{y}_1+\overline{y}_2+\overline{y}_3) 
				\rangle_{V, V'}
	$.
	The corresponding equations follows for $(C_{i,n})_{n \in \N}$, 
	$i \in \{2, 3, 4, 5, 6\}$
	analogously. Therefore we obtain that
	\begin{align}
		\nonumber
				&4\langle 
					\overline{y}_1, 
					C_P (\overline{y}_2+\overline{y}_3) 
				\rangle_{V, V'} \\
		\nonumber
			={}&
					\sup_{C \in P}
						(\langle 
							(\overline{y}_1+\overline{y}_2+\overline{y}_3), 
							C(\overline{y}_1+\overline{y}_2+\overline{y}_3) 
						\rangle_{V, V'})
					-\sup_{C \in P}
						(\langle 
							(\overline{y}_1-\overline{y}_2-\overline{y}_3), 
							C(\overline{y}_1-\overline{y}_2-\overline{y}_3) 
						\rangle_{V, V'}) \\
		\nonumber
			={}&
					\lim_{n \to \infty} (
						\langle
							(\overline{y}_1+\overline{y}_2+\overline{y}_3), 
							C_{0,n} (\overline{y}_1+\overline{y}_2+\overline{y}_3) 
						\rangle_{V, V'}
					- \langle 
							(\overline{y}_1-\overline{y}_2-\overline{y}_3), 
							C_{0,n} (\overline{y}_1-\overline{y}_2-\overline{y}_3) 
						\rangle_{V, V'}
					) \\
		\nonumber
			={}&
					\lim_{n \to \infty} (
						4\langle 
							\overline{y}_1, 
							C_{0,n} (\overline{y}_2+\overline{y}_3) 
						\rangle_{V, V'}
					)
			=
					\lim_{n \to \infty} (
						4 \langle 
								\overline{y}_1, C_{0,n} \overline{y}_2 
							\rangle_{V, V'}
						+4\langle 
								\overline{y}_1, C_{0,n} \overline{y}_3 
							\rangle_{V, V'}
					) \\
			={}&
					\lim_{n \to \infty} (
						\langle 
							(\overline{y}_1+\overline{y}_2), 
							C_{0,n} (\overline{y}_1+\overline{y}_2) 
						\rangle_{V, V'}
					- \langle 
							(\overline{y}_1-\overline{y}_2), 
							C_{0,n} (\overline{y}_1-\overline{y}_2) 
						\rangle_{V, V'} \\
		\nonumber
					&+ \langle 
							(\overline{y}_1+\overline{y}_3), 
							C_{0,n} (\overline{y}_1+\overline{y}_3) 
						\rangle_{V, V'}
					-\langle 
							(\overline{y}_1-\overline{y}_3), 
							C_{0,n} (\overline{y}_1-\overline{y}_3) 
						\rangle_{V, V'}
					) \\
		\nonumber
			={}&
					\sup_{C \in P} 
						(\langle 
							(\overline{y}_1+\overline{y}_2), 
							C(\overline{y}_1+\overline{y}_2)
						\rangle_{V, V'})
					-\sup_{C \in P}
						(\langle 
							(\overline{y}_1-\overline{y}_2), 
							C (\overline{y}_1-\overline{y}_2) 
						\rangle_{V, V'}) \\
		\nonumber
					&+\sup_{C \in P}
						(\langle 
							(\overline{y}_1+\overline{y}_3), 
							C (\overline{y}_1+\overline{y}_3) 
						\rangle_{V, V'}) 
					-\sup_{C \in P}
						(\langle
							(\overline{y}_1-\overline{y}_3),
							C (\overline{y}_1-\overline{y}_3) 
						\rangle_{V, V'}) \\
		\nonumber
			={}&
				4\langle \overline{y}_1, C_P \overline{y}_2 \rangle_{V, V'}
				+ 4\langle \overline{y}_1, C_P \overline{y}_3 \rangle_{V, V'}.
	\end{align}
	Since $\overline{y}_1, \overline{y}_2, \overline{y}_3$ were arbitrary,
	this shows that $C_P \in S$ and that $P$ has a maximal element.
	Zorn's lemma now implies the existence of
	a maximal element $A \in S$. 
	Next denote by $\tilde{H}_2 \subseteq H$ the set satisfying that
	$\tilde{H}_2 = \{ x_2 \in H \colon (\forall x_1 \in H_1 \colon g(A,x_1,x_2) =0\}$.
	Then we obtain from \eqref{eq: g bounded from below on H_1} and
	from Lemma \ref{lem: orthogonal set} that $\tilde{H}_2$
	is a with respect to the $\| \cdot \|_H$-norm
	closed linear subspace such that $H_1 \cap \tilde{H}_2 = \{0\}$
	and such that $H_1 + \tilde{H}_2 =H$.
	Combining this with the fact that 
	$H_1 = \ker((D_{\mathbb{H}, \mathbb{V}} f)(x_0))$ and the fact that
	$(D_{\mathbb{H}, \mathbb{V}} f)(x_0)$ is surjective
	yields that
	$(D_{\mathbb{H}, \mathbb{V}} f)(x_0) |_{\tilde{H}_2}$ is bijective. 
	Next we show that for all $x \in \tilde{H}_2$ it holds that
	$g(A,x,x) = 0$.
	Therefore denote for every $y \in \tilde{H}_2$ by
	$\hat{A}_y \in \mathbb{S}_{\mathbb{V}, \mathbb{V}'}$ 
	the function satisfying for all $x_1, x_2 \in H$
	that
	\begin{equation}
	\label{eq: def of hat A}
			\langle 
				\hat{A}_y (D_{\mathbb{H}, \mathbb{V}} f)(x_0) (x_1),
				(D_{\mathbb{H}, \mathbb{V}} f)(x_0) (x_2)
			\rangle_{V', V}
		=
			\frac{g(A, x_1, y) \cdot g(A, x_2, y)} {g(A,y,y)+1}.
	\end{equation}
	Note that it follows 
	from the fact that $H_1 = \ker((D_{\mathbb{H}, \mathbb{V}} f)(x_0))$,
	the fact that for all
	$y \in \tilde{H}_2$ and all $x \in H_1$ it holds that
	$g(A,x,y)=0$, and from the fact that 
	$
		(D_{\mathbb{H}, \mathbb{V}} f)(x_0)
	$
	is surjective
	that for every $y \in \tilde{H}_2$ it holds that
	$\hat{A}_y$ is well defined.
	Moreover, we get with \eqref{eq: def of g},
	\eqref{eq: def of hat A},
	the fact that $H_1 = \ker((D_{\mathbb{H}, \mathbb{V}} f)(x_0))$, 
	the fact that for all $x \in H_1$ and all $z \in H_2$ it holds that
	$g(A,x,z)=0$
	and with the fact for all $x \in H$ it holds that
	$
		g(A,x,x) \geq 0
	$
	that for all $x_1 \in H_1$ and all $x_2, y \in \tilde{H}_2$ it holds
	that
	\begin{equation}
		\begin{split}
				&g(A+\hat{A}_y,x_1 +x_2 , x_1 +x_2) \\
			={}
				&g(A, x_1 +x_2, x_1 +x_2)
				- \langle 
						\hat{A}_y (D_{\mathbb{H}, \mathbb{V}} f)(x_0) (x_1 + x_2),
						(D_{\mathbb{H}, \mathbb{V}} f)(x_0) (x_1 + x_2)
					\rangle_{V', V} \\
			={}&
				g(A,x_1, x_1)
				+ g(A,x_2, x_2)
				- \langle 
						\hat{A}_y (D_{\mathbb{H}, \mathbb{V}} f)(x_0) (x_2),
						(D_{\mathbb{H}, \mathbb{V}} f)(x_0) (x_2)
					\rangle_{V', V} \\
			={}&
				g(A,x_1, x_1)
				+ g(A,x_2, x_2)
				- \frac{g(A, x_2, y) \cdot g(A, x_2, y)} {g(A,y,y)+1} \\
			\geq{}
				&g(A,x_1, x_1)
				+ g(A,x_2, x_2)
				- \frac{g(A, x_2, x_2) \cdot g(A, y, y)} {g(A,y,y)+1} \\
			\geq{}&
				g(A,x_1, x_1)
				+ g(A,x_2, x_2)
				- g(A, x_2, x_2)
			\geq
				0.
		\end{split}
	\end{equation}
	Thus we get for all $y \in \tilde{H}_2$ that
	$A+\hat{A}_y \in S$ and this with the fact that for all $y \in \tilde{H}_2$
	it holds that
	$\hat{A}_y \geq 0$ and the fact that $A \in S$ is a maximal element
	shows for every $y \in \tilde{H}_2$ that $\hat{A}_y = 0$.
	Taking now
	$x_1=x_2= y$ 
	in \eqref{eq: def of hat A}
	yields for all $y \in \tilde{H}_2$
	that
	$g(A, y, y) = 0$.
	Next note that the fact that 
	$(D_{\mathbb{H}, \mathbb{V}} f)(x_0)|_{\tilde{H}_2} \colon \tilde{H}_2 \to V$
	is bijective together with 
	the inverse function theorem implies that
	there exist a relative to $\{x_0 \} + \tilde{H}_2$ open set 
	$O_2 \subseteq O \cap (\{x_0 \} + \tilde{H}_2)$ 
	and an open set
	$U_2 \subseteq U$ such that 
	$x_0 \in O_2$,
	and that
	$(f|_{O_2})^{-1} \in \C_{\mathbb{V}, \mathbb{H}}^2(U_2, O_2)$.
	Therefore we have that
	\begin{equation}
	\label{eq: frak f boundedness}
		\begin{split}
				\| 
					D_{\mathbb{V},\mathbb{H}} (f|_{O_2})^{-1} (f(x_0))
				\|^2_{L(\mathbb{V},\mathbb{H})}
			\geq{}&
				\limsup_{U_2 \ni y \to f(x_0)}
					\frac 
						{ \| (f|_{O_2})^{-1}(y) - (f|_{O_2})^{-1}(f(x_0)) \|^2_H}
						{ \| y - f(x_0) \|^2_V} \\
			\geq{}&
				\limsup_{O_2 \ni x \to x_0}
					\frac 
						{ \| x - x_0 \|^2_H}
						{ \| f(x) - f(x_0) \|^2_V}.
		\end{split}
	\end{equation}
	Again denote by $r_1, r_2 \in \mathbb{M}(O, U)$ the functions
	satisfying \eqref{eq: r2 convergence},
	\eqref{eq: taylor f}, and \eqref{eq: def of r_1}.
	Then it follows from \eqref{eq: frak f boundedness},
	\eqref{eq: def of r_1},
	and from
	\eqref{eq: r2 convergence}
	that
	\begin{align}
		\nonumber
				&\limsup_{O_2 \ni x \to x_0} \left(
					\frac{r_2(x)}{\|f(x) - f(x_0)\|^2_V} 
					+\frac{r_1(x)}{\|f(x) - f(x_0)\|_V}
				\right) \\
		\nonumber
			\leq{}&
				\limsup_{O_2 \ni x \to x_0} \left(
					\frac
						{
							r_2(x) \cdot
							\| 
								D_{\mathbb{V},\mathbb{H}} (f|_{O_2})^{-1} (f(x_0))
							\|^2_{L(\mathbb{V},\mathbb{H})}
						}
						{\|x - x_0\|^2_H} 
					+\frac
						{
							r_1(x) \cdot
							\| 
								D_{\mathbb{V},\mathbb{H}} (f|_{O_2})^{-1} (f(x_0))
							\|_{L(\mathbb{V},\mathbb{H})}
						}
						{\|x - x_0\|_H}
				\right)	\\
	\label{eq: r1 + r2 estimates}
			={}& 0.
	\end{align}
	Furthermore,
	the fact that 
	$f|_{O_2} \colon O_2 \to U_2$ is bijective,
	\eqref{eq: taylor f}, \eqref{eq: def of r_1},
	\eqref{eq: r1 + r2 estimates},
	\eqref{eq: def of p},
	\eqref{eq: def of g},
	the fact that for all $x \in O_2$ it holds that
	$x-x_0 \in \tilde{H}_2$ and thus
	$g(A,x-x_0, x-x_0)=0$,
	\eqref{eq: frak f boundedness},
	and 
	$(\tilde{p}, \tilde{A}) \in (J^2_{\mathbb{H}, +} (u \circ f)) (x_0)$
	yields
	that
	\begin{align}
		\nonumber
				&\limsup_{U \ni y \to f(x_0)} \bigg(
						\tfrac
							{
								u(y) - u(f(x_0)) 
								- \langle 
										p, y - f(x_0)
									\rangle_{V', V}
								-\frac 12 
									\langle 
										y- f(x_0), 
										A (y- f(x_0) ) 
									\rangle_{V, V'}
							}
							{ \| y - f(x_0) \|^2_V} 
				\bigg) \\
		\nonumber
			={}&
				\limsup_{U_2 \ni y \to f(x_0)} \bigg(
						\tfrac
							{
								u(y) - u(f(x_0)) 
								- \langle 
										p, y - f(x_0)
									\rangle_{V', V}
								-\frac 12 
									\langle 
										y- f(x_0), 
										A (y- f(x_0) ) 
									\rangle_{V, V'}
							}
							{ \| y - f(x_0) \|^2_V} 
				\bigg) \\
		\nonumber
			={}&
				\limsup_{O_2 \ni x \to x_0} \bigg(
						\tfrac
							{
								u(f(x)) - u(f(x_0)) 
								- \langle 
										p, f(x) - f(x_0)
									\rangle_{V', V}
								-\frac 12 
									\langle 
										f(x)- f(x_0), 
										A (f(x)- f(x_0) ) 
									\rangle_{V, V'}
							}
							{ \| f(x) - f(x_0) \|^2_V} 
				\bigg) \\
		\nonumber
			={}&
				\limsup_{O_2 \ni x \to x_0} \bigg( 
						\tfrac
							{
								u(f(x)) - u(f(x_0)) 
								- \big \langle 
										p, 
										\big( (D_{\mathbb{H}, \mathbb{V}}f)(x_0) \big)(x- x_0)
										+\frac 12 
										\big(
											\big( (D^2_{\mathbb{H}, \mathbb{V}}f)(x_0) \big)(x-x_0) 
										\big)
											(x-x_0)
										+r_2(x)
									\big \rangle_{V', V}
							}
							{ \| f(x) - f(x_0) \|^2_V} \\
					& \qquad \qquad \qquad
						-\tfrac{
								\frac 12 
									\big \langle 
										\big( (D_{\mathbb{H}, \mathbb{V}}f)(x_0) \big)(x- x_0)
										+ r_1(x), 
										A \big(
											\big( (D_{\mathbb{H}, \mathbb{V}}f)(x_0) \big)(x- x_0)+ r_1(x)
										\big)
									\big \rangle_{V, V'}
							}
							{ \| f(x) - f(x_0) \|^2_V}
				\bigg) \\
		\nonumber
			={}&
				\limsup_{O_2 \ni x \to x_0} \bigg(
						\tfrac
							{
								u(f(x)) - u(f(x_0)) 
								- \big \langle 
										p, \big( (D_{\mathbb{H}, \mathbb{V}}f)(x_0) \big)(x- x_0)
								+ \frac 12 
									\big(\big( (D^2_{\mathbb{H}, \mathbb{V}}f)(x_0) \big)(x-x_0) \big)(x-x_0)
								\big \rangle_{V', V}
							}
							{ \| f(x) - f(x_0) \|^2_V} \\
		\nonumber
				& \qquad \qquad \qquad
						-\tfrac{
								\frac 12 
									\big \langle 
										\big( (D_{\mathbb{H}, \mathbb{V}}f)(x_0) \big)(x- x_0), 
										A \big(
											\big( (D_{\mathbb{H}, \mathbb{V}}f)(x_0) \big)(x- x_0)
										\big)
									\big \rangle_{V,V'}
							}
							{ \| f(x) - f(x_0) \|^2_V} 
				\bigg) \\
		\nonumber
			={}&
				\limsup_{O_2 \ni x \to x_0} \bigg( 
						\tfrac
							{
								u(f(x)) - u(f(x_0)) 
								- \langle 
										\tilde{p}, (x- x_0)
									\rangle_{H', H}
								+\frac 12 g(A,  (x- x_0), (x- x_0))
								- \frac 12 \langle  (x- x_0), \tilde{A}  (x- x_0) \rangle_{H, H'}
							}
							{ \| x - x_0 \|^2_H} \\
		\nonumber
					& \qquad \qquad \qquad
						-\tfrac
							{
								-\frac \lambda2 \|\pi^H_{H_1} x- \pi^H_{H_1} x_0\|^2_H
							}
							{ \| f(x) - f(x_0) \|^2_V} 
				\bigg) \\
		\nonumber
			\leq{}&
				\limsup_{O_2 \ni x \to x_0} \bigg(
						\tfrac
							{
								u(f(x)) - u(f(x_0)) 
								- \langle 
										\tilde{p}, (x- x_0)
									\rangle_{H', H}
								- \frac 12 \langle  (x- x_0), \tilde{A}  (x- x_0) \rangle_{H, H'}
							}
							{ \| x - x_0 \|^2_H} 
						\cdot
						\tfrac 
							{ \| x - x_0 \|^2_H}
							{ \| f(x) - f(x_0) \|^2_V}
				\bigg) 
			\leq
				0.
	\end{align}
	This shows that $(p, A) \in (J_{\mathbb{V}, +}^2 u) (f(x_0))$.
	Moreover, let $B \in \mathbb{S}_{\mathbb{H}, \mathbb{H}'}$ be the function 
	satisfying for all
	$x_1, x_2 \in H$ that 
	$
			\langle x_1, B x_2 \rangle_{H, H'} 
		= 
			g(A, x_1, x_2) - \lambda \langle \pi^H_{H_1} x_1, x_2 \rangle_H 
	$.
	Then 
	we have that $B \geq - \lambda I_\mathbb{H} \pi^H_{H_1}$ and it follows from 
	\eqref{eq: def of g} that for all $x_1, x_2 \in H$ it holds that 
	\begin{align}
	\nonumber
			&\big \langle 
				A \big( (D_{\mathbb{H}, \mathbb{V}} f)(x_0) x_1 \big), 
				(D_{\mathbb{H}, \mathbb{V}} f)(x_0) x_2 
			\big \rangle_{V', V}
			+ \big \langle
					p, \big( \big( (D^2_{\mathbb{H}, \mathbb{V}} f)(x_0) \big) (x_1) \big) (x_2)
				\big \rangle_{V', V}
			+ \langle x_1, Bx_2 \rangle_{H, H'} \\
		&={}
			\langle x_1, \tilde{A} x_2 \rangle_{H, H'}.
	\end{align}
	This shows the second part of the assertion and therefore completes the proof
	of Lemma \ref{lem: chain rule}.
\end{proof}
\section{Lifting viscosity solutions to bigger Hilbert spaces}
The chain rule for semijets allows us to
lift viscosity solutions to bigger spaces.
To the best of our knowledge, lifting of viscosity solutions to bigger
Hilbert spaces is new.
\begin{prop}[Lifting of classical viscosity solutions to
	viscosity solutions in bigger Hilbert spaces]
\label{prop: u proj is viscosity solution}
	Assume the setting in Section \ref {ssec: Setting H X without t}
	and assume that
	$
			\{ 
				y \in X \cap O 
					\colon (D_{\mathbb{X}} (h|_X)) (y) \in D(E_{\mathbb{X}', \mathbb{H}'})  
			\}
		=
			\{ 
				y \in X \cap O \colon  
				~ (J^{2}_{\mathbb{H}, -}h) (y) \neq \emptyset 
			\}
	$,
	that for all $x \in O \cap X$ it holds that
	$
			\underline{h}_\mathbb{H}^W(x)
		=
			h(x)
	$,
	and that for all $x \in O \backslash X$ it holds that
	$
		h(x) = \infty
	$,
	let $\mathbb{V}=(V, \langle \cdot, \cdot \rangle_V, \| \cdot \|_V)$
	be a real Hilbert space,
	let $\mathbb{V}' = (V', \| \cdot \|_{H'})$ be its dual space,
	let $f \in \C_{\mathbb{H}, \mathbb{V}}(O, V)$ satisfy that
	$f|_X \in \C^2_{\mathbb{X}, \mathbb{V}}(O \cap X, V)$
	and that for all 
	$
		x_0 \in X \cap O
	$ 
	it holds that
	$(D_{\mathbb{X}, \mathbb{V}} (f|_X))(x_0) \in \C_{\mathbb{H},\mathbb{V}}(X, V)$
	and that
	$(D_{\mathbb{X}, \mathbb{V}} (f|_X))(x_0) \colon X \to V$ is surjective,
	denote for every $x \in X$ by $K_x^\perp$ 
	the orthogonal set of $\ker((D_{\mathbb{X}, \mathbb{V}} (f|_X))(x))$ satisfying that
	$
			K_x^\perp 
		= 
			\{
				y \in X \colon 
				(
					\forall v \in \ker((D_{\mathbb{X}, \mathbb{V}} (f|_X))(x)) 
						\colon \langle y, v \rangle_X = 0
				) 
			\}
	$,
	let 
	$
		F \colon 
			W 
			\times \R \times H' \times \mathbb{S}_{\mathbb{X}, \mathbb{X}'} \to \R
	$ 
	and
	$
		\tilde{F} \colon 
			f(O) 
			\times \R \times V' \times \mathbb{S}_{\mathbb{V}, \mathbb{V}'} \to \R
	$ 
	be degenerate elliptic functions
	such that for all 
	$x_0 \in W$, $r \in \R$, $\tilde{p} \in H'$, 
	$\tilde{A} \in \mathbb{S}_{\mathbb{X}, \mathbb{X}'}$,
	$p \in V'$, $A \in \mathbb{S}_{\mathbb{V}, \mathbb{V}'}$,
	and all
	$(A_\lambda)_{\lambda \in (0, \infty)} \subseteq \mathbb{S}_{\mathbb{H}, \mathbb{H}'}$,
	with 
	\begin{align}
			\begin{split}
					&\forall x_1 \in X \colon
					\langle p, (D_{\mathbb{X}, \mathbb{V}} f)(x_0)x_1 \rangle_{V', V} 
				= 
					\langle \tilde{p}|_X, x_1 \rangle_{X', X} \qquad \textrm{ and with} \\
					&\forall x_1,x_2 \in X \colon
					\big \langle 
						A \big( (D_{\mathbb{X}, \mathbb{V}} f)(x_0)  x_1 \big), 
						(D_{\mathbb{X}, \mathbb{V}} f)(x_0) x_2 
					\big \rangle_{V', V}
					+ \big \langle
							p, \big( \big( (D^2_{\mathbb{X}, \mathbb{V}} f)(x_0) \big) (x_1) \big) (x_2)
						\big \rangle_{V', V} \\
			& \qquad \qquad \quad
				=
					\langle x_1, \tilde{A} x_2 \rangle_{X, X'}
			\end{split}
		\end{align}
	it holds that
	\begin{equation}
	\label{eq: F tilde F relation}
			\tilde{F}(f(x_0),r,p, A)
		=
			F(
				x_0,
				r,
				\tilde{p},
				\tilde{A}
			),
	\end{equation}
	that
	\begin{equation}
			\liminf_{\lambda \downarrow 0}
				F(
					x_0,
					r,
					\tilde{p}, 
					(A_\lambda|_X)|_X 
					- \lambda I_{\mathbb{X}} \pi_{K_{x_0}^\perp}^X
				)
		\leq
			\limsup_{\lambda \downarrow 0} F(x_0,r,\tilde{p}, (A_\lambda|_X)|_X),
	\end{equation}
	and that
	\begin{equation}
			\limsup_{\lambda \downarrow 0} 
				F(
					x_0,
					r,
					\tilde{p}, 
					(A_\lambda|_X)|_X + \lambda I_{\mathbb{X}} \pi_{K_{x_0}^\perp}^X
				)
		\geq
			\liminf_{\lambda \downarrow 0} F(x_0,r,\tilde{p}, (A_\lambda|_X)|_X),
	\end{equation}
	assume that $\tilde{F}$ is continuous with respect to the
	$
		\| \cdot \|_{V \times \R \times V' \times L(\mathbb{V}, \mathbb{V}')}
	$-norm
	and let $u \in \C_{\mathbb{V}} (f(O), \R)$
	be a classical viscosity solution of
	$
		\tilde{F}=0
	$.
	Then
	$
		u \circ f
	$
	is a viscosity solution of
	$
		F=0
	$
	relative to  $(h, \mathbb{H}, \mathbb{X})$.
\end{prop}
\begin{proof}[Proof of Proposition~\ref{prop: u proj is viscosity solution}]
	If 
	$ 
			\{
				x \in O \colon
					(J_{\mathbb{H}, +}^2 (u \circ f)^{-,W}_{\mathbb{H}, \delta, h})(x) \neq \emptyset 
			\}
		= \emptyset	
	$, 
  then the assertion is trivial. 
  So for the rest of 
  the proof, we assume that 
  $ 
			\{
				x \in O \colon
					(J_{\mathbb{H}, +}^2 (u \circ f)^{-,W}_{\mathbb{H}, \delta, h})(x) \neq \emptyset 
			\}
		\neq \emptyset	
	$.
	Fix $x_0 \in O$, $(\tilde{p}, \tilde{A}) \in H' \times \mathbb{S}_{\mathbb{H}, \mathbb{H}'}$,
	$\delta \in (0, \infty)$ such that
	$
		(\tilde{p}, \tilde{A}) \in 
			(J_{\mathbb{H}, +}^2 (u \circ f)^{-,W}_{\mathbb{H}, \delta, h})(x_0).
	$
	Note that the continuity 
	with respect to the $\| \cdot \|_{\mathbb{H}}$-norm 
	of $u$ and $f$ together with
	the assumption that for all $x \in X \cap O$ it holds that
	$
			\underline{h}_\mathbb{H}^W(x)
		=
			h(x)
	$
	shows that
	\begin{equation}
			(u \circ f)^{-,W}_{\mathbb{H}, \delta, h} 
		=
			\overline{(u \circ f - \delta h)}_\mathbb{H}^W
		= 
			u \circ f - \delta \underline{h}_\mathbb{H}^W
		=
			u \circ f - \delta h.
	\end{equation}
	Furthermore, the assumption that for all $x \in O \backslash X$ it holds that
	$
		h(x) = \infty
	$
	implies that $x_0 \in O \cap X$
	and this together with 
	the assumption 
	that $\mathbb{X}$ is embedded continuously in $\mathbb{H}$
	and the assumption that
	$
		(\tilde{p}, \tilde{A}) \in 
			(J_{\mathbb{H}, +}^2 (u \circ f)^{-,W}_{\mathbb{H}, \delta, h})(x_0)
	$
	ensures that
	$	
		(\tilde{p}|_X, (\tilde{A}|_X)|_X) \in 
			(J_{\mathbb{X}, +}^2 (u \circ f|_X - \delta h|_X))(x_0)
	$.
	Moreover, the assumption that
	$h|_X \in \C_{\mathbb{X}}^2(X \cap O,\R)$
	yields that
	$	
		(
			\tilde{p}|_X + \delta (D_{\mathbb{X}}(h|_X))(x_0), 
			(\tilde{A}|_X)|_X + \delta (D^2_{\mathbb{X}}(h|_X))(x_0)
		)
			\in (J_{\mathbb{X}, +}^2 (u \circ f|_X))(x_0)
	$.
	Thus 
	Lemma \ref{lem: chain rule}
		(with 
		$\mathbb{H} \leftarrow \mathbb{X}$,
		$U \leftarrow f (O)$,
		$\tilde{p} \leftarrow \tilde{p}|_X + \delta (D_{\mathbb{X}}(h|_X))(x_0)$,
		and with
		$\tilde{A} \leftarrow (\tilde{A}|_X)|_X + \delta (D_{\mathbb{X}}(h|_X))(x_0)$)
	implies that for every $\lambda \in (0,\infty)$ 
	there exist
	$B_\lambda \in \mathbb{S}_{\mathbb{X}, \mathbb{X}'}$ and
	$(p, A_\lambda ) \in (J_{\mathbb{V}, +}^2 u) (f(x_0))$
	such that
	$B_\lambda \geq -\lambda I_\mathbb{X} \pi^X_{K_{x_0}^\perp}$ and that
	for all $x_1, x_2 \in X$ it holds that
	\begin{align}
	\label{eq: def of p lem}
				&\langle p, (D_{\mathbb{X}, \mathbb{V}} \, (f|_X))(x_0)x_1 \rangle_{V', V} 
			= 
				\langle 
					\tilde{p}|_X + \delta (D_{\mathbb{X}}(h|_X))(x_0), x_1 
				\rangle_{X', X} \qquad \textrm{ and that} \\
		\nonumber
				&\big \langle 
					A_\lambda \big((D_{\mathbb{X}, \mathbb{V}} \, (f|_X))(x_0)  x_1 \big), 
					(D_{\mathbb{X}, \mathbb{V}} \, (f|_X))(x_0) x_2 
				\big \rangle_{V', V}
				+ \big \langle
						p, 
						\big( \big( 
								(D^2_{\mathbb{X}, \mathbb{V}} \, (f|_X))(x_0) 
						\big) (x_1) \big) (x_2)
					\big \rangle_{V', V} \\
		\label{eq: def of A lem}
				&+ \langle x_1, B_\lambda x_2 \rangle_{X, X'}
			=
				\langle 
					x_1, \big( (\tilde{A}|_X)|_X + \delta (D^2_{\mathbb{X}}(h|_X))(x_0) \big)x_2 
				\rangle_{X, X'}.
	\end{align}
	Note that \eqref{eq: def of p lem} ensures that $p$ does not depend on
	$\lambda$. 
	Moreover, we get from \eqref{eq: def of p lem} and the
	assumption 
	$(D_{\mathbb{X}, \mathbb{V}} (f|_X))(x_0) \in \C_{\mathbb{H},\mathbb{V}}(X, V)$
	that
	$(D_{\mathbb{X}}(h|_X))(x_0) \in D(E_{\mathbb{X}', \mathbb{H}'})$
	and thus it holds that $x_0 \in W$.
	Combining now 
	$x_0 \in W$,
	\eqref{eq: def F+ delta h u},
	\eqref{eq: def F+ delta h},
	the fact that for all $\lambda \in (0, \infty)$ it holds that
	$B_\lambda \geq -\lambda I_{\mathbb{X}} \pi^X_{K_{x_0}^\perp}$
	and the assumption that F is degenerate elliptic,
	the assumption that for all $x \in W$, $r \in \R$, $q \in H'$, 
	and all
	$(A_\lambda)_{\lambda \in (0, \infty)} \subseteq \mathbb{S}_{\mathbb{H}, \mathbb{H}'}$
	it holds that
	$
			\liminf_{\lambda \downarrow 0} 
				F(x,r,q, (A_\lambda|_X)|_X - \lambda I_{\mathbb{X}} \pi^X_{V^\perp})
		\leq
			\limsup_{\lambda \downarrow 0} 
				F(x,r,q, (A_\lambda|_X)|_X)
	$
	yields 
	that
	\begin{equation}
		\begin{split}
				&F^+_{\mathbb{H}, \mathbb{X}, \delta, h, u \circ f}
				\big (	
					x_0, 
					(u \circ f)^{-,W}_{\mathbb{H}, \delta, h}(x_0), 
					\tilde{p},
					\tilde{A}
				\big ) 
			\leq
				F^+_{\mathbb{H}, \mathbb{X}, \delta, h}
				\big(	
					x_0, 
					u( f (x_0))- \delta h(x_0), 
					\tilde{p},
					\tilde{A}
				\big) \\
			={}&
				F
				\big(	
					x_0, 
					u( f (x_0))- \delta h(x_0) + \delta h(x_0), 
					\tilde{p} + \delta  E_{\mathbb{X}', \mathbb{H}'} ((D_{\mathbb{X}} (h|_X))(x_0)),
					(\tilde{A}|_X) |_X + \delta  (D^2_{\mathbb{X}} (h|_X))(x_0)
				\big) \\
			\leq{}&
				\liminf_{ \lambda \downarrow 0} 
					F
					\big(	
						x_0, 
						u( f (x_0)), 
						\tilde{p} + \delta  E_{\mathbb{X}', \mathbb{H}'} ((D_{\mathbb{X}} (h|_X))(x_0)),
						(\tilde{A}|_X) |_X + \delta  (D^2_{\mathbb{X}} (h|_X))(x_0) -B_\lambda 
		\\ & \qquad \qquad
						-\lambda I_\mathbb{X} \pi^X_{K_{x_0}^\perp}
					\big) \\
			\leq{}&
				\limsup_{ \lambda \downarrow 0} 
					F
					\big(	
						x_0, 
						u( f (x_0)), 
						\tilde{p} + \delta  E_{\mathbb{X}', \mathbb{H}'} ((D_{\mathbb{X}} (h|_X))(x_0)),
						(\tilde{A}|_X) |_X + \delta  (D^2_{\mathbb{X}} (h|_X))(x_0) -B_\lambda
					\big).
		\end{split}
	\end{equation}
	Therefore it follows from \eqref{eq: F tilde F relation},
	\eqref{eq: def of p lem}, \eqref{eq: def of A lem},
	the assumption that u is a classical viscosity subsolution of $\tilde{F} = 0$,
	the fact that for all $\lambda \in (0, \infty)$ it holds that
	$(p, A_\lambda ) \in (J_{\mathbb{V}, +}^2 u) (f(x_0))$,
	and from the continuity of $\tilde{F}$
	that
	\begin{equation}
		\begin{split}
				&F^+_{\mathbb{H}, \mathbb{X}, \delta, h, u \circ f}
				(	
					x_0, 
					(u \circ f)^{-,W}_{\mathbb{H}, \delta, h}(x_0), 
					\tilde{p},
					\tilde{A}
				)\\
			\leq{}&
				\limsup_{ \lambda \downarrow 0} 
					F
					\big(	
						x_0, 
						u( f (x_0)), 
						\tilde{p} + \delta  E_{\mathbb{X}', \mathbb{H}'} ((D_{\mathbb{X}} (h|_X))(x_0)),
						(\tilde{A}|_X) |_X + \delta  (D^2_{\mathbb{X}} (h|_X))(x_0) -B_\lambda
					\big)  \\
			={}&
				\limsup_{\lambda \downarrow 0} 
					\tilde{F}(f(x_0),u( f(x_0)),p, A_\lambda)
			=
				\limsup_{\lambda \downarrow 0} 
					\tilde{F}_{\mathbb{V},\mathbb{V},1,0,u}(f(x_0),u( f (x_0)),p, A_\lambda)
			\leq
				0.
		\end{split}
	\end{equation}
	Thus Corollary \ref{cor:semijets_equivalence} shows that
	$u \circ \pi_V$ is a viscosity subsolution of 
	$ F  = 0 $ 
	relative to $(h, \mathbb{H}, \mathbb{X})$.
	In the same way, it can be shown
  that $ u \circ \pi_V$ is a viscosity supersolution
  of $ F = 0 $ relative to $(h, \mathbb{H}, \mathbb{X})$ and we thereby obtain
  that $ u \circ \pi_V$ is a viscosity solution
  of $ F = 0 $ relative to $(h, \mathbb{H}, \mathbb{X})$.
	This
	completes the proof
	of Proposition \ref{prop: u proj is viscosity solution}.
\end{proof}
\begin{corollary}[Lifting classical viscosity solutions to
	classical viscosity solutions in bigger Hilbert spaces]
\label{prop: u proj is classical viscosity solution}
	Assume the setting in Section \ref {ssec: Setting H X without t},
	let $\mathbb{V}=(V, \langle \cdot, \cdot \rangle_V, \| \cdot \|_V)$
	be a real Hilbert space,
	let $\mathbb{V}' = (V', \| \cdot \|_{H'})$ be its dual space,
	let $f \in \C^2_{\mathbb{H}, \mathbb{V}}(O, V)$, satisfy 
	for all $x_0 \in O$ that
	$(D_{\mathbb{H}, \mathbb{V}} f)(x_0) \colon X \to V$ is surjective,
	denote for every $x \in H$ by $K_x^\perp$ 
	the orthogonal set of $\ker((D_{\mathbb{H}, \mathbb{V}} f)(x))$ satisfying that
	$
			K_x^\perp 
		= 
			\{
				y \in H \colon 
				(
					\forall v \in \ker((D_{\mathbb{H}, \mathbb{V}} f)(x)) 
						\colon \langle y, v \rangle_H = 0
				) 
			\}
	$,
	let 
	$
		F \colon 
			O
			\times \R \times H' \times \mathbb{S}_{\mathbb{H}, \mathbb{H}'} \to \R
	$ 
	and
	$
		\tilde{F} \colon 
			f(O) 
			\times \R \times V' \times \mathbb{S}_{\mathbb{V}, \mathbb{V}'} \to \R
	$ 
	be degenerate elliptic functions
	such that for all 
	$x_0 \in O$, $r \in \R$, $\tilde{p} \in H'$, 
	$\tilde{A} \in \mathbb{S}_{\mathbb{H}, \mathbb{H}'}$,
	$p \in V'$, $A \in \mathbb{S}_{\mathbb{V}, \mathbb{V}'}$,
	and all
	$(A_\lambda)_{\lambda \in (0, \infty)} \subseteq \mathbb{S}_{\mathbb{H}, \mathbb{H}'}$,
	with 
	\begin{align}
		\nonumber
					&\forall x_1 \in H \colon
					\langle p, (D_{\mathbb{H}, \mathbb{V}} f)(x_0)x_1 \rangle_{V', V} 
				= 
					\langle \tilde{p}|_X, x_1 \rangle_{H', H} \qquad \textrm{ and with} \\
			\nonumber
					&\forall x_1,x_2 \in H \colon
					\big \langle 
						A \big( (D_{\mathbb{H}, \mathbb{V}} f)(x_0)  x_1 \big), 
						(D_{\mathbb{H}, \mathbb{V}} f)(x_0) x_2 
					\big \rangle_{V', V}
					+ \big \langle
							p, \big( \big( (D^2_{\mathbb{H}, \mathbb{V}} f)(x_0) \big) (x_1) \big) (x_2)
						\big \rangle_{V', V} \\
		& \qquad \qquad \quad
				=
					\langle x_1, \tilde{A} x_2 \rangle_{H, H'}
		\end{align}
	it holds that
	\begin{align}
		\nonumber
			\tilde{F}(f(x_0),r,p, A)
		={}
			&F(
				x_0,
				r,
				\tilde{p},
				\tilde{A}
			), \\
			\liminf_{\lambda \downarrow 0}
				F(x,r,\tilde{p}, A_\lambda - \lambda I_{\mathbb{H}} \pi_{K_{x_0}^\perp}^H)
		\leq{}
			&\limsup_{\lambda \downarrow 0} F(x,r,\tilde{p}, A_\lambda), \\
	\nonumber
	\textrm{ and that } \quad
			\limsup_{\lambda \downarrow 0} 
				F(x,r,\tilde{p}, A_\lambda + \lambda I_{\mathbb{H}} \pi_{K_{x_0}^\perp}^H)
		\geq{}
			&\liminf_{\lambda \downarrow 0} F(x,r,\tilde{p}, A_\lambda),
	\end{align}
	let
	$\tilde{F}$ be continuous with respect to the
	$
		\| \cdot \|_{V \times \R \times V' \times L(\mathbb{V}, \mathbb{V}')}
	$-norm
	and let $u \in \C_{\mathbb{V}} (f(O), \R)$
	be a classical viscosity solution of
	$
		\tilde{F}=0
	$.
	Then
	$
		u \circ f
	$
	is a classical viscosity solution of
	$
		F=0
	$.
	\end{corollary}
	\begin{proof}
		The proof follows directly from Proposition \ref{prop: u proj is viscosity solution}
		 (with $\mathbb{X} \leftarrow \mathbb{H}$ and $h \leftarrow 0$).
	\end{proof}
\chapter{Uniqueness of viscosity 
solutions of parabolic PDEs}
	The main result of this chapter is a comparison result
	for viscosity solutions of parabolic PDEs (Theorem \ref{l:comparison.viscosity.solution})
	and its application to Kolmogorov equations
	(Corollary \ref{cor:uniqueness2})
	under suitable assumptions. 
	This generalizes Lemma 4.13 and Corollary 4.14 in
	Hairer, Hutzenthaler \& Jentzen \cite{HairerHutzenthalerJentzen2015}
	to separable Hilbert spaces.
	Moreover, in contrast to Theorem 6.1 in Ishii \cite{Ishii1993}
	we don't have a global Lipschitz assumption 
	and we allow
	unbounded
	viscosity solution as long as they satisfy a suitable growth condition
	(often polynomial growth).

	In Section \ref{ssec: Construction of semijets}
	we show several results on constructing suitable
	semijets. 
	In particular Lemma \ref{l: suitable semijets finite dim}
	is  a small improvement on Theorem 3.2 in
	Crandall, Ishii \citationand~Lions~\cite{CrandallIshiiLions1992}.
	Lemma \ref{l: suitable semijets matrix} then generalizes
	Lemma 4 in Lions \cite{Lions1989} to arbitrary linear functions
	and 
	Corollary \ref{cor: suitable semijets function} generalizes
	this result to general $\C^2$-functions.
	Corollary \ref{cor: suitable semijets special case} is a special case of
	Corollary \ref{cor: suitable semijets function} and will be needed in the proof
	of Lemma \ref{l:technical.lemma.uniqueness}.
	In the Section \ref{sec: uniqueness 2nd order PDE}
	we prove our main theorem (Theorem \ref{l:comparison.viscosity.solution}).
	To achieve this, we establish a comparison result for bounded 
	viscosity solutions vanishing at infinity
	(Lemma \ref{l:technical.lemma.uniqueness}).  
	Using a scaling argument (Lemma \ref{c:technical.lemma.uniqueness})
	we can transfer this result to  viscosity solutions satisfying a certain growth condition (Theorem \ref{l:comparison.viscosity.solution}).
	In Section \ref{sec: Uniqueness of Kolmogorov equations}
	we apply Theorem \ref{l:comparison.viscosity.solution} to 
	Kolmogorov equations (Corollary \ref{cor:uniqueness2}). 
	This chapter is based on Section 4.3 in
	Hairer, Hutzenthaler \& Jentzen \cite{HairerHutzenthalerJentzen2015} and
	Section 6 in Ishii \cite{Ishii1993}. 
	
	Throughout this chapter we use the notation 
	from Section \ref{ssec: Notation with t} below.
\section{Setting}
\label{ssec: Setting H X with t}
	Throughout this chapter the following setting is frequently used.
	Let $ T \in (0,\infty) $, 
  let $\mathbb{H} = (H, \langle \cdot, \cdot \rangle_H, \| \cdot \|_H)$ and 
	$\mathbb{X} = (X, \langle \cdot, \cdot \rangle_X, \| \cdot \|_X)$  
	be real separable Hilbert spaces  
	such that $X \subseteq H$, $\mathbb{X}$ is embedded continuously in $\mathbb{H}$, and
	$X$ is dense in $H$ with respect to the $\| \cdot \|_H$-norm,
	let
	$\mathbb{H}' = (H', \| \cdot \|_{H'})$,
	$\mathbb{X}' = (X, \langle \cdot, \cdot \rangle_{X'}, \| \cdot \|_{X})$,
	$(\R \times \mathbb{H})' = ((\R \times H)', \| \cdot \|_{(\R \times H)'})$,
	$(\R \times \mathbb{X})' = ((\R \times X)', \| \cdot \|_{(\R \times X)'})$,
	and
	$((\R \times \mathbb{H})^2)' = (((\R \times H)^2)', \| \cdot \|_{((\R \times H)^2)'})$
	be the corresponding dual spaces,
	let 
  $ O \subseteq H $ 
  be an open set,
	let $ O_n \subseteq O $, $ n \in \N $
	be the $\mathbb{H}$-bounded and $\mathbb{H}$-closed sets satisfying for all $n \in \N$ that, 
	$
			O_n
		=
			\left\{
				x \in O \colon
				\dist_{\mathbb{H}}(x, H \backslash O)  
				\ge \tfrac{ 1 }{ n }
				\text{ and }
				\left\| x \right\|_H \le n
			\right\},
	$
	and denote for all $n \in \N$ by $ O_n^c = O \backslash O_n $
	the complement of $ O_n $ in $ O $,
	let $W \subseteq (0,T) \times O$ satisfy 
	that 
	$
			W
	$
	is dense in $(0,T) \times O$ with respect to the $\| \cdot \|_{\R \times H}$-norm, 
	let $h \colon [0,T] \times O \to [0, \infty] $ be a 
	with respect to the $\| \cdot \|_{\R \times \mathbb{H}}$-norm
	lower semicontinuous 
	function
	with the property that 
	$
		h |_{(0,T) \times X} \in \C_{\R \times \mathbb{X}}^2((0,T) \times (X \cap O)) 
	$
	and that
	$
			W
		=
			\{ 
				z \in (0,T) \times (X \cap O) \colon 
					(D_{\mathbb{X}} (h|_{(0,T) \times X})) (z) 
						\in D(E_{\mathbb{X}', \mathbb{H}'}), 
					~ (J^{2}_{\R \times \mathbb{H}, -}h) (z) \neq \emptyset 
			\}.
	$
\section{Notation}
\label{ssec: Notation with t}
Assume the setting in Section \ref{ssec: Setting H X with t}
and let
$G \colon W \times \R \times H' \times \mathbb{S}_{\mathbb{X}, \mathbb{X}'} \to \R$
then we will in the following always consider the function
$
	F \colon 
		W 
		\times \R 
		\times (\R \times H)' 
		\times \mathbb{S}_{\R \times \mathbb{X}, (\R \times \mathbb{X})'} 
		\to \R
$
satisfying for all
$
	((t,x),r,(q,p),A, \tilde{A}) 
		\in W 
		\times \R 
		\times (\R \times H)' 
		\times \mathbb{S}_{\R \times \mathbb{X}, (\R \times \mathbb{X})'} 	
		\times \mathbb{S}_{\mathbb{X}, \mathbb{X}'} 
$
with 
$\forall y \in X \colon \tilde{A} y = \pi^{(\R \times X)'}_{2} A (0,y)$
that
\begin{equation}
		F((t,x),r,(q,p),A)
	= 
		I^{-1}_{\R}q- G((t,x), r, p, \tilde{A}),
\end{equation}
and the corresponding equation
\begin{equation}
	\begin{split}
			&F(
				(t,x),
				u(t,x),
				(D_{\R \times \mathbb{H}} u)(t,x) ,
				((D^2_{\R \times \mathbb{H}} u)(t,x)|_{\R \times X}) |_{\R \times X}
			) \\
		={}& 
			\tfrac{\partial }{\partial t} u (t,x) 
			- G((t,x), u(t,x),(D_{\mathbb{H}} u)(t,x) ,((D^2_{\mathbb{H}} u)(t,x)|_X ) |_X)
		=
			0.
	\end{split}
\end{equation}
Therefore we will use the following notation, which is more suitable for 
parabolic equations
throughout this chapter.
More precisely  
a function
  $
    G \colon 
    W 
    \times\R\times H'
    \times \mathbb{S}_{\mathbb{X}, \mathbb{X}'}
    \to\R
  $
is here called 
{\it degenerate elliptic} if
for all 
$ (t, x) \in W$,
$ r \in \R $,
$ p \in H' $
and all 
$ A, B \in \mathbb{S}_{\mathbb{X}, \mathbb{X}'} $
with
$ A \leq B $ it holds that
$
  G((t, x), r, p, A) \leq G((t, x), r, p, B)
$
(see, e.g.,
inequality~(1.2)
in Appendix~C
in Peng~\cite{Peng2010}
and compare also
with Section~\ref{ssec:Definition 
of viscosity solutions}
above).
Moreover, for $\delta \in (0, \infty)$
we will denote by 
$
	d^{+, W}_{\mathbb{H}, \delta, h, u},
	~d^{-, W}_{\mathbb{H}, \delta, h, u} \colon 
		((0,T) \times O \times \R \times H' \times \mathbb{S}_{\mathbb{H}, \mathbb{H}'})^2 
			\to [0, \infty  ]
$ 
the functions
satisfying for all 
$
		(\xi, \eta)
	= 
		((t,x,r,p,B),(v,y,s,q,C)) 
			\in ((0,T) \times O \times \R \times H' \times \mathbb{S}_{\mathbb{H}, \mathbb{H}'})^2
$ 
that
\begin{equation}
\label{eq:def of d delta}
			d^{+, W}_{\mathbb{H}, \delta, h, u}(\xi, \eta)
		=
			\begin{cases}
				&\|x-y\|_H \vee |t - v| \vee |r-s| \vee \|p-q\|_{H'} 
				\vee \| B-C \|_{L(\mathbb{H},\mathbb{H}')} \\
				&\vee |u^{+, W}_{\R \times \mathbb{H}, \delta, h}(t, x)-
					u^{+, W}_{\R \times \mathbb{H}, \delta, h}(v, y)| \\
				& \qquad \qquad
					\text{if }  u^{+, W}_{\R \times \mathbb{H}, \delta, h}(t, x),
						u^{+, W}_{\R \times \mathbb{H}, \delta, h}(v, y) \in \R\\
				&\infty \\ 
				& \qquad \qquad 
					\text{if } u^{+, W}_{\R \times \mathbb{H}, \delta, h}(t, x),
						u^{+, W}_{\R \times \mathbb{H}, \delta, h}(v, y) \notin \R\\
			\end{cases} 
	\end{equation}
	and that 
	\begin{equation}
	\label{eq:def of d delta-}
			d^{-, W}_{\mathbb{H}, \delta, h, u}(\xi, \eta)
		=
			\begin{cases}
				&\|x-y\|_H \vee |t - v| \vee |r-s| \vee \|p-q\|_{H'} 
				\vee \| B-C \|_{L(\mathbb{H},\mathbb{H}')} \\
				&\vee |u^{-, W}_{\R \times \mathbb{H}, \delta, h}(t, x)-
					u^{-, W}_{\R \times \mathbb{H}, \delta, h}(v, y)| \\
				& \qquad \qquad
					\text{if }  u^{-, W}_{\R \times \mathbb{H}, \delta, h}(t, x),
						u^{-, W}_{\R \times \mathbb{H}, \delta, h}(v, y) \in \R\\
				&\infty \\ 
				& \qquad \qquad 
					\text{if } u^{-, W}_{\R \times \mathbb{H}, \delta, h}(t, x),
						u^{-, W}_{\R \times \mathbb{H}, \delta, h}(v, y) \notin \R\\
			\end{cases} 
\end{equation}
and by 
$G^+_{\mathbb{H}, \mathbb{X}, \delta, h}$, 
$
	G^-_{\mathbb{H}, \mathbb{X}, \delta, h} 
	\colon W \times \R \times H' \times \mathbb{S}_{\mathbb{H}, \mathbb{H}'} \to \R
$ 
the functions satisfying for all
$((t,x), r, p, A) \in W \times \R \times H' \times \mathbb{S}_{\mathbb{H}, \mathbb{H}'}$
that
\begin{equation}
	\begin{split}
			&G^+_{\mathbb{H}, \mathbb{X}, \delta, h}
				((t,x), r, p, A) \\
		={}&
			G\Big(
				(t,x), 
				r + \delta h(t,x), 
				p +\delta E_{\mathbb{X}', \mathbb{H}'} 
					\big(
						D_{\mathbb{X}} (h(t,x) |_{(0,T) \times (X \cap O)})
					\big)(t,x) , \\
			&\qquad
				(A|_X) |_X 
				+\delta \big(
					D^2_{\mathbb{X}} (h |_{(0,T) \times (X \cap O)}) (t,x)
				\big)
			\Big)
	\end{split}
\end{equation}
and that 
\begin{equation}
	\begin{split}
			&G^-_{\mathbb{H}, \mathbb{X}, \delta, h}
				((t,x), r, p, A) \\
		={}&
			G\Big(
				(t,x), 
				r - \delta h(t,x), 
				p -\delta E_{\mathbb{X}', \mathbb{H}'} 
					\big(
						D_{\mathbb{X}} (h(t,x) |_{(0,T) \times (X \cap O)})
					\big)(t,x) , 
		\\ & \qquad
				(A|_X) |_X 
				-\delta \big(
					D^2_{\mathbb{X}} (h |_{(0,T) \times (X \cap O)}) (t,x)
				\big)
			\Big)
	\end{split}
\end{equation}
(compare this also with the Section~\ref{ssec:Definition 
of viscosity solutions} above).
\section{Construction of semijets}
\label{ssec: Construction of semijets}
First we recall a 
minor improvement 
(since we don’t assume that $u_i$ is upper
semicontinuous)
of Theorem 3.2 (with $k \leftarrow 2$) in
Crandall, Ishii \citationand~Lions~\cite{CrandallIshiiLions1992}.
It
will be used in
the proof of 
Lemma~\ref{l: suitable semijets matrix} below.
\begin{lemma}[Construction of suitable semijets]
\label{l: suitable semijets finite dim}
	Let for all $i \in \{ 1, 2 \}$
	$N_i \in \N$,
	$\mathcal{O}_i$ be a locally compact subset of $\R^{N_i}$,
	and let for all $i \in \{ 1, 2 \}$
	$u_i \colon \mathcal{O}_i \to [-\infty, \infty)$,
	let $N \in \N$ satisfy that
	$N = N_1 + N_2$,
	let $\mathcal{O}$ be the set satisfying that
	$\mathcal{O} = \mathcal{O}_1 \times \mathcal{O}_2$,
	let $\varphi \in \C_{\R^N}^2(\mathcal{O}, \R)$,
	and let 
	$\hat{x}= (\hat{x}_1, \hat{x}_2) \in \mathcal{O}$
	be a local maximum of
	$O \ni (x_1, x_2) \to u_1(x_1) + u_2(x_2) - \varphi(x_1, x_2) \in \R$.
	Then for all $i \in \{1, 2\}$ and all $\eps > 0$ 
	there exists $X_i \in \mathbb{S}_{\R^{N_i}, (\R^{N_i})'}$ such that 
	$(\pi^{\mathcal{O}}_i D_{\R^{N}}\varphi(\hat{x}), X_i) \in (\hat{J}^2_{\R^{N_i}, +} u_i) (x_i)$
	and that
	\begin{equation}
			- \left( \frac{1}{\eps} + \| (D_{\R^N}^2 \varphi)(\hat{x}) \| \right) I
		\leq
			\left( \begin{array}{ccc}
				X_1 & 0 \\
				0 & X_2
			\end{array} \right) 
		\leq
			(D_{\R^N}^2 \varphi)(\hat{x}) + \eps ((D_{\R^N}^2 \varphi)(\hat{x}))^2.
	\end{equation}
\end{lemma}
\begin{proof}[Proof of Lemma \ref{l: suitable semijets finite dim}]
	The proof follows analogously to Theorem 3.2 
	in Crandall, Ishii \& Lions \cite{CrandallIshiiLions1992} but instead of the
	upper semicontinuity in the last inequality (the one after (A.11) on page 62)
	we use
	\begin{equation}
		\begin{split}
				u_1(0) 
			={}& 
				0 
			= 
				\lim_{n \to \infty} 
					\Big \langle 
						\Big ( \xi_n + \frac{q_n}{\lambda}, 0 \Big ),
						A \Big ( \xi_n + \frac{q_n}{\lambda}, 0 \Big )
					\Big \rangle_{\R^N, (\R^N)'}
			\geq
				\limsup_{n \to \infty} \Big (
					u_1\Big ( \xi_n + \frac{q_n}{\lambda} \Big )
					+ u_2(0)
				\Big ) \\
			={}&
				\limsup_{n \to \infty}
					u_1 \Big ( \xi_n + \frac{q_n}{\lambda} \Big )
			=
				\limsup_{n \to \infty} \Big (
					\hat{u}_1(\xi_n) + \frac{1}{2 \lambda} |q_n|^2
				\Big )
			=
				\hat{u}_1(0)
			\geq
				u_1(0).
		\end{split}
	\end{equation}
	For $u_2$ this follows analogously, 
	which completes the proof of Lemma \ref{l: suitable semijets finite dim}.
\end{proof}
The next lemma generalizes
	Lemma 4 in Lions \cite{Lions1989} to arbitrary linear functions.
\begin{lemma}[Construction of suitable semijets for strict maxima wrt.~to a linear operator]
	\label{l: suitable semijets matrix}
	Let $\mathbb{H} = (H, \langle \cdot, \cdot \rangle_H, \| \cdot \|_H)$
	be a real Hilbert space,
	let 
	$
			\mathbb{V} 
		= 
			(
				V, 
				\langle \cdot, \cdot \rangle_H|_{V^2}, 
				\| \cdot \|_H |_{V}
			)
	$
	be a finite-dimensional linear subspace of $H$,
	let 
	$\mathbb{H}'=(H,\| \cdot \|_{H'})$,
	$(\mathbb{H}^2)'= ((H^2)',\| \cdot \|_{(H^2)'})$,
	$\mathbb{V}'=(V,\| \cdot \|_{V'})$,
	$(\mathbb{V}^2)'= ((V^2)',\| \cdot \|_{(V^2)'})$,
	be the dual spaces,
	let $u_1$, $u_2 \colon H \to [- \infty, \infty)$ be  
	upper semicontinuous
	with respect to the $\| \cdot \|_{\mathbb{H}}$-norm,
	let 
	$
		A \in \mathbb{S}_{\mathbb{H}^2, (\mathbb{H}^2)'}
	$,
	let $(\overline{z}_1$, $\overline{z}_2) \in H^2$ 
	be a strict $\mathbb{H}^2$-maximum of
	$
		H^2 \ni (z_1,z_2) \to u_1(z_1) + u_2(z_2) 
		- \frac 12 \langle (z_1,z_2) , A (z_1, z_2) \rangle_{H^2, (H^2)'} 
			\in [- \infty, \infty)
	$, 
	denote by $E_1$, $E_2 \colon H \to H^2$ the canonical embeddings
	satisfying for all $x \in H$ that
	$E_1(x)=(x,0)$ and that $E_2(x)=(0,x)$,
	and let
	$\eps, \delta, \lambda \in (0, \infty) $
	satisfy that
	\begin{equation}
		\begin{split}
			\lambda \geq{}
				&\frac {1}{\delta}
					\left(
						\|
							\pi^{H'}_{V'} \pi^{(H^2)'}_1 A E_2 \pi^{H}_{V^\perp} 
						\|_{L(\mathbb{H},\mathbb{H}' )}^2 
						\vee \| 
							\pi^{H'}_{(V^\perp)'} \pi^{(H^2)'}_1 A E_2 \pi^{H}_{V}
						\|_{L(\mathbb{H},\mathbb{H}' )}^2
					\right) \\
				&+ \|
						\pi^{H'}_{(V^\perp)'} \pi^{(H^2)'}_1 A E_2 \pi^{H}_{V^\perp} 
					\|_{L(\mathbb{H},\mathbb{H}' )}.
		\end{split}
	\end{equation}
	Then there exist operators 
	$X_{1}$, $X_{2} \in \mathbb{S}_{\mathbb{H}, \mathbb{H}'}$ 
	satisfying that
	$X_{1} = \pi^{H'}_{V'} X_{1} \pi^{H}_V $, $X_{2}=\pi^{H'}_{V'} X_{2} \pi^{H}_V$, that
	\begin{equation}
			-\Big (
				\frac {1}{\eps} 
				+ \| A \|_{L(\mathbb{H}^2, (\mathbb{H}^2)')} 
			\Big ) 
			I_{\mathbb{H}^2}
		\leq
			\left(
				\begin{array}{cc}
					X_{1}
					&
					0
					\\
					0
					&
					X_{2}
				\end{array}
			\right)
		\leq
			\Big ( 
				\pi^{(H^2)'}_{(V^2)'} A \pi^{H^2}_{V^2} 
				+ \eps \big( 
						\pi^{(H^2)'}_{(V^2)'} A \pi^{H^2}_{V^2} \big
				)^2 
			\Big ),
	\end{equation}
	and satisfying that for all $i \in \{ 1, 2\}$ it holds that 
	\begin{equation}
			\left( 
				\pi^{(H^2)'}_i A (\overline{z}_1,\overline{z}_2), 
				X_{i} 
				+ \pi^{(H^2)'}_i A E_i -  \pi^{H'}_{V'} \pi^{(H^2)'}_i A E_i \pi^{H}_{V}   
				+ \delta I_{\mathbb{H}} \pi^{H}_V + \lambda I_{\mathbb{H}} \pi^{H}_{V^\perp}
			\right) 
		\in (\hat{J}^2_{ \mathbb{H}, + } u_i) (\overline{z}_i). 
	\end{equation}	
\end{lemma}
\begin{proof}[Proof of Lemma~\ref{l: suitable semijets matrix}]
	Let $\overline{x}_1$, $\overline{x}_2 \in V$
	be the elements satisfying that
	$\overline{x}_1= \pi^{H}_V (\overline{z}_1)$ and that
	$\overline{x}_2= \pi^{H}_V (\overline{z}_2)$,
	let 
	$\overline{y}_1$, $\overline{y}_2 \in V^\perp$
	be the elements satisfying that
	$\overline{y}_1= \pi^{H}_{V^\perp} (\overline{z}_1)$ and that
	$\overline{y}_2= \pi^{H}_{V^\perp} (\overline{z}_2)$,
	let
	$\hat{z} \in H^2$ be the element satisfying that
	$\hat{z}=\frac{(\overline{z}_1,\overline{z}_2)}{2}$, and let
	$
		B \in \mathbb{S}_{\mathbb{H}^2, (\mathbb{H}^2)'}
	$ 
	be the function satisfying that
	\begin{equation}
	\label{eq: def of B semi}
		B =
			\left (
				\begin{array}{cc}
					-I_{\mathbb{H}} (\delta \pi^{H}_V + \lambda \pi^{H}_{V^\perp})
					& \pi^{(H^2)'}_1 A E_2 
					- \pi^{H'}_{V'} \pi^{(H^2)'}_1 A E_2 \pi^{H}_{V}  \\
					\pi^{(H^2)'}_2 A E_1 
					- \pi^{H'}_{V'} \pi^{(H^2)'}_2 A E_1 \pi^{H}_{V}
					& -I_{\mathbb{H}} (\delta \pi^{H}_V + \lambda \pi^{H}_{V^\perp})
				\end{array}
			\right).
	\end{equation}
	Moreover, denote by $\Phi \colon H^2 \to  [- \infty, \infty)$ 
	the function satisfying
	for all $z_1$, $z_2 \in H$ that
	\begin{equation}
		\Phi(z_1,z_2)= u_1(z_1) + u_2(z_2) 
			- \frac 12 \langle (z_1,z_2) , A (z_1, z_2) \rangle_{H^2, (H^2)'},
	\end{equation}
	denote by $\phi_1$, $\phi_2 \colon H \to  [- \infty, \infty)$ the functions satisfying for all
	$x \in H$ and all $i \in \{1,2\}$ that
	\begin{equation}
	\label{eq: def of phit}
		\begin{split}
				\phi_i(x) 
			={}& 
				u_i(x) 
					- \tfrac 12 \big \langle x , 
						(\pi^{(H^2)'}_i A E_i 
						-\pi^{H'}_{V'} \pi^{(H^2)'}_i A E_i \pi^{H }_V) x \big \rangle_{H, H'}  \\
					&+\Big \langle 
						\big ( \hat{z} -E_i (x) \big ), 
						 B  \big (\hat{z} - E_i (x) \big ) 
					\Big \rangle_{H^2, (H^2)'}    
					+ \tfrac {\delta}{2} \|\pi^{H}_{V} (x)\|_H^2  
					+ \tfrac {\lambda}{2} \|\pi^{H}_{V^\perp} (x)\|_H^2 \\			
		\end{split}
	\end{equation}
	denote by $\tilde{u}_1,$ $\tilde{u}_2 \colon V \to  [- \infty, \infty]$ 
	the functions satisfying for all $x \in V$ that
	\begin{equation}
	\label{eq: def of ut}
		\tilde{u}_1(x) = \sup_{y \in V^\perp} \phi_1(x+y), \qquad
		\tilde{u}_2(x) = \sup_{y \in V^\perp} \phi_2(x+y),
	\end{equation}
	and denote by $\tilde{\Phi} \colon V^2 \to  [- \infty, \infty]$
	the function satisfying
	for all $x_1$, $x_2 \in V$ that
	\begin{equation}
			\tilde{\Phi}(x_{1},x_{2})
		=\tilde{u}_1(x_1) + \tilde{u}_2(x_{2}) 
			- \tfrac 12 \langle (x_1,x_2), A (x_1,x_2) \rangle_{H^2, (H^2)'}.
	\end{equation}
	Furthermore, using the symmetry of B 
	and \eqref{eq: def of B semi} 
	we derive that for every $x_1$, $x_2 \in V$
	and for every $y_1$, $y_2 \in V^\perp$ it holds that 
	\begin{align}
		\nonumber
				&\langle (x_1+y_1,x_2+y_2), B (x_1+y_1,x_2+y_2) \rangle_{H^2, (H^2)'} \\
		\nonumber
			={}& 
				-\delta \langle x_1+y_1, x_1 \rangle_{H} 
				-\delta \langle x_2+y_2, x_2 \rangle_{H}
				-\lambda \langle x_1+y_1, y_1 \rangle_{H} 
				-\lambda \langle x_2+y_2, y_2 \rangle_{H} \\
		\nonumber
				&-2 \Big \langle 
						x_1+y_1, 
						\big ((\pi^{H'}_{V'} 
						+ \pi^{H'}_{(V^\perp)'}) \pi^{(H^2)'}_1 A 
								E_2 (\pi^{H}_V + \pi^{H}_{V^\perp}) 
						- \pi^{H'}_{V'} \pi^{(H^2)'}_1 A E_2 \pi^{H}_V \big ) (x_2+y_2) \Big 
					\rangle_{H, H'}  \\
		\nonumber
			={}&
				-\delta \| x_1 \|_H^2 -\delta \| x_2 \|_H^2
				-\lambda \| y_1 \|_H^2 -\lambda \| y_2\|_H^2 
				-2 \Big \langle 
					x_1+y_1, 
					\big (
						\pi^{H'}_{V'} \pi^{(H^2)'}_1 A E_2 \pi^{H}_{V^\perp} 
		\\
		\begin{split}
			& \qquad \qquad
						+ \pi^{H'}_{(V^\perp)'} \pi^{(H^2)'}_1 A E_2 \pi^{H}_V 
						+ \pi^{H'}_{(V^\perp)'} \pi^{(H^2)'}_1 A E_2 \pi^{H}_{V^\perp} 
					\big ) (x_2 + y_2)
				\Big \rangle_{H, H'} \\
			={}&
				-\delta \| x_1 \|_H^2 -\delta \| x_2\|_H^2
				-\lambda \| y_1 \|_H^2 -\lambda \| y_2\|_H^2
				-2 \Big \langle 
						x_1, 
						\pi^{H'}_{V'} \pi^{(H^2)'}_1 A E_2 \pi^{H}_{V^\perp} y_2 \Big 
					\rangle_{H, H'}
		\end{split} \\
		\nonumber
				&-2\Big \langle 
						y_1, 
						\pi^{H'}_{(V^\perp)'} \pi^{(H^2)'}_1 A E_2 \pi^{H}_V x_2 
					\Big \rangle_{H, H'} 
				-2 
				\Big \langle 
					y_1, \pi^{H'}_{(V^\perp)'} \pi^{(H^2)'}_1 A E_2 \pi^{H}_{V^\perp} y_2 
				\Big \rangle_{H, H'} \\
			\nonumber
			\leq{}&
				-\delta \| x_1 \|_H^2 -\delta \| x_2\|_H^2
				-\lambda \| y_1 \|_H^2 -\lambda \| y_2\|_H^2
				+2\| 
						\pi^{H'}_{V'} \pi^{(H^2)'}_1 A E_2 \pi^{H}_{V^\perp} 
					\|_{L(\mathbb{H}, \mathbb{H}')} 
					\cdot \|y_2\|_H \, \| x_1 \|_H \\
		\nonumber
			&+2 \| 
							\pi^{H'}_{(V^\perp)'} \pi^{(H^2)'}_1 A E_2 \pi^{H }_V 
						\|_{L(\mathbb{H}, \mathbb{H}')} 
					\cdot \|x_2 \|_H \, \| y_1\|_H 
		\\& \nonumber
				+2\| 
						\pi^{H'}_{(V^\perp)'} \pi^{(H^2)'}_1 A E_2 \pi^{H}_{V^\perp} 
					\|_{L(\mathbb{H}, \mathbb{H}')} 
					\cdot \| y_1 \|_H  \, \| y_2 \|_H
	\end{align}
	Hence we obtain that for every $x_1$, $x_2 \in V$
	and for every $y_1$, $y_2 \in V^\perp$ it holds that
	\begin{align}
	\nonumber
				&\langle (x_1+y_1,x_2+y_2), B (x_1+y_1,x_2+y_2) \rangle_{H^2, (H^2)'} \\ \nonumber
			\leq{}&
				-\left (\sqrt{\delta} \cdot \| x_1 \|_H 
					- \frac
							{
								\| 
									\pi^{H'}_{V'} \pi^{(H \times H)'}_1 A E_2 \pi^{H }_{V^\perp}
								\|_{L(\mathbb{H}, \mathbb{H}')}
							}
							{\sqrt{\delta}} 
						\| y_2 \|_{H}
				\right)^2		
		\\& \nonumber
				-\left( \sqrt{\delta} \cdot \| x_2 \|_{H} 
					- \frac
							{
								\| 
									\pi^{H'}_{(V^\perp)'} \pi^{(H^2)'}_1 A E_2 \pi^{H }_V
								\|_{L(\mathbb{H}, \mathbb{H}')}
							}
							{\sqrt{\delta}} 
						\| y_1 \|_{H} 
				\right)^2 \\ 
				&-\| 
						\pi^{H'}_{(V^\perp)'} \pi^{(H^2)'}_1 A E_2 \pi^{H}_{V^\perp} 
					\|_{L(\mathbb{H}, \mathbb{H}')} 
					\left ( \| y_1 \|_H - \| y_2 \|_H \right )^2 \\ \nonumber
				&- \left (\lambda
						-\frac
							{
								\| 
									\pi^{H'}_{(V^\perp)'} \pi^{(H^2)'}_1 A E_2 \pi^{H}_V 
								\|_{L(\mathbb{H}, \mathbb{H}')}^2}
							{\delta}
						- \| 
								\pi^{H'}_{(V^\perp)'} \pi^{(H^2)'}_1 A E_2 \pi^{H}_{V^\perp} 
							\|_{L(\mathbb{H}, \mathbb{H}')} 
					\right ) \| y_1 \|_H^2 \\ \nonumber
				&-\left( \lambda 
						-\frac
							{
								\| 
									\pi^{H'}_{V'} \pi^{(H^2)'}_1 A E_2 \pi^{H}_{V^\perp} 
								\|_{L(\mathbb{H}, \mathbb{H}')}^2
							}
							{\delta}
						- \| 
								\pi^{H'}_{(V^\perp)'} \pi^{(H^2)'}_1 A E_2 \pi^{H}_{V^\perp} 
							\|_{L(\mathbb{H}, \mathbb{H}')}
					\right )\| y_2\|_H^2 
			\leq 
				0. 
	\end{align}
	Thus $B$
	is negative semi-definite. 
	Moreover, we have for all $z_1$, $z_2 \in H$ that
	\begin{align}
		\nonumber
				&-\langle 
						(z_1,0), 
						\big( 
							A - \pi^{(H^2)'}_{(V^2)'} A \pi^{H^2}_{V^2}
						\big) (0,z_2) 
					\rangle_{H^2, (H^2)'} 
				+ 2 \left \langle 
							\left ( \hat{z} - \tfrac{(z_1,z_2)}{2} \right),
							B \left ( \hat{z} - \tfrac{(z_1,z_2)}{2} \right) 
						\right \rangle_{H^2, (H^2)'}  \\
		\nonumber
			={}&
				- \langle (z_1,0), B (0,z_2) \rangle_{H^2, (H^2)'}  
				+ 2 \left \langle 
							\hat{z}, B \, \hat{z} 
						\right \rangle_{H^2, (H^2)'} 
				- 4 \left \langle 
							\tfrac{(z_1,z_2)}{2} , B \, \hat{z}
						\right \rangle_{H^2, (H^2)'} \\
		\nonumber
				&+ 2 \left \langle 
							\tfrac{(z_1,z_2)}{2} , B \left( \tfrac{(z_1,z_2)}{2} \right ) 
						\right \rangle_{H^2, (H^2)'}    \\
		\nonumber
			={}&
				\left \langle 
					\hat{z}, B \, \hat{z} 
				\right \rangle_{H^2, (H^2)'} 
				- 2 \left \langle 
							(z_1,0) , B \, \hat{z} 
						\right \rangle_{H^2,(H^2)'}
				+ \big \langle (z_1,0), B (z_1,0) \big \rangle_{H^2 (H^2)'}\\
		\nonumber
				&+ \left \langle 
						\hat{z}, B \, \hat{z} 
					\right \rangle_{H^2, (H^2)'} 
				- 2 \left \langle 
							(0,z_2) , B \, \hat{z} 
						\right \rangle_{H^2, (H^2)'}
				+ \big \langle (0,z_2), B (0,z_2) \big \rangle_{H^2, (H^2)'} \\
		\label{eq: parallel equality} 
				&+ \tfrac 12 
					\big \langle 
						(z_1,0)+(0,z_2) , B \big ( (z_1,0)+(0,z_2) \big ) 
					\big \rangle_{H^2, (H^2)'} 
				- \langle (z_1,0), B (0,z_2) \rangle_{H^2, (H^2)'} \\
		\nonumber
				&- \big \langle (z_1,0), B (z_1,0) \big \rangle_{H^2, (H^2)'}
				- \big \langle (0, z_2), B (0, z_2) \big \rangle_{H^2, (H^2)'}\\
		\nonumber
			={}&
				\big \langle 
					(\hat{z} -(z_1,0)), B (\hat{z} -(z_1,0)) 
				\big \rangle_{H^2, (H^2)'}
				+\big \langle 
					(\hat{z} -(0,z_2)), B (\hat{z} -(0,z_2)) 
				\big \rangle_{H^2, (H^2)'} \\
		\nonumber
				&- \tfrac 12 \big \langle 
						(z_1,0) , B (z_1,0)
					\big \rangle_{H^2, (H^2)'} 
				- \tfrac 12 \big \langle 
						(0,z_2) , B (0,z_2) \big 
					\rangle_{H^2, (H^2)'}  \\
		\nonumber
			={}&
				\big \langle 
					(\hat{z} -(z_1,0)), B (\hat{z} -(z_1,0)) 
				\big \rangle_{H^2, (H^2)'}
				+\big \langle 
					(\hat{z} -(0,z_2)), B (\hat{z} -(0,z_2)) 
				\big \rangle_{H^2, (H^2)'} \\
		\nonumber
				&+\tfrac {\delta}{2} \|\pi^{H}_V (z_1)\|_H^2 
				+ \tfrac {\lambda}{2} \|\pi^{H}_{V^\perp} (z_1)\|_H^2 
				+\tfrac {\delta}{2} \|\pi^{H}_V (z_2)\|_H^2 
				+ \tfrac {\lambda}{2} \|\pi^{H}_{V^\perp} (z_2)\|_H^2.
	\end{align}
	This implies
	that for all $x_1$, $x_2 \in V$ and all $y_1, y_2 \in V^\perp$ it holds that
\begin{align}
		\nonumber
			&\Phi(\overline{z}_1,\overline{z}_2)
		\geq{}
			\Phi(x_1+y_1,x_2+y_2) \\
		\nonumber
		\geq{}&
			\Phi(x_1+y_1,x_2+y_2)
			+ 2 \Big \langle 
						\Big ( \hat{z} - \tfrac{(x_1+y_1,x_2+y_2)}{2} \Big),
						B \Big ( \hat{z} - \tfrac{(x_1+y_1,x_2+y_2)}{2} \Big) 
					\Big \rangle_{H^2, (H^2)'} \\
		\nonumber
		={}& 
			u_1(x_1+y_1) + u_2(x_2+y_2) 
			+ 2 \Big \langle 
						\Big ( \hat{z} - \tfrac{(x_1+y_1,x_2+y_2)}{2} \Big),
						B \Big ( \hat{z} - \tfrac{(x_1+y_1,x_2+y_2)}{2} \Big)
					\Big \rangle_{H^2, (H^2)'}  \\
			\nonumber
			&- \tfrac 12 
				\big \langle 
					(x_1+y_1,x_2+y_2) , A (x_1+y_1,x_2+y_2) 
				\big \rangle_{H^2, (H^2)'} \\	
		\nonumber
		={}&
			u_1(x_1+y_1) + u_2(x_2+y_2) 
			+ 2 \Big \langle 
						\Big ( \hat{z} - \tfrac{(x_1+y_1,x_2+y_2)}{2} \Big),
						B \Big ( \hat{z} - \tfrac{(x_1+y_1,x_2+y_2)}{2} \Big) 
					\Big \rangle_{H^2, (H^2)'}  \\
			\nonumber
			&- \tfrac 12 
				\big \langle 
					(x_1+y_1,x_2+y_2) , 
					(A-\pi^{(H^2)'}_{(V^2)'} A \pi^{H^2}_{V^2}) 
						(x_1+y_1,x_2+y_2) 
				\big \rangle_{H^2, (H^2)'} \\
			\nonumber
			&- \tfrac 12 
				\big \langle 
					(x_1+y_1,x_2+y_2) , 
					(\pi^{(H^2)'}_{(V^2)'} A \pi^{H^2}_{V^2}) 
						(x_1+y_1,x_2+y_2) 
				\big \rangle_{H^2, (H^2)'} \\
		\nonumber
		={}&
			u_1(x_1+y_1) + u_2(x_2+y_2) 
			+ 2 \Big \langle 
						\Big ( \hat{z} - \tfrac{(x_1+y_1,x_2+y_2)}{2} \Big),
						B \Big ( \hat{z} - \tfrac{(x_1+y_1,x_2+y_2)}{2} \Big) 
					\Big \rangle_{H^2, (H^2)'}  \\
		\label{eq: phi and phit}
			&- \tfrac 12 
				\big \langle 
					(x_1+y_1,0) , 
					(A-\pi^{(H^2)'}_{(V^2)'} A \pi^{H^2}_{V^2}) (x_1+y_1,0) 
				\big \rangle_{H^2, (H^2)'} \\
			\nonumber
			&	- \tfrac 12 
				\big \langle 
					(0,x_2+y_2) , 
					(A-\pi^{(H^2)'}_{(V^2)'} A \pi^{H^2}_{V^2}) (0,x_2+y_2) 
				\big \rangle_{H^2, (H^2)'} \\  
			\nonumber
			& -\big \langle 
					(x_1+y_1,0) , 
					(A-\pi^{(H^2)'}_{(V^2)'} A \pi^{H^2}_{V^2}) (0,x_2+y_2) 
				\big \rangle_{H^2, (H^2)'} 
			- \tfrac 12 
				\big \langle 
					(x_1,x_2) , A  (x_1,x_2) 
				\big \rangle_{H^2, (H^2)'} \\
		\nonumber
		={}&
			u_1(x_1+y_1) + u_2(x_2+y_2)  
			+\Big \langle 
					\Big (\hat{z} -\big (x_1+y_1,0 \big ) \Big ), 
					B  \Big (\hat{z} -\big (x_1+y_1,0 \big ) \Big ) 
				\Big \rangle_{H^2, (H^2)'} \\
		\nonumber
			& + \Big \langle 
					\Big (\hat{z} -\big (0,x_2+y_2 \big ) \Big ), 
					B \Big (\hat{z} -\big (0,x_2+y_2 \big ) \Big ) 
				\Big \rangle_{H^2, (H^2)'}
			+ \tfrac {\delta}{2} \|x_1\|_H^2 + \tfrac {\lambda}{2} \|y_1\|_H^2
			+ \tfrac {\delta}{2} \|x_2\|_H^2 \\
	\nonumber
			&+ \tfrac {\lambda}{2} \|y_2\|_H^2
			- \tfrac 12 
				\big \langle 
					(x_1+y_1) , 
					(
						\pi^{(H^2)'}_1 A E_1
						-\pi^{H'}_{V'} \pi^{(H^2)'}_1 A E_1 \pi^{H}_{V}
					) (x_1+y_1) 
				\big \rangle_{H^2, (H^2)'} \\ 
		\nonumber
			& - \tfrac 12 
				\big \langle 
					(x_2+y_2) , 
					(
						\pi^{(H^2)'}_2 A E_2
						-\pi^{H'}_{V'} \pi^{(H^2)'}_2 A E_2 \pi^{H}_{V}
					) (x_2+y_2) 
				\big \rangle_{H^2, (H^2)'} \\
		\nonumber
			&- \tfrac 12 
				\big \langle 
					(x_1,x_2) , A  (x_1,x_2) 
				\big \rangle_{H^2, (H^2)'} \\
		\nonumber
		={}&
			\phi_1(x_1+y_1) + \phi_2(x_2+y_2)
			- \tfrac 12 
				\big \langle 
					(x_1,x_2) , A  (x_1,x_2) 
				\big \rangle_{H^2, (H^2)'}.
	\end{align}
	From (\ref{eq: phi and phit}) it follows that
	$\tilde{u}_1$ and $\tilde{u}_2$ are well-defined 
	and that it holds 
	\begin{equation}
	\label{eq: same max}
		\begin{split}
				\tilde{\Phi}(\overline{x}_1,\overline{x}_2)
			={}&
				\tilde{u}_1(\overline{x}_1)+\tilde{u}_2(\overline{x}_2)
				-\frac 12 
					\langle 
						(\overline{x}_1,\overline{x}_2), A (\overline{x}_1,\overline{x}_2) 
					\rangle_{H^2, (H^2)'} \\
			={}&
				\sup_{y_1 \in V^\perp} \sup_{y_2 \in V^\perp} 
				\big (
					\phi_1(\overline{x}_1+y_1) + \phi_2(\overline{x}_2+y_2) 
					-\frac 12 
						\langle 
							(\overline{x}_1,\overline{x}_2), A (\overline{x}_1,\overline{x}_2) 
						\rangle_{H^2, (H^2)'} 
				\big ) \\
			\geq{}&
				\phi_1(\overline{x}_1+\overline{y}_1) 
				+ \phi_2(\overline{x}_2+\overline{y}_2) 
					-\frac 12 
						\langle 
							(\overline{x}_1,\overline{x}_2), A (\overline{x}_1,\overline{x}_2)
						\rangle_{H^2, (H^2)'}
			= \Phi(\overline{z}_1,\overline{z}_2).
		\end{split}
	\end{equation} 
	Thus we have that
	\begin{equation}
	\label{eq: max of phi}
			\tilde{u}_1(\overline{x}_1)
		=
			\phi_1(\overline{z}_1),
		\qquad
			\tilde{u}_2(\overline{x}_2)
		=
			\phi_2(\overline{z}_2)
	\end{equation}
	and that
	$(\overline{x}_1, \overline{x}_2)$ is a global maximum of
	$\tilde{\Phi}$.
	Lemma \ref{l: suitable semijets finite dim}
		(with $\varphi \leftarrow \tilde{\Phi}$)
	implies the existence of
	$\tilde{X}_{1}, \tilde{X}_{2} \in \mathbb{S}_{\mathbb{V}, \mathbb{V}'}$
	(which we fix for the rest of the proof)
	such that it holds that
	\begin{equation}
	\label{eq: matrix inequality for Xt}
			-\Big(
				\frac {1}{\eps} 
				+ \big \| 
						\pi^{(H^2)'}_{(V^2)'} A \big |_{V^2} 
					\big \|_{L(\mathbb{V}^2, (\mathbb{V}^2)')} 
			\Big) 
				I_{\mathbb{V}^2}
		\leq
			\left(
				\begin{array}{cc}
					\tilde{X}_{1}
					&
					0
					\\
					0
					&
					\tilde{X}_{2}
				\end{array}
			\right) 
		\leq
			\Big ( \pi^{(H^2)'}_{(V^2)'} A \big |_{V^2} 
				+\eps \big( 
						\pi^{(H^2)'}_{(V^2)'} A \big |_{V^2} 
				\big )^2 
			\Big),   
	\end{equation}
	\begin{equation}
		\label{eq: semijets hat}
		(
			(\pi^{H'}_{V'} \pi^{(H^2)'}_1 A (\overline{x}_1,\overline{x}_2)) |_V, 
			\tilde{X}_{1}
		) 
			\in (\hat{J}^2_{\mathbb{V}, + } \tilde{u}_1)(\overline{x}_1), \quad
		(
			(\pi^{H'}_{V'} \pi^{(H^2)'}_2 A (\overline{x}_1,\overline{x}_2))|_V, 
			\tilde{X}_{2}
		) 
			\in (\hat{J}^2_{\mathbb{V}, + } \tilde{u}_2)(\overline{x}_2).
	\end{equation}
	From (\ref{eq: semijets hat}) it follows that
	there exist sequences $(x_{1,k})_{k \in \N}$, $(x_{2,k})_{k \in \N}$ $\subseteq V$,
	$(\xi_{1,k})_{k \in \N}, (\xi_{2,k})_{k \in \N}$ $\subseteq V'$ 
	and 
	$
		(\tilde{X}_{1,k})_{k \in \N}, (\tilde{X}_{2,k})_{k \in \N} 
			\subseteq \mathbb{S}_{\mathbb{V}, \mathbb{V}'}
	$ 
	(which we fix for the rest of the proof)
	such that
	\begin{align}
	\label{eq: k convergence for x}
		&\lim_{k \to \infty}
			\| x_{1,k} - \overline{x}_{1} \|_H  = 0, \qquad 
		\lim_{k \to \infty}
			\| x_{2,k} - \overline{x}_{2} \|_H = 0,                  \\ 
	\label{eq: k convergence for ut}
		&\lim_{k \to \infty}
			\tilde{u}_1(x_{1,k}) 
			= \tilde{u}_1(\overline{x}_{1}), \qquad 
		\lim_{k \to \infty}
			\tilde{u}_2(x_{2,k}) 
			= \tilde{u}_2(\overline{x}_{2}),                    \\
	\label{eq: k convergence for xi1 in V}
		&\lim_{k \to \infty}
			\| 
				\xi_{1,k} 
				- (\pi^{H'}_{V'} \pi^{(H^2)'}_1 A (\overline{x}_1,\overline{x}_2))|_V 
			\|_{V'} = 0, \qquad \\
	\label{eq: k convergence for xi2 in V}
		&\lim_{k \to \infty}
			\| 
				\xi_{2,k} 
				- (\pi^{H'}_{V'} \pi^{(H^2)'}_2 A (\overline{x}_1,\overline{x}_2))|_V 
			\|_{V'} = 0,\\
		&\lim_{k \to \infty}
	\label{eq: k convergence for Xt}
			\| \tilde{X}_{1,k} - \tilde{X}_{1} \|_{L(\mathbb{V}, \mathbb{V})} = 0, \qquad
		\lim_{k \to \infty}
			\| \tilde{X}_{2,k} - \tilde{X}_{2} \|_{L(\mathbb{V}, \mathbb{V})} = 0, \\
	\label{eq: semijets}
		&(\xi_{1,k},\tilde{X}_{1,k}) \in (J^2_{\mathbb{V}, + } \tilde{u}_1)(x_{1,k}), \quad
		(\xi_{2,k},\tilde{X}_{2,k}) \in (J^2_{\mathbb{V}, + } \tilde{u}_2)(x_{2,k}).
	\end{align}
	Now denote by $X_{1}$, $X_{2} \in \mathbb{S}_{\mathbb{H}, \mathbb{H}'}$ 
	the functions satisfying
	for all $x_1$, $x_2 \in H$ that
	\begin{equation}
		\begin{split}
				&\langle X_{1} (x_1), x_2 \rangle_{H', H} 
			= 
				\langle \tilde{X}_{1} (\pi^{H}_V (x_1)), \pi^{H}_V (x_2) \rangle_{V', V}, \\
				&\langle X_{2} (x_1), x_2 \rangle_{H', H} 
			= 
				\langle \tilde{X}_{2} (\pi^{H}_V (x_1)), \pi^{H}_V (x_2) \rangle_{V', V}
		\end{split}
	\end{equation} 
	and for all $k \in \N$ denote by 
	$X_{1,k}$, $X_{2,k} \in \mathbb{S}_{\mathbb{H}, \mathbb{H}'}$ 
	the functions satisfying
	for all $x_1$, $x_2 \in H$ that
	\begin{equation}
	\label{eq: Def of X}
		\begin{split}
				&\langle X_{1,k} (x_1), x_2 \rangle_{H', H} 
			= 
				\langle \tilde{X}_{1,k} (\pi^{H}_V (x_1)), \pi^{H}_V (x_2) \rangle_{V', V}, \\
				&\langle X_{2,k} (x_1), x_2 \rangle_{H', H} 
			= 
				\langle \tilde{X}_{2,k} (\pi^{H}_V (x_1)), \pi^{H}_V (x_2) \rangle_{V', V}.
		\end{split}
	\end{equation}
	Then (\ref{eq: Def of X}), 
	(\ref{eq: k convergence for Xt}), 
	(\ref{eq: matrix inequality for Xt}) and 
	$
			\| 
				\pi^{(H^2)'}_{(V^2)'} A \big |_{V^2} 
			\|_{L(\mathbb{V}^2, (\mathbb{V}^2)')}
		\leq 
			\| A \|_{L(\mathbb{H}^2, (\mathbb{H}^2)')}
	$
	imply
	\begin{align}
	\label{eq: X k finite dimensional}
		&X_{1,k} = \pi^{H'}_{V'} X_{1,k} \pi^{H}_V, 
		\qquad X_{2,k}= \pi^{H'}_{V'} X_{2,k} \pi^{H}_V, \\
	\label{eq: k convergence for X}
		&\lim_{k \to \infty}
			\|
				X_{1,k} - X_{1}
			\|_{L(\mathbb{H}, \mathbb{H}')}
		= 0, 
	\qquad
		\lim_{k \to \infty}
			\|
				X_{2,k} - X_{2}
			\|_{L(\mathbb{H}, \mathbb{H}')}
		= 0 \\
	\label{eq: matrix inequality}
		&-\left(
			\frac {1}{\eps} 
			+ \| A \|_{L(\mathbb{H}^2, (\mathbb{H}^2)')} 
		\right) I_{\mathbb{H}^2}
		\leq
			\left(
				\begin{array}{cc}
					X_{1}
					&
					0
					\\
					0
					&
					X_{2}
				\end{array}
			\right)
		\leq
			\left( 
				\pi^{(H^2)'}_{(V^2)'} A \pi^{H^2}_{V^2} 
				+ \eps (\pi^{(H^2)'}_{(V^2)'} A \pi^{H^2}_{V ^2})^2 
			\right).
	\end{align}
	Moreover, \eqref{eq: k convergence for xi1 in V},
	\eqref{eq: k convergence for xi2 in V}, and the definition of
	$\pi^{H'}_{V'}$ show for all $i \in \{1,2\}$ that
	\begin{align}
	\label{eq: k convergence for xi}
		\begin{split}
				&\lim_{k \to \infty}
					\| 
						\xi_{i,k} \pi^{H}_V 
						- \pi^{H'}_{V'} \pi^{(H^2)'}_i A (\overline{x}_1,\overline{x}_2) 
					\|_{H'} 
			=
				\lim_{k \to \infty}
					\| 
						\xi_{i,k} \pi^{H}_V 
						- \pi^{H'}_{V'} \pi^{(H^2)'}_i A (\overline{x}_1,\overline{x}_2) \pi^{H}_V 
					\|_{H'} 
			= 
				0.
		\end{split}
	\end{align}
	To conclude the theorem we pull back the expansion implied 
	by~(\ref{eq: semijets}) 
	on $\tilde{u}_1$ or $\tilde{u}_2$ to the level of $u_1$ and $u_2$.  
	We will do this only for 
	$\tilde{u}_1$ since for $\tilde{u}_2$ it follows exactly the same way.
	From (\ref{eq: semijets}) 
	and Lemma \ref{lem:semijets}
	we deduce for all $k \in \N$ the existence of
	$\varphi_{k} \in C_{\mathbb{V}}^2(V,\R)$ 
	satisfying for all $x \in V$ that
	\begin{equation}
	\label{eq: def of varphi}
		\tilde{u}_1(x) \leq \varphi_{k}(x),
		\quad \varphi_{k}(x_{1,k}) = \tilde{u}_1(x_{1,k}), 
		\quad (D_{\mathbb{V}}\varphi_{k})(x_{1,k}) = \xi_{1,k},
		\quad (D_{\mathbb{V}}^2\varphi_{k})(x_{1,k}) = (X_{1,k} |_V)|_V.
	\end{equation}
	Without loss of generality we may assume that for all $k \in \N$ it holds that
	$x_{1,k}$ is a strict $\mathbb{V}$-maximum of 
	$V \ni x \to \tilde{u}_1(x) - \varphi_{k}(x) \in  [- \infty, \infty)$ 
	and that for all $k \in \N$ it holds that
	$\varphi_{k}$ grows at least not slower than 
	$
		V \ni x \to 
			\frac 12 \|A\|_{L(\mathbb{H}^2, (\mathbb{H}^2)')} 
				\, \| x \|_H^2 \in \R
	$ 
	at infinity
	in the sense that there exists an $M_k \in (0, \infty)$ such that
	for all $x \in H$ with $\| x \|_H \geq M_k$ it holds that
	\begin{equation}
	\label{eq: behaviour of varphi}
			\varphi_k(x) 
		\geq 
			\frac 12 \|
				A
			\|_{L(\mathbb{H}^2, (\mathbb{H}^2)')} 
				\, \|x\|_H^2,
	\end{equation}
	otherwise take a function $\eta \in C_{\mathbb{V}}^2(V,[0,1])$ satisfying for all
	$x \in V$ with $\| x \|_H \leq 1$ that $\eta(x)=1$ and for all
	$x \in V$ with $\| x \|_H \geq 2$ that $\eta(x)=0$
	and replace $\varphi_k$ by the function
	$
		V \ni x \to 
			\eta(x-x_{1,k}) (\varphi_k(x) + \|x-x_{1,k}\|_H^4)
			+(1-\eta(x-x_{1,k})) 
				( \frac 12 \|
						A
					\|_{L(\mathbb{H}^2, (\mathbb{H}^2)')} 
						\, \|x\|_H^2 
					+ (\varphi_k(x))^2 + 1) \in \R
	$.
	From (\ref{eq: behaviour of varphi}) we obtain 
	for all $k \in \N$ and all $x \in H$ that 
	\begin{align}
		\begin{split}
				&- \1_{ \{ \| \pi^{H}_V (x) \|_H \geq M_k \} } \varphi_k(\pi^{H}_V (x)) 
				- \tfrac 12 \| A \|_{L(\mathbb{H}^2, (\mathbb{H}^2)')} 
						\, M_k^2 \\
			\leq{}&
				- \1_{ \{ \| \pi^{H}_V (x) \|_H \geq M_k \} } \varphi_k(\pi^{H}_V (x)) 
				- \1_{ \{ \| \pi^{H}_V (x) \|_H < M_k \} } 
					\tfrac 12 \| A \|_{L(\mathbb{H}^2, (\mathbb{H}^2)')} 
						\, \| \pi^{H}_V (x) \|_H^2 \\
			\leq{}&
				- \tfrac 12 
					\| A \|_{L(\mathbb{H}^2, (\mathbb{H}^2)')}  
						\, \| \pi^{H}_V (x) \|_H^2 
			\leq
				- \tfrac 12 \langle 
						\pi^{H}_V (x) , \pi^{(H^2)'}_1 A E_1 \pi^{H}_V (x) 
					\rangle_{H, H'} \\
			={}&
				- \tfrac 12 \langle 
						x , \pi^{H'}_{V'} \pi^{(H^2)'}_1 A E_1 \pi^{H}_V x 
					\rangle_{H, H'}.  
		\end{split}
	\end{align}
	Therefore the fact that $(\overline{z}_1,\overline{z}_2)$ is a
	maximum of  $\Phi$ ensures for all $k \in \N$ and all $x \in H$ that 
	\begin{equation}
	\label{eq: control first term of phi}
		\begin{split}
				&u_1(x)
					- \frac 12 \langle 
							x,
							(
								\pi^{(H^2)'}_1 A E_1
								- \pi^{H'}_{V'} \pi^{(H^2)'}_1 A E_1 \pi^{H}_V
							) x 
						\rangle_{H, H'} 
					- \1_{ \{ \| \pi^{H}_V (x) \|_H \geq M_k \} } \varphi_k(\pi^{H}_V (x)) \\
					&- \frac 12 \| A \|_{L(\mathbb{H}^2, (\mathbb{H}^2)')} 
							\, M_k^2 \\
			\leq{}&
				u_1(x)
					- \frac 12 \langle 
							x,
							(
								\pi^{(H^2)'}_1 A E_1
								- \pi^{H'}_{V'} \pi^{(H^2)'}_1 A E_1 \pi^{H}_V
							) x 
						\rangle_{H, H'} 
					- \frac 12 \langle 
							x , \pi^{H'}_{V'} \pi^{(H^2)'}_1 A E_1 \pi^{H}_V x 
						\rangle_{H, H'}  \\
			={}&
				u_1(x)
					- \frac 12 \langle (x, 0) , 
						A (x,0) \rangle_{H^2, (H^2)'} 
			=
				\Phi(x,0) - u_2(0)
			\leq
				\Phi(\overline{z}_1,\overline{z}_2) - u_2(0).
		\end{split}
	\end{equation}
	Thus equation (\ref{eq: def of phit})
	and inequality (\ref{eq: control first term of phi}) imply 
	for all $k \in \N$ and all $x \in H$ that 
	\begin{equation}
	\label{eq: behavior at infty without lin}
		\begin{split}
				&\phi_1(x) - \varphi_{k}(\pi^{H}_V (x))  \\
			={}& 
				u_1(x)
				- \tfrac 12 \langle 
						x, 
						(	
							\pi^{(H^2)'}_1 A E_1
							- \pi^{H'}_{V'} \pi^{(H^2)'}_1 A E_1 \pi^{H}_V
						) x 
					\rangle_{H, H'}  
				+ \Big \langle 
						\big (\hat{z} - (x,0  ) \big ), 
						B  \big (\hat{z} (x,0  ) \big ) 
					\Big \rangle_{H^2, (H ^2)'} \\
				&+ \tfrac {\delta}{2} \|\pi^{H}_V (x) \|_H^2 
				+ \tfrac {\lambda}{2} \|\pi^{H}_{V^\perp} (x) \|_H^2
				- \1_{ \{ \| \pi^{H}_V (x) \|_H \geq M_k \} } \varphi_k(\pi^{H}_V (x))
				- \1_{ \{ \| \pi^{H}_V (x) \|_H < M_k \} } \varphi_k(\pi^{H}_V (x)) \\
				&- \frac 12 \| 
						A 
					\|_{L(\mathbb{H}^2, (\mathbb{H}^2)')}
						\, M_k^2 
				+ \frac 12 \| 
						A 
					\|_{L(\mathbb{H}^2, (\mathbb{H}^2)')} 
						\, M_k^2\\
			\leq{}&
				\Phi(\overline{z}_1,\overline{z}_2) -u_2(0) 
				+\tfrac {\delta}{2} \|\pi^{H}_V (x)\|_H^2 
				+ \tfrac {\lambda}{2} \|\pi^{H}_{V^\perp} (x)\|_H^2
				+ \langle (x,0), B (x,0) \rangle_{H^2, (H^2)'}
				+ \langle \hat{z}, B \, \hat{z} \rangle_{H^2, (H^2)'}	\\
				&-2 \langle \hat{z}, B (x,0) \rangle_{H^2, (H^2)'}  
				- \1_{ \{ \| \pi^{H}_V (x) \|_H < M_k \} } \varphi_k(\pi^{H}_V (x))
				+ \tfrac 12 \| 
						A 
					\|_{L(\mathbb{H}^2, (\mathbb{H}^2)')} \, M_k^2 \\
			\leq{}&
				\Phi(\overline{z}_1,\overline{z}_2) -u_2(0) 
				-\tfrac {\delta}{2} \|\pi^{H}_V (x)\|_H^2 
				- \tfrac {\lambda}{2} \|\pi^{H}_{V^\perp} (x)\|_H^2
				+\| 
						B 
				\|_{L(\mathbb{H}^2, (\mathbb{H}^2)')}
						\, \| \hat{z} \|_{H^2}^2 \\
				&+ 2 \|
						B 
					\|_{L(\mathbb{H}^2, (\mathbb{H}^2)')} 
						\, \| \hat{z} \|_{H^2} \, \|x\|_H 
				- \1_{ \{ \| \pi^{H}_V (x) \|_H < M_k \} } \varphi_k(\pi^{H}_V (x))
				+ \tfrac 12 \| 
						A 
					\|_{L(\mathbb{H}^2, (\mathbb{H}^2)')}
						\, M_k^2.
		\end{split}
	\end{equation}
	From (\ref{eq: behavior at infty without lin}) and
	$\delta$, $\lambda >0$ 
	it follows then that for all $k \in \N$
	it holds that 
	\begin{equation}
	\label{eq: convergence at infty whithout lin}
			\lim_{r \to \infty} \sup 
			\big \{ 
				\phi_1(x) - \varphi_{k}(\pi^{H}_V (x)) 
					\colon x \in H,~ \|x\|_H \geq r
			\big \}
		= -\infty.
	\end{equation}
	In a finite-dimensional Hilbert space, 
	(\ref{eq: convergence at infty whithout lin}) would be enough to conclude 
	the existence of a global maximum but in an infinite-dimensional Hilbert space
	we need a perturbation argument to get the existence of a maximum.
	Since $u_1$ is upper semicontinuous 
	with respect to the $\| \cdot \|_{\mathbb{H}}$-norm
	and $A, B \in \mathbb{S}_{\mathbb{H} \times \mathbb{H}, (\mathbb{H} \times \mathbb{H})'}$
	it follows from
	(\ref{eq: def of phit}) that $\phi_1$ is upper semicontinuous 
	with respect to the $\| \cdot \|_{\mathbb{H}}$-norm
	and this together with
	the fact that $\varphi_k$ is continuous 
	with respect to the $\| \cdot \|_{\mathbb{H}}$-norm
	implies that also
	$H \ni x \to \phi_1(x) - \varphi_{k}(\pi^{H}_V (x))$ is upper semicontinuous
	with respect to the $\| \cdot \|_{\mathbb{H}}$-norm.
	Furthermore, it follows from the fact that for all $r \in (0, \infty)$ it holds that
	$\{ x \in H \colon \| x\|_H \leq r \}$ is convex, $\mathbb{H}$-bounded, and
	closed that it is weakly compact. Therefore it follows from page $3$
	in Stegall \cite{Stegall1978} that for all $r \in (0, \infty)$
	$\{ x \in H \colon \| x\|_H \leq r \}$ is also an RNP set 
	and this together with
	the theorem starting on page $4$ in Stegall \cite{Stegall1978} implies that for all  
	$r \in (0, \infty)$, $k \in \N$, and all $\gamma \in (0,1)$
	there exist $p_{\gamma, k,r} \in H'$ and 
	$\hat{x}_{\gamma, k,r} \in H$  
	such that $\| p_{\gamma, k,r} \|_{H'} < \gamma$ 
	and $\hat{x}_{ \gamma, k,r}$ is a strict $\mathbb{H}$-maximum of
	\begin{equation}
	\label{eq: strikt max on Br}
		\{ x \in H \colon \| x\|_H \leq r \} \ni x \to \phi_1(x) - \varphi_{k}(\pi^{H}_V (x))
		+ \langle p_{ \gamma, k,r}, x \rangle_{H', H} \in \R.
	\end{equation}
	Furthermore from (\ref{eq: behavior at infty without lin}) 
	it follows that
	for all $r \in (0, \infty)$
  and for all $k \in \N$
	it holds that
	\begin{equation}
		\label{eq: behavior at infty with lin}
		\begin{split}
				\lim_{r \to \infty} 
					\sup_{\gamma \in (0,1)} 
					\sup
					\big \{
						\phi_1(x) - \varphi_{k}(\pi^{H}_V (x)) 
						+ \langle p_{\gamma, k,r}, x \rangle_{H', H}
							\colon x \in H,~ \|x\|_H \geq r
					\big \}
			= 
				-\infty.
		\end{split}
	\end{equation}
	Thus for all $k \in \N$ there exists an $\hat{r}_k \in (0, \infty)$ such  that 
	for all $\gamma \in (0,1)$ and
	all $x \in H$ with $\| x\|_H \geq \hat{r}_k$ it holds that
	\begin{equation}
	\label{eq: phi-varphi + lin bounded at infty}
			\phi_1(x) - \varphi_{k}(\pi^{H}_V (x)) 
			+ \langle p_{\gamma, k,\hat{r}_k}, x \rangle_{H', H} 
		\leq
			\phi_1(0_H) - \varphi_{k}(0_V)-1.
	\end{equation}
	Moreover, equation (\ref{eq: convergence at infty whithout lin}) 
	implies that for all $k \in \N$ 
	there exists an $r_k \in (\hat{r}_k,\infty)$ such that
	for all $x \in H$ with $\|x\|_H \geq r_k$ it holds that
	\begin{equation}
	\label{eq: phi-varphi bounded at infty}
		 \phi_1(x) - \varphi_{k}(\pi^{H}_V (x)) \leq  \phi_1(x_{1,k}) - \varphi_{k}(x_{1,k})-1.
	\end{equation}
	Then from (\ref{eq: phi-varphi bounded at infty}) it follows that for all 
	$x \in V$ and for all $k \in \N$ it holds that
	\begin{equation}
	\label{eq: sup is less}
			\sup_{\substack{y \in V^\perp \\ \|y \|_H \geq r_k}}
			\Big [
				\phi_1(x_{1,k}+y) - \varphi_{k}(x_{1,k})
			\Big ]
		\leq
			\phi_1(x_{1,k}) - \varphi_{k}(x_{1,k})-1
		\leq
			\sup_{\substack{y \in V^\perp \\ \|y \|_H \leq r_k}}
			\Big [
				\phi_1(x_{1,k}+y) - \varphi_{k}(x_{1,k})
			\Big ]
	\end{equation}
	and this ensures that for all $k \in \N$ it holds that
	\begin{equation}
	\label{eq: sup is equal}
			\sup_{\substack{y \in V^\perp \\ \|y \|_H \leq r_k}}
				\Big [
					\phi_1(x_{1,k}+y) - \varphi_{k}(x_{1,k})
				\Big]
		=
				\sup_{\substack{y \in V^\perp}}
				\Big [
					\phi_1(x_{1,k}+y) - \varphi_{k}(x_{1,k})
				\Big].
	\end{equation}
	For every $\gamma \in (0,1)$ and every $k \in \N$ 
	denote by $p_{\gamma, k}$, $\hat{x}_{\gamma, k} \in H$ the
	elements satisfying that
	$p_{\gamma, k}=p_{\gamma, k,r_k}$ and that
	$\hat{x}_{\gamma, k}=\hat{x}_{\gamma, k,r_k}$.
	Then from (\ref{eq: strikt max on Br}) and (\ref{eq: phi-varphi + lin bounded at infty})
	it follows that for all $k \in \N$, for all $\gamma \in (0,1)$, and
	for all $x \in H$ with $\|x \|_H \geq r_k$ it holds that
	\begin{equation}
	\label{eq: strikt max outside Br}
		\begin{split}
			&\phi_1(x) - \varphi_{k}(\pi^{H}_V (x)) + \langle p_{\gamma, k}, x \rangle_{H', H} 
		\leq{}
			\phi_1(0_H) - \varphi_{k}(0_V)-1 \\
		\leq{} 
			&\phi_1(\hat{x}_{\gamma, k}) 
			- \varphi_{k}(\pi^{H}_V (\hat{x}_{\gamma, k})) 
			+ \langle p_{\gamma, k}, x \rangle_{H', H}-1.
		\end{split}
	\end{equation}
	Combining (\ref{eq: strikt max on Br}) and (\ref{eq: strikt max outside Br}) we get
	for all $k \in \N$ and all $\gamma \in (0,1)$
	that
	$\hat{x}_{\gamma, k}$ is also a strict $\mathbb{H}$-maximum of
	\begin{equation}
	\label{eq: max pertubation}
		\begin{split}
			H \ni x \to \phi_1(x) -\varphi_{k}(\pi^{H}_V (x)) 
			+ \langle p_{\gamma, k}, x \rangle_{H', H} \in \R
		\end{split}	
	\end{equation}
	and that 
	\begin{equation}
	\label{eq: boundedness of x(gamma,k)}
		\sup_{\gamma \in (0,1)} \| \hat{x}_{\gamma, k} \|_H \leq r_k.
	\end{equation}
	Moreover 
	equation (\ref{eq: def of ut}) 
	implies that for all $k \in \N$ and all $\gamma \in (0,1)$
	it holds that
	\begin{equation}
	\label{eq: u-phi inequality}
		\begin{split}
				&\tilde{u}_1(\pi^{H}_V (\hat{x}_{\gamma, k}))
				-\varphi_{k}(\pi^{H}_V (\hat{x}_{\gamma, k})) 
			=
				\sup_{\tilde{y} \in V^\perp} 
					\big[ \phi_1(\pi^{H}_V (\hat{x}_{\gamma, k})+\tilde{y}) \big]
				-\varphi_{k}(\pi^{H}_V (\hat{x}_{\gamma, k})) \\
			\geq{}&
				\phi_1(\hat{x}_{\gamma, k})
				-\varphi_{k}(\pi^{H}_V (\hat{x}_{\gamma, k})).
		\end{split}
	\end{equation}
	Hence for all $k \in \N$, $\gamma \in (0,1)$,
	and all $y \in V^\perp$ with $ \| y \|_H \leq r_k$ 
	(\ref{eq: max pertubation}) and (\ref{eq: boundedness of x(gamma,k)}) imply that
	\begin{equation}
	\label{eq: phi convergence for hatx}
		\begin{split}
				&\phi_1(\hat{x}_{\gamma, k})
				-\varphi_{k}(\pi^{H}_V (\hat{x}_{\gamma, k})) \\
			={}
				&\phi_1(\hat{x}_{\gamma, k})
				-\varphi_{k}(\pi^{H}_V (\hat{x}_{\gamma, k})) 
				+\langle p_{\gamma, k}, \hat{x}_{\gamma, k} \rangle_{H', H}
				-\langle p_{\gamma, k}, \hat{x}_{\gamma, k} \rangle_{H', H} \\
			\geq{}&
				\phi_1(x_{1, k}+y)
				-\varphi_{k}(x_{1,k})
				+\langle p_{\gamma, k}, x_{1, k}+y \rangle_{H', H}
				-\langle p_{\gamma, k}, \hat{x}_{\gamma, k} \rangle_{H', H} \\
			\geq{}&
				\phi_1(x_{1, k}+y)
				-\varphi_{k}(x_{1,k})
				-\gamma (\| x_{1, k}+y \|_H + \| \hat{x}_{\gamma, k} \|_H) \\
			\geq{}& 
				\phi_1(x_{1, k}+y)
				-\varphi_{k}(x_{1,k})
				-\gamma (\| x_{1, k} \|_H + 2 r_k).
		\end{split}
	\end{equation}
	Maximizing over $\{ y \in V^\perp : \| y \|_H \leq r_k \} $ together with
	(\ref{eq: def of ut}) and
	(\ref{eq: sup is equal})
	yields that for all $k \in \N$ and
	for all $\gamma \in (0,1)$ it holds that
	\begin{equation}
		\label{eq: phi-u convergence for hatx}
		\begin{split}
				&\phi_1(\hat{x}_{\gamma,k})
				- \varphi_k(\pi^{H}_V (\hat{x}_{\gamma, k}))
			\geq
				\sup_{\substack{y \in V^\perp \\ \|y \| \leq r_k}} 
				\Big [
					\phi_1(x_{1, k}+y)
					-\varphi_{k}(x_{1,k})- \gamma (\| x_{1, k} \|_H + 2 r_k)
				\Big ] \\
			={}&
				\sup_{y \in V^\perp} 
				\Big [
					 \phi_1(x_{1, k}+y)
					-\varphi_{k}(x_{1,k})- \gamma (\| x_{1, k} \|_H + 2 r_k)
				\Big ] \\
			={} 
				&\tilde{u}_1(x_{1,k})
				-\varphi_k(x_{1,k}) - \gamma (\| x_{1, k} \|_H + 2 r_k).
		\end{split}
	\end{equation}
	Thus it follows from
	(\ref{eq: u-phi inequality}) and
	(\ref{eq: phi-u convergence for hatx}) that
	for all $k \in \N$ it holds that
	\begin{equation}
		\begin{split}
				&\liminf_{\gamma \to 0} 
				\Big(
					\tilde{u}_1(\pi^{H}_V (\hat{x}_{\gamma, k}))-\varphi_k(\pi^{H}_V (\hat{x}_{\gamma, k}))
				\Big)
			\geq
				\liminf_{\gamma \to 0}
				\Big(
					\phi_1(\hat{x}_{\gamma,k})
						- \varphi_k(\pi^{H}_V (\hat{x}_{\gamma, k}))
				\Big) \\
			\geq{}
				&\tilde{u}_1(x_{1,k}) - \varphi_k(x_{1,k})
		\end{split}
	\end{equation}
	and since $x_{1,k}$ is a strict $\mathbb{V}$-maximum of 
	$V \ni x \to \tilde{u}_1(x) - \varphi_{k}(x) \in \R$ we have
	for all $k \in \N$ that
	\begin{align}
	\label{eq: gamma convergence for ut}
			&\lim_{\gamma \to 0} 
			\Big(
				\tilde{u}_1(\pi^{H}_V (\hat{x}_{\gamma, k}))
				-\varphi_k(\pi^{H}_V (\hat{x}_{\gamma, k}))
			\Big)
		=
			\tilde{u}_1(x_{1,k}) - \varphi_k(x_{1,k}) \\
	\label{eq: gamma convergence for x}
			&\lim_{\gamma \to 0} 
				\pi^{H}_V (\hat{x}_{\gamma, k}) = x_{1,k}.
	\end{align}
	From (\ref{eq: u-phi inequality}) and
	(\ref{eq: phi-u convergence for hatx})
	it follows then that
	for all $\gamma \in (0,1)$ and
	all $k \in \N$ it holds that
	\begin{equation}
	\label{eq: gamma inequality for u}
		\begin{split}
				&\tilde{u}_1(\pi^{H}_V (\hat{x}_{\gamma, k}))
				-\varphi_{k}(\pi^{H}_V (\hat{x}_{\gamma, k})) 
			\geq
				\phi_1(\hat{x}_{\gamma, k})
				-\varphi_{k}(\pi^{H}_V (\hat{x}_{\gamma, k})) \\
			\geq{} 
				&\tilde{u}_1(x_{1,k})
				-\varphi_k(x_{1,k}) - \gamma (\| x_{1, k} \|_H + 2 r_k).
		\end{split}
	\end{equation}
	Hence we get with 
	(\ref{eq: gamma convergence for ut}), 
	(\ref{eq: gamma convergence for x}), and with
	\eqref{eq: gamma inequality for u} that
	for all $k \in \N$ it holds that
	\begin{equation}
		\begin{split}
				\tilde{u}_1(x_{1,k})
			=
				\lim_{\gamma \to 0} (
					\tilde{u}_1(\pi^{H}_V(\hat{x}_{\gamma, k}))
				)
			\geq
				\limsup_{\gamma \to 0} (
					\phi_1(\hat{x}_{\gamma, k})
				)
			\geq
				\liminf_{\gamma \to 0} (
					\phi_1(\hat{x}_{\gamma, k})
				)
			\geq
				\tilde{u}_1(x_{1,k})
		\end{split}
	\end{equation}
	and this shows that
	\begin{equation}
	\label{eq: gamma convergence for u}
			\lim_{\gamma \to 0} \phi_1(\hat{x}_{\gamma, k}) 				                 
		=
			\tilde{u}_1(x_{1,k}).
	\end{equation}
	Now with   
	(\ref{eq: phi and phit}),
	(\ref{eq: max of phi}),
	(\ref{eq: gamma convergence for x}),
	and (\ref{eq: gamma convergence for u})  we derive 
	for all $k \in \N$  that 
	\begin{equation}
	\label{eq: xhat max of Phit}
		\begin{split}
				&\Phi(\overline{z}_1, \overline{z}_2)
			\geq{}
				\limsup_{\gamma \to 0} 
					\Phi(\hat{x}_{\gamma, k} , \overline{z}_2)  \\
			\geq{}& 
				\limsup_{\gamma \to 0} 
				\bigg ( 
					\Phi(\hat{x}_{\gamma, k} ,  \overline{z}_2) 
					+ 2 \left \langle 
						\left ( \hat{z} - \frac{(\hat{x}_{\gamma, k}, \overline{z}_2)}{2} \right),
						B \left ( \hat{z} - \frac{(\hat{x}_{\gamma, k}, \overline{z}_2)}{2} \right) 
					\right \rangle_{H^2, (H^2)'} 
				\bigg ) \\
			={}&
				\limsup_{\gamma \to 0} 
				\Big (
					\phi_1(\hat{x}_{\gamma, k}) + \phi_2( \overline{z}_2) 
					- \frac 12 
						\big \langle 
							(\pi^{H}_V (\hat{x}_{\gamma, k}) ,\overline{x}_2), 
							A (\pi^{H}_V (\hat{x}_{\gamma, k}) ,\overline{x}_2) 
						\big \rangle_{H^2, (H^2)'}
				\Big )	\\
			={}& 
				\tilde{u}_1(x_{1,k}) + \tilde{u}_2(\overline{x}_2)
				- \frac 12 
					\big \langle 
						(x_{1,k},\overline{x}_2), A (x_{1,k},\overline{x}_2) 
					\big \rangle_{H^2, (H^2)'}. 
		\end{split}
	\end{equation}
	Then letting k go to infinity (\ref{eq: xhat max of Phit}) together with
	(\ref{eq: k convergence for x}) and (\ref{eq: k convergence for ut}) implies that
	\begin{equation}
	\label{eq: Phi convergence for xhat}
		\begin{split}
				&\Phi(\overline{z}_1, \overline{z}_2) 
			\geq{} 
				\lim_{k \to \infty} \limsup_{\gamma \to 0} 
					\Phi(\hat{x}_{\gamma, k}, \overline{z}_2) \\
			\geq{}&
				\lim_{k \to \infty} 
				\Big (
					\tilde{u}_1(x_{1,k}) + \tilde{u}_2(\overline{x}_2)
					- \frac 12 \big 
						\langle 
							(x_{1,k},\overline{x}_2), A (x_{1,k},\overline{x}_2) 
						\big \rangle_{H^2, (H^2)'} 
				\Big ) \\
			={}&
				\tilde{u}_1(\overline{x}_1) + \tilde{u}_2(\overline{x}_2)
				- \frac 12 
					\big \langle 
						(\overline{x}_1,\overline{x}_2), A (\overline{x}_1,\overline{x}_2) 
					\big \rangle_{H^2, (H^2)'}
			= 
				\tilde{\Phi}(\overline{x}_1, \overline{x}_2) 
			= 
				\Phi(\overline{z}_1, \overline{z}_2).
		\end{split}
	\end{equation}
	From (\ref{eq: Phi convergence for xhat}) it follows then that
	\begin{equation}
	\label{eq: gamma, k convergence for Phi}
			\lim_{k \to \infty} \limsup_{\gamma \to 0}
					\Phi(\hat{x}_{\gamma, k} , \overline{z}_2)  
		=
			\Phi(\overline{z}_1, \overline{z}_2).
	\end{equation}
	In addition we get from 
	\eqref{eq: X k finite dimensional},
	(\ref{eq: def of varphi}), and
	(\ref{eq: gamma convergence for x})
	for all $k \in \N$ that
	\begin{align}
	\nonumber
				&\lim_{\gamma \to 0} \varphi_k(\pi^{H}_V(\hat{x}_{\gamma,k}))
			= 
				 \varphi_k(x_{1,k})
			=
				 \tilde{u}_1(x_{1,k})  \\
				&\lim_{\gamma \to 0} 
					(D_{\mathbb{H}} (\varphi_k \pi^{H}_V))(\hat{x}_{\gamma,k})
			= 
				(D_{\mathbb{H}}(\varphi_k \pi^{H}_V))(x_{1,k})
			=
				\xi_{1,k} \pi^{H}_V \\ \nonumber
				&\lim_{\gamma \to 0} (D^2_{\mathbb{H}} (\varphi_k \pi^{H}_V))(\hat{x}_{\gamma,k})
			= 
				(D^2_{\mathbb{H}} (\varphi_k \pi^{H}_V)) (x_{1,k})
			=
				\pi^{H'}_{V'} X_{1,k} \pi^{H}_V
			= 
				X_{1,k}.
	\end{align}
	Thus for all $k \in \N$ there exist a $C_k \in (0,1)$ such that
	for all $\gamma \in (0,C_k)$ it holds that
	\begin{equation}
	\label{eq: inequality ensuring convergence}
		\begin{split}
				&| \varphi_k(\pi^{H}_V(\hat{x}_{\gamma,k}))- \tilde{u}_1(x_{1,k})| 
			\leq \frac 1k,\\
				&\| 
					(D_{\mathbb{H}} (\varphi_k \pi^{H}_V))(\hat{x}_{\gamma,k}) -
					\xi_{1,k} \pi^{H}_V 
				\|_{H'} 
			\leq \frac 1k,
			\qquad \textrm{ and that}\\
				&\| 
					(D_{\mathbb{H}}^2 (\varphi_k \pi^{H}_V))(\hat{x}_{\gamma,k}) - X_{1,k}
				\|_{L(\mathbb{H}, \mathbb{H}')} 
			\leq \frac 1k.
		\end{split}	
	\end{equation}
	From (\ref{eq: gamma, k convergence for Phi}) it follows now
	that for all $k \in \N$ there exist a $\gamma_k \in (0,C_k)$
	such that
	\begin{align}
	\label{eq: k convergence for Phi}
			&\lim_{k \to \infty}
					\Phi(\hat{x}_{\gamma_k, k} , \overline{z}_2)  
		=
			\Phi(\overline{z}_1, \overline{z}_2).
	\end{align}
	With the assumption that $(\overline{z}_1,\overline{z}_2)$ 
	is a strict $\mathbb{H} \times \mathbb{H}$-maximum of $\Phi$ we get
	\begin{equation}
	\label{eq: gamma, k convergence for xhat}
			\lim_{k \to \infty}
				\hat{x}_{\gamma_k, k} = \overline{z}_1.
	\end{equation}
	Now (\ref{eq: k convergence for Phi})
	and (\ref{eq: gamma, k convergence for xhat}) imply that
	\begin{equation}
		\begin{split}
				&u_1(\overline{z}_1) + u_2(\overline{z}_2)
				-\tfrac 12 
				\big \langle 
					(\overline{z}_1,\overline{z}_2), A (\overline{z}_1,\overline{z}_2)
				\big \rangle_{H^2, (H^2)'} 
			=
				\Phi(\overline{z}_1,\overline{z}_2)
			= 
				\lim_{k \to \infty}
					\Phi(\hat{x}_{\gamma_k, k} , \overline{z}_2) \\
			={}&
				\lim_{k \to \infty}
				\Big (
					u_1(\hat{x}_{\gamma_k, k}) + u_2( \overline{z}_2)
					-\tfrac 12 
					\big \langle 
						(\hat{x}_{\gamma_k, k}, \overline{z}_2), 
						A (\hat{x}_{\gamma_k, k}, \overline{z}_2)
					\big \rangle_{H^2, (H^2)'} 
				\Big ) \\
			={}&
				\lim_{k \to \infty}
				\Big (
					u_1(\hat{x}_{\gamma_k, k}) + u_2( \overline{z}_2)
					-\tfrac 12 
					\big \langle 
						(\overline{z}_1,\overline{z}_2), A (\overline{z}_1,\overline{z}_2)
					\big \rangle_{H^2, (H^2)'}
				\Big )
		\end{split}
	\end{equation}
	and thus
	\begin{equation}
	\label{eq:  gamma, k convergence for u}
		\lim_{k \to \infty}
		u_1 (\hat{x}_{\gamma_k, k}) = u(\overline{z}_1).
	\end{equation}
	Moreover,
	(\ref{eq: k convergence for ut}),
	(\ref{eq: k convergence for X}),
	(\ref{eq: k convergence for xi}),
	(\ref{eq: inequality ensuring convergence}), and 
	the fact that
	that $\forall k \in \N \colon \gamma_k < C_k$
	ensure that
	\begin{equation}
	\label{eq: final convergence}
		\begin{split}
				&\lim_{k \to \infty} 
					\varphi_k(\pi^{H}_V(\hat{x}_{\gamma_k,k})) 
			= 
				\lim_{k \to \infty} \left(
					\varphi_k(\pi^{H}_V(\hat{x}_{\gamma_k,k})) 
					- \tilde{u}_1(x_{1,k})
					+ \tilde{u}_1(x_{1,k})
				\right)
			=
				\tilde{u}_1(\overline{x}_1),  \\
				&\lim_{k \to \infty} 
					\left (
						(D_{\mathbb{H}} (\varphi_k \pi^{H}_V)) (\hat{x}_{\gamma_k,k})
					\right)
			=
				\lim_{k \to \infty} 
					\left(
						(D_{\mathbb{H}} (\varphi_k \pi^{H}_V )) (\hat{x}_{\gamma_k,k}) 
						- \xi_{1,k} \pi^{H}_V + \xi_{1,k} \pi^{H}_V
					\right) \\
				& \qquad \qquad
			=
				\pi^{H'}_{V'} \pi^{(H \times H)'}_1 A(\overline{x}_1, \overline{x}_2), 
				\\&
				\lim_{k \to \infty} 
					(D^2_{\mathbb{H}} (\varphi_k \pi^{H}_V) (\hat{x}_{\gamma_k,k})
			=
				\lim_{k \to \infty} \left(
					(D^2_{\mathbb{H}} (\varphi_k \pi^{H}_V)) (\hat{x}_{\gamma_k,k})
					-X_{1,k} + X_{1,k}
				\right)
			=
				X_1.
		\end{split}
	\end{equation}
	On the other hand we obtain from (\ref{eq: def of phit}) and
	from (\ref{eq: max pertubation}) that
	for all $k \in \N$ it holds that
	\begin{align}
	\nonumber
			\big( 
				&(D_{\mathbb{H}} (\varphi_{k} \pi^{H}_V)) (\hat{x}_{\gamma_k, k})
				- p_{\gamma_k, k},
				(D_{\mathbb{H}}^2 (\varphi_{k} \pi^{H}_V))(\hat{x}_{\gamma_k, k})
			\big) 
			\in (J^2_{\mathbb{H}, + } \phi_1)(\hat{x}_{\gamma_k, k}) \\ \nonumber
		\Leftrightarrow
			\big( 
				&D_{\mathbb{H}} (\varphi_{k} \pi^{H}_V) (\hat{x}_{\gamma_k, k})
				+ \pi^{(H^2)'}_1 A E_1(\hat{x}_{\gamma_k, k}) 
				- \pi^{H'}_{V'} \pi^{(H^2)'}_1 A E_1(\pi^{H}_V (\hat{x}_{\gamma_k, k})) \\ \nonumber
				&- 2 \pi^{(H^2)'}_{1} B((\hat{x}_{\gamma_k, k},0)-\hat{z}) 
				- \delta I_{\mathbb{H}} \pi^{H}_V (\hat{x}_{\gamma_k, k}) 
				- \lambda I_{\mathbb{H}} \pi^{H}_{V^\perp} (\hat{x}_{\gamma_k, k} )
				- p_{\gamma_k, k}, 
				(D_{\mathbb{H}}^2 (\varphi_{k} \pi^{H}_V))(\hat{x}_{\gamma_k, k}) \\ 
	\label{eq: semijet}
				&+ \pi^{(H^2)'}_1 A E_1 -\pi^{H'}_{V'} \pi^{(H^2)'}_1 A E_1 \pi^{H}_V
				- 2 \pi^{(H^2)'}_{1} B E_{1}
				- \delta I_{\mathbb{H}} \pi^{H}_{V}
				- \lambda I_{\mathbb{H}} \pi^{H}_{V^\perp}
			\big) 
			\in (J^2_{\mathbb{H}, + } u_1)(\hat{x}_{\gamma_k, k}).
	\end{align}
	Letting k go to infinity we deduce from
	\eqref{eq: gamma, k convergence for xhat}
	and \eqref{eq: final convergence} that	
	\begin{equation}
		\begin{split}
			&\lim_{k \to \infty} \big(
					D_{\mathbb{H}} (\varphi_{k} \pi^{H}_V) (\hat{x}_{\gamma_k, k})
					+ \pi^{(H^2)'}_1 A (\hat{x}_{\gamma_k, k},0) 
					- \pi^{H'}_{V'} \pi^{(H^2)'}_1 A (\pi^{H}_V (\hat{x}_{\gamma_k, k},0))  \\
					&- 2 \pi^{(H^2)'}_{1} B((\hat{x}_{\gamma_k, k},0)-\hat{z}) 
					- \delta I_{\mathbb{H}} \pi^{H}_V (\hat{x}_{\gamma_k, k}) 
					- \lambda I_{\mathbb{H}} \pi^{H}_{V^\perp} (\hat{x}_{\gamma_k, k}) 
					- p_{\gamma_k, k}
				\big)\\ 
			={}	
				&\pi^{H'}_{V'} \pi^{(H^2)'}_1 A (\overline{x}_1,\overline{x}_2)
				+ \pi^{(H^2)'}_1 A(\overline{z}_1,0) 
				- \pi^{H'}_{V'} \pi^{(H^2)'}_1 A (\overline{x}_1,0)
				- \pi^{(H^2)'}_{1} B(\overline{z}_1,-\overline{z}_2) \\
				&- \delta I_{\mathbb{H}} \overline{x}_1 
				- \lambda I_{\mathbb{H}} \overline{y}_1 \\
			={}	
				&\pi^{H'}_{V'} \pi^{(H^2)'}_1 A (0,\overline{x}_2)
				+ \pi^{(H^2)'}_1 A(\overline{z}_1,0) 
				+ \delta I_{\mathbb{H}} \overline{x}_1
				+ \lambda I_{\mathbb{H}} \overline{y}_1 \\
				&+ (
						\pi^{(H^2)'}_1 A E_2 
						-\pi^{H'}_{V'} \pi^{(H^2)'}_1 A E_2 \pi^{H}_V
					) \overline{z}_2 
				- \delta I_{\mathbb{H}} \overline{x}_1 
				- \lambda I_{\mathbb{H}} \overline{y}_1 \\
			={}
				&\pi^{(H^2)'}_1 A (\overline{z}_1,\overline{z}_2)
		\end{split}
	\end{equation}
	and that
	\begin{align}
	\nonumber
			\lim_{k \to \infty}
				\big( 
					&D_{\mathbb{H}}^2 (\varphi_{k} \pi^{H}_V)\hat{x}_{\gamma_k, k}
					+ \pi^{(H^2)'}_1 A E_1 -\pi^{H'}_{V'} \pi^{(H^2)'}_1 A E_1 \pi^{H}_V
					- 2 \pi^{(H^2)'}_{1} B E_{1}
					- \delta I_{\mathbb{H}} \pi^{H}_{V}
					- \lambda I_{\mathbb{H}} \pi^{H}_{V^\perp}
				\big)   \\ 
		\nonumber
			=	
				&X_1 
				+ \pi^{(H^2)'}_1 A E^{(H^2)'}_1 
				-\pi^{H'}_{V'} \pi^{(H^2)'}_1 A E_1 \pi^{H}_V
				+ 2 \delta I_{\mathbb{H}} \pi^{H}_{V}
				+ 2 \lambda I_{\mathbb{H}} \pi^{H}_{V^\perp}
				- \delta I_{\mathbb{H}} \pi^{H}_{V}
				- \lambda I_{\mathbb{H}} \pi^{H}_{V^\perp}\\
			=	    
				&X_1 
				+ \pi^{(H^2)'}_1 A E_1 -\pi^{H'}_{V'} \pi^{(H^2)'}_1 A E_1 \pi^{H}_V
				+ \delta I_{\mathbb{H}} \pi^{H}_{V}
				+ \lambda I_{\mathbb{H}} \pi^{H}_{V^\perp}\\
		\nonumber
			=
				&X_1 
				+ \pi^{(H^2)'}_1 A E_1 -\pi^{H'}_{V'} \pi^{(H^2)'}_1 A E_1 \pi^{H}_V
					+ \delta I_{\mathbb{H}} \pi^{H}_{V}
					+ \lambda I_{\mathbb{H}} \pi^{H}_{V^\perp}.
	\end{align}
	This together with 
	(\ref{eq: gamma, k convergence for xhat}),
	(\ref{eq:  gamma, k convergence for u}), and
	(\ref{eq: semijet}) shows that
	\begin{equation}
			\big ( 
				\pi^{(H^2)'}_1 A (\overline{z}_1,\overline{z}_2), 
				X_1 
				+ \pi^{(H^2)'}_1 A E_1 -\pi^{H'}_{V'} \pi^{(H^2)'}_1 A E_1 \pi^{H}_V
				+ \delta I_{\mathbb{H}} \pi^{H}_{V}
				+ \lambda I_{\mathbb{H}} \pi^{H}_{V^\perp}
			\big )
		\in (\hat{J}^2_{\mathbb{H}, +} u_1)(\overline{z}_1)
	\end{equation}
	which completes the proof of Lemma \ref{l: suitable semijets matrix}.

\end{proof}

For completeness we also state the result for general 
$\C^2$ functions. It follows from 
Lemma \ref{l: suitable semijets matrix}
and a standard translation technique
(see, e.g., First and Second Reduction on page 57 in 
 Crandall, Ishii \& Lions \cite{CrandallIshiiLions1992})
\begin{corollary}[Construction of suitable semijets for strict maxima wrt.~to a $\C^2$ function]
\label{cor: suitable semijets function}
	Let $\mathbb{H} = (H, \langle \cdot, \cdot \rangle_H, \| \cdot \|_H)$,
	be a real Hilbert space,
	let $\mathbb{V} = (V, \langle \cdot, \cdot \rangle_H, \| \cdot \|_H)$ 
	be a finite-dimensional linear subspace of $\mathbb{H}$,
	let
	$\mathbb{H}' = (H', \| \cdot \|_{H'})$,
	$\mathbb{V}' = (V', \| \cdot \|_{V'})$,
	$(\mathbb{H}^2)' = ((H^2)', \| \cdot \|_{(H^2)'})$,
	$(\mathbb{V}^2)' = ((V^2)', \| \cdot \|_{(V^2)'})$
	be the dual spaces,
	let $O,U \subseteq H$ be open sets,
	let $u_1 \colon O \to  [- \infty, \infty)$, 
	$u_2 \colon U \to  [- \infty, \infty)$ be upper semicontinuous
	with respect to the $\| \cdot \|_{\mathbb{H}}$-norm,
	let $\Psi \in \C_{\mathbb{H}^2}^2(O \times U, \R)$,
	let $(\overline{z}_1$, $\overline{z}_2) \in O \times U$ 
	be a strict $\mathbb{H}^2$-maximum of
	$
		O \times U \ni (z_1,z_2) \to u_1(z_1) + u_2(z_2) -\Psi(z_1, z_2) \in  [- \infty, \infty)
	$, 
	denote by $E_1$, $E_2 \colon H \to H^2$ the canonical embeddings
	satisfying for all $x \in H$ that
	$E_1(x)=(x,0)$ and that $E_2(x)=(0,x)$,
	and let $\eps, \delta, \lambda \in (0, \infty) $
	satisfy that
	\begin{equation}
		\begin{split}
			\lambda >{}
				&\tfrac {1}{\delta}
					\Big(
						\|
							\pi^{H'}_{V'} \pi^{(H^2)'}_1 
								((D^2_{\mathbb{H}^2} \Psi) (\overline{z}_1, \overline{z}_2))  
									E_2 \pi^{H}_{V^\perp} 
						\|_{L(\mathbb{H},\mathbb{H}' )}^2 \\
						& \qquad
						\vee 
						\| 
							\pi^{H'}_{(V^\perp)'} \pi^{(H^2)'}_1 
								((D^2_{\mathbb{H}^2} \Psi) (\overline{z}_1, \overline{z}_2)) 
									E_2 \pi^{H}_{V} 
						\|_{L(\mathbb{H},\mathbb{H}' )}^2
					\Big) \\
				&+\|
					\pi^{H'}_{(V^\perp)'} \pi^{(H^2)'}_1 
						((D^2_{\mathbb{H}^2} \Psi) (\overline{z}_1, \overline{z}_2)) 
							E_2 \pi^{H}_{V^\perp} 
				\|_{L(\mathbb{H},\mathbb{H}' )}.
		\end{split}
	\end{equation}
	Then
	there exist operators 
	$X_{1}$, $X_{2} \in \mathbb{S}_{\mathbb{H}, \mathbb{H}'}$ 
	satisfying that
	$X_{1} = \pi^{H'}_{V'} X_{1} \pi^{H}_V $, 
	$X_{2} = \pi^{H '}_{V'} X_{2} \pi^{H}_V$, 
	\begin{equation}
		\begin{split}
				&-\left(
					\tfrac {1}{\eps} 
					+ \| 
							((D^2_{\mathbb{H}^2} \Psi) (\overline{z}_1, \overline{z}_2)) 
						\|_{L(\mathbb{H}^2, (\mathbb{H}^2)')} 
				\right) 
				I_{\mathbb{H}^2}
			\leq
				\left(
					\begin{array}{cc}
						X_{1}
						&
						0
						\\
						0
						&
						X_{2}
					\end{array}
				\right) \\
			\leq{}&
				\left( 
					\pi^{(H^2)'}_{(V^2)'} 
						((D^2_{\mathbb{H}^2} \Psi) (\overline{z}_1, \overline{z}_2))
							\pi^{H^2}_{V^2} 
					+ \eps (\pi^{(H^2)'}_{(V^2)'} 
						((D^2_{\mathbb{H}^2} \Psi) (\overline{z}_1, \overline{z}_2))
							\pi^{H^2}_{V^2})^2 
				\right),
		\end{split}
	\end{equation}
	and satisfying for all $i \in \{1,2\}$ that 
	\begin{align}
		\nonumber
			&\Big( 
				\pi^{(H^2)'}_i ((D_{\mathbb{H}^2} \Psi) 
					(\overline{z}_1, \overline{z}_2)), 
				X_{i} 
				+ \pi^{(H^2)'}_i 
					((D^2_{\mathbb{H}^2} \Psi) (\overline{z}_1, \overline{z}_2)) E_i \\
			& \qquad \qquad \qquad
				- \pi^{H'}_{V'} \pi^{(H^2)'}_i 
					((D^2_{\mathbb{H}^2} \Psi) (\overline{z}_1, \overline{z}_2))
						E_i \pi^{H}_{V} 
				+ \delta I_{\mathbb{H}} \pi^{H}_V + \lambda I_{\mathbb{H}} \pi^{H}_{V^\perp}
			\Big)
		\in (\hat{J}^2_{ \mathbb{H}, + } u_i) (\overline{z}_i). 
	\end{align}
\end{corollary}
\begin{proof}
	Let
	$
		\tilde{\Psi} \in \C^2_{\mathbb{H}^2}(H^2, \R)
	$
	be the function satisfying for all $z \in H^2$ that
	\begin{equation}
			\tilde{\Psi}(z)
		=
			\begin{cases}
				\Psi(z+ (\overline{z}_1, \overline{z}_2))
				- \Psi(\overline{z}_1, \overline{z}_2)
				-	\langle
						(D_{\mathbb{H}^2} \Psi)(\overline{z}_1, \overline{z}_2),
						z
					\rangle_{(H^2)', H^2}, 
						 &\textrm{ if } z+ (\overline{z}_1, \overline{z}_2) \in O \times U \\
				-\infty, &\textrm{ if } z+ (\overline{z}_1, \overline{z}_2) \notin O \times U,
			\end{cases}
	\end{equation}
	let
	$\tilde{u}_1$, $\tilde{u}_2 \colon H \to[- \infty, \infty)$ 
	be the function satisfying for all $z \in H$ that
	\begin{equation}
			\tilde{u}_1(z) 
		= 
			\begin{cases}
				u_1(z+ \overline{z}_1)
				- u_1(\overline{z}_1)
				-	\langle
						(D_{\mathbb{H}^2} \Psi)(\overline{z}_1, \overline{z}_2),
						(z,0)
					\rangle_{(H^2)', H^2},
					\quad &\textrm{ if } z+ \overline{z}_1 \in O \\
				-\infty, &\textrm{ if } z+ \overline{z}_1 \notin O,
			\end{cases}
	\end{equation}
	and that
	\begin{equation}
			\tilde{u}_2(z) 
		= 
			\begin{cases}
				u_2(z+ \overline{z}_2)
				- u_2(\overline{z}_2)
				-	\langle
						(D_{\mathbb{H}^2} \Psi)(\overline{z}_1, \overline{z}_2),
						(0,z)
					\rangle_{(H^2)', H^2},
					\quad &\textrm{ if } z+ \overline{z}_2 \in U \\
				-\infty, &\textrm{ if } z+ \overline{z}_2 \notin U,
			\end{cases}
	\end{equation}
	denote by $\overline{\eta} \in (0,1]$ a real number satisfying that
	\begin{align}
	\nonumber
				\overline{\eta}
			\leq{}& 
				\tfrac{\delta}{1+\delta} \Big(
					\lambda 
					-\tfrac {1}{\delta}
						\Big(
							\|
								\pi^{H'}_{V'} \pi^{(H^2)'}_1 
									((D^2_{\mathbb{H}^2} \Psi)
										(\overline{z}_1, \overline{z}_2))  
											E_2 \pi^{H}_{V^\perp} 
							\|_{L(\mathbb{H},\mathbb{H}' )}^2 \\
					& \qquad
							\vee 
							\| 
								\pi^{H'}_{(V^\perp)'} \pi^{(H^2)'}_1 
									((D^2_{\mathbb{H}^2} \Psi) 
										(\overline{z}_1, \overline{z}_2)) 
											E_2 \pi^{H}_{V} 
							\|_{L(\mathbb{H},\mathbb{H}' )}^2
						\Big) \\ \nonumber
					&-\|
						\pi^{H'}_{(V^\perp)'} \pi^{(H^2)'}_1 
							((D^2_{\mathbb{H}^2} \Psi) (\overline{z}_1, \overline{z}_2)) 
								E_2 \pi^{H}_{V^\perp} 
					\|_{L(\mathbb{H},\mathbb{H}' )}
				\Big),
	\end{align}
	and denote for every $\eta \in [0,\infty)$ by 
	$A_\eta \in \mathbb{S}_{\mathbb{H} \times \mathbb{H},(\mathbb{H} \times \mathbb{H})'}$ 
	the operator
	satisfying that
	\begin{equation}
	\label{eq: def A eta}
			A_\eta
		=
			(D^2_{\mathbb{H}^2} \Psi) (\overline{z}_1,\overline{z}_2) 
				+ \eta I_{\mathbb{H}^2}.
	\end{equation}
	Then it holds that 
	$\tilde{u}_1(0) = \tilde{u}_2(0) = \tilde{\Psi}(0,0) = 0$, 
	$(D_{\mathbb{H}^2} \tilde{\Psi}) (0,0) = (0,0)$,
	$
			(D^2_{\mathbb{H}^2} \tilde{\Psi}) (0,0)
		= 
			(D^2_{\mathbb{H}^2} \Psi) (\overline{z}_1,\overline{z}_2)
	$,
	and that $(0,0)$ is a strict $\mathbb{H}^2$-maximum of
	$
		H^2 \ni (z_1,z_2) \to 
			\tilde{u}_1(z_1) + \tilde{u}_2(z_2) -\tilde{\Psi}(z_1, z_2) \in  [- \infty, \infty)
	$.
	Moreover, the fact that
	$\tilde{\Psi} \in \C^2_{\mathbb{H}^2}(O \times U, \R)$
	shows that for all $\eta \in [0,\infty)$ it holds that 
	\begin{equation}
		\begin{split}
				&\lim_{H^2 \ni z \to 0} \left(
					\frac
						{
							\tilde{\Psi}(z) - \tilde{\Psi}(0) 
							- \langle
									(D_{\mathbb{H}^2} \tilde{\Psi}) (0),
									z
								\rangle_{(H^2)', H^2}
							- \langle
									(D^2_{\mathbb{H}^2} \tilde{\Psi}) (0) z,
									z
								\rangle_{(H^2)', H^2}
						}
						{\|z\|^2_{H^2}}
				\right) \\
			={}&
				\lim_{H^2 \ni z \to 0} \left(
					\frac
						{
							\tilde{\Psi}(z)
							- \langle
									(A_\eta - \eta I_{\mathbb{H}^2}) z,
									z
							\rangle_{(H^2)', H^2}
						}
						{\|z\|^2_{H^2}}
				\right)
			= 0.
		\end{split}
	\end{equation}
	Thus for every $\eta \in (0,\infty)$ there exists an $r_{\eta} \in (0, \infty)$
	such that for all $z \in (H^2)\backslash\{ 0 \}$
	with $\| z\|^2_{H^2} \leq 2 r_{\eta}$
	it holds that
	$
			\tilde{\Psi} (z) 
		<
			\langle
				A_\eta z,
				z
			\rangle_{(H^2)', H^2}.
	$
	Now denote for every $\eta \in (0,\infty)$ by 
	$\tilde{u}_{1,\eta}$, $\tilde{u}_{2,\eta} \colon H \to [-\infty, \infty)$
	the functions satisfying for all $i \in \{1, 2\}$
	and all $z \in H$ that
	\begin{equation}
		\begin{split}
				\tilde{u}_{i, \eta}(z) 
			= 
				\begin{cases}
					\tilde{u}_{i}(z) & \textrm{ if } \|z\|_H^2 < r_\eta \\
					-\infty & \textrm{ if } \|z\|_H^2 \geq r_\eta.
				\end{cases} 
		\end{split}
	\end{equation}
	Then we have for all $\eta \in (0,\infty)$ 
	that $\tilde{u}_{1, \eta}$, $\tilde{u}_{2, \eta}$
	are upper semicontinuous
	with respect to the $\| \cdot \|_{\mathbb{H}}$-norm
	and that $(0,0)$ is a strict $\mathbb{H}^2$-maximum of
	\begin{equation}
		H^2 \ni (z_1, z_2) \to 
			\tilde{u}_{1, \eta}(z_1) + \tilde{u}_{2, \eta}(z_2) 
			- \langle
					A_\eta (z_1, z_2),
					(z_1, z_2)
				\rangle_{(H^2)', H^2}.
	\end{equation} 
	Moreover, for all $\eta \in (0,\overline{\eta}]$ it holds that
	\begin{align}
	\nonumber
				\lambda 
			\geq{}&
				\tfrac {1}{\delta}
				\Big(
					\|
						\pi^{H'}_{V'} \pi^{(H^2)'}_1 
							A
								E_2 \pi^{H}_{V^\perp} 
					\|_{L(\mathbb{H},\mathbb{H}' )}^2
					\vee 
					\| 
						\pi^{H'}_{(V^\perp)'} \pi^{(H^2)'}_1 
							A
								E_2 \pi^{H}_{V} 
					\|_{L(\mathbb{H},\mathbb{H}' )}^2
				\Big) \\ \nonumber
				&+\|
					\pi^{H'}_{(V^\perp)'} \pi^{(H^2)'}_1 
						A
							E_2 \pi^{H}_{V^\perp} 
				\|_{L(\mathbb{H},\mathbb{H}' )}
				+\eta(\tfrac{1+\delta}{\delta})\\
		\begin{split}
			\geq{}&
				\tfrac {1}{\delta}
				\Big(
					\|
						\pi^{H'}_{V'} \pi^{(H^2)'}_1 
							A
								E_2 \pi^{H}_{V^\perp} 
					\|_{L(\mathbb{H},\mathbb{H}' )}^2
					\vee 
					\| 
						\pi^{H'}_{(V^\perp)'} \pi^{(H^2)'}_1 
							A
								E_2 \pi^{H}_{V} 
					\|_{L(\mathbb{H},\mathbb{H}' )}^2
				\Big) \\
				&+\tfrac{\eta^2}{\delta}
				+\|
					\pi^{H'}_{(V^\perp)'} \pi^{(H^2)'}_1 
						A
							E_2 \pi^{H}_{V^\perp} 
				\|_{L(\mathbb{H},\mathbb{H}' )}
				+\eta
		\end {split} \\ \nonumber
			\geq{}&
				\tfrac {1}{\delta}
				\Big(
					\|
						\pi^{H'}_{V'} \pi^{(H^2)'}_1 
							A_\eta  
								E_2 \pi^{H}_{V^\perp} 
					\|_{L(\mathbb{H},\mathbb{H}' )}^2 
					\vee 
					\| 
						\pi^{H'}_{(V^\perp)'} \pi^{(H^2)'}_1 
							A_\eta 
								E_2 \pi^{H}_{V} 
					\|_{L(\mathbb{H},\mathbb{H}' )}^2
				\Big) \\ \nonumber
				&-\|
					\pi^{H'}_{(V^\perp)'} \pi^{(H^2)'}_1 
						A_\eta
							E_2 \pi^{H}_{V^\perp} 
				\|_{L(\mathbb{H},\mathbb{H}' )}.
	\end{align}
	Thus for every $\eta \in (0,\overline{\eta}]$
	Lemma \ref{l: suitable semijets matrix}
		(with
		$
				A 
			\leftarrow 
				A_\eta
		$,
		$u_1 \leftarrow \tilde{u}_{1, \eta}$, 
		$u_2 \leftarrow \tilde{u}_{2, \eta}$,
		$(\overline{z}_1, \overline{z}_2) \leftarrow (0, 0)$)
	yields that for every
	$\eta \in (0,\overline{\eta}]$ 
	there exist operators 
	$X_{1, \eta}$, $X_{2, \eta} \in \mathbb{S}_{\mathbb{H}, \mathbb{H}'}$ 
	satisfying that
	$X_{1, \eta} = \pi^{H'}_{V'} X_{1, \eta} \pi^{H}_V $,
	$X_{2, \eta} = \pi^{H'}_{V'} X_{2, \eta} \pi^{H}_V$,
	\begin{equation}
		\label{eq: X bound}
			-\left(
				\frac {1}{\eps} 
				+ \| A_\eta \|_{L(\mathbb{H}^2, (\mathbb{H}^2)')} 
			\right) 
			I_{\mathbb{H}^2}
		\leq
			\left(
				\begin{array}{cc}
					X_{1, \eta}
					&
					0
					\\
					0
					&
					X_{2, \eta}
				\end{array}
			\right)
		\leq
			\left( 
				\pi^{(H^2)'}_{(V^2)'} A_\eta \pi^{H^2}_{V^2} 
				+ \eps (\pi^{(H^2)'}_{(V^2)'} A_\eta \pi^{H^2}_{V^2})^2 
			\right),
	\end{equation}
	and satisfying for all $i \in \{1, 2\}$ that
	\begin{equation}
	\label{eq: semijets u eta}
		\begin{split}
				&\left( 
					0, 
					X_{i, \eta} 
					+ \pi^{(H^2)'}_i A_\eta E_i 
					-  \pi^{H'}_{V'} \pi^{(H^2)'}_i A_\eta E_i \pi^{H}_{V}   
					+ \delta I_{\mathbb{H}} \pi^{H}_V + \lambda I_{\mathbb{H}} \pi^{H}_{V^\perp}
				\right)
			\in (\hat{J}^2_{ \mathbb{H}, + } \tilde{u}_{i, \eta}) (0). 
		\end{split}
	\end{equation}	
	In addition, it holds for all $\eta \in (0,\overline{\eta}]$
	and all $i \in \{1, 2\}$
	that
	\begin{equation}
	\label{eq: semijets u and u eta}
			(\hat{J}^2_{ \mathbb{H}, + } \tilde{u}_{i, \eta}) (0)
		=
			(\hat{J}^2_{ \mathbb{H}, + } \tilde{u}_{i}) (0)
		=
			(\hat{J}^2_{ \mathbb{H}, + } u_i) (\overline{z}_i) 
			- \big \{ \big(
					\pi^{(H^2)'}_i (D_{\mathbb{H}^2} \Psi)
						(\overline{z}_1, \overline{z}_2), 
					0
				\big ) \big \}.
	\end{equation}
	Combining \eqref{eq: semijets u eta} and
	\eqref{eq: semijets u and u eta}
	shows then that for all $\eta \in (0,\overline{\eta}]$ and all $i \in \{1, 2\}$
	it holds that
	\begin{equation}
	\label{eq: semijets u}
		\begin{split}
			&\left( 
					\pi^{(H^2)'}_i (D_{\mathbb{H}^2} \Psi)
						(\overline{z}_1, \overline{z}_2), 
					X_{i, \eta} + \pi^{(H^2)'}_i A_\eta E_i 
					-  \pi^{H'}_{V'} \pi^{(H^2)'}_i A_\eta E_i \pi^{H}_{V}   
					+ \delta I_{\mathbb{H}} \pi^{H}_V 
					+ \lambda I_{\mathbb{H}} \pi^{H}_{V^\perp}
				\right) 
			\\ & \qquad 
			\in (\hat{J}^2_{ \mathbb{H}, + } u_i) (\overline{z}_i).
		\end{split}
	\end{equation}
	To conclude the theorem we need to let $\eta$ go to $0$ and ensure
	that the limits 
	$ \lim_{\eta \to 0} X_{1, \eta}$ and $ \lim_{\eta \to 0} X_{2, \eta}$
	exists. 
	To this end, note that 
	\eqref{eq: def A eta} and \eqref{eq: X bound} show that
	$(X_{1,\eta})_{\eta \in (0,\overline{\eta}]}$ 
	and $(X_{2,\eta})_{\eta \in (0,\overline{\eta}]}$ 
	are uniformly $\mathbb{L}(\mathbb{H}, \mathbb{H}')$-bounded
	and combining this with the fact that the Banach space
	$
		(
			\{
				A \in \mathbb{S}_{\mathbb{H}, \mathbb{H}'} \colon A=  \pi^{H'}_{V'} A \pi^{H}_V 
			\},
			\| \cdot \|_{L(\mathbb{H}, \mathbb{H}')}
		)
	$
	is isomorphic to the finite-dimensional Banach space
	$ (\mathbb{S}_{\mathbb{V}, \mathbb{V}'}, \| \cdot \|_{L(\mathbb{V}, \mathbb{V}')})$
	yields that there exist a sequence $(\eta_i)_{i \in \N} \subseteq (0,\overline{\eta}]$ 
	and operators 
	$X_1$, 
	$
		X_2 \in 
			\{
				A \in \mathbb{S}_{\mathbb{H}, \mathbb{H}'} \colon A = \pi^{H'}_{V'} A \pi^{H}_V 
			\}
	$
	such that
	$\lim_{i \to \infty} \eta_{i} = 0$, 
	$\lim_{i \to \infty} \| X_{1, \eta_i} - X_1 \|_{L(\mathbb{H}, \mathbb{H}')} = 0$, and that 
	$\lim_{i \to \infty} \| X_{2, \eta_i} - X_2 \|_{L(\mathbb{H}, \mathbb{H}')} = 0$.
	This together with \eqref{eq: semijets u}, the fact that
	$
		\lim_{i \to \infty} 
			\| 
				A_{\eta_i} 
				- (D^2_{\mathbb{H}^2} \Psi) (\overline{z}_1,\overline{z}_2) 
			\|_{L(\mathbb{H}^2, (\mathbb{H}^2)')}
		= 0
	$,
	and the fact that $\hat{J}^2_{ \mathbb{H}, +}$
	is closed under limits
	then completes the proof of Corollary \ref{cor: suitable semijets function}.
\end{proof}
The next corollary is an application of Corollary 
\ref{cor: suitable semijets function} to special functions.
It will will be used in the proof of 
Lemma \ref{l:technical.lemma.uniqueness} below.  
\begin{corollary}[Construction of suitable semijets for strict maxima 
	wrt.~to a $C^2$ function with suitable 2nd derivative]
\label{cor: suitable semijets special case}
	Let $\mathbb{H} = (H, \langle \cdot, \cdot \rangle_H, \| \cdot \|_H)$
	be a real Hilbert space,
	let $\mathbb{V} = (V, \langle \cdot, \cdot \rangle_H, \| \cdot \|_H) \subseteq \mathbb{H}$ 
	be a finite-dimensional linear subspace of $\mathbb{H}$,
	let $\mathbb{H}'=(H', \|\cdot \|_{H'})$,
	$\mathbb{V}'=(V', \|\cdot \|_{V'})$,
	$(\mathbb{H}^2)'=((H^2)', \|\cdot \|_{(H^2)'})$,
	$(\mathbb{V}^2)'=((V^2)', \|\cdot \|_{(V^2)'})$,
	be the dual spaces,
	let $O,U \subseteq H$ be open sets,
	let $u_1\colon O \to  [- \infty, \infty)$, 
	$u_2 \colon U \to  [- \infty, \infty)$ be upper semicontinuous
	with respect to the $\| \cdot \|_{\mathbb{H}}$-norm,
	let $\Psi \in \C_{\mathbb{H}^2}^2(O \times U, \R)$,
	let $(\overline{z}_1$, $\overline{z}_2) \in O \times U$ be a strict 
	$\mathbb{H}^2$-maximum of
	$
		H \times H \ni (z_1,z_2) \to u_1(z_1) + u_2(z_2) -\Psi(z_1, z_2) \in  [- \infty, \infty)
	$,
	denote by $E_1$, $E_2 \colon H \to H^2$ the canonical embeddings
	satisfying for all $x \in H$ that
	$E_1(x)=(x,0)$ and that $E_2(x)=(0,x)$, 
	assume that
	$
			\pi^{H'}_{V'} \pi^{(H^2)'}_1 
				((D^2_{\mathbb{H}^2} \Psi) (\overline{z}_1, \overline{z}_2))  
					E_2 \pi^{H}_{V^\perp}
		= 0
	$, that
	$
			\pi^{H'}_{(V^\perp)'} \pi^{(H^2)'}_1 
				((D^2_{\mathbb{H}^2} \Psi) (\overline{z}_1, \overline{z}_2)) 
					E_2 \pi^{H}_{V}
		= 0
	$, that
	$
			\pi^{H'}_{(V^\perp)'} \pi^{(H^2)'}_1 
				((D^2_{\mathbb{H}^2} \Psi) (\overline{z}_1, \overline{z}_2)) 
					E_1 \pi^{H}_{V}
		= 0
	$, and that
	$
			\pi^{H'}_{(V^\perp)'} \pi^{(H^2)'}_2 
				((D^2_{\mathbb{H}^2} \Psi) (\overline{z}_1, \overline{z}_2)) 
					E_2 \pi^{H}_{V}
		= 0,
	$
	and let 
	$\eps$, $\lambda \in (0, \infty) $
	satisfy that
	\begin{equation}
		\begin{split}
			\lambda >
				&\max_{i \in \{1, 2\}}
					\left(
					\|
					\pi^{H'}_{(V^\perp)'} \pi^{(H^2)'}_i 
						((D^2_{\mathbb{H}^2} \Psi) (\overline{z}_1, \overline{z}_2)) 
							E_i \pi^{H}_{V^\perp} 
					\|_{L(\mathbb{H},\mathbb{H}' )}
				\right) \\
				&+\|
					\pi^{H'}_{(V^\perp)'} \pi^{(H^2)'}_1 
						((D^2_{\mathbb{H}^2} \Psi) (\overline{z}_1, \overline{z}_2)) 
							E_2 \pi^{H}_{V^\perp} 
				\|_{L(\mathbb{H},\mathbb{H}' )}.
		\end{split}
	\end{equation}
	Then there exist operators 
	$X_{1}$, $X_{2} \in \mathbb{S}_{\mathbb{H}, \mathbb{H}'}$ 
	satisfying that
	$X_{1} = \pi^{H'}_{V'} X_{1} \pi^{H}_V $, 
	$X_{2} = \pi^{H'}_{V'} X_{2} \pi^{H}_V$,
	\begin{equation}
	\label{eq: X delta bound}
		\begin{split}
				&-\left(
					\frac {1}{\eps} 
					+ \| 
							((D^2_{\mathbb{H}^2} \Psi) (\overline{z}_1, \overline{z}_2)) 
						\|_{L(\mathbb{H}^2, (\mathbb{H}^2)')} 
				\right) 
				I_{\mathbb{H}^2}
			\leq
				\left(
					\begin{array}{cc}
						X_{1}
						&
						0
						\\
						0
						&
						X_{2}
					\end{array}
				\right) \\
			\leq{}&
				\left( 
					\pi^{(H^2)'}_{(V^2)'} 
						((D^2_{\mathbb{H}^2} \Psi) (\overline{z}_1, \overline{z}_2))
							\pi^{H^2}_{V^2} 
					+ \eps (\pi^{(H^2)'}_{(V ^2)'} 
						((D^2_{\mathbb{H}^2} \Psi) (\overline{z}_1, \overline{z}_2))
							\pi^{H^2}_{V^2})^2 
				\right),
		\end{split}
	\end{equation}
	and satisfying for all $i \in \{1,2\}$ that 
	\begin{align}
			&\left( 
				\pi^{(H^2)'}_i ((D_{\mathbb{H}^2} \Psi)
					(\overline{z}_1, \overline{z}_2)), 
				X_{i} + \lambda I_{\mathbb{H}} \pi^{H}_{V^\perp}
			\right) 
		\in (\hat{J}^2_{ \mathbb{H}, + } u_i) (\overline{z}_i). 
	\end{align}
\end{corollary}

\begin{proof}[Proof of Corollary \ref{cor: suitable semijets special case}]
	For the rest of the proof let $\tilde{\lambda} \in (0, \infty)$ be the real number 
	satisfying that
	\begin{equation}
			\tilde{\lambda} 
		= 
			\lambda
			-\max_{i \in \{1, 2\} } \big(
				\|
				\pi^{H'}_{(V^\perp)'} \pi^{(H^2)'}_i 
					((D^2_{\mathbb{H}^2} \Psi) (\overline{z}_1, \overline{z}_2)) 
						E_i \pi^{H}_{V^\perp} 
				\|_{L(\mathbb{H},\mathbb{H}' )}
			\big),
	\end{equation}
	Then it follows from Corollary \ref{cor: suitable semijets function} 
	and from the assumptions that
	for every $\delta \in (0, \infty)$
	there exist operators 
	$X_{1,\delta}$, $X_{2,\delta} \in \mathbb{S}_{\mathbb{H}, \mathbb{H}'}$ 
	satisfying that
	$X_{1,\delta} = \pi^{H'}_{V'} X_{1,\delta} \pi^{H}_V$, 
	$X_{2,\delta} = \pi^{H'}_{V'} X_{2,\delta} \pi^{H}_V$,
	\begin{equation}
		\begin{split}
				&-\left(
					\frac {1}{\eps} 
					+ \| 
							((D^2_{\mathbb{H}^2} \Psi) (\overline{z}_1, \overline{z}_2)) 
						\|_{L(\mathbb{H}^2, (\mathbb{H}^2)')} 
				\right) 
				I_{\mathbb{H}^2}
			\leq
				\left(
					\begin{array}{cc}
						X_{1,\delta}
						&
						0
						\\
						0
						&
						X_{2,\delta}
					\end{array}
				\right) \\
			\leq{}&
				\left( 
					\pi^{(H^2)'}_{(V^2)'} 
						((D^2_{\mathbb{H}^2} \Psi) (\overline{z}_1, \overline{z}_2))
							\pi^{H^2}_{V^2} 
					+ \eps (\pi^{(H^2)'}_{(V^2)'} 
						((D^2_{\mathbb{H}^2} \Psi) (\overline{z}_1, \overline{z}_2))
							\pi^{H^2}_{V^2})^2 
				\right),
		\end{split}
	\end{equation}
	and satisfying that for all $i \in \{1, 2\}$ 
	it holds that 
	\begin{align}
		\nonumber
			&\Big( 
				\pi^{(H^2)'}_i ((D_{\mathbb{H}^2} \Psi) 
					(\overline{z}_1, \overline{z}_2)), 
				X_{i,\delta} 
				+ \pi^{(H^2)'}_i 
					((D^2_{\mathbb{H}^2} \Psi) (\overline{z}_1, \overline{z}_2)) E_i \\
			\label{eq: lies in semijet u1}
			& \qquad \qquad
				- \pi^{H'}_{V'} \pi^{(H^2)'}_i 
					((D^2_{\mathbb{H}^2} \Psi) (\overline{z}_1, \overline{z}_2))
						E_i \pi^{H}_{V}   
				+ \delta I_{\mathbb{H}} \pi^{H}_V 
				+ \tilde{\lambda} I_{\mathbb{H}} \pi^{H}_{V^\perp}
			\Big)
		\in (\hat{J}^2_{ \mathbb{H}, + } u_i) (\overline{z}_i). 
	\end{align}
	Moreover, we obtain from the assumptions
	and from the fact that for all $i \in \{1, 2\}$
	it holds that
	\begin{equation}
			\pi^{H'}_{V'} \pi^{(H^2)'}_i 
					((D^2_{\mathbb{H}^2} \Psi) (\overline{z}_1, \overline{z}_2)) 
						E_i \pi^{H}_{V^\perp}
		=
			(\pi^{H'}_{(V^\perp)'} \pi^{(H^2)'}_i 
					((D^2_{\mathbb{H}^2} \Psi) (\overline{z}_1, \overline{z}_2)) 
						E_i \pi^{H}_{V} )^*
	\end{equation}
	that for all $i \in \{1, 2\}$
	it holds that
	\begin{align}
		\nonumber
				&\pi^{(H^2)'}_i 
					((D^2_{\mathbb{H}^2} \Psi) (\overline{z}_1, \overline{z}_2)) 
						E_i
				-\pi^{H'}_{V'} \pi^{(H^2)'}_i 
					((D^2_{\mathbb{H}^2} \Psi) (\overline{z}_1, \overline{z}_2)) 
						E_i \pi^{H}_{V} \\
		\nonumber
			={}&
				(\pi^{H'}_{V'} + \pi^{H'}_{(V^\perp)'}) \pi^{(H^2)'}_i 
					((D^2_{\mathbb{H}^2} \Psi) (\overline{z}_1, \overline{z}_2)) 
						E_i (\pi^{H}_{V} + \pi^{H}_{V^\perp})
				-\pi^{H'}_{V'} \pi^{(H^2)'}_i 
					((D^2_{\mathbb{H}^2} \Psi) (\overline{z}_1, \overline{z}_2)) 
						E_i \pi^{H}_{V} \\
		\nonumber
			={}&
				\pi^{H'}_{V'} \pi^{(H^2)'}_i 
					((D^2_{\mathbb{H}^2} \Psi) (\overline{z}_1, \overline{z}_2)) 
						E_i \pi^{H}_{V^\perp}
				+\pi^{H'}_{(V^\perp)'} \pi^{(H^2)'}_i 
					((D^2_{\mathbb{H}^2} \Psi) (\overline{z}_1, \overline{z}_2)) 
						E_i \pi^{H}_{V}  \\
			\nonumber
				&+\pi^{H'}_{(V^\perp)'} \pi^{(H^2)'}_i 
					((D^2_{\mathbb{H}^2} \Psi) (\overline{z}_1, \overline{z}_2)) 
						E_i \pi^{H}_{V^\perp} \\
			\label{eq: operator difference}
			={}&
				\pi^{H'}_{(V^\perp)'} \pi^{(H^2)'}_i 
					((D^2_{\mathbb{H}^2} \Psi) (\overline{z}_1, \overline{z}_2)) 
						E_i \pi^{H}_{V^\perp}.
	\end{align}
	Thus it follows from \eqref{eq: operator difference} that
	for all $i \in \{1, 2\}$
	it holds that
	\begin{equation}
		\begin{split}	
				&\pi^{(H^2)'}_i 
					((D^2_{\mathbb{H}^2} \Psi) (\overline{z}_1, \overline{z}_2)) 
						E_i
				-\pi^{H'}_{V'} \pi^{(H^2)'}_i 
					((D^2_{\mathbb{H}^2} \Psi) (\overline{z}_1, \overline{z}_2)) 
						E_i \pi^{H}_{V}
				+ \tilde{\lambda} I_{\mathbb{H}} \pi^{H}_{V^\perp} \\
			\leq{}&
				\|
					\pi^{H'}_{(V^\perp)'} \pi^{(H \times H)'}_i 
						((D^2_{\mathbb{H}^2} \Psi) (\overline{z}_1, \overline{z}_2)) 
							E_i \pi^{H}_{V^\perp}
				\|_{L(\mathbb{H}, \mathbb{H}')} I_{\mathbb{H}} \pi^{H}_{V^\perp}
				+ \tilde{\lambda} I_{\mathbb{H}} \pi^{H}_{V^\perp}
			\leq
				\lambda I_{\mathbb{H}} \pi^{H}_{V^\perp}.
		\end{split}
	\end{equation}
	Combining this with 
	\eqref{eq: lies in semijet u1} 
	and with the fact that
	$\hat{J}^2_{ \mathbb{H}, + }$ is monotone in the second argument
	then shows 
	for all $i \in \{1,2\}$ and all $\delta \in (0,\infty)$ that
	\begin{align}
				&\left( 
					\pi^{(H^2)'}_i ((D_{\mathbb{H}^2} \Psi) 
						(\overline{z}_1, \overline{z}_2)), 
					X_{i,\delta}   
					+ \delta I_{\mathbb{H}} \pi^{H}_V 
					+ \lambda I_{\mathbb{H}} \pi^{H}_{V^\perp}
				\right) 
			\in (\hat{J}^2_{ \mathbb{H}, + } u_i) (\overline{z}_i). 
	\end{align}
	To this end, note that 
	\eqref{eq: X delta bound} shows that
	$(X_{1,\delta})_{\delta \in (0,\infty)}$ 
	and $(X_{2,\delta})_{\delta \in (0,\infty)}$ 
	are uniformly $\mathbb{L}(\mathbb{H}, \mathbb{H}')$-bounded
	and combining this with the fact that the Banach space
	$
		(
			\{
				A \in \mathbb{S}_{\mathbb{H}, \mathbb{H}'} \colon A=  \pi^{H'}_{V'} A \pi^{H}_V 
			\},
			\| \cdot \|_{L(\mathbb{H}, \mathbb{H}')}
		)
	$
	is isomorphic to the finite-dimensional Banach space
	$ (\mathbb{S}_{\mathbb{V}, \mathbb{V}'}, \| \cdot \|_{L(\mathbb{V}, \mathbb{V}')})$
	yields that there exist a sequence $(\delta_i)_{i \in \N} \subseteq (0,\infty)$ 
	and operators 
	$X_1$, 
	$
		X_2 \in 
			\{
				A \in \mathbb{S}_{\mathbb{H}, \mathbb{H}'} \colon A = \pi^{H'}_{V'} A \pi^{H}_V 
			\}
	$
	such that
	$\lim_{i \to \infty} \delta_{i} = 0$, 
	$\lim_{i \to \infty} \| X_{1, \delta_i} - X_1 \|_{L(\mathbb{H}, \mathbb{H}')} = 0$, 
	and that 
	$\lim_{i \to \infty} \| X_{2, \delta_i} - X_2 \|_{L(\mathbb{H}, \mathbb{H}')} = 0$.
	This together with the fact that $\hat{J}^2_{ \mathbb{H}, +}$
	is closed under limits
	then completes the proof of Corollary \ref{cor: suitable semijets special case}.
\end{proof}
\section{Uniqueness of viscosity solutions of nonlinear 2nd order PDEs}
\label{sec: uniqueness 2nd order PDE}
The next lemma establishes a comparison result for
bounded viscosity solutions
vanishing at infinity. It generalizes
Lemma 4.10 in Hairer, Hutzenthaler \& Jentzen 
\cite{HairerHutzenthalerJentzen2015}
(which assumed finite-dimensional Hilbert spaces
and used the classical notion of viscosity solutions)
and in contrast to Theorem 6.1 in Ishii \cite{Ishii1993}
we don’t have a one sided global uniform continuity condition.
\begin{lemma}[A domination result for viscosity subsolutions] 
  \label{l:technical.lemma.uniqueness}
	Assume the setting in Section \ref{ssec: Setting H X with t},
	assume that $O$ is convex,
	and assume that for all $R \in (0, \infty)$ 
	there exists a $\Lambda_R \in (0, \infty)$ such that
	for all $t \in (0,T)$ and all $x \in X \cap O$ with $\|x \|_H \leq R$ it holds that
	\begin{equation}
	\label{eq: h_t bound}
		\left | \frac{\partial}{\partial t} h(t,x) \right | \leq \Lambda_R \cdot h(t,x).
	\end{equation}
	Moreover,
  let
  $
    G_1, G_2 
    \colon 
    W
		\times \R \times H' \times \mathbb{S}_{\mathbb{X}, \mathbb{X}'} \to \R
  $
  be degenerate elliptic functions
  and let
  $ 
    u_1, u_2 \colon [0,T] \times O \to \R \cup \{- \infty \}
  $
  be bounded from above  
	on every $\R \times \mathbb{H}$-bounded subset of $[0,T] \times O$ such 
  that for every
  $ i \in \{ 1, 2 \} $
  it holds that 
  $ u_i|_{ (0,T) \times O } $ 
  is  
  a viscosity subsolution of
  \begin{equation} \label{eq:parabolic.equation.i}
    \tfrac{ \partial }{ \partial t }
    u_i(t,x) -
    G_i\big( (t, x), u_i(t,x), (D_\mathbb{H} u_i)(t,x),
      ((D^2_\mathbb{H} u_i)(t,x) |_X) |_X
    \big) = 0
  \end{equation}
	for 
  $ (t,x) \in (0,T) \times O $ relative to 
	$(h, \R \times \mathbb{H}, \R \times \mathbb{X})$.
	Furthermore, assume 
	that for all $ R \in (0,\infty)$,
	$\delta \in (0,1]$,
	and all $ i \in \{ 1, 2 \} $
	it holds that
	\begin{equation}
	\label{eq: uniformly bounded at p}
		\begin{split}
			\lim_{\eps \to 0} \Big [
				\sup \{
					& |(G_i)_{\mathbb{H}, \mathbb{X}, \delta, h}^{+}((t,x),r,p,A)
						-(G_i)_{\mathbb{H}, \mathbb{X}, \delta, h}^{+}((t,x),r,q,A)| 
					\colon
					~ (t, x) \in W,
					~ r \in \R, \\
					~ &p, q \in H', 
					~ A \in \mathbb{S}_{\mathbb{H}, \mathbb{H}'}, 
					~\max \{ h(t,x), \|x\|_H, |r|, \|p\|_{H'}, \| A \|_{L(\mathbb{H}, \mathbb{H}')} \} 
						\leq R, \\
					~&\|p-q\|_{H'} \leq \eps, 
					~(G_i)_{\mathbb{H}, \mathbb{X}, \delta, h}^{+}((t,x),r,p,A) \geq -R\}
			\Big ] = 0,
		\end{split}
	\end{equation}
	assume that there exist an increasing sequence of finite-dimensional
	linear subspaces 
	$H_1 \subseteq H_2 \subseteq \ldots \subseteq H$
	and a function 
	$
		m \colon 
			(0, \infty) \times (0,T) \times O \times \R \times H' \times (0,1] \to (0, \infty)
	$
	satisfying for all
	$R \in (0,\infty)$, $t \in (0,T)$, $x \in O$, $r \in \R$, $p \in H'$,  
	$\delta \in (0,1]$, and all $ i \in \{ 1, 2 \} $
	that
	$\cup_{N =1}^\infty H_N$ 
	is dense in $H$ with respect to the $\| \cdot \|_H$-norm
	and that
	\begin{equation}
	\label{eq: uniformly bounded in dimension}
		\begin{split}
			\lim_{N \to \infty} \big [
				\sup \{ 
					&(G_i)_{\mathbb{H}, \mathbb{X}, \delta, h}^{+}
						((\tau, \xi), \nu, \rho, A + \alpha I_\mathbb{H} \pi^{H}_{H^\perp_N})
					-(G_i)_{\mathbb{H}, \mathbb{X}, \delta, h}^{+}
						((\tau, \xi), \nu, \rho, A) \colon 
				~ \alpha \in (0, R), \\
				~ &\nu \in \R, 
				~ (\tau, \xi) \in W,
				~ \rho \in H', 
				~ A \in \mathbb{S}_{\mathbb{H},\mathbb{H}'}, 
				~ h(\tau, \xi) \leq R, 
				~ \| A \|_{L(\mathbb{H},\mathbb{H}')} \leq R, \\
				~ &| t-\tau | \vee \|x- \xi\|_H \vee |r -\nu| \vee \|p - \rho\|_{H'} 
						\leq m(R,t,x,r,p,\delta), \\
				~ &(G_i)_{\mathbb{H}, \mathbb{X}, \delta, h}^{+}
						((\tau, \xi), \nu, \rho, A) \geq -R \}
			\big ]
				\leq 0,
		\end{split}
	\end{equation}
  assume that 
	there exist sequences $(\tilde{\beta}_n)_{n \in \N} \subseteq (0, \infty)$
	and $ (\tilde{\delta}_n)_{n \in \N} \in (0,1]$
	such that
	$
		\lim_{n \to \infty} \tilde{\delta}_n =0 
	$
	and such that
	for all 
  $ i \in \{ 1, 2 \} $ and all
  $ 
    ( 
      (t_i^{ (n) }, x_i^{ (n) }), r_i^{ (n) }, A_i^{ (n) } 
    )_{n \in \N}
    \in
    W \times \R \times \mathbb{S}_{\mathbb{H}, \mathbb{H}'} 
  $,
  satisfying that
  $
    w-\lim_{ n \to \infty }
    ( t_1^{ (n) }, x_1^{ (n) } )
    \in [ 0, T ) \times O
  $, that
  $
    \lim_{ n \to \infty }
    \big(
      \sqrt{ n } 
      \| 
         x_2^{ (n) }
        -
        x_{ 1 }^{ (n) }
      \|_{H}
    \big)
    = 0
  $,
	  that
  $
    \lim_{ n \to \infty }
    \Big(
      \sqrt{ \tilde{\beta}_n } 
      | 
         t_2^{ (n) }
        -
        t_{ 1 }^{ (n) }
      |
    \Big)
    = 0
  $,
	that
  $
    \lim_{ n \to \infty }
    \big(
      \tilde{\delta}_n 
        \big( h(t_{ 1 }^{ (n) }, x_{ 1 }^{ (n) })
        +
        h(t_{ 2 }^{ (n) }, x_{ 2 }^{ (n) }) \big)
    \big)
    = 0
  $,
	that
  $
    \lim_{ n \to \infty }
	$
	$
    \sum_{ i = 1 }^2 
    r_i^{ (n) } 
    > 0
  $,
  that 
  $
    \sup_{ n \in \N } 
    \sum_{ i = 1 }^2
    | r_i^{ (n) } |
    < \infty
  $,
  and that for all
  $
		R \in (0, \infty)
	$
	and all
	$
		(z^{(n)}_1)_{n \in \N},
	$
	$
		(z^{(n)}_2)_{n \in \N} 
			\subseteq \{z \in H \colon ~ \| z \|_H \leq R\}
	$
	with
	$
		\limsup_{n \to \infty}
    \sum_{ i = 1 }^2
    \langle z_i^{ (n) }, A_i^{(n)} z_i^{ (n) }
    \rangle_{H, H'}
  \leq 
		\limsup_{n \to \infty} (
			3 \| z_2^{ (n) } 
	$
	$
		- z^{ (n) }_{ 1 } \|_H^2
		)
  $
	it holds that
  \begin{equation}  
  \label{eq:00assumption}
    \limsup_{ 
      n \to \infty 
    }
    \left[
    \sum_{ i = 1 }^2 
      (G_i)_{\mathbb{H}, \mathbb{X},\tilde{\delta}_n, h}^+\left(
         (t_i^{(n)}, x_i^{(n)}),
         r_i^{ (n) } ,
         n \, I_\mathbb{H}
           \big( 
             x_i^{ (n) } -
             x_{ 3 - i }^{ (n) }
           \big),
         n \,
         A_i^{ (n) } 
      \right) 
    \right]
    \leq 0,
  \end{equation}
  assume that for all $R \in (0,\infty)$ it holds that
	\begin{equation}
	\label{eq: uniformly bounded at 0}
		\begin{split}
				\lim_{r \downarrow 0} \lim_{\eps \downarrow 0}	
				\sup \Bigg \{
					\sum_{ i = 1 }^2 u_i(t_i,x_i) 
					\colon
					~&(t_1, x_1), (t_2, x_2) \in W,
					~h(t_1, x_1) \vee h(t_2, x_2) \leq R, \\
					~&\|x_1-x_2\|_H \leq r,
					~t_1 \vee t_2 \leq \eps 
				\Bigg \}
			\leq 
				0,
		\end{split}
	\end{equation}
  and assume that for all $i \in \{ 1, 2\}$ it holds that
  \begin{equation}    
  \label{eq:attains.maximum}
			\lim_{ n \to \infty }
				\sup_{
						(t,x) \in 
						([0,T] \times O_n^c) \cap W
						}
				u_i(t,x)
			\leq 0.
  \end{equation}
  Then 
	for all $\tilde{T} \in (0,T)$ 
	and all $R \in (0,\infty)$ it holds that
  \begin{equation}
		\begin{split}
			\lim_{r \downarrow 0} \lim_{\eps \downarrow 0}
				\sup \Bigg \{
					&\sum_{ i = 1 }^2 u_i(t_i,x_i) \colon 
					~(t_1, x_1), (t_2, x_2) \in W, 
					~ t_1 \vee t_2 \leq \tilde{T}, \\
					~&h(t_1, x_1) \vee h(t_2, x_2) \leq R, 
					~\|x_1-x_2\|_H \leq r, 
					~|t_1-t_2| \leq \eps 
				\Bigg \}
			\leq
				0.
		\end{split}
  \end{equation}
\end{lemma}

\begin{proof}[Proof
of
Lemma~\ref{l:technical.lemma.uniqueness}]
  If $ O = \emptyset $, 
  then the assertion is trivial. 
  So for the rest of 
  the proof, we assume that 
  $ O \neq \emptyset $.
  We will show
  that for all 
  $ 
    \mu \in (0,1] 
  $
	and all 
	$ 
    R \in (0,\infty) 
  $
	it holds that
  \begin{equation}
	\label{eq: assertion with extra mu term}
	\begin{split}
		\lim_{r \downarrow 0} \lim_{\eps \downarrow 0}
			\sup \Bigg \{
				&\sum_{ i = 1 }^2 \left( 
					u_i(t_i,x_i) - \tfrac{ \mu }{ (T - t_i) } 
				\right ) \colon
				~(t_1, x_1), (t_2, x_2) \in W, \\
				&~h(t_1, x_1) \vee h(t_2, x_2) \leq R, 
				~\|x_1-x_2\|_H \leq r, 
				~|t_1-t_2| \leq \eps 
			\Bigg \}
    \leq
			0.
	\end{split}
  \end{equation}
  Letting $\mu \to 0 $ 
  then will yield that for all $\tilde{T} \in (0,T)$ 
	and all $ R \in (0,\infty)$ it holds that 
  \begin{equation}
		\begin{split}
			\lim_{r \downarrow 0} \lim_{\eps \downarrow 0}
				\sup \Bigg \{
					&\sum_{ i = 1 }^2 u_i(t_i,x_i) \colon 
					~(t_1, x_1), (t_2, x_2) \in W,
					~ t_1 \vee t_2 \leq \tilde{T}, \\
					&~h(t_1, x_1) \vee h(t_2, x_2) \leq R, 
					~\|x_1-x_2\|_H \leq r, 
					~|t_1-t_2| \leq \eps 
				\Bigg \}
			\leq
				0.
		\end{split}
  \end{equation}
  For the rest of the proof we thus fix 
  $ \mu \in (0,1] $.
  In a first step of this proof, 
  we modify the problem.
  More precisely,
  denote by
	$W_n$, $n \in \N$, the sets satisfying for all $n \in \N$ that
	$W_n = ([0,T] \times O_n) \cap W$, by
	$W_n^c$, $n \in \N$, the sets satisfying for all $n \in \N$ that
	$W_n^c = ([0,T] \times O_n^c) \cap W$, by
  $
    \tilde{u}_1,
    \tilde{u}_2
    \colon
    [0,T] \times O
    \to [ - \infty, \infty )
  $
	the functions satisfying for all
	$i \in \{1, 2\}$, $t \in [0,T]$, and all $x \in O$
  that
  $
    \tilde{u}_i(t,x)
    =
    u_i( t, x ) - 
    \tfrac{ \mu }{ (T - t) }
  $,
  and by
  $
    \tilde{G}_1,
    \tilde{G}_2 \colon 
    W 
    \times \R \times H'
    \times \mathbb{S}_{\mathbb{X},\mathbb{X}'}
		\to \R
  $
	the functions satisfying for all
	$i \in \{1, 2\}$ and all $( (t, x), r, p, A) \in 
    W 
    \times \R \times H'
    \times \mathbb{S}_{\mathbb{X},\mathbb{X}'}
	$
  that
  \begin{equation}
  \label{eq:def_tildeG}
			\tilde{G}_i\left( (t, x), r, p, A
			\right) 
			={}
			G_i\!\left(
				(t, x),
				r + \tfrac{ \mu }{ (T - t) }
				,
				p ,
				A
			\right)
			-
			\tfrac{ \mu }{ (T - t)^2 }.
  \end{equation}
  Then it holds 
  for every $ i \in \{ 1, 2 \} $
  that
  $ 
    \tilde{u}_i|_{ (0,T) \times O } 
  $ 
  is
	a viscosity subsolution of
  \begin{equation} 
  \label{eq:parabolic.equation.tilde.ui}
    \tfrac{ \partial }{ \partial t }
    \tilde{u}_i(t,x) -
    \tilde{G}_i\big( (t, x), \tilde{u}_i(t,x), (D_\mathbb{H} \tilde{u}_i)(t,x),
      (D^2_\mathbb{H}\tilde{u}_i)(t,x)
    \big) = 0
  \end{equation}
	for 
  $ (t,x) \in (0,T) \times O $ relative to $(h, \R \times \mathbb{H}, \R \times \mathbb{X})$.
  It remains to prove that 
	for all $R \in (0,\infty)$ it holds that
  \begin{equation}
		\begin{split}
			\lim_{r \downarrow 0} \lim_{\eps \downarrow 0}
				\sup \Bigg \{ 
					\sum_{ i = 1 }^2 \tilde{u}_i(t_i,x_i) \colon 
					~&(t_1, x_1), (t_2, x_2) \in W,
					~h(t_1, x_1) \vee h(t_2, x_2) \leq R, \\
					~&\|x_1-x_2\|_H \leq r, 
					~|t_1-t_2| \leq \eps 
				\Bigg \}
			\leq
				0.
		\end{split}
  \end{equation}
  Aiming at a contradiction, let
	$S_0 \in (-\infty, \infty]$ be
  the extended real number satisfying that
  \begin{equation}
	\label{eq: def of S0}
		\begin{split}
			S_{0}
			=
				\lim_{R \to \infty} \lim_{r \downarrow 0} \lim_{\eps \downarrow 0}
					\sup \Bigg \{
						&\sum_{ i = 1 }^2 \tilde{u}_i(t_1,x_1) \colon 
						~(t_1, x_1), (t_2, x_2) \in W, \\
						~&h(t_1, x_1) \vee h(t_2, x_2) \leq R, 
						~\|x_1-x_2\|_H \leq r, 
						~|t_1-t_2| \leq \eps 
					\Bigg \}		
		\end{split}
  \end{equation}
  and assume that 
  $
    S_{0} 
    \in
    (0,\infty]
  $.
  Assumption~\eqref{eq:attains.maximum} 
  then implies that
  there exists a natural
  number
  $ n_0 \in \N $ 
  such that $W_{n_0-1}$ is non-empty
  and such that
	for all 
  $
   z \in 
   W_{n_0-1}^c
  $
	it holds that 
  \begin{equation}
  \label{eq:sum_tildeu_outside}
			\tilde{u}_1( z ) \vee \tilde{u}_2(z)
		\leq
			u_1( z ) \vee u_2(z)
		\leq 
			\min( 1, \tfrac{ S_0 }{ 4 } ). 
  \end{equation}
	Denote now by $K$ the set satisfying $K = O_{n_0}$ that
  and note that the
  functions 
	$u_1$, $u_2$
  are bounded from above
  on the $\R \times \mathbb{H}$-bounded set
  $ [0,T] \times K$. 
	Thus there exists an $L \in (0, \infty)$
	such that for all $z \in [0,T] \times K$ it holds that
	\begin{equation}
	\label{eq: u bound on K}
		\tilde{u}_1(z) \vee \tilde{u}_2(z) \leq u_1(z) \vee u_2(z) \leq L.
	\end{equation}
  Combining this with \eqref{eq:sum_tildeu_outside}
  proves that
  $
    S_{0} < \infty
  $
  and we thus get
  $
    S_{0} \in ( 0, \infty ) 
  $.
  For several
  $ n \in \N $ and suitable $p_n \in ((\R \times H)^2)'$ and $\beta_n \in (0, \infty)$
  we will apply 
  Corollary~\ref{cor: suitable semijets special case}
	 \big(with 
		$
			\CO \leftarrow (0,T) \times O
		$,
		$
			\varepsilon \leftarrow \frac{ 1 }{ n }
		$,
		$
				\Psi
			\leftarrow
				\big(
					((0,T) \times O)^2 \ni \big( (t_1,x_1), (t_2,x_2) \big) \to
						n \| x_1 - x_2 \|_H^2
						+ \beta_n |t_1 - t_2|^2
						+\langle 
							p_n, 
							( (t_1,x_1), (t_2,x_2) ) 
						\rangle_{((\R \times H)^2)', (\R \times H)^2} 
					\in \R
				\big)
		$\big)
  below.
  For this we 
  now check the
  assumptions
  of Corollary~\ref{cor: suitable semijets special case}.
	Denote by
  $
    \eta \colon 
    ( [0,T] \times K )^2 \times [0, \infty)
    \to [-\infty,\infty)
  $
  the function satisfying
  for all 
  $
    z_1, z_2 \in [0,T] \times K
  $
	and all $\delta \in [0, \infty)$
	that
  $
    \eta( z_1, z_2, \delta)
    =
    \sum_{ i = 1 }^2
    (\tilde{u}_i)_{\R \times \mathbb{H}, \delta, h} ^{-,W} ( z_i ).
  $
  Note that we obtain from \eqref{eq: u bound on K} that
	$(\tilde{u}_{1})_{\R \times \mathbb{H}, \delta, h} ^{-,W}$ 
	and $(\tilde{u}_{2})_{\R \times \mathbb{H}, \delta, h} ^{-,W} $
	are bounded from above on the set $[0,T] \times K$ 
	and thus we have
	for every 
  $ \alpha, \beta, \delta \in (0,\infty) $
  that the function
  $
    ( [0,T] \times K )^2
    \ni 
    \big( (t_1,x_1), (t_2,x_2) \big)
    \mapsto
    \eta \big(
       \big( (t_1,x_1), (t_2,x_2) \big), \delta
    \big)
    -
    \frac \alpha 2 \| x_1 - x_2 \|^2_H
		- 
		\frac \beta 2 |t_1 -t_2|^2
    \in
    [ - \infty, \infty )
  $
	is bounded from above.
	Furthermore, it follows from Lemma 5.1.4 in Kato \cite{Kato1980},
	the Banach-Saks Theorem, and from
	the fact that $[0,T] \times K$ is convex, $\R \times \mathbb{H}$-bounded, and
	$\R \times \mathbb{H}$-closed that it is weakly compact. Therefore it follows from
	page $3$ in Stegall \cite{Stegall1978} that $[0,T] \times K$ is also an RNP set 
	and this together with
	the theorem starting on page $4$ in Stegall \cite{Stegall1978} implies that for all  
	$\alpha, \beta, \gamma, \delta \in (0, \infty)$ there exists a
	$
		p_{\alpha, \beta, \gamma, \delta} =
		\big ( 
			(q_{1,\alpha, \beta, \gamma, \delta}, p_{1,\alpha, \beta, \gamma, \delta}),
			(q_{2,\alpha, \beta, \gamma, \delta}, p_{2,\alpha, \beta, \gamma, \delta}) 
		\big ) 
		\in ((\R \times H)^2)'
	$ 
	satisfying for all
	$\alpha, \beta, \gamma, \delta \in (0, \infty)$ 
	that
	$ 
			\| 
				p_{\alpha, \beta, \gamma, \delta}
			\|_{((\R \times H)^2)'} 
		< \gamma 
	$
	and that
	\begin{equation}
		\begin{split}
				&([0,T] \times K)^2 \ni \big( (t_1,x_1), (t_2,x_2) \big)
			\mapsto 
				\eta \big( \big( (t_1,x_1), (t_2,x_2) \big), \delta \big) 
				- \tfrac \alpha2 \|x_1-x_2\|_H^2 \\
			& \qquad \qquad \quad
				- \tfrac \beta2 |t_1 -t_2|^2  
				+ \big \langle 
						p_{\alpha, \beta, \gamma, \delta}, \big( (t_1,x_1), (t_2,x_2) \big) 
					\big \rangle_{((\R \times H)^2)', (\R \times H)^2} 
			\in [-\infty,\infty)
		\end{split}
	\end{equation}
  attains a strict $\mathbb{H}$-maximum
  \begin{equation}
	\label{eq: def of S alpha beta gamma delta}
		\begin{split}
			S_{ \alpha, \beta, \gamma, \delta}
			={} 
			&\sup_{ 
				\big( (t_1,x_1), (t_2,x_2) \big)
				\in 
				( [0,T] \times K )^2
			}
			\Big(
				\eta \big( \big( (t_1,x_1), (t_2,x_2) \big), \delta \big)
				- \tfrac \alpha 2 \| x_1 - x_2 \|_H^2 \\
		&\qquad
				- \tfrac \beta 2 | t_1 - t_2|^2 
				+ \langle 
						p_{\alpha, \beta, \gamma, \delta}, \big( (t_1,x_1), (t_2,x_2) \big) 
					\rangle_{((\R \times H)^2)', (\R \times H)^2}
			\Big)
			< \infty
		\end{split}
  \end{equation}
  in a point
  $
    \underline{z}^{ ( \alpha, \beta, \gamma, \delta) }
    =
    \big(
      \big( 
        t^{ ( \alpha, \beta, \gamma, \delta) }_1, 
        x^{ ( \alpha, \beta, \gamma, \delta) }_1
      \big) ,
      \big( 
        t^{ ( \alpha, \beta, \gamma, \delta) }_2, 
        x^{ ( \alpha, \beta, \gamma, \delta) }_2     
      \big)
    \big)
    \in ( [0,T] \times K )^2
  $.
	In addition, from \eqref{eq: def of S0} 
	and from $S_0 < \infty$
	it follows that
	for all $\eps, r \in (0, \infty)$ 
	and all $l \in (0, \infty]$
	there exist  
	$(t_{l, r, \eps},x_{l, r, \eps})$, $(s_{l, r, \eps},y_{l, r, \eps}) \in W$ and 
	$\tilde{R}_{l} \in (0, \infty)$ such that
	$
			\tilde{u}_1(t_{l, r, \eps}, x_{l, r, \eps})
			+\tilde{u}_2(s_{l, r, \eps}, y_{l, r, \eps})
		\geq 
			S_0 - l
	$,
	that 
	$
			h(t_{l, r, \eps}, x_{l, r, \eps}) \vee h(s_{l, r, \eps}, y_{l, r, \eps}) 
		\leq \tilde{R}_{l}
	$,
	that
	$
			\|x_{l, r, \eps}-y_{l, r, \eps}\|^2_{H} 
		\leq 
			r
	$,
	and that 
	$
			| t_{l, r, \eps} - s_{l, r, \eps} |^2
		\leq 
			\eps
	$.
	Furthermore, it follows from \eqref{eq:sum_tildeu_outside} and 
	from the definition of $O_n$ that for all 
	$\eps \in (0, \infty)$,
	$r \in \left( 0, \frac{1}{n_0(n_0-1)} \right)$,
	and all $l \in (0, \frac{S_0}{2})$ it holds that
	$x_{l, r, \eps}$, $y_{l, r, \eps} \in K$.
	Combining this with the fact
	that for all $x \in K$ it holds that
	$\|x\|_H \leq n_0$
	we get that 
	for all
	$l \in (0, \frac{S_0}{2})$,
	$\alpha \in [n_0, \infty)$, 
	$\beta \in (0, \infty)$,
	$r \in \left( 0, \frac {2l} {\alpha} \right]$,
	$\delta \in \big( 0, \frac {l} {\tilde{R}_{l}} \big]$,
	$\gamma \in (0, l]$,
	and all $\eps \in (0, \frac {2l}{\beta}]$
	it holds that
	\begin{equation}
	\label{eq: S lower bound}
		\begin{split}
					S_{\alpha, \beta, \gamma, \delta} 
				={}& 
					\sup_{\big ( (t,x), (s, y) \big) \in ( [0,T] \times K )^2}
					\Big(
						\eta \big( \big ( (t,x), (s, y) \big), \delta \big)
						- \tfrac \alpha 2 \| x -y \|_H^2
						- \tfrac \beta 2 |t-s|^2 \\
					&\qquad \qquad \qquad \qquad \qquad \qquad
							+ \big \langle 
								p_{\alpha, \beta, \gamma, \delta}, \big ( (t,x), (s, y) \big) 
							\big \rangle_{((\R \times H)^2)', (\R \times H)^2}
					\Big) \\
				\geq{}&
					\overline{(\tilde{u}_1-\delta h)}^{W}_{\mathbb{H}}
						(t_{l, r, \eps} , x_{l, r, \eps} )
					+\overline{(\tilde{u}_2-\delta h)}^{W}_{\mathbb{H}}
						(s_{l, r, \eps} , y_{l, r, \eps} )
					- \tfrac \alpha 2
							\| x_{l, r, \eps} - y_{l, r, \eps} \|_H^2 \\
					&- \tfrac \beta 2
							| t_{l, r, \eps} - s_{l, r, \eps} |^2 
					+ \big \langle 
							p_{\alpha, \beta, \gamma, \delta}, 
							\big( 
								(t_{l, r, \eps} , x_{l, r, \eps} ), 
								(s_{l, r, \eps} , y_{l, r, \eps}) 
							\big)
						\big \rangle_{((\R \times H)^2)', (\R \times H)^2} \\
				\geq{}&
					\tilde{u}_1(t_{l, r, \eps} , x_{l, r, \eps} )
					- \delta h (t_{l, r, \eps} , x_{l, r, \eps} )
					+ \tilde{u}_2(s_{l, r, \eps} , y_{l, r, \eps} )
					- \delta h (s_{l, r, \eps} , y_{l, r, \eps} )
					- \tfrac {\alpha r}{2} \\
					&- \tfrac {\beta \eps}{2} 
					- \gamma 
						\big \|
							\big( (t_{l, r, \eps} , x_{l, r, \eps} ), (s_{l, r, \eps} , y_{l, r, \eps}) \big )
						\big \|_{(\R \times H)^2} \\
				\geq{}&
					S_0 - l 
					-2 \tfrac {l} {\tilde{R}_{l}} \tilde{R}_{l} 
					- \tfrac {2\alpha  l}{2 \alpha}
					- \tfrac {2\beta l}{2 \beta}
					- 2l \cdot (n_0 + T) 
				=
					S_0 - l \cdot (5+ 2 (n_0 + T)).
		\end{split}
	\end{equation}
	Now denote for all $\eps$, $r$, $R \in (0, \infty)$
	by $\Omega_{1, \eps, r ,R} \subseteq W_{n_0}^2$ 
	the set satisfying that
	\begin{equation}
		\begin{split}
				\Omega_{1, \eps, r ,R} 
			= 
				\Big \{
					\big( (t,x),(s,y) \big) \in W_{n_0}^2 \colon
							&|t-s|^2 < \eps, 
							~\|x - y\|_H^2 < r,
							~h(t, x) \vee h(s, y) < R, \\
							~&t^2 \wedge s^2 \geq \eps
				\Big \},
		\end{split}
	\end{equation}
	by $\Omega_{2, \eps, r ,R} \subseteq W_{n_0}^2$ 
	the set satisfying that
	\begin{equation}
		\begin{split}
			\Omega_{2, \eps, r ,R} 
			= 
				\Big \{
					\big( (t,x),(s,y) \big) \in W_{n_0}^2 \colon
							&|t-s|^2 < \eps, 
							~\|x - y\|_H^2 < r,
							~h(t, x) \vee h(s, y) < R, \\
							~&t^2 \wedge s^2 < \eps
				\Big \},
		\end{split}
	\end{equation}
	by $\Omega_{3, \eps, r ,R} \subseteq W_{n_0}^2$ 
	the set satisfying that
	\begin{equation}
			\Omega_{3, \eps, r ,R} 
		= 
			\left \{
				\big( (t,x),(s,y) \big) \in W_{n_0}^2 \colon
						|t-s|^2 \geq \eps
			\right \},
	\end{equation}
	by $\Omega_{4, \eps, r ,R} \subseteq W_{n_0}^2$ 
	the set satisfying that
	\begin{equation}
			\Omega_{4, \eps, r ,R} 
		= 
			\left \{
				\big( (t,x),(s,y) \big) \in W_{n_0}^2 \colon 
						~\|x - y\|_H^2 \geq r
			\right \},
	\end{equation}
	and by $\Omega_{5, \eps, r ,R} \subseteq W_{n_0}^2$ 
	the set satisfying that
	\begin{equation}
			\Omega_{5, \eps, r ,R} 
		= 
			\left \{
				\big( (t,x),(s,y) \big) \in W_{n_0}^2 \colon
						~h(t, x) \vee h(s, y) \geq R
			\right \}.
	\end{equation}
	From \eqref{eq: def of S alpha beta gamma delta},
	the fact that for all $\eps$, $r$, $R \in (0, \infty)$ it holds that
	$\bigcup_{i=1}^{5} \Omega_{i, \eps, r, R} =  W_{n_0}^2$,
	the fact that norm and scalar product are continuous,
	from the density of $W \in O$,
	and from the definition of $\eta$
	it follows that
	for all
	$\alpha$, $\beta$, $\gamma$, $\delta$, $\eps$, $r$, $R \in (0, \infty)$
	it holds that
	\begin{equation}
	\label{eq: splitting S}
		\begin{split}
				S_{\alpha, \beta, \gamma, \delta}
			=
				\max \bigg \{
					&\sup_{\big( (t,x),(s,y) \big) \in \Omega_{i, \eps , r, R}}
					\Big\{
						\eta \big( \big( (t,x),(s, y) \big), \delta \big)
						- \tfrac \alpha 2 \|x - y\|_H^2
						- \tfrac \beta 2 |t-s|^2 \\
						&+ \big \langle 
								p_{\alpha, \beta, \gamma, \delta}, \big( (t,x),(s, y) \big) 
							\big \rangle_{((\R \times H)^2)', (\R \times H)^2}
					\Big \}
					\colon i \in \{1, 2, 3, 4, 5\} 
			\bigg \}.
			\end{split}
	\end{equation}
	Furthermore, 
	we get 
	from \eqref{eq: uniformly bounded at 0},
	\eqref{eq: def of S0},
	from $S_0 \in \R$
	and from $W_{n_0} \subseteq W$ that for all
	$l \in (0, \infty]$ and all $R \in (0, \infty)$
	there exist
	$\tilde{r}_{l,R}, \tilde{\eps}_{l, R} \in (0, \infty)$ such that
	$
			\sup \{
				\tilde{u}_1(t,x) + \tilde{u}_2(s,y) \colon 
				(t,x), 
	$
	$				
					(s,y) \in W_{n_0},
				~h(t, x) \vee h(s, y) < R, 
				~\|x-y\|^2_H < \tilde{r}_{l,R},
				~|t-s|^2 < \tilde{\eps}_{l, R} 
			\}
		\leq
			S_0 + l
	$
	and that
	$
			\sup \{
				\tilde{u}_1(t,x) + \tilde{u}_2(s,y) \colon
				(t,x), (s,y) \in W_{n_0},
				~h(t, x) \vee h(s, y) < R, 
				~\|x-y\|^2_H < \tilde{r}_{l,R},
				~t^2 \vee s^2 < 4 \tilde{\eps}_{l, R} 
			\}
		\leq
			l
	$.
	Thus we obtain from the fact that for all $(t,x) \in W_{n_0}$
	it holds that $\| x \|_H \leq n_0$
	and the fact that for all $(t,x) \in (0, T) \times O$
	it holds that $h(t, x) \geq 0$
	that for all 
	$l \in (0, \infty]$,
	$\delta,
	R,
	\alpha,
	\beta \in (0 , \infty )$,
	and all $\gamma \in ( 0, l]$ it holds that 
	\begin{equation}
	\label{eq: upper bound W1}
		\begin{split}
				&\sup_{
					\big( (t,x), (s,y) \big) \in \Omega_{1, \tilde{\eps}_{l,R}, \tilde{r}_{l,R}, R}
				} 
				\Big \{
					\eta \big( \big( (t,x),(s, y) \big), \delta \big)
					- \tfrac \alpha 2 \|x - y\|_H^2
					- \tfrac \beta 2 |t-s|^2 \\
				& \qquad\qquad\qquad\qquad
					+ \big \langle 
							p_{\alpha, \beta, \gamma, \delta}, \big( (t,x),(s, y) \big) 
					\big \rangle_{((\R \times H)^2)', (\R \times H)^2} 
			\Big \} \\
		\leq{}&
				\sup \Big \{
					\eta \big( \big( (t,x),(s, y) \big), 0 \big)
					+ \gamma \big \| 
							\big( (t,x),(s, y) \big) 
					\big \|_{(\R \times H)^2} \colon 
						(t, x), (s, y) \in W_{n_0}, \\
					&\qquad \qquad\qquad \qquad
						~h(t, x) \vee h(s, y) < R, 
						~\|x-y\|^2_H < \tilde{r}_{l,R},
						~|t-s|^2 < \tilde{\eps}_{l, R} 
				\Big \} \\
			\leq{}&
				S_0 + l \cdot (1 + 2 (n_0 + T)).
		\end{split}
	\end{equation}
	Moreover, we get from
	the fact that for all $(t,x) \in W_{n_0}$
	it holds that $\| x \|_H \leq n_0$,
	from
	the fact that for all $(t,x) \in (0, T) \times O$
	it holds that $h(t, x) \geq 0$
	and from the fact that for all $t$, $s \in (0,T)$ it holds that
	$
			t^2 \vee s^2 
		= 
			((t \vee s) - (t \wedge s) + (t \wedge s))^2
		\leq
			2 (t-s)^2 + 2 (t^2 \wedge s^2)
	$
	that for all
	$l \in (0, \infty]$,
	$\delta,
	R,
	\alpha,
	\beta \in (0 , \infty )$,
	and all $\gamma \in ( 0, l]$
	it holds that
	\begin{equation}
	\label{eq: upper bound W2}
		\begin{split}
				&\sup_{
					\big( (t,x), (s,y) \big) \in \Omega_{2, \tilde{\eps}_{l,R}, \tilde{r}_{l,R}, R}
				}
				\Big \{
					\eta \big( \big( (t,x),(s, y) \big), \delta \big)
					- \tfrac \alpha 2 \|x - y\|_H^2
					- \tfrac \beta 2 |t-s|^2 \\
				& \qquad\qquad
					+ \big \langle 
							p_{\alpha, \beta, \gamma, \delta}, \big( (t,x),(s, y) \big) 
					\big \rangle_{((\R \times H)^2)', (\R \times H)^2}
				\Big \} \\
			\leq{}&
				\sup \Big \{
					\eta \big( \big( (t,x),(s, y) \big), 0 \big)
					+ \gamma \big \| 
							\big( (t,x),(s, y) \big) 
					\big \|_{(\R \times H)^2} \colon 
						(t, x), (s, y) \in W_{n_0}, \\
						&\qquad \qquad
						~h(t, x) \vee h(s, y) < R,
						~\|x-y\|^2_H < \tilde{r}_{l,R}, 
						~|t-s|^2 < \tilde{\eps}_{l, R},
						~t^2 \wedge s^2 < \tilde{\eps}_{l, R} 
				\Big \} \\
			\leq{}&
				\sup \Big \{
					\tilde{u}_1 (t,x) + \tilde{u}_2 (s,y)
					+ 2l \cdot (n_0 +T)\colon 
						(t, x), (s, y) \in W_{n_0},
						~h(t, x) \vee h(s, y) < R, \\
					&\qquad \qquad
						~\|x-y\|^2_H < \tilde{r}_{l,R},
						~t^2 \vee s^2 < 4\tilde{\eps}_{l, R}
				\Big \} \\
			\leq{}&
				l \cdot (1 + 2(n_0 +T)).
		\end{split}
	\end{equation}
	In addition, \eqref{eq: u bound on K},
	the fact that for all $(t,x) \in W_{n_0}$
	it holds that $\| x \|_H \leq n_0$
	and the fact that for all $(t,x) \in (0, T) \times O$
	it holds that $h(t, x) \geq 0$ implies that 
	for all
	$l \in (0, \infty]$,
	$\delta,
	R,
	\alpha \in ( 0, \infty )$,
	$\beta \in \big[ \frac{4L}{\tilde{\eps}_{l,R}}, \infty \big )$
	and all $\gamma \in ( 0, l]$
	it holds that
	\begin{equation}
	\label{eq: upper bound W3}
		\begin{split}
				&\sup_{
					\big( (t,x), (s,y) \big) \in \Omega_{3, \tilde{\eps}_{l,R}, \tilde{r}_{l,R}, R}
				}
				\Big \{
					\eta \big( \big( (t,x),(s, y) \big), \delta \big)
					- \tfrac \alpha 2 \|x - y\|_H^2
					- \tfrac \beta 2 |t-s|^2 \\
				&\qquad \qquad\qquad \qquad\qquad \qquad
					+ \big \langle 
							p_{\alpha, \beta, \gamma, \delta}, \big( (t,x),(s, y) \big) 
					\big \rangle_{((\R \times H)^2)', (\R \times H)^2}
				\Big \} \\
			\leq{}&
				\sup \Big \{
					\eta \big( \big( (t,x),(s, y) \big), 0 \big)
					- \tfrac \beta 2 |t-s|^2
					+ \gamma \big \| 
							\big( (t,x),(s, y) \big) 
					\big \|_{(\R \times H)^2} \colon \\
					&\qquad \qquad\qquad \qquad\qquad \qquad
						(t, x), (s, y) \in W_{n_0}, 
						~|t-s|^2 \geq \tilde{\eps}_{l, R} 
				\Big \} \\
			\leq{}&
				\sup \Big\{
					2 L - \tfrac \beta 2 \tilde{\eps}_{l, R}
					+ \gamma \big \| 
							\big( (t,x),(s, y) \big) 
					\big \|_{(\R \times H)^2} \colon 
						(t, x), (s, y) \in W_{n_0}
				\Big \} \\
			\leq{}&
				2l \cdot(n_0 +T).
		\end{split}
	\end{equation}
	Similarly we have for all
	$l \in (0, \infty]$,
	$\delta, \beta,
	R \in (0,\infty)$,
	$\alpha \in \big[ \frac{4L}{\tilde{r}_{l,R}}, \infty \big )$,
	and all $\gamma \in ( 0, l]$
	that
	\begin{align}
	\nonumber
				&\sup_{
					\big( (t,x), (s,y) \big) \in \Omega_{4, \tilde{\eps}_{l,R}, \tilde{r}_{l,R}, R}
				}
				\Big \{
					\eta \big( \big( (t,x),(s, y) \big), \delta \big)
					- \tfrac \alpha 2 \|x - y\|_H^2
					- \tfrac \beta 2 |t-s|^2 \\ \nonumber
			&\qquad \qquad\qquad \qquad\qquad \qquad\qquad
					+ \big \langle 
							p_{\alpha, \beta, \gamma, \delta}, \big( (t,x),(s, y) \big) 
					\big \rangle_{((\R \times H)^2)', (\R \times H)^2}
				\Big \} \\
	\label{eq: upper bound W4}
			\leq{}&
				\sup \Big \{
					\eta \big( \big( (t,x),(s, y) \big), 0 \big)
					- \tfrac \alpha 2 \|x - y\|_H^2
					+ \gamma \big \| 
							\big( (t,x),(s, y) \big) 
					\big \|_{(\R \times H)^2} \colon \\ \nonumber
					&\qquad \qquad\qquad \qquad\qquad \qquad\qquad
						(t, x), (s, y) \in W_{n_0},
						~\|x-y\|_H^2 \geq \tilde{r}_{l,R} 
				\Big \} \\ \nonumber
			\leq{}&
				\sup \Big \{
					2 L - \tfrac \alpha 2 \tilde{r}_{l,R}
					+ \gamma \big \| 
							\big( (t,x),(s, y) \big) 
					\big \|_{(\R \times H)^2} \colon 
						(t, x), (s, y) \in W_{n_0}
				\Big \} \\ \nonumber
			\leq{}&
				2l \cdot (n_0 +T)
	\end{align}
	and finally we obtain with
	\eqref{eq: u bound on K},
	the fact that for all $(t,x) \in W_{n_0}$
	it holds that $\| x \|_H \leq n_0$, 
	and with the lower semicontinuity of $h$
	with respect to the $\| \cdot \|_{\R \times \mathbb{H}}$-norm that 
	for all
	$l \in (0, \infty]$,
	$\alpha, \beta, \delta \in (0, \infty)$,
	$R \in \big[ \frac{L}{\delta},\infty \big)$,
	and all $\gamma \in (0, l]$
	it holds that
	\begin{equation}
	\label{eq: upper bound W5}
		\begin{split}
				&\sup_{
					\big( (t,x), (s,y) \big) \in \Omega_{5, \tilde{\eps}_{l,R}, \tilde{r}_{l,R}, R}
				}
				\Big \{
					\eta \big( \big( (t,x),(s, y) \big), \delta \big)
					- \tfrac \alpha 2 \|x - y\|_H^2
					- \tfrac \beta 2 |t-s|^2 \\
				&\qquad \qquad\qquad \qquad\qquad \qquad \qquad
					+ \big \langle 
							p_{\alpha, \beta, \gamma, \delta}, \big( (t,x),(s, y) \big) 
					\big \rangle_{((\R \times H)^2)', (\R \times H)^2}
				\Big \} \\
			\leq{}&
				\sup \Big \{
					\eta \big( \big( (t,x),(s, y) \big), 0 \big)
					- \delta (h(t, x) + h(s, y))
					+ \gamma \big \| 
							\big( (t,x),(s, y) \big) 
					\big \|_{(\R \times H)^2} \colon \\
				&\qquad \qquad\qquad \qquad\qquad \qquad\qquad
						(t, x), (s, y) \in W_{n_0},
						~h(t, x) \vee h(s, y) \geq R
				\Big \} \\
			\leq{}&
				\sup \Big \{
					2 L - 2 \delta R
					+ \gamma \big \| 
							\big( (t,x),(s, y) \big) 
					\big \|_{(\R \times H)^2} \colon 
						(t, x), (s, y) \in W_{n_0}
				\Big \} \\
			\leq{}&
				2l \cdot (n_0 +T).
		\end{split}
	\end{equation}
	Combining
	\eqref{eq: splitting S}, \eqref{eq: upper bound W1},
	\eqref{eq: upper bound W2}, \eqref{eq: upper bound W3}, 
	\eqref{eq: upper bound W4}, and \eqref{eq: upper bound W5} then results
	for all 
	$l \in (0, \infty]$,
	$\delta \in (0, \infty)$,
	$R \in \big [\frac{L}{\delta},\infty \big)$,
	$\alpha \in \big[ \frac{4L}{\tilde{r}_{l,R}}, \infty \big )$,
	$\beta \in \big[ \frac{4L}{\tilde{\eps}_{l,R}}, \infty \big )$,
	and all $\gamma \in ( 0, l]$
	in
	\begin{equation}
	\label{eq: S upper bound}
			S_{\alpha, \beta, \gamma, \delta}
		\leq
			S_0 + l (1+ 2 (n_0 +T)).
	\end{equation}		
	In addition \eqref{eq: S lower bound} yields for all 
	$l \in \big(0, \frac {S_0} {2(5+2(n_0+T))} \big)$,
	$\alpha \in \big[ n_0, \infty \big )$,
	$\delta \in \big( 0, \frac {l} {\tilde{R}_{l}} \big]$,
	$\beta \in ( 0, \infty )$,
	and all $\gamma \in ( 0, l ]$ that
	\begin{equation}
		\begin{split}
				l ( 1+ 2 (n_0 +T))
			<
				\tfrac{S_0}{2}
			< 
				S_0 - l (5+2(n_0 +T))
			\leq
				S_{\alpha, \beta, \gamma, \delta}.
		\end{split}
	\end{equation}
	Combining this with \eqref{eq: splitting S},
	\eqref{eq: upper bound W2}, \eqref{eq: upper bound W3}, 
	\eqref{eq: upper bound W4}, \eqref{eq: upper bound W5}, and
	the fact that for all 
	$\alpha$, $\beta$, $\gamma$, $\delta \in (0, \infty)$
	it holds that
	\begin{equation}
		\begin{split}
			S_{\alpha, \beta, \gamma, \delta}
		={}
			&\eta ( \underline{z}^{(\alpha, \beta, \gamma, \delta)}, \delta )
						- \tfrac \alpha 2 \|
								x^{(\alpha, \beta, \gamma, \delta)}_1 - x^{(\alpha, \beta, \gamma, \delta)}_2
							\|_H^2
						- \tfrac \beta 2 |
								t^{(\alpha, \beta, \gamma, \delta)}_1-t^{(\alpha, \beta, \gamma, \delta)}_2
							|^2 \\
						&+ \big \langle 
								p_{\alpha, \beta, \gamma, \delta}, 
								\underline{z}^{(\alpha, \beta, \gamma, \delta)} 
						\big \rangle_{((\R \times H)^2)', (\R \times H)^2}
		\end{split}
	\end{equation}
	shows that for all
	$l \in \big(0, \frac {S_0} {2(5+2(n_0+T))} \big)$,
	$R \in [\frac{2L}{\delta},\infty)$,
	$\alpha \in \big[ \frac{4L}{\tilde{r}_{l,R}} \vee n_0, \infty \big )$,
	$\delta \in \big( 0, \frac {l} {\tilde{R}_{l}} \big]$,
	$\beta \in \big[ \frac{4L}{\tilde{\eps}_{l,R}}, \infty \big )$,
	and all $\gamma \in (0, l]$ it holds that 
	$
			\underline{z}^{(\alpha, \beta, \gamma, \delta)} 
		\in 
			\Omega_{1, \tilde{\eps}_{l,R}, \tilde{r}_{l,R}, R}
	$
	and
	therefore that
	\begin{equation}
	\label{eq: t bigger 0}
			\big( t^{(\alpha, \beta, \gamma, \delta)}_{1} \big)^2 
			\wedge \big( t^{(\alpha, \beta, \gamma, \delta)}_{2} \big)^2
		\geq \tilde{\eps}_{l, R} >0.
	\end{equation}
	Moreover,
	from \eqref{eq: uniformly bounded at p} 
	it follows that for all $R \in (0, \infty)$ and all
	$\delta \in (0,1]$,
	there exists a constant $C_{R, \delta} \in (0, \infty)$ such that 
	for all 
	$(t, x) \in W$, $r \in \R$, $p$, $q \in H'$, 
	$A \in \mathbb{S}_{\mathbb{H}, \mathbb{H}'}$,
	and all $i \in \{1, 2\} $ satisfying that
	$\max \{ h(t,x), \|x\|_H, |r|, \|p\|_{H'}, \| A \|_{L(\mathbb{H},\mathbb{H}')} \} \leq R$,
	$(G_i)^{+}_{\mathbb{H}, \mathbb{X}, \delta, h} \left( (t,x),r,p,A \right) \geq -R$,
	and satisfying that $\|p - q\|_{H'} \leq 4C_{R, \delta}$
	it holds that
	\begin{equation}
	\label{eq: diff small in p}
		\begin{split}
			|	(G_i)^{+}_{\mathbb{H}, \mathbb{X}, \delta, h} \left( (t,x),r,p,A \right)
				-(G_i)^{+}_{\mathbb{H}, \mathbb{X}, \delta, h} 
					\left( (t,x),r,q,A \right) | \leq \tfrac 1R.
		\end{split}
	\end{equation}
	Next denote 
	by 
	$
		(l_\alpha, \delta_\alpha, \beta_\alpha, K_\alpha)_{\alpha \in (0,\infty)}
		\subseteq (0, \infty]^4,
	$
	$
		(\gamma_\alpha)_{\alpha \in (0,\infty)}
		\subseteq (0, \infty],
	$
	by
	$(\underline{z}^{(\alpha)})_{\alpha \in (0,\infty)} \subseteq ([0,T] \times K)^2$, 
	by
	$
			(t^{ (\alpha) }_1)_{\alpha \in (0,\infty)}
	$,
	$
			(t^{ (\alpha) }_2)_{\alpha \in (0,\infty)} 
		\subseteq [0,T]
	$,
	by
	$
			(x^{ (\alpha) }_1, x^{ (\alpha) }_2)_{\alpha \in (0,\infty)}
		\subseteq K^2
	$, 
	by
	$(p_{ \alpha })_{\alpha \in (0,\infty)} \subseteq ((\R \times H)^2)'$,
	$
			(q_{ 1, \alpha }, q_{2, \alpha })_{\alpha \in (0,\infty)}
		\subseteq (\R')^2
	$,
	by
	$
			(p_{ 1, \alpha })_{\alpha \in (0,\infty)}, 
			(p_{ 2, \alpha })_{\alpha \in (0,\infty)} 
		\subseteq H'
	$, and by
	$(S_{ \alpha })_{\alpha \in (0,\infty)} \subseteq \R$
	the sequences
	satisfying for all $\alpha \in (0, \infty)$ 
	that
	$
			l_\alpha 
		= 
			\inf
				\{ l \in [\alpha, \infty) \colon 
					~ \alpha \leq \frac {l}{\tilde{R}_l}
				\} 
	$, 
	that
	$
			\delta_\alpha 
		= 
			\inf (
				\{ \delta \in [\tilde{\delta}_ { \lceil \alpha \rceil}, \infty) \colon 
					~ \alpha \geq \frac {8L}{\tilde{r}_{l_{\delta},2L \delta^{-1}}}
				\} \cup \{ 1\}
			)
	$,
	that
	$
			\beta_\alpha 
		= 
			\frac{8 L}{\tilde{\eps}_{l_{\delta_\alpha},2L \delta_\alpha^{-1}}} 
			\vee \alpha
			\vee \tilde{\beta}_{\lceil \alpha \rceil}
			\vee \1_{(1,\infty)}(\alpha) \,
			\beta_{(\alpha-1) \vee \frac 12}
	$,
	that
	\begin{equation}
	\label{eq: def of K alpha}
		\begin{split}
			K_\alpha 
		=
			&\left( \tfrac {1}{\delta_\alpha} \right)
			\vee \left( l_{\delta_\alpha} \cdot (5+ 4 (n_0 + T)) + \tfrac{1}{\alpha}+L \right)
			\vee \left( 2\alpha \cdot n_0 + \tfrac 1\alpha + l_{\delta_\alpha} \right) \\
			&\vee \left(  4\alpha +2 \beta_\alpha + \tfrac 1\alpha \right)\vee \left( \tfrac {2}{\alpha}+ l_{\delta_\alpha} +1 + \Lambda_{n_0} \right),
		\end{split}
	\end{equation}
	that 
	$
			\gamma_\alpha 
		= 
			l_{\delta_{\alpha}} 
			\wedge C_{K_\alpha, \delta_\alpha} 
	$,
	that
	$
			\underline{z}^{(\alpha)} 
		= 
			\underline{z}^{(\alpha, \beta_\alpha, \gamma_\alpha, \delta_\alpha)}
	$,
	that
	$t^{(\alpha)}_1 = t^{(\alpha, \beta_\alpha, \gamma_\alpha, \delta_\alpha)}_1$,
	that
	$t^{(\alpha)}_2 = t^{(\alpha, \beta_\alpha, \gamma_\alpha, \delta_\alpha)}_2$,
	that
	$x^{(\alpha)}_1 = x^{(\alpha, \beta_\alpha, \gamma_\alpha, \delta_\alpha)}_1$, 
	that
	$x^{(\alpha)}_2 = x^{(\alpha, \beta_\alpha, \gamma_\alpha, \delta_\alpha)}_2$, 
	that
	$p_{ \alpha } = p_{ \alpha, \beta_\alpha, \gamma_\alpha, \delta_\alpha}$, 
	that
	$p_{1, \alpha} = p_{1, \alpha, \beta_\alpha, \gamma_\alpha, \delta_\alpha}$,
	that
	$p_{2, \alpha} = p_{2, \alpha, \beta_\alpha, \gamma_\alpha, \delta_\alpha}$,
	that
	$q_{1, \alpha} = q_{1, \alpha, \beta_\alpha, \gamma_\alpha, \delta_\alpha}$, 
	that
	$q_{2, \alpha} = q_{2, \alpha, \beta_\alpha, \gamma_\alpha, \delta_\alpha}$, and that
	\begin{equation}
	\label{eq: def of Sn}
		\begin{split}
				S_{ \alpha } 
			={}& 
					(\tilde{u}_1)_{\R \times \mathbb{H}, \delta_\alpha, h} ^{-,W}
						( t^{(\alpha)}_{1},x^{(\alpha)}_{1} )
					+(\tilde{u}_2)_{\R \times \mathbb{H}, \delta_\alpha, h} ^{-,W} 
						(t^{(\alpha)}_{2},x^{(\alpha)}_{2} )
						- \tfrac \alpha 2 \| x^{(\alpha)}_{1} - x^{(\alpha)}_{2} \|_H^2 \\
						&- \tfrac {\beta_{ \alpha }} {2} | t^{(\alpha)}_{1} - t^{(\alpha)}_{2}|^2 
					+ \langle 
							p_{\alpha}
							\big( 
								(t^{(\alpha)}_{1},x^{(\alpha)}_{1}), 
								(t^{(\alpha)}_{2},x^{(\alpha)}_{2}) 
							\big) 
						\rangle_{((\R \times H)^2)', (\R \times H)^2}.
		\end{split}
	\end{equation}
	Furthermore denote for every $n \in \N$ by $\Psi_n \colon ((0, T) \times O)^2 \to \R$
	the function satisfying for all $n \in \N$
	and all $\big( (t_1, x_1), (t_2, x_2) \big) \in ((0, T) \times O)^2$ that
	\begin{equation}
		\begin{split}
				\Psi_n \big( (t_1, x_1), (t_2, x_2) \big)
			=
				\tfrac n2 \| x_1 - x_2 \|_H^2
				+ \tfrac {\beta_n}{2} |t_1 - t_2|^2
				+ \langle 
					p_n, ((t_1, x_1), (t_2, x_2)) \, 
				\rangle_{((\R \times H)^2)', (\R \times H)^2}.
		\end{split}
	\end{equation}
	Note that \eqref{eq: S lower bound}, \eqref{eq: S upper bound},
	and the fact that 
	$\lim_{\alpha \to \infty} \delta_\alpha =  \lim_{\delta \to 0} l_\delta =0$
	imply that $\lim_{\alpha \to \infty} S_{\alpha}$ $= S_0$.
  Moreover the set
	$([0,T]\times K)^2$ is closed, $(\R \times \mathbb{H})^2$-bounded, 
	and convex and therefore 
  the Banach-Saks Theorem and Lemma $5.1.4$ in
	Kato \cite{Kato1980} show that
	it is also weakly compact.
	Thus there exists a strictly increasing sequence
	$
    (n_j)_{j \in \N} \subseteq \N 
  $
	and a point 
  $
			\hat{\underline{z}}
		=
			\big ((\hat{t}_1, \hat{x}_1), (\hat{t}_2, \hat{x}_2) \big )
			\in
			([0,T] \times K)^2
  $
	such that
	$
			w-\lim_{ j \to \infty }
			\underline{z}^{(n_j)}
    =
			\hat{\underline{z}} 
  $.
  Clearly, 
  $ 
    \tilde{u}_i(T,x) = - \infty 
  $
  for all 
  $ x \in K $ and all 
  $ i \in \{ 1, 2 \} $
  implies 
  that   
	for all 
  $\alpha \in (0,\infty) $
	it holds that
  $
    t^{ ( \alpha) }_1, t^{ ( \alpha ) }_2
    \in [0, T).
	$
  In addition, observe that if
  $
    \big( \hat{t}_1, \hat{t}_2 \big) 
    \in 
    [0,T]^2 \backslash [0,T-\frac{\mu}{2L})^2
  $,
  then 
	\eqref{eq: u bound on K},
  $\lim_{\alpha \to \infty} S_{\alpha} = S_0$,
	the fact that for all $(t,x) \in [0, T] \times O$
	it holds that $h(t, x) \geq 0$,
	the fact that for all $\alpha \in (0, \infty)$
	it holds that $\|p_\alpha\|_{((\R \times H)^2)'} \leq l_{\delta_\alpha}$
	and the fact that
	$\lim_{\alpha \to \infty} l_{\delta_\alpha} =0$
  imply that
  \begin{align}
	\nonumber
				0
			<
				S_0
		  ={}&
				\lim_{ \alpha \to \infty }
				S_{ \alpha}
			=
				\lim_{ j \to \infty }
				S_{ n_j} \\ \nonumber
			={}& 
				\lim_{ j \to \infty }
				\Big(
					\eta \big( \underline{z}^{ ( n_j ) }, \delta_{ n_j} \big)
					-\tfrac { n_j}{2} \|
						x_1^{( n_j)}
						- x_2^{( n_j)} 
					\|_{H}^2
					- \tfrac {\beta_ {n_j}}{2} | t_1^{( n_j)} - t_2^{( n_j)} |^2
				\\ & \nonumber \qquad \qquad
					+ \langle 
							p_{n_j}, \underline{z}^{ ( n_j ) } 
						\rangle_{((\R \times H)^2)', (\R \times H)^2}
				\Big) \\
			\leq{}&
				\limsup_{ j \to \infty }
				\left(
					\eta\big( \underline{z}^{ ( n_j ) }, 0 \big)
					+ l_{\delta_{n_j}} \cdot \|\underline{z}^{ ( n_j ) }\|_{(\R \times H)^2}
				\right) \\ \nonumber
			\leq{}
				&\left(
					\sum_{ i = 1 }^2
					\left[
						\sup_{
							z \in [0,T] \times K
						}
						u_i( z )
					\right]
				\right)
				- 
				\liminf_{ j \to \infty }
				\left(
					\tfrac{\mu}{T-t_1^{(n_j)}}
					+
					\tfrac{\mu}{T-t_2^{(n_j)}} 
				\right)\\ \nonumber
			\leq{}&
				2L-2L = 0
  \end{align}
  and this contradiction shows that
  $ 
    \big( \hat{t}_1, \hat{t}_2 \big) 
    \in [ 0, T-\frac{\mu}{2L} ]^2
  $. 
	Now denote by 
	$(N_n)_{n  \in \N} \subseteq \N$ the natural numbers satisfying for all $n \in \N$ that 
	\begin{equation}
	\label{eq: def of N}
		\begin{split}
				N_n 
			={} 
				&\max_{i \in \{1, 2\} } \bigg \{ 
					\inf \bigg \{ 
						N \in \N \colon  
							\sup \Big \{
								(G_i)_{\mathbb{H}, \mathbb{X}, \delta_n, h}^{+}
									((\tau, \xi), \nu, \rho, A + 3\alpha I_{\mathbb{H}} \pi^{H}_{H^\perp_m}) \\
								&-(G_i)_{\mathbb{H}, \mathbb{X}, \delta_n, h}^{+}
									((\tau, \xi), \nu, \rho, A) \colon 
								m \in \{N, N+1, \ldots\}, 
								~ \alpha \in (0, K_n/3), 
								~\nu \in \R, \\
								~ &(\tau, \xi) \in W, 
								~ \rho \in H', 						
								~ A \in \mathbb{S}_{\mathbb{H}, \mathbb{H}'}, 
								~ \| A \|_{L(\mathbb{H},\mathbb{H}')} \leq K_n, 
								~ h(\tau, \xi) \leq K_n,\\
								~ &| t_i^{(n)}-\tau | \vee \|x_i^{(n)}- \xi\|_H 
									\vee	|
													(u_i)_{\R \times \mathbb{H}, \delta_n, h}^{-,W}
														\big( t_i^{ (n) }, x_i^{ (n) } \big)  
														-\nu
												| 
									\vee \|n I_{\mathbb{H}} \, (x^{ (n) }_i - x^{ (n) }_{3-i}) - \rho\|_{H'} \\ 
									&\leq 
										m(
											K_n,
											t_i^{ (n) },
											x_i^{ (n) },
											(u_i)_{\R \times \mathbb{H}, \delta_n, h}^{-,W}
												\big( t_i^{ (n) }, x_i^{ (n) } \big),
											n I_{\mathbb{H}} \, (x^{ (n) }_i - x^{ (n) }_{3-i}),
											\delta_n
										), \\
								~ &(G_i)_{\mathbb{H}, \mathbb{X}, \delta_n, h}^{+}
									((\tau, \xi), \nu, \rho, A) \geq -K_n
							\Big \} \leq \frac {1}{K_n} 
					\bigg \}
				\bigg \}.
		\end{split}
	\end{equation}
	Note that \eqref{eq: uniformly bounded in dimension} ensures that
	$(N_n)_{n  \in \N} \subseteq \N$ is well defined.
	Moreover, denote
	by $(B_n)_{n \in \N} \in \mathbb{S}_{\R \times \mathbb{H}, (\R \times \mathbb{H})'}$ 
	the operators satisfying
	for all $n \in \N$ that
	$
			B_n
		=
			\left ( \begin{array}{cc}
				\frac{\beta_n}{n} I_{\R} & 0\\
				0 & I_{\mathbb{H}} 
			\end{array} \right),
	$
	by
	$\mathcal{N}_0$ the set with the property that
	$
			\mathcal{N}_0
		=
			\{n \in \N \colon ~ 
				\underline{z}_n \in ((0,T) \times \overset{\circ}{K})^2,
				~ l_{\delta_{n}} < \infty,
				~ \delta_n < 1,
				~ t^{(n)}_1 \vee t^{(n)}_2 \leq T-\frac{\mu}{3L} - \frac{1}{\beta_n}
			\} \cap
			\{n_j \in \N \colon ~ j \in \N \},
	$
	and by $E_1$, $E_2 \colon \R \times H \to (\R \times H)^2$
	the operator satisfying for all 
	$z \in \R \times H$ that
	$E_1(z) = (z,0)$ and that
	$E_2(z) = (0,z)$.
	Note  that 
	$\lim_{\alpha \to \infty} \delta_\alpha =\lim_{\delta \to 0} l_\delta =0$,
	$\lim_{\alpha \to \infty} S_\alpha =S_0$,
	$\lim_{n \to \infty} \frac {1}{\beta_n} = 0$,
	$\big( \hat{t}_1, \hat{t}_2 \big) \in [ 0, T-\frac{\mu}{2L} )^2$,
	\eqref{eq:sum_tildeu_outside}, and
	\eqref{eq: t bigger 0}
	ensure that $\#\mathcal{N}_0= \infty$.
	Furthermore, observe that for all $n \in \N$ and all
	closed linear subspaces $V \subseteq H$ it holds that
	\begin{equation}
		\begin{split}
				&\pi^{(\R \times H)'}_{(\R \times V)'} \pi^{((\R \times H)^2)'}_1 
					((D^2_{(\R \times \mathbb{H})^2} \Psi_n) (\overline{z}_1, \overline{z}_2)) 
						E_2 \pi^{\R \times H}_{(\R \times V)^\perp} \\
			={}
				&\pi^{(\R \times H)'}_{(\R \times V)'} \pi^{((\R \times H)^2)'}_1 
					\left( \begin{array}{cc}
						n B_n & - n B_n \\
						- n B_n & n B_n \\
					\end{array} \right)
						E_2 \pi^{\R \times H}_{(\R \times V)^\perp} \\
				={}& 
					\pi^{(\R \times H)'}_{(\R \times V)'}
						(- n B_n)
						 \pi^{\R \times H}_{(\R \times V)^\perp} 
				=
					\pi^{(\R \times H)'}_{(\R \times V)'}
						\left( \begin{array}{cc}
							-\beta_n I_{\R} & 0 \\
							0 & - n I_\mathbb{H} \\
						\end{array} \right)
							\pi^{\R \times H}_{(\R \times V)^\perp} \\
				={}
					&\left( \begin{array}{cc}
						-\beta_n I_{\R} \cdot 0 & 0 \\
						0 & - \pi^{H'}_{V'} n I_\mathbb{H} \pi^{H}_{V^\perp}\\
					\end{array} \right)
				= 0,
		\end{split}
	\end{equation}
	and that
	\begin{align}
	\nonumber
				&\left \|
					\pi^{(\R \times H)'}_{((\R \times V)^\perp)'} \pi^{((\R \times H)^2)'}_1 
						((D^2_{(\R \times \mathbb{H})^2} \Psi_n) (\overline{z}_1, \overline{z}_2)) 
							E_1 \pi^{\R \times H}_{(\R \times V)^\perp}
				\right \|_{L(\R \times \mathbb{H}, (\R \times \mathbb{H})')} \\ \nonumber
			={}&
				\left \|
					\pi^{(\R \times H)'}_{((\R \times V)^\perp)'} \pi^{((\R \times H)^2)'}_1
						\left( \begin{array}{cc}
							n B_n & -n B_n \\
							-n B_n & n B_n \\
						\end{array} \right)
							E_1 \pi^{\R \times H}_{(\R \times V)^\perp} 
				\right \|_{L(\R \times \mathbb{H}, (\R \times \mathbb{H})')} \\
			={} 
				&\left \|
					\pi^{(\R \times H)'}_{((\R \times V)^\perp)'}
						n B_n \pi^{\R \times H}_{(\R \times V)^\perp} 
				\right \|_{L(\R \times \mathbb{H}, (\R \times \mathbb{H})')} \\ \nonumber
			={}&
				\left \|
					\pi^{(\R \times H)'}_{((\R \times V)^\perp)'}
						\left( \begin{array}{cc}
							\beta_n I_{\R} & 0 \\
							0 & n I_\mathbb{H} \\
						\end{array} \right)
							\pi^{\R \times H}_{(\R \times V)^\perp} 
				\right \|_{L(\R \times \mathbb{H}, (\R \times \mathbb{H})')} \\ \nonumber
			={}&
				\left \|
					\left( \begin{array}{cc}
						0 \cdot \beta_n I_{\R} \cdot 0 & 0 \\
						0 & \pi^{H'}_{(V^\perp)'} n I_\mathbb{H} \pi^{H}_{V^\perp} \\
					\end{array} \right)
				\right \|_{L(\R \times \mathbb{H}, (\R \times \mathbb{H})')}
			= n.
	\end{align}
	Analogously it follows for all $n \in \N$ and all
	closed linear subspaces $V \subseteq H$ that
	\begin{equation}
		\begin{split}
				&\pi^{(\R \times H)'}_{((\R \times V)^\perp)'} \pi^{((\R \times H)^2)'}_1 
					((D^2_{(\R \times \mathbb{H})^2} \Psi_n) (\overline{z}_1, \overline{z}_2)) 
						E_2 \pi^{\R \times H}_{\R \times V} \\
			={}
				&\pi^{(\R \times H)'}_{((\R \times V)^\perp)'} \pi^{((\R \times H)^2)'}_1 
					((D^2_{(\R \times \mathbb{H})^2} \Psi_n) (\overline{z}_1, \overline{z}_2)) 
						E_1 \pi^{\R \times H}_{\R \times V} \\
			={}&
				\pi^{(\R \times H)'}_{((\R \times V)^\perp)'} \pi^{((\R \times H)^2)'}_2 
					((D^2_{(\R \times \mathbb{H})^2} \Psi_n) (\overline{z}_1, \overline{z}_2)) 
						E_2 \pi^{\R \times H}_{\R \times V}
			= 0
		\end{split}
	\end{equation}
	and that
	\begin{equation}
		\begin{split}
				&\left \|
					\pi^{(\R \times H)'}_{((\R \times V)^\perp)'} \pi^{((\R \times H)^2)'}_2 
						((D^2_{(\R \times \mathbb{H})^2} \Psi_n) (\overline{z}_1, \overline{z}_2)) 
							E_2 \pi^{\R \times H}_{(\R \times V)^\perp}
				\right \|_{L(\R \times \mathbb{H}, (\R \times \mathbb{H})')} \\
			={}&
				\left \|
					\pi^{(\R \times H)'}_{((\R \times V)^\perp)'} \pi^{((\R \times H)^2)'}_1 
						((D^2_{(\R \times \mathbb{H})^2} \Psi_n) (\overline{z}_1, \overline{z}_2)) 
							E_2 \pi^{\R \times H}_{(\R \times V)^\perp}
				\right \|_{L(\R \times \mathbb{H}, (\R \times \mathbb{H})')}
			= n.
		\end{split}
	\end{equation}
  Then
	Corollary~\ref{cor: suitable semijets special case}
		(with
			$
				O \leftarrow (0,T) \times O
			$,
			$
				U \leftarrow (0,T) \times O
			$,
			$
				V \leftarrow \R \times H_{N_n}
			$,
			$ 
				\eps \leftarrow \tfrac{ 1 }{ n } 
			$,
			$
				\lambda \leftarrow 3n
			$,
			$
					u_1 
				\leftarrow 
					(\tilde{u}_1)_{\R \times \mathbb{H}, \delta_n, h}^{-,W}|_{ (0,T) \times O }
			$,
			$
					u_2
				\leftarrow
					(\tilde{u}_2)_{\R \times \mathbb{H}, \delta_n, h}^{-,W}|_{ (0,T) \times O }
			$, 
			$
					\Psi 
				\leftarrow
					\Psi_n
			$,
			and 
			with the local maximum point
			$
					(\overline{z}_1,\overline{z}_2)
				\leftarrow
					\underline{z}^{ (n) }
			$)
  yields that for all
  $
    n \in
    \mathcal{N}_0
  $ 
	there exist linear operators 
	$A^{(n)}_1$, $A^{(n)}_2 \in \mathbb{S}_{\R \times \mathbb{H}, (\R \times \mathbb{H})'}$,       
  such that 
  \begin{align}
			&A^{(n)}_1
		= 
			\pi^{(\R \times H)'}_{(\R \times H_{N_n})'} A^{(n)}_1 
				\pi^{\R \times H}_{\R \times H_{N_n}}, 
		\qquad
			A^{(n)}_2
		= 
			\pi^{(\R \times H)'}_{(\R \times H_{N_n})'} A^{(n)}_2 
				\pi^{\R \times H}_{\R \times H_{N_n}}, \\
  \label{eq:liesinJplus1}
	\begin{split}
  &
    \Big(  
			I_{\R \times \mathbb{H}}
				\big( \beta_n (t^{(n)}_1-t^{(n)}_2), n (x^{(n)}_1-x^{(n)}_2) \big)
			+ (q_{1, n}, p_{1, n}), 
      n A_1^{ (n) } + 3n I_{\R \times \mathbb{H}} \, \pi^{H}_{H^\perp_{N_n}}
    \Big)
	\\ & \qquad\qquad\qquad
    \in 
    \big( 
      \hat{J}^2_{ \R \times \mathbb{H}, + } 
				(\tilde{u}_1)_{\R \times \mathbb{H},\delta_n, h}^{-,W}
    \big)( 
      t^{ (n) }_1, 
      x^{ (n) }_1
    ),
	\end{split}
	\\
	\label{eq:liesinJplus2}
	\begin{split}
    &\Big(
			I_{\R \times \mathbb{H}}
				\big( \beta_n (t^{(n)}_2-t^{(n)}_1), n (x^{(n)}_2-x^{(n)}_1) \big)
			+ (q_{2, n}, p_{2, n}), 
      n A_2^{ (n) } + 3n I_{\R \times \mathbb{H}} \, \pi^{H}_{H^\perp_{N_n}}
    \Big)
	\\ & \qquad\qquad\qquad
    \in 
    \big( 
      \hat{J}^2_{\R \times \mathbb{H}, + }
				(\tilde{u}_2)_{\R \times \mathbb{H}, \delta_n, h}^{-,W}
    \big)( 
      t^{ (n) }_2, 
      x^{ (n) }_2
    ),
	\end{split}
	\end{align}
	and that
	\begin{align}
    &-
    \left[
      n
      +
			n
      \left \|
				\left ( \begin{array}{cc}
					B_n & -B_n \\
					-B_n & B_n
				\end{array} \right)
      \right \|_{L\left( (\R \times \mathbb{H})^2, ((\R \times \mathbb{H})^2)' \right)}
    \right]
    I_{(\R \times \mathbb{H})^2}
  \leq
    \left(
      \begin{array}{cc}
          n A_1^{ (n) }
        & 
          0
        \\
          0
        &
          n A_2^{ (n) } 
      \end{array}
    \right)
  \\ & 
    \leq
    n \,
				\left ( \begin{array}{cc}
					B_n & -B_n \\
					-B_n & B_n
				\end{array} \right)
      +
      \tfrac{ 1 }{ n } 
      \,
      \left[ 
        n 
        \,
				\left ( \begin{array}{cc}
					B_n & -B_n \\
					-B_n & B_n
				\end{array} \right)
      \right]^2 
    .
    \nonumber
  \end{align}
	Combining this with 
  the fact that for all $n \in \mathcal{N}_0$ it holds that
  \begin{equation}
			\left \|
				\left ( \begin{array}{cc}
					B_n & -B_n \\
					-B_n & B_n
				\end{array} \right)
			\right \|_{L((\R \times \mathbb{H})^2, ((\R \times \mathbb{H})^2)')}
		=
			2 \left \|
				B_n
			\right \|_{L(\R \times \mathbb{H}, (\R \times \mathbb{H})')}
		=
			2 \frac{\beta_n}{n}
  \end{equation}
  then implies for all $ n \in \mathcal{N}_0 $ that
  \begin{equation}
  \label{eq:matrix_A_est}
		\begin{split}
			&-\left( 1 + 2 \frac{\beta_n}{n} \right) I_{(\R \times \mathbb{H})^2}
		\leq
			\left(
				\begin{array}{cc}
						A_1^{ (n) }
					&
						0
					\\
						0
					&
						A_2^{  (n) } 
				\end{array}
			\right) 
		\leq
			\left(
				\begin{array}{cc}
					B_n + 2B_n^2  
					& -B_n - 2B_n^2 \\
					-B_n - 2B_n^2  
					& B_n + 2B_n^2
				\end{array}
			\right).
		\end{split}
  \end{equation}
	Moreover, Corollary~\ref{cor:semijets_equivalence2} together 
  with \eqref{eq:liesinJplus1}, \eqref{eq:liesinJplus2}, and with
  and the fact for every
  $ i \in \{ 1, 2 \} $ it holds that 
  $
    \tilde{u}_i|_{(0, T) \times O}
  $
  is a viscosity subsolution of \eqref{eq:parabolic.equation.tilde.ui}
	relative to $(h, \R \times \mathbb{H}, \R \times \mathbb{X})$
  then proves that
	all
	$n \in \mathcal{N}_0$ there exist
	$
		\big( 
			((\tau^{(n)}_i, \xi^{(n)}_i), \nu^{(n)}_i, 
				\sigma^{(n)}_i, \rho^{(n)}_i, \mathfrak{A}^{(n)}_i )
		\big)_{i \in \{1,2\}}
		\subseteq W \times \R \times \R' \times H' 
			\times \mathbb{S}_{\R \times H,(\R \times H)'}
	$ 
	such that
	for all 
	$i \in  \{1,2\}$ it holds that
	\begin{equation}
	\label{eq: distance help points}
		\begin{split}
				& \big |t^{(n)}_i- \tau^{(n)}_i \big | 
				\vee \big \| x^{(n)}_i-\xi^{(n)}_i \big \|_H 
				\vee \big | 
							(\tilde{u}_i)_{\R \times \mathbb{H}, \delta_n, h}^{-,W} 
								(t^{(n)}_i, x^{(n)}_i)
							-(\tilde{u}_i)_{\R \times \mathbb{H}, \delta_n, h}^{-, W} 
								(\tau^{(n)}_i, \xi^{(n)}_i) 
						\big | \\
				&\vee \big | 
							\nu^{(n)}_i
							-(\tilde{u}_i)_{\R \times \mathbb{H}, \delta_n, h}^{-, W} 
								(t^{(n)}_i, x^{(n)}_i) 
						\big | 
				\vee
					\big \| n \mathfrak{A}^{(n)}_i 
								- n A_i^{ (n) } - 3n I_{\R \times \mathbb{H}} \, \pi^{H}_{H^\perp_{N_n}} 
						\big \|_{L(\R \times \mathbb{H}, (\R \times \mathbb{H})')} \\
				&\vee \big \|
						(\sigma^{(n)}_i, \rho^{(n)}_i) 
						- I_{\R \times \mathbb{H}}
							\big( \beta_n (t^{(n)}_i-t^{(n)}_{3-i}), n (x^{(n)}_i-x^{(n)}_{3-i}) \big)
						- (q_{i, n}, p_{i, n}) 
					\big \|_{(\R \times H)'} \\
			\leq
				&\tfrac {1}{\beta_n} \wedge \tfrac{C_{K_n, \delta_n}}{n}
				\wedge 
									\tfrac {1}{2n} m(K_n, t^{(n)}_i, x^{(n)}_i, 
									(u_i)_{\R \times \mathbb{H}, \delta_n, h}^{-, W} (t^{(n)}_i, x^{(n)}_i),
									n I_{\mathbb{H}} (x^{(n)}_i-x^{(n)}_{3-i}), \delta_n)
		\end{split}
	\end{equation}
	and that
	\begin{equation}
	\label{eq: equation close to 0 at the help points}
		\begin{split}
				& I^{-1}_{\R} \sigma^{(n)}_i 
				+ \delta_n \cdot \tfrac{\partial}{\partial t} h(\tau^{(n)}_i,\xi^{(n)}_i)
				- (\tilde{G}_i)_{\mathbb{H}, \mathbb{X}, \delta_n, h}^{+}\!\left(
						(\tau^{(n)}_i, \xi^{(n)}_i),
						\nu^{(n)}_i,
						\rho^{(n)}_i,
						n \pi^{(\R \times H)'}_{2} \mathfrak{A}^{(n)}_i \big |_{H}
				\right)
			\leq \frac 1n.
		\end{split}
	\end{equation}
	Furthermore, it follows from
	\eqref{eq: distance help points},
	\eqref{eq:matrix_A_est},
	the fact that for all $n \in \mathcal{N}_0$ it holds that
	$\pi^{(\R \times H)'}_{2} B_n \big | _H = I_{\mathbb{H}}$,
	the fact that for $n \in \mathcal{N}_0$ it holds that
	$\beta_n \geq n$,
	and from
	\eqref{eq: def of K alpha}
	that for all $n \in \mathcal{N}_0$
	and all $l \in \{ 1, 2 \}$
	it holds that
	\begin{equation}
	\label{eq: A bound}
		\begin{split}
				\|
					n \pi^{(\R \times H)'}_{2} \mathfrak{A}^{(n)}_l \big |_{H} 
				\|_{L(\mathbb{H}, \mathbb{H}')}
			\leq{}&
				\|
					n \pi^{(\R \times H)'}_{2} A^{(n)}_l \big | _H 
					+3n I_{\mathbb{H}} \, \pi^{H}_{H^\perp_{N_n}} 
				\|_{L(\mathbb{H}, \mathbb{H}')}
				+ \tfrac 1n \\
			\leq{}
				&n  \left( \left( 1 + 2 \tfrac{\beta_n}{n} \right) \vee 3 \right )
				+3n+\tfrac 1n 
			={}
				4n +2 \beta_n + \tfrac 1n
			\leq
				K_n.
		\end{split}
	\end{equation}
	and that
	\begin{equation}
	\label{eq: A -3n bound}
		\begin{split}
				&\|
					n \pi^{(\R \times H)'}_{2} \mathfrak{A}^{(n)}_l \big |_{H}
					-3n I_{\mathbb{H}} \, \pi^{H}_{H^\perp_{N_n}}
				\|_{L(\mathbb{H}, \mathbb{H}')}
			\leq{}
				\|
					n \pi^{(\R \times H)'}_{2} A^{(n)}_l \big | _H 
				\|_{L(\mathbb{H}, \mathbb{H}')}
				+ \tfrac 1n \\
			\leq{}
				&n  \left( \left( 1 + 2 \tfrac{\beta_n}{n} \right) \vee 3 \right )
				+\tfrac 1n 
			=
				n +2 \beta_n + \tfrac 1n
			\leq
				K_n.
		\end{split}
	\end{equation}
	Moreover, it follows from
	\eqref{eq: distance help points},
	the lower semicontinuity of $h$
	with respect to the $\| \cdot \|_{\R \times \mathbb{H}}$-norm,
	the fact that for all $x \in K$ it holds that $\|x \|_H \leq n_0$,
	the fact that for all $n \in \mathcal{N}_0$ it holds that
	$\beta_n  \geq n$,
	the fact that 
	for all $n \in \mathcal{N}_0$ it holds that 
	$\|p_{n} \|_{((\R \times H)^2)'} \leq \gamma_n \leq l_{\delta_n}$,
	and from the fact that for all $a, b, c \in H$ it holds that
	$
			\|a-b-c\|_H^2
		= 
			\frac 34 \|a\|_H^2 - 6\|b\|_H^2 -6\|c\|_H^2 
			+\|\frac 12a -2b-2c\|_H^2
			+3\|b-c\|_H^2
		\geq 
			\frac 34 \|a\|_H^2 - 6\|b\|_H^2 -6\|c\|_H^2
	$
	that for all
	$n \in \mathcal{N}_0$ it holds that
	\begin{align}
		\label{eq: S_alpha and S_alpha/2}
				S_{n}
			={}&
				\eta \big( \underline{z}^{(n)}, \delta_n \big)
				- \tfrac {n}{2} \big \|
						x_1^{(n)} - x_2^{(n)} 
					\big \|_{H}^2
				- \tfrac {\beta_n}{2} \big |
						t_1^{(n)} - t_2^{(n)} 
					\big |^2
				+ \big \langle
						p_{n}, \underline{z}^{(n)}
					\big \rangle_{((\R \times H )^2)', (\R \times H)^2} 
		\\ \nonumber \leq{}&	
				(\tilde{u}_1)_{\R \times \mathbb{H}, \delta_n, h}^{-, W}
					\big( t_1^{(n)}, x_1^{(n)} \big)
				+ (\tilde{u}_2)_{\R \times \mathbb{H}, \delta_n, h}^{-, W}
					\big( t_2^{(n)}, x_2^{(n)} \big) 
		\\ \nonumber &
				- \tfrac {n}{2}
						\big \| \xi_1^{(n)}-\xi_2^{(n)}
						- (\xi_1^{(n)} -x_1^{(n)}) 
						- (x_2^{(n)}-\xi_2^{(n)}) \big \|^2_H 
		 \\ \nonumber &
				- \tfrac {\beta_n}{2}\big ( 
						\big | \tau_1^{(n)}-\tau_2^{(n)} \big |
						- \big | t_1^{(n)}-\tau_1^{(n)} \big | 
						- \big | t_2^{(n)}-\tau_2^{(n)} \big |
					\big )^2
				+ \gamma_n \big \| 
						\underline{z}^{(n)} 
					\big \|_{(\R \times H)^2} 
			\\ \nonumber \leq{}&	
				(\tilde{u}_1)_{\R \times \mathbb{H}, \delta_n, h}^{-, W}
					\big( \tau_1^{(n)}, \xi_1^{(n)} \big)
				+ (\tilde{u}_2)_{\R \times \mathbb{H}, \delta_n, h}^{-, W}
					\big( \tau_2^{(n)}, \xi_2^{(n)} \big)
				+ \tfrac {2}{\beta_n} 
				- \tfrac {3n}{8} \big \|
						\xi_1^{(n)} - \xi_2^{(n)} 
					\big \|_{H}^2
			\\ \nonumber &
			  - \tfrac {3\beta_n}{8} \big |
						\tau_1^{(n)} - \tau_2^{(n)} 
					\big |^2
				+\tfrac {6n}{\beta_n^2} 
				+ \tfrac {6}{\beta_n} 
				+ 2l_{\delta_n} (n_0+T) 
			\\ \nonumber \leq{}&	
				(\tilde{u}_1)_{\R \times \mathbb{H}, \delta_n/2, h}^{-, W}
					\big( \tau_1^{(n)}, \xi_1^{(n)} \big)
				+ (\tilde{u}_2)_{\R \times \mathbb{H}, \delta_n/2, h}^{-, W}
					\big( \tau_2^{(n)}, \xi_2^{(n)} \big) 
				- \tfrac {\delta_n}{2} \big(
						h \big( \tau_1^{(n)}, \xi_1^{(n)} \big)
						+ h \big( \tau_2^{(n)}, \xi_2^{(n)} \big) \big) 
			\\ \nonumber &
				- \tfrac {n}{4} \big \|
						\xi_1^{(n)} - \xi_2^{(n)} 
					\big \|_{H}^2 
				- \tfrac {n}{8} \big \|
						\xi_1^{(n)} - \xi_2^{(n)} 
					\big \|_{H}^2 
				- \tfrac {\beta_n}{4} \big |
						\tau_1^{(n)} - \tau_2^{(n)} 
					\big |^2
				- \tfrac {\beta_n}{8} \big |
						\tau_1^{(n)} - \tau_2^{(n)} 
					\big |^2 
				+ \tfrac {14}{n} 
			\\ \nonumber &
				+ 2l_{\delta_n} (n_0+T) 
				+ \tfrac{l_{\delta_n}}{2} \big \| 
							\big ( (\tau_1^{(n)}, \xi_1^{(n)}), (\tau_2^{(n)}, \xi_2^{(n)}) \big) 
						\big \|_{(\R \times H)^2} 
		 \\ \nonumber &
					+ \big \langle
						p_{n/2, \beta_n/2, \gamma_n/2, \delta_n/2}, 
						\big ( (\tau_1^{(n)}, \xi_1^{(n)}) (\tau_2^{(n)}, \xi_2^{(n)}) \big)
					\big \rangle_{((\R \times H)^2)', (\R \times H)^2}
			\\ \nonumber \leq{}&	
				S_{n/2, \beta_n/2, \gamma_n/2, \delta_n/2}
				- \tfrac {\delta_n}{2} \big(
						h \big( \tau_1^{(n)}, \xi_1^{(n)} \big)
						+ h \big( \tau_2^{(n)}, \xi_2^{(n)} \big) \big)
				- \tfrac {n}{8} \big \|
						\xi_1^{(n)} - \xi_2^{(n)}
					\big \|_{H}^2 
		\\ \nonumber &
				- \tfrac {\beta_n}{8} \big |
						\tau_1^{(n)} - \tau_2^{(n)} 
					\big |^2 
				+ \tfrac {14}{n}
				+ 3l_{\delta_n} (n_0+T).
	\end{align}
	Hence we obtain from 
	\eqref{eq: S lower bound}, 
	\eqref{eq: S_alpha and S_alpha/2},
	\eqref{eq: S upper bound},
	and from the fact that for all 
	$n \in \mathcal{N}_0$ it holds that
	$l_n < \infty$,
	$n \geq  \frac{8 L}{\tilde{r}_{l_{\delta_n},2L \delta_n^{-1}}}$,
	$\delta_n \leq \frac {l_{\delta_n}}{\tilde{R}_{l_{\delta_n}}}$,
	and that
	$\beta_n \geq \frac{8 L}{\tilde{\eps}_{l_{\delta_n},2L \delta_n^{-1}}}$
	that for all 
	$n \in \mathcal{N}_0$
	it holds that
	\begin{equation}
	\label{eq: eq for norm convergence}
		\begin{split}
				&S_0 - l_{\delta_n} \cdot (5+ 2 (n_0 + T))
			\leq
				S_{n} \\
			\leq{}&
				S_{n/2, \beta_n/2, \gamma_n/2, \delta_n/2}
				- \tfrac {\delta_n}{2} (
						h(\tau_1^{(n)}, \xi_1^{(n)})
						+ h(\tau_2^{(n)}, \xi_2^{(n)}) )
				- \tfrac {n}{8} \|
						\xi_1^{(n)}
						- \xi_2^{(n)}
					\|_{H}^2 \\
				&- \tfrac {\beta_n}{8} |
						\tau_1^{(n)} - \tau_2^{(n)} 
					|^2
				+ 3l_{\delta_n} (n_0+T) 
				+ \tfrac {14}{n}\\
			\leq{}&
				S_0 + l_{\delta_n} (1+ 5 (n_0 +T)) + \tfrac {14}{n}
				- \tfrac {\delta_n}{2} (
						h(\tau_1^{(n)}, \xi_1^{(n)})
						+ h(\tau_2^{(n)}, \xi_2^{(n)}) )
				- \tfrac {n}{8} \|
						\xi_1^{(n)}
						- \xi_2^{(n)}
					\|_{H}^2 \\
				&- \tfrac {\beta_n}{8} |
						\tau_1^{(n)} - \tau_2^{(n)} 
					|^2.
		\end{split}
	\end{equation} 
	Taking now the limit 
	$\limsup_{\mathcal{N}_0 \ni n \to \infty}$ 
	in \eqref{eq: eq for norm convergence} 
	together with the fact that
	$\lim_{n \to \infty} l_{\delta_n} = 0$ 
	shows then
	\begin{equation}
		\begin{split}
				&\limsup_{\mathcal{N}_0 \ni n \to \infty} \left [
					\tfrac {\delta_n}{2} (h(\tau_1^{(n)}, \xi_1^{(n)}) + h(\tau_2^{(n)}, \xi_2^{(n)}) )
					+\tfrac {n}{8} \|
						\xi_1^{(n)} - \xi_2^{(n)} 
					\|_{H}^2
					+ \tfrac {\beta_n}{8} | \tau_1^{(n)} - \tau_2^{(n)} |^2
				\right ] \\
			\leq{}&
				\limsup_{\mathcal{N}_0 \ni n \to \infty} \left [
					l_{\delta_n} (6+ 7 (n_0 +T)) + \tfrac {14}{n}
				\right ]
			=
				0
		\end{split}
	\end{equation}
	and this yields that
	\begin{equation}
	\label{eq: norm convergence}
		\begin{split}
				&\lim_{\mathcal{N}_0 \ni n \to \infty} \left [
					\delta_n \big( h(\tau_1^{(n)}, \xi_1^{(n)}) + h(\tau_2^{(n)}, \xi_2^{(n)}) \big)
				\right ]
			=
				\lim_{\mathcal{N}_0 \ni n \to \infty} \left [
					n \|
						\xi_1^{(n)}
							- \xi_2^{(n)} 
					\|_{H}^2
				\right ] \\
			={}
				&\lim_{\mathcal{N}_0 \ni n \to \infty} \left [
					\beta_n | \tau_1^{(n)} - \tau_2^{(n)} |^2
				\right ]
			=
				0.
		\end{split}
	\end{equation}
	Furthermore, \eqref{eq: h_t bound},
	\eqref{eq: distance help points},
	\eqref{eq: norm convergence},
	and the fact that for all $n \in \mathcal{N}_0$ it holds that
	$\| \xi^{(n)}_1 \|_H \vee \| \xi^{(n)}_2 \|_H \leq n_0$ imply that
	\begin{align}
		\nonumber
				&\limsup_{\mathcal{N}_1 \ni n \to \infty} 
					\left |
						\sum_{i=1}^2 \left( 
							I_{\R}^{-1} \sigma_i^{(n)} 
							+ \delta_n \cdot \tfrac{\partial}{\partial t} h(\tau^{(n)}_i, \xi^{(n)}_i)
						\right) 
						-\tfrac{2}{n} 
					\right | \\ \nonumber
			\leq{}&
				\limsup_{\mathcal{N}_1 \ni n \to \infty} \Big(
					\left| I_{\R}^{-1} \sigma_1^{(n)} - \beta_n (t_1^{(n)}- t_2^{(n)}) \right |
						+ \left | I_{\R}^{-1} \sigma_2^{(n)} - \beta_n (t_2^{(n)}- t_1^{(n)}) \right |  
		\\ \label{eq: sigma convergence}
					& \qquad\qquad\qquad 
					+  \Lambda_{n_0} \delta_n \cdot  
						\left( 
							h(\tau^{(n)}_1, \xi^{(n)}_1) + h(\tau^{(n)}_2, \xi^{(n)}_2)
						\right)
					+  \tfrac{2}{n}
				\Big) \\ \nonumber
			\leq{}&
				\limsup_{\mathcal{N}_1 \ni n \to \infty} \left(
					\tfrac {2}{\beta_n} + |q_{1,n}| + |q_{2,n}|
					+ \Lambda_{n_0} \delta_n \cdot
						\left( 
							h(\tau^{(n)}_1, \xi^{(n)}_1) + h(\tau^{(n)}_2, \xi^{(n)}_2)
						\right)
					+  \tfrac{2}{n}  
				\right) \\ \nonumber
			\leq{}&
				\limsup_{\mathcal{N}_1 \ni n \to \infty} \left(
					\tfrac {2}{\beta_n} + 2 \gamma_n
					+ \Lambda_{n_0} \delta_n \cdot 
						\left( 
							h(\tau^{(n)}_1, \xi^{(n)}_1) + h(\tau^{(n)}_2, \xi^{(n)}_2)
						\right)
					+  \tfrac{2}{n} 
				\right)
			=
				0.
	\end{align}	
	Denote now by $\mathcal{N}_1 \subseteq \N$ the set satisfying that
	$
			\mathcal{N}_1
		=
			\{  
				n \in \mathcal{N}_0 \colon
				h(\tau_1^{(n)}, \xi_1^{(n)}) \vee h(\tau_2^{(n)}, \xi_2^{(n)}) \leq \frac {1}{\delta_n},
				~|t_1^{(n)} - t_2^{(n)} | \leq \frac {1}{\beta_n},
				~n \geq \tfrac{9L^2}{\mu}
			\}.
	$
	Then \eqref{eq: norm convergence}
	together with $\# \mathcal{N}_0=\infty$ implies that
	$\# \mathcal{N}_1 = \infty$.
	In addition, we get from 
	\eqref{eq: equation close to 0 at the help points},
	\eqref{eq: distance help points},
	\eqref{eq: h_t bound},
	\eqref{eq: def of K alpha},
	the fact that for all $x \in K$ it holds that
	$\| x \|_H \leq n_0$,
	the fact that for all $n \in \mathcal{N}_1$ and $i \in \{1,2\}$ it holds that
	$|I^{-1}_{\R} q_{i,n}|\leq \|p_n\|_{H'} \leq l_{\delta_n}$,
	the fact that for all $n \in \mathcal{N}_1$ it holds that
	$|t_1^{(n)} - t_{2}^{(n)}| \leq \frac {1}{\beta_n} \leq \frac {1}{n}$,
	and from the fact that for all $n \in \mathcal{N}_1$ 
	and all $i \in \{1, 2\}$ it holds that
	$h(\tau^{(n)}_i, \xi^{(n)}_i) \leq \frac {1}{\delta_n}$
	that for all
	$n \in \mathcal{N}_1$ and all $i \in \{1, 2\}$ it holds that
	\begin{equation}
	\label{eq: G bound}
		\begin{split}
					&(G_i)_{\mathbb{H}, \mathbb{X}, \delta_n, h}^{+}\!\left(
						(\tau^{(n)}_i, \xi^{(n)}_i),
						\nu^{(n)}_i + \tfrac{\mu}{T-\tau^{(n)}_i},
						\rho^{(n)}_i,
						n \mathfrak{B}^{(n)}_i
					\right) \\
			\geq{} 
				&-\tfrac 1n + \tfrac {\mu}{(T-\tau_i^{(n)})^2}+ I^{-1}_{\R} \sigma^{(n)}_i 
				+ \delta_n \cdot \tfrac{\partial}{\partial t} h(\tau^{(n)}_i,\xi^{(n)}_i) \\
			\geq{}&
				-\tfrac 1n - \beta_n |t_1^{(n)} - t_{2}^{(n)}| 
				- |I^{-1}_{\R} q_{i,n}| - \tfrac {1}{\beta_n}
				- \Lambda_{n_0} \delta_n \cdot h(\tau^{(n)}_1,\xi^{(n)}_1)
			\geq 
				- \tfrac 2n-l_{\delta_n} -1 - \Lambda_{n_0} \\
			\geq{}
				&- K_n.
		\end{split}
	\end{equation}
	Moreover, it follows from 
	\eqref{eq: u bound on K},
	\eqref{eq: S lower  bound},
	\eqref{eq: def of Sn},
	\eqref{eq: distance help points}, and from
	the fact that $S_0 \geq 0$
	that for all $n \in \mathcal{N}_1$
	it holds that
	\begin{align}
	\nonumber
				&-l_{\delta_n} \cdot (5+ 4 (n_0 + T)) - \tfrac{1}{n} -L \\ \nonumber
			\leq{}&
				S_0 - l_{\delta_n} \cdot (5+ 2 (n_0 + T)) 
				- (\tilde{u}_2)_{\R \times \mathbb{H}, \delta_n, h}^{-, W} (t^{(n)}_2, x^{(n)}_2)
				-2l_{\delta_n}(n_0+T) 
				-\tfrac 1n\\
		\begin{split}
			\leq{}&
				S_n 
				- (\tilde{u}_2)_{\R \times \mathbb{H}, \delta_n, h}^{-, W} (t^{(n)}_2, x^{(n)}_2) 
				+ \tfrac n2 \| x^{(n)}_1 - x^{(n)}_2 \|^2_H
				+ \tfrac {\beta_n}{2} | t^{(n)}_1 - t^{(n)}_2 |^2 \\
				&-\langle 
					p_{n}, 
					\big( 
						(t^{(n)}_{1},x^{(n)}_{1}), 
						(t^{(n)}_{2},x^{(n)}_{2}) 
					\big) 
				\rangle_{((\R \times H)^2)', (\R \times H)^2}
				-\tfrac 1n 
		\end{split} \\ \nonumber
			\leq{}&
				(\tilde{u}_1)_{\R \times \mathbb{H}, \delta_n, h}^{-, W} (t^{(n)}_1, x^{(n)}_1)
				- \tfrac {1}{\beta_n} \\ \nonumber
			\leq
				\nu_1^{(n)}
			\leq{}&
				(\tilde{u}_1)_{\R \times \mathbb{H}, \delta_n, h}^{-, W} (t^{(n)}_1, x^{(n)}_1) 
				+ \tfrac {1}{\beta_n}
			\leq
				L + \tfrac 1n
			\leq
				l_{\delta_n} \cdot (5+ 4 (n_0 + T)) + \tfrac{1}{n}+L.
	\end{align}
	Similar we get for all $n \in \mathcal{N}_1$ that
	\begin{equation}
			-l_{\delta_n} \cdot (5+ 4 (n_0 + T)) - \tfrac{1}{n}-L
		\leq
			\nu_2^{(n)}
		\leq
			l_{\delta_n} \cdot (5+ 4 (n_0 + T)) + \tfrac{1}{n}+L
	\end{equation}
	and this together with 
	the fact that for all $n \in \mathcal{N}_1$
	it holds that $t_1^{(n)} \vee t_2^{(n)} \leq T - \frac{\mu}{3L}$
	and with
	\eqref{eq: def of K alpha},
	proves that for all $n \in \mathcal{N}_1$ it holds that
	\begin{equation}
	\label{eq: nu bound}
			\max_{ \{ i \in 1,2 \}} 
				\Big | \nu^{(n)}_i + \tfrac{\mu}{T-\tau^{(n)}_i} \Big | 
		\leq 
			L + l_n + \tfrac{4(n_0+T+1)}{n} + \tfrac{\mu}{T-(T-\frac{\mu}{3L})}
		=
			4L + l_n + \tfrac{4(n_0+T+1)}{n}
		\leq
			K_n.
	\end{equation}
	In addition, 
	\eqref{eq: distance help points},
	the fact that for all $n \in \mathcal{N}_1$ it holds that
	$\| p_{1, n} \|_{H'} + \| p_{2, n}\|_{ H'} \leq C_{K_n, \delta_n} \wedge l_{\delta_n}$,
	the fact that for all $n \in \mathcal{N}_1$ it holds that
	$\beta_n \geq n$,
	and
	\eqref{eq: def of K alpha}
	show that for all 
	$i \in \{ 1, 2 \}$ and all $n \in \mathcal{N}_1$ it holds that
	\begin{equation}
	\label{eq: p bound}
				\big \|\rho^{(n)}_i \big \|_{H'} 
			\leq 
				n \big \|x^{(n)}_1 - x^{(n)}_2 \big \|_H 
				+ \big \| p_{i, n} \big \|_{H'} + \frac{1}{\beta_n}
			\leq 
				2 n \cdot n_0 + \frac{1}{n} + l_{\delta_n}
			\leq K_n 
	\end{equation}
	and that
	\begin{align}
	\nonumber
				&\big \|
					\rho^{(n)}_i - n I_{\mathbb{H}} (\xi^{(n)}_i - \xi^{(n)}_{3-i} ) 
				\big \|_{H'} \\ \label{eq: p diff bound}
			\leq{}
				&\big \|\rho^{(n)}_i - n I_{\mathbb{H}} (x^{(n)}_i - x^{(n)}_{3-i} ) \big \|_{H'} 
				+ \big \|
						n (x^{(n)}_i - x^{(n)}_{3-i} ) - n (\xi^{(n)}_i - \xi^{(n)}_{3-i} ) 
					\big \|_H  \\ \nonumber
			\leq{}&
				\frac{C_{K_n, \delta_n}}{n} + \big \| p_{i, n} \big \|_{H'} 
				+ n \big \| x^{(n)}_1 - \xi^{(n)}_1 \big \|_H
				+ n \big \| x^{(n)}_2 - \xi^{(n)}_2 \big \|_H
			\leq
				\frac{C_{K_n, \delta_n}}{n} + 3 C_{K_n, \delta_n} \\ \nonumber
			\leq{} 
				&4 C_{K_n, \delta_n}.
	\end{align}
	Then \eqref{eq: diff small in p}
		(with 
		$R \leftarrow K_n$, $\delta \leftarrow \delta_n$, 
		$(t,x) \leftarrow (\tau_i^{(n)},\xi_i^{(n)})$
		$r \leftarrow \nu_i^{(n)}+\tfrac {\mu}{T-\tau_i^{(n)}}$,
		$p \leftarrow \rho_i^{(n)}$, 
		and with
		$A \leftarrow n \pi^{(\R \times H)'}_{2} \mathfrak{A}^{(n)}_i \big |_{H}$)
	together with
	the fact that for all $n \in \mathcal{N}_1$ it holds that
	$
			h(\tau^{(n)}_1, \xi^{(n)}_1) \vee h(\tau^{(n)}_2, \xi^{(n)}_2) 
		\leq 
			\frac {1}{\delta_n} \leq K_n
	$,
	the fact that for all $n \in \mathcal{N}_1$ it holds that
	$
			\|\xi^{(n)}_1 \|_H \vee \|\xi^{(n)}_2 \|_H \leq n_0 +1
	$,
	\eqref{eq: nu bound},
	\eqref{eq: p bound},
	\eqref{eq: A bound},
	\eqref{eq: G bound},
	\eqref{eq: p diff bound},
  \eqref{eq: diff small in p},
	and with the fact that for all $n \in \mathcal{N}_1$ it holds that
	$\frac{1}{K_n} \leq \frac {1}{n}$
	implies that
	for all $n \in \mathcal{N}_1$
	and all $i \in \{1,2\}$ it holds that
	\begin{equation}
	\label{eq: G diff small in p}
		\begin{split}
				&(G_i)_{\mathbb{H}, \mathbb{X}, \delta_n, h}^+\!\left(
					(\tau^{(n)}_i, \xi^{(n)}_i),
					\nu^{(n)}_i + \tfrac {\mu}{T-\tau^{(n)}_i},
					n I_{\mathbb{H}} \, (\xi^{(n)}_i -\xi^{(n)}_i),
					n \pi^{(\R \times H)'}_{2} \mathfrak{A}^{(n)}_i \big |_{H} 
				\right) + \tfrac 1n \\
			\geq{}&
				(G_1)_{\mathbb{H}, \mathbb{X}, \delta_n, h}^+\!\left(
					(\tau^{(n)}_1, \xi^{(n)}_1),
					\nu^{(n)}_1 + \tfrac {\mu}{T-\tau^{(n)}_1},
					\rho^{(n)}_1,
					n \pi^{(\R \times H)'}_{2} \mathfrak{A}^{(n)}_i \big |_{H} 
				\right)
			\geq
				- \tfrac 3n- 1 - \Lambda_{n_0}.
		\end{split}
	\end{equation}
	Note that \eqref{eq: distance help points}
	ensures that for all $n \in \N$ it holds that
	\begin{equation}
	\label{eq: xi x diff small}
		\begin{split}
				&\| 
					n I_{\mathbb{H}} \, (\xi^{(n)}_1 -\xi^{(n)}_2) 
					- n I_{\mathbb{H}} \, (x^{(n)}_1 -x^{(n)}_2)
				\|_{H'}	
			\leq
				n \| \xi^{(n)}_1- x^{(n)}_1 \|_H
				+ n \| \xi^{(n)}_2- x^{(n)}_2 \|_H \\
			\leq{}&
				\min \left \{ 
					m(K_n, t^{(n)}_i, x^{(n)}_i, 
					(u_i)_{\R \times \mathbb{H}, \delta_n, h}^{-, W} (t^{(n)}_i, x^{(n)}_i),
					n I_{\mathbb{H}} (x^{(n)}_i-x^{(n)}_{3-i}), \delta_n)
					\colon i \in \{ 1, 2\}
				\right \}.
		\end{split}
	\end{equation}
	Moreover,
	\eqref{eq: distance help points},
	the fact that 
	for all $n \in \mathcal{N}_1$ and all $i \in \{1,2\}$
	it holds that
	$t^{(n)}_i \vee \tau_i^{(n)} \leq T-\tfrac{3L}{\mu}$,
	and the fact that for all $n \in \mathcal{N}_1$ it holds that
	$n \geq \tfrac{9L^2}{\mu}$
	imply that for all $n \in \mathcal{N}_1$ 
	and all $i \in \{1,2\}$ it holds that
	\begin{equation}
	\label{eq: nu u + mu term diff small}
		\begin{split}
				&\Big | 
					(u_i)_{\R \times \mathbb{H}, \delta_n, h}^{-,W}
						\big( t_i^{ (n) }, x_i^{ (n) } \big)
					+\tfrac{\mu}{T-t^{(n)}_i}
					-\nu^{(n)}_i - \tfrac{\mu}{T-\tau^{(n)}_i} 
				\Big | \\
			\leq{}
				&\Big | 
					(u_i)_{\R \times \mathbb{H}, \delta_n, h}^{-,W}
						\big( t_i^{ (n) }, x_i^{ (n) } \big)
					-\nu^{(n)}_i 
				\Big |
				+\Big |
					\tfrac{\mu \tau^{(n)}_i - \mu t^{(n)}_i}{(T-t^{(n)}_i)(T-\tau^{(n)}_i)}
				\Big | \\
			\leq{}&
				\nicefrac 12 \cdot
					m(K_n, t^{(n)}_i, x^{(n)}_i, 
						(u_i)_{\R \times \mathbb{H}, \delta_n, h}^{-, W} (t^{(n)}_i, x^{(n)}_i),
						n I_{\mathbb{H}} (x^{(n)}_i-x^{(n)}_{3-i}), \delta_n)
				+\Big |
					\mu \, |\tau^{(n)}_i - t^{(n)}_i| \cdot \tfrac{(3L)^2}{\mu^2}
				\Big | \\
			\leq{}&
				m(K_n, t^{(n)}_i, x^{(n)}_i, 
						(u_i)_{\R \times \mathbb{H}, \delta_n, h}^{-, W} (t^{(n)}_i, x^{(n)}_i),
						n I_{\mathbb{H}} (x^{(n)}_i-x^{(n)}_{3-i}), \delta_n).\\
		\end{split}
	\end{equation}
	Thus  
	\eqref{eq: equation close to 0 at the help points}
	and
	\eqref{eq: G diff small in p}
	together with
	\eqref{eq: def of N}
	\eqref{eq: distance help points},
	\eqref{eq: nu u + mu term diff small},
	\eqref{eq: xi x diff small},
	\eqref{eq: A -3n bound},
	the fact that for all $n \in \mathcal{N}_1$ it holds that
	$
			h(\tau^{(n)}_1, \xi^{(n)}_1) \vee h(\tau^{(n)}_2, \xi^{(n)}_2) 
		\leq 
			\frac {1}{\delta_n} \leq K_n
	$,
	\eqref{eq: G bound},
	\eqref{eq: def of K alpha},
	and with
	the fact that for all $n \in \mathcal{N}_1$ it holds that
	$\frac{1}{K_n} \leq \frac {1}{n}$
	proves that
	for all
	$ n \in \mathcal{N}_1 $
	it holds that
	\begin{equation}
	\label{eq: G sum for all n}
		\begin{split}
				&\sum_{i=1}^2 \left( 
					I_{\R}^{-1} \sigma_i^{(n)} 
					+ \delta_n \cdot \tfrac{\partial}{\partial t} h(\tau^{(n)}_i,\xi^{(n)}_i)
					+ \tfrac {\mu}{(T-\tau_i^{(n)})^2} 
				\right)
			 -\frac{2}{n} \\
			\leq{}&
				\sum_{i=1}^2 \left(
					(G_i)_{\mathbb{H}, \mathbb{X}, \delta_n, h}^+ \left( 
						(\tau_i^{(n)}, \xi_i^{(n)}), 
						\nu_i^{(n)} + \tfrac {\mu}{T-\tau_i^{(n)}},
						\rho_i^{(n)},
						n \pi^{(\R \times H)'}_{2} \mathfrak{A}^{(n)}_i \big |_{H} 
					\right)
				\right) \\
			\leq{}&
				\sum_{i=1}^2 \left(
					(G_i)_{\mathbb{H}, \mathbb{X}, \delta_n, h}^+ \left( 
						(\tau_i^{(n)}, \xi_i^{(n)}), 
						\nu_i^{(n)} + \tfrac {\mu}{T-\tau_i^{(n)}},
						n I_{\mathbb{H}} \, (\xi_i^{(n)} - \xi_{3-i}^{(n)}),
						n \pi^{(\R \times H)'}_{2} \mathfrak{A}^{(n)}_i \big |_{H} 
					\right)
				\right)
				+ \frac 2n \\
			\leq{}&
				\sum_{i=1}^2 \bigg(
					(G_i)_{\mathbb{H}, \mathbb{X}, \delta_n, h}^+ \bigg( 
						(\tau_i^{(n)}, \xi_i^{(n)}), 
						\nu_i^{(n)} + \tfrac {\mu}{T-\tau_i^{(n)}},
						n I_{\mathbb{H}} \, (\xi_i^{(n)} - \xi_{3-i}^{(n)}),
						n \pi^{(\R \times H)'}_{2} \mathfrak{A}^{(n)}_i \big |_{H}  
		\\ & \qquad\qquad\qquad\qquad				
						- 3n I_{\mathbb{H}} \, \pi^{H}_{H^\perp_{N_n}}
					\bigg)
				\bigg)
				+ \frac 4n.
		\end{split}
	\end{equation}
	Moreover, 
	\eqref{eq: nu bound} ensures that
	$
			\sup_{n \in \N} \sup_{i \in \{1,2\}} \Big(
				\nu^{ (n) }_i + \tfrac{ \mu }{ T - \tau_i^{ (n) } }
			\Big)
		< 
			\infty
	$ 
	and
	\eqref{eq: distance help points}
	together with 
	$\lim_{\mathcal{N}_1 \ni n \to \infty} S_n = S_0$ 
	and with the fact that 
	$
			\lim_{\mathcal{N}_1 \ni n \to \infty} (\tau_1^{(n)}, \tau_2^{(n)})
		= (\hat {t}_1, \hat {t}_1)
	$
	shows that
	\begin{equation}
		\begin{split}
				&\lim_{\mathcal{N}_1 \ni n \to \infty}
					\sum_{i=1}^2 \left(
						\nu^{ (n) }_i + \tfrac{ \mu }{T - \tau_i^{ (n) } }
					\right) \\
			={}
				&\lim_{\mathcal{N}_1 \ni n \to \infty}
					\sum_{i=1}^2 \left(
						(\tilde{u}_i)_{\R \times \mathbb{H}, \delta_n, h}^{-,W} (t^{(n)}_i, x^{(n)}_i)
						+ \tfrac{ \mu }{  T - \tau_i^{ (n) }  }
						+ \frac {1}{\beta_n}
					\right) 
			=
				S_0
				+ \tfrac{2 \mu }{ T - \hat{t}_1 }
		\end{split}
	\end{equation}
	and thus we have that
	$
		\lim_{n \to \infty}
			\sum_{i=1}^2 \left(
				\nu^{ (n) }_i + \tfrac{ \mu }{ T - \tau_i^{ (n) } }
			\right)
		=
			S_0
			+ \tfrac{2 \mu }{ T - \hat{t}_1 }
		> 0.
	$
	In addition, we get with \eqref{eq: distance help points},
	with \eqref{eq:matrix_A_est},
	and with the fact that for all
	$n \in \mathcal{N}_1$ it holds that 
	$\pi^{(\R \times H)'}_{2} B_n \big |_H = I_{\mathbb{H}}$
	that for all sequences 
	$( z_1^{(n)} )_{n \in \N}$, 
	$( z_2^{(n)} )_{n \in \N} \subseteq H$
	satisfying that $\sup_{n \in \N} \|z_1^{(n)}\|_H + \|z_2^{(n)}\|_H < \infty$
	it holds that
	\begin{align}
	\label{eq: A condition}
				&\limsup_{\mathcal{N}_1 \ni n \to \infty} \left( 
					\sum_{ i = 1 }^2
						\left \langle 
							z_i^{ (n) },  
							\left( 
								\pi^{(\R \times H)'}_{2} \mathfrak{A}^{(n)}_i \big |_{H} 
								- 3I_{\mathbb{H}} \pi^{H}_{H^\perp_{N_n}} 
							\right) z_i^{ (n) } 
						\right \rangle_{H, H'}
				\right) \\ \nonumber
			\leq{}& 
				\limsup_{\mathcal{N}_1 \ni n \to \infty} \left( 
					\sum_{ i = 1 }^2
						\left \langle 
							z_i^{ (n) },  
							\left( \pi^{(\R \times H)'}_{2} A^{ (n) }_i \big |_H \right) z_i^{ (n) } 
						\right \rangle_{H, H'}
					+ \frac{\|z_1^{(n)}\|^2_H + \|z_2^{(n)}\|^2_H}{\beta_n}					
				\right) \\ \nonumber
			\leq{}&
				\limsup_{\mathcal{N}_1 \ni n \to \infty} \Bigg( 
						\bigg \langle 
							(z_1^{ (n) }, z_2^{ (n) }),  
							\left(
								\begin{array}{cc}
									\pi^{(\R \times H)'}_{2} (B_n|_H + 2(B_n|_H)^2)  
									& \pi^{(\R \times H)'}_{2} (-B_n|_H - 2(B_n|_H)^2) \\
									\pi^{(\R \times H)'}_{2} (-B_n|_H - 2(B_n|_H)^2)  
									& \pi^{(\R \times H)'}_{2} (B_n|_H + 2(B_n|_H)^2)
								\end{array}
							\right) 
			\\&\nonumber \qquad\qquad\qquad\qquad
							(z_1^{ (n) }, z_2^{ (n) }) 
						\bigg \rangle_{H^2, (H^2)'}					
				\Bigg)
				 \\
			\leq{}& \nonumber
				\limsup_{\mathcal{N}_1 \ni n \to \infty} \Bigg( 
						\bigg \langle 
							(z_1^{ (n) }, z_2^{ (n) }),  
							\left(
								\begin{array}{cc}
									3 I_{\mathbb{H}}  
									& -3 I_{\mathbb{H}} \\
									-3 I_{\mathbb{H}}  
									& 3 I_{\mathbb{H}}
								\end{array}
							\right) (z_1^{ (n) }, z_2^{ (n) }) 
						\bigg \rangle_{H^2, (H^2)'}					
				\Bigg) \\ \nonumber
			\leq{}
				&\limsup_{\mathcal{N}_1 \ni n \to \infty} \left(
					3 \| z_1^{ (n) } - z_2^{ (n) } \|_H^2
				\right).
	\end{align}
	Finally 
	\eqref{eq: sigma convergence},
	$
			\lim_{\mathcal{N}_1 \ni n \to \infty} (\tau^{(n)}_1, \tau^{(n)}_2) 
		= 
			(\hat{t}_1, \hat{t}_1)
	$,
	\eqref{eq: G sum for all n},
	\eqref{eq:00assumption}
		(with
		$(t^{(n)}_i)_{n\in \N} \leftarrow (\tau_i^{(n)})_{n\in \N} $,
		$(x^{(n)}_i)_{n\in \N}  \leftarrow (\xi_i^{(n)})_{n\in \N} $,
		$(r_i^{(n)})_{n\in \N}  \leftarrow (\nu_i^{(n)} + \frac {\mu}{T-\tau_i^{(n)}})_{n\in \N} $,
		and with
		$
				(A_i^{(n)})_{n\in \N} 
			\leftarrow
				\big(
					\pi^{(\R \times H)'}_{2} \mathfrak{A}^{(n)}_i \big |_{H} 
					- 3I_{\mathbb{H}} \, \pi^{H}_{H^\perp_{N_n}}
				\big)_{n\in \N} 
		$)	
	and \eqref{eq: A condition}
	yield that 
	\begin{align}
				&0
			<
				\frac {2\mu}{(T-\hat{t}_1)^2} 
			=
				\limsup_{\mathcal{N}_1 \ni n \to \infty} \left(
					\sum_{i=1}^2 \left( 
							I_{\R}^{-1} \sigma_i^{(n)} 
							+ \delta_n \cdot \tfrac{\partial}{\partial t} h(\tau^{(n)}_i,\xi^{(n)}_i)
							+ \tfrac {\mu}{(T-\tau_i^{(n)})^2}
						\right)
					 -\frac{2}{n}
				\right )\\ \nonumber
			\leq{}&
				\limsup_{\mathcal{N}_1 \ni n \to \infty} \bigg(
					\sum_{i=1}^2 \bigg(
						(G_i)_{\mathbb{H}, \mathbb{X}, \delta_n, h}^+ \bigg( 
							(\tau_i^{(n)}, \xi_i^{(n)}), 
							\nu_i^{(n)} + \tfrac {\mu}{T-\tau_i^{(n)}},
							n I_{\mathbb{H}} (\xi_i^{(n)} - \xi_{3-i}^{(n)}),
							n \pi^{(\R \times H)'}_{2} \mathfrak{A}^{(n)}_i \big |_{H} 
		\\ \nonumber & \qquad\qquad\qquad
							- 3n I_{\mathbb{H}} \, \pi^{H}_{H^\perp_{N_n}}
						\bigg)
					\bigg) 
					+ \frac 4n
				\bigg) \\ \nonumber
			\leq{} &0.
	\end{align}
  This contradiction implies that 
  $ S_0 \leq 0 $.
  This proves \eqref{eq: assertion with extra mu term} and finishes the proof of 
  Lemma~\ref{l:technical.lemma.uniqueness}.
\end{proof}
As a direct consequence
of Lemma~\ref{lem:sign_changes} and
Lemma~\ref{l:technical.lemma.uniqueness}
we get the following result. 
Its proof is clear and therefore omitted.
\begin{corollary}[A comparison result
for viscosity sub- and supersolutions] 
  \label{c:technical.lemma.uniqueness}
	Assume the setting in Section \ref{ssec: Setting H X with t},
	assume that $O$ is convex,
	and assume that for all $R \in (0, \infty)$ 
	there exists a $\Lambda_R \in (0, \infty)$ such that
	for all $t \in (0,T)$ and all $x \in X \cap O$ with $\|x \|_H \leq R$ it holds that
	\begin{equation}
	\label{eq: h_t bound COR}
		\left | \frac{\partial}{\partial t} h(t,x) \right | \leq \Lambda_R \cdot h(t,x).
	\end{equation}
	Moreover,
  let
  $ 
    u_1 \colon [0,T]\times O \to \R \cup \{ -\infty \}
  $
  be bounded from above  
	on every $\R \times \mathbb{H}$-bounded subset of $[0,T] \times O$,
	let
  $ 
    u_2 \colon [0,T]\times O \to \R \cup \{ \infty \}
  $
  be bounded from below  
	on every $\R \times \mathbb{H}$-bounded subset of $[0,T] \times O$,
  let
  $
    G
    \colon 
    W \times \R \times H' \times \mathbb{S}_{\mathbb{X}, \mathbb{X}'} \to \R
  $
  be a degenerate elliptic function
  and assume
	that
  $ u_1|_{ (0,T) \times O } $ 
  is
  a viscosity subsolution of
   \begin{equation} 
  \label{eq:parabolic.equation.i_COR}
    \tfrac{ \partial }{ \partial t }
    u(t,x) -
    G\big( t, x, u(t,x), (D_{\mathbb{H}} u)(t,x),
      ((D^2_{\mathbb{H}} u)(t,x) |_X) |_X
    \big) = 0
  \end{equation}
	for 
  $ 
    (t,x) \in (0,T) \times O 
  $
	relative to $(h, \R \times \mathbb{H}, \R \times \mathbb{X})$
  and that $ u_2|_{ (0,T) \times O } $
  is a viscosity supersolution of
  \eqref{eq:parabolic.equation.i_COR} relative to 
	$(h, \R \times \mathbb{H}, \R \times \mathbb{X})$.
	Furthermore, assume 
	that for all $ R \in (0,\infty)$ and all
	$\delta \in (0,1]$
	it holds that
	\begin{equation}
	\label{eq: uniformly bounded at p COR+}
		\begin{split}
			&\lim_{\eps \to 0} \Big [
				\sup \{ 
					|G_{\mathbb{H}, \mathbb{X}, \delta, h}^{+}((t,x),r,p,A )
					-G_{\mathbb{H}, \mathbb{X}, \delta, h}^{+}((t,x),r,q,A )| \colon
					~ (t, x) \in W,
					~ r \in \R, 
			\\& \qquad \qquad	
					~ p, q \in H', 
					~ A \in \mathbb{S}_{\mathbb{H}, \mathbb{H}'}, 
					~\|p-q\|_{H'} \leq \eps, 
					~G_{\mathbb{H}, \mathbb{X}, \delta, h}^{+}((t,x),r,p,A) \geq -R,			
		\\& \qquad \qquad
				\max \{ h(t, x), \|x\|_H, |r|, \|p\|_{H'}, \| A \|_{L(\mathbb{H}, \mathbb{H}')} \} \leq R\}
			\Big ] = 0
	\end{split}
	\end{equation}
	and that
	\begin{equation}
	\label{eq: uniformly bounded at p COR-}
		\begin{split}
			&\lim_{\eps \to 0} \Big [
				\sup \{
					|G_{\mathbb{H}, \mathbb{X}, \delta, h}^{-}((t,x),r,p,A)
					-G_{\mathbb{H}, \mathbb{X}, \delta, h}^{-}((t,x),r,q,A)| \colon
					~ (t, x) \in W, 
					~ r \in \R,
		\\& \qquad \qquad
					~ p, q \in H',
					~ A \in \mathbb{S}_{\mathbb{H}, \mathbb{H}'},
					~\|p-q\|_{H'} \leq \eps,  
					~G_{\mathbb{H}, \mathbb{X}, \delta, h}^{-}((t,x),r,p,A)\leq R,
		\\& \qquad \qquad
		\max \{ h(t, x), \|x\|_H, |r|, \|p\|_{H'}, \| A \|_{L(\mathbb{H}, \mathbb{H}')} \} 
						\leq R\}
			\Big ] = 0,
		\end{split}
	\end{equation}
	assume that there exist an increasing sequence of
	finite-dimensional linear subspaces 
	$H_1 \subseteq H_2 \subseteq \ldots \subseteq H$
	and a function 
	$m \colon (0, \infty) \times (0,T) \times O \times \R \times H' \times (0,1] \to (0, \infty)$
	satisfying for all
	$R \in (0,\infty)$, $t \in (0,T)$, $x \in O$, $r \in \R$, $p \in H'$,  
	and all $\delta \in (0,1]$
	that
	$\cup_{N =1}^\infty H_N$ is dense in $H$ with respect to the $\| \cdot \|_H$-norm,
	that
	\begin{equation}
	\label{eq: uniformly bounded in dimension COR+}
		\begin{split}
			\lim_{N \to \infty} \big [
				\sup
				\{ &G_{\mathbb{H}, \mathbb{X}, \delta, h}^{+}
						((\tau, \xi), \nu, \rho, A + \alpha I_{\mathbb{H}} \pi^{H}_{H^\perp_N})
					-G_{\mathbb{H}, \mathbb{X}, \delta, h}^{+}((\tau, \xi), \nu, \rho, A) \colon 
				~ \alpha \in (0, R),\\
				~ &\nu \in \R, 
				~ (\tau, \xi) \in W,
				~ \rho \in H', 
				~ A \in \mathbb{S}_{\mathbb{H},\mathbb{H}'}, 
				~ h(\tau, \xi) \leq R,
				~ \| A \|_{L(\mathbb{H},\mathbb{H}')} \leq R,\\
				~ &| t-\tau | \vee \|x- \xi\|_H \vee |r -\nu| \vee \|p - \rho\|_{H'} 
						\leq m(R,t,x,r,p,\delta), \\
				~ &G_{\mathbb{H}, \mathbb{X}, \delta, h}^{+}((\tau, \xi), \nu, \rho, A) \geq -R \}
			\big ]
				\leq 0,
		\end{split}
	\end{equation}
	and that
	\begin{equation}
	\label{eq: uniformly bounded in dimension COR-}
		\begin{split}
			\lim_{N \to \infty} \big [
				\sup
				\{ &G_{\mathbb{H}, \mathbb{X}, \delta, h}^{-}((\tau, \xi), \nu, \rho, A)
					-G_{\mathbb{H}, \mathbb{X}, \delta, h}^{-}
						((\tau, \xi), \nu, \rho, A- \alpha I_{\mathbb{H}} \pi^{H}_{H^\perp_N}) \colon 
				~ \alpha \in (0, R), \\
				~ &\nu \in \R,
				~ (\tau, \xi) \in W,
				~ \rho \in H',																													
				~ A \in \mathbb{S}_{\mathbb{H},\mathbb{H}'},
				~ h(\tau, \xi) \leq R,
				~ \| A \|_{L(\mathbb{H},\mathbb{H}')} \leq R, \\
				~ &| t-\tau | \vee \|x- \xi\|_H \vee |r -\nu| \vee \|p - \rho\|_{H'} 
						\leq m(R,t,x,r,p,\delta), \\
				~ &G_{\mathbb{H}, \mathbb{X}, \delta, h}^{-}((\tau, \xi), \nu, \rho, A) \leq R \}
			\big ]
				\leq 0,
		\end{split}
	\end{equation}
	and assume that
	there exists a sequence
	$(\tilde{\beta}_n, \tilde{\delta}_n)_{n \in \N} \subseteq (0, \infty)\times (0,1]$
	such that
	$
		\lim_{n \to \infty} \tilde{\delta}_n =0 
	$
	and such that
	for all 
	$  n \in \N $, 
  $ 
    ( 
      (t_n, x_n), r_n, A_n 
    )
	$
	and all
	$
    ( 
      (\hat{t}_n, \hat{x}_n), 
      \hat{r}_n, 
      \hat{A}_n 
    )
    \in
    W \times \R \times \mathbb{S}_{\mathbb{H}, \mathbb{H}'} 
  $ 
  satisfying that
  $
    w-\lim_{ n \to \infty }
    ( t_n, x_n )
    \in [ 0, T ) \times O
  $, 
	that
  $
    \lim_{ n \to \infty }
    \big(
      \sqrt{ n } 
      \| 
         x_n
        -
        \hat{x}_{ n }
      \|_{H}
    \big)
    = 0
  $,
	  that
  $
    \lim_{ n \to \infty }
    \Big(
      \sqrt{ \tilde{\beta}_n } 
      | 
         t_n
        -
        \hat{t}_{ n }
      |
    \Big)
    = 0
  $,
	that
  $
    \lim_{ n \to \infty }
    \big(
      \tilde{\delta}_n ( 
				h(t_{ n }, x_{ n })
        + h(\hat{t}_{ n }, \hat{x}_{ n })
			)
    \big)
    = 0
  $, 
  that
  $
    0 <
    \lim_{ n \to \infty }
    \left( 
      r_n - \hat{r}_n
    \right)
    \leq
    \sup_{ n \in \N } 
    ( 
      | r_n |
      +
      | \hat{r}_n |
    )
    < \infty
  $,
  and that for all
  $
		R \in (0, \infty)
	$
	and all
	$
		(z^{(n)})_{n \in \N}, (\hat{z}^{(n)})_{n \in \N} 
			\subseteq \{z \in H \colon ~ \| z \|_H \leq R\}
	$
   with
	$
		\limsup_{n \to \infty}
    \langle z^{ (n) }, A^{(n)} z^{ (n) } \rangle_{H, H'}
		-
		\langle \hat{z}^{ (n) }, \hat{A}^{(n)} \hat{z}^{ (n) } \rangle_{H, H'}
  \leq 
		\limsup_{n \to \infty}
    3
    \| z^{ (n) } - \hat{z}^{ (n) }
    \|_H^2
  $
	it holds that
  \begin{equation}  
  \label{eq:00assumption_COR}
		\begin{split}
			\limsup_{ 
				n \to \infty 
			}
			\Big[
				&G_{\mathbb{H}, \mathbb{X}, \tilde{\delta}_n, h}^+\big(
					 (t_n, x_n) ,
					 r_n ,
					 n I_{\mathbb{H}} \,
					 (
						 x_n - 
						 \hat{x}_n
					 ) ,
					 n \,
					 A_n
				\big)
		\\ & 
				-G_{\mathbb{H}, \mathbb{X}, \tilde{\delta}_n, h}^-\big(
					(\hat{t}_n, \hat{x}_n) ,
					 \hat{r}_n ,
					 n I_{\mathbb{H}} \,
					 (
						 x_n -
						 \hat{x}_n 
					 ) ,
					 n \,
					 \hat{A}_n
				\big) 
			\Big]
			\leq 0. 
		\end{split}
  \end{equation}
	 In addition, assume that 
	for all $R \in (0,\infty)$ it holds that
	\begin{equation}
	\label{eq: uniformly bounded at 0 COR}
		\begin{split}
				\lim_{r \downarrow 0} \lim_{\eps \downarrow 0}	
				\sup \bigg \{
					u_1(t,x) -  u_2(\hat{t}, \hat{x})
					\colon
					~&(t, x), (\hat{t}, \hat{x}) \in W,
					~h(t, x) \vee h(\hat{t}, \hat{x}) \leq R, \\
					~&\| x - \hat{x}\|_H \leq r,
					~t \vee \hat{t} \leq \eps 
				\bigg \}
			\leq 
				0
		\end{split}
	\end{equation}
  and assume that
  \begin{equation}    
  \label{eq:attains.maximum_COR}
    \lim_{ n \to \infty }
      \sup_{
          (t,x) \in 
          ([0,T] \times O_n^c) \cap W
          }
			u_1(t,x) \vee (- u_2(t,x))
    \leq 0.
  \end{equation}
	Then 
	for all $\tilde{T} \in (0,T)$
	and all $R \in (0,\infty)$ it holds that
  \begin{equation}
		\begin{split}
			\lim_{r \downarrow 0} \lim_{\eps \downarrow 0}	
			\sup \bigg \{
				&u_1(t,x) -  u_2(\hat{t}, \hat{x})
				\colon
				~(t, x), (\hat{t}, \hat{x}) \in W, 
				~t \vee \hat{t} \leq \tilde{T}, \\
				~&h(t, x) \vee h(\hat{t}, \hat{x}) \leq R, 
				~\| x - \hat{x}\|_H \leq r,
				~ |t - \hat{t}| \leq \eps 
			\bigg \} 
    \leq 
			0.
		\end{split}
  \end{equation}
\end{corollary}
Next we state a scaling result for viscosity solutions.
Lemma \ref{l:quotient.subsolution.smooth}  generalizes  Lemma 4.12 in 
Hairer, Hutzenthaler \& Jentzen \cite{HairerHutzenthalerJentzen2015}
to the general notion of viscosity solutions and general
separable Hilbert spaces and will be used in the proof of
Theorem \ref{l:comparison.viscosity.solution} below.
\begin{lemma}[Scaling
of viscosity sub- and supersolutions] 
\label{l:quotient.subsolution.smooth}
	Assume the setting in Section \ref{ssec: Setting H X with t}.
	Moreover,
	let
  $
    V \in
    \C_{\R \times \mathbb{H}}^2( 
      (0,T) \times O, (0,\infty)
    )
  $
	be Lipschitz continuous with respect to the $\R \times \mathbb{H}$-norm
	on every $\R \times \mathbb{H}$-bounded subset of
	$(0,T) \times O$,
  let
  $
    G
    \colon 
    W \times \R \times H' \times \mathbb{S}_{\mathbb{X}, \mathbb{X}'} \to \R
  $
  be a degenerate elliptic function,
  let 
  $ u \colon (0,T) \times O \to \R \cup \{ -\infty, \infty\}$
	be
	bounded from above (below) 
	on $\R \times \mathbb{H}$-bounded subsets of $(0,T) \times O$ and
  be a viscosity subsolution (supersolution)
  of \eqref{eq:parabolic.equation.i_COR} relative to 
	$(Vh, \R \times \mathbb{H}, \R \times \mathbb{X})$
  and let
  $
    \tilde{G} \colon 
    W \times \R \times H' \times \mathbb{S}_{\mathbb{X},\mathbb{X}'}
    \to \R
  $
  be a function satisfying
	for all 
  $
    ((t, x), r, p, A)
    \in
		W
		\times \R \times H' 
    \times \mathbb{S}_{\mathbb{X}, \mathbb{X}'}
  $ that 
  \begin{align}  
	\nonumber
			&\tilde{G}((t,x),r,p,A)
	\\ \nonumber
		={}&
			\tfrac{ 1 }{ V(t,x) }
			\,
			G\Big( (t, x), 
				r \, V(t,x),
				p \, V(t,x)
				+ 
				r \, (D_{\mathbb{H}} V)(t,x),
				A \, V(t,x)
				+
				 p |_X \otimes
					( D_{\mathbb{H}} V )(t,x) |_X
			\\ 
	\label{eq:Gtilde}
			& \qquad\qquad\qquad
				+
				(D_{\mathbb{H}} V)(t,x) |_X
				\otimes p |_X
				+
				r \, 
				(( D^2_{\mathbb{H}} V)(t,x) |_X) |_X
			\Big) 
			-
			r \,
			\tfrac{
				\frac{ \partial }{\partial t } V(t,x)}
				{V(t,x)}.
  \end{align}
  Then 
  $ \tilde{G} $ is
  degenerate elliptic and
  the function 
  $
    \tilde{u} \colon (0,T)
    \times O \to \R
  $
  satisfying
	for all 
  $ (t,x) \in (0,T) \times O $ 
	that
  $
    \tilde{u}(t,x)
    =
    \tfrac{ u(t,x) }{ V(t,x) }
  $
  is
	a viscosity subsolution (supersolution) of
  \begin{equation}   
  \label{eq:parabolic.equation.tilde}
    \tfrac{ \partial }{
      \partial t
	  }  
    \tilde{u}( t, x )
    - \tilde{G}\big(
      (t, x), \tilde{u}(t,x),
      (D_{\mathbb{H}} \tilde{u})(t,x),
      ((D^2_{\mathbb{H}} \tilde{u})(t,x)|_X ) |_X
    \big)
    = 0
  \end{equation}
	for $ (t,x) \in (0,T) \times O $ relative to
	$(h, \R \times \mathbb{H}, \R \times \mathbb{X})$.
\end{lemma}
\begin{proof}[Proof
of 
Lemma~\ref{l:quotient.subsolution.smooth}]
  We proof Lemma~\ref{l:quotient.subsolution.smooth}
  in the case where $ u $ is a viscosity subsolution
  of \eqref{eq:parabolic.equation.i_COR} 
	relative to $(V h, \R \times \mathbb{H}, \R \times \mathbb{X})$.
  The case where $ u $ is a viscosity supersolution 
  of \eqref{eq:parabolic.equation.i_COR} 
	relative to 
	$(V h, \R \times \mathbb{H}, \R \times \mathbb{X})$ follows 
  analogously.
  We thus assume in the following that $ u $ is a viscosity 
  subsolution of \eqref{eq:parabolic.equation.i_COR} 
	relative to 
	$(V h, \R \times \mathbb{H}, \R \times \mathbb{X})$.
	First note that Proposition \ref{prop: h and Vh} 
	shows that $Vh$ fulfills the assumption of 
	Definition \ref{d:viscosity.solution}.
  In the next step
	assume 
  that there exist
  a vector
  $ (t,x) \in (0,T) \times O $,
	a value
	$\delta \in (0, \infty)$,
  and 
  a function
  $
    \phi \in
    \C_{\R \times \mathbb{H}}^2( (0,T) \times O, \R)
  $
  satisfying that
  $ 
    \phi( t, x ) = 
    \tilde{u}_{\R \times \mathbb{H}, \delta, h}^{-, W}(t,x)
  $ 
  and that
  $ 
    \phi \geq \tilde{u}_{\R \times \mathbb{H}, \delta, h}^{-, W}(t,x)
  $.
  Then the function 
  $
    (0,T) \times O
    \ni
    (s,y) \mapsto 
    \phi(s,y) \, V(s,y) \in \R
  $
  is in 
  $
    \C_{\R \times \mathbb{H}}^{ 2 }( (0,T) \times O, \R)
  $ 
  and satisfies that
  $
			\phi(t, x) \, V(t,x) 
		= 
			\tilde{u}_{\R \times \mathbb{H}, \delta, h}^{-, W}(t,x) \, V(t,x) 
		=
			\overline{(\tilde{u} - \delta h)}_{\R \times \mathbb{H}}^{W} (t,x) \, V(t,x)
		=
			\overline{(u - \delta V h)}_{\R \times \mathbb{H}}^{W} (t,x)
  $ 
  and that
  $
			\phi \cdot V 
		\geq 
			\tilde{u}_{\R \times \mathbb{H}, \delta, h}^{-, W} \cdot V 
		= 
			\overline{\tilde{u} -\delta h}_{\R \times \mathbb{H}}^{W} \cdot V 
		=
			\overline{u -\delta V h}_{\R \times \mathbb{H}}^{W}
  $.
  As $u$ is a viscosity 
  subsolution 
  of~\eqref{eq:parabolic.equation.i_COR} relative to 
	$(V h, \R \times \mathbb{H}, \R \times \mathbb{X})$,
  we get that
  \begin{equation}
		\label{eq: u sub relative to Vh I}
		\begin{split}
				&\lim_{\eps \to 0} \inf \Big \{
					\sigma
					+ \delta \tfrac{ \partial }{ \partial t } (Vh) (\tau, \xi)
					-G _{\mathbb{H}, \mathbb{X}, \delta, V h}^+  
						( (\tau, \xi), \nu, \rho, \mathfrak{A}) \colon 
					(\tau, \xi) \in W,
					~\nu, \sigma \in \R,
					~\rho \in H',
			\\ & \qquad \qquad
					~\mathfrak{A} \in \mathbb{S}_{\mathbb{H},\mathbb{H}'},
					~ |\sigma - \tfrac{ \partial }{ \partial t } (V \phi) (t,x) | \leq \eps, 
						d^{-,W}_{\mathbb{H}, \delta, V h, u} \Big(
							( \tau, \xi, \nu, \rho, \mathfrak{A}), \\
				 & \qquad \qquad \qquad
							\Big( 
								t, x, \phi(t,x) \, V(t,x),
								\big( D_{\mathbb{H}}( \phi V) \big)(t,x),
								\big( D^2_{\mathbb{H}}( \phi V) \big)(t,x)
							\Big)
						\Big) \leq \eps
				\Big \}
			\leq 0.
		\end{split}     
  \end{equation}
	Moreover, we have that
	\begin{align}
	\nonumber
			&\lim_{\eps \to 0} \inf \Big \{
				\delta \tfrac{ \partial }{ \partial t } (Vh) (\tau, \xi)
				-G_{\mathbb{H}, \mathbb{X}, \delta, V h}^+ 
						((\tau, \xi), \nu, \rho, \mathfrak{A}) \colon
					~ (\tau, \xi) \in W,
					~\nu \in \R,
					~\rho \in H', \\ \nonumber
					~ & \qquad 
					~\mathfrak{A} \in \mathbb{S}_{\mathbb{H},\mathbb{H}'}, 
						~d^{-,W}_{\mathbb{H}, \delta, V h, u} \Big(
							( \tau, \xi, \nu, \rho, \mathfrak{A}),
					\\ \nonumber & \qquad \qquad 
							\Big( 
								t, x, \phi(t,x) \, V(t,x),
								\big( D_{\mathbb{H}}( \phi V) \big)(t,x),
								\big( D^2_{\mathbb{H}}( \phi V) \big)(t,x)
							\Big)
						\Big) \leq \eps
				\Big \} \\ \nonumber
			={}&
				\lim_{\eps \to 0} \inf \Big \{
					\delta \tfrac{ \partial }{ \partial t } (Vh) (\tau, \xi)
					-G \Big( (\tau, \xi), 
						\nu + \delta V(\tau, \xi) h(\tau, \xi) , 
						\rho
						+ \delta V(\tau, \xi) 
		\\ \nonumber & \qquad
								\cdot E_{\mathbb{X}', \mathbb{H}'}
									((D_{\mathbb{X}} (h|_{(0,T) \times (O \cap X)})) (\tau, \xi)) 
						+ \delta	(D_{\mathbb{H}} V)(\tau, \xi) \cdot h(\tau, \xi), 
						(\mathfrak{A}|_X) |_X 
		\\ \nonumber & \qquad 
							+ \delta ((D^2_{\mathbb{H}} V)(\tau, \xi)|_X) |_X \cdot h (\tau, \xi) 
							+	\delta V(\tau, \xi) 
									\cdot (D^2_{\mathbb{X}} (h|_{(0,T) \times (O \cap X)})) (\tau, \xi) \\ \nonumber
						~ & \qquad 
							+ \delta (D_{\mathbb{X}} (h|_{(0,T) \times (O \cap X)})) (\tau, \xi) 
									\otimes (D_{\mathbb{H}} V)(\tau, \xi)|_X 
							+ \delta (D_{\mathbb{H}} V)(\tau, \xi) |_X
		\\ \nonumber & \qquad 
									\otimes (D_{\mathbb{X}} (h |_{(0,T) \times (O \cap X)})) (\tau, \xi) 
					\Big ) \colon 
						(\tau, \xi) \in W,
						~\nu \in \R,
						~\rho \in H',
						~\mathfrak{A} \in \mathbb{S}_{\mathbb{H},\mathbb{H}'}, \\ \nonumber
						~ & \qquad
							d^{-, W}_{\mathbb{H}, \delta, V h, u} \Big(
								( \tau, \xi, \nu, \rho, \mathfrak{A}),
								\Big( 
									t, x, \phi(t,x) \, V(t,x),
									(D_{\mathbb{H}} V)(t, x) \cdot \phi(t, x)
		\\ \nonumber & \qquad \qquad
										+ V(t, x) \cdot (D_{\mathbb{H}} \phi) (t, x ), 
									(D^2_{\mathbb{H}} V)(t, x) \cdot \phi (t, x)
							+	V(t, x) \cdot (D^2_{\mathbb{H}} \phi) (t, x)
			\\ \nonumber & \qquad \qquad  
							+ (D_{\mathbb{H}} \phi) (t, x) \otimes (D_{\mathbb{H}} V)(t, x)
							+  (D_{\mathbb{H}} V)(t, x) \otimes (D_{\mathbb{H}} \phi) (t, x)
								\Big)
							\Big) \leq \eps
					\Big \} \\ \nonumber
			={}&
				\lim_{\eps \to 0} \inf \Big \{
					\delta \tfrac{ \partial }{ \partial t } (Vh) (\tau, \xi)
					- G\Big( (\tau, \xi), 
						\nu \cdot V(\tau, \xi) + \delta V(\tau, \xi) h(\tau, \xi) , 
						(D_{\mathbb{H}} V)(\tau, \xi) \cdot \nu
				~\\  & \nonumber \qquad 
							+ V(\tau, \xi) \cdot \rho
							+ \delta	(D_{\mathbb{H}} V)(\tau, \xi) \cdot h(\tau, \xi) 
							+ \delta V(\tau, \xi)
									\cdot E_{\mathbb{X}', \mathbb{H}'}
										((D_{\mathbb{X}} (h|_{(0,T) \times (O \cap X)})) (\tau, \xi)), 
				~\\ \nonumber & \qquad 
						((D^2_{\mathbb{H}} V)(\tau, \xi)|_X) |_X \cdot \nu
							+	V(\tau, \xi) \cdot (\mathfrak{A}|_X) |_X 
							+ \rho|_X \otimes (D_{\mathbb{H}} V)(\tau, \xi)|_X
							+ (D_{\mathbb{H}} V)(\tau, \xi) |_X 
					\\\nonumber & \qquad 	
							\otimes \rho|_X 
							+ \delta ((D^2_{\mathbb{H}} V)(\tau, \xi)|_X) |_X \cdot h (\tau, \xi) 
							+ \delta	V(\tau, \xi) 
									\cdot (D^2_{\mathbb{X}} (h|_{(0,T) \times (O \cap X)})) (\tau, \xi) 
						~ \\ \nonumber & \qquad
							+ \delta (D_{\mathbb{X}} (h|_{(0,T) \times (O \cap X)})) (\tau, \xi) 
									\otimes (D_{\mathbb{H}} V)(\tau, \xi)|_X 
							+ \delta (D_{\mathbb{H}} V)(\tau, \xi) |_X
					\\ \nonumber & \qquad
									\otimes (D_{\mathbb{X}} (h|_{(0,T) \times (O \cap X)})) (\tau, \xi)
					\Big ) \colon 
						(\tau, \xi) \in W,
						~\nu \in \R,
						~\rho \in H',
						~\mathfrak{A} \in \mathbb{S}_{\mathbb{H},\mathbb{H}'},\\ 
	\label{eq: u sub relative to Vh II}
						~ & \qquad 
						~	d^{-, W}_{\mathbb{H}, \delta, Vh, u} \Big(
								( \tau, \xi, \nu, \rho, \mathfrak{A}),
								\Big( 
									t, x, \phi(t,x),
									(D_{\mathbb{H}} \phi) (t, x ),
									(D^2_{\mathbb{H}} \phi)(t, x)
								\Big)
							\Big) \leq \eps
				\Big \}.
	\end{align}
	In addition, the Lipschitz continuity of $V$ 
	with respect to the $\| \cdot \|_{\R \times \mathbb{H}}$-norm
	on every $\R \times \mathbb{H}$-bounded subset of
	$(0,T) \times O$ implies that for all $R \in (0, \infty)$
	there exist a $C_R \in (1, \infty)$
	such that for all 
	$z_1$, $z_2 \in \{ z \in (0,T) \times O \colon \|z \|_{\R \times H} \leq R \}$
	it holds that
	$V(z_1) \leq C_R$ and that $|V(z_1) - V(z_2)| \leq C_R \|z_1-z_2\|_{\R \times H}$. 
	This yields that
	for all $R \in (0, \infty)$ and all 
	$z_1$, $z_2 \in \{ z \in (0, T) \times O \colon \|z \|_{\R \times H} \leq R \}$ 
	it holds that
	\begin{align}
	\nonumber
				&|\overline{(u-\delta Vh)}^{W}_{\R \times \mathbb{H}} (z_1) 
				- \overline{(u-\delta Vh)}^{W}_{\R \times \mathbb{H}}(z_2)| \\ 
			={}
				&|\overline{(\tilde{u}-\delta h)}^{W}_{\R \times \mathbb{H}}(z_1) V(z_1) 
				- \overline{(\tilde{u}-\delta h)}^{W}_{\R \times \mathbb{H}}(z_2) V(z_2)| \\ \nonumber
			\leq{}&
				|\overline{(\tilde{u}-\delta h)}^{W}_{\R \times \mathbb{H}}(z_1) 
				- \overline{(\tilde{u}-\delta h)}^{W}_{\R \times \mathbb{H}}(z_2)| V(z_1)
				+ |\overline{(\tilde{u}-\delta h)}^{W}_{\R \times \mathbb{H}}(z_2)|
						\cdot |V(z_1) - V(z_2)| \\ \nonumber
			\leq{}&
				|\overline{(\tilde{u}-\delta h)}^{W}_{\R \times \mathbb{H}} (z_1) 
				- \overline{(\tilde{u}-\delta h)}^{W}_{\R \times \mathbb{H}}(z_2)| \, C_R
				+ |\overline{(\tilde{u}-\delta h)}^{W}_{\R \times \mathbb{H}}(z_2)| 
					\cdot \|z_1 - z_2 \|_{\R \times H} \cdot C_R.
	\end{align}
	Combining this with \eqref{eq:def of d delta-} shows that for all
	$R \in (0, \infty)$,
	$z_1$, $z_2 \in \{ z \in (0,T) \times O \colon \|z \|_{\R \times H} \leq R \}$, 
	$r_1$, $r_2 \in \R$,
	$p_1$, $p_2 \in H'$ and all $A_1$, $A_2 \in \mathbb{S}_{\mathbb{H},\mathbb{H}'}$ 
	it holds that
	\begin{align}
					&d^{-,W}_{\mathbb{H}, \delta, V h, u} \big(
					( z_1, r_1, p_1, A_1),
					( z_2, r_2, p_2, A_2)
				\big) \\ \nonumber
			\leq{}&
				d^{-, W}_{\mathbb{H}, \delta, h, \tilde{u}} \big(
					( z_1, r_1, p_1, A_1),
					( z_2, r_2, p_2, A_2)
				\big) \\ \nonumber
				&\vee \left( 
					|\overline{(\tilde{u}-\delta h)}^{W}_{\R \times \mathbb{H}}(z_1) 
					- \overline{(\tilde{u}-\delta h)}^{W}_{\R \times \mathbb{H}}(z_2)| \, C_R
					+ |\overline{(\tilde{u}-\delta h)}^{W}_{\R \times \mathbb{H}}(z_2)| 
					\cdot \|z_1 - z_2 \|_{\R \times H} \cdot C_R
				\right) \\ \nonumber
			\leq{}&
				d^{-, W}_{\mathbb{H}, \delta, h, \tilde{u}} \big(
					( z_1, r_1, p_1, A_1),
					( z_2, r_2, p_2, A_2)
				\big) 
				\vee \Big( C_R \, d^{-, W}_{\mathbb{H}, \delta, h, \tilde{u}} \big(
					( z_1, r_1, p_1, A_1),
					( z_2, r_2, p_2, A_2)
				\big)
		\\ \nonumber &
				+ C_R \,|\overline{(\tilde{u}-\delta h)}^{W}_{\R \times \mathbb{H}}(z_2)|
					\cdot d^{-, W}_{\mathbb{H}, \delta, h, \tilde{u}} \big(
						( z_1, r_1, p_1, A_1),
						( z_2, r_2, p_2, A_2)
					\big)
				\Big) \\ \nonumber
			={}&
				C_R \, d^{-, W}_{\mathbb{H}, \delta, h, \tilde{u}} \big(
					( z_1, r_1, p_1, A_1),
					( z_2, r_2, p_2, A_2)
				\big)
				 (|\overline{(\tilde{u}-\delta h)}^{W}_{\R \times \mathbb{H}}(z_2)| +1 ).
	\end{align}
	Thus replacing the metric $d^{-,W}_{\mathbb{H}, \delta, V h, u}$ by the metric
	$d^{-,W}_{\mathbb{H}, \delta, h, \tilde{u}}$ in \eqref{eq: u sub relative to Vh II} can
	only decrease the infimum. This results in
	\begin{align}
	\nonumber
				&\lim_{\eps \to 0} \inf \Big \{
					\delta \tfrac{ \partial }{ \partial t } (Vh) (\tau, \xi)
					- G\Big( (\tau, \xi), 
						\nu \cdot V(\tau, \xi) + \delta V(\tau, \xi) h(\tau, \xi) , 
						(D_{\mathbb{H}} V)(\tau, \xi) \cdot \nu
							+ V(\tau, \xi) \cdot \rho
				\\ & \nonumber \qquad \qquad	
							+ \delta	(D_{\mathbb{H}} V)(\tau, \xi) \cdot h(\tau, \xi) 
							+ \delta V(\tau, \xi)
									\cdot E_{\mathbb{X}', \mathbb{H}'}
										((D_{\mathbb{X}} (h|_{(0,T) \times (O \cap X)})) (\tau, \xi)), 
				~ \\&\nonumber \qquad \qquad 
						((D^2_{\mathbb{H}} V)(\tau, \xi)|_X) |_X \cdot \nu
							+	V(\tau, \xi) \cdot (\mathfrak{A}|_X) |_X 
							+ \rho|_X \otimes (D_{\mathbb{H}} V)(\tau, \xi)|_X 
							+ (D_{\mathbb{H}} V)(\tau, \xi) |_X 
					\\& \nonumber \qquad \qquad
							\otimes \rho|_X 
							+ \delta ((D^2_{\mathbb{H}} V)(\tau, \xi)|_X) |_X \cdot h (\tau, \xi) 
							+ \delta	V(\tau, \xi) 
									\cdot (D^2_{\mathbb{X}} (h|_{(0,T) \times (O \cap X)})) (\tau, \xi) 
					~ \\& \nonumber \qquad \qquad 
							+ \delta (D_{\mathbb{X}} (h|_{(0,T) \times (O \cap X)})) (\tau, \xi) 
									\otimes (D_{\mathbb{H}} V)(\tau, \xi)|_X 
							+ \delta (D_{\mathbb{H}} V)(\tau, \xi) |_X
					\\& \nonumber \qquad \qquad
									\otimes (D_{\mathbb{X}} (h|_{(0,T) \times (O \cap X)})) (\tau, \xi)
					\Big ) \colon 
						(\tau, \xi) \in W,
						~\nu \in \R,
						~\rho \in H',
						~\mathfrak{A} \in \mathbb{S}_{\mathbb{H},\mathbb{H}'}, \\
						~ & \nonumber \qquad \qquad
						~	d^{-, W}_{\mathbb{H}, \delta, Vh, u} \Big(
								( \tau, \xi, \nu, \rho, \mathfrak{A}),
								\Big( 
									t, x, \phi(t,x),
									(D_{\mathbb{H}} \phi) (t, x ),
									(D^2_{\mathbb{H}} \phi)(t, x)
								\Big)
							\Big) \leq \eps
				\Big \} \\ \nonumber
			\geq{}&
				\lim_{\eps \to 0} \inf \Big \{
					\delta \tfrac{ \partial }{ \partial t } (Vh) (\tau, \xi)
					- \tilde{G} \Big( (\tau, \xi), 
						\nu + \delta h(\tau, \xi), 
						\rho 
						+ \delta E_{\mathbb{X}', \mathbb{H}'} 
							((D_{\mathbb{X}} (h|_{(0,T) \times (O \cap X)})) (\tau, \xi)),
			\\& \nonumber \qquad \qquad
						(\mathfrak{A}|_X) |_X 
						+ \delta (D^2_{\mathbb{X}} (h |_{(0,T) \times (O \cap X)})) (\tau, \xi)
					\Big ) 
				\cdot V(\tau,\xi) 
				- (\nu + \delta h(\tau, \xi))
					\frac{ \partial }{\partial t } V(\tau,\xi) \colon 
				~ \\& \nonumber \qquad \qquad 
						(\tau, \xi) \in W, 
						~\nu \in \R,
						~\rho \in H', 
						~\mathfrak{A} \in \mathbb{S}_{\mathbb{H},\mathbb{H}'},
			\\& \nonumber \qquad \qquad
						~	d^{-, W}_{\mathbb{H}, \delta, h, \tilde{u}} \Big(
								( \tau, \xi, \nu, \rho, \mathfrak{A}),
								\Big( 
									t, x, \phi(t,x),
									(D_{\mathbb{H}} \phi) (t, x ),
									(D^2_{\mathbb{H}} \phi)(t, x)
								\Big)
							\Big) \leq \eps
				\Big \} \\ \nonumber
			\geq{}&
				\lim_{\eps \to 0} \inf \Big \{
					\delta V(\tau, \xi) (\tfrac{ \partial }{ \partial t } h) (\tau, \xi)
					- \tilde{G}^{+}_{\mathbb{H}, \mathbb{X}, \delta, h} \left( 
						(\tau, \xi), 
						\nu, 
						\rho,
						\mathfrak{A}
					\right ) 
				\cdot V(\tau,\xi) \colon
						(\tau, \xi) \in W, 
						~\nu \in \R,
						~\rho \in H',
				\\& \nonumber \qquad \qquad
						~\mathfrak{A} \in \mathbb{S}_{\mathbb{H},\mathbb{H}'}, 
						~	d^{-, W}_{\mathbb{H}, \delta, h, \tilde{u}} \Big(
								( \tau, \xi, \nu, \rho, \mathfrak{A}),
								\Big( 
									t, x, \phi(t,x),
									(D_{\mathbb{H}} \phi) (t, x ),
									(D^2_{\mathbb{H}} \phi)(t, x)
								\Big)
							\Big) \leq \eps
				\Big \} 
			\\& 
	\label{eq: u sub relative to Vh III}
				- \phi(t,x) 
					\tfrac{ \partial }{\partial t } V(t,x).
	\end{align} 
	Now \eqref{eq: u sub relative to Vh I}, \eqref{eq: u sub relative to Vh II}
	and \eqref{eq: u sub relative to Vh III}
  imply
	\begin{align}
	\label{eq:Gtildephi}
				&0  
			\geq{}
				\lim_{\eps \to 0} \inf \Big \{
					\sigma
						+ \delta \tfrac{ \partial }{ \partial t } (Vh) (\tau, \xi)
						-G_{\mathbb{H}, \mathbb{X}, \delta, V h}^+  
							( (\tau, \xi), \nu, \rho, \mathfrak{A}) \colon 
						(\tau, \xi) \in W,
						~\nu, \sigma \in \R,
						~\rho \in H',
				\\ \nonumber & \qquad \qquad
						~ |\sigma - \tfrac{ \partial }{ \partial t } (V \phi) (t,x) | \leq \eps, 
						~\mathfrak{A} \in \mathbb{S}_{\mathbb{H},\mathbb{H}'},
							~d^{-,W}_{\mathbb{H}, \delta, V h, u} \Big(
								( \tau, \xi, \nu, \rho, \mathfrak{A}),
				\\ \nonumber & \qquad \qquad\qquad
								\Big( 
									t, x, \phi(t,x) \, V(t,x),
									\big( D_{\mathbb{H}}( \phi V) \big)(t,x),
									\big( D^2_{\mathbb{H}}( \phi V) \big)(t,x)
								\Big)
							\Big) \leq \eps
					\Big \} \\ \nonumber
				\geq{}&
					\lim_{\eps \to 0} \inf \Big \{
						\delta V(\tau, \xi) (\tfrac{ \partial }{ \partial t } h) (\tau, \xi)
						- \tilde{G}^+_{\mathbb{H}, \mathbb{X}, \delta, h} \big( 
							(\tau, \xi), 
							\nu, 
							\rho,
							\mathfrak{A}
						\big ) 
					\cdot V(\tau,\xi) \colon
							(\tau, \xi) \in W, 
							~\nu \in \R,
							~\rho \in H',
					\\ \nonumber & \qquad \qquad
							~\mathfrak{A} \in \mathbb{S}_{\mathbb{H},\mathbb{H}'},
							~	d^{+, W}_{\mathbb{H}, \delta, h, \tilde{u}} \Big(
									( \tau, \xi, \nu, \rho, \mathfrak{A}),
									\Big( 
										t, x, \phi(t,x),
										(D_{\mathbb{H}} \phi) (t, x ),
										(D^2_{\mathbb{H}} \phi)(t, x)
									\Big)
								\Big) \leq \eps
					\Big \} 
			\\ \nonumber & 
					+ \tfrac{ \partial }{ \partial t } (V \phi) (t,x)
					- \phi(t,x) 
						\tfrac{ \partial }{\partial t } V(t,x) \\ \nonumber
			={}&
				\lim_{\eps \to 0} \inf \Big \{
					V(\tau, \xi)
						(\sigma + \delta (\tfrac{ \partial }{ \partial t } h) (\tau, \xi))
					- \tilde{G}_{\mathbb{H}, \mathbb{X}, \delta, h} \left( 
						(\tau, \xi), 
						\nu, 
						\rho,
						\mathfrak{A}
					\right ) 
					\cdot V(\tau,\xi) \colon
						(\tau, \xi) \in W, 
						~\nu \in \R,
			\\ \nonumber & \qquad \qquad
						~\rho \in H',
						~\mathfrak{A} \in \mathbb{S}_{\mathbb{H},\mathbb{H}'}, 
						~ |\sigma - (\tfrac{ \partial }{ \partial t } \phi) (t,x) | \leq \eps, 
				~ \\ \nonumber& \qquad \qquad 
						~	d^{-, W}_{\mathbb{H}, \delta, h, \tilde{u}} \Big(
								( \tau, \xi, \nu, \rho, \mathfrak{A}),
								\Big( 
									t, x, \phi(t,x),
									(D_{\mathbb{H}} \phi) (t, x ),
									(D^2_{\mathbb{H}} \phi)(t, x)
								\Big)
							\Big) \leq \eps
					\Big \} \\ \nonumber
				={}&
					V(t, x) \cdot 
				\Big(
						W \times \R \times (\R \times H)' \times \mathbb{S}_{\mathbb{H},\mathbb{H}'}
						\ni ((\tau,\xi),\nu,(\rho_1,\rho_2),\mathfrak{A})
				\\ \nonumber & \qquad \qquad\qquad
						\to 
						I_{\R}^{-1}\rho_1 
						-\tilde{G} ( 
							(\tau, \xi), \nu,
							\rho_2,
							\mathfrak{A}
						) \in \R
					\Big)^+_{\R \times \mathbb{H}, \R \times \mathbb{X}, \delta, h, \tilde{u}}  
				\\ \nonumber & \qquad \qquad\qquad\qquad 
						\Big(
						t,x, 
						\phi(t,x),
						(\tfrac{ \partial }{ \partial t } \phi) (t,x),
						(D_{\mathbb{H}} \phi) (t,x ), 
						(D^2_{\mathbb{H}} \phi)(t,x)
					\Big).
	\end{align}
  This proves 
  inequality~\eqref{eq:Gtildephi}
  for all 
	$ (t,x) \in (0,T) \times O $, 
	$\delta \in (0, \infty)$,
	and all
  $
    \phi \in
    \{
      \psi \in
      \C_{\R \times \mathbb{H}}^{ 2 }( (0,T) \times O, \R )
      \colon
      \psi( t, x ) = 
				\tilde{u}_{\R \times \mathbb{H}, \delta, h}^{-, W}(t,x)
      \text{ and }
      \psi \geq \tilde{u}_{\R \times \mathbb{H}, \delta, h}^{-, W}(t,x)
    \}
  $.
  Therefore,
  $ \tilde{u} $ is a viscosity 
  subsolution
  of~\eqref{eq:parabolic.equation.tilde}
  and the proof of 
  Lemma~\ref{l:quotient.subsolution.smooth}
  is completed.
\end{proof}
Next we state our main theorem of this section.
For the proof we use a suitable classical supersolution
and Lemma \ref{l:quotient.subsolution.smooth} to rescale our equation.
After this transformation we can apply 
Lemma \ref{l:technical.lemma.uniqueness}.
This theorem generalizes
Lemma 4.13 in Hairer, Hutzenthaler \& Jentzen \cite{HairerHutzenthalerJentzen2015} to the general notion
of viscosity solutions and arbitrary separable Hilbert spaces.
\begin{theorem}[A comparison result for viscosity 
sub- and supersolutions]
\label{l:comparison.viscosity.solution}
	Assume the setting in Section \ref{ssec: Setting H X with t},
	assume that $O$ is convex,
	and assume that for all $R \in (0, \infty)$ there exists a 
	$\Lambda_R \in (0, \infty)$ such that
	for all $t \in (0,T)$ and all $x \in X \cap O$ with $\|x \|_H \leq R$ it holds that
	\begin{equation}
	\label{eq: h_t bound LEM}
		\left | \frac{\partial}{\partial t} h(t,x) \right | \leq \Lambda_R \cdot h(t,x).
	\end{equation}
	Moreover,
	let
  $ 
    u_1 \colon [0,T]\times O \to \R \cup \{ -\infty \}
  $
  be bounded from above  
	on every $\R \times \mathbb{H}$-bounded subset of $[0,T] \times O$,
	let
  $ 
    u_2 \colon [0,T]\times O \to \R \cup \{ \infty \}
  $
  be bounded from below  
	on every $\R \times \mathbb{H}$-bounded subset of $[0,T] \times O$,
	let
  $
    V \in
    \C_{\R \times \mathbb{H}}^2( 
      [0,T] \times O, (1,\infty)
    )
  $
	satisfy that
	$D_{\mathbb{H}}  V$ 
	is $\mathbb{H}'$-bounded 
	on every $\R \times \mathbb{H}$-bounded subset
	of $[0,T] \times O$,
	that
	$
		\frac{\partial} {\partial t} V
	$
	is uniformly continuous 
	with respect to the $\| \cdot \|_{\mathbb{H}}$-norm
	on every $\R \times \mathbb{H}$-bounded subset
	of $[0,T] \times O$,
	and that
	$D^2_{\mathbb{H}}  V$ 
	is uniformly continuous 
	with respect to the $\| \cdot \|_{\mathbb{H}}$-norm and the
	$\| \cdot \|_{L(\mathbb{H}, \mathbb{H}')}$-norm
	on every $\R \times \mathbb{H}$-bounded subset
	of $[0,T] \times O$,
  let
  $
    G
    \colon 
    W  \times \R \times H' \times \mathbb{S}_{\mathbb{X}, \mathbb{X}'} \to \R
  $
  be a degenerate elliptic function
  and assume
	that
  $ u_1|_{ (0,T) \times O } $ 
  is
  a viscosity subsolution of
  \begin{equation} 
  \label{eq:parabolic.equation.unbounded.case}
    \tfrac{ \partial }{ \partial t }
    u(t,x) -
    G\big( (t, x), u(t,x), (D_{\mathbb{H}} u)(t,x),
      ((D^2_{\mathbb{H}} u)(t,x)|_X)|_X
    \big) = 0
  \end{equation}
	for 
  $ (t,x) \in (0,T) \times O $ relative to 
	$(Vh, \R \times \mathbb{H}, \R \times \mathbb{X})$,
	that
  $
    u_2|_{ (0,T) \times O }
  $
  is a viscosity supersolution
  of~\eqref{eq:parabolic.equation.unbounded.case} relative to 
	$(V h, \R \times \mathbb{H}, \R \times \mathbb{X})$,  
  and that
  for every $ r \in ( 0, \infty ) $
  it holds that
  $ r V|_{(0,T)\times O}$
  is a classical supersolution
  of~\eqref{eq:parabolic.equation.unbounded.case}
	i.e., suppose that for all $ (t,x) \in W$ it holds that
  \begin{equation}
		(\tfrac{\partial}{\partial t} r V)(t,x)
    -G\big(
      (t, x), rV(t,x), (D_{\mathbb{H}} \, r V)(t, x), 
      ((D_{\mathbb{H}}^2 \, r V)(t, x) |_X) |_X
    \big)
    \geq 0.
  \end{equation}
	Furthermore, assume 
	that for all $ R \in (0,\infty)$ and all
	$\delta \in (0,1]$
	it holds that
	\begin{equation}
	\label{eq: uniformly bounded at p and A+}
		\begin{split}
			\lim_{\eps \to 0} \Big [
				\sup \{
					&|G_{\mathbb{H}, \mathbb{X}, \delta, Vh}^{+}((t,x),r,p,A)
					-G_{\mathbb{H}, \mathbb{X}, \delta, Vh}^{+}((t,x),r,q,B)| \colon
					(t, x) \in W, 
					~ r \in \R, \\
					~ &p, q \in H',
					~A, B \in \mathbb{S}_{\mathbb{H}, \mathbb{H}'}, 
						\max \{ h(t, x), \|x\|_H, |r|, \|p\|_{H'}, \| A \|_{L(\mathbb{H}, \mathbb{H}')} \} 
							\leq R, \\
					~&\|p-q\|_{H'} \leq \eps, 
					~\|A -B \|_{L(\mathbb{H}, \mathbb{H}')} \leq \eps, 
						G_{\mathbb{H}, \mathbb{X}, \delta, Vh}^{+}((t,x),r,p,A) \geq -R\}
			\Big ] = 0 
		\end{split}
	\end{equation}
	and that
	\begin{equation}
	\label{eq: uniformly bounded at p and A-}
		\begin{split}
			\lim_{\eps \to 0} \Big [
				\sup \{
					&|G_{\mathbb{H}, \mathbb{X}, \delta, Vh}^{-}((t,x),r,p,A)
					-G_{\mathbb{H}, \mathbb{X}, \delta, Vh}^{-}((t,x),r,q,B)| \colon
					(t, x) \in W,
					~r \in \R, \\
					~ &p, q \in H', 
					~ A, B \in \mathbb{S}_{\mathbb{H}, \mathbb{H}'}, 
						\max \{ h(t, x), \|x\|_H, |r|, \|p\|_{H'}, \| A \|_{L(\mathbb{H}, \mathbb{H}')} \}
							\leq R, \\
					~&\|p-q\|_{H'} \leq \eps, 
					~\|A -B \|_{L(\mathbb{H}, \mathbb{H}')} \leq \eps, 
					~G_{\mathbb{H}, \mathbb{X}, \delta, Vh}^{-}((t,x),r,p,A) \leq R\}
			\Big ] = 0,
		\end{split}
	\end{equation}
	assume that there exist an increasing sequence of finite-dimensional linear subspaces 
	$H_1 \subseteq H_2 \subseteq \ldots \subseteq H$
	and a function 
	$
		m \colon (0, \infty) \times (0,T) \times O \times \R \times H' \times (0,1] \to (0, \infty)
	$
	satisfying for all
	$R \in (0,\infty)$, $t \in (0,T)$, $x \in O$, $r \in \R$, $p \in H'$,  
	and all $\delta \in (0,1]$
	that
	$\cup_{N =1}^\infty H_N$ is dense in $H$ with respect to the $\| \cdot \|_H$-norm,
	that 
	\begin{equation}
	\label{eq: uniformly bounded in dimension+ LEM}
		\begin{split}
			\lim_{N \to \infty} \big [
				\sup \{ 
					&G_{\mathbb{H}, \mathbb{X}, \delta, Vh}^{+}
						((\tau, \xi), \nu, \rho, A + \alpha I_{\mathbb{H}} \pi^{H}_{H^\perp_N})
					-G_{\mathbb{H}, \mathbb{X}, \delta, Vh}^{+}((\tau, \xi), \nu, \rho, A) \colon 
					\alpha \in (0, R), \\
					~ &\nu \in \R, 
					~ (\tau, \xi) \in W,
					~ \rho \in H', 
					~ A \in \mathbb{S}_{\mathbb{H},\mathbb{H}'}, 
					~ h(\tau, \xi) \leq R, 
					~ \| A \|_{L(\mathbb{H},\mathbb{H}')} \leq R, \\
					~& | t-\tau | \vee \|x- \xi\|_H \vee |r -\nu| \vee \|p - \rho\|_{H'} 
							\leq m(R,t,x,r,p,\delta), \\
					~ &G_{\mathbb{H}, \mathbb{X}, \delta, Vh}^{+}
							((\tau, \xi), \nu, \rho, A) 
						\geq -R 
				\}
			\big ]
				\leq 0, 
		\end{split}
	\end{equation}
	and that
	\begin{equation}
	\label{eq: uniformly bounded in dimension- LEM}
		\begin{split}
			\lim_{N \to \infty} \big [
				\sup \{ 
					&G_{\mathbb{H}, \mathbb{X}, \delta, Vh}^{-}((\tau, \xi), \nu, \rho, A)
					-G_{\mathbb{H}, \mathbb{X}, \delta, Vh}^{-}
						((\tau, \xi), \nu, \rho, A- \alpha I_{\mathbb{H}} \pi^{H}_{H^\perp_N}) \colon 
					\alpha \in (0, R), \\
					~ &\nu \in \R,
					~ (\tau, \xi) \in W,
					~ \rho \in H', 
					~ A \in \mathbb{S}_{\mathbb{H},\mathbb{H}'}, 
					~ h(\tau, \xi) \leq R,
					~ \| A \|_{L(\mathbb{H},\mathbb{H}')} \leq R, \\
					~&| t-\tau | \vee \|x- \xi\|_H \vee |r -\nu| \vee \|p - \rho\|_{H'} 
							\leq m(R,t,x,r,p,\delta), \\
					~ &G_{\mathbb{H}, \mathbb{X}, \delta, Vh}^{-}
							((\tau, \xi), \nu, \rho,\tilde{A}) 
						\leq R 
				\}
			\big ]
				\leq 0,
		\end{split}
	\end{equation}
	and assume that there exists a sequence
	$(\beta_n, \delta_n)_{n \in \N} \subseteq (0, \infty) \times (0,1]$
	such that
	$
    \lim_{ n \to \infty }
			\delta_n
    = 0
  $
	and such that
	for all
  $
    ( (t_n, x_n), r_n, 
			A_n, B_n ),
    ( (\hat{t}_n, \hat{x}_n), \hat{r}_n, 
			\hat{A}_n, \hat{B}_n )
    \in 
		W
		\times \R
    \times \mathbb{S}_{\mathbb{H}, \mathbb{H}'}
    \times \mathbb{S}_{\mathbb{H}, \mathbb{H}'}
  $,
	$ 
			\tilde{A}_n
		\in 
			\mathbb{S}_{\mathbb{H}, \mathbb{H}'}
	$,
	$n \in \N$,
  satisfying
  that
  $
		w-\lim_{n \to \infty}  (t_n,x_n) \in [0,T) \times O
	$,
	that
	$
    \lim_{ n \to \infty }
    \big(
      \sqrt{ n }
      \,
      \|
        x_n - \hat{x}_n 
      \|_{H}
    \big)
    = 0
  $,
	that
	$
    \lim_{ n \to \infty }
    \big(
      \sqrt{ \beta_n }
      \,
      |
				t_n - \hat{t}_n
      |
    \big)
    = 0
  $,
	that
	$
    \lim_{ n \to \infty }
    \big(
      \delta_n h(t_n, x_n) \vee \delta_n h(\hat{t}_n, \hat{x}_n)
    \big)
    = 0
  $,
  that 
  $
    0 < 
    \liminf_{ n \to \infty }
    ( r_n - \hat{r}_n )
    \leq
    \sup_{ n \in \N } 
    ( 
      |r_n| 
      + 
      |\hat{r}_n|
    ) 
    < \infty
  $, 
  that
  $
     A_n
		=
     n 
     (I_\mathbb{H} \left( 
       x_n - \hat{x}_n
     \right)) 
		\otimes ( D_\mathbb{H} V )( t_n , x_n )
     +
     (
       D_\mathbb{H} 
       V
     )( t_n, x_n )
		\otimes
     \left( 
			 n I_\mathbb{H} \left( 
       x_n - \hat{x}_n 
     \right) \right)
     +
     r_n \, 
     ( D^2_\mathbb{H} V)( t_n, x_n ) 
  $,
  that
  $
   \hat{A}_n
		=
     n 
     (I_\mathbb{H} \left( 
       x_n - \hat{x}_n
     \right))
     \otimes ( D_\mathbb{H} V )( \hat{t}_n , \hat{x}_n )
     +
     (
       D_\mathbb{H} 
       V
     )( \hat{t}_n , \hat{x}_n )
     \otimes 
     \left(n I_\mathbb{H} \left( 
       x_n - \hat{x}_n 
     \right) \right)
     +
     \hat{r}_n \, 
     ( D^2_\mathbb{H} V)( \hat{t}_n , \hat{x}_n )
	$,
	that
	$
    \tilde{A}_n
		=
     (r_n-\hat{r}_n) \, 
     ( D^2_\mathbb{H} V)( t_n , x_n )
	$,
	and that for all
	$
		R \in (0, \infty)
	$
	and for all
	$
		(z^{(n)})_{n \in \N}, (\hat{z}^{(n)})_{n \in \N} 
			\subseteq \{z \in H \colon ~ \| z \|_H \leq R\}
	$
	with
	\begin{equation}
		\limsup_{n \to \infty} \Big(
			\langle z^{ (n) }, B_n z^{ (n) } \rangle_{H, H'}
			- \langle \hat{z}^{ (n) }, \hat{B}_n \hat{z}^{ (n) } \rangle_{H, H'}
		\Big)
  \leq 
		\limsup_{n \to \infty}
    3
    \| z^{ (n) } - \hat{z}^{ (n) }
    \|_H^2
  \end{equation}
	it holds that
  \begin{equation}
  \label{eq:comparison.viscosity.solution.assumption}
		\begin{split}
			&\limsup_{ n \to \infty }
			\bigg (
				\tfrac{1}{V( t_n, x_n )} \, \,
				G_{\mathbb{H}, \mathbb{X}, \delta_n, Vh}^{+}\Big( (t_n, x_n), r_n, 
					n I_{\mathbb{H}} \left( x_n - \hat{x}_n \right) V( t_n, x_n )
					+ r_n \, ( D_{\mathbb{H}} V )( t_n, x_n ), 
			\\ &\qquad \qquad
					A_n + n B_n \, V( t_n, x_n )
				\Big)
			-
				\tfrac{1}{V( \hat{t}_n, \hat{x}_n )} \, \,
				G_{\mathbb{H}, \mathbb{X}, \delta_n, Vh}^{-}\Big( (\hat{t}_n, \hat{x}_n), \hat{r}_n, 
					n I_{\mathbb{H}} \left( x_n - \hat{x}_n \right) V( \hat{t}_n, \hat{x}_n )
			\\ & \qquad \qquad
					+ \hat{r}_n \, ( D_{\mathbb{H}} V )( \hat{t}_n, \hat{x}_n ), 
					\hat{A}_n + n \hat{B}_n \, V( \hat{t}_n, \hat{x}_n ) \Big)
			\bigg ) \\
		\leq{}& 
			\limsup_{ n \to \infty } \left (
				\frac{
					G( 
						(t_n, x_n), 
						r_n- \hat{r}_n, 
						(r_n- \hat{r}_n) ( D_{\mathbb{H}} V )(t_n, x_n), 
						(\tilde{A}_n |_X) |_X 
					)
				}{
					V( t_n, x_n )
				}
			\right ).
		\end{split}
  \end{equation}  
	In addition, assume that for all $R,\tilde{R} \in (0, \infty)$ it holds that
	\begin{equation}
	\label{eq: uniformly bounded at 0 with V}
		\begin{split}
			\lim_{r \downarrow 0} \lim_{\eps \downarrow 0}	
			\sup \Bigg \{
				&u_1(t,x) -  u_2(\hat{t}, \hat{x})
				\colon
				~(t, x), (\hat{t}, \hat{x}) \in W,
				~h(t, x) \vee h(\hat{t}, \hat{x}) \leq R, 
				~\|x \|_H \leq \tilde{R}, \\
				~&\| x - \hat{x}\|_H \leq r,
				~t \vee \hat{t} \leq \eps 
			\Bigg \}
    \leq 
			0
		\end{split}
	\end{equation}
  and assume that
  \begin{equation}  
  \label{eq:attains.maximum.difference2.comparison}
    \lim_{ n \to \infty }
    \Bigg[
      \sup_{
        \substack{
          (t,x) \in ((0,T) \times O_n^c) \cap W       }
      } 
			\left(
        \frac{ u_1(t, x)}{ V(t,x) } \vee \frac{-u_2(t,x) }{ V(t,x) }
		\right)
    \Bigg]
    \leq 0 .
  \end{equation}
 Then 
	for all $\tilde{T} \in (0,T)$ and all
	$R,\tilde{R} \in (0, \infty)$ it holds that
  \begin{equation}
		\begin{split}
			\lim_{r \downarrow 0} \lim_{\eps \downarrow 0}
				\sup \Bigg \{
					u_1(t,x) - u_2(\hat{t} , \hat{x}) \colon 
					&~(t, x), (\hat{t}, \hat{x}) \in W,
					~h(t, x) \vee h(\hat{t}, \hat{x}) \leq R,
					~\|x\|_H \leq \tilde{R}, \\
					~ &t, \hat{t} \leq \tilde{T}, 
					~\|x- \hat{x}\|_H \leq r, 
					~|t-\hat{t}| \leq \eps 
				\Bigg \}
			\leq
				0.
		\end{split}
  \end{equation}
\end{theorem}
\begin{proof}[Proof
of Theorem~\ref{l:comparison.viscosity.solution}]
	First note that Proposition \ref{prop: h and Vh} 
	shows that $(Vh)|_{(0,T) \times O}$ fulfills the assumption of 
	Definition \ref{d:viscosity.solution}.
	Now 
  let 
  $
    \tilde{u}_1,
    \tilde{u}_2
    \colon [0,T] \times O \to \R \cup \{-\infty, \infty\}
  $
  and let
  $
    \tilde{G} \colon
    W
    \times \R \times H'
    \times \mathbb{S}_{\mathbb{X},\mathbb{X}'}
    \to \R
  $
  be the functions with the property that
	for all 
  $
    (t,x) \in [0,T] \times O
  $
	it holds that
  $
    \tilde{u}_1(t,x) = 
    \frac{ u_1(t,x) }{ V(t,x) }
  $
  and that 
  $
    \tilde{u}_2(t,x) = 
    \frac{ u_2(t,x) }{ V(t,x) }
  $ 
  and
	for all 
  $
    ((t, x), r, p, A)
    \in
    W \times \R \times H' 
    \times \mathbb{S}_{\mathbb{X}, \mathbb{X}'}
  $ 
	it holds that
  \begin{equation}  
  \label{eq:Gtilde1}
		\begin{split}
    &\tilde{G}( t, x, r, p, A) \\
    ={}
    &\tfrac{ 1 }{ V(t,x) }
    \,
    G\Big( t, x, 
      r \, V(t,x),
      p \, V(t,x)
      + 
      r \, (D_{\mathbb{H}} V)(t,x),
      A \, V(t,x)
      + p|_X \otimes ( D_{\mathbb{H}} V )(t,x) |_X
				\\ & \qquad \qquad 
      + (D_{\mathbb{H}} V)(t,x) |_X \otimes p |_X
      + r \, (( D^2_{\mathbb{H}} V)(t,x) |_X) |_X
    \Big) 
    -
    r \,
    \tfrac{
      \frac{ \partial }{
        \partial t 
      }
      V(t,x)
    }{
      V(t,x)
    }.
		\end{split}
  \end{equation}
  Lemma~\ref{l:quotient.subsolution.smooth}
  then
  ensures that 
  $ \tilde{G} $
  is degenerate elliptic,
  that 
  $
    \tilde{u}_1|_{ (0,T) \times O }
  $ 
  is
	a viscosity subsolution of
  \begin{equation} 
  \label{eq:parabolic.equation.unbounded.case.i_PROOF}
    \tfrac{ \partial }{ \partial t }
    u(t,x) -
    \tilde{G}\big( 
      t, x, u(t,x), (D_{\mathbb{H}} u)( t, x ),
      ((D^2_{\mathbb{H}} \, u)(t,x) |_X ) |_X
    \big) = 0
  \end{equation}
	for 
  $
    (t,x) \in (0,T) \times O 
  $ 
	relative to 
	$(h, \R \times \mathbb{H}, \R \times \mathbb{X})$
  and that
  $
    \tilde{u}_2|_{ (0,T) \times O }
  $  
  is viscosity supersolution of
  \eqref{eq:parabolic.equation.unbounded.case.i_PROOF} 
	relative to $(h, \R \times \mathbb{H}, \R \times \mathbb{X})$.
  Below we will finish this proof 
  by an application 
  of Corollary~\ref{c:technical.lemma.uniqueness}
  with $ \tilde{u}_1 $, $ \tilde{u}_2 $ 
  and $ \tilde{G} $.
  For this we now check the assumptions
  of Corollary~\ref{c:technical.lemma.uniqueness}.
  First, observe that 
	the assumption that
	$D^2_{\mathbb{H}}  (V|_{ (0,T) \times O })$
	is uniformly continuous with respect to the 
	$\| \cdot \|_{\mathbb{H}}$ and the
	$\| \cdot \|_{L(\mathbb{H},\mathbb{H}')}$-norm
	on $\R \times \mathbb{H}$-bounded subsets 
	of $(0,T) \times O$,
	that $\frac {\partial} {\partial t} V$
	is uniformly continuous with respect to the 
	$\| \cdot \|_{\mathbb{H}}$-norm
	on $\R \times \mathbb{H}$-bounded subsets 
	of $(0,T) \times O$, and that
	$D_{\mathbb{H}}  (V|_{ (0,T) \times O })$
	is $\mathbb{H}'$-bounded 
	on $\R \times \mathbb{H}$-bounded subsets 
	of $(0,T) \times O$
	imply
	that $D^2_{\mathbb{H}}  (V|_{ (0,T) \times O })$
	is $\mathbb{L}(\mathbb{H},\mathbb{H}')$-bounded
	on $\R \times \mathbb{H}$-bounded subsets of $(0,T) \times O$,
	that $\frac {\partial} {\partial t} V$
	is bounded on $\R \times \mathbb{H}$-bounded subsets of $(0,T) \times O$,
	and that $V$
	is Lipschitz continuous with respect to the $\| \cdot \|_{\R \times \mathbb{H}}$-norm
	and therefore also bounded
	on $\R \times \mathbb{H}$-bounded subsets
	of $(0,T) \times O$. 
	This together with the assumption that 
	$u_1 \vee -u_2$ is bounded from above on 
	$\R \times \mathbb{H}$-bounded subsets of $(0,T) \times O$
	shows that
	for all $R \in (0, \infty)$ there exists a $C_R \in (1, \infty)$
	such that for all $(t,x)$, $(\tau, \xi) \in (0,T) \times O$ with 
	$\|x\|_H \vee \|\xi\|_H \leq R$
	it holds that 
	\begin{equation}
	\label{eq: V bound}
		\begin{split}
				&|V(t,x)| 
				\vee \|(D_{\mathbb{H}} V) (t,x) \|_{H'} 
				\vee \| (D^2_{\mathbb{H}} V)(t,x) \|_{L(\mathbb{H},\mathbb{H}')}
				\vee | \tfrac {\partial} {\partial t} V(t,x) |
				\vee u_1(t,x)
				\vee -u_2(t,x)
			\leq C_R, \\ 
				&|V(t,x)-V(\tau, \xi)| 
			\leq 
				C_R \|(t,x) - (\tau, \xi)\|_{\R \times H} 
		\end{split}
	\end{equation} 
	Combining this with \eqref{eq: h_t bound LEM} yields then that
	for all $R \in (0, \infty)$
	and all $(t,x) \in (0, T) \times (X \cap O)$ with $\|x \|_H \leq R$
	it holds that
	\begin{equation}
			|\tfrac{\partial}{\partial t} (h V) (t,x)|
		\leq
			|(\tfrac{\partial}{\partial t} h)(t,x) \cdot V (t,x)|
			+ |h(t,x) \cdot (\tfrac{\partial}{\partial t} V)(t,x)|
		\leq
			(\Lambda_R \cdot C_R + C_R) h(t,x).
	\end{equation}
	Next, note
	\eqref{eq:attains.maximum.difference2.comparison} ensures that
	\eqref{eq:attains.maximum_COR} is fulfilled.
	Moreover, \eqref{eq: uniformly bounded at 0 with V}
	and the assumption $V > 1$
	imply that for  all $n \in \{2, 3, 4, \ldots \}$
	and for all $R \in (0, \infty)$
	it holds that 
	\begin{align}
	\nonumber
			&\lim_{r \downarrow 0} \lim_{\eps \downarrow 0}
				\sup \bigg \{
					\tilde{u}_1(t,x) - \tilde{u}_2(\hat{t} , \hat{x}) \colon 
					~(t, x), (\hat{t}, \hat{x}) \in W,
					~h(t, x) \vee h(\hat{t}, \hat{x}) \leq R,
			\\  \nonumber & \qquad\qquad\qquad
					~\|x- \hat{x}\|_H \leq r,
					~t \vee \hat{t} \leq \eps 
				\bigg \} \\ \nonumber
			={}&
				\lim_{r \downarrow 0} \lim_{\eps \downarrow 0}
				\sup \bigg \{
					\frac {u_1(t,x)}{V(t,x)} 
					- \frac {u_2(\hat{t} , \hat{x})}{V(\hat{t} , \hat{x})} \colon 
					~(t, x), (\hat{t}, \hat{x}) \in W,
					~h(t, x) \vee h(\hat{t}, \hat{x}) \leq R,
		\\ \nonumber & \qquad\qquad\qquad
					~\|x- \hat{x}\|_H \leq r,
					~t \vee \hat{t} \leq \eps 
				\bigg \} \\ \nonumber
			={}&
				\lim_{r \downarrow 0} \lim_{\eps \downarrow 0}
				\sup \bigg \{
					\frac {u_1(t,x)}{V(t,x)} 
					- \frac {u_2(\hat{t} , \hat{x})}{V(\hat{t} , \hat{x})} \colon 
					~(t, x), (\hat{t}, \hat{x}) \in W,
					~h(t, x) \vee h(\hat{t}, \hat{x}) \leq R,
			\\ \nonumber & \qquad\qquad\qquad
					~\|x\|_H \leq n,
					~\|x- \hat{x}\|_H \leq r, 
					~t \vee \hat{t} \leq \eps 
				\bigg \} \\
				&\vee
				\lim_{r \downarrow 0} \lim_{\eps \downarrow 0}
				\sup \bigg \{
					\frac {u_1(t,x)}{V(t,x)} 
					- \frac {u_2(\hat{t} , \hat{x})}{V(\hat{t} , \hat{x})} \colon 
					~(t, x), (\hat{t}, \hat{x}) \in W,
					~h(t, x) \vee h(\hat{t}, \hat{x}) \leq R,
			\\ \nonumber & \qquad\qquad\qquad
					~\|x\|_H \geq n,
					~\|x- \hat{x}\|_H \leq r,
					~t \vee \hat{t} \leq \eps 
				\bigg \} \\ \nonumber
			\leq{}&
				\lim_{r \downarrow 0} \lim_{\eps \downarrow 0}
				\sup \bigg \{
					u_1(t,x) - u_2(\hat{t} , \hat{x}) \colon 
					~(t, x), (\hat{t}, \hat{x}) \in W,
					~h(t, x) \vee h(\hat{t}, \hat{x}) \leq R,
			\\ \nonumber & \qquad\qquad\qquad
					~\|x\|_H \leq n,
					~\|x- \hat{x}\|_H \leq r, 
					~t \vee \hat{t} \leq \eps 
				\bigg \} \\ \nonumber
				&\vee
				\lim_{r \downarrow 0} \lim_{\eps \downarrow 0}
				\sup \bigg \{
					\frac {u_1(t,x)}{V(t,x)} 
					- \frac {u_2(\hat{t} , \hat{x})}{V(\hat{t} , \hat{x})} \colon 
					~(t, x), (\hat{t}, \hat{x}) \in ((0,T) \times O_{ \lfloor n-r \rfloor}^{c}) \cap W
				\bigg \} 
				\vee 0 \\ \nonumber
			\leq{}&
				0 \vee \left( 2\sup \bigg \{
					\frac {u_1(t,x)}{V(t,x)} 
					\vee \frac {-u_2(\hat{t} , \hat{x})}{V(\hat{t} , \hat{x})} \colon 
					~(t, x), (\hat{t}, \hat{x}) \in ((0,T) \times O_{n-1}^{c}) \cap W
				\bigg \} \right). 
	\end{align}
	Letting $n$ tend to infinity together with 
	\eqref{eq:attains.maximum.difference2.comparison} then shows that
	for all $R \in (0, \infty)$
	it holds that 
	\begin{align}
		\nonumber
				&\lim_{r \downarrow 0} \lim_{\eps \downarrow 0}
					\sup \Big \{
						\tilde{u}_1(t,x) - \tilde{u}_2(\hat{t} , \hat{x}) \colon 
						~(t, x), (\hat{t}, \hat{x}) \in W,
						~h(t, x) \vee h(\hat{t}, \hat{x}) \leq R,
			\\ & \qquad\qquad\qquad	
						~\|x- \hat{x}\|_H \leq r,
						~t \vee \hat{t} \leq \eps 
					\Big \} \\ \nonumber
			\leq{}&
				0 \vee \lim_{n \to \infty} \Big( 2\sup \Big \{
					\tfrac {u_1(t,x)}{V(t,x)} 
					\vee \tfrac {-u_2(\hat{t} , \hat{x})}{V(\hat{t} , \hat{x})} \colon 
					~(t, x), (\hat{t}, \hat{x}) \in ((0,T) \times O_{n-1}^{c}) \cap W
				\Big \} \Big)
			= 0
	\end{align}
	and this shows \eqref{eq: uniformly bounded at 0 COR}.
Furthermore, we have
for all $\delta \in (0, \infty)$ and all 
$((t,x),r,p,A) \in W \times \R \times H' \times \mathbb{S}_{\mathbb{H},\mathbb{H}'}$ that
\begin{align}
\label{eq: tilde G + and G +}
			&\tilde{G}_{\mathbb{H}, \mathbb{X}, \delta, h}^{+}((t,x),r,p,A) \\ \nonumber
		={}&
			\tilde{G}(
								(t,x),
								r +  \delta h(t,x),
								p
								+ \delta E_{\mathbb{X}', \mathbb{H}'} 
										\big( D_{\mathbb{X}}(h|_{(0,T) \times (O \cap X)}) \big) (t,x),
								(A|_X)|_X 
		\\ \nonumber & \qquad\quad
								+ \delta (D_{\mathbb{X}}^2(h|_{(0,T) \times (O \cap X)})) (t, x) \\ \nonumber
		={}&
			\tfrac{1}{V(t,x)} \Big( 
						G \Big(
							(t,x),
							V(t,x) r + \delta V(t,x) \cdot h(t,x),
							V(t,x) p 
							+ \delta V(t,x) 
			\\ \nonumber & \qquad\quad
							\cdot E_{\mathbb{X}', \mathbb{H}'} 
										\big(D_{\mathbb{X}}(h|_{(0,T) \times (O \cap X)}) \big) (t,x) 
								+(r + \delta h(t,x)) \cdot (D_{\mathbb{H}} V) (t,x), 
							V(t,x) \cdot (A |_X)|_X 
			\\ \nonumber & \qquad\quad
								+ \delta V(t,x) \cdot (D_{\mathbb{X}}^2(h|_{(0,T) \times (O \cap X)})) (t, x) 
								+ (D_{\mathbb{H}} V) (t,x)|_X 
			\\ \nonumber & \qquad\quad
									\otimes (p|_X + \delta (D_{\mathbb{X}}(h|_{(0,T) \times (O \cap X)}))(t,x)
								+  (p|_X + \delta (D_{\mathbb{X}}(h|_{(0,T) \times (O \cap X)}))(t,x)) 
			\\ \nonumber & \qquad\quad
										\otimes (D_{\mathbb{H}} V) (t,x) |_X
								+ (r + \delta h(t,x)) \cdot ((D^2_{\mathbb{H}} \, V) (t,x) |_X)|_X
						\Big) \Big) 
						- (r + \delta h(t,x)) \cdot \tfrac{ \frac{\partial}{\partial t} V(t,x)}{V(t,x)} \\ \nonumber
		={}&
			\tfrac{1}{V(t,x)} \Big(
				G^+_{\mathbb{H}, \mathbb{X}, \delta, Vh} \Big(
							(t,x),
							V(t,x) r,
							V(t,x) p + r \cdot (D_{\mathbb{H}} V) (t,x), 
							V(t,x) \cdot A + (D_{\mathbb{H}} V) (t,x)  \otimes p \\ \nonumber
						& \qquad \qquad
								+ p \otimes (D_{\mathbb{H}} V) (t,x) 			 
								+ r \cdot (D^2_{\mathbb{H}} \, V) (t,x)  
						\Big)
				- (r + \delta h(t,x)) \cdot \tfrac{\partial}{\partial t} V(t,x)
			\Big).
\end{align}
Analogously we obtain for all $\delta \in (0, \infty)$ and all 
$((t,x),r,p,A) \in W \times \R \times H' \times \mathbb{S}_{\mathbb{H},\mathbb{H}'}$ that
\begin{equation}
\label{eq: tilde G - and G -}
	\begin{split}
			&\tilde{G}_{\mathbb{H}, \mathbb{X}, \delta, h}^{-}((t,x),r,p,A ) \\
		={}&
			\tfrac{1}{V(t,x)} \Big(
				G^-_{\mathbb{H}, \mathbb{X}, \delta, Vh} \Big(
							(t,x),
							V(t,x) r,
							V(t,x) p + r \cdot (D_{\mathbb{H}} V) (t,x), 
							V(t,x) \cdot A 
							+ (D_{\mathbb{H}} V) (t,x) \otimes p
			\\ & \qquad \quad
							+ p \otimes (D_{\mathbb{H}} V) (t,x) 
							+ r \cdot (D^2_{\mathbb{H}} \, V) (t,x)
						\Big)
				- (r - \delta h(t,x)) \cdot \tfrac{\partial}{\partial t} V(t,x)
			\Big).
	\end{split}
\end{equation}
In addition,
\eqref{eq: V bound} yields that 
for all $((t,x),r,p,A) \in W \times \R \times H' \times \mathbb{S}_{\mathbb{H},\mathbb{H}'}$,
$\delta \in (0, 1]$
and all $R \in (0, \infty)$
with
$
	h(t, x) \vee \|x\|_H \vee |r| \vee \|p\|_{H'} \vee \| A \|_{L(\mathbb{H}, \mathbb{H}')}  
		\leq R
$ 
it holds that
\begin{align}
\nonumber
			&h(t, x) \leq R \leq 4 C_R \cdot R, \quad
			\|x\|_H \leq R \leq 4 C_R \cdot R, \quad
			|V(t,x) r| \leq C_R \cdot R \leq 4 C_R \cdot R, \\ \nonumber
			&\| V(t,x) p + r \cdot (D_{\mathbb{H}} V) (t,x)\|_{H'}
				\leq 2 C_R \cdot R \leq 4 C_R \cdot R, \\ \nonumber
			&\| V(t,x) \cdot A + (D_{\mathbb{H}} V) (t,x) \otimes p
					+ p \otimes (D_{\mathbb{H}} V) (t,x) 
					+ r \cdot (D^2_{\mathbb{H}} \, V) (t,x) 
			\|_{L(\mathbb{H},\mathbb{H}')}
				\leq 4 C_R \cdot R, \\
	\label{eq: different bound}
			& (r+ \delta h(t,x)) \cdot \tfrac {\partial} {\partial t} V(t,x)
				\geq - 2C_R \cdot R.
\end{align}
We thus get from \eqref{eq: uniformly bounded at p and A+}, 
\eqref{eq: tilde G + and G +},
from \eqref{eq: different bound},
and from the assumption that $V >1$ that
for all $\delta \in (0,1]$ 
and all $R \in (0, \infty)$ it holds that 
\begin{align}
				&\lim_{\eps \to 0} \Big [
					\sup \{|\tilde{G}_{\mathbb{H}, \mathbb{X}, \delta, h}^{+}((t,x),r,p,A)
						-\tilde{G}_{\mathbb{H}, \mathbb{X}, \delta, h}^{+}((t,x),r,q,A)| \colon
						(t, x) \in W,
						~ r \in \R, 
						~ p, q \in H',
				\\& \nonumber \qquad \qquad 
						~ A \in \mathbb{S}_{\mathbb{H}, \mathbb{H}'},  
						~\max \{ h(t, x), \|x\|_H, |r|, \|p\|_{H'}, \| A \|_{L(\mathbb{H}, \mathbb{H}')} \} 
							\leq R,
						~\|p-q\|_{H'} \leq \eps,
				\\& \nonumber \qquad \qquad 
						~\tilde{G}_{\mathbb{H}, \mathbb{X}, \delta, h}^{+}((t,x),r,p,A) \geq -R\}
				\Big ] \\ \nonumber
			={}&
				\lim_{\eps \to 0} \Big [
					\sup \Big \{ \tfrac{1}{V(t,x)} \Big|
						G^+_{\mathbb{H}, \mathbb{X}, \delta, Vh} \Big(
							(t,x),
							V(t,x) r,
							V(t,x) p + r \cdot (D_{\mathbb{H}} V) (t,x), 
							V(t,x) \cdot A
			\\ \nonumber & \qquad 
							+ (D_{\mathbb{H}} V) (t,x) \otimes p 
								+ p \otimes (D_{\mathbb{H}} V) (t,x)
								+ r \cdot (D^2_{\mathbb{H}} \, V) (t,x)
						\Big) 
						- (r + \delta h(t,x)) \cdot \tfrac{\partial}{\partial t} V(t,x) \\ \nonumber
					& \quad
						-G^+_{\mathbb{H}, \mathbb{X}, \delta, Vh} \Big(
							(t,x),
							V(t,x) r,
							V(t,x) q + r \cdot (D_{\mathbb{H}} V) (t,x), 
							V(t,x) \cdot A + (D_{\mathbb{H}} V) (t,x) \otimes q
					\\ \nonumber & \qquad 
								+ q \otimes (D_{\mathbb{H}} V) (t,x) 
								+ r \cdot (D^2_{\mathbb{H}} \, V) (t,x)
						\Big) 
						+ (r + \delta h(t,x)) \cdot \tfrac{\partial}{\partial t} V(t,x) \Big | \colon 
						(t, x) \in W,
						~ r \in \R, 
				\\ \nonumber & \qquad 
						~ p, q \in H',
						~ A \in \mathbb{S}_{\mathbb{H}, \mathbb{H}'}, 
						~\max \{ h(t, x), \|x\|_H, |r|, \|p\|_{H'}, \| A \|_{L(\mathbb{H}, \mathbb{H}')} \}
							\leq R,
						~\|p-q\|_{H'} \leq \eps, \\ \nonumber
						& \qquad 
						~\tfrac {1}{V(t,x)} \Big(
							G^+_{\mathbb{H}, \mathbb{X}, \delta, Vh} \Big(
								(t,x),
								V(t,x) r,
								V(t,x) p + r \cdot (D_\mathbb{H} V) (t,x), 
						\\ \nonumber & \qquad
								V(t,x) \cdot A + (D_{\mathbb{H}} V) (t,x) \otimes p 
									+ p \otimes (D_{\mathbb{H}} V) (t,x) 			
									+ r \cdot (D^2_{\mathbb{H}} \, V) (t,x) 
							\Big) 
						\\ \nonumber  & \qquad 
							- (r +  \delta h(t,x)) \cdot \tfrac{\partial}{\partial t} V(t,x) 
						\Big) \geq -R
					\Big\}
				\Big ] \\ \nonumber
			\leq{}&
				\lim_{\eps \to 0} \Big [
					\sup \Big \{ |
						G^+_{\mathbb{H}, \mathbb{X}, \delta, Vh} ((t, x), r, p, A)  
						-G^+_{\mathbb{H}, \mathbb{X}, \delta, Vh} ((t, x), r, q, B) | \colon 
						(t, x) \in W,
						~ r \in \R, 
						~ p, q \in H' \! ,
				\\ \nonumber & \qquad 
						~ A, B \in \mathbb{S}_{\mathbb{H}, \mathbb{H}'}, 
						~\max \{ h(t, x), \|x\|_H, |r|, \|p\|_{H'}, \| A \|_{L(\mathbb{H}, \mathbb{H}')} \} 
							\leq 4 C_R \, R,
						~\|p-q\|_{H'} \leq C_R \, \eps,
				\\ \nonumber & \qquad  
						~\|A-B\|_{L(\mathbb{H},\mathbb{H}')} \leq 2 C_R \cdot \eps, 
						~G_{\mathbb{H}, \mathbb{X}, \delta, Vh}^{+}((t,x),r,p,A) \geq -3 C_R \cdot R 
					\Big\}
				\Big ] \vee 0= 0.
\end{align}
and thus \eqref{eq: uniformly bounded at p COR+} is fulfilled.
Analogously it follows \eqref{eq: uniformly bounded at p COR-}.
	Moreover, the assumption that $D^2_\mathbb{H} \, V$ is 
	continuous
	with respect to the
	$\| \cdot \|_{\R \times H}$ and the
	$\| \cdot \|_{L(\mathbb{H},\mathbb{H}')}$-norm
	ensures
	that there exist a function
	$K \colon (0, \infty) \times (0,T) \times O \times \R \times H' \times (0,1] \to (0, \infty)$
	such that for all 
	$(R,t,x,r,p,\delta) \in  (0,T) \times O \times \R \times H' \times (0,1]$,
	and all $(\tau, \xi) \in  (0,T) \times O$ with
	$|t - \tau| \vee \| \xi - x\|_{H} \leq K(R,t,x,r,p,\delta)$
	it holds that
	\begin{equation}
	\label{eq: def of K}
		\begin{split}
				&\|
					(D^2_{\mathbb{H}} \, V) (\tau, \xi) -(D^2_{\mathbb{H}} \, V) (t, x) 
				\|_{L(\mathbb{H},\mathbb{H}')}
			\leq
				 \frac{m(C_R (2R +2 \| p\|_{H'} +|r| + 3),t,x,r,p,\delta)}{2 (|r| +1)}.
		\end{split}	
	\end{equation}
	Now 
	denote
	by
	$
		\tilde{m} \colon 
			(0, \infty) 
			\times (0,T) 
			\times O 
			\times \R 
			\times H' 
			\times (0,1] 
			\to (0, \infty)
	$
	the function satisfying for all 
	$(R,t,x,r,p,\delta) \in (0,\infty) \times (0,T) \times O \times \R \times H' \times (0,1]$
	that
	\begin{equation}
	\label{eq: def of tilde m}
			\tilde{m}(R,t,x,r,p,\delta)
		=
			\frac
				{m(C_R (2R +2 \| p\|_{H'} +|r| + 3),t,x,r,p,\delta)}
				{2 C_R (2R+ 4\|p\|_{H'} + 2|r| + 7)}
			\wedge K(R,t,x,r,p,\delta)
			\wedge 1.	
	\end{equation}
	Combining then \eqref{eq: V bound}, \eqref{eq: def of K}, 
	and \eqref{eq: def of tilde m} yields that for all 
	$R \in (0, \infty)$,
	$\delta \in (0,1]$,
	$(t,x,r,p,A) \in (0,T) \times O \times \R \times H' \times L(\mathbb{H},\mathbb{H}')$,
	and all
	$((\tau, \xi), \nu, \rho) \in W \times \R \times H'$ satisfying that
	$\|x\|_H +1 \leq R$, 
	$ h(\tau, \xi) \leq R $,
	$\| A \|_{L(\mathbb{H},\mathbb{H}')} \leq R$, and that
	$ 
		| t-\tau | \vee \|x- \xi\|_H \vee |r -\nu| \vee \|p - \rho\|_{H'} 
				\leq \tilde{m}(R,t,x,r,p,\delta)
	$
	it holds that
	\begin{equation}
	\label{eq: new A bound}
		\begin{split}
				&\big\| 
					V(\tau, \xi) \cdot A 
					+ (D_{\mathbb{H}} V) (\tau, \xi) \otimes \rho 
					+ \rho \otimes (D_{\mathbb{H}} V) (\tau, \xi) 
					+ \nu \cdot (D^2_{\mathbb{H}} \, V) (\tau, \xi) 
				\big\|_{L(\mathbb{H},\mathbb{H}')} \\
			\leq{}&
				C_{R} \cdot R
				+ C_R (\|p\|_{H'}+1)
				+ C_R (\|p\|_{H'}+1)
				+ C_R (|r|+1)
			\leq
				C_R (2R+ |r|+ 2 \|p\|_{H'}+3),
		\end{split}
	\end{equation}
	that
	\begin{equation}
	\label{eq: new G bound}
		\begin{split}
				&(r + \delta h(\tau, \xi)) \cdot \tfrac{\partial}{\partial t} V(\tau, \xi)
			\leq
				(|r| + R) \cdot C_R
			\leq
				C_R (R+ |r|+ 2 \|p\|_{H'}+3),
		\end{split}
	\end{equation}
	that
	\begin{equation}
	\label{eq: new nu close to new r}
		\begin{split}
				&|\nu \cdot V(\tau, \xi) - r V(t,x)|
			\leq
				|\nu -r | \cdot V(t, x) + |\nu| \cdot |V(\tau, \xi) - V(t,x)| \\
			\leq{}&
				\tilde{m}(R,t,x,r,p,\delta) \cdot C_R 
				+ ( |r| + 1 ) 
					\cdot C_R \|(\tau, \xi) - (t, x) \|_{\R \times H} \\
			\leq{}&
				\tilde{m}(R,t,x,r,p,\delta) \cdot C_R (1 +2 |r| +2 ) 
			\leq
				m(C_R (2R +2 \| p\|_{H'} +|r| + 3),t,x,r,p,\delta),
		\end{split}
	\end{equation}
	that
	\begin{equation}
	\label{eq: new rho close to new p}
		\begin{split}
				&\|V(\tau, \xi) \rho + \nu \cdot (D_{\mathbb{H}} V) (\tau, \xi)
					- V(t, x) p - r \cdot (D_{\mathbb{H}} V) (t, x) \|_{H'} \\
			\leq{}&
				|V(\tau, \xi) - V(t, x)| \cdot \| \rho \|_{H'}
				+ V(t, x) \| \rho - p \|_{H'}
				+ |\nu| \cdot \| (D_{\mathbb{H}} V) (\tau, \xi) - (D_{\mathbb{H}} V) (t, x) \|_{H'} \\
				&+ \| (D_{\mathbb{H}} V) (t, x) \|_{H'} \cdot | \nu - r| \\
			\leq{}&
				C_R \|(\tau, \xi) - (t, x) \|_{\R \times H} 
					\cdot ( \|p\|_{H'} + 1 )
				+ C_R \cdot \tilde{m}(R,t,x,r,p,\delta) \\
				&+ ( |r| + 1 ) 
					\cdot C_R \|(\tau, \xi) - (t, x) \|_{\R \times H} 
				+ \tilde{m}(R,t,x,r,p,\delta) \cdot C_R \\
			\leq{}&
				2 \tilde{m}(R,t,x,r,p,\delta) \cdot C_R ( 
					3+ \|p \|_{H'} + |r|
				)
			\leq
				m(C_R (2R +2 \| p\|_{H'} +|r| + 3),t,x,r,p,\delta),
		\end{split}
	\end{equation}
	that
	\begin{align}
	\nonumber
				& \big \|
					(D_\mathbb{H} V) (\tau, \xi) \otimes \rho
					+ \rho \otimes (D_\mathbb{H} V) (\tau, \xi)  
					- (D_\mathbb{H} V) (t,x) \otimes p
					- p \otimes (D_\mathbb{H} V) (t,x)
				\big \|_{L(\mathbb{H},\mathbb{H}')} \\ \nonumber
			\leq{}& 
				\|(D_\mathbb{H} V) (\tau, \xi) - (D_\mathbb{H} V) (t, x) \|_{H'} 
					\cdot \|\rho \|_{H'} 
				+ \|(D_\mathbb{H} V) (t, x) \|_{H'} \cdot \| \rho - p \|_{H'} \\
	\label{eq: new A close to new tilde A1}
				&+ \|\rho \|_{H'} \cdot \|(D_\mathbb{H} V) (\tau, \xi) - (D_\mathbb{H} V) (t, x) \|_{H'} 
				+ \| \rho - p \|_{H'} \cdot \|(D_\mathbb{H} V) (t, x) \|_{H'} \\ \nonumber
			\leq{}&
				C_R \|(\tau, \xi) - (t, x) \|_{\R \times H} 
					\cdot  ( \|p\|_{H'} + 1 ) 
				+ C_R \cdot \tilde{m}(R,t,x,r,p,\delta) \\ \nonumber
				&+  ( \|p\|_{H'} + 1 )
					\cdot C_R \|(\tau, \xi) - (t, x) \|_{\R \times H}
				+ \tilde{m}(R,t,x,r,p,\delta) \cdot C_R \\ \nonumber
			\leq{}&
				C_R \cdot \tilde{m}(R,t,x,r,p,\delta) ( 
					4 \|p\|_{H'} +6
				),
	\end{align}
	and that
	\begin{equation}
	\label{eq: new A close to new tilde A2}
		\begin{split}
				& \big \| V(\tau, \xi) \cdot A
				+ \nu \cdot (D^2_\mathbb{H} \, V) (\tau, \xi) 
				-V(t,x) \cdot A 
				- r \cdot (D^2_\mathbb{H} \, V) (t, x) \big \|_{L(\mathbb{H},\mathbb{H}')} \\
			\leq{}& 
				|V(\tau, \xi) - V(t, x) | \cdot \|A\|_{L(\mathbb{H},\mathbb{H}')} 
				+ \| (D^2_\mathbb{H} \, V) (t, x) \|_{L(\mathbb{H},\mathbb{H}')} \cdot | \nu - r| \\
				&+ |\nu| \cdot \| 
						(D^2_\mathbb{H} \, V) (\tau, \xi) - (D^2_\mathbb{H} \, V) (t, x) 
					\|_{L(\mathbb{H},\mathbb{H}')} \\
			\leq{}&
				C_R \|(\tau, \xi) - (t, x) \|_{\R \times H} \cdot \|A\|_{L(\mathbb{H},\mathbb{H}')} 
				+ \tilde{m}(R,t,x,r,p,\delta) \cdot C_R  \\
				&+ ( |r| + 1 ) \cdot \| 
						(D^2_\mathbb{H} \, V) (\tau, \xi) - (D^2_\mathbb{H} \, V) (t, x) 
					\|_{L(\mathbb{H},\mathbb{H}')} \\
			\leq{}&
				C_R \cdot \tilde{m}(R,t,x,r,p,\delta) ( 
					2 R + 1
				) 
				+ ( |r| + 1 ) 
					\cdot \tfrac{m(C_R (2R +2 \| p\|_{H'} +|r| + 3),t,x,r,p,\delta)}{2(|r| + 1)}. 
		\end{split}
	\end{equation}
	Combining \eqref{eq: new A close to new tilde A1},
	\eqref{eq: new A close to new tilde A2}, 
	and \eqref{eq: def of tilde m} proves
	\begin{align}
	\label{eq: new A close to new tilde A}
				&\big \| 
					V(\tau, \xi) \cdot A
					+(D_\mathbb{H} V) (\tau, \xi) \otimes \rho
					+ \rho \otimes (D_\mathbb{H} V) (\tau, \xi)  
					+ \nu \cdot (D^2_\mathbb{H} \, V) (\tau, \xi) \\ \nonumber
					&-V(t,x) \cdot A
					- (D_\mathbb{H} V) (t,x) \otimes p
					- p \otimes (D_\mathbb{H} V) (t,x)
					- r \cdot (D^2_\mathbb{H} \, V) (t, x)
				\big \|_{L(\mathbb{H},\mathbb{H}')} \\ \nonumber
		\leq{}
			&\big \| V(\tau, \xi) \cdot A
				+ \nu \cdot (D^2_\mathbb{H} \, V) (\tau, \xi) 
				-V(t,x) \cdot A 
				- r \cdot (D^2_\mathbb{H} \, V) (t, x) \big \|_{L(\mathbb{H},\mathbb{H}')} \\ \nonumber
			&+\big \|
					(D_\mathbb{H} V) (\tau, \xi) \otimes \rho
					+ \rho \otimes (D_\mathbb{H} V) (\tau, \xi)  
					- (D_\mathbb{H} V) (t,x) \otimes p
					- p \otimes (D_\mathbb{H} V) (t,x)
			\big \|_{L(\mathbb{H},\mathbb{H}')} \\ \nonumber
		\leq{}&
				C_R \cdot \tilde{m}(R,t,x,r,p,\delta) ( 
					2 R + 4 \|p\|_{H'} +7
				) 
				+\tfrac{m(C_R (2R +2 \| p\|_{H'} +|r| + 3),t,x,r,p,\delta)}{2} \\ \nonumber
			\leq{}
				&\tfrac{m(C_R (2R +2 \| p\|_{H'} +|r| + 3),t,x,r,p,\delta)}{2} 
				+ \tfrac{m(C_R (2R +2 \| p\|_{H'} +|r| + 3),t,x,r,p,\delta)}{2} \\ \nonumber
			\leq{}&
				m(C_R (2R +2 \| p\|_{H'} +|r| + 3),t,x,r,p,\delta).
	\end{align}
	Thus 
	\eqref{eq: uniformly bounded in dimension+ LEM},
	\eqref{eq: tilde G + and G +},
	\eqref{eq: def of tilde m},
	\eqref{eq: new A bound},
	\eqref{eq: new G bound},
	\eqref{eq: new nu close to new r},
	\eqref{eq: new rho close to new p},
	\eqref{eq: new A close to new tilde A},
	and
	the assumption $V>1$
	shows that for all 
	$t \in (0,T)$, $x \in O$, $r \in \R$, $p \in H'$, $R \in (\|x\|_H+1,\infty)$ 
	and all $\delta \in (0,1]$ it holds that 
	\begin{align}
	\label{eq: uniformly bounded in dimension with R bound}
				&\lim_{N \to \infty} \big [
					\sup
					\{ \tilde{G}_{\mathbb{H}, \mathbb{X}, \delta, h}^{+}
						((\tau, \xi), \nu, \rho, A + \alpha I_{\mathbb{H}} \pi^{H}_{H^\perp_N})
						-\tilde{G}_{\mathbb{H}, \mathbb{X}, \delta, h}^{+}
							((\tau, \xi), \nu, \rho, A) \colon 
					\alpha \in (0, R),
					~ \nu \in \R, 
			~ \\ \nonumber & \quad 
					~ (\tau, \xi) \in W,
					~ \rho \in H', 
					~ A \in \mathbb{S}_{\mathbb{H},\mathbb{H}'}, 
					~ h(\tau, \xi) \leq R,
					~ \| A \|_{L(\mathbb{H},\mathbb{H}')} \leq R, 
			~ \\ \nonumber & \quad 
					~ | t-\tau | \vee \|x- \xi\|_H \vee |r -\nu| \vee \|p - \rho\|_{H'} 
							\leq \tilde{m}(R,t,x,r,p,\delta), 
			~ \\\nonumber & \quad
					~ \tilde{G}_{\mathbb{H}, \mathbb{X}, \delta, h}^{+}
							((\tau, \xi), \nu, \rho, A) \geq -R \}
				\big ] \\ \nonumber
			={}&
				\lim_{N \to \infty} \Big [
					\sup
					\Big \{ \tfrac{1}{V(\tau, \xi)} \Big(
						G^+_{\mathbb{H}, \mathbb{X}, \delta, Vh} \Big(
							(\tau, \xi),
							V(\tau, \xi) \nu,
							V(\tau, \xi) \rho + \nu \cdot (D_{\mathbb{H}} V) (\tau, \xi), 
							V(\tau, \xi)
					\\ \nonumber & \quad
							\cdot (A + \alpha I_{\mathbb{H}} \pi^{H}_{H_N^\perp}) 
								+ (D_{\mathbb{H}} V) (\tau, \xi) \otimes \rho 
								+ \rho \otimes (D_{\mathbb{H}} V) (\tau, \xi)
								+ \nu \cdot (D^2_{\mathbb{H}} \, V) (\tau, \xi) 
						\Big)
				\\ \nonumber & \quad
						- (\nu + \delta h(\tau, \xi)) 	
						\cdot \tfrac{\partial}{\partial t} V(\tau, \xi)
					\Big) 
						-\tfrac{1}{V(\tau, \xi)} \Big(
						G^+_{\mathbb{H}, \mathbb{X}, \delta, Vh} \Big(
							(\tau, \xi),
							V(\tau, \xi) \nu,
							V(\tau, \xi) \rho 
					\\ \nonumber & \quad	
							+ \nu \cdot (D_{\mathbb{H}} V) (\tau, \xi), 
							V(\tau, \xi) \cdot A+ (D_{\mathbb{H}} V) (\tau, \xi) \otimes \rho 
								+ \rho \otimes (D_{\mathbb{H}} V) (\tau, \xi)
								+ \nu \cdot (D^2_{\mathbb{H}} \, V) (\tau, \xi)
						\Big)
				\\ \nonumber & \quad
						- (\nu + \delta h(\tau, \xi)) \cdot \tfrac{\partial}{\partial t} V(\tau, \xi)
					\Big) \colon 
					\alpha \in (0, R),
					~ \nu \in \R, 
					~ (\tau, \xi) \in W, 
					~ \rho \in H',
					~ A \in \mathbb{S}_{\mathbb{H},\mathbb{H}'}, 
			\\\nonumber & \quad
					~ h(\tau, \xi) \leq R,
					~ \| A \|_{L(\mathbb{H},\mathbb{H}')} \leq R,
					~ | t-\tau | \vee \|x- \xi\|_H \vee |r -\nu| \vee \|p - \rho\|_{H'} 
							\leq \tilde{m}(R,t,x,r,p,\delta), \\ \nonumber 
				& \quad 
					~ \tfrac{1}{V(\tau, \xi)} \Big(
				G^+_{\mathbb{H}, \mathbb{X}, \delta, Vh} \Big(
							(\tau, \xi),
							V(\tau, \xi) r,
							V(\tau, \xi) p + r \cdot (D_{\mathbb{H}} V) (\tau, \xi), 
							V(\tau, \xi) \cdot A  
							+ (D_{\mathbb{H}} V) (\tau, \xi) \otimes p
				\\ \nonumber & \quad
							+ p \otimes (D_{\mathbb{H}} V) (\tau, \xi) 
								+ r \cdot (D^2_{\mathbb{H}} \, V) (\tau, \xi) 
						\Big)
				- (r + \delta h(\tau, \xi)) \cdot \tfrac{\partial}{\partial t} V(\tau, \xi)
			\Big) \geq -R 
				\Big \} \Big ] \\ \nonumber
			\leq{}&
				\lim_{N \to \infty} \Big [
					\sup
					\Big \{ \tfrac{1}{V(\tau, \xi)} \Big(
						G^+_{\mathbb{H}, \mathbb{X}, \delta, Vh} (
							(\tau, \xi),
							\nu,
							\rho, 
							A+ \alpha I_{\mathbb{H}} \pi^{H}_{H_N^\perp} 
						)
						-G^+_{\mathbb{H}, \mathbb{X}, \delta, Vh} (
							(\tau, \xi),
							\nu,
							\rho, 
							A
						) \Big) \colon
			\\ \nonumber & \quad 
					\alpha \in (0, C_R\cdot R),
					~ \nu \in \R, 
					~ (\tau, \xi) \in W,
					~ \rho \in H',
					~ A \in \mathbb{S}_{\mathbb{H}, \mathbb{H}'}, 
					~ h(\tau, \xi) \leq R, 
	 		\\ \nonumber & \quad 
					~ \| A \|_{L(\mathbb{H}, \mathbb{H}')} \leq C_R (2R +2 \| p\|_{H'} +|r| + 3), 
					~ | t-\tau | \vee \|x- \xi\|_H 
						\vee |V(t,x) r -\nu| 
			\\ \nonumber & \quad
						\vee \|V(t, x) p + r \cdot (D_{\mathbb{H}} V) (t, x)  - \rho\|_{H'} 
							\leq m(C_R (2R +2 \| p\|_{H'} +|r| + 3),t,x,r,p,\delta), 
			\\ \nonumber & \quad
				~G^+_{\mathbb{H}, \mathbb{X}, \delta, Vh} \Big(
							(\tau, \xi),
							\nu,
							\rho, 
							A
						\Big)
				\geq - C_R (2R +2 \| p\|_{H'} +|r| + 3)
				\Big \} \Big ]
				\\ \nonumber
			\leq{}&
				0.
	\end{align}
	Obviously \eqref{eq: uniformly bounded in dimension with R bound} then also holds for all
	$t \in (0,T)$, $x \in O$, $r \in \R$, $p \in H'$, $R \in (0,\infty)$ 
	and all $\delta \in (0,1]$.
	and this implies \eqref{eq: uniformly bounded in dimension COR+}.
	Analogously it follows \eqref{eq: uniformly bounded in dimension COR-}.
  It remains to verify \eqref{eq:00assumption_COR}.
	Therefore note that without loss of generality we can assume that
	$\beta_n \geq n$ since otherwise we can replace $\beta_n$ by $\beta_n \vee n$.
  Now let
  $ 
    ( 
      (t_n, x_n), r_n, A_n 
    )
    ,
    ( 
      (\hat{t}_n, \hat{x}_n), 
      \hat{r}_n, 
      \hat{A}_n 
    )
    \in
    W \times \R \times \mathbb{S}_{\mathbb{H}, \mathbb{H}'} 
  $,
  $  n \in \N $,
  be sequences 
  satisfying that
  $
    w-\lim_{ n \to \infty }
    ( t_n, x_n )
    \in [ 0, T ) \times O
  $,
  that
  $
    \lim_{ n \to \infty }
    \big(
      \sqrt{ n } 
      \,
      \|
        x_n - \hat{x}_n 
      \|_{H}
    \big)
    = 0
  $,
	that
  $
    \lim_{ n \to \infty }
    \big(
      \sqrt{ \beta_n } 
      \,
      |
        t_n - \hat{t}_n 
      |
    \big)
    = 0
  $,
  that
  $
    \lim_{ n \to \infty }
    \big(
      \tilde{\delta}_n ( 
				h(t_{n}, x_{n})
        + h(\hat{t}_{n}, \hat{x}_{n})
			)
    \big)
    = 0
  $,
  that
  $
    0 <
    \lim_{ n \to \infty }
    \left( 
      r_n - \hat{r}_n
    \right)
    \leq
    \sup_{ n \in \N } 
    ( 
      | r_n |
      +
      | \hat{r}_n |
    )
    < \infty
  $,
  and that for all
	$R \in (0, \infty)$,
	and all
	$
		(z^{(n)})_{n \in \N}, (\hat{z}^{(n)})_{n \in \N} 
		\subseteq \{z \in H \colon ~ \| z \|_H \leq R\}
	$
	it holds that 
	\begin{equation}
			\limsup_{n \to \infty} \left(
			\langle z^{ (n) }, A^{(n)} z_i^{ (n) } \rangle_{H, H'}
			-
			\langle \hat{z}^{ (n) }, \hat{A}^{(n)} \hat{z}_i^{ (n) } \rangle_{H, H'}
			\right)
		\leq 
			\limsup_{n \to \infty}
			3
			\| z^{ (n) } - \hat{z}^{ (n) }
			\|_H^2.
  \end{equation}
  To verify \eqref{eq:00assumption_COR},
  we will apply assumption~\eqref{eq:comparison.viscosity.solution.assumption}.
  For this we denote by
  $ 
    \big( 
      ({\bf t}_n, {\bf x}_n),
      {\bf r}_n,
      {\bf p}_n
    \big)
    \in
    W \times 
    \R \times 
    H' 
  $,
	by
	$ 
    \big( 
      ({\bf \hat{t}}_n, {\bf \hat{x}}_n),
      {\bf \hat{r}}_n,
      {\bf \hat{p}}_n,
    \big)
    \in
    W \times 
    \R \times 
    H' 
  $,
	by
	$
		{\bf A}_n,
    {\bf B}_n,
		{\bf \hat{A}}_n,
    {\bf \hat{B}}_n
		\in
    \mathbb{S}_{\mathbb{H}, \mathbb{H}'}
	$,
	and by
	$
		({\bf \tilde{r}}_n, {\bf \tilde{p}}_n, {\bf \tilde{A}}_n) 
			\in \R \times H' \times \mathbb{S}_{\mathbb{H}, \mathbb{H}'} 
	$,
  $ n \in \N $,
  the elements satisfying for all $n \in \N$
  that
  $
    ({\bf t}_n,
     {\bf x}_n,
     {\bf r}_n
    )
    = (t_n,
        x_n,
        r_n V( t_n, x_n )
       )
  $, that
  $
    ({\bf \hat{t}}_n,
     {\bf \hat{x}}_n,
     {\bf \hat{r}}_n
    )
    = (\hat{t}_n,
        \hat{x}_n,
        \hat{r}_n V( \hat{t}_n, \hat{x}_n )
       )
  $,
	that
  $
    {\bf B}_n =
    A_n
  $,
	that
  $
    {\bf \hat{B}}_n =
    \hat{A}_n
  $,
	and that
  \begin{align}
	  {\bf \tilde{r}}_n
    &=
		(r_n-\hat{r}_n) V(t_n,x_n), 
	\\
    {\bf p}_n
    &=
    n I_\mathbb{H} \left( x_n - \hat{x}_n \right) 
    V( t_n, x_n )
    + 
    r_n \, ( D_\mathbb{H} V )(t_n, x_n ) ,
  \\
    {\bf \hat{p}}_n
    &=
    n I_\mathbb{H} \left( x_n - \hat{x}_n \right) 
    V( \hat{t}_n, \hat{x}_n )
    + 
    \hat{r}_n \, 
    ( D_\mathbb{H} V )( \hat{t}_n, \hat{x}_n )
    ,
  \\
		{\bf \tilde{p}}_n
    &=
    (r_n-\hat{r}_n) \, 
    ( D_\mathbb{H} V )( t_n, x_n )
    ,
  \\
	\begin{split}
			{\bf A}_n
		& =
			 n 
			 (I_\mathbb{H} \left( 
				 x_n - \hat{x}_n
			 \right)) 
			\otimes ( D_\mathbb{H} V )( t_n , x_n )
			 +
			 (
				 D_\mathbb{H} 
				 V
			 )( t_n, x_n )
			\otimes
			 \left( 
				 n I_\mathbb{H} \left( 
				 x_n - \hat{x}_n 
			 \right) \right) 
		\\ &
			 +
			 r_n \, 
			 ( D^2_\mathbb{H} V)( t_n, x_n ) 
			 ,
	\end{split}
  \\
	\begin{split}
			{\bf \hat{A}}_n
		& =
			 n 
			 (I_\mathbb{H} \left( 
				 x_n - \hat{x}_n
			 \right))
			 \otimes ( D_\mathbb{H} V )( \hat{t}_n , \hat{x}_n )
			 +
			 (
				 D_\mathbb{H} 
				 V
			 )( \hat{t}_n , \hat{x}_n )
			 \otimes 
			 \left(n I_\mathbb{H} \left( 
				 x_n - \hat{x}_n 
			 \right) \right)
		\\ &
			 +
			 \hat{r}_n \, 
			 ( D^2_\mathbb{H} V)( \hat{t}_n , \hat{x}_n )
			,
	\end{split}
	\\ 
    {\bf \tilde{A}}_n
  & =
     (r_n-\hat{r}_n) \, 
     ( D^2_\mathbb{H} V)( t_n , x_n ). 
  \end{align}
  Moreover the assumption
	$V >1$ 
	together with the fact that
	 $ 
    \lim_{ n \to \infty }
    \big(
      \sqrt{ n } \,
      \|
        x_n 
        -
        \hat{x}_n 
      \|_H
    \big)
		+ \lim_{ n \to \infty }
    \big(
      \sqrt{ \beta_n } \,
      |
        t_n 
        -
        \hat{t}_n 
      |
    \big)
    = 0
  $,
	the Lipschitz continuity  of $ V $ 
	with respect to the $\| \cdot \|_{\R \times \mathbb{H}}$-norm
	on $\R \times \mathbb{H}$-bounded subsets of $(0,T) \times O$,
  and with
  the fact 
  $
    0 < 
    \lim_{ n \to \infty }
    \left( 
      r_n - \hat{r}_n 
    \right) 
  \leq
    \sup_{ n \in \N }
    \big(
      | r_n |
      +
      | \hat{r}_n |
    \big)
    < \infty
  $
  then imply that
  \begin{align}
	\nonumber
				0
			<{}&
				\lim_{ n \to \infty }
				\left( 
					r_n - \hat{r}_n 
				\right) 
				\inf_{(t,x) \in (0,T) \times O} V(t,x) \\ 
	\label{eq:lim_r}
			\leq{}
				&\liminf_{ n \to \infty }
				\big(
					(r_n-\hat{r}_n) V( t_n, x_n )
					- \hat{r}_n (V( \hat{t}_n, \hat{x}_n ) -V( t_n, x_n ))
				\big) \\ \nonumber
			={}&
				\liminf_{ n \to \infty }
				\big(
					r_n V( t_n, x_n )
					- \hat{r}_n V( \hat{t}_n, \hat{x}_n )
				\big)
			=
				\liminf_{ n \to \infty }
				\big(
					{\bf r}_n
					- {\bf \hat{r}}_n
				\big)
			\leq
				\sup_{ n \in \N }
				\left(
					| {\bf r}_n | + | {\bf \hat{r}}_n |
				\right) \\ \nonumber
			={}&
				\sup_{ n \in \N }
				\left(
					|r_n| |V(t_n,x_n)| + | \hat{r}_n | |V(\hat{t}_n,\hat{x}_n)|
				\right)
			<
				\infty.
  \end{align}
Thus
assumption \eqref{eq:comparison.viscosity.solution.assumption}
shows that
  \begin{equation} 
		\begin{split}
				&\limsup_{ n \to \infty }
				\left(
					\tfrac{
						G_{\mathbb{H}, \mathbb{X}, \delta_n, Vh}^+( 
							({\bf t}_n, {\bf x}_n), 
							{\bf r}_n, 
							{\bf p}_n, 
							{\bf A}_n
							+ n {\bf B}_n
							V( {\bf t}_n, {\bf x}_n ) 
						)
					}
					{V( {\bf t}_n, {\bf x}_n ) }
					- \tfrac{
						G_{\mathbb{H}, \mathbb{X}, \delta_n, Vh}^-( 
							({\bf \hat{t}}_n, {\bf \hat{x}}_n), 
							{\bf \hat{r}}_n, 
							{\bf \hat{p}}_n, 
							{\bf \hat{A}}_n
							+ n {\bf \hat{B}}_n
							V( {\bf \hat{t}}_n, {\bf \hat{x}}_n ) 
						)
					}
					{ V( {\bf \hat{t}}_n, {\bf \hat{x}}_n )}
				\right) \\
			\leq{}& 
				\limsup_{ n \to \infty } \left (
					\tfrac
						{G( (t_n, x_n), {\bf \tilde{r}}_n, {\bf \tilde{p}}_n, ({\bf \tilde{A}}_n|_X) |_X)}
						{V( t_n, x_n ) }
				\right ).
		\end{split}
  \end{equation}  
	In addition, \eqref{eq: tilde G + and G +}, \eqref{eq: tilde G - and G -},
	$\lim_{n \to \infty} \delta_n (h(t_n, x_n) \vee h(\hat{t}_n, \hat{x}_n)) = 0$
  the uniform continuity 
	with respect to the $\| \cdot \|_{\R \times \mathbb{H}}$-norm
	of $\frac{\partial}{\partial t} V$
	and of $V$ 
	on 
	$\R \times \mathbb{H}$-bounded subsets of $(0,T) \times O$,
	the definition of $ \tilde{G} $,
	and $\lim_{n \to \infty} r_n - \hat{r}_n > 0$ hence 
  imply that
  \begin{align}  
  &
    \limsup_{ 
      n \to \infty 
    }
    \Big(
       \tilde{G}_{\mathbb{H}, \mathbb{X}, \delta_n, h}^+\big(
         (t_n , x_n) ,
         r_n ,
         n \, I_\mathbb{H}
         (
           x_n - 
           \hat{x}_n
         ) ,
         n \,
         A_n
       \big) 
	\\ \nonumber & \qquad\qquad
       -
       \tilde{G}_{\mathbb{H}, \mathbb{X}, \delta_n, h}^-\big(
         (\hat{t}_n, \hat{x}_n),
         \hat{r}_n ,
         n \, I_\mathbb{H}
         (
           x_n -
           \hat{x}_n 
         ) ,
         n \,
         \hat{A}_n
       \big) 
    \Big)
  \\ \nonumber ={}&
    \limsup_{ 
      n \to \infty 
    }
    \Big(
      \tfrac{
        G_{\mathbb{H}, \mathbb{X}, \delta_n, Vh}^+( 	
          (t_n, x_n), {\bf r}_n ,  
          {\bf p}_n ,
          {\bf A}_n + n {\bf B}_n  V( t_n , x_n )
        )
        -
        (r_n+ \delta_n h(t_n, x_n))
        \frac{ \partial }{ \partial t } V( t_n, x_n )
      }{
        V( t_n, x_n )
      } \\ \nonumber
     & \qquad \qquad \qquad
			-
      \tfrac{
        G_{\mathbb{H}, \mathbb{X}, \delta_n, Vh}^-( 
          (\hat{t}_n, \hat{x}_n), {\bf \hat{r}}_n ,  
          {\bf \hat{p}}_n ,
          {\bf \hat{A}}_n 
          +
          n {\bf \hat{B}}_n V( \hat{t}_n, \hat{x}_n )
        )
        -
        (\hat{r}_n- \delta_n h(\hat{t}_n, \hat{x}_n))
        \frac{ \partial }{ \partial t } 
        V( \hat{t}_n, \hat{x}_n )
      }{
        V( \hat{t}_n, \hat{x}_n )
      }
    \Big)
  \\ \nonumber \leq{}&
    \limsup_{ n \to \infty } \Big (
			\tfrac
				{G( (t_n, x_n), {\bf \tilde{r}}_n, {\bf \tilde{p}}_n, ({\bf \tilde{A}}_n |_X) |_X )}
				{V( t_n, x_n ) }
			- \tfrac{(r_n-\hat{r}_n) \tfrac{ \partial }{ \partial t } V( t_n, x_n ) }
					{ V( t_n, x_n ) }
			+\hat{r}_n \left (
				\tfrac {\frac{ \partial }{ \partial t } V( \hat{t}_n, \hat{x}_n )}
					{V( \hat{t}_n, \hat{x}_n )}
				- \tfrac{ \frac{ \partial }{ \partial t } V( t_n, x_n )}
						{ V( t_n, x_n ) }
			\right ) \\ \nonumber
		& \qquad \qquad \qquad
			- \delta_n h( t_n, x_n )
					\tfrac{\frac{ \partial }{ \partial t } V( t_n, x_n )} {V( t_n, x_n )}
			- \delta_n h( \hat{t}_n, \hat{x}_n )
					\tfrac
						{\frac{ \partial }{ \partial t } V( \hat{t}_n, \hat{x}_n )} 
						{V( \hat{t}_n, \hat{x}_n )}
		\Big )
  \\  \nonumber
    ={}&
			\limsup_{ n \to \infty } \Big(
				\tfrac{
							-\frac{ \partial }{ \partial t }
								(r_n - \hat{r}_n)V( t_n, x_n )
				}
				{	V( t_n, x_n ) }
		\\ \nonumber & \qquad\qquad
				+\tfrac{
								G( 
								(t_n, x_n), 
								(r_n-\hat{r}_n) V( t_n, x_n ), 
								(r_n-\hat{r}_n) ( D_\mathbb{H} \tilde{V} )( t_n, x_n ),
								(r_n-\hat{r}_n) (( D^2_\mathbb{H} \tilde{V} )( t_n, x_n ) |_X) |_X
							)
				}
				{	V( t_n, x_n ) }
			\Big)
    \leq 0
  \end{align}  
  as $\lim_{n \to \infty} r_n-\hat{r}_n > 0$ and
	$(r_n-\hat{r}_n) V $ is by assumption for all $n \in \N$ with
	$r_n-\hat{r}_n>0$ a classical supersolution
  of \eqref{eq:parabolic.equation.unbounded.case}.
  We can thus apply Corollary~\ref{c:technical.lemma.uniqueness}
  to obtain that
	for all $\tilde{T} \in (0,T)$  
	it holds that
  \begin{equation}
	\label{eq: equation LEM 4.10}
		\begin{split}
			\lim_{R \to \infty} \lim_{r \downarrow 0} \lim_{\eps \downarrow 0}
				\sup \bigg \{
					&\tilde{u}_1(t,x) - \tilde{u}_2(\hat{t} , \hat{x}) \colon 
					~(t, x), (\hat{t}, \hat{x}) \in W,
					~ t, \hat{t} \leq \tilde{T},
					~h(t, x) \vee h(\hat{t}, \hat{x}) \leq R,
			\\
					~&\|x- \hat{x}\|_H \leq r,
					~|t-\hat{t}| \leq \eps 
				\bigg \}
			\leq
				0.
		\end{split}
  \end{equation}
	In the last step we need to show that \eqref{eq: equation LEM 4.10} implies 
	Theorem~\ref{l:comparison.viscosity.solution}.
	Therefore, note that \eqref{eq: V bound}
	shows that for all
	$\tilde{R} \in (0, \infty)$
	$x$, $\hat{x} \in O$ and all $t$, $\hat{t} \in (0, T)$
	with 
	$\|x\|_H \leq \tilde{R}$,
	$\|x- \hat{x}\|_H \leq \tilde{R}$, and with
	$\frac{1}{V(\hat{t}, \hat{x})} - \frac{1}{V(t, x)} \geq 0$ it holds that
	\begin{equation}
		\begin{split}
				u_1(t,x) - u_2(\hat{t}, \hat{x})
			\leq{}& 
				\left ( \frac{C_{2\tilde{R}}}{V(\hat{t}, \hat{x})} 
					(u_1(t,x) - u_2(\hat{t}, \hat{x})) \right)
				\vee 0 \\
			\leq{}
				&C_{2\tilde{R}} \left ( 
					\frac{u_1(t,x)}{V(t, x)}
					-\frac{u_2(\hat{t}, \hat{x})}{V(\hat{t}, \hat{x})}
					-\frac{u_1(t,x)}{V(t, x)}
					+\frac{u_1(t,x)}{V(\hat{t}, \hat{x})}
				\right)
				\vee 0 \\
			\leq{}&
				\left(
					C_{2\tilde{R}} (\tilde{u}_1(t,x)- \tilde{u}_2(\hat{t}, \hat{x}))
					+ C_{2\tilde{R}} \left( \frac{1}{V(\hat{t}, \hat{x})} - \frac{1}{V(t, x)} \right)
				\right)
				\vee 0.
		\end{split}
	\end{equation}
	On the other hand, \eqref{eq: V bound} implies that for all
	$\tilde{R} \in (0, \infty)$
	$x$, $\hat{x} \in O$ and all $t$, $\hat{t} \in (0, T)$
	with 
	$\|x\|_H \leq \tilde{R}$,
	$\|x- \hat{x}\|_H \leq \tilde{R}$, and with
	$\frac{1}{V(t, x)} - \frac{1}{V(\hat{t}, \hat{x})} \geq 0$ it holds that
	\begin{align}
		\nonumber
				u_1(t,x) - u_2(\hat{t}, \hat{x})
			\leq{}& 
				\left ( \frac{C_{2\tilde{R}}}{V(t, x)} (u_1(t,x) - u_2(\hat{t}, \hat{x})) \right)
				\vee 0 \\
			\leq{}
				&C_{2\tilde{R}} \left ( 
					\frac{u_1(t,x)}{V(t, x)}
					-\frac{u_2(\hat{t}, \hat{x})}{V(\hat{t}, \hat{x})}
					+\frac{u_2(\hat{t}, \hat{x})}{V(\hat{t}, \hat{x})}
					-\frac{u_2(\hat{t}, \hat{x})}{V(t, x)}
				\right)
				\vee 0 \\ \nonumber
			\leq{}&
				\left(
					C_{2\tilde{R}} (\tilde{u}_1(t,x)- \tilde{u}_2(\hat{t}, \hat{x}))
					+ C_{2\tilde{R}} \left( \frac{1}{V(t, x)} - \frac{1}{V(\hat{t}, \hat{x})} \right)
				\right)
				\vee 0.
	\end{align}
	Combining this inequalities yields then that for all
	$\tilde{R} \in (0, \infty)$
	$x$, $\hat{x} \in O$ and all $t$, $\hat{t} \in (0, T)$
	with 
	$\|x\|_H \leq \tilde{R}$ and with
	$\|x- \hat{x}\|_H \leq \tilde{R}$
	it holds that
	\begin{equation}
	\label{eq: u_1 - u_2 bound}
		\begin{split}
				u_1(t,x) - u_2(\hat{t}, \hat{x})
			\leq 
				\left(
					C_{2\tilde{R}} (\tilde{u}_1(t,x)- \tilde{u}_2(\hat{t}, \hat{x}))
					+ C_{2\tilde{R}} \left| \frac{1}{V(t, x)} - \frac{1}{V(\hat{t}, \hat{x})} \right|
				\right)
				\vee 0.
		\end{split}
	\end{equation}
	Finally we obtain with \eqref{eq: equation LEM 4.10},
	\eqref{eq: u_1 - u_2 bound}, and with
	the fact that $\frac 1V$ is uniformly continuous 
	with respect to the $\| \cdot \|_{\R \times \mathbb{H}}$-norm
	on $\R \times \mathbb{H}$-bounded subsets of $(0,T) \times O$
	that for all $R, \tilde{R} \in (0, \infty)$ and all $\tilde{T} \in (0, T)$ it holds that
	\begin{align}
				&\lim_{r \downarrow 0} \lim_{\eps \downarrow 0}
					\sup \bigg \{
						u_1(t,x) - u_2(\hat{t} , \hat{x}) \colon 
						(t, x), (\hat{t}, \hat{x}) \in W,
						~t, \hat{t} \leq \tilde{T},
						~h(t, x) \vee h(\hat{t}, \hat{x}) \leq R,
				~\\ \nonumber & \qquad \qquad \quad
						~\|x\|_H \leq \tilde{R},
						~\|x- \hat{x}\|_H \leq r, 
						~|t-\hat{t}| \leq \eps 
					\bigg \} \\ \nonumber
			\leq{}&
				\lim_{r \downarrow 0} \lim_{\eps \downarrow 0}
					\sup \bigg \{
						\bigg (
							C_{2\tilde{R}} (\tilde{u}_1(t,x)- \tilde{u}_2(\hat{t}, \hat{x}))
							+ C_{2\tilde{R}} \left| \frac{1}{V(t, x)} - \frac{1}{V(\hat{t}, \hat{x})} \right|
						\bigg ) \colon 
						(t, x), (\hat{t}, \hat{x}) \in W,
					~\\ \nonumber & \qquad \qquad \quad
						~t, \hat{t} \leq \tilde{T}, 
						~h(t, x) \vee h(\hat{t}, \hat{x}) \leq R,
						~\|x\|_H \leq \tilde{R},
						~\|x- \hat{x}\|_H \leq r, 
						~|t-\hat{t}| \leq \eps 
					\bigg \} \vee 0\\ \nonumber
			\leq{}&
				\Bigg( 
					C_{2\tilde{R}} \lim_{r \downarrow 0} \lim_{\eps \downarrow 0}
						\sup \bigg \{
							 \tilde{u}_1(t,x)- \tilde{u}_2(\hat{t}, \hat{x}) \colon 
							(t, x), (\hat{t}, \hat{x}) \in W,
							~t, \hat{t} \leq \tilde{T}, 
							~h(t, x) \vee h(\hat{t}, \hat{x}) \leq R,
						~\\ \nonumber & \qquad \qquad \quad
							~\|x\|_H \leq \tilde{R}, 
							~\|x- \hat{x}\|_H \leq r,
							~|t-\hat{t}| \leq \eps 
						\bigg \} \\ \nonumber
					&+
					C_{2\tilde{R}} \lim_{r \downarrow 0} \lim_{\eps \downarrow 0}
						\sup \bigg \{
							\left| \frac{1}{V(t, x)} - \frac{1}{V(\hat{t}, \hat{x})} \right| \colon 
							(t, x), (\hat{t}, \hat{x}) \in W,
							~t, \hat{t} \leq \tilde{T},
							~h(t, x) \vee h(\hat{t}, \hat{x}) \leq R,
					~\\ \nonumber & \qquad \qquad \quad
							~\|x\|_H \leq \tilde{R}, 
							~\|x- \hat{x}\|_H \leq r,
							~|t-\hat{t}| \leq \eps 
						\bigg \} 
				\Bigg)  \vee 0 \\ \nonumber
			\leq{}& 0.
	\end{align}
  This finishes the proof
  of Theorem~\ref{l:comparison.viscosity.solution}.
\end{proof}
\section{Uniqueness of viscosity solutions of Kolmogorov equations}
\label{sec: Uniqueness of Kolmogorov equations}
Next we apply Theorem 
\ref{l:comparison.viscosity.solution}
to Kolmogorov equations of SPDEs.
This generalizes  Corollary 4.14 in 
Hairer, Hutzenthaler \& Jentzen \cite{HairerHutzenthalerJentzen2015}
(which assumed finite-dimensional Hilbert spaces and
classical viscosity solutions).
\begin{corollary}[Uniqueness
of viscosity solutions of Kolmogorov type equations]
\label{cor:uniqueness2}
  Let 
  $ T $, $ \theta \in (0, \infty)$,
	let $\mathbb{U}=(U, \langle \cdot, \cdot \rangle_{U}, \| \cdot \|_U)$
	and $\mathbb{H}=(H, \langle \cdot, \cdot \rangle_{H}, \| \cdot \|_H)$
	be real separable Hilbert spaces,  
	let
	$\mathbb{H}'=(H', \| \cdot \|_{H'})$ be the dual space,
	let $(e_i)_{i \in \N} \subseteq H$ be an orthonormal basis of $\mathbb{H}$,
	let $\lambda \colon \N \to (- \infty, 0)$ be a function with 
	$\lim_{n \to \infty} \lambda_n = - \infty$,
	let $A \colon D(A) \subseteq H \to H$ be the linear operator such that
	$
			D(A) 
		= 
			\{ 
				x \in H \colon 
					\sum_{n =1}^\infty |\lambda_n \langle e_n, x \rangle_H|^2 < \infty 
			\}
	$
	and such that for all $x \in D(A)$ it holds that
	$
		Av = \sum_{n =1}^\infty \lambda_n \langle e_n, x \rangle_H e_n
	$,
	let $\mathbb{H}_r=(H_r, \langle \cdot, \cdot \rangle_{H_r}, \| \cdot \|_{H_r})$,
	$r \in \R$, be a family of interpolation spaces associated to $-A$,
	(see e.g., Definition 3.6.30 in Jentzen \cite{Jentzen2015}),
	and let $\mathbb{H}'_r=(H'_r, \| \cdot \|_{H'_r})$,
	$r \in \R$,
	be the corresponding dual spaces.
	By abuse of notation we will also denote by
		$A$ and  by $\| \cdot \|_{H_r}$, $r \in \R$,
		the extended functions
		$
			A \colon \bigcup_{i=1}^{\infty} H_{-i} \to \bigcup_{i=1}^{\infty} H_{-i}
		$ 
		and 
		$
			\| \cdot \|_{H_r} \colon \bigcup_{i=1}^{\infty} H_{-i} \to [0, \infty]
		$,
		$r \in \R$
		satisfying for all
		$r \in \R$, $x \in H_{r}$,
		and all $y \in \bigcup_{i=1}^{\infty} H_{-i}$
		that
		\begin{equation}
					\|y\|_{H_r}
				=
					\begin{cases} 
						\|y\|_{H_{r}} 
							& \textrm{ if } z \in H_{r} \\
						\infty 
							& \textrm{ if }z \notin H_{r}
					\end{cases}
		\end{equation}
		and that
		\begin{equation}
				(A(x) = y) 
			\Leftrightarrow 
				(
						\lim_{\eps \downarrow 0} \sup \{
							\|A(\xi) - y\|_{H_{r-1}} \colon \xi \in D(A), ~\|x- \xi\|_{H_{r}} \leq \eps 
						\}
					=
						0
				),
		\end{equation}
	let
  $ O \subseteq H $ 
  be an open and convex set,
	let $ O_n, O^c_n \subseteq O $, $ n \in \N $
	be the sets satisfying for all $n \in \N$ that
	$
			O_n
		=
			\left\{
				x \in O \colon
				\dist_{\mathbb{H}}(x, H \backslash O)  
				\ge \tfrac{ 1 }{ n }
				\text{ and }
				\left\| x \right\|_H \le n
			\right\},
	$
	and that $ O_n^c = O \backslash O_n $,
	let
  $
    \varphi \in \C_{\R \times \mathbb{H}} ((0,T) \times O, \R )
  $,
	let $K \colon [0, \infty) \to (0, \infty)$ be increasing,
	let $\omega \colon [0, \infty) \to [0, \infty)$
	satisfy that 
		$\lim_{x \downarrow 0} \omega (x) = 0$
	and that
		$\sup_{x \in [0,\infty)} \omega(x) < \infty$,
	let $\alpha_1, \alpha_2, \beta_1, \beta_2 \in (0, \infty)$,
	$\var \in [1/2,\infty)$,
	$\kappa\in (-1/2, 0]$,
	let
	$ 
    v \colon (0,T) \times O \to \R  
  $,
  $ 
    F \colon 
    (0,T) \times O \to H_{-1}
  $,
  and let
  $ 
    B \colon 
    [0,T] \times O \to 
    HS(\mathbb{U}, \mathbb{H}_{\kappa})
  $
	satisfy that
	for all  
	$t, \hat{t} \in (0,T)$ 
	and all $x, \hat{x} \in H_{2\var}$ 
	it holds that
	\begin{align}
	\label{eq: H bound for F, B and v}
			&\| F (t,x)\|_{H_{\var-1}}
		\leq 
			K(\|x\|_{H_\var}), \quad
			\| B(t,x) \|_{HS(\mathbb{U},\mathbb{H}_\kappa)}
		\leq
			K(\|x\|_{H}), \quad
			|v(t,x)|
		\leq K(\|x\|_{H}) \\
	\label{eq: H_3/2 bound for F}
			&\| F (t,x)\|_{H_{\var - 1/2}}
		\leq 
			K(\|x\|_{H}) (\|x\|_{H_{\var +1/2 - \alpha_1}} +1),
	\\
	\label{eq: H_3/2 bound for B}
			&\| B(t,x) \|^2_{HS(\mathbb{U},\mathbb{H}_\var)} 
		\leq
			K(\|x\|_{H}) (\|x\|^2_{H_{\var + 1/2- \beta_1}} + 1), \\
	\label{eq: Lipschitz continuity of F}
			&\|F(t,x) - F(\hat{t},\hat{x})\|_{H_{-1/2}}
		\leq
			K(\|x\|_{H} \vee \|\hat{x}\|_{H}) 
				(\|x-\hat{x}\|_{H_{1/2 - \alpha_2}} + |t - \hat{t}|), \\
	\label{eq: Lipschitz continuity of B}
			&\|B(t,x) - B(\hat{t},\hat{x})\|_{HS(\mathbb{U}, \mathbb{H})}
		\leq
			K(\|x\|_{H} \vee \|\hat{x}\|_{H}) 
				(\|x-\hat{x}\|_{H_{1/2 - \beta_2}} + |t - \hat{t}|), \\
	\label{eq: uniform continuity of v}
			&|v(t,x) - v(\hat{t},\hat{x})|
		\leq
			K(\|x\|^2_{H_\var} \vee \|\hat{x}\|^2_{H_\var}) \cdot
			(\|x- \hat{x}\|^2_{H_{1/2}} + \omega(\|x- \hat{x}\|^2_{H_{1/2}} + |t -\hat{t}|^2)),
	\end{align}
	let 
	$
		\| \cdot \|_{\boldsymbol{L}(\mathbb{H}_{\kappa},\mathbb{H}'_{\kappa})} 
			\colon L(\mathbb{H}, \mathbb{H}') \to [0, \infty]
	$ 
	and
	$
		\| \cdot \|_{\boldsymbol{L}(\mathbb{H}_{\kappa},\mathbb{H}_{-\kappa})} 
			\colon L(\mathbb{H}, \mathbb{H}) \to [0, \infty]
	$
	be the extended norms satisfying for all $C \in L(\mathbb{H}, \mathbb{H}')$
	that 
	\begin{equation}
	\label{eq: def of bold L}
			\| I^{-1}_{\mathbb{H}} C \|_{\boldsymbol{L}(\mathbb{H}_{\kappa},\mathbb{H}_{-\kappa})}
		=
			\| C \|_{\boldsymbol{L}(\mathbb{H}_{\kappa},\mathbb{H}'_{\kappa})}
		=
			\sup_{x,y \in H \backslash \{0 \}} 
				\frac{|\langle x, Cy \rangle_{H, H'}|}{\|x\|_{H_{\kappa}}\|y\|_{H_{\kappa}}},
	\end{equation}
  and
	let $ V \in \C_{\mathbb{H}}^2( O, (1,\infty) ) $, 
	satisfy
	that $D_{\mathbb{H}}^2 \, V \colon H \to L(\mathbb{H}, \mathbb{H}')$ 
	is uniformly continuous with respect to the
	$\| \cdot \|_{\mathbb{H}}$ and the
	$\| \cdot \|_{\boldsymbol{L}(\mathbb{H}_{\kappa},\mathbb{H}'_{\kappa})}$-norm
	on all $\mathbb{H}$-bounded subsets of $O$,
	that for all $ (t,x) \in (0,T) \times O \cap H_{2 \varz} $ it holds that
  \begin{equation}
  \label{eq:assumption_V}
		\begin{split}
			&v(t,x) \, V(x)
			+
			\langle
				F(t,x) + A(x), (D_\mathbb{H} V)(x)
			\rangle_{H, H'}
	\\ & \qquad \qquad\qquad\qquad
			+\nicefrac 12
			\langle 
				B(t,x),
				I_{\mathbb{H}_\var}^{-1} (D_\mathbb{H}^2 \, V)(x) \, B(t,x)
			\rangle_{HS(\mathbb{U},\mathbb{H}_\var)} 
			\leq \theta \cdot V(x),
		\end{split}
  \end{equation}
	that for all  
	$x$, $\hat{x} \in H_{2 \var}$ 
	it holds that
	\begin{align}
	\label{eq: V continuous with respect to H1}
			&\bigg \|
				\frac{I_\mathbb{H}^{-1} (D_\mathbb{H} V)(x)}{V(x)} 
				- \frac{I_\mathbb{H}^{-1} (D_\mathbb{H} V)(\hat{x})}{V(\hat{x})}
			\bigg \|_{H_{1-\var}} 
		\leq{}
			K(\|x\|^2_{H_\var} \vee \|\hat{x}\|^2_{H_\var}) 
			\, (\|x- \hat{x}\|^2_{H_{1/2}} + \omega(\|x- \hat{x}\|^2_{H_{1/2}})), \\
	\label{eq: V Lipschitz with respect to H -gamma}
			&\bigg \|
				\frac
					{I_{\mathbb{H}}^{-1} ( D_\mathbb{H} V)( x )}
					{V (x )}
				-\frac
					{I_{\mathbb{H}}^{-1} (D_\mathbb{H} V )( \hat{x} )}
					{V( \hat{x} )}
			\bigg \|_{H_{-\kappa}}
		\leq 
			K(\|x\|_{H} \vee \| \hat{x} \|_H) \cdot \|x - \hat{x}\|_{H_{1/2}}, \\
		\label{eq: V bounded with respect to H}
			&\left \|
				I_\mathbb{H}^{-1}(D_\mathbb{H} V)(x) 
			\right \|_{H_{-\kappa}}
		\leq
			K(\|x\|_{H}),
	\end{align}
	and	that 
	\begin{equation}
	\label{eq: V bounded with respect to H 1/2}
			\left \|
				I_\mathbb{H}^{-1}(D_\mathbb{H} V)(x) 
			\right \|_{H_{1/2}}
		\leq
			K(\|x\|^2_{H_\var}).
	\end{equation}
  Then there exists 
  at most one
	with respect to the
	$\| \cdot \|_{\R \times H}$-norm
  continuous function
  $ 
    u 
    \colon [0,T] \times O \to \R 
  $
	which is 
	bounded on $\mathbb{H}$-bounded subsets of $(0,T) \times O$,
	which fulfills for all $R, \tilde{R} \in (0, \infty)$ that
	\begin{equation}
		\begin{split}
				\lim_{r \downarrow 0} \lim_{\eps \downarrow 0}	
				\sup \Bigg \{
					&u(t,x) -  \varphi(\hat{t}, \hat{x})
					\colon
					~x, \hat{x} \in H_{2 \var} \cap O,
					~t, \hat{t} \in (0,T),
					~\tfrac{e^{-\theta t} \|x\|_{H_\var}^2}{V(x)} 
							\vee \tfrac{e^{-\theta \hat{t}} \|\hat{x}\|_{H_\var}^2}{V(\hat{x})} \leq R, \\
					~&\|x \|_H \leq \tilde{R}, 
					~\| x - \hat{x}\|_H \leq r,
					~t \vee \hat{t} \leq \eps 
				\Bigg \}
			= 
				0,
		\end{split}
	\end{equation}
  which fulfills
  $
    \lim_{ n \to \infty }
      \sup_{
        \substack{
          (t,x) \in [0,T] \times
          O_n^c
        }
      }
      \frac{
        | u(t, x) |
      }{
        V(x)
      }
  $
  $
    = 0 
  $,
  and which fulfills that
  $ 
    u|_{ (0,T) \times O } 
  $
  is a viscosity solution of
  \begin{equation} 
  \label{eq:second-order.PDE}
		\begin{split}
			&\tfrac{ \partial }{ \partial t }
			u(t,x) -
			v(t,x) \, u(t,x)
			-
			\big\langle
				F(t,x) + A(x), (D_\mathbb{H} u)(t,x)
			\big\rangle_{H, H'} 
			\\ & \qquad \qquad\qquad \qquad
			-\big\langle 
				B(t,x),
				I_{\mathbb{H}_\var}^{-1} \big( 
					(D^2_\mathbb{H} \, u)(t, x) \, B(t,x) 
				\big) \big |_{H_\var}
			\big \rangle_{HS(\mathbb{U},\mathbb{H}_\var)} 
			= 0
		\end{split}
  \end{equation}
	for $ (t,x) \in (0,T) \times O $ relative to 
	$
		((0,T) \times O \ni (t,x) \to \|x \|^2_{H_\var} \in \R \cup \{\infty \},
		\R \times \mathbb{H}, \R \times \mathbb{H}_\var)
	$.
\end{corollary}
\begin{proof}[Proof 
of Corollary~\ref{cor:uniqueness2}]
	First
  denote for every $N \in \N$ 
	by $V_{N} = \Span_{\mathbb{H}} (\{ e_1, ... , e_N\})$ 
	the linear span of the first $N$ basis vectors,
	by $W \subseteq (0,T) \times O$ the set satisfying that
	$
			W
		=
			(0,T) \times (H_{2\var} \cap O)
	$,
  by
  $
    G \colon W \times
    \R \times H' \times
    \mathbb{S}_{\mathbb{H}_{\var}, \mathbb{H}_{\var}'}
    \to \R 
  $
  the function satisfying for all
	$
		((t, x), r, p, C) 
			\in W \times R \times H' \times \mathbb{S}_{\mathbb{H}_\var, \mathbb{H}_\var'}
	$ 
	that
  \begin{equation}
    G( (t, x), r, p, C)
    =
    v(t,x) r
    +
    \left\langle
      F(t,x) +A(x), p
    \right\rangle_{H, H' }
    +
    \langle
      B(t,x) ,
      I_{\mathbb{H}_\var}^{-1} C \, B(t,x)  
    \rangle_{HS(\mathbb{U}, \mathbb{H}_\var)},
  \end{equation}
	by $\tilde{V} \colon [0,T] \times O \to (1, \infty)$
	the function satisfying for all $(t,x) \in [0,T] \times O$ that
	\begin{equation}
	\label{eq: def of tilde V}
		\tilde{V}(t,x)= e^{\theta t} \cdot V (x),
	\end{equation}
	and by $h \colon [0,T] \times O \to \R \cup \{ \infty \}$ the function
	satisfying for all $t \in [0,T]$ and all $x \in O$ that
	\begin{equation}
	\label{eq: def of h}
			h(t,x) 
		=  
				\frac {\|x\|^2_{H_\var}}{\tilde{V}(t,x)}. 
	\end{equation}
	Note that we obtain for all $(t,x) \in W$ that
	\begin{equation}
	\label{eq: h derivative}
		\begin{split}
				&E_{\mathbb{H}_\var', \mathbb{H}'} \big(
					(D_{\mathbb{H}_\var} (\tilde{V}h|_{\R \times H_\var}))(t,x)
				\big)
			=
				2 \cdot E_{\mathbb{H}_\var', \mathbb{H}'} \big( I_{\mathbb{H}_\var}x \big)
			=
			 2 \cdot E_{\mathbb{H}_\var', \mathbb{H}'} \big(
				H_\var \ni y \to \langle y , x \rangle_{H_\var} \in \R
			\big) \\
			={}&
			 2 \cdot E_{\mathbb{H}_\var', \mathbb{H}'} \big(
				H_\var \ni y \to \langle y , (-A)^{2\var} x \rangle_{H} \in \R
			\big)
			=
				2 \cdot I_{\mathbb{H}} ((-A)^{2\var}x).	
		\end{split}
	\end{equation}
	In addition, the assumption
	that $D_{\mathbb{H}}^2 \, V \colon H \to L(\mathbb{H}, \mathbb{H}')$ 
	is uniformly continuous with respect to the
	$\| \cdot \|_{\mathbb{H}}$ and the
	$\| \cdot \|_{\boldsymbol{L}(\mathbb{H}_{\kappa}, \mathbb{H}_{\kappa}')}$-norm
	on $\mathbb{H}$-bounded subsets of $O$ implies that
	$D_{\mathbb{H}}^2 \, V \colon H \to L(\mathbb{H}, \mathbb{H}')$ 
	is also uniformly continuous with respect to the
	$\| \cdot \|_{\mathbb{H}}$ and the
	$\| \cdot \|_{L(\mathbb{H}, \mathbb{H}')}$-norm
	on $\mathbb{H}$-bounded subsets of $O$
	and therefore $V$ is Lipschitz continuous
	with respect to the
	$\| \cdot \|_{\mathbb{H}}$-norm
	on $\mathbb{H}$-bounded subsets of $O$.
	Moreover,
	for all $x \in H_{\var}$ and all $x_0 \in H_{2 \var}$ it holds that
	$
			\|x \|_{H_{\var}}^2 
		= 
			\|x_0\|_{H_{\var}}^2 + 2 \langle x_0, x-x_0 \rangle_{H_{\var}} 
			+ \langle x-x_0, x-x_0 \rangle_{H_{\var}}
		\geq
			\|x_0\|_{H_{\var}}^2 
			+2\langle I_{\mathbb{H}_{\var}} 
					(x_0), x-x_0 
				\rangle_{H'_{\var},H_{\var}}.
	$
	This shows that 
	$
			\{ 
				y \in H_\var \colon 
					(J^{2}_{\mathbb{H}, -} \| \cdot \|^2_{H_\var}) (y) \neq \emptyset
			\} 
		=
			\{ 
				y \in H_\var \colon 
					I_{\mathbb{H}_{\var}} (y) \in D(E_{\mathbb{H}_\var', \mathbb{H}'})
			\}
		=
			H_{2\var}
	$.
	Combining this with \eqref{eq: h derivative} yields that
	$
			\{ 
				z \in (0,T) \times (O \cap H_{\var}) \colon 
					(D_{\R \times \mathbb{H}_\var} (h\tilde{V} |_{\R \times H_\var}))(z) \in
					D(E_{\R \times \mathbb{H}_\var',\R \times \mathbb{H}'}),
					~(J^{2}_{\mathbb{H}, -} (h \tilde{V})) (z) \neq \emptyset
			\} 
		=
			W
	$
	and this together with
	Proposition \ref{prop: h and Vh} 
	implies that 
	$
			\{ 
				z \in (0,T) \times (O \cap H_{\var}) \colon 
					(D_{\R \times \mathbb{H}_\var} (h |_{\R \times H_\var}))(z) \in
					D(E_{\R \times \mathbb{H}_\var',\R \times \mathbb{H}'}),
					~(J^{2}_{\mathbb{H}, -} h) (z) \neq \emptyset
			\} 
		=
			W
	$. 
	Next let 
  $
    u_1,
    u_2
    \colon [0,T] \times O \to \R
  $ 
  be two 
	with respect to the $\| \cdot \|_{\R \times H}$-norm continuous
  functions such that
	$u_1$, $u_2$
	are bounded
	on $\R \times \mathbb{H}$-bounded subsets of $ [0,T] \times O$,
	such that for all $R, \tilde{R} \in  (0, \infty)$ and $i \in \{1, 2 \}$
  it holds that
	\begin{equation}
	\label{eq: uniformly bounded at 0 with phi}
		\begin{split}
				\lim_{r \downarrow 0} \lim_{\eps \downarrow 0}	
				\sup \Bigg \{
					&u_i(t,x) -  \varphi(\hat{t}, \hat{x})
					\colon
					~x, \hat{x} \in H_{2\var} \cap O,
					~t, \hat{t} \in (0,T),
					~\tfrac{\|x\|_{H_{\var}}^2}{V(x)} 
						\vee \tfrac{\|\hat{x}\|_{H_{\var}}^2}{V(\hat{x})}  \leq R, \\
					&~\|x \|_H \leq \tilde{R}, 
					~\| x - \hat{x}\|_H \leq r,
					~t \vee \hat{t} \leq \eps 
				\Bigg \}
			= 
				0,
		\end{split}
	\end{equation}
  such that
  $
    \lim_{ n \to \infty }
      \sup_{
        \substack{
          (t,x) \in (0,T) \times
          O_n^c
        }
      }
      \frac{
        | u_1(t, x) |
        +
        | u_2(t, x) |
      }{
        V(x)
      }
  $
  $
    = 0 
  $
  and such that
  $
    u_1|_{ (0,T) \times O }
  $ 
  and
  $
    u_2|_{ (0,T) \times O }
  $ 
  are viscosity solutions 
  of~\eqref{eq:second-order.PDE}  relative to 
	$
		(
			h \tilde{V},
			\R \times \mathbb{H},
			\R \times \mathbb{H}_{\var}
		)
	$,
	let
	$\mathcal{I} \subseteq \N$,
	and let
	$(\tilde{e}_i)_{i \in \mathcal{I}}$ be an orthonormal basis
	of $\mathbb{U}$.
	We show Corollary~\ref{cor:uniqueness2} by applying 
  Theorem~\ref{l:comparison.viscosity.solution}
	with $\mathbb{X} \leftarrow \mathbb{H}_{\var}$
	and with $V \leftarrow \tilde{V}$.
	Note that we have for all $x \in H$ that
	\begin{equation}
		\begin{split}
				&(-A)^{2\var} I^{-1}_{\mathbb{H}_\var} \big( 
					(H \ni y \to \langle x, y \rangle_{H} \in \R)|_{H_\var}
				\big) 
			=
				(-A)^{2\var} I^{-1}_{\mathbb{H}_\var} \big( 
					H_\var \ni y \to \langle (-A)^{-2\var} x, y \rangle_{H_\var} \in \R
				\big) 
			=
				x \\
			=
				&I^{-1}_{\mathbb{H}} \big( 
					H \ni y \to \langle x, y \rangle_{H} \in \R
				\big)
		\end{split}
	\end{equation}
	and from this it follows that for all
	$(t,x) \in W$ and all 
	$C \in L(\mathbb{H},\mathbb{H}')$ it holds that
	\begin{equation}
	\label{eq: different HS norm}
		\begin{split}
				&\big \langle 
					B(t,x), I_{\mathbb{H}_\var}^{-1} \big( (C B(t,x))|_{H_\var} \big)
				\big \rangle_{HS(\mathbb{U},\mathbb{H}_\var)}
			=
				\sum_{ i\in \mathcal{I}}
					\big \langle 
						B(t,x) \tilde{e}_i, 
						I_{\mathbb{H}_\var}^{-1} \big( (C B(t,x) \tilde{e}_i)|_{H_\var} \big)
					\big \rangle_{H_\var} \\
			={}&
				\sum_{ i\in \mathcal{I}}
					\big\langle
						(-A)^{\var} B(t,x) \tilde{e}_i, 
						(-A)^{\var} I_{\mathbb{H}_\var}^{-1} \big( 
							(C B(t,x) \tilde{e}_i)|_{H_\var} 
						\big)
					\big\rangle_{H} \\
			={}&
				\sum_{ i\in \mathcal{I}} 
					\big \langle 
						B(t,x) \tilde{e}_i, 
						(-A)^{2\var} I_{\mathbb{H}_\var}^{-1} \big(
							(C B(t,x) \tilde{e}_i)|_{H_\var} 
						\big)
					\big \rangle_{H} \\
			={}&
				\sum_{ i\in \mathcal{I}} 
					\langle B(t,x) \tilde{e}_i, I_\mathbb{H}^{-1} C B(t,x) \tilde{e}_i \rangle_{H} 
			=
				\langle B(t,x), I_\mathbb{H}^{-1} C B(t,x) \rangle_{HS(\mathbb{U},\mathbb{H})}.
		\end{split}
	\end{equation}
	Furthermore, observe that for all $\eps \in (0, \infty)$
	there exist an $M_\eps \in \R$ such that
	for all $x \in H_{2 \var}$ it holds that
	\begin{equation}
	\label{eq: H1/2- norm bounds}
			\|x\|_{H_{1/2 - \alpha_2}}
			\vee \|x\|_{H_{1/2 - \beta_2}}
			\vee \|x\|_{H_{- \kappa}}
		\leq
			\eps \|x\|_{H_{1/2}} + M_\eps\|x\|_{H}
	\end{equation}
	and that
	\begin{equation}
	\label{eq: H3/2- norm bounds}
			\|x\|_{H_{\var + 1/2 - \alpha_1}}
			\vee \|x\|_{H_{\var + 1/2 - \beta_1}}
			\vee \|x\|_{H_{\var}}
		\leq
			\eps \|x\|_{H_{\var + 1/2}} + M_\eps\|x\|_{H}.
	\end{equation}
	In addition we get from \eqref{eq: h derivative},
	\eqref{eq: H_3/2 bound for F},
	the monotonicity of $K$, and from
	\eqref{eq: H3/2- norm bounds}
	that
	for all $(t, x) \in W$ and all $R, \eps, \delta \in (0, \infty)$
	with $\|x\|_H \leq R$
	it holds that
	\begin{equation}
	\label{eq: delta F equation}
    \begin{split}
				&
						\langle 
							F(t, x) + A(x),
							\delta 
							E_{\mathbb{H}_\var', \mathbb{H}'} \big(
								(D_{\mathbb{H}_\var} (\tilde{V}h|_{\R \times H_\var}))(t,x)
							\big)
						\rangle_{H, H'} \\
						&2\langle 
							F(t, x) + A(x),
							\delta I_{\mathbb{H}} (-A)^{2\var} (x)
						\rangle_{H, H'}
			\\ ={}&
						2\langle 
							F(t, x),
							\delta (-A)^{2\var} (x)
						\rangle_{H}
						+2\langle 
							A(x),
							\delta (-A)^{2\var} (x)
						\rangle_{H} \\
			\leq{}
						&2\delta \| F(t, x) \|_{H_{\var-1/2}} 
						\cdot \|x \|_{H_{1/2 + \var}}
						- 2\delta \| x \|^2_{H_{1/2+\var}}
			\\ \leq{}&
				2\delta \left(
						K(R) \cdot (\| x\|_{H_{\var + 1/2- \alpha_1}} + 1) 
							\|x\|_{H_{\var + 1/2}} -\|x\|^2_{H_{\var + 1/2}}
				\right)
			\\ \leq{}&
				2\delta \left(
						K(R) \cdot (\eps \, \|x\|_{H_{\var + 1/2}} 
						+ M_\eps \cdot \| x\|_{H} + 1) \|x\|_{H_{\var + 1/2}} 
						-\|x\|^2_{H_{\var + 1/2}}
				\right)
			\\ \leq{}&
				2\delta \left(
						K(R) \, (M_\eps \cdot R + 1) \|x\|_{H_{\var + 1/2}} 
						- (1- \eps \cdot K(R)) \|x\|^2_{H_{\var + 1/2}}
				\right).
		\end{split}
	\end{equation}
	Moreover, it follows from 
	\eqref{eq: def of h},
	\eqref{eq: H_3/2 bound for B},
	the monotonicity of $K$,
	\eqref{eq: H3/2- norm bounds},
	and from the fact that for all 
	$a,b \in \R$
	it holds that $(a+b)^2 \leq 2a^2 + 2b^2$
	that for all
	$(t, x) \in W$ and all $R, \eps, \delta \in (0, \infty)$
	with $\|x\|_H \leq R$
	it holds that
	\begin{equation} 
	\label{eq: delta B equation}
    \begin{split}
				&
						\langle 
							B(t, x), 
							\delta I_{\mathbb{H}_\var}^{-1}(D^2_{\mathbb{H}_\var} (h\tilde{V})(t, x))  
								\, B(t, x)
						\rangle_{HS(\mathbb{U},\mathbb{H}_\var)}
			=
						\langle 
							B(t, x), 
							\delta I_{\mathbb{H}_\var}^{-1} \, 2 I_{\mathbb{H}_\var}
								\, B(t, x)
						\rangle_{HS(\mathbb{U},\mathbb{H}_\var)}
			\\ ={}&
				2\delta \,
						\| B(t, x) \|^2_{HS(\mathbb{U},\mathbb{H}_\var)}
			\leq
				2\delta \,
						K(R) (\|x\|^2_{H_{\var + 1/2-\beta_1}} + 1)
		\\ \leq{}&
				2\delta \,
						K(R) (
							(\sqrt{\eps} \, \|x\|_{H_{\var + 1/2}} 
							+ M_{\sqrt{\eps}} \cdot \|x\|_{H})^2 + 1
						)
			\leq
				2\delta \,
						K(R) (
							2\eps \, \|x\|^2_{H_{\var + 1/2}} 
							+ 2 M^2_{\sqrt{\eps}} \cdot \|x\|^2_{H} + 1
						)
		\\ \leq{}&
				2\delta \,
						K(R) (2\eps \, \|x\|^2_{H_{\var + 1/2}} + 2 M^2_{\sqrt{\eps}} \cdot R^2 + 1).
		\end{split}
	\end{equation}
	Furthermore, \eqref{eq: different HS norm},
	\eqref{eq: H_3/2 bound for F},
	\eqref{eq: H_3/2 bound for B},
	the monotonicity of $K$,
	\eqref{eq: H3/2- norm bounds},
	and the fact that for all 
	$a,b \in \R$
	it holds that $(a+b)^2 \leq 2a^2 + 2b^2$
	imply
	that for all
	$(t, x) \in W$, $R, \eps \in (0, \infty)$,
	$p \in H'$, and all $C \in \mathbb{S}_{\mathbb{H}, \mathbb{H}'}$
	with $\|x\|_H \vee \|p\|_{H'} \vee \|C\|_{L(\mathbb{H}, \mathbb{H}')} \leq R$
	it holds that
	\begin{align}
	\nonumber
				& \big \langle F(t,x) + A(x), p \big \rangle_{H,H'} 
				+\big \langle 
					B(t,x), I^{-1}_{\mathbb{H}_\var} \big( (C B(t,x))|_{H_\var} \big)
				\big \rangle_{HS(\mathbb{U},\mathbb{H}_\var)} \\ \nonumber
			={}&
				\big \langle F(t,x) + A(x), p \big \rangle_{H,H'} 
				+\big \langle 
					B(t,x), I^{-1}_{\mathbb{H}} C B(t,x)
				\big \rangle_{HS(\mathbb{U},\mathbb{H})} \\ \nonumber
			\leq{}
				&(\|F(t,x)\|_{H} +\|x\|_{H_1}) \|p\|_{H'}
				+\|B(t,x)\|^2_{HS(\mathbb{U},\mathbb{H})}
				\|C\|_{L(\mathbb{H}, \mathbb{H}')} \\ \nonumber
			\leq{}&
				(
					\sup_{n \in \N} (|\lambda_n|^{1/2 - \var}) \|F(t,x)\|_{H_{\var -1/2}} 
					+\sup_{n \in \N} (|\lambda_n|^{1/2 - \var}) \|x\|_{H_{\var + 1/2}}
				) \cdot R
			\\  \nonumber &
				+\sup_{n \in \N} (|\lambda_n|^{- 2\var}) 
					\|B(t,x)\|^2_{HS(\mathbb{U},\mathbb{H}_\var)} \cdot R\\
	\label{eq: F + B equation}
		\begin{split}
			\leq{}&
				(
					\sup_{n \in \N} (|\lambda_n|^{1/2 - \var})
						K(R) (\|x\|_{H_{\var +1/2 -\alpha_1}} +1) 
					+\sup_{n \in \N} (|\lambda_n|^{1/2 - \var}) \, \|x\|_{H_{\var + 1/2}}
				) \cdot R 
			\\ &
				+\sup_{n \in \N} (|\lambda_n|^{- 2\var}) 
					K(R) (\|x\|^2_{H_{\var +1/2 -\beta_1}} +1) \cdot R 
		\end{split} \\ \nonumber
			\leq{}&
				(
					\sup_{n \in \N} (|\lambda_n|^{1/2 - \var})
						K(R) (\eps \|x\|_{H_{\var +1/2}} + M_\eps \|x\|_{H} +1) 
					+\sup_{n \in \N} (|\lambda_n|^{1/2 - \var}) \, \|x\|_{H_{\var + 1/2}}
				) \cdot R \\ \nonumber
				&+\sup_{n \in \N} (|\lambda_n|^{- 2\var}) 
					K(R) (
						(
							\sqrt{\eps} \|x\|_{H_{\var +1/2}}
							+  M_{\sqrt{\eps}} \|x\|_{H}
						)^2 +1
					) 
						\cdot R \\ \nonumber
			\leq{}& 
				\|x\|_{H_{\var + 1/2}}
					R \cdot
					\sup_{n \in \N} (|\lambda_n|^{1/2 - \var})
					(\eps \, K(R) +1)
				+\sup_{n \in \N} (|\lambda_n|^{1/2 - \var}) \,
					K(R) \cdot R \, (M_\eps \, R +1) \\ \nonumber
				&+\sup_{n \in \N} (|\lambda_n|^{- 2\var}) \,
					K(R) \cdot R \, (
						2 \eps \|x\|^2_{H_{\var +1/2}}
						+  2M^2_{\sqrt{\eps}} \, R^2 +1
					). 
	\end{align}
	Finally, we obtain 
	from \eqref{eq: H bound for F, B and v},
	the monotonicity of $K$ ,
	\eqref{eq: H3/2- norm bounds},
	and the fact that for all 
	$a,b \in \R$
	it holds that $(a+b)^2 \leq 2a^2 + 2b^2$
	that for all 
	$(t, x) \in W$ and all $R, \eps, \delta \in (0, \infty)$, $r \in \R$
	with $\|x\|_H \vee |r| \leq R$
	it holds that
	\begin{equation}
	\label{eq: v equation}
		\begin{split}
				&v(t,x) (r + \delta \|x\|^2_{H_\var})
			\leq
				K(R) \, (
					R +\delta (\sqrt{\eps} \|x\|_{H_{\var+ 1/2}} + M_{\sqrt{\eps}} \|x\|_H)^2
				) \\
			\leq{}
				&K(R) \, (
					R +2 \eps \cdot \delta \, \|x\|^2_{H_{\var+ 1/2}} 
					+2 \delta \, M^2_{\sqrt{\eps}} R^2
				).
		\end{split}
	\end{equation}
	Combining now \eqref{eq: delta F equation},
	\eqref{eq: delta B equation}, \eqref{eq: F + B equation},
	\eqref{eq: v equation} yields
	that for all 
	$(t, x) \in W$, $R, \eps, \delta \in (0, \infty)$, $r, K_1, K_2, K_3 \in \R$,
	$p \in H'$, and all $C \in \mathbb{S}_{\mathbb{H}, \mathbb{H}'}$
	with
	\begin{align}
	\nonumber
				K_1 
			={}&
				8 \delta \, K(R)
				+2 \sup_{n \in \N} (|\lambda_n|^{-2 \var}) \, K(R) \cdot R,\\ 
				K_2
			={}&
				2 \delta \, K(R) \, (M_\eps \cdot R +1)
				+ \sup(|\lambda_n|^{1/2 - \var}) \cdot R \, (\eps \, K(R) + 1), \\ \nonumber
				K_3 
			={}&
				2 \delta \, K(R) \, (2 M^2_{\sqrt{\eps}} \cdot R^2 +1)
				+ \sup_{n \in \N} (|\lambda_n|^{1/2 - \var}) K(R) \cdot (M_\eps R + 1) \\ \nonumber
				&+ \sup_{n \in \N} (|\lambda_n|^{- 2\var}) \, K(R) \cdot R \,
					(2 M^2_{\sqrt{\eps}}R^2 +1) 
				+K(R) (R + 2 \delta M^2_{\sqrt{\eps}} R^2),
	\end{align}
	with
	$\delta - \eps K_1 >0,$
	and with
	$
		\|x\|_H \vee |r| \vee \|p\|_{H'} \vee \|C\|_{L(\mathbb{H}, \mathbb{H}')} \leq R
	$
	it holds that
	\begin{align}
				&G_{\mathbb{H}, \mathbb{H}_\var, \delta, \tilde{V}h}^{+}((t,x),r,p,C)
			=
				G( 
					(t, x), 
					r + \delta \|x\|^2_{H_\var} , 
					p + \delta I_\mathbb{H} ((-A)^{2\var} x),
					(C|_{H_\var})|_{H_\var} + \delta I_{\mathbb{H}_\var}
				) \\ \nonumber
			={}&
				v(t,x) (r + \delta \|x\|^2_{H_\var})
				+
				\big\langle
					F(t,x) + A(x), p + \delta I_\mathbb{H} ((-A)^{2\var} x)
				\big\rangle_{H, H' }
				\\ \nonumber &+
				\langle
					B(t,x) ,
					I_{\mathbb{H}_\var}^{-1} \big(
						(C \, B(t,x))|_{H_\var}  + \delta I_{\mathbb{H}_\var} \, B(t,x)  
					\big)
				\rangle_{HS(\mathbb{U}, \mathbb{H}_\var)} \\ \nonumber
			\leq
			&- \|x\|^2_{H_{\var + 1/2}} 
			(
				2\delta
				- \eps 
				(
					2 \delta \, K(R)
					+ 4 \delta \, K(R)
					+2 \sup_{n \in \N} (|\lambda_n|^{-2 \var}) \, K(R) \cdot R
					+ 2 \delta \, K(R)
				)
			)\\ \nonumber
			&+\|x\|_{H_{\var + 1/2}} 
			(
				2 \delta \, K(R) \, (M_\eps \cdot R +1)
				+ \sup(|\lambda_n|^{1/2 - \var}) (\eps \, K(R) + R)
			) \\ \nonumber
			&+2 \delta \, K(R) \, (2 M^2_{\sqrt{\eps}} \cdot R^2 +1)
			+ \sup_{n \in \N} (|\lambda_n|^{1/2 - \var}) K(R) \cdot (M_\eps R + 1) \\ \nonumber
			&+ \sup_{n \in \N} (|\lambda_n|^{- 2\var}) \, K(R) \cdot R \,
				(2 M^2_{\sqrt{\eps}}R^2 +1) 
			+K(R) (R + 2 \delta M^2_{\sqrt{\eps}} R^2) \\ \nonumber
		\leq{}&
			- \|x\|^2_{H_{\var + 1/2}} 
			(
				2\delta
				- \eps K_1
			)
			+\|x\|_{H_{\var + 1/2}} \cdot K_2 
			+K_3 \\ \nonumber
		\leq{}&
			-\delta \|x\|^2_{H_{\var + 1/2}} 
			- (\delta- \eps K_1)
			(
				\|x\|_{H_{\var + 1/2}} + \tfrac{K_2}{2\delta- 2\eps K_1}
			)^2
			+\tfrac{K^2_2}{4\delta- 4\eps K_1}
			+K_3 \\ \nonumber
		\leq{}&
			-\delta \|x\|^2_{H_{\var + 1/2}}
			+\tfrac{K^2_2}{4\delta- 4\eps K_1}
			+K_3. 
	\end{align}
	Therefore we derive that for all 
	$R, \delta \in (0, \infty)$
	there exist an $L_{R, \delta} \in (0,\infty)$ such that for all
	$(t, x) \in W$, $r \in \R$,
	$p \in H'$, and all $C \in \mathbb{S}_{\mathbb{H}, \mathbb{H}'}$
	with
	$
		\|x\|_H \vee |r| \vee \|p\|_{H'} \vee \|C\|_{L(\mathbb{H}, \mathbb{H}')} \leq R
	$
	and with
	$G_{\mathbb{H}, \mathbb{H}_\var, \delta, \tilde{V}h}^{+}((t,x),r,p,C) \geq -R$
	it holds that
	\begin{equation}
	\label{eq: x H_1 bound}
		\|x\|_{\var + 1/2} \leq L_{R, \delta}.
	\end{equation}
	Moreover, we get from
	\eqref{eq: different HS norm},
	\eqref{eq: H_3/2 bound for F},
	and from \eqref{eq: H_3/2 bound for B}
	that for all 
	$(t, x) \in W$, $r \in \R$, $\delta \in (0, \infty)$, 
	$p$, $q \in H'$, $C$, $\mathfrak{C} \in \mathbb{S}_{\mathbb{H},\mathbb{H}'}$ 
	with 
	$\| C - \mathfrak{C} \|_{L(\mathbb{H}, \mathbb{H}')} \vee \|p - q\|_{H'} \leq \eps$ 
	it holds that
	\begin{equation}
	\label{eq: diff in p & A}
		\begin{split}
				&|G_{\mathbb{H}, \mathbb{H}_\var, \delta, \tilde{V}h}^+
					((t,x),r,p,C)-
				G_{\mathbb{H}, \mathbb{H}_\var, \delta, \tilde{V}h}^+
					((t,x),r,q,\mathfrak{C})| \\
			\leq{}&
				| \langle F (t,x) + A(x), p-q \rangle_{H,H'} | 
				+| \langle 
					B(t,x),
					I_{\mathbb{H}_\var}^{-1} \big( ((C-\mathfrak{C}) \, B(t,x))|_{H_\var} \big)
				\rangle_{HS(\mathbb{U},\mathbb{H}_\var)}  | \\
			\leq{}& 
				\|p - q\|_{H'} \| F (t,x) \|_{H} 
				+ \|p - q\|_{H'} \| A (x) \|_{H} 
				+| \langle 
					B(t,x),
					I_{\mathbb{H}}^{-1}(C-\mathfrak{C}) \, B(t,x)
				\rangle_{HS(\mathbb{U},\mathbb{H})}  | \\
			\leq{}&
				\eps \| F (t,x) \|_{H} 
				+ \eps \| x \|_{H_1} 
				+ \| C-\mathfrak{C} \|_{L(\mathbb{H}, \mathbb{H}')} 
					\| B(t,x) \|^2_{HS(\mathbb{U},\mathbb{H})} \\
			\leq{}&
				\eps \big(
					\sup_{n \in \N} (|\lambda_n|^{1/2 - \var}) \| F (t,x) \|_{H_{\var- 1/2}} 
					+\sup_{n \in \N} (|\lambda_n|^{1/2 - \var}) \|x\|_{H_{\var + 1/2}} \\
					&+\sup_{n \in \N} (|\lambda_n|^{- 2\var}) 
						\| B(t,x) \|^2_{HS(\mathbb{U},\mathbb{H}_\var )}
				\big) \\
			\leq{}&
				\eps \big(
					K(\|x\|_H) \,
						\sup_{n \in \N} (|\lambda_n|^{1/2 - \var}) 
						(\| x \|_{H_{\var+ 1/2 - \alpha_1}} +1) 
					+\sup_{n \in \N} (|\lambda_n|^{1/2 - \var}) \|x\|_{H_{\var + 1/2}} \\
					&+K(\|x\|_H) \,
						\sup_{n \in \N} (|\lambda_n|^{- 2\var}) 
						(\|x \|^2_{H_{\var+1/2- \beta_1}} +1)
				\big). 
		\end{split}
	\end{equation}
	Thus
	\eqref{eq: x H_1 bound},
	\eqref{eq: diff in p & A},
	and the monotonicity of $K$
	show that
	for all $\delta$, $R \in (0, \infty)$ it holds that
	\begin{equation}
	\label{eq: uniformly cont in p and A +}
		\begin{split}
				\lim_{\eps \downarrow 0} \Big [
					\sup \{& |G_{\mathbb{H}, \mathbb{H}_\var, \delta, \tilde{V}h}^{+}
						((t,x),r,p,C)
						-G_{\mathbb{H}, \mathbb{H}_\var, \delta, \tilde{V}h}^{+}
							((t,x),r,q,\mathfrak{C})| 
						\colon
						~ (t, x) \in W,
						~ r \in \R, \\
						~ &p, q \in H',
						~ C, \mathfrak{C} \in \mathbb{S}_{\mathbb{H}, \mathbb{H}'},
						~\max \{ 
								h(t, x), \|x\|_H, |r|, \|p\|_{H'}, \| C \|_{L(\mathbb{H}, \mathbb{H}')} 
							\} \leq R, \\
						~&\|p-q\|_{H'} \leq \eps,
						~\|C -\mathfrak{C} \|_{L(\mathbb{H},\mathbb{H}')} \leq \eps, 
						~G_{\mathbb{H}, \mathbb{H}_\var, \delta, \tilde{V}h}^{+}
							((t,x),r,p,C) \geq -R
					\}
				\Big ] \\
			\leq
				\lim_{\eps \downarrow 0} \Big [
					\sup \big \{& \eps \big(
						K(\|x\|_H) \,
							\sup_{n \in \N} (|\lambda_n|^{1/2 - \var}) 
							(\| x \|_{H_{\var+ 1/2 - \alpha_1}} +1) 
						+\sup_{n \in \N} (|\lambda_n|^{1/2 - \var}) \|x\|_{H_{\var + 1/2}} \\
						&+K(\|x\|_H) \,
							\sup_{n \in \N} (|\lambda_n|^{- 2\var}) 
							(\|x \|^2_{H_{\var+1/2- \beta_1}} +1)
					\big) \colon
						~ (t, x) \in W, \\
						~&\|x\|_{H_{\var + 1/2}} \leq L_{R, \delta},
						~\|x\|_H \leq R \big\}
				\Big ] \\
			\leq
				\lim_{\eps \downarrow 0} \Big [
					\eps \sup \big \{&
						K(R) \,
							\sup_{n \in \N} (|\lambda_n|^{1/2 - \var})
							\big(
								\sup_{n \in \N} (|\lambda_n|^{-\alpha_1})  \cdot\| x \|_{H_{\var+ 1/2}} 
								+1
							\big) 
						+\sup_{n \in \N} (|\lambda_n|^{1/2 - \var}) \|x\|_{H_{\var + 1/2}} \\
						&+K(R) \,
							\sup_{n \in \N} (|\lambda_n|^{- 2\var}) 
							\big(
								\sup_{n \in \N} (|\lambda_n|^{-\beta_1}) \|x \|^2_{H_{\var+1/2}} 
								+1
							\big)
					\colon 
						(t, x) \in W,\\
						~&\|x\|_{H_{\var + 1/2}} \leq L_{R, \delta},
				\Big ] 
			=
				0,
		\end{split}
	\end{equation}
	which implies \eqref{eq: uniformly bounded at p and A+}. 
	Analogously it follows \eqref{eq: uniformly bounded at p and A-}.
	Next fix $R \in (0,\infty)$, $t \in (0,T)$, $x \in O$, $r \in \R$, $p \in H'$,  
	and $\delta \in (0,1]$ and
	denote by $M \in (0, \infty)$ the value satisfying that
	$
			M 
		= 
			(\|x\|_H \vee |r| \vee \|p\|_{H'} \vee \|C\|_{L(\mathbb{H}, \mathbb{H}')} \vee R) +1
	$.
	In addition, let $(\tau_{N}, \xi_{N})_{N \in \N} \subseteq W$,
	be a sequence satisfying for all $N \in \N$ that
	$
			\|\xi_{N}\|_{H_{\var + 1/2}} 
		\leq 
			L_{M, \delta},
	$
	that
	$
		\|x - \xi_{N} \|_{H} \leq \dist_{\mathbb{H}}(x, H \backslash O) /2,
	$
	and that
	\begin{align}
	\nonumber
				&\lim_{N \to \infty} 
					\langle 
						B(\tau_N, \xi_N), \pi^{H}_{V_N^\perp} B(\tau_N, \xi_N) 
					\rangle_{HS(\mathbb{U},\mathbb{H})}\\
	\label{eq: def of sequence tau, xi}
			={}&
				\limsup_{N \to \infty} \big [
						\sup
						\{ 
							\langle 
								B(\tau, \xi) , \pi^{H}_{V_{N}^\perp} B(\tau, \xi)
							\rangle_{HS(\mathbb{U}, \mathbb{H})} \colon 
						(\tau, \xi) \in W,
						~ \|\xi\|_{H_{\var + 1/2}} \leq L_{M, \delta},
			\\ \nonumber & \qquad \qquad
						~ \|x - \xi \|_{H} \leq \dist_{\mathbb{H}}(x, H \backslash O) /2
					\big ].
	\end{align}
	Then it follows from the fact that
	$\mathbb{H}_{\var + 1/2}$ is compactly embedded in $\mathbb{H}_{1/2 - \beta_2}$ that the set
	$
		\{ y \in O \colon \|y\|_{H_{\var + 1/2}} \leq L_{M, \delta}, 
			~\|x-y \|_H \leq \dist_{\mathbb{H}}(x, H \backslash O) /2 \}
	$ 
	is compact with respect to the $\| \cdot \|_{H_{1/2- \beta_2}}$-norm and 
	therefore there exist a subsequence 
	$	
			(\tau_{N_i}, \xi_{N_i})_{i \in \N}
		\subseteq
			(\tau_{N}, \xi_{N})_{N \in \N}
	$
	and a $(\hat{\tau}, \hat{\xi}) \in [0,T] \times O$
	satisfying that
	$
		\lim_{i \to \infty} 
			\|(\tau_{N_i}- \hat{\tau}, \xi_{N_i} - \hat{\xi}) \|_{\R \times H_{1/2-\beta_2}} 
		= 0
	$.
	Thus,
	we have with 
	\eqref{eq: def of h},
	\eqref{eq: x H_1 bound},
	\eqref{eq: different HS norm}, 
	\eqref{eq: def of sequence tau, xi},
	and with
	\eqref{eq: Lipschitz continuity of B} 
	that 
	\begin{align}
				&\lim_{N \to \infty} \big [
					\sup
					\{ G_{\mathbb{H}, \mathbb{H}_\var, \delta, \tilde{V}h}^{+}
							((\tau, \xi), \nu, \rho, C + \alpha I_\mathbb{H} \pi^{H}_{V^\perp_N})
						-G_{\mathbb{H}, \mathbb{H}_\var, \delta, \tilde{V}h}^{+}((\tau, \xi), \nu, \rho, C) 
					\colon 
					\alpha \in (0, R), 
			~\\ \nonumber &\qquad\qquad 
					~ \nu \in \R, 
					~ (\tau, \xi) \in W,
					~ \rho \in H', 
					~ C \in \mathbb{S}_{\mathbb{H},\mathbb{H}'},
					~ h(\tau, \xi) \leq R, 
					~ \| C \|_{L(\mathbb{H},\mathbb{H}')} \leq R,  
			~\\ \nonumber &\qquad\qquad
					~ | t-\tau | \vee \|x- \xi\|_H \vee |r -\nu| \vee \|p - \rho\|_{H'} 
							\leq (\dist_{\mathbb{H}}(x, H \backslash O) /2) \vee 1, 
			~ \\ \nonumber &\qquad\qquad
					~ G_{\mathbb{H}, \mathbb{H}_\var, \delta, \tilde{V}h}^{+}((\tau, \xi), \nu, \rho, C) 
							\geq -R \}
				\big ] \\ \nonumber
			\leq{}&
				\limsup_{N \to \infty} \big [
					\sup
					\{ 
						\langle 
							B(\tau, \xi) , 
							\alpha I_{\mathbb{H}_\var}^{-1} \big( 
								(I_\mathbb{H} \pi^{H}_{V_{N}^\perp} B(\tau, \xi)) |_{H_\var}
							\big)
						\rangle_{HS(\mathbb{U}, \mathbb{H}_\var)} \colon 
					~ \alpha \in (0, R), 
					~ (\tau, \xi) \in W,
			~\\ \nonumber &\qquad\qquad 
					~\|x - \xi \|_H \leq (\dist_{\mathbb{H}}(x, H \backslash O) /2), 
					~ \|\xi\|_{\var + 1/2} \leq L_{M,\delta}
				\big ] \\ \nonumber
			\leq{}&
				R \limsup_{N \to \infty} \big [
					\sup
					\{ 
						\langle 
							B(\tau, \xi) , \pi^{H}_{V_{N}^\perp} B(\tau, \xi)
						\rangle_{HS(\mathbb{U}, \mathbb{H})} \colon  
					~ (\tau, \xi) \in W,
					~\|x - \xi \|_H \leq \dist_{\mathbb{H}}(x, H \backslash O) /2,
			\\ \nonumber &\qquad\qquad 
					~ \| \xi \|_{H_{\var + 1/2}} \leq  L_{M,\delta}
				\big ] \\ \nonumber
			={}&
				R \lim_{N \to \infty} 
					\langle 
						B(\tau_N, \xi_N), \pi^{H}_{V_N^\perp} B(\tau_N, \xi_N) 
					\rangle_{HS(\mathbb{U},\mathbb{H})} 
			=
				R \lim_{i \to \infty} 
					\langle 
						B(\tau_{N_i}, \xi_{N_i}), \pi^{H}_{V_{N_i}^\perp} B(\tau_{N_i}, \xi_{N_i}) 
					\rangle_{HS(\mathbb{U},\mathbb{H})} \\ \nonumber
			={}&
				R \lim_{i \to \infty} \Big(
					\langle 
						B(\tau_{N_i}, \xi_{N_i}), 
						\pi^{H}_{V_{N_i}^\perp} (B(\tau_{N_i}, \xi_{N_i}) - B(\hat{\tau}, \hat{x}))
					\rangle_{HS(\mathbb{U},\mathbb{H})}
			\\ \nonumber&\qquad\qquad 
					+ \langle 
						B(\tau_{N_i}, \xi_{N_i}), 
						\pi^{H}_{V_{N_i}^\perp} B(\hat{\tau}, \hat{x})
					\rangle_{HS(\mathbb{U},\mathbb{H})} 
				\Big)\\ \nonumber
			\leq{}&
				R \lim_{i \to \infty} \Big(
					\| B(\tau_{N_i}, \xi_{N_i}) - B(\hat{\tau}, \hat{x}) \|_{HS(\mathbb{U},\mathbb{H})} 
						\| 
							\pi^{H}_{V_{N_i}^\perp} B(\tau_{N_i}, \xi_{N_i}) 
						\|_{HS(\mathbb{U},\mathbb{H})}
			\\ \nonumber&\qquad\qquad 
					+ \| B(\tau_{N_i}, \xi_{N_i}) \|_{HS(\mathbb{U},\mathbb{H})} 
						\| 
							\pi^{H}_{V_{N_i}^\perp} B(\hat{\tau}, \hat{x}) 
						\|_{HS(\mathbb{U},\mathbb{H})} 
				\Big)\\ \nonumber
			={}&
				\lim_{i \to \infty} \left [
					R \, \left ( 
						\sum_{j=N_i+1}^\infty \| [B(\hat{\tau}, \hat{x}) ]^{*} e_j \|^2_{U} 
					\right )^{\frac 12} \,
					\| B(\tau_{N_i}, \xi_{N_i}) \|_{HS(\mathbb{U}, \mathbb{H})}
				\right ]  = 0,
	\end{align}
	which implies \eqref{eq: uniformly bounded in dimension+ LEM}.
	Analogously it follows \eqref{eq: uniformly bounded in dimension- LEM}.
  To this end we now verify 
  \eqref{eq:comparison.viscosity.solution.assumption}.
  For this let
	$\delta_n \in (0, 1]$,
  $
    ((t_n, x_n), r_n, p_n, C_n, \mathfrak{C}_n),
	$
	$
    ((\hat{t}_n, \hat{x}_n), \hat{r}_n, \hat{p}_n, \hat{C}_n, \hat{\mathfrak{C}}_n )
    \in
    W \times \R \times H' \times
    \mathbb{S}_{\mathbb{H}, \mathbb{H}'} \times \mathbb{S}_{\mathbb{H}, \mathbb{H}'}
  $,
	$
		(
		\tilde{p}_n, \tilde{C}_n) 
		\in 
		H' \times \mathbb{S}_{\mathbb{H}, \mathbb{H}'}
	$,
  $  
    n \in \N 
  $,
  satisfy that 
  $
    w-\lim_{ n \to \infty }( t_n, x_n ) \in [0,T) \times O
  $, that
	\begin{align}
    &\lim_{ n \to \infty }
    \big( 
      \sqrt{ n }
      \, 
      \|
        x_n  
        -
        \hat{x}_n
      \|_{H}
    \big)
  +
    \lim_{ n \to \infty }
    \big( 
      \sqrt{ n }
      \, 
      | t_n - \hat{t}_n|
    \big)
    = 0, \\
	\label{eq: def of delta}
			&\delta_n
		=
			\inf \Big( 
				\Big\{ 
					\delta \in(\nicefrac 1n, \infty) \colon 
							(K(\tfrac 1 \delta))^2+\tfrac {1}{\delta^2}+1
						\leq
							\big(
								(\sup_{m \in [n,\infty)}\omega(m^{- \nicefrac 14}))+n^{- \nicefrac 12}
							\big)^{-\nicefrac 12}
				\Big\}
				\cup \{1\}
			\Big), \\
    &\lim_{ n \to \infty }
    \big( 
      \delta_n
      \, 
      (h(t_n,x_n) \vee h(\hat{t}_n, \hat{x}_n)
    \big)
    = 0 
    <
    \liminf_{ n \to \infty } (r_n - \hat{r}_n )
    \leq
    \sup_{ n \in \N } 
    ( | r_n | + | \hat{r}_n | )
    <
    \infty, \\
	  &p_n
    =
    n I_\mathbb{H} \left( x_n - \hat{x}_n \right) 
    \tilde{V}( t_n, x_n )
    + 
    r_n \, ( D_\mathbb{H} \tilde{V} )(t_n, x_n ), \\
    &\hat{p}_n
    =
    n I_\mathbb{H} \left( x_n - \hat{x}_n \right) 
    \tilde{V}( \hat{t}_n, \hat{x}_n )
    + 
    \hat{r}_n \, 
    ( D_\mathbb{H} \tilde{V} )( \hat{t}_n, \hat{x}_n ), \\
		&\tilde{p}_n
    =
    (r_n - \hat{r}_n) \, 
    ( D_\mathbb{H} \tilde{V} )( t_n, x_n ), \\
		\begin{split}
     &C_n
		=
     n 
     (I_\mathbb{H} \left( 
       x_n - \hat{x}_n
     \right)) 
		\otimes ( D_\mathbb{H} \tilde{V} )( t_n , x_n )
     +
     (
       D_\mathbb{H} 
       \tilde{V}
     )( t_n, x_n )
		\otimes
     \left( 
			 n I_\mathbb{H} \left( 
       x_n - \hat{x}_n 
     \right) \right)
     \\& \qquad+
     r_n \, 
     ( D^2_\mathbb{H} \tilde{V})( t_n, x_n ), 
	\end{split}\\
	\begin{split}
   &\hat{C}_n
		=
     n 
     (I_\mathbb{H} \left( 
       x_n - \hat{x}_n
     \right))
     \otimes ( D_\mathbb{H} \tilde{V} )( \hat{t}_n , \hat{x}_n )
     +
     (
       D_\mathbb{H} 
       \tilde{V}
     )( \hat{t}_n , \hat{x}_n )
     \otimes 
     \left(n I_\mathbb{H} \left( 
       x_n - \hat{x}_n 
     \right) \right) \\
     & \qquad+
     \hat{r}_n \, 
     ( D^2_\mathbb{H} \tilde{V})( \hat{t}_n , \hat{x}_n ), 
	\end{split}\\
    &\tilde{C}_n
		=
     (r_n-\hat{r}_n) \, 
     ( D^2_\mathbb{H} \tilde{V})( t_n , x_n ), 
	\end{align}
  and that
	for all
	$
		R \in (0, \infty)
	$,
	and for all
	$
		(z^{(n)})_{n \in \N}, (\hat{z}^{(n)})_{n \in \N} 
			\subseteq \{z \in H \colon ~ \| z \|_H \leq R\}
	$
	it holds that
	\begin{equation}
		\limsup_{n \to \infty} \left(
			\langle z^{ (n) }, \mathfrak{C}_n z^{ (n) } \rangle_{H, H'}
			- \langle \hat{z}^{ (n) }, \hat{\mathfrak{C}}_n \hat{z}^{ (n) } \rangle_{H, H'}
		\right)
  \leq 
		\limsup_{n \to \infty}
    3
    \| z^{ (n) } - \hat{z}^{ (n) } \|_H^2.
	\end{equation}
	First note, that the fact that $(x_n)_{n \in \N}$ and $(\hat{x}_n)_{n \in \N}$ 
	are weakly convergent, the fact that
	$
			\lim_{ n \to \infty }
			\big( 
				\sqrt{ n }
				\, 
				| t_n - \hat{t}_n|
			\big)
    = 
			\lim_{ n \to \infty }
			\big( 
				\sqrt{ n }
				\, 
				\| x_n - \hat{x}_n\|_{H}
			\big)
    = 0,
	$
	and the fact that
	$
		\sup_{ n \in \N } 
			( | r_n | + | \hat{r}_n | )
    <
    \infty
	$
	imply that there exists a $R \in (0, \infty)$ such that for all $n \in \N$
	it holds that
	\begin{equation}
	\label{eq: R bound}
			\|x_n\|_{H} \vee \|\hat{x}_n\|_{H} \vee| r_n | \vee | \hat{r}_n | 
			\vee (\sqrt{n} |t_n - \hat{t}_n|) \vee (\sqrt{n} \|x_n - \hat{x}_n\|_H)
		\leq R.
	\end{equation}
	Furthermore, it follows from the monotonicity of $K$ 
	and from \eqref{eq: def of delta} that
	for all $n \in \N$ 
	with 
	$\delta_n (\|x_n\|_{H_\var}^2 \vee \|\hat{x}_n\|_{H_\var}^2) \leq 1$
	and with $\delta_n < 1$
	it holds that
	\begin{align}
	\label{eq: x vartheta bound 1}
		\begin{split}
				&\Big(
					\big(K(\|x_n\|^2_{H_\var} \vee \|\hat{x}_n\|^2_{H_\var}) \big)^2
					\vee K(\|x_n\|^2_{H_\var} \vee \|\hat{x}_n\|^2_{H_\var})
			\\ &
					\vee \big(
						K(\|x_n\|^2_{H_\var} \vee \|\hat{x}_n\|^2_{H_\var}) 
						(\|x_n\|^2_{H_\var} \vee \|\hat{x}_n\|^2_{H_\var})
					\big)
				\Big) 
				\cdot (
					\|x_n-\hat{x}_n\|^2_{H_{1/2}}
					+\omega(\|x_n-\hat{x}_n\|^2_{H_{1/2}})
				) \\ 
			\leq{}&
				((K(\tfrac{1}{\delta_n}))^2+\tfrac{1}{\delta^2_n}+1)
				\cdot	(
					(\sup_{m \in [n,\infty)}\omega(m^{- \nicefrac 14} ))
					+n^{- \nicefrac 12}
				)^{\nicefrac 12}
				\cdot	(
					(\sup_{m \in [n,\infty)}\omega( m^{- \nicefrac 14} ))
					+n^{- \nicefrac 12}
				)^{-\nicefrac 12} \\
				&\cdot \Big(
					\|x_n-\hat{x}_n\|^2_{H_{1/2}}
					+\big( \sup_{m \in [n,\infty)} \omega( m^{- \nicefrac 14} ) \big)
					+\big( \sup_{m \in [0,\infty)} \omega(m ) \big)
						\cdot n^{\nicefrac 14} \, \|x_n-\hat{x}_n\|^2_{H_{1/2}}
				\Big)\\
			\leq{}&
					n^{\nicefrac 14} \|x_n-\hat{x}_n\|^2_{H_{1/2}}
					+((\sup_{m \in [n,\infty)}\omega(m^{- \nicefrac 14} ))^{\nicefrac 12}
					+\big( \sup_{m \in [0,\infty)} \omega(m ) \big)
						\cdot n^{\nicefrac 12} \, \|x_n-\hat{x}_n\|^2_{H_{1/2}},
		\end{split}
	\end{align}
	that
	\begin{align}
	\nonumber
				&K(\|x_n\|^2_{H_\var} \vee \|\hat{x}_n\|^2_{H_\var}) \\
	\label{eq: x vartheta bound 2}
			\leq{}&
				((K(\tfrac{1}{\delta_n}))^2+\tfrac{1}{\delta^2_n}+1)
				\cdot	(
					(\sup_{m \in [n,\infty)}\omega(m^{- \nicefrac 14} ))
					+n^{- \nicefrac 12}
				)^{\nicefrac 12}
				\cdot	(
					(\sup_{m \in [n,\infty)}\omega( m^{- \nicefrac 14} ))
					+n^{- \nicefrac 12}
				)^{-\nicefrac 12} \\ \nonumber
			\leq{}&
				n^{\nicefrac 14},
	\end{align}
	and that
	\begin{align}
	\nonumber
				&K(\|x_n\|^2_{H_\var} \vee \|\hat{x}_n\|^2_{H_\var}) 
				\cdot (
					\|x_n-\hat{x}_n\|^2_{H_{1/2}}
					+\omega(\|x_n-\hat{x}_n\|^2_{H_{1/2}} +|t_n - \hat{t}_n|^2)
				) \\ 
	\label{eq: x vartheta bound 3}
			\leq{}
				&(
					n^{\nicefrac 14} \|x_n-\hat{x}_n\|^2_{H_{1/2}}
					+((\sup_{m \in [n,\infty)}\omega(m^{- \nicefrac 14} ))^{\nicefrac 12}
					+\big( \sup_{m \in [0,\infty)} \omega(m ) \big)
			\\ \nonumber &
						\cdot n^{\nicefrac 12} \, (\|x_n-\hat{x}_n\|^2_{H_{1/2}} +|t_n - \hat{t}_n|^2)
				).
	\end{align}
  In addition we get with 
	\eqref{eq: H bound for F, B and v},
	the monotonicity of $K$,
	and
	the fact that 
	$   
		\lim_{ n \to \infty }
    \big( 
      \delta_n
      \, 
      (h(t_n,x_x) \vee h(\hat{t}_n, \hat{x}_n)
    \big)
    = 0
	$ that
	\begin{equation}  
	\label{eq: delta v limit}
				\limsup_{ n \to \infty } \big(\delta_n (
					v(t_n, x_n)  h(t_n,x_n)
					+ v(\hat{t}_n, \hat{x}_n) h(\hat{t}_n,\hat{x}_n)
				)\big) 
			\leq{} 
				K(R) \limsup_{ n \to \infty } \big(
					\delta_n (
						h(t_n,x_n)
						+h(\hat{t}_n,\hat{x}_n)
					)
				\big)
			= 0.
	\end{equation}
	Moreover,
	combining \eqref{eq: delta F equation} and \eqref{eq: delta B equation}
	together with the fact that 
	$   
		\lim_{ n \to \infty }
      \delta_n
    = 0
	$,
	and the fact that $\tilde{V} > 1$,
	yields for all $\eps \in (0, \tfrac{1}{3 K(R)})$ that
	\begin{align}
	\label{eq: delta F and B limit}
				&\limsup_{n \to \infty} \Big(
					\tfrac
					{
						\langle 
							F(t_n, x_n) + A(x_n),
							\delta_n 
							E_{\mathbb{H}_\var', \mathbb{H}'} (
								(D_{\mathbb{H}_\var} (\tilde{V}h|_{\R \times H_\var}))(t_n,x_n)
							)
						\rangle_{H, H'}
					}
					{\tilde{V}(t_n, x_n)}
				\\ \nonumber  & \qquad \qquad
				+\tfrac
						{
							\big \langle 
								B(t_n, x_n), 
								\delta_n I_{\mathbb{H}_\var}^{-1} \big(
									(D^2_{\mathbb{H}_\var} (h\tilde{V})) (t_n, x_n) \, B(t_n, x_n)
								\big)
							\big \rangle_{HS(\mathbb{U},\mathbb{H}_\var)}
						}
						{ \tilde{V}(t_n,x_n) }
				\Big) \\ \nonumber
			\leq{}&
				\limsup_{n \to \infty} \Big(
					2\delta_n
					\, \tfrac	
						{
							-(1- 3 \eps \, K(R)) \|x_n\|^2_{H_{\var + 1/2}}
							+K(R) \, (M_\eps \cdot R + 1) \|x_n\|_{H_{\var + 1/2}}
						} 
						{ \tilde{V}(t_n,x_n) }
						+2 \delta_n \,
						\tfrac	
						{
							K(R) \, (2 M^2_{\sqrt{\eps}} \cdot R^2 + 1)
						} 
						{ \tilde{V}(t_n,x_n) } 
				\Big) \\ \nonumber
			\leq{}&
				\limsup_{n \to \infty} \Big(
					2 \delta_n \,
					\tfrac	
						{
							-(1 - 3 \eps \, K(R)) 
							\big( 
								\|x_n\|_{H_{\var + 1/2}} 
								- K(R) \cdot\frac{ M_{\eps} \cdot R +1 }{2(1 - 3 \eps \, K(R))} 
							\big)^2
						} 
						{ \tilde{V}(t_n, x_n) }
					+	2 \delta_n \, 
					\tfrac	
						{
							(K(R))^2 \cdot\frac{ (M_{\eps} \cdot R +1)^2 }{4(1 - 3 \eps \, K(R))} 
							+ K(R) \, (2 M^2_{\sqrt{\eps}} \cdot R^2 + 1)
						} 
						{ \tilde{V}(t_n, x_n) } 
				\bigg) \\ \nonumber
			\leq{}&
				\limsup_{n \to \infty} \Big(
					2 \delta_n \,
					\tfrac	
						{
							(K(R))^2 \cdot\frac{ (M_{\eps} \cdot R +1)^2 }{4(1 - 3 \eps \, K(R))} 
							+ K(R) \, (2 M^2_{\sqrt{\eps}} \cdot R^2 + 1)
						} 
						{ \tilde{V}(t_n, x_n) } 
				\Big)
			= 0.
	\end{align}
	Analogously it follow that
	\begin{equation}
	\label{eq: delta F and B limit hat}
		\begin{split}
				&\limsup_{n \to \infty} \Big(
					\tfrac
						{
							\langle 
								F(\hat{t}_n,\hat{x}_n) + A(\hat{x}_n),
								\delta_n 
								E_{\mathbb{H}_\var', \mathbb{H}'} (
									(D_{\mathbb{H}_\var} (\tilde{V}h|_{\R \times H_\var}))
										(\hat{t}_n,\hat{x}_n)
							)
							\rangle_{H, H'}
						}
						{ \tilde{V}(\hat{t}_n,\hat{x}_n) } \\
				& \qquad \qquad
					+ \tfrac
						{
							\big \langle 
								B( \hat{t}_n, \hat{x}_n ), 
								\delta_n I_{\mathbb{H}_\var}^{-1} \big(
									(D^2_{\mathbb{H}_\var} (h\tilde{V})) ( \hat{t}_n, \hat{x}_n ) 
										\, B( \hat{t}_n, \hat{x}_n )
								\big)
							\big \rangle_{HS(\mathbb{U},\mathbb{H}_\var)}
						}
						{ \tilde{V}(\hat{t}_n,\hat{x}_n) }
				\Big) 
			\leq 0.
		\end{split}
	\end{equation}
	In addition,
	we get from
	$\tilde{V}> 1$
	\eqref{eq: H bound for F, B and v},
	\eqref{eq: uniform continuity of v},
	\eqref{eq: R bound},
	and from
	the monotonicity of $K$
	that for all $n \in \N$ it holds that
	\begin{equation}
	\label{eq: v difference inequality}
		\begin{split}
				&\tfrac{v(t_n,x_n)}{\tilde{V}(t_n, x_n)} r_n
					-\tfrac{v(\hat{t}_n,\hat{x}_n)}{\tilde{V}(\hat{t}_n,\hat{x}_n)} \hat{r}_n 
					-\tfrac{v(t_n,x_n)}{\tilde{V}(t_n, x_n)} (r_n -\hat{r}_n) \\
			\leq{}
				&\Big |
					\tfrac{v(t_n,x_n) - v(\hat{t}_n,\hat{x}_n)}{\tilde{V}(t_n, x_n)} \hat{r}_n
				\Big|
				+\Big|
					v(\hat{t}_n,\hat{x}_n) \hat{r}_n
					\big (
						\tfrac{1} {\tilde{V}(t_n,x_n)}
						-\tfrac {1}{\tilde{V}(\hat{t}_n,\hat{x}_n)}
					\big )
				\Big | \\
			\leq{}&
				\left |
					v(t_n,x_n) - v(\hat{t}_n,\hat{x}_n)
				\right|
				\cdot |\hat{r}_n|
				+ |\hat{r}_n| \,
				|
					v(\hat{t}_n,\hat{x}_n)
				|
				\cdot
				\Big |
					\tfrac
						{
							\tilde{V}(\hat{t}_n,\hat{x}_n)
							-\tilde{V}(t_n,x_n)
						}
						{
						\tilde{V}(\hat{t}_n,\hat{x}_n)
						\cdot \tilde{V}(t_n,x_n)
						}
				\Big | \\
			\leq{}&
				K(\|x\|^2_{H_\var} \vee \|\hat{x}\|^2_{H_\var}) \cdot
					(\|x- \hat{x}\|^2_{H_{1/2}} + \omega(\|x- \hat{x}\|^2_{H_{1/2}} + |t -\hat{t}|^2)) 
					\, R 
		\\& \qquad
				+ R \cdot K(R) 
				\,
				\big |
					\tilde{V}(\hat{t}_n,\hat{x}_n)
					-\tilde{V}(t_n,x_n)
				\big |. 
		\end{split}
	\end{equation}
	Moreover,
	the definition of $(p_n)_{n \in \N}$ and $(\hat{p}_n)_{n \in \N}$,
	\eqref{eq: def of tilde V},
	\eqref{eq: V continuous with respect to H1},
	\eqref{eq: V bounded with respect to H 1/2},
	and
	\eqref{eq: R bound}
	implies that
	for all $n \in \N$ it holds that
	\begin{align}
	\label{eq: A difference inequality}
				&\tfrac{\langle A(x_n), p_n \rangle_{H, H'}}
					{\tilde{V}(t_n, x_n)}
				- \tfrac{\langle A(\hat{x}_n), \hat{p}_n \rangle_{H, H'}}
					{\tilde{V}(\hat{t}_n, \hat{x}_n)} 
				- \tfrac {\left \langle A(x_n), \tilde{p}_n \right \rangle_{H, H'}}
					{\tilde{V}(t_n, x_n)}
			\\ \nonumber ={}&
				\langle A(x_n), n I_\mathbb{H} (x_n - \hat{x}_n) \rangle_{H, H'} 
				-\langle A(\hat{x}_n), n I_\mathbb{H} (x_n - \hat{x}_n) \rangle_{H, H'} 
				+\Big \langle 
					A(x_n), \tfrac{r_n \, ( D_\mathbb{H} \tilde{V} )(t_n, x_n )} {\tilde{V}(t_n, x_n)}
				\Big \rangle_{H, H'} 
			\\ \nonumber &
				-\Big \langle 
					A(\hat{x}_n), 
					\tfrac{\hat{r}_n \, ( D_\mathbb{H} \tilde{V} )(\hat{t}_n, \hat{x}_n)} 
						{\tilde{V}(\hat{t}_n, \hat{x}_n)}
				\Big \rangle_{H, H'}	
				-\Big \langle 
					A(x_n), 
					\tfrac{(r_n - \hat{r}_n)\, ( D_\mathbb{H} \tilde{V} )(t_n, x_n )} {\tilde{V}(t_n, x_n)}
				\Big \rangle_{H, H'} 
			\\ \nonumber ={}&
				\langle A(x_n-\hat{x}_n), n (x_n - \hat{x}_n) \rangle_{H}
				+\Big \langle 
					A(x_n - \hat{x}_n), 
					\tfrac{\hat{r}_n \, ( D_\mathbb{H} V )(\hat{x}_n)} 
						{V(\hat{x}_n)}
				\Big \rangle_{H, H'} 
		\\ \nonumber &
				+\Big \langle 
					A(x_n), 
					\tfrac{\hat{r}_n \, ( D_\mathbb{H} V )(x_n )} {V(x_n)}
					-\tfrac{\hat{r}_n \, ( D_\mathbb{H} V )(\hat{x}_n)} 
						{V(\hat{x}_n)}
				\Big \rangle_{H, H'}
			\\ \nonumber \leq{}&
				-n \|x_n - \hat{x}_n \|^2_{H_{1/2}} 
				+\|x_n - \hat{x}_n \|_{H_{1/2}} \cdot
					\tfrac{|\hat{r}_n| \, \|I_\mathbb{H}^{-1}( D_\mathbb{H} V )(\hat{x}_n) \|_{H_{1/2}}} 
						{|V(\hat{x}_n)|}
			\\ \nonumber &
				+ |\hat{r}_n| \, \| x_n \|_{H_{\var}} 
					\Big \|
						\tfrac{( D_\mathbb{H} V )(x_n )} {V(x_n)} 
						-\tfrac{( D_\mathbb{H} V )(\hat{x}_n)} 
						{V(\hat{x}_n)}
					\Big \|_{H_{1-\var}}
			\\ \nonumber \leq{}&
				-n \|x_n - \hat{x}_n \|^2_{H_{1/2}}
				+R \cdot K(\| \hat{x}_n \|^2_{H_{\var}})
					\|x_n - \hat{x}_n \|_{H_{1/2}} 
			\\ \nonumber &
				+ R \cdot K(\| x_n \|^2_{H_{\var}} \vee \| \hat{x}_n \|^2_{H_{\var}}) 
					\cdot \| x_n \|_{H_{\var}} \,
					\big (
						\| x_n - \hat{x}_n \|^2_{H_{1/2}} 
						+ \omega(\| x_n - \hat{x}_n \|^2_{H_{1/2}})
					\big ).
	\end{align}
	Furthermore, we obtain from
	the definition of $(p_n)_{n \in \N}$ and $(\hat{p}_n)_{n \in \N}$,
	and from
	\eqref{eq: def of tilde V}
	that for all $n \in \N$ 
	it holds that
	\begin{align}
			&\frac{\langle F(t_n, x_n), p_n \rangle_{H, H'}}
						{\tilde{V}(t_n, x_n)}
					- \frac{\langle F(\hat{t}_n, \hat{x}_n), \hat{p}_n \rangle_{H, H'}}
						{\tilde{V}(\hat{t}_n, \hat{x}_n)} 
					- \frac {\left \langle F(t_n, x_n), \tilde{p}_n \right \rangle_{H, H'}}
						{\tilde{V}(t_n, x_n)}
				\\ \nonumber ={} &
					\langle F(t_n, x_n), n I_\mathbb{H} (x_n - \hat{x}_n) \rangle_{H, H'} 
					-\langle 
						F(\hat{t}_n, \hat{x}_n), n I_\mathbb{H} (x_n - \hat{x}_n)
					\rangle_{H, H'} 
				\\ \nonumber &
					+\Big \langle 
						F(t_n, x_n), 
						\frac{r_n \, ( D_\mathbb{H} \tilde{V} )(t_n, x_n )} {\tilde{V}(t_n, x_n)}
					\Big \rangle_{H, H'}		
					-\Big \langle 
						F(\hat{t}_n, \hat{x}_n), 
						\frac{\hat{r}_n \, ( D_\mathbb{H} \tilde{V} )(\hat{t}_n, \hat{x}_n)} 
							{\tilde{V}(\hat{t}_n, \hat{x}_n)}
					\Big \rangle_{H, H'}	
				\\ \nonumber &
					-\Big \langle 
						F(t_n, x_n), 
						\frac
							{(r_n - \hat{r}_n)\, ( D_\mathbb{H} \tilde{V} )(t_n, x_n )} 
							{\tilde{V}(t_n, x_n)}
					\Big \rangle_{H, H'}
			\\ \nonumber ={}&
				\langle F(t_n, x_n)-F(\hat{t}_n, \hat{x}_n), n (x_n - \hat{x}_n) \rangle_{H}
				+\Big \langle 
					F(t_n, x_n), 
					\frac{\hat{r}_n \, ( D_\mathbb{H} V )(x_n )} {V(x_n)}
					-\frac{\hat{r}_n \, ( D_\mathbb{H} V )(\hat{x}_n)} 
						{V(\hat{x}_n)}
				\Big \rangle_{H, H'}
			\\ \nonumber &
				+\Big \langle 
					F(t_n, x_n) - F(\hat{t}_n, \hat{x}_n), 
					\frac{\hat{r}_n \, ( D_\mathbb{H} V )(\hat{x}_n)} 
						{V(\hat{x}_n)}
				\Big \rangle_{H, H'}	
			\\ \nonumber \leq{}&
				n \|F(t_n, x_n) - F(\hat{t}_n, \hat{x}_n) \|_{H_{-1/2}} 
					\cdot \|x_n - \hat{x}_n \|_{H_{1/2}} 
			\\ \nonumber &
				+\|F(t_n, x_n) - F(\hat{t}_n, \hat{x}_n) \|_{H_{-1/2}} \cdot
					\frac
						{|\hat{r}_n| \, \|I_\mathbb{H}^{-1}( D_\mathbb{H} V )(\hat{x}_n) \|_{H_{1/2}}} 
						{|V(\hat{x}_n)|} 
			\\ \nonumber &
				+ |\hat{r}_n| \, \| F(t_n, x_n) \|_{H_{\var-1}} 
					\cdot \left \|
						\frac{( D_\mathbb{H} V )(x_n )} {V(x_n)} 
						-\frac{( D_\mathbb{H} V )(\hat{x}_n)} 
						{V(\hat{x}_n)}
					\right \|_{H_{1-\var}} 
	\end{align}
	and this together with
	\eqref{eq: H bound for F, B and v},
	\eqref{eq: Lipschitz continuity of F},
	\eqref{eq: V continuous with respect to H1},
	\eqref{eq: V bounded with respect to H 1/2}, 
	\eqref{eq: R bound},
	the monotonicity of $K$,
	$V > 1$,
	with \eqref{eq: H1/2- norm bounds}
	and with 
	\eqref{eq: x vartheta bound 2}
	shows for all
	$n \in \N$ and all $\eps \in (0, \infty)$ that
	\begin{align}
	\nonumber
				&\frac{\langle F(t_n, x_n), p_n \rangle_{H, H'}}
						{\tilde{V}(t_n, x_n)}
					- \frac{\langle F(\hat{t}_n, \hat{x}_n), \hat{p}_n \rangle_{H, H'}}
						{\tilde{V}(\hat{t}_n, \hat{x}_n)} 
					- \frac {\left \langle F(t_n, x_n), \tilde{p}_n \right \rangle_{H, H'}}
						{\tilde{V}(t_n, x_n)}\\ \nonumber
			\leq{}& 
				n \, K (R) \, (\|x_n - \hat{x}_n \|_{H_{1/2 - \alpha_2}} + |t_n- \hat{t}_n|)
					\cdot \|x_n - \hat{x}_n\|_{H_{1/2}} \\ \nonumber
				&+R \cdot K(R) \cdot K(\| \hat{x}_n \|^2_{H_{\var}}) \, 
					(\|x_n - \hat{x}_n \|_{H_{1/2 - \alpha_2}} + |t_n- \hat{t}_n|)\\ \nonumber
				&+ R \cdot (K(\| x_n \|^2_{H_{\var}} \vee \| \hat{x}_n \|^2_{H_{\var}}))^2
					\, (\|x_n-\hat{x}_n\|^2_{H_{1/2}}+\omega(\|x_n-\hat{x}_n\|^2_{H_{1/2}}))\\ 
	\label{eq: F difference inequality}	
			\leq{}&
				n \, K (R) \, 
				(
					\eps \|x_n - \hat{x}_n \|_{H_{1/2}} 
					+ M_\eps \|x_n - \hat{x}_n \|_{H} 
					+ |t_n- \hat{t}_n|
				)
					\cdot \|x_n - \hat{x}_n\|_{H_{1/2}} \\ \nonumber
				&+R \cdot K(R) \cdot K(\| \hat{x}_n \|^2_{H_{\var}})\, 
					(
						\eps \|x_n - \hat{x}_n \|_{H_{1/2}} 
						+ M_\eps \, \|x_n - \hat{x}_n\|_H + |t_n- \hat{t}_n|
					) \\ \nonumber
				&+ R \cdot (K(\| x_n \|^2_{H_{\var}} \vee \| \hat{x}_n \|^2_{H_{\var}}))^2
					\, (\|x_n-\hat{x}_n\|^2_{H_{1/2}}+\omega(\|x_n-\hat{x}_n\|^2_{H_{1/2}}))\\ \nonumber
			\leq{}&
				n \, K (R) \, 
				(
					\eps \|x_n - \hat{x}_n \|_{H_{1/2}} 
					+ M_\eps \|x_n - \hat{x}_n \|_{H} 
					+ |t_n- \hat{t}_n|
				)
					\cdot \|x_n - \hat{x}_n\|_{H_{1/2}} \\ \nonumber
				&+R \cdot K(R) \cdot n^{\nicefrac 14}\, 
					(
						\eps \|x_n - \hat{x}_n \|_{H_{1/2}} 
						+(M_\eps + 1) \, \tfrac{R}{\sqrt{n}}
					) \\ \nonumber
				&+ R \cdot (K(\| x_n \|^2_{H_{\var}} \vee \| \hat{x}_n \|^2_{H_{\var}}))^2
					\, (\|x_n-\hat{x}_n\|^2_{H_{1/2}}+\omega(\|x_n-\hat{x}_n\|^2_{H_{1/2}}))
	\end{align}
	Moreover, it follows from \eqref{eq: def of tilde V} and from $V>1$
	that for all $(t,x) \in (0,T) \times O$ it holds that
	$\tilde{V}(t,x) > 1$ and thus
	we get from \eqref{eq: different HS norm} that
	for all $n \in \N$ it holds that
	\begin{align}
	\label{eq: B difference}
				&\frac
					{
						\big \langle
							B( t_n, x_n ) ,
							I_{\mathbb{H}_\var}^{-1} \big(
								((C_n - \tilde{C}_n) \, B( t_n, x_n ))|_{H_\var}
							\big)
						\big \rangle_{HS(\mathbb{U}, \mathbb{H}_\var)} 					
					}
					{ \tilde{V}( t_n, x_n ) } \\ \nonumber
				&-\frac
					{
						\big \langle
							B( \hat{t}_n, \hat{x}_n ) ,
							I_{\mathbb{H}_\var}^{-1} \big(
								(\hat{C}_n \, B( \hat{t}_n, \hat{x}_n ))|_{H_\var}
							\big)
						\big \rangle_{HS(\mathbb{U}, \mathbb{H}_\var)}
					}
					{ \tilde{V}( \hat{t}_n, \hat{x}_n ) } \\ \nonumber
			={}&
				\frac
					{
						\langle
							B( t_n, x_n ) ,
							I_{\mathbb{H}}^{-1} 
								(C_n - \tilde{C}_n) \, B( t_n, x_n )
						\rangle_{HS(\mathbb{U}, \mathbb{H})} 					
					}
					{ \tilde{V}( t_n, x_n ) } 
				-\frac
					{
						\langle
							B( \hat{t}_n, \hat{x}_n ) ,
							I_{\mathbb{H}}^{-1} \hat{C}_n \, B( \hat{t}_n, \hat{x}_n )
						\rangle_{HS(\mathbb{U}, \mathbb{H})}
					}
					{ \tilde{V}( \hat{t}_n, \hat{x}_n ) }  \\ \nonumber
			={}&
				\frac
					{
						\langle
							B( t_n, x_n ) ,
							(-A)^{-2 \kappa} I_{\mathbb{H}}^{-1} 
								(C_n - \tilde{C}_n) \, B( t_n, x_n )
						\rangle_{HS(\mathbb{U}, \mathbb{H}_{\kappa})} 					
					}
					{ \tilde{V}( t_n, x_n ) } \\ \nonumber
				&-\frac
					{
						\langle
							B( \hat{t}_n, \hat{x}_n ) ,
							(-A)^{-2 \kappa} I_{\mathbb{H}}^{-1} \hat{C}_n \, B( \hat{t}_n, \hat{x}_n )
						\rangle_{HS(\mathbb{U}, \mathbb{H}_{\kappa})}
					}
					{ \tilde{V}( \hat{t}_n, \hat{x}_n ) }  \\ \nonumber
			={}&
				\bigg \langle
							B( t_n, x_n ),
							(-A)^{-2 \kappa} I_{\mathbb{H}}^{-1}
								\bigg (
									\frac{C_n}{\tilde{V}( t_n, x_n )} 
									-\frac{\hat{C}_n}{\tilde{V}( \hat{t}_n, \hat{x}_n )} 
									-\frac{\tilde{C}_n}{\tilde{V}( t_n, x_n )}
								\bigg )
								\, B(t_n, x_n )
				\bigg \rangle_{HS(\mathbb{U}, \mathbb{H}_{\kappa})} \\ \nonumber
				&+\frac
					{
						\big \langle
							B( t_n, x_n ) - B( \hat{t}_n, \hat{x}_n ),
							(-A)^{-2 \kappa} I_{\mathbb{H}}^{-1}
								\hat{C}_n \, \big( B( t_n, x_n ) + B( \hat{t}_n, \hat{x}_n ) \big)
						\big \rangle_{HS(\mathbb{U}, \mathbb{H}_{\kappa})} 					
					}
					{ \tilde{V}( \hat{t}_n, \hat{x}_n ) } \\ \nonumber
			\leq{}& 
				\|
						B( t_n, x_n )
					\|^2_{HS(\mathbb{U}, \mathbb{H}_{\kappa})}
					\bigg \|
						I_{\mathbb{H}}^{-1} 
							\bigg(
								\frac{C_n}{\tilde{V}( t_n, x_n )} 
								-\frac{\hat{C}_n}{\tilde{V}( \hat{t}_n, \hat{x}_n )} 
								-\frac{\tilde{C}_n}{\tilde{V}( t_n, x_n )}
							\bigg ) 
					\bigg \|_{\boldsymbol{L}(\mathbb{H}_{\kappa}, \mathbb{H}_{-\kappa})}\\ \nonumber
				&+\|
						B( t_n, x_n ) - B( \hat{t}_n, \hat{x}_n )
					\|_{HS(\mathbb{U}, \mathbb{H}_{\kappa})}
					\| 
						I_{\mathbb{H}}^{-1} \hat{C}_n
					\|_{\boldsymbol{L}(\mathbb{H}_{\kappa}, \mathbb{H}_{-\kappa})}	\\  \nonumber
					&\cdot
					\big(
						\|
							B( t_n, x_n )
						\|_{HS(\mathbb{U}, \mathbb{H}_{\kappa})}
						+
						\|
							B( \hat{t}_n, \hat{x}_n )
						\|_{HS(\mathbb{U}, \mathbb{H}_{\kappa})}
					\big).
	\end{align}
	In addition 
	the definition of $\hat{C}_n$, $n \in \N$,
	\eqref{eq: def of tilde V},
	\eqref{eq: def of bold L},
	\eqref{eq: V bounded with respect to H},
	and
	\eqref{eq: R bound} 
	imply that for all $n \in \N$ it holds that
	\begin{align}
	\label{eq: C hat bound}
				& \| 
					I_{\mathbb{H}}^{-1} \hat{C}_n 
				\|_{\boldsymbol{L}(\mathbb{H}_{\kappa}, \mathbb{H}_{-\kappa})} \\ \nonumber
			={}&
				\big \|  
					I_{\mathbb{H}}^{-1} \big(
						n (I_\mathbb{H} \left( x_n - \hat{x}_n\right))
							\otimes ( D_\mathbb{H} \tilde{V} )( \hat{t}_n , \hat{x}_n ) 
						+
						n (D_\mathbb{H} \tilde{V} )( \hat{t}_n , \hat{x}_n )
							\otimes \left(n I_\mathbb{H} \left( x_n - \hat{x}_n \right) \right) \\ \nonumber
						&+
						\hat{r}_n 
							\, ( D^2_\mathbb{H} \tilde{V})( \hat{t}_n , \hat{x}_n )
					\big)
				\big \|_{\boldsymbol{L}(\mathbb{H}_{\kappa}, \mathbb{H}_{-\kappa})} \\ \nonumber
			\leq{}&
				\|  
					I_{\mathbb{H}}^{-1} \big(
						n (I_\mathbb{H} \left( x_n - \hat{x}_n\right))
							\otimes ( D_\mathbb{H} V )( \hat{x}_n )
					\big)
				\|_{\boldsymbol{L}(\mathbb{H}_{\kappa}, \mathbb{H}_{-\kappa})} \\ \nonumber
				&+
				\|	
					I_{\mathbb{H}}^{-1} \big(
						(D_\mathbb{H} V )( \hat{x}_n )
							\otimes \left(n I_\mathbb{H} \left( x_n - \hat{x}_n \right) \right)
					\big)
				\|_{\boldsymbol{L}(\mathbb{H}_{\kappa}, \mathbb{H}_{-\kappa})} 
				+
				\|
					\hat{r}_n
						\, I_{\mathbb{H}}^{-1} ( D^2_\mathbb{H} V)(\hat{x}_n ) 
				\|_{\boldsymbol{L}(\mathbb{H}_{\kappa}, \mathbb{H}_{-\kappa})} \\ \nonumber
			={}&
				\sup_{x,y \in H \backslash \{0\}} 
					\frac
						{
							|
								n \langle  x_n - \hat{x}_n, x \rangle_H \cdot
									\langle 
										I_{\mathbb{H}}^{-1} ( D_\mathbb{H} V )( \hat{x}_n ), y 
									\rangle_{H}
							|
						}
						{	
							\|x\|_{H_\kappa} \|y\|_{H_\kappa}
						}
				+ \\ \nonumber
				&\sup_{x,y \in H \backslash \{0\}} 
					\frac
						{
							|
								n \langle  x_n - \hat{x}_n, y \rangle_H \cdot
									\langle 
										I_{\mathbb{H}}^{-1} ( D_\mathbb{H} V )( \hat{x}_n ), x 
									\rangle_{H}
							|
						}
						{	
							\|x\|_{H_\kappa} \|y\|_{H_\kappa}
						} 
				+
				\|
					\hat{r}_n
						\, I_{\mathbb{H}}^{-1} ( D^2_\mathbb{H} V)(\hat{x}_n ) 
				\|_{\boldsymbol{L}(\mathbb{H}_{\kappa}, \mathbb{H}_{-\kappa})} \\ \nonumber
			={}&
				2 n \|  
						x_n - \hat{x}_n
				\|_{H_{-\kappa}}
				\cdot	
				\|
					I_{\mathbb{H}}^{-1} ( D_\mathbb{H} V )( \hat{x}_n )
				\|_{H_{-\kappa}}
				+
				|\hat{r}_n|
				\cdot
				\|
					( D^2_\mathbb{H} V)( \hat{x}_n ) 
				\|_{\boldsymbol{L}(\mathbb{H}_{\kappa}, \mathbb{H}'_{\kappa})} \\ \nonumber
			\leq{}&
				2 n (\sup_{m \in \N} \nicefrac {1} {|\lambda_m|})^{\nicefrac 12 + \kappa} 
					\|x_n- \hat{x}_n\|_{H_{1/2}}
					\cdot
					K(R)
				+R \,
					\|
						( D^2_\mathbb{H} V)( \hat{x}_n ) 
					\|_{\boldsymbol{L}(\mathbb{H}_{\kappa}, \mathbb{H}'_{\kappa})} 
	\end{align}
	and 
	the definition of $C_n$, $\hat{C}_n$, and of $\tilde{C}_n$, $n \in \N$,
	\eqref{eq: def of tilde V},
	\eqref{eq: V Lipschitz with respect to H -gamma},
	\eqref{eq: R bound}, 
	\eqref{eq: H1/2- norm bounds},
	and 
	imply that for all $n \in \N$ and all $\eps \in (0, \infty)$ it holds that
	\begin{align}
	\label{eq: C difference bound}
				&\left \|
					I_{\mathbb{H}}^{-1} 
						\left(
							\frac{C_n}{\tilde{V}( t_n, x_n )} 
							-\frac{\hat{C}_n}{\tilde{V}( \hat{t}_n, \hat{x}_n )} 
							-\frac{\tilde{C}_n}{\tilde{V}( t_n, x_n )}
						\right ) 
				\right \|_{\boldsymbol{L}(\mathbb{H}_{\kappa}, \mathbb{H}_{-\kappa})} \\ \nonumber
			={}&
				\big \|
					I_{\mathbb{H}}^{-1} 
					\big(
						n (I_\mathbb{H} \left( x_n - \hat{x}_n \right) )
							\otimes 
								\tfrac
									{(D_\mathbb{H} \tilde{V} )( t_n , x_n )}
									{\tilde{V}( t_n, x_n )}
						+ \tfrac{(D_\mathbb{H} \tilde{V} )( t_n , x_n )} {\tilde{V}( t_n, x_n )}
							\otimes \left( n I_\mathbb{H} \left( x_n - \hat{x}_n \right) \right) \\ \nonumber
						&+
						r_n 
							\, \tfrac
								{(D^2_\mathbb{H} \tilde{V} )( t_n , x_n )} 
								{\tilde{V}( t_n , x_n )}
						-n (I_\mathbb{H} \left( x_n - \hat{x}_n \right) )
							\otimes 
								\tfrac
									{(D_\mathbb{H} \tilde{V} )( \hat{t}_n , \hat{x}_n )}
									{\tilde{V}( \hat{t}_n , \hat{x}_n )} 
						- \tfrac
								{(D_\mathbb{H} \tilde{V} )( \hat{t}_n , \hat{x}_n )}
								{\tilde{V}( \hat{t}_n , \hat{x}_n )}
							\otimes \left(n I_\mathbb{H} \left( x_n - \hat{x}_n \right) \right) \\ \nonumber
						&-
						\hat{r}_n 
							\, \tfrac
								{(D^2_\mathbb{H} \tilde{V} )( \hat{t}_n , \hat{x}_n )}
								{\tilde{V}( \hat{t}_n , \hat{x}_n )}
						- (r_n - \hat{r}_n) 
							\, \tfrac
								{( D^2_\mathbb{H} \tilde{V})( t_n , x_n )}
								{\tilde{V}( t_n , x_n )}			
					\big )  
				\big \|_{\boldsymbol{L}(\mathbb{H}_{\kappa}, \mathbb{H}_{-\kappa})} \\ \nonumber
			={}&
				\Big \|
					I_{\mathbb{H}}^{-1} 
					\Big(
						n (I_\mathbb{H} \left( x_n - \hat{x}_n \right) )
							\otimes 
							\Big (
								\tfrac{(D_\mathbb{H} V )( x_n )}{V(  x_n )}
								-\tfrac
									{(D_\mathbb{H} V )( \hat{x}_n )}
									{V(\hat{x}_n )}
							\Big ) \\ \nonumber
						&+\Big(
							\tfrac{(D_\mathbb{H} V )(x_n )}{V (x_n )}
							-\tfrac
								{(D_\mathbb{H} V )(\hat{x}_n )}
								{V( \hat{x}_n )}
						\Big)
							\otimes \Big( n I_\mathbb{H} \Big( x_n - \hat{x}_n \Big) \Big) 
					+
						\hat{r}_n 
							\,\Big(
								\tfrac
									{( D^2_\mathbb{H} V)(x_n )}
									{V(x_n )}
								-\tfrac
									{(D^2_\mathbb{H} V )( \hat{x}_n )}
									{V( \hat{x}_n )}
							\Big) 
					\Big )
				\Big \|_{\boldsymbol{L}(\mathbb{H}_{\kappa}, \mathbb{H}_{-\kappa})} \\ \nonumber
			\leq{}&
				2 n \|x - x_n\|_{H_{-\kappa}} \cdot
					\big \|
						\tfrac
							{I_{\mathbb{H}}^{-1} ( D_\mathbb{H} V)( x_n )}
							{V (x_n )}
						-\tfrac
							{I_{\mathbb{H}}^{-1} (D_\mathbb{H} V )( \hat{x}_n )}
							{V( \hat{x}_n )}
					\big \|_{H_{-\kappa}}
					+|\hat{r}_n| \cdot
						\big \|
								\tfrac
									{( D^2_\mathbb{H} V)(x_n )}
									{V( x_n )}
								-\tfrac
									{(D^2_\mathbb{H} V )( \hat{x}_n )}
									{V( \hat{x}_n )}
						\big \|_{\boldsymbol{L}(\mathbb{H}_{\kappa}, \mathbb{H}'_{\kappa})} \\ \nonumber
			\leq{}&
				2 n \|x - x_n\|_{H_{-\kappa}}
				\cdot
				K(R) \, \|x_n- \hat{x}_n\|_{H_{1/2}} 
				+R \,
					\big \|
							\tfrac
								{( D^2_\mathbb{H} V)(x_n )}
								{V( x_n )}
							-\tfrac
								{(D^2_\mathbb{H} V )( \hat{x}_n )}
								{V( \hat{x}_n )}
					\big \|_{\boldsymbol{L}(\mathbb{H}_{\kappa}, \mathbb{H}'_{\kappa})} \\ \nonumber
			\leq{}&
				2 \sqrt{n} (
					\eps \cdot \sqrt{n} \, \|x_n- \hat{x}_n\|_{H_{1/2}}
					+ M_\eps \cdot \sqrt{n} \, \|x_n- \hat{x}_n\|_{H} 
				)
				\cdot
				K(R) \, \|x_n- \hat{x}_n\|_{H_{1/2}} 
			\\ \nonumber &
				+R \,
					\big \|
							\tfrac
								{( D^2_\mathbb{H} V)(x_n )}
								{V( x_n )}
							-\tfrac
								{(D^2_\mathbb{H} V )( \hat{x}_n )}
								{V( \hat{x}_n )}
					\big \|_{\boldsymbol{L}(\mathbb{H}_{\kappa}, \mathbb{H}'_{\kappa})} \\ \nonumber
			\leq{}&
				2 \eps \, n
				\, K(R) \cdot
				\|x_n- \hat{x}_n\|^2_{H_{1/2}}
				+
				2\sqrt{n} \,
				R \cdot K(R) \, 
				M_\eps \,
				\|x_n- \hat{x}_n\|_{H_{1/2}} \\ \nonumber
				&+R \,
					\big \|
							\tfrac
								{( D^2_\mathbb{H} V)(x_n )}
								{V( x_n )}
							-\tfrac
								{(D^2_\mathbb{H} V )( \hat{x}_n )}
								{V( \hat{x}_n )}
					\big \|_{\boldsymbol{L}(\mathbb{H}_{\kappa}, \mathbb{H}'_{\kappa})}.
	\end{align}
	Moreover, \eqref{eq: Lipschitz continuity of B},
	\eqref{eq: R bound},
	the monotonicity of $K$,
	and
	\eqref{eq: H1/2- norm bounds}
	yield  that for all $\eps \in (0, \infty)$
	and all $n \in \N$ it holds that
	\begin{equation}
	\label{eq: B difference bound}
		\begin{split}
				&\|
					B( t_n, x_n ) - B( \hat{t}_n, \hat{x}_n )
				\|_{HS(\mathbb{U}, \mathbb{H}_{\kappa})}
			\leq 
				(\sup_{m \in \N} \nicefrac {1} {|\lambda_m|})^{\nicefrac 12 + \kappa}
				\|
					B( t_n, x_n ) - B( \hat{t}_n, \hat{x}_n )
				\|_{HS(\mathbb{U}, \mathbb{H})} \\
			\leq{}&
				K(R) \cdot
				(\sup_{m \in \N} \nicefrac {1} {|\lambda_m|})^{\nicefrac 12 + \kappa}
				(
					\|x_n - \hat{x}_n\|_{H_{1/2 - \beta_2}}
					+ |t_n - \hat{t}_n|
				) \\
			\leq{}&
				K(R) \cdot
				(\sup_{m \in \N} \nicefrac {1} {|\lambda_m|})^{\nicefrac 12 + \kappa}
				(
					\eps \, \|x_n - \hat{x}_n\|_{H_{1/2}}
					+M_\eps \|x_n - \hat{x}_n\|_{H}
					+ |t_n - \hat{t}_n|
				)\\
			\leq{}&
				K(R) \cdot
				(\sup_{m \in \N} \nicefrac {1} {|\lambda_m|})^{\nicefrac 12 + \kappa}
				(
					\eps \, \|x_n - \hat{x}_n\|_{H_{1/2}}
					+ M_\eps \cdot \tfrac {R}{\sqrt{n}}
					+ |t_n - \hat{t}_n|
				).
		\end{split}
	\end{equation}
	Hence we obtain from \eqref{eq: H bound for F, B and v},
	\eqref{eq: B difference}, 
	\eqref{eq: C hat bound},
	\eqref{eq: C difference bound}, 
	and from \eqref{eq: B difference bound}
	that for all $n \in \N$ and all $\eps \in (0, \infty)$ it holds that 
	\begin{align}
	\label{eq: B1 inequality for small H1/2 diff}
				&\frac
					{
						\big \langle
							B( t_n, x_n ) ,
							I_{\mathbb{H}_\var}^{-1} \big(
								((C_n - \tilde{C}_n) \, B( t_n, x_n ))|_{H_\var}
							\big)
						\big \rangle_{HS(\mathbb{U}, \mathbb{H}_\var)} 					
					}
					{ \tilde{V}( t_n, x_n ) } \\ \nonumber
				&-\frac
					{
						\big \langle
							B( \hat{t}_n, \hat{x}_n ) ,
							I_{\mathbb{H}_\var}^{-1} \big(
								(\hat{C}_n \, B( \hat{t}_n, \hat{x}_n ))|_{H_\var}
							\big)
						\big \rangle_{HS(\mathbb{U}, \mathbb{H}_\var)}
					}
					{ \tilde{V}( \hat{t}_n, \hat{x}_n ) } \\ \nonumber
			\leq{}&
				(K(R))^2 \, \Big(
					2 n
					\|x_n- \hat{x}_n\|_{H_{-\kappa}}
					\cdot
					K(R) \, \|x_n- \hat{x}_n\|_{H_{1/2}} 
					+R \,
						\big \|
								\tfrac
									{( D^2_\mathbb{H} V)(x_n )}
									{V( x_n )}
								-\tfrac
									{(D^2_\mathbb{H} V )( \hat{x}_n )}
									{V( \hat{x}_n )}
						\big \|_{\boldsymbol{L}(\mathbb{H}_{\kappa}, \mathbb{H}'_{\kappa})}
				\Big) \\ \nonumber
				&+2(K(R))^2 \cdot
				(\sup_{m \in \N} \nicefrac {1} {|\lambda_m|})^{\nicefrac 12 + \kappa}
				\|x_n - \hat{x}_n\|_{H_{1/2 - \beta_2}} 
				\\  \nonumber
				& \cdot 
				\Big(
					2 n (\sup_{m \in \N} \nicefrac {1} {|\lambda_m|})^{1 + 2\kappa} 
					\|x_n- \hat{x}_n\|_{H_{1/2}}
					\cdot K(R)
					+R \,
						\|
							( D^2_\mathbb{H} V)( \hat{x}_n ) 
						\|_{\boldsymbol{L}(\mathbb{H}_{\kappa}, \mathbb{H}'_{\kappa})} 
				\Big)
	\end{align}
	and that
	\begin{align}
	\nonumber
				&\frac
					{
						\big \langle
							B( t_n, x_n ) ,
							I_{\mathbb{H}_\var}^{-1} \big(
								((C_n - \tilde{C}_n) \, B( t_n, x_n ))|_{H_\var}
							\big)
						\big \rangle_{HS(\mathbb{U}, \mathbb{H}_\var)} 					
					}
					{ \tilde{V}( t_n, x_n ) } 
			\\ \nonumber &
				-\frac
					{
						\big \langle
							B( \hat{t}_n, \hat{x}_n ) ,
							I_{\mathbb{H}_\var}^{-1} \big(
								(\hat{C}_n \, B( \hat{t}_n, \hat{x}_n ))|_{H_\var}
							\big)
						\big \rangle_{HS(\mathbb{U}, \mathbb{H}_\var)}
					}
					{ \tilde{V}( \hat{t}_n, \hat{x}_n ) } 
			\\ \label{eq: B1 inequality for big H1/2 diff} \leq{}&
				(K(R))^2 \big(
					2 \eps \, n
					\, K(R) \cdot
					\|x_n- \hat{x}_n\|^2_{H_{1/2}}
					+
					2\sqrt{n} \,
					R \cdot K(R) \, 
					M_\eps \,
					\|x_n- \hat{x}_n\|_{H_{1/2}}
				\\ \nonumber &
					+R \,
						\big \|
								\tfrac
									{( D^2_\mathbb{H} V)(x_n )}
									{V( x_n )}
								-\tfrac
									{(D^2_\mathbb{H} V )( \hat{x}_n )}
									{V( \hat{x}_n )}
						\big \|_{\boldsymbol{L}(\mathbb{H}_{\kappa}, \mathbb{H}'_{\kappa})}
				\big) 
		\\ \nonumber &
			+2 (K(R))^2 \cdot
				(\sup_{m \in \N} \nicefrac {1} {|\lambda_m|})^{\nicefrac 12 + \kappa}
				(
					\eps \, \|x_n - \hat{x}_n\|_{H_{1/2}}
					+ M_\eps \cdot \tfrac {R}{\sqrt{n}}
				) 
				\, \sqrt{n} 
		\\ \nonumber & \cdot 
				\big(
					2 \sqrt{n} (\sup_{m \in \N} \nicefrac {1} {|\lambda_m|})^{1 + 2\kappa} 
					\|x_n- \hat{x}_n\|_{H_{1/2}}
					\cdot K(R)
					+\tfrac{R}{\sqrt{n}} \,
					\|
						( D^2_\mathbb{H} V)( \hat{x}_n ) 
					\|_{\boldsymbol{L}(\mathbb{H}_{\kappa}, \mathbb{H}'_{\kappa})}
				\big). 
	\end{align}
	Furthermore, \eqref{eq: different HS norm},
	the properties of $\mathfrak{C}$ and $\hat{\mathfrak{C}}$,
	\eqref{eq: Lipschitz continuity of B},
	\eqref{eq: R bound},
	the monotonicity of $K$,
	\eqref{eq: H1/2- norm bounds},
	and
	the fact that for all $a, b \in \R$ it holds that
	$(a+b)^2 \leq 2a^2+2b^2$
	show that for all $n \in \N$ and all $\eps \in (0, \infty)$ it holds that 
	\begin{align}
			\nonumber
				&n
					\big[
						\big \langle
							B( t_n, x_n ), 
							I_{\mathbb{H}_\var}^{-1} \big(
								(\mathfrak{C}_n \, B( t_n, x_n ) )|_{H_\var}
							\big)
						\big \rangle_{HS(\mathbb{U}, \mathbb{H}_\var)} \\ \nonumber
						&-	
						\big \langle
							B( \hat{t}_n, \hat{x}_n ),
							I_{\mathbb{H}_\var}^{-1} \big(
								(\hat{\mathfrak{C}}_n \, B( \hat{t}_n, \hat{x}_n ))|_{H_\var}
							\big)
						\big \rangle_{HS(\mathbb{U}, \mathbb{H}_\var)}
					\big] \\ \nonumber
				\leq{}&
					n
					\sum_{ i \in \mathcal{I}} \left[
						\left \langle
							B( t_n, x_n ) \, \tilde{e}_i ,
							I_\mathbb{H}^{-1} \mathfrak{C}_n \, B( t_n, x_n ) \, \tilde{e}_i
						\right \rangle_{H}
						-
						\left \langle
							B( \hat{t}_n, \hat{x}_n ) \, \tilde{e}_i ,
							I_\mathbb{H}^{-1} \hat{\mathfrak{C}}_n \, 
								B( \hat{t}_n, \hat{x}_n ) \, \tilde{e}_i
						\right \rangle_{H}
					\right] \\ \label{eq: B2 inequality}
				\leq{}&
					n
						\sum_{ i\in \mathcal{I}} 
						3 \,
						\big\|
							B( t_n, x_n ) \, \tilde{e}_i
							-
							B( \hat{t}_n, \hat{x}_n ) \, \tilde{e}_i 
						\big\|_{H}^2 
				= 
					3
						n
						\,
						\|
							B( t_n, x_n ) 
							-
							B( \hat{t}_n, \hat{x}_n )  
						\|_{
							HS( \mathbb{U}, \mathbb{H} )
						}^2 \\ \nonumber
				\leq{}&
					3K(R) \cdot n \, 
						(\|x_n -\hat{x_n}\|_{H_{1/2 - \beta_2}} + |t_n- \hat{t}_n| )^2 \\ \nonumber
				\leq{}&
					3K(R) \cdot n \,
					(
						\sqrt{\eps} \|x_n -\hat{x_n}\|_{H_{1/2}} 
						+ M_{\sqrt{\eps}} \, \|x_n -\hat{x_n}\|_{H} 
						+ |t_n- \hat{t}_n| 
					)^2 \\ \nonumber
				\leq{}&
					6 \eps \, n \, K(R) \cdot \|x_n -\hat{x}_n\|^2_{H_{1/2}}
					+ 6K(R) \, (M_{\sqrt{\eps}} \cdot R + R)^2. 
	\end{align}	
	Therefore
	\eqref{eq: v difference inequality},
	\eqref{eq: A difference inequality},
	\eqref{eq: F difference inequality},
	\eqref{eq: B1 inequality for big H1/2 diff},
	\eqref{eq: B2 inequality},
	\eqref{eq: x vartheta bound 1},
	\eqref{eq: x vartheta bound 2},
	\eqref{eq: x vartheta bound 3},
	\eqref{eq: R bound},
	together with  the fact that 
	$\tilde{V}$ is bounded on $\R \times \mathbb{H}$- bounded subsets of $(0,T) \times O$,
	the fact that
	$
		D_{\mathbb{H}}^2 V
	$
	is $\boldsymbol{L}(\mathbb{H}_{\kappa}, \mathbb{H}'_{\kappa})$-bounded 
	on $\mathbb{H}$- bounded subsets of $O$,
	the fact that for all $n \in \N$ it holds that $n^{\nicefrac 14} \leq \sqrt{n} + 1$,
	and the fact that $\sup_{x \in [0,\infty)} \omega (x) < \infty$
	imply that 
	there exists
	$R_1 \in (1,\infty)$ such that
	for all
	$\eps \in (0, \nicefrac {1}{(4R_1)})$ 
	and all
	$ n \in \N \cap (64(R_1)^2, \infty)$ 
	with 
	$\delta_n (\|x_n\|_{H_\var}^2 \vee \|\hat{x}_n\|_{H_\var}^2) \leq 1$
	and with $\delta_n < 1$
	there exist
	$R_{2, \eps}, R_{3, \eps} \in (0,\infty)$
	such that 
	\begin{align}
	\label{eq: diff bound for big H1/2 difference}
					&\tfrac{v(t_n,x_n)}{\tilde{V}(t_n, x_n)} r_n
						-\tfrac{v(\hat{t}_n,\hat{x}_n)}{\tilde{V}(\hat{t}_n,\hat{x}_n)} \hat{r}_n
					+\tfrac{\langle A(x_n), p_n \rangle_{H, H'}}
						{\tilde{V}(t_n, x_n)}
					-\tfrac{\langle A(\hat{x}_n), \hat{p}_n \rangle_{H, H'}}
						{\tilde{V}(\hat{t}_n, \hat{x}_n)} 
					+\tfrac{\langle F(t_n, x_n), p_n \rangle_{H, H'}}
						{\tilde{V}(t_n, x_n)} 
					-\tfrac{\langle F(\hat{t}_n, \hat{x}_n), \hat{p}_n \rangle_{H, H'}}
						{\tilde{V}(\hat{t}_n, \hat{x}_n)} \\ \nonumber
					&+\tfrac
							{
								\big \langle
									B( t_n, x_n ) ,
									I_{\mathbb{H}_\var}^{-1} \big(
										(C_n \, B( t_n, x_n ))|_{H_\var}
									\big)
								\big \rangle_{HS(\mathbb{U}, \mathbb{H}_\var)} 					
							}
							{ \tilde{V}( t_n, x_n ) } 
						-\tfrac
							{
								\big \langle
									B( \hat{t}_n, \hat{x}_n ) ,
									I_{\mathbb{H}_\var}^{-1} \big(
										(\hat{C}_n \, B( \hat{t}_n, \hat{x}_n ))|_{H_\var}
									\big)
								\big \rangle_{HS(\mathbb{U}, \mathbb{H}_\var)}
							}
							{ \tilde{V}( \hat{t}_n, \hat{x}_n ) } 
				\\ \nonumber &
					+n
							\Big[
								\big \langle
									B( t_n, x_n ), 
									I_{\mathbb{H}_\var}^{-1} \big(
										(\mathfrak{C}_n \, B( t_n, x_n ) )|_{H_\var}
									\big)
								\big \rangle_{HS(\mathbb{U}, \mathbb{H}_\var)} 
					\\ \nonumber &
								-	
								\big \langle
									B( \hat{t}_n, \hat{x}_n ),
									I_{\mathbb{H}_\var}^{-1} \big(
										(\hat{\mathfrak{C}}_n \, B( \hat{t}_n, \hat{x}_n ))|_{H_\var}
									\big)
								\big \rangle_{HS(\mathbb{U}, \mathbb{H}_\var)}
							\Big] 
					-\tfrac{v(t_n,x_n)}{\tilde{V}(t_n, x_n)} (r_n -\hat{r}_n) 
				\\ \nonumber &
					-\tfrac {\left \langle A(x_n), \tilde{p}_n \right \rangle_{H, H'}}
						{\tilde{V}(t_n, x_n)} 
					-\tfrac {\left \langle F(t_n, x_n), \tilde{p}_n \right \rangle_{H, H'}}
						{\tilde{V}(t_n, x_n)} 
				\\ \nonumber &
					-\Big \langle
									B( t_n, x_n ),
									I_{\mathbb{H}_\var}^{-1} 
										\Big( \big(
											\tfrac{\tilde{C}_n}{\tilde{V}( t_n, x_n )}
												\, B( t_n, x_n )
										\big) \Big|_{H_\var} \Big )
							\Big \rangle_{HS(\mathbb{U}, \mathbb{H}_{\var})} 
				\\ \nonumber
			\leq{}&
				-n \|x_n -\hat{x}_n\|^2_{H_{1/2}}
				+\|x_n -\hat{x}_n\|^2_{H_{1/2}} R_{1} (1 +\sqrt{n} + n \eps) 
				+\|x_n -\hat{x}_n\|_{H_{1/2}} R_{2 ,\eps} (
					1 + \sqrt{n}
				) 
				+R_{3, \eps} \\ \nonumber
			={}&
				- \nicefrac n2 \|x_n -\hat{x}_n\|^2_{H_{1/2}}
					-(n \, (\nicefrac 12- R_1 \eps) - \sqrt{n} R_1 -R_1) \,
					\big (
						\|x_n -\hat{x}_n\|_{H_{1/2}} 
						- \tfrac
							{
								R_{2 ,\eps} (1 + \sqrt{n})
							}
							{
								(n \, (1- 2 R_1 \eps) - 2\sqrt{n} R_1- 2R_1)
							}
					\big )^2 
			\\ \nonumber & 
				+(R_{2 ,\eps})^2 \tfrac
					{
						(
							1 + \sqrt{n}
						)^2
					}
					{2 (n \, (1- 2 R_1 \eps) - 2\sqrt{n} R_1- 2R_1)} 
				+R_{3, \eps} \\ \nonumber
			\leq{}&
				- \nicefrac n2 \|x_n -\hat{x}_n\|^2_{H_{1/2}}
				+(R_{2 ,\eps})^2 \frac
					{
							2 + 2n
					}
					{2 (\tfrac{n}{4}- 2R_1)} 
				+R_{3, \eps} \\ \nonumber
			\leq{}
				&- \nicefrac n2 \|x_n -\hat{x}_n\|^2_{H_{1/2}}
				+12(R_{2 ,\eps})^2
				+R_{3, \eps} 
	\end{align}
	and this 
	together with  the fact that 
	$
			\lim_{n \to \infty} \delta_n (\|x_n\|_{H_\var}^2 \vee \|\hat{x}_n\|_{H_\var}^2) 
		=
			\lim_{n \to \infty} \delta_n
		= 0
	$
	yields that
	\begin{align}
	\label{eq: limit bound for big H1/2 diff}
				&\limsup_{L \to \infty} \limsup_{n \to \infty} \Bigg(
					\1_{n \| x_n - \hat{x}_n \|^2_{H_{1/2}} > L } \bigg(
						\frac{v(t_n,x_n)}{\tilde{V}(t_n, x_n)} r_n
							-\frac{v(\hat{t}_n,\hat{x}_n)}{\tilde{V}(\hat{t}_n,\hat{x}_n)} \hat{r}_n
						+\frac{\langle A(x_n), p_n \rangle_{H, H'}}
							{\tilde{V}(t_n, x_n)}
				~\\ \nonumber &\qquad 
						-\frac{\langle A(\hat{x}_n), \hat{p}_n \rangle_{H, H'}}
							{\tilde{V}(\hat{t}_n, \hat{x}_n)} 
						+\frac{\langle F(t_n, x_n), p_n \rangle_{H, H'}}
							{\tilde{V}(t_n, x_n)} 
						-\frac{\langle F(\hat{t}_n, \hat{x}_n), \hat{p}_n \rangle_{H, H'}}
							{\tilde{V}(\hat{t}_n, \hat{x}_n)}
					~\\ \nonumber &\qquad 
						+\frac
								{
									\langle
										B( t_n, x_n ) ,
										I_{\mathbb{H}_\var}^{-1}  \big(
											(C_n \, B( t_n, x_n )) |_{H_\var}
										\big)
									\rangle_{HS(\mathbb{U}, \mathbb{H}_\var)} 					
								}
								{ \tilde{V}( t_n, x_n ) } 
					\\ \nonumber & \qquad
						-\frac
							{
								\langle
									B( \hat{t}_n, \hat{x}_n ) ,
									I_{\mathbb{H}_\var}^{-1} \big(
										(\hat{C}_n \, B( \hat{t}_n, \hat{x}_n )) |_{H_\var}
									\big)
								\rangle_{HS(\mathbb{U}, \mathbb{H}_\var)}
							}
							{ \tilde{V}( \hat{t}_n, \hat{x}_n ) } 
						-\frac{v(t_n,x_n)}{\tilde{V}(t_n, x_n)} (r_n -\hat{r}_n) \\ \nonumber
				& \qquad
						+n
							\Big[
								\langle
									B( t_n, x_n ), 
									I_{\mathbb{H}_\var}^{-1} \big(
										(\mathfrak{C}_n \, B( t_n, x_n ))|_{H_\var}
									\big)
								\rangle_{HS(\mathbb{U}, \mathbb{H}_\var)}
						~\\ \nonumber &\qquad 
								-	
								\langle
									B( \hat{t}_n, \hat{x}_n ),
									I_{\mathbb{H}_\var}^{-1} \big(
										(\hat{\mathfrak{C}}_n \, B( \hat{t}_n, \hat{x}_n ))|_{H_\var}
									\big)
								\rangle_{HS(\mathbb{U}, \mathbb{H}_\var)}
							\Big] 
						-\frac {\left \langle A(x_n), \tilde{p}_n \right \rangle_{H, H'}}
							{\tilde{V}(t_n, x_n)}
					\\ \nonumber
					& \qquad
						-\frac {\left \langle F(t_n, x_n), \tilde{p}_n \right \rangle_{H, H'}}
							{\tilde{V}(t_n, x_n)} 
						-\bigg \langle
										B( t_n, x_n ),
										I_{\mathbb{H}_\var}^{-1} 
										\Big(
											\Big(
												\frac{\tilde{C}_n}{\tilde{V}( t_n, x_n )} \, B( t_n, x_n )
											\Big) \Big |_{H_\var}
										\Big)
								\bigg \rangle_{HS(\mathbb{U}, \mathbb{H}_{\var})} 
						\bigg)
					\Bigg)
			\leq 0.
	\end{align}
	In addition, it follows from
	\eqref{eq: v difference inequality},
	\eqref{eq: x vartheta bound 3},
	\eqref{eq: def of tilde V},
	the fact that 
	$
			\lim_{n \to \infty} (\sqrt{n} \|x_n - \hat{x}_n\|_H)
		=
			\lim_{n \to \infty} ( \sqrt{n} |t_n - \hat{t}_n|)
		=
			\lim_{n \to \infty} ( \delta_n (\|x_n\|^2_{H_\var} \vee \|\hat{x}_n\|^2_{H_\var}))
		=
			\lim_{n \to \infty} \delta_n
		=
			\lim_{n \to \infty} \omega(\nicefrac 1n )
		= 0,
	$
	the fact that 
	$
		\sup_{x \in (0,\infty)} \omega(x) < \infty,
	$
	and from the fact that
	$V$ is Lipschitz continuous with respect to the $\| \cdot \|_{\mathbb{H}}$-norm
	on $\mathbb{H}$-bounded subsets of $O$
	that for all $L \in (0, \infty)$ it holds that
	\begin{align}
	\label{eq: v limit bound for small H1/2 diff}
				&\limsup_{n \to \infty} \bigg(
					\1_{n \| x_n - \hat{x}_n \|^2_{H_{1/2}} \leq L } \bigg(
						\frac{v(t_n,x_n)}{\tilde{V}(t_n, x_n)} r_n
						-\frac{v(\hat{t}_n,\hat{x}_n)}{\tilde{V}(\hat{t}_n,\hat{x}_n)} \hat{r}_n 
						-\frac{v(t_n,x_n)}{\tilde{V}(t_n, x_n)} (r_n -\hat{r}_n)
					\bigg)
				\bigg)\\ \nonumber
			\leq{}&
				\limsup_{n \to \infty} \bigg(
					\1_{n \| x_n - \hat{x}_n \|^2_{H_{1/2}} \leq L } \bigg(
						K(\|x\|^2_{H_\var} \vee \|\hat{x}\|^2_{H_\var}) \, R 
						\cdot
							(\|x- \hat{x}\|^2_{H_{1/2}} + \omega(\|x- \hat{x}\|^2_{H_{1/2}} + |t -\hat{t}|^2))
				\\ \nonumber & \qquad
						+ R \cdot K(R) 
						\,
						\left |
							\tilde{V}(\hat{t}_n,\hat{x}_n)
							-\tilde{V}(t_n,x_n)
						\right |
					\bigg)
				\bigg) \\ \nonumber
			\leq{}&
				\limsup_{n \to \infty} \bigg(
					\1_{n \| x_n - \hat{x}_n \|^2_{H_{1/2}} \leq L } \bigg(
						\big(
							n^{\nicefrac 14} \|x_n-\hat{x}_n\|^2_{H_{1/2}}
							+((\sup_{m \in [n,\infty)}\omega(m^{- \nicefrac 14} ))^{\nicefrac 12}
				\\ \nonumber & \qquad
							+\big( \sup_{m \in [0,\infty)} \omega(m ) \big)
								\cdot n^{\nicefrac 12} \, (\|x_n-\hat{x}_n\|^2_{H_{1/2}}+|t_n-\hat{t}_n|^2)
						\big)
							\, R
				\\ \nonumber & \qquad
						+ R \cdot K(R) 
						\,
						\left |
							\tilde{V}(\hat{t}_n,\hat{x}_n)
							-\tilde{V}(t_n,x_n)
						\right |
					\bigg)
				\bigg)
			=
				0.
	\end{align}
	Moreover, we get from 
	\eqref{eq: A difference inequality},
	\eqref{eq: F difference inequality},
	\eqref{eq: x vartheta bound 1},
	\eqref{eq: x vartheta bound 2},
	the fact that 
	$
			\lim_{n \to \infty} (\sqrt{n} \|x_n - \hat{x}_n\|_H)
		=
			0
	$, that
	$
			\lim_{n \to \infty} ( \sqrt{n} |t_n - \hat{t}_n|)
		= 0
	$, that
	$
			\lim_{n \to \infty} ( \delta_n (\|x_n\|^2_{H_\var} \vee \|\hat{x}_n\|^2_{H_\var}))
		= 0
	$, that
	$
			\lim_{n \to \infty} \delta_n
		= 0
	$, that
	$
			\lim_{n \to \infty} \omega(\nicefrac 1n )
		= 0
	$,
	the fact that 
	$
		\sup_{x \in (0,\infty)} \omega(x) < \infty,
	$
	and from the fact that 
	$
		\forall (z_n)_{n \in \N} \subseteq H_2 \colon 
			(
				(\lim_{n \to \infty} \|z_n\|_H = 0) 
				\wedge (\limsup_{n \to \infty} \|z_n\|_{H_{1/2}} < \infty)
			)
			\Rightarrow
			(\lim_{n \to \infty} \|z_n\|_{H_{1/2- \alpha_2}} = 0)
	$
	that for all $L \in (0, \infty)$ it holds that
	\begin{align}
	\nonumber
				&\limsup_{n \to \infty} \bigg(
					\1_{n \| x_n - \hat{x}_n \|^2_{H_{1/2}} \leq L } \bigg(
						\frac{\langle A(x_n), p_n \rangle_{H, H'}}
							{\tilde{V}(t_n, x_n)}
						-\frac{\langle A(\hat{x}_n), \hat{p}_n \rangle_{H, H'}}
							{\tilde{V}(\hat{t}_n, \hat{x}_n)} 
						-\frac{\langle A(x_n), \tilde{p}_n \rangle_{H, H'}}
							{\tilde{V}(t_n, x_n)}\\ \nonumber
					& \qquad 
						+\frac{\langle F(t_n, x_n), p_n \rangle_{H, H'}}
							{\tilde{V}(t_n, x_n)} 
						-\frac{\langle F(\hat{t}_n, \hat{x}_n), \hat{p}_n \rangle_{H, H'}}
							{\tilde{V}(\hat{t}_n, \hat{x}_n)}
						-\frac{\langle F(t_n, x_n), \tilde{p}_n \rangle_{H, H'}}
							{\tilde{V}(t_n, x_n)}
						\bigg)						
					\bigg)\\ \nonumber
			\leq{}&
				\limsup_{n \to \infty} \bigg(
					\1_{n \| x_n - \hat{x}_n \|^2_{H_{1/2}} \leq L } \bigg(
						-n \|x_n - \hat{x}_n \|^2_{H_{1/2}}
						+R \cdot K(\| \hat{x}_n \|^2_{H_{\var}})
							\|x_n - \hat{x}_n \|_{H_{1/2}} 
				\\ \nonumber & \quad
						+ R \cdot K(\| x_n \|^2_{H_{\var}} \vee \| \hat{x}_n \|^2_{H_{\var}}) 
							\cdot \| x_n \|_{H_{\var}} \,
							\big (
								\| x_n - \hat{x}_n \|^2_{H_{1/2}} 
								+ \omega(\| x_n - \hat{x}_n \|^2_{H_{1/2}})
							\big )
				\\\nonumber  & \quad
						+n \, K (R) \, (\|x_n - \hat{x}_n \|_{H_{1/2 - \alpha_2}} + |t_n- \hat{t}_n|)
							\cdot \|x_n - \hat{x}_n\|_{H_{1/2}}
				\\
	\label{eq: A+F limit bound for small H1/2 diff}
				& \quad
						+R \cdot K (R) \cdot K(\| \hat{x}_n \|^2_{H_{\var}}) \, 
							(
									\|x_n - \hat{x}_n \|_{H_{1/2-\alpha_2}} 
								+ |t_n- \hat{t}_n|
							)
				\\ \nonumber & \quad
						+ R \cdot (K(\| x_n \|^2_{H_{\var}} \vee \| \hat{x}_n \|^2_{H_{\var}}))^2
							\, (\|x_n-\hat{x}_n\|^2_{H_{1/2}}+\omega(\|x_n-\hat{x}_n\|^2_{H_{1/2}}))
					\bigg)
				\bigg) \\ \nonumber
			\leq{}&
				\limsup_{n \to \infty} \bigg(
					\1_{n \| x_n - \hat{x}_n \|^2_{H_{1/2}} \leq L } \bigg(
						-n \|x_n - \hat{x}_n \|^2_{H_{1/2}}
						+R \cdot n^{\nicefrac 14} \|x_n-\hat{x}_n\|_{H_{1/2}} 
				\\ \nonumber & \quad
						+R \cdot K(R) \, n^{\nicefrac 14} 
							(\|x_n-\hat{x}_n\|_{H_{1/2-\alpha_2}} +|t_n- \hat{t}_n|)
						+ \big (
							n^{\nicefrac 14} \|x_n-\hat{x}_n\|^2_{H_{1/2}}
				\\ \nonumber & \quad
							+((\sup_{m \in [n,\infty)}\omega(m^{- \nicefrac 14} ))^{\nicefrac 12}
							+\big( \sup_{m \in [0,\infty)} \omega(m ) \big)
								\cdot n^{\nicefrac 12} \, \|x_n-\hat{x}_n\|_{H_{1/2}}
						\big )
						\cdot 
							2R	
				\\ \nonumber & \quad
						+n \, K (R) \, (\|x_n - \hat{x}_n \|_{H_{1/2 - \alpha_2}} + |t_n- \hat{t}_n|)
							\cdot \|x_n - \hat{x}_n\|_{H_{1/2}}
					\bigg)
				\bigg) 
			\leq{} 0.
	\end{align}
	Furthermore, \eqref{eq: B1 inequality for small H1/2 diff},
	the fact that 
	$
			\lim_{n \to \infty} (\sqrt{n} \|x_n - \hat{x}_n\|_H)
		=
			\lim_{n \to \infty} ( \sqrt{n} |t_n - \hat{t}_n|)
		=0
	$,
	the fact that $D_{\mathbb{H}}^2 \, V $ 
	is uniformly continuous with respect to the
	$\| \cdot \|_{\mathbb{H}}$ and the
	$\| \cdot \|_{\boldsymbol{L}(\mathbb{H}_{\kappa}, \mathbb{H}_{\kappa}')}$-norm
	on all $\mathbb{H}$-bounded subsets of $O$,
	and from the fact that 
	$
		\forall (z_n)_{n \in \N} \subseteq H_2 \colon 
			(
				(\lim_{n \to \infty} \|z_n\|_H = 0) 
				\wedge (\limsup_{n \to \infty} \|z_n\|_{H_{1/2}} < \infty)
			)
			\Rightarrow
			(\lim_{n \to \infty} (\|z_n\|_{H_{1/2- \beta_2}} + \|z_n\|_{H_{-\kappa}}) = 0)
	$
	show
	that for all $L \in (0, \infty)$ it holds that
	\begin{align}
	\label{eq: B1 limit bound for small H1/2 diff}
				&\limsup_{n \to \infty} \bigg(
					\1_{n \| x_n - \hat{x}_n \|^2_{H_{1/2}} \leq L } \bigg(
						\frac
								{
									\big \langle
										B( t_n, x_n ) ,
										I_{\mathbb{H}_\var}^{-1} \big(
											((C_n - \tilde{C}_n) \, B( t_n, x_n ) )|_{H_\var}
										\big)
									\big \rangle_{HS(\mathbb{U}, \mathbb{H}_\var)} 					
								}
								{ \tilde{V}( t_n, x_n ) } \\ \nonumber
					& \qquad \qquad
						-\frac
							{
								\langle
									B( \hat{t}_n, \hat{x}_n ) ,
									I_{\mathbb{H}_\var}^{-1} \big(
										(\hat{C}_n \, B( \hat{t}_n, \hat{x}_n ))|_{H_\var}
									\big)
								\rangle_{HS(\mathbb{U}, \mathbb{H}_\var)}
							}
							{ \tilde{V}( \hat{t}_n, \hat{x}_n ) } 
					\bigg)
				\bigg)\\ \nonumber
			\leq{}&
				\limsup_{n \to \infty} \bigg(
					\1_{n \| x_n - \hat{x}_n \|^2_{H_{1/2}} \leq L } \bigg(
						K^2(R) \, 
						\Big(
							2 
							\sqrt{n} \,
							\|x_n- \hat{x}_n\|_{H_{-\kappa}}
							\cdot
							K(R) \, \sqrt{n} \, \|x_n- \hat{x}_n\|_{H_{1/2}} \\ \nonumber
				& \quad
							+R \,
								\big \|
									\tfrac
										{( D^2_\mathbb{H} V)(x_n )}
										{V( x_n )}
									-\tfrac
										{(D^2_\mathbb{H} V )( \hat{x}_n )}
										{V( \hat{x}_n )}
								\big \|_{\boldsymbol{L}(\mathbb{H}_{\kappa}, \mathbb{H}'_{\kappa})}
						\Big) 
						+2K^2(R) \cdot
						(\sup_{m \in \N} \nicefrac {1} {|\lambda_m|})^{\nicefrac 12 + \kappa}
						\sqrt{n} \,
						\|x_n - \hat{x}_n\|_{H_{1/2 - \beta_2}} \\ \nonumber
			& \quad
						\cdot 
						\Big(
							2 (\sup_{m \in \N} \nicefrac {1} {|\lambda_m|})^{1 + 2\kappa} 
							\sqrt{n} \,
							\|x_n- \hat{x}_n\|_{H_{1/2}}
							\cdot K(R)
							+\frac {R}{\sqrt{n}} \,
							\|
								( D^2_\mathbb{H} V)( \hat{x}_n ) 
							\|_{\boldsymbol{L}(\mathbb{H}_{\kappa}, \mathbb{H}'_{\kappa})} 
						\Big)
					\Bigg)
				\Bigg) \\ \nonumber
			={}&
				0
	\end{align}
	and \eqref{eq: B2 inequality} together with
	the fact that 
	$
			\lim_{n \to \infty} (\sqrt{n} \|x_n - \hat{x}_n\|_H)
		=
			\lim_{n \to \infty} ( \sqrt{n} |t_n - \hat{t}_n|)
		=0
	$
	and with the fact that 
	$
		\forall (z_n)_{n \in \N} \subseteq H_2 \colon 
			(
				(\lim_{n \to \infty} \|z_n\|_H = 0) 
				\wedge (\limsup_{n \to \infty} \|z_n\|_{H_{1/2}} < \infty)
			)
			\Rightarrow
			(\lim_{n \to \infty} \|z_n\|_{H_{1/2- \beta_2}} = 0)
	$
	shows that for all $L \in (0, \infty)$ it holds that
	\begin{align}
	\nonumber
				&\limsup_{n \to \infty} \bigg(
					\1_{\| x_n - \hat{x}_n \|^2_{H_{1/2}} \leq L } \bigg(
						n
							\Big[
								\big \langle
									B( t_n, x_n ), 
									I_{\mathbb{H}_\var}^{-1} \big(
										(\mathfrak{C}_n \, B( t_n, x_n ))|_{H_\var} 
									\big)
								\big \rangle_{HS(\mathbb{U}, \mathbb{H}_\var)} \\ 
	\label{eq: B2 limit bound for small H1/2 diff}
					& \qquad \qquad
								-	
								\big \langle
									B( \hat{t}_n, \hat{x}_n ),
									I_{\mathbb{H}_\var}^{-1} \big(
										(\hat{\mathfrak{C}}_n \, B( \hat{t}_n, \hat{x}_n ))|_{H_\var}
									\big)
								\big \rangle_{HS(\mathbb{U}, \mathbb{H}_\var)}
							\Big] 
						\bigg)
					\bigg)\\ \nonumber
			\leq{}&
				\limsup_{n \to \infty} \bigg(
					\1_{\| x_n - \hat{x}_n \|^2_{H_{1/2}} \leq L } \bigg(
						3K(R) \, 
						( 
							\sqrt{n} \, \|x_n -\hat{x_n}\|_{H_{1/2 - \beta_2}} 
							+ \sqrt{n} \, |t_n- \hat{t}_n| 
						)^2
					\bigg)
				\bigg)
			=
				0.
	\end{align}
	Combining 
	\eqref{eq: v limit bound for small H1/2 diff},
	\eqref{eq: A+F limit bound for small H1/2 diff},
	\eqref{eq: B1 limit bound for small H1/2 diff},
	\eqref{eq: B2 limit bound for small H1/2 diff}
	implies that for all $L \in  (0, \infty)$ it holds that
	\begin{align}
			&\limsup_{n \to \infty} \Bigg(
					\1_{n \| x_n - \hat{x}_n \|^2_{H_{1/2}} \leq L } \Bigg(
						\frac{v(t_n,x_n)}{\tilde{V}(t_n, x_n)} r_n
							-\frac{v(\hat{t}_n,\hat{x}_n)}{\tilde{V}(\hat{t}_n,\hat{x}_n)} \hat{r}_n
						+\frac{\langle A(x_n), p_n \rangle_{H, H'}}
							{\tilde{V}(t_n, x_n)}
				\\ \nonumber & \qquad
						-\frac{\langle A(\hat{x}_n), \hat{p}_n \rangle_{H, H'}}
							{\tilde{V}(\hat{t}_n, \hat{x}_n)}
						+\frac{\langle F(t_n, x_n), p_n \rangle_{H, H'}}
							{\tilde{V}(t_n, x_n)} 
						-\frac{\langle F(\hat{t}_n, \hat{x}_n), \hat{p}_n \rangle_{H, H'}}
							{\tilde{V}(\hat{t}_n, \hat{x}_n)}
				\\\nonumber & \qquad
						+\frac
								{
									\big \langle
										B( t_n, x_n ) ,
										I_{\mathbb{H}_\var}^{-1} \big(
											(C_n \, B( t_n, x_n ))|_{H_\var}
										\big)
									\big \rangle_{HS(\mathbb{U}, \mathbb{H}_\var)} 					
								}
								{ \tilde{V}( t_n, x_n ) } 
					\\ \nonumber & \qquad
						-\frac
							{
								\big \langle
									B( \hat{t}_n, \hat{x}_n ) ,
									I_{\mathbb{H}_\var}^{-1} \big( 
										(\hat{C}_n \, B( \hat{t}_n, \hat{x}_n ))|_{H_\var}
									\big)
								\big \rangle_{HS(\mathbb{U}, \mathbb{H}_\var)}
							}
							{ \tilde{V}( \hat{t}_n, \hat{x}_n ) } 
						-\frac{v(t_n,x_n)}{\tilde{V}(t_n, x_n)} (r_n -\hat{r}_n)\\ \nonumber
				& \qquad
						+n
							\Big[
								\big \langle
									B( t_n, x_n ), 
									I_{\mathbb{H}_\var}^{-1} \big(
										(\mathfrak{C}_n \, B( t_n, x_n ) )|_{H_\var}
									\big)
								\big \rangle_{HS(\mathbb{U}, \mathbb{H}_\var)}
					\\ \nonumber & \qquad
								-	
								\big \langle
									B( \hat{t}_n, \hat{x}_n ),
									I_{\mathbb{H}_\var}^{-1} \big(
										(\hat{\mathfrak{C}}_n \, B( \hat{t}_n, \hat{x}_n ))|_{H_\var}
									\big)
								\big \rangle_{HS(\mathbb{U}, \mathbb{H}_\var)}
							\Big] 
						-\frac {\left \langle A(x_n), \tilde{p}_n \right \rangle_{H, H'}}
							{\tilde{V}(t_n, x_n)}
				\\ \nonumber & \qquad
						-\frac {\left \langle F(t_n, x_n), \tilde{p}_n \right \rangle_{H, H'}}
							{\tilde{V}(t_n, x_n)}
						-\bigg \langle
										B( t_n, x_n ),
										I_{\mathbb{H}_\var}^{-1} 
											\Big(
												\Big(
													\frac{\tilde{C}_n}{\tilde{V}( t_n, x_n )} \, B( t_n, x_n )
												\Big) \Big |_{H_\var}
											\Big )
								\bigg \rangle_{HS(\mathbb{U}, \mathbb{H}_{\var})} 
						\Bigg)
					\Bigg) \\ \nonumber
			&\leq 0
	\end{align}
	and this together with \eqref{eq: limit bound for big H1/2 diff}
	yields that
	\begin{align}
	\nonumber
			&\limsup_{n \to \infty} \Bigg(
					\frac{v(t_n,x_n)}{\tilde{V}(t_n, x_n)} r_n
						-\frac{v(\hat{t}_n,\hat{x}_n)}{\tilde{V}(\hat{t}_n,\hat{x}_n)} \hat{r}_n
					+\frac{\langle A(x_n), p_n \rangle_{H, H'}}
						{\tilde{V}(t_n, x_n)}
					-\frac{\langle A(\hat{x}_n), \hat{p}_n \rangle_{H, H'}}
						{\tilde{V}(\hat{t}_n, \hat{x}_n)} 
				\\ \nonumber & \qquad
					+\frac{\langle F(t_n, x_n), p_n \rangle_{H, H'}}
						{\tilde{V}(t_n, x_n)} 
					-\frac{\langle F(\hat{t}_n, \hat{x}_n), \hat{p}_n \rangle_{H, H'}}
						{\tilde{V}(\hat{t}_n, \hat{x}_n)}
			\\ \nonumber& \qquad
					+\frac
							{
								\big \langle
									B( t_n, x_n ) ,
									I_{\mathbb{H}_\var}^{-1} \big(
										(C_n \, B( t_n, x_n )) |_{H_\var}
									\big)
								\big \rangle_{HS(\mathbb{U}, \mathbb{H}_\var)} 					
							}
							{ \tilde{V}( t_n, x_n ) } 
			\\ \label{eq: limit bound}  & \qquad
					-\frac
						{
							\big \langle
								B( \hat{t}_n, \hat{x}_n ) ,
								I_{\mathbb{H}_\var}^{-1} \big(
									(\hat{C}_n \, B( \hat{t}_n, \hat{x}_n )) |_{H_\var}
								\big )
							\big \rangle_{HS(\mathbb{U}, \mathbb{H}_\var)}
						}
						{ \tilde{V}( \hat{t}_n, \hat{x}_n ) } 
					-\frac{v(t_n,x_n)}{\tilde{V}(t_n, x_n)} (r_n -\hat{r}_n)\\ \nonumber
			& \qquad
					+n
						\Big[
							\big \langle
								B( t_n, x_n ), 
								I_{\mathbb{H}_\var}^{-1} \big(
									(\mathfrak{C}_n \, B( t_n, x_n )) |_{H_\var}
								\big)
							\big \rangle_{HS(\mathbb{U}, \mathbb{H}_\var)}
					\\ \nonumber & \qquad
							-	
							\big \langle
								B( \hat{t}_n, \hat{x}_n ),
								I_{\mathbb{H}_\var}^{-1} \big(
									(\hat{\mathfrak{C}}_n \, B( \hat{t}_n, \hat{x}_n )) |_{H_\var}
								\big)
							\big \rangle_{HS(\mathbb{U}, \mathbb{H}_\var)}
						\Big] 
					-\frac {\left \langle A(x_n), \tilde{p}_n \right \rangle_{H, H'}}
						{\tilde{V}(t_n, x_n)}
				\\ \nonumber & \qquad
					-\frac {\left \langle F(t_n, x_n), \tilde{p}_n \right \rangle_{H, H'}}
						{\tilde{V}(t_n, x_n)} 
					-\bigg \langle
									B( t_n, x_n ),
									I_{\mathbb{H}_\var}^{-1} 
										\Big(
											\Big(
												\frac{\tilde{C}_n}{\tilde{V}( t_n, x_n )} \, B( t_n, x_n )
											\Big) \Big |_{H_\var}
										\Big)
							\bigg \rangle_{HS(\mathbb{U}, \mathbb{H}_{\var})} 
				\Bigg) \\ \nonumber
			&\leq 0.
	\end{align}
	Thus we obtain from
	\eqref{eq: different HS norm},
	\eqref{eq: delta v limit}, 
	\eqref{eq: delta F and B limit},
	\eqref{eq: delta F and B limit hat},
	and from \eqref{eq: limit bound} that
  \begin{align}  
				&\limsup_{ n \to \infty } \Big(
					\tfrac{ 1 }{ \tilde{V}( t_n , x_n ) }
					G_{\mathbb{H}, \mathbb{H}_\var, \delta_n, \tilde{V}h}^+\big( 
						(t_n, x_n), r_n, p_n, C_n+ n \mathfrak{C}_n \tilde{V}( t_n, x_n ) 
					\big) \\ \nonumber
				&\qquad -
					\tfrac{ 1 }{ \tilde{V}( \hat{t}_n , \hat{x}_n ) }
					G_{\mathbb{H}, \mathbb{H}_\var, \delta_n, \tilde{V}h}^-\big( 
						(\hat{t}_n, \hat{x}_n), \hat{r}_n, \hat{p}_n, 
						\hat{C}_n + n \hat{\mathfrak{C}}_n \tilde{V}( \hat{t}_n, \hat{x}_n ) 
					\big)
				\Big) \\ \nonumber
			={}&
				\limsup_{ n \to \infty } \bigg(
					\tfrac 
						{v( t_n, x_n ) (r_n + \delta_n \tilde{V}(t_n,x_n) h(t_n, x_n))}
						{\tilde{V}( t_n, x_n )}
					- 
					\tfrac
						{
							v( \hat{t}_n, \hat{x}_n ) 
							(\hat{r}_n - \delta_n \tilde{V}(\hat{t}_n, \hat{x}_n) h(\hat{t}_n, \hat{x}_n))
						}
						{\tilde{V}( \hat{t}_n, \hat{x}_n ) }
				\\ \nonumber
				& \qquad +
					\tfrac
						{
							\left\langle
								F( t_n, x_n) + A(x_n), 
								(p_n + \delta_n (D_{\mathbb{H}_\var} (h\tilde{V})) ( t_n, x_n))
							\right\rangle_{H, H'}
						}
						{\tilde{V}( t_n, x_n ) }
			\\ \nonumber & \qquad
					-
					\tfrac
						{
							\left\langle
								F( \hat{t}_n, \hat{x}_n ) + A( \hat{x}_n ) , 
									(\hat{p}_n - 
										\delta_n (D_{\mathbb{H}_\var} (h\tilde{V})) ( \hat{t}_n, \hat{x}_n ))
							\right\rangle_{H, H'}
						}
						{ \tilde{V}( \hat{t}_n , \hat{x}_n ) }
				\\  \nonumber
				& \qquad +
					\tfrac
						{
							\big \langle
								B( t_n, x_n ) ,
								I_{\mathbb{H}_\var}^{-1} \big(
									(C_n B( t_n, x_n ))|_{H_\var} 
									+ \delta_n (D_{\mathbb{H}_\var}^2 (h\tilde{V})) ( t_n, x_n ) 
										\, B( t_n, x_n )
								\big)
							\big \rangle_{HS(\mathbb{U}, \mathbb{H}_\var)} 
						}
						{ \tilde{V}( t_n, x_n ) } \\ \nonumber
				& \qquad -
					\tfrac
						{
							\big \langle
								B( \hat{t}_n, \hat{x}_n ) ,
								I_{\mathbb{H}_\var}^{-1} \big(
									(\hat{C}_n B( \hat{t}_n, \hat{x}_n )) |_{H_\var}
									- \delta_n (D_{\mathbb{H}_\var}^2 ( h\tilde{V})) ( \hat{t}_n, \hat{x}_n ) 
										\, B( \hat{t}_n, \hat{x}_n )
								\big)
							\big\rangle_{HS(\mathbb{U}, \mathbb{H}_\var)}
						}
						{ \tilde{V}( \hat{t}_n, \hat{x}_n ) }
				\\  \nonumber
				& \qquad + 
					n
					\Big[
						\big \langle
							B( t_n, x_n ), 
							I_{\mathbb{H}_\var}^{-1} \big(
								(\mathfrak{C}_n \, B( t_n, x_n ))|_{H_\var}
							\big)
						\big \rangle_{HS(\mathbb{U}, \mathbb{H}_\var)}
			\\ \nonumber & \qquad
						-
						\big \langle
							B( \hat{t}_n, \hat{x}_n ),
							I_{\mathbb{H}_\var}^{-1} \big(
								(\hat{\mathfrak{C}}_n \, B( \hat{t}_n, \hat{x}_n ))|_{H_\var}
							\big)
						\big \rangle_{HS(\mathbb{U}, \mathbb{H}_\var)}
					\Big]
				\bigg) 
				\\  \nonumber
			\leq{}&
				\limsup_{ n \to \infty } \Big(
					\tfrac{ v( t_n, x_n ) ( r_n - \hat{r}_n )}
						{\tilde{V}( t_n, x_n )}
					+\tfrac
						{
							\left\langle
								F( t_n, x_n), \tilde{p}_n
							\right\rangle_{H, H'}
						}
						{ \tilde{V}( t_n, x_n ) }
					+\tfrac {\left \langle A(x_n), \tilde{p}_n \right \rangle_{H, H'}}
						{\tilde{V}(t_n, x_n)} 
					+\big \langle
						\tfrac { B( t_n, x_n ) } { \tilde{V}( t_n, x_n )},
						I_\mathbb{H}^{-1} \tilde{C}_n \, B( t_n, x_n ) 
					\big \rangle_{HS(\mathbb{U}, \mathbb{H})}
				\Big) \\ \nonumber
			={}&  
				\limsup_{n \to \infty} \left (
					\tfrac{
						G( t_n, x_n, (r_n - \hat{r}_n), \tilde{p}_n, (\tilde{C}_n|_{H_\var})|_{H_\var} )
					}{
						\tilde{V}( t_n, x_n )
					} 
				\right ).
\end{align}
 This shows assumption~\eqref{eq:comparison.viscosity.solution.assumption}.
  Furthermore,
  by assumption,
  $
    u_1|_{ (0,T) \times O }
  $ 
  is a viscosity subsolution
  of~\eqref{eq:second-order.PDE}
	relative to
	$
		(h\tilde{V}, \R \times \mathbb{H}, \R \times \mathbb{H}_\var)
	$
  and
  $
    u_2|_{(0,T)\times O}
  $ 
  is a viscosity supersolution
  of~\eqref{eq:second-order.PDE}
	relative to
	$
		(h\tilde{V}, \R \times \mathbb{H}, \R \times \mathbb{H}_\var)
	$.
  Moreover,
  \eqref{eq:assumption_V} shows 
  for every $ r \in (0,\infty) $
  that the function 
	$r \tilde{V}$
  is a classical supersolution
  of~\eqref{eq:second-order.PDE}.
  In addition, observe that \eqref{eq:attains.maximum.difference2.comparison}
  follows from 
  $
    \lim_{ n \to \infty }
      \sup_{
        \substack{
          (t,x) \in [0,T] \times
          O_n^c
        }
      }
      \frac{
        | u_1(t, x) |
        +
        | u_2(t, x) |
      }{
        V(x)
      }
  $
  $
    = 0 
  $
	and \eqref{eq: uniformly bounded at 0 with V} follows from
	\eqref{eq: uniformly bounded at 0 with phi}.
  Consequently,
  Theorem~\ref{l:comparison.viscosity.solution}
		(with $\tilde{T} \leftarrow t$ and with
		$\tilde{R} \leftarrow \| x \|_H$)
	implies that
	for all $(t, x) \in W$ 
	and all $R \in (0,\infty)$ it holds that 
	 \begin{equation}
		\begin{split}
				0
			\geq{}&
				\lim_{r \downarrow 0} \lim_{\eps \downarrow 0}
					\sup \Bigg \{
						u_1(\tilde{t},\tilde{x}) - u_2(\hat{t} , \hat{x}) \colon 
						~(\tilde{t}, \tilde{x}), (\hat{t}, \hat{x}) \in W,
						~\tilde{t}, \hat{t} \leq t,
						~h(\tilde{t}, \tilde{x}) \vee h(\hat{t}, \hat{x}) \leq R,
				\\& \qquad \qquad \qquad
						~\|\tilde{x} \|_H \leq \| x \|_H, 
						~\|\tilde{x}- \hat{x}\|_H \leq r,
						~|\tilde{t}-\hat{t}| \leq \eps 
					\Bigg \} \\
			\geq{}&
				u_1(t,x) - u_2(t , x).
		\end{split}
  \end{equation}
	Hence, we obtain with the assumption that $u_1$ and $u_2$ are continuous 
	with respect to the $\| \cdot \|_{\R \times H}$-norm and 
	the fact that $W$ is dense in $[0,T] \times O$ 
	with respect to the $\| \cdot \|_{\R \times H}$-norm
	that $u_1 \leq u_2$.
  Repeating these arguments with $u_1$ and $u_2$ interchanged
  finally shows that $u_2\leq u_1$ so that $u_1=u_2$.
  This proves uniqueness and finishes the proof of
  Corollary~\ref{cor:uniqueness2}.
\end{proof}

The next remark gives sufficient conditions on $V$ to satisfy the 
assumptions in Corollary \ref{cor:uniqueness2}.
	\begin{remark}[Sufficient conditions for V]
	\label{rem: assumption on V}
		Let 
		$
				\mathbb{H}
			=
				( 
					H,
					\langle \cdot, \cdot \rangle_H, 
					\left\| \cdot \right\|_H
				)
		$
		be a real separable Hilbert space
		and
		$
				\mathbb{H}'
			=
				( 
					H',
					\left\| \cdot \right\|_{H'}
				)
		$
		its dual space,
		let $ (e_i)_{i \in \N} \subseteq H $ be an orthonormal basis of $ \mathbb{H} $,
		let $\lambda \colon \N \to (- \infty, 0)$ be a function with 
		$\lim_{n \to \infty} \lambda_n = -\infty$,
		let
		$
			A \colon D(A) \subseteq H \rightarrow H
		$ 
		be the linear operator such that
		\begin{equation}
			D(A) 
			= 
				\biggl\{ 
					v \in H 
					\colon
					\sum^{\infty}_{ i = 1 } 
			\left| 
				\lambda_i 
				\langle e_i, v \rangle_H 
			\right|^2
			< \infty
				\biggr\}
		\end{equation}
		and such that for all $ v \in D(A) $ it holds that
		\begin{equation}
				Av
			 =
				\sum^{\infty}_{ i =1 } 
				\lambda_i \langle e_i, v \rangle_H \, e_i,
		\end{equation}
		let 
		$ 
				\mathbb{H}_r
			=
				( H_r , \left< \cdot , \cdot \right>_{ H_r }, \left\| \cdot \right\|_{ H_r } ) 
		$,
		$ r \in \R $,
		be a family of interpolation spaces associated with
		$ - A $
		(see, e.g., Definition~3.6.30 in Jentzen \cite{Jentzen2015}),
		and let 
		$ 
				\mathbb{H}'_r
			=
				( H'_r, \left\| \cdot \right\|_{ H'_r } ) 
		$
		be the corresponding dual spaces.
		By abuse of notation we will also denote by
		$A$ and  by $\| \cdot \|_{H_r}$, $r \in \R$
		the extended operators 
		$
			A \colon \bigcup_{i=1}^{\infty} H_{-i} \to \bigcup_{i=1}^{\infty} H_{-i}
		$ 
		and 
		$
			\| \cdot \|_{H_r} \colon \bigcup_{i=1}^{\infty} H_{-i} \to [0, \infty]
		$,
		$r \in \R$
		satisfying for all
		$r \in \R$, $x \in H_{r}$, $y \in H_{r-1}$ 
		and all $z \in \bigcup_{i=1}^{\infty} H_{-i}$
		that
		\begin{equation}
				(A(x) = y) 
			\Leftrightarrow 
				(
						\lim_{\eps \downarrow 0} \sup \{
							\|A(\xi) - y\|_{H_{r-1}} \colon \xi \in H_1, ~\|x- \xi\|_{H_{r}} \leq \eps 
						\}
					=
						0
				)
		\end{equation}
		and that
		\begin{equation}
					\|z\|_{H_r}
				=
					\begin{cases} 
						\|z\|_{H_{r}} 
							& \textrm{ if } z \in H_{r} \\
						\infty 
							& \textrm{ if }z \notin H_{r}
					\end{cases}.
		\end{equation}
		Let
		$\vartheta \in [\nicefrac {1}{2} , \infty)$,
		$\kappa \in (-\nicefrac 12, 0]$,
		$l \in [-\kappa, \infty)$,
		let 
		$
			\| \cdot \|_{\boldsymbol{L}(\mathbb{H}_{\kappa},\mathbb{H}'_{\kappa})} 
				\colon L(\mathbb{H}, \mathbb{H}') \to [0, \infty],
		$
		and
		$
			\| \cdot \|_{\boldsymbol{L}(\mathbb{H}_{\kappa},\mathbb{H}_{-\kappa})} 
				\colon L(\mathbb{H}, \mathbb{H}) \to [0, \infty]
		$
		be the extended norms satisfying for all $C \in L(\mathbb{H}, \mathbb{H}')$
		that 
		\begin{equation}
				\| I^{-1}_{\mathbb{H}} C \|_{\boldsymbol{L}(\mathbb{H}_{\kappa},\mathbb{H}_{-\kappa})}
			=
				\| C \|_{\boldsymbol{L}(\mathbb{H}_{\kappa},\mathbb{H}'_{\kappa})}
			=
				\sup_{x,y \in H \backslash \{0 \}} 
					\frac{|\langle Cx, y \rangle_{H, H'}|}{\|x\|_{H_{\kappa}} \|y\|_{H_{\kappa}}},
		\end{equation}
		let $f \in \C^1([0,\infty), (1,\infty))$
		satisfy that $f|_{(0,\infty)} \in \C^2((0,\infty), (1,\infty))$
		and that $\lim_{t \downarrow 0} (t f''(t))= 0$,
		and let
		$V \colon H \to \R$ satisfy for all $ x \in H$ that
		$V(x) = f(\nicefrac 12 \|(-A)^{-l} x\|^2_H)$.
		Then it holds that
		$V \in \C^2_{\mathbb{H}}(H,(1, \infty))$,
		that $D_{\mathbb{H}}^2 \, V \colon H \to L(\mathbb{H}, \mathbb{H}')$ 
		is uniformly continuous with respect to the
		$\| \cdot \|_{\mathbb{H}}$ and the
		$\| \cdot \|_{\boldsymbol{L}(\mathbb{H}_{\kappa}, \mathbb{H}_{\kappa}')}$-norm
		on all $\mathbb{H}$-bounded subsets of $H$,
		and that there exist functions
		$K \colon [0, \infty) \to (0, \infty)$ and
		$\omega \colon [0, \infty) \to [0, \infty)$
		such that $K$ is increasing,
		such that
			$\lim_{x \downarrow 0} \omega (x) = 0$,
		such that
			$\sup_{x \in [0,\infty)} \omega(x) < \infty$,
		and such that for all 
		$x$, $\hat{x} \in H_{2 \var}$ 
		it holds that
		\begin{equation}
		\label{eq: V continuous with respect to H1 rem}
				\bigg \|
					\frac{I_\mathbb{H}^{-1} (D_\mathbb{H} V)(x)}{V(x)} 
					- \frac{I_\mathbb{H}^{-1} (D_\mathbb{H} V)(\hat{x})}{V(\hat{x})}
				\bigg\|_{H_{1-\var}}
			\leq
				K(\|x\|^2_{H_\var} \vee \|\hat{x}\|^2_{H_\var}) 
				\, (\|x- \hat{x}\|^2_{H_{1/2}} + \omega(\|x- \hat{x}\|^2_{H_{1/2}})), 
		\end{equation}
		\begin{align}
		\label{eq: V continuous with respect to H kappa rem}
				&\left \|
					\frac
						{I_{\mathbb{H}}^{-1} ( D_\mathbb{H} V)( x )}
						{V (x )}
					-\frac
						{I_{\mathbb{H}}^{-1} (D_\mathbb{H} V )( \hat{x} )}
						{V( \hat{x} )}
				\right \|_{H_{-\kappa}}
			\leq 
				K(\|x\|_{H} \vee \| \hat{x} \|_H) \cdot \|x - \hat{x}\|_{H_{1/2}}, \\
		\label{eq: V bounded with respect to H kappa rem}
				&\left \|
					I_\mathbb{H}^{-1}(D_\mathbb{H} V)(x) 
				\right \|_{H_{-\kappa}}
			\leq
				K(\|x\|_{H}),
		\end{align}
		and	that 
		\begin{equation}
		\label{eq: V bounded with respect to H 1/2 rem}
				\left \|
					I_\mathbb{H}^{-1}(D_\mathbb{H} V)(x) 
				\right \|_{H_{1/2}}
			\leq
				K(\|x\|^2_{H_\var}).
		\end{equation}
	\end{remark}
	
	\begin{proof}
		First note that we get
		from the assumption
		$l \geq -\kappa$ 
		that
		\begin{equation}
		\label{eq: 2nd derivative of the l-norm}
			\begin{split}
					&\| 
						D^2_\mathbb{H} (\nicefrac 12 \| \cdot \|^2_{H_{-l}}) 
					\|_{\boldsymbol{L}(\mathbb{H}_{\kappa},\mathbb{H}_{\kappa}')}
				=
					\|
						(-A)^{-2l}
					\|_{\boldsymbol{L}(\mathbb{H}_{\kappa},\mathbb{H}_{\kappa}')}
				=
					\sup_{x, y \in H \backslash \{0\} }
						\tfrac
							{\langle (-A)^{-2l} x, y \rangle_{H',H}}
							{\|x\|_{H_{\kappa}} \|y\|_{H_{\kappa}}} \\
				=
					&\sup_{x, y \in H_{-\kappa} \backslash \{0\}}
						\tfrac
							{\langle (-A)^{-2l} (-A)^{-\kappa}x, (-A)^{-\kappa} y \rangle_{H',H}}
							{\|(-A)^{-\kappa} x \|_{H_{\kappa}} \|(-A)^{-\kappa} y \|_{H_{\kappa}}}
				=
					\sup_{x, y \in H \backslash \{0\}}
						\tfrac
							{\langle  (-A)^{-2l-2\kappa} x, y \rangle_{H',H}}
							{\| x \|_{H} \|y \|_H} \\
				\leq{}
					&\sup_{n \in \N} (|\lambda_n|^{-2l-2\kappa})
				<
					\infty.
			\end{split}
		\end{equation}	
		Thus we obtain
		that
		\begin{equation}
		\label{eq: 2nd derivative continuous wrt kappa}
				D^2_\mathbb{H} (\nicefrac 12 \| \cdot \|^2_{H_{-l}})
			\in \C_{\mathbb{H},\boldsymbol{L}(\mathbb{H}_{\kappa},\mathbb{H}_{\kappa}')}
						(H, L(\mathbb{H},\mathbb{H}')).
		\end{equation}
		Combining this with $f|_{(0,\infty)} \in \C^2((0,\infty), (1,\infty))$
		shows then that
		\begin{equation}
		\label{eq: V diffbar outside 0}
				D_\mathbb{H}^2 (V|_{H\backslash \{0\}}) 
			\in \C_{\mathbb{H},\boldsymbol{L}(\mathbb{H}_{\kappa},\mathbb{H}_{\kappa}')}
						(H\backslash \{0\}, L(\mathbb{H},\mathbb{H}')). 
		\end{equation}
		Next note that we have for all $x \in H \backslash \{ 0 \}$ 
		that 
		\begin{equation}
			\begin{split}
				&(D^2_\mathbb{H} V)(x) \\
			={} 
				&f''(\nicefrac 12 \|(-A)^{-l} x\|^2_H) \, 
					I_\mathbb{H} ((-A)^{-2l} x) \otimes I_\mathbb{H} ((-A)^{-2l} x)
				+ f'(\nicefrac 12 \|(-A)^{-l} x\|^2_H) \, I_\mathbb{H} \, (-A)^{-2l}. 
			\end{split}
		\end{equation}
		Moreover, it holds that 
		\begin{align}
			\nonumber
					&\lim_{H \backslash \{0 \} \ni x \to 0}
						\frac
							{|V(x)-V(0) - \nicefrac 12 \, f'(0) \cdot \langle (-A)^{-2l} x, x \rangle_H|}
							{\|x\|^2_H} \\
				\leq{}&
					\lim_{H \backslash \{0 \} \ni x \to 0}
						\frac
							{
								|
									f(\nicefrac 12 \|(-A)^{-l} x\|^2_H)-f(0) 
									- \nicefrac 12 \, f'(0) \cdot \| x \|^2_{H_{-l}}
								|
							}
							{(\inf_{n \in \N} |\lambda_n|)^{l} \, \|x\|^2_{H_{-l}}} \\  \nonumber
				={}
					&\lim_{t \downarrow 0}
						\frac
							{
								|
									f(t)-f(0) 
									- f'(0) \cdot t
								|
							}
							{2 (\inf_{n \in \N} |\lambda_n|)^{l} \, t} 
				=
					0
		\end{align}
		and this yields that 
		\begin{equation}
		\label{eq: V diffbar in 0}
			(D^2_\mathbb{H} V)(0) = f'(0) I_\mathbb{H} \, (-A)^{-2l}.
		\end{equation}
		In addition, it follows from the assumptions
		$\lim_{t \downarrow 0} t f''(t)=0$
		and $l + \kappa \geq 0$
		that
		\begin{align}
					&\limsup_{\eps \downarrow 0} \Big \{
						\big\|
							f''(\nicefrac 12 \|(-A)^{-l} x\|^2_H) \, 
							I_\mathbb{H} ((-A)^{-2l} x) \otimes I_\mathbb{H} ((-A)^{-2l} x)
						\big \|_{\boldsymbol{L}(\mathbb{H}_{\kappa}, \mathbb{H}_{\kappa}')}
						\colon 
				\\ \nonumber & \qquad	\qquad	
						x \in H \backslash \{0\}, ~\|x\|_{H} \leq \eps 
					\Big \} \\ \nonumber
				={}&
					\limsup_{\eps \downarrow 0} \bigg \{
						\sup_{y,z \in H \backslash \{0\}}
							f''(\nicefrac 12 \| x\|^2_{H_{-l}}) \, 
							\frac
								{
									\langle (-A)^{-2l} x, y \rangle_H 
									\cdot \langle (-A)^{-2l} x, z \rangle_H
								}
								{
									\|y\|_{H_{\kappa}} \|z\|_{H_{\kappa}}
								}
						\colon x \in H \backslash \{0\}, 
				\\ \nonumber & \qquad	\qquad
						\|x\|_{H} \leq \eps 
					\bigg \} \\ \nonumber
				={}&
					\limsup_{\eps \downarrow 0} \bigg \{
						f''(\nicefrac 12 \|x\|^2_{H_{-l}}) \,
						\sup_{y,z \in H_{-\kappa} \backslash \{0\}}	 
							\frac
								{
									\langle (-A)^{-2l} x, (-A)^{-\kappa} y \rangle_H 
									\cdot \langle (-A)^{-2l} x, (-A)^{-\kappa} z \rangle_H
								}
								{
									\|y\|_{H} \|z\|_{H}
								}
						\colon 
				\\ \nonumber & \qquad	\qquad
						x \in H \backslash \{0\}, ~\|x\|_{H} \leq \eps 
					\bigg \} \\ \nonumber
				={}&
					\limsup_{\eps \downarrow 0} \bigg \{
							f''(\nicefrac 12 \|x\|^2_{H_{-l}}) \, 
							\| (-A)^{-l} (-A)^{-l-\kappa} x\|^2_H 
						\colon 
						x \in H \backslash \{0\}, ~\|x\|_{H} \leq \eps 
					\bigg \} \\ \nonumber
				\leq{}&
					\limsup_{\eps \downarrow 0} \left \{
							f''(\nicefrac 12 \|x\|^2_{H_{-l}}) \, 
							\sup_{n \in \N} (|\lambda_n|^{-l-\kappa} )
							\| (-A)^{-l} x\|^2_H 
						\colon x \in H \backslash \{0\}, ~\|x\|_{H} \leq \eps 
					\right \} \\  \nonumber
				\leq{}&
					\limsup_{\eps \downarrow 0} \left \{
							\sup_{n \in \N} (|\lambda_n|^{-l-\kappa} ) \,
							|f''(\nicefrac 12 \|x\|^2_{H_{-l}})| \cdot
							\|  x\|^2_{H_{-l}}
						\colon x \in H \backslash \{0\}, ~\|x\|_{H_{-l}} \leq \eps 
					\right \} 
				= 0.
		\end{align}
		Combining this with \eqref{eq: V diffbar outside 0} 
		and with \eqref{eq: V diffbar in 0} ensures that 
		$
				D_\mathbb{H}^2 \, V 
			\in \C_{\mathbb{H},\boldsymbol{L}(\mathbb{H}_{\kappa},\mathbb{H}_{\kappa}')}
						(H, L(\mathbb{H},\mathbb{H}')).
		$
		In particular the assumption $\kappa \leq 0$ verifies that
		$
				D_\mathbb{H}^2 \, V
			\in \C_{\mathbb{H},L(\mathbb{H},\mathbb{H}')}
						(H, L(\mathbb{H},\mathbb{H}'))
		$
		and therefore
		$ V \in \C^2_{\mathbb{H}} (H, \R).$
		In addition, note that the fact that 
		$D^2_\mathbb{H} (\nicefrac 12 \| \cdot \|^2_{H_{-l}})$ is constant shows
		that $D^2_\mathbb{H} (\nicefrac 12 \| \cdot \|^2_{H_{-l}})$ is also 
		uniform continuous with respect to the $\| \cdot \|_H$-norm and the
		$\|\cdot\|_{\boldsymbol{L}(\mathbb{H}_\kappa, \mathbb{H}'_\kappa)}$-norm
		and the
		fact that $(f|_{(0, \infty)})''$ is continuous implies for all $r, R \in (0, \infty)$ 
		that $(f|_{(0, \infty)})''$ is uniform continuous on 
		the compact set $[r, R]$.
		Combining this with the fact that the function
		$H \ni x \to \nicefrac 12 \| x \|^2_{H_{-l}} \in \R$
		is bounded on $\mathbb{H}$-bounded subsets of $H$
		yields that $D^2_\mathbb{H} V$ is for all $r \in (0, \infty)$
		uniform continuous with respect to the $\| \cdot \|_H$-norm and the
		$\|\cdot\|_{\boldsymbol{L}(\mathbb{H}_\kappa, \mathbb{H}'_\kappa)}$-norm
		on $\mathbb{H}$-bounded subsets of 
		$ \{ x \in H \colon r \leq \|x\|_H\}$.
		This means that for all $\eps, r, R \in (0, \infty)$ there exist a 
		$\delta_{\eps,r,R}$ such that for all 
		$x,y \in \{ z \in H \colon r \leq \|z\|_H \leq R\}$ 
		with $\|x-y\|_H \leq \delta_{\eps, r, R}$ it holds that
		\begin{equation}
		\label{eq: D^2 V uniform continuous outside 0}
				\|
					(D^2_\mathbb{H} V) (x) - (D^2_\mathbb{H} V) (y)
				\|_{\boldsymbol{L}(\mathbb{H}_{\kappa},\mathbb{H}_{\kappa}')}
			<	
				\eps.
		\end{equation}
		On the other hand it follows from
		\eqref{eq: 2nd derivative continuous wrt kappa} that for all
		$\eps \in (0, \infty)$ there exist a $\delta_\eps$ such that
		for all $x \in \{ z \in H \colon \|z\|_H \leq \delta_\eps\}$ it holds that
		\begin{equation}
				\|
					(D^2_\mathbb{H} V) (x) - (D^2_\mathbb{H} V) (0)
				\|_{\boldsymbol{L}(\mathbb{H}_{\kappa},\mathbb{H}_{\kappa}')}
			<	
				\tfrac \eps 2
		\end{equation}
		and this ensures that for all 
		$x, y \in \{ z \in H \colon \|z\|_H \leq \delta_\eps\}$ it holds that
		\begin{equation}
		\label{eq: D^2 V uniform continuous in 0}
			\begin{split}
				&\|
					(D^2_\mathbb{H} V) (x) - (D^2_\mathbb{H} V) (y)
				\|_{\boldsymbol{L}(\mathbb{H}_{\kappa},\mathbb{H}_{\kappa}')} \\
			\leq{}	
				&\|
					(D^2_\mathbb{H} V) (x) - (D^2_\mathbb{H} V) (0)
				\|_{\boldsymbol{L}(\mathbb{H}_{\kappa},\mathbb{H}_{\kappa}')}
				+
				\|
					(D^2_\mathbb{H} V) (y) - (D^2_\mathbb{H} V) (0)
				\|_{\boldsymbol{L}(\mathbb{H}_{\kappa},\mathbb{H}_{\kappa}')}
			<	
				\eps.
			\end{split}
		\end{equation}
		Finally note that we have for all
		$\delta \in (0, \infty)$ and all $x, y \in H$
		with $\|x-y\|_H \leq \tfrac {\delta}{2}$ that
		\begin{align}
		\label{eq: x,y both small or big}
			\begin{split}
					&(\|x\|_H \leq \tfrac \delta 2) \vee
					(\|y\|_H \leq \tfrac \delta 2) \vee
					(\min\{ \|x\|_H, \|y\|_H \} > \tfrac \delta 2) \\
				\Rightarrow
					&(
						(\|x\|_H \leq \tfrac \delta 2) 
						\wedge (\|y\|_H \leq \|x\|_H + \|x-y\|_H \leq \delta)
					) 
				\\ &
					\vee
					(
						(\|y\|_H \leq \tfrac \delta 2) 
						\wedge (\|x\|_H \leq \|y\|_H + \|x-y\|_H \leq \delta)
					) 
					\vee 
					(\min\{ \|x\|_H, \|y\|_H \} > \tfrac \delta 2) \\
				\Rightarrow
					&(\max \{\|x\|_H, \|y\|_H \} \leq \delta)
					\vee
					(\min\{ \|x\|_H, \|y\|_H \} \geq \tfrac \delta 2).
			\end{split}
		\end{align}
		Combining now \eqref{eq: D^2 V uniform continuous outside 0},
		\eqref{eq: D^2 V uniform continuous in 0},
		\eqref{eq: x,y both small or big} implies that for all
		$(\eps, R) \in (0, \infty)^2$ and all
		$(x,y) \in H^2$ with
		$
				\|x-y\|_H 
			\leq 
				 \tfrac {\delta_\eps}{2} \wedge \delta_{\eps, \delta_{\eps}/ 2, R}
		$
		and with $\|x\|_H \vee \|y\|_H \leq R$
		it holds that
		\begin{equation}
				\|
					(D^2_\mathbb{H} V) (x) - (D^2_\mathbb{H} V) (y)
				\|_{\boldsymbol{L}(\mathbb{H}_{\kappa},\mathbb{H}_{\kappa}')}
			<	
				\eps
		\end{equation}
		and this verifies that $(D^2_\mathbb{H} V)$ is uniform continuous with
		respect to the $\| \cdot \|_H$ and the
		$\| \cdot \|_{\boldsymbol{L}(\mathbb{H}_\kappa, \mathbb{H}_\kappa')}$-norm
		on $\mathbb{H}$-bounded subsets of $H$.
		Next note that we have for all $x \in H$ that
		\begin{equation}
		\label{eq: derivative of V}
				(D_\mathbb{H} V)(x)
			=
				f'(\nicefrac 12 \|(-A)^{-l} x\|^2_H) \, I_{\mathbb{H}} \, ((-A)^{-2l} x).
		\end{equation}
		Moreover, the fact that $f'$ and $f$ are continuous 
		and the fact that for all $R \in (0, \infty)$
		it holds that $[0, R]$ is compact implies that
		there exist an increasing function $\tilde{K} \colon [0, \infty) \to (0, \infty)$
		and a function $\tilde{\omega} \colon [0, \infty) \to [0, \infty)$
		with $\lim_{t \downarrow 0} \tilde{\omega}(t) =0$
		and with $\sup_{t \in (0,\infty)} \tilde{\omega}(t) <\infty$
		such that
		for all $t, \hat{t} \in [0, \infty)$ it holds that
		\begin{equation}
		\label{eq: f and f' inequalities}
			\begin{split}
					&|f(\tfrac t2 ) - f(\tfrac {\hat{t}}{2})| 
				\leq 
					\tilde{K}( t \vee \hat{t}) 
						\cdot |t - \hat{t}|, \qquad 
					f(\tfrac t2) \vee |f'(\tfrac t2)| 
				\leq
					\tilde{K}( t)
			\end{split}
		\end{equation}
		and such that
		for all $t, \hat{t} \in [0, 2]$ it holds that
		\begin{equation}
		\label{eq: f' uniform cont}
			\begin{split}
					|f'(\tfrac t2) - f'(\tfrac {\hat{t}} {2})| 
				\leq 
					\tilde{\omega}(|t - \hat{t}|).
			\end{split}
		\end{equation}
		This implies that for all $x,\hat{x} \in H$ with
		$\|x\|^2_{H_{-l}} \vee\| \hat{x}\|^2_{H_{-l}} \leq 2$ it holds that
		\begin{equation}
		\label{eq: f' uniform cont for small x}
					|
						f'(\nicefrac 12 \|(-A)^{-l} x\|^2_H) 
						-f'(\nicefrac 12 \|(-A)^{-l} \hat{x}\|^2_H)
					|
				\leq
					\tilde{\omega}(
						\|x\|_{H_{-l}}^2
						-\|\hat{x}\|_{H_{-l}}^2
					)
				\leq 
					\sup_{t\in [0,4]} \omega (t \cdot \|x-\hat{x}\|_{H_{-l}}).
		\end{equation}
		In addition, \eqref{eq: f and f' inequalities}
		shows that for all $x, \hat{x} \in H$
		with $\|x\|^2_{H_{-l}} \vee \| \hat{x}\|^2_{H_{-l}} \geq 2$
		and with $ \|x\|^2_{H_{-l}} \wedge \| \hat{x}\|^2_{H_{-l}} \leq 1$ 
		it holds that
		\begin{equation}
		\label{eq: f' uniform cont for small x and big x}
			\begin{split}
					&|
						f'(\nicefrac 12 \|(-A)^{-l} x\|^2_{H_{-l}}) 
						-f'(\nicefrac 12 \|(-A)^{-l} \hat{x}\|^2_{H_{-l}})
					|
				\leq
					(\tilde{K}(\|x\|^2_{H_{-l}}) +\tilde{K}(\|\hat{x}\|^2_{H_{-l}})) \\
				\leq{} 
					&(\tilde{K}(\|x\|^2_{H_{-l}}) +\tilde{K}(\|\hat{x}\|^2_{H_{-l}}))
					\cdot \|x-\hat{x}\|_{H_{-l}}.
			\end{split}
		\end{equation}
		Furthermore,we get for all $x, \hat{x} \in H$
		with with $ \|x\|^2_{H_{-l}} \wedge \| \hat{x}\|^2_{H_{-l}} \geq 1$
		that
		\begin{equation}
		\label{eq: f' uniform cont for big x}
					|
						f'(\nicefrac 12 \|(-A)^{-l} x\|^2_H) 
						-f'(\nicefrac 12 \|(-A)^{-l} \hat{x}\|^2_H)
					|
				\leq
					(\sup_{t \in [1,\|x\|_{H_{-l}}^2 \vee \|\hat{x}\|_{H_{-l}}^2]} |f''(\tfrac t2)|)
					\cdot \|x-\hat{x}\|_{H_{-l}}.
		\end{equation}
		Combining \eqref{eq: f' uniform cont for small x},
		\eqref{eq: f' uniform cont for small x and big x},
		and \eqref{eq: f' uniform cont for big x}
		shows that for all $x,\hat{x} \in H$ it holds that
		\begin{align}
		\nonumber
					&\|
						f'(\nicefrac 12 \|(-A)^{-l} x\|^2_H) \, (-A)^{-2l} x
						-f'(\nicefrac 12 \|(-A)^{-l} \hat{x}\|^2_H) \, (-A)^{-2l} \hat{x}
					\|_{H_{1-\vartheta}}  \\ \nonumber
				\leq{}&
					\|
						(
							f'(\nicefrac 12 \|(-A)^{-l} x\|^2_H)
							-f'(\nicefrac 12 \|(-A)^{-l} \hat{x}\|^2_H)
						)
						\, (-A)^{-2l} x
					\|_{H_{1-\vartheta}}
			\\ \nonumber &
					+\|
						f'(\nicefrac 12 \|(-A)^{-l} \hat{x}\|^2_H) \, (-A)^{-2l} (x-\hat{x})
					\|_{H_{1-\vartheta}} \\ \nonumber
				\leq{}&
					|
							f'(\nicefrac 12 \|(-A)^{-l} x\|^2_H)
							-f'(\nicefrac 12 \|(-A)^{-l} \hat{x}\|^2_H)
						|
						\cdot \| x \|_{H_{1-\vartheta-2l}}
			\\ 
	\label{eq: f' x uniform cont}
					&
						+|f'(\nicefrac 12 \|(-A)^{-l} \hat{x}\|^2_H)| 
						\cdot \|x-\hat{x}\|_{H_{1-\vartheta-2l}} \\ \nonumber
				\leq{}&
					\big(
						(
							\sup_{t \in [1,\|x\|_{H_{-l}}^2 \vee \|\hat{x}\|_{H_{-l}}^2]} 
								|f''(\tfrac t2)|
						)
						+\tilde{K}(\|x\|^2_{H_{-l}}) + \tilde{K}(\|\hat{x}\|^2_{H_{-l}})
						+1
					\big)
					\cdot \| x \|_{H_{1-\vartheta-2l}} \\ \nonumber
					\cdot 
					&\big(
						\|x-\hat{x}\|^2_{H_{-l}}
						+(
							\1_{\|x-\hat{x}\|^2_{H_{-l}} \leq 1} \|x-\hat{x}\|_{H_{-l}}
							+\sup_{t\in [0,4]} \tilde{\omega}(t \cdot \|x-\hat{x}\|_{H_{-l}})
						)
					\big)
		\\  \nonumber &
					+|f'(\nicefrac 12 \|(-A)^{-l} \hat{x}\|^2_H)| 
						\cdot \|x-\hat{x}\|_{H_{1-\vartheta-2l}}.
		\end{align}
		Moreover, we have for all $(x,\hat{x}) \in H^2 \backslash \{(0,0)\}$
		that
		\begin{align}
		\label{eq: f' Lipschitz cont}
					&\|
						f'(\nicefrac 12 \|(-A)^{-l} x\|^2_H) \, (-A)^{-2l} x
						-f'(\nicefrac 12 \|(-A)^{-l} \hat{x}\|^2_H) \, (-A)^{-2l} \hat{x}
					\|_{H_{-\kappa}}  
				\\ \nonumber \leq{}&
					\int_{0}^{1} 
						\Big \|
							f''(\nicefrac 12 \|(-A)^{-l} (tx +(1-t) \hat{x})\|^2_H)
								\cdot \big \langle
									(-A)^{-l} (tx +(1-t) \hat{x}),
									(-A)^{-l} (x - \hat{x})
								\big \rangle_H
				\\ \nonumber & \qquad
									\cdot (-A)^{-2l} (tx +(1-t) \hat{x})
							+f'(\nicefrac 12 \|(-A)^{-l} (tx +(1-t) \hat{x})\|^2_H)
								\cdot (-A)^{-2l} (x-\hat{x})
						\Big \|_{H_{-\kappa}}
					\ud t 
				\\ \nonumber \leq{}&
					\int_{0}^{1} 
						\Big \|
							f''(\nicefrac 12 \|(-A)^{-l} (tx +(1-t) \hat{x})\|^2_H)
									\cdot \|tx +(1-t) \hat{x} \|_{H_{-l}} \cdot
									\|x - \hat{x} \|_{H_{-l}}
				\\ \nonumber & \qquad
									\cdot (-A)^{-2l} (tx +(1-t) \hat{x})
						\Big \|_{H_{-\kappa}}
				\\ \nonumber & \qquad
						+\Big \|
							f'(\nicefrac 12 \|(-A)^{-l} (tx +(1-t) \hat{x})\|^2_H)
								\cdot (-A)^{-2l} (x-\hat{x})
						\Big \|_{H_{-\kappa}}
					\ud t 
				\\ \nonumber \leq{}&
					\sup_{t\in (0,1)} \Big(
							f''(\nicefrac 12 \|(-A)^{-l} (tx +(1-t) \hat{x})\|^2_H)
								\cdot \|tx +(1-t) \hat{x} \|_{H_{-l}} \cdot
								\|tx +(1-t) \hat{x} \|_{H_{-2l-\kappa}}
				\\ \nonumber & \qquad
								\cdot \|x - \hat{x} \|_{H_{-l}}
						+f'(\nicefrac 12 \|(-A)^{-l} (tx +(1-t) \hat{x})\|^2_H)
								\cdot \|x - \hat{x} \|_{H_{-2l-\kappa}}
					\Big) 
				\\ \nonumber \leq{}&
					\sup_{t\in (0,\|x\|^2_{H_{-l}} \vee \|\hat{x}\|^2_{H_{-l}} ]} \big(
							|f''(\tfrac t2)|
							\cdot t
					\big)
								\inf_{n \in \N} ( |\lambda_n|^{-l-\kappa})
								\cdot \|x - \hat{x} \|_{H_{-l}}
				\\ \nonumber & \qquad
						+\sup_{t\in (0,\|x\|^2_{H_{-l}} \vee \|\hat{x}\|^2_{H_{-l}} ]} \big(
							|f'(\tfrac t2)|
						\big)
							\cdot \|x - \hat{x} \|_{H_{-2l-\kappa}}. 
		\end{align}
		In addition,
		\eqref{eq: derivative of V},
		the assumption $f > 1$,
		and \eqref{eq: f and f' inequalities}
		imply for all $r \in (-\infty,\infty)$
		and all $x, \hat{x} \in H$
		that
		\begin{align}
		\label{eq: DV /V Difference}
					&\left \|
						\frac{I_\mathbb{H}^{-1} (D_\mathbb{H} V)(x)}{V(x)} 
						-\frac{I_\mathbb{H}^{-1} (D_\mathbb{H} V)(\hat{x})}{V(\hat{x})} 
					\right \|_{H_r} \\ \nonumber
				={}&
					\left \|
						\frac
							{
								I_\mathbb{H}^{-1} f'(\nicefrac 12 \|(-A)^{-l} x\|^2_H) \, 
								I_{\mathbb{H}} \, ((-A)^{-2l} x)
							}
							{f(\nicefrac 12 \|(-A)^{-l} x\|^2_H)} 
						-\frac
							{
								I_\mathbb{H}^{-1} f'(\nicefrac 12 \|(-A)^{-l} \hat{x}\|^2_H) \, 
								I_{\mathbb{H}} \, ((-A)^{-2l} \hat{x})
							}
							{f(\nicefrac 12 \|(-A)^{-l} \hat{x}\|^2_H)} 
					\right \|_{H_r} \\ \nonumber
				\leq{}&
					\frac
						{
							\|
								f'(\nicefrac 12 \|(-A)^{-l} x\|^2_H) \, (-A)^{-2l} x
								-f'(\nicefrac 12 \|(-A)^{-l} \hat{x}\|^2_H) \, 
									(-A)^{-2l} \hat{x}
							\|_{H_r}
						}
						{f(\nicefrac 12 \|(-A)^{-l} x\|^2_H)} \\ \nonumber
					&+
						\|(-A)^{-2l} \hat{x}\|_{H_{r}}
						\Big|
							\frac
								{1}
								{f(\nicefrac 12 \|(-A)^{-l} x\|^2_H)} 
							-\frac
								{1}
								{f(\nicefrac 12 \|(-A)^{-l} \hat{x}\|^2_H)} 
						\Big|\\ \nonumber
				\leq{}&
						\|
							f'(\nicefrac 12 \|(-A)^{-l} x\|^2_H) \, (-A)^{-2l} x
							-f'(\nicefrac 12 \|(-A)^{-l} \hat{x}\|^2_H) \, 
								(-A)^{-2l} \hat{x}
						\|_{H_r} \\ \nonumber
					&+
						\|\hat{x}\|_{H_{r-2l}} \,
						|
							f(\nicefrac 12 \|(-A)^{-l} \hat{x}\|^2_H)
							-f(\nicefrac 12 \|(-A)^{-l} x\|^2_H)
						|\\ \nonumber
				\leq{}&
						\|
							f'(\nicefrac 12 \|(-A)^{-l} x\|^2_H) \, (-A)^{-2l} x
							-f'(\nicefrac 12 \|(-A)^{-l} \hat{x}\|^2_H) \, 
								(-A)^{-2l} \hat{x}
						\|_{H_r} \\ \nonumber
					&+
						\|\hat{x}\|_{H_{r-2l}} \,
						\tilde{K}(\|x\|^2_{H_{-l}} \vee \|x\|_{H_{-l}}^2)
						\cdot \big |
							\|(-A)^{-l} x\|^2_H
							-\|(-A)^{-l} \hat{x}\|^2_H
						\big|\\ \nonumber
				\leq{}&
						\|
							f'(\nicefrac 12 \|(-A)^{-l} x\|^2_H) \, (-A)^{-2l} x
							-f'(\nicefrac 12 \|(-A)^{-l} \hat{x}\|^2_H) \, 
								(-A)^{-2l} \hat{x}
						\|_{H_r} \\ \nonumber
					&+
						\|\hat{x}\|_{H_{r-2l}} \,
						\tilde{K}(\|x\|^2_{H_{-l}} \vee \|x\|_{H_{-l}}^2)
						\cdot
							(\|x\|_{H_{-l}} + \|\hat{x}\|_{H_{-l}}) \,
							\|x-\hat{x}\|_{H_{-l}}.
		\end{align}
		Combining \eqref{eq: DV /V Difference} 
			(with $r \leftarrow 1-\vartheta$)
		with \eqref{eq: f' x uniform cont},
		the fact that for all $r,s \in \R$ and all $x \in H$ it holds that
		$\|x\|_{H_s} \leq \sup_{n \in \N} (|\lambda_n|^{s-r} ) \|x\|_{H_r}$,
		the fact that $1-\vartheta-2l \leq \vartheta$,
		the fact that $-l \leq \nicefrac 12 \leq \vartheta$,
		the fact that $f'$ is continuous on $[0,\infty)$,
		and with the fact that $f''$ is continuous on $[1,\infty)$,
		verifies \eqref{eq: V continuous with respect to H1 rem}.
		Moreover, \eqref{eq: DV /V Difference} 
			(with $r \leftarrow -\kappa$),
		\eqref{eq: f' Lipschitz cont},
		the fact that for all $r,s \in \R$ and all $x \in H$ it holds that
		$\|x\|_{H_s} \leq \sup_{n \in \N} (|\lambda_n|^{s-r} ) \|x\|_{H_r}$,
		the fact that
		$-\kappa-2l \leq -l-\kappa \leq 0 \leq \nicefrac 12$,
		the fact that
		$-l \leq 0 \leq \nicefrac 12$,
		the fact that $f'$ is continuous on $[0,\infty)$,
		the fact that $f''$ is continuous on $(0,\infty)$,
		and the fact that $\lim_{t \downarrow 0} f''(t) t =0$
		proves \eqref{eq: V continuous with respect to H kappa rem}.
		In addition, note that we obtain from \eqref{eq: derivative of V},
		\eqref{eq: f and f' inequalities},
		and
		the monotonicity of $\tilde{K}$ that
		\begin{equation}
			\begin{split}
					&\left \|
						I_\mathbb{H}^{-1}(D_\mathbb{H} V)(x) 
					\right \|_{H_{-\kappa}}
				=
					\left \|
						I_\mathbb{H}^{-1} f'(\nicefrac 12 \|(-A)^{-l} x\|^2_H) \, 
							I_{\mathbb{H}} \, ((-A)^{-2l} x)
					\right \|_{H_{-\kappa}} \\
				={}
					&|f'(\nicefrac 12 \|(-A)^{-l} x\|^2_H)|
					\left \| 
							((-A)^{-2l} x)
					\right \|_{H_{-\kappa}} 
				={}
					|f'(\nicefrac 12 \| x\|^2_{H_{-l}})|
					\left \| 
						x
					\right \|_{H_{-\kappa-2l}} \\
				\leq{}
					&\tilde{K}(\| x\|^2_{H_{-l}})
					\sup_{n \in \N} (|\lambda_n|^{-\kappa - 2l})
					\left \| 
						x
					\right \|_{H}
				\leq
					\tilde{K}(\sup_{n \in \N} (|\lambda_n|^{-l}) \| x\|^2_{H})
					\sup_{n \in \N} (|\lambda_n|^{-\kappa - 2l})
					\left \| 
						x
					\right \|_{H}
			\end{split}
		\end{equation}
		and this together with $-l \leq 0$
		and with $-\kappa -2l \leq 0$
		implies \eqref{eq: V bounded with respect to H kappa rem},
		Furthermore,
		\eqref{eq: derivative of V},
		\eqref{eq: f and f' inequalities},
		and
		the monotonicity of $\tilde{K}$ yield
		that for all $x \in H$ it holds that
		\begin{align}
					&\left \|
						I_\mathbb{H}^{-1}(D_\mathbb{H} V)(x) 
					\right \|_{H_{1/2}}
				=
					\left \|
						I_\mathbb{H}^{-1} f'(\nicefrac 12 \|(-A)^{-l} x\|^2_H) \, 
							I_{\mathbb{H}} \, ((-A)^{-2l} x)
					\right \|_{H_{1/2}} \\ \nonumber
				={}
					&|f'(\nicefrac 12 \|(-A)^{-l} x\|^2_H)|
					\left \| 
							((-A)^{-2l} x)
					\right \|_{H_{1/2}} 
				={}
					|f'(\nicefrac 12 \| x\|^2_{H_{-l}})|
					\left \| 
						x
					\right \|_{H_{1/2-2l}} \\ \nonumber
				\leq{}
					&\tilde{K}(\| x\|^2_{H_{-l}})
					\sup_{n \in \N} (|\lambda_n|^{\nicefrac 12 - 2l-\vartheta})
					\left \| 
						x
					\right \|_{H_\vartheta} 
				\leq{}
					\tilde{K}(\sup_{n \in \N} (|\lambda_n|^{-l-\vartheta}) \| x\|^2_{H_\vartheta})
					\sup_{n \in \N} (|\lambda_n|^{\nicefrac 12 - 2l -\vartheta})
					\left \| 
						x
					\right \|_{H_\vartheta}.
		\end{align}
		Combining this, $-l-\vartheta \leq 0$,
		and $\nicefrac 12 -2l -\vartheta \leq 0$
		implies \eqref{eq: V bounded with respect to H 1/2 rem},
		which finishes the proof of Remark \ref{rem: assumption on V}.
	\end{proof}
\chapter{Existence of viscosity solutions}
Existence of viscosity solutions of 2nd-order PDEs in infinite-dimensional
Hilbert spaces has been established e.g.\@ in Lions \cite[Theorem 3]{Lions1989}
and in Ishii \cite[Theorem 4.2 and Theorem 6.3]{Ishii1993}
with Perron's method. In this chapter we prove existence of viscosity solutions 
using the stability of limits (Lemma \ref{l:limits.of.viscosity.solutions}). 
This allows us to get a representation
of the solution. Moreover, we do not need to assume that the initial function is bounded. 

In Section \ref{sec: Convergence of Galerkin approximations}
we recall several results of Cox et al.\@ \cite[Section 4]{CoxHutzenthalerJentzenVanNeervenWelti2020},
which we slightly generalized to fit our purposes.
In particular we prove the convergence of Galerkin approximations for SPDEs in probability
uniformly on bounded sets
(Theorem \ref{t: uniform convergence}). 
In Section \ref{sec: Existence of viscosity solutions} 
we prove Theorem \ref{thm: existence}, which shows
the existence
of viscosity solutions
of Kolmogorov equations of SPDEs under suitable assumptions
and that the solutions can be represented by the solution of the
corresponding SPDEs. For its proof we use Theorem 4.16 in 
Hairer, Hutzenthaler \& Jentzen \cite{HairerHutzenthalerJentzen2015}
to show that the Kolmogorov equation of every
Galerkin approximation has a finite-dimensional classical viscosity solution.
Then we apply Proposition \ref{prop: u proj is viscosity solution}
to extend these to viscosity solutions of a common infinite-dimensional
Hilbert space. To conclude the result we then use the convergence
of the Galerkin approximations (Theorem \ref{t: uniform convergence})
and the stability of viscosity solutions under limits
(Lemma \ref{l:limits.of.viscosity.solutions}).
This chapter is based on 
Hairer, Hutzenthaler \& Jentzen \cite[Section 4.4]{HairerHutzenthalerJentzen2015} and
Cox et al.\@ \cite[Section 4]{CoxHutzenthalerJentzenVanNeervenWelti2020}.
	
\section{Setting}
\label{ssec:setting}
Throughout this chapter the following setting
is frequently used.
Let 
$
		\mathbb{H}
	=
		( 
			H,
			\langle \cdot, \cdot \rangle_H, 
			\left\| \cdot \right\|_H
		)
$
and
$
		\mathbb{U}
	=
		( 
			U ,
			\langle \cdot, \cdot \rangle_U , 
			\left\| \cdot \right\|_U
		)
$
be real separable $\R$-Hilbert spaces,
let
$
		\mathbb{H}'
	=
		( 
			H',
			\left\| \cdot \right\|_{H'}
		)
$
and
$
		\mathbb{U}'
	=
		( 
			U',
			\left\| \cdot \right\|_{U'}
		)
$
be the corresponding dual spaces,
let $ (e_i)_{i \in \N} \subseteq H $ be an orthonormal basis of $ \mathbb{H} $,
let 
$ 
	T
		\in (0,\infty) 
$,
let 
$ ( \Omega, \mathcal{F}, \P ) $
be a probability space with a normal filtration 
$ ( \mathbb{F}_t )_{ t \in [0,T] } $, 
let 
$
  ( W_t )_{ t \in [0,T] } 
$ 
be an $ \operatorname{Id}_U $-cylindrical $ ( \mathbb{F}_t )_{ t \in [0,T] } $-Wiener
process,
let $\lambda \colon \N \to (- \infty, 0)$ be a function with 
$\lim_{n \to \infty} \lambda_n = -\infty$,
let
$
  A \colon D(A) \subseteq H \rightarrow H
$ 
be the linear operator such that
\begin{equation}
  D(A) 
  = 
    \biggl\{ 
      v \in H 
      \colon
      \sum^{\infty}_{ i = 1 } 
	\left| 
	  \lambda_i 
	  \langle e_i, v \rangle_H 
	\right|^2
	< \infty
    \biggr\}
\end{equation}
and such that for all $ v \in D(A) $ it holds that
\begin{equation}
    Av
   =
    \sum^{\infty}_{ i =1 } 
    \lambda_i \langle e_i, v \rangle_H \, e_i,
\end{equation}
let 
$ 
		\mathbb{H}_r
	=
		( H_r , \left< \cdot , \cdot \right>_{ H_r }, \left\| \cdot \right\|_{ H_r } ) 
$,
$ r \in \R $,
be a family of interpolation spaces associated with
$ - A $
(see, e.g., Definition~3.6.30 in Jentzen \cite{Jentzen2015})
and 
$ 
		\mathbb{H}'_r
	=
		( H'_r , \left\| \cdot \right\|_{ H'_r } ) 
$,
$r \in \R$,
be the corresponding dual spaces.
Then we can extend the operator $A$ and
the norms $\| \cdot \|_{H_r}$, $r \in \R$,
to an operator 
$
	A \colon \bigcup_{i=1}^{\infty} H_{-i} \to \bigcup_{i=1}^{\infty} H_{-i}
$ 
and 
semi-norms
$
	\| \cdot \|_{H_r} \colon \bigcup_{i=1}^{\infty} H_{-i} \to [0, \infty]
$,
$r \in \R$,
such that for all
$r \in \R$, $x \in H_{r}$ and all $y \in \bigcup_{i=1}^{\infty} H_{-i}$ 
it holds that
\begin{equation}
			\|y\|_{H_r}
		=
			\begin{cases} 
				\|y\|_{H_{r}} 
					& \textrm{ if } y \in H_{r} \\
				\infty 
					& \textrm{ if } y \notin H_{r}
			\end{cases}
\end{equation}
and that
\begin{equation}
		(A(x) = y) 
	\Leftrightarrow 
		(
				\lim_{\eps \downarrow 0} \sup \{
					\|A(\xi) - y\|_{H_{r-1}} \colon \xi \in D(A), ~\|x- \xi\|_{H_{r}} \leq \eps 
				\}
			=
				0
		).
\end{equation}
Moreover,
let $\mathscr{E}_r \colon [0, \infty) \to [0, \infty)$,
$r \in (0, \infty)$,
be the functions satisfying for all $r \in (0, \infty)$ and all
$x \in [0, \infty)$ that
\begin{equation}
		\mathscr{E}_r(x)
	= 
		\left(
			1+
			\sum^\infty_{n=1}
				\frac
					{(x^{2}\Gamma(r))^n}
					{\Gamma(nr+1)}
		\right)^{\nicefrac 12}
\end{equation}
(cf.\@ Chapter~7 in Henry~\cite{Henry1981} and, e.g., Definition~1.3.1 
in Jentzen \cite{Jentzen2015}),
let 
$ \gamma \in \R $, $ \alpha \in [0,1) $, 
$ \beta \in [ 0, \nicefrac{ 1 }{ 2 } ) $, 
$ \chi \in [ \beta, \nicefrac{ 1 }{ 2 } ) $,
$ 
		F \in 
			\C_{
				\mathbb{H}_\gamma, 
				\mathbb{H}_{\gamma - \alpha}
			}( 
				H_{ \gamma }, H_{ \gamma - \alpha } 
			) 
$,
$ 
	B \in 
		\C_{
			\mathbb{H}_{\gamma}, 
			\mathbb{HS}(\mathbb{U}, \mathbb{H}_{\gamma - \beta})
		}( 
			H_{ \gamma }, HS( \mathbb{U}, \mathbb{H}_{ \gamma - \beta } ) 
		) 
$
satisfy that for all $\mathbb{H}_\gamma$-bounded sets
$ E \subseteq H_{ \gamma }$
it holds that
\begin{equation}
\label{eq: Lipschitz bound for F and B}
  | (F|_E) |_{
    \C^1( E, \left\| \cdot \right\|_{  H_{ \gamma - \alpha }  } )
  }
  +
  | (B|_E) |_{
    \C^1( E, \left\| \cdot \right\|_{  HS( \mathbb{U}, \mathbb{H}_{ \gamma - \beta } )  } )
  }
  < \infty,
\end{equation}
let $\mathfrak{A}_N \subseteq \{e_i, i \in \N\}$, $N \in \N$, be
finite non-empty sets satisfying that 
$\cup_{N \in \N} \mathfrak{A}_N = \{e_i, i \in \N\}$
and that
$
	\forall i,j \in \N \colon
		(\mathfrak{A}_i \cap \mathfrak{A}_j) = \mathfrak{A}_{i \wedge j}
$,
let $ V_N \subseteq H_\gamma$,
$ N \in \N_0 $,
be the sets
satisfying that
$
  V_0 = H_\gamma
$
and that
for all $N \in \N$ it holds that
$
	V_N = \overline{\Span_\mathbb{H} (\mathfrak{A}_N)}_{\mathbb{H}}
$,
let $ \mathbb{V}_N$,
$ N \in \N$,
be the real linear subspaces of $\mathbb{H}$
satisfying for all $N \in \N$ that
$
		\mathbb{V}_N 
	= 
		(
			V_N, 
			\langle \cdot , \cdot \rangle_{H} \, |_{V_N \times V_N}, 
			\| \cdot \|_{H} \,|_{V_N}
		)
$,
let 
$ \mathscr{P}_N \in L( U ) $,
$ N \in \N_0 $,
and let 
$
  X^{N,x} \colon [0,T] \times \Omega \rightarrow H_{ \gamma }
$,
$ N \in \N_0 $,
$x \in H_{\gamma}$,
be
$
  (\mathbb{F}_t )_{ t \in [0,T] }
$-adapted stochastic processes 
with continuous sample paths satisfying
for all $ N \in \N_0 $, $ t \in [0,T] $, and all $x \in H_\gamma$
 $\P$-a.s.\@ that
\begin{equation}
\label{eq:solution}
		X_t^{N,x}
  = 
		e^{ t A } \pi^{H_\gamma}_{V_N} x
		+
		\int_0^t
			e^{ ( t - s ) A }
			\pi^{H_{\gamma-\alpha}}_{V_N} F( X_s^{N,x}  )
		\ds
		+
		\int_0^t
			e^{ ( t - s ) A }
			\pi^{H_{\gamma-\beta}}_{V_N} B( X_s^{N,x}  ) \mathscr{P}_N
		\dWs.
\end{equation}
	\section{Convergence of Galerkin approximations}
	\label{sec: Convergence of Galerkin approximations}
	The following
	Lemma \ref{l: Lp estimate for difference} is a minor improvement on
	Lemma 4.2 in Cox et al.\@ \cite{CoxHutzenthalerJentzenVanNeervenWelti2020}
	since we allow different (but deterministic) starting points on the left hand side of
	\eqref{eq: Lp bound for different starting points}.
	The proof follows closely the proof of
	Lemma 4.2 in Cox et al.\@ \cite{CoxHutzenthalerJentzenVanNeervenWelti2020}
	and in some parts literally.
	\begin{lemma}
		[$L^p$-approximation bound for the Galerkin approximation for SPDEs with globally Lipschitz continuous coefficients]
		\label{l: Lp estimate for difference}
		Assume the setting in 
		Section~\ref{ssec:setting},
		let 
		$ p \in [2,\infty) $,
		$ 
			\eta \in \big[ \max\{ \alpha, 2 \beta \}, 1 \big)
		$,
		$x, \hat{x} \in H_\gamma $,
		$ N \in \N_0 $,
		and assume that 
		\begin{equation}
			| F |_{
				\C^1( H_{ \gamma }, \left\| \cdot \right\|_{  H_{ \gamma - \alpha }  } )
			}
			+
			| B |_{
				\C^1( H_{ \gamma }, \left\| \cdot \right\|_{  HS( \U, \H_{ \gamma - \beta } )  } )
			}
			< \infty.
		\end{equation}
		Then 
		\begin{equation}
		\label{eq: Lp bound for different starting points}
			\begin{split}
				&
					\sup_{ t \in [0,T] }
					\left\|
						X^{0, \hat{x}}_t - X^{N, x}_t
					\right\|_{ 
						L^p( 
							\P ; 
							\left\| \cdot \right\|_{  H_{ \gamma }  }
						) 
					} \\
					\leq{}
					&\Bigg[
						\sqrt{ 2 }
						\sup_{ t \in [0,T] }
						\left\| ( \id_{H} - \pi^{H_\gamma}_{V_N} ) X_t^{0, \hat{x}} \right\|_{
							L^p( \P ; \left\| \cdot \right\|_{  H_{ \gamma }  } )
						}
						+
						\sqrt{2} \,
						\|
							\pi^{H_\gamma}_{V_N} (\hat{x} - x)
						\|_{ H_{ \gamma } }
						+
						\tfrac{
							T^{ 1 / 2 - \chi }
							\sqrt{ p \, ( p - 1 ) } 
						}{
							\sqrt{ 1 - 2 \chi }
						}
				\\ &
							\left(
								1 +
								\sup_{ t \in [0,T] }
								\| X^{N, x}_t \|_{
									L^p( \P ; \left\| \cdot \right\|_{  H_{ \gamma }  } )
								}
							\right)
							\left(
								\sup_{ v \in H_{ \gamma } }
								\frac{
									\| B( v ) ( \mathscr{P}_0 - \mathscr{P}_N ) \|_{
										HS( \U, \H_{ \gamma - \chi } )
									}
								}{
									1 + \| v \|_{ H_{ \gamma } }
								}
							\right)
					\Bigg]
				\\ & 
					\cdot 
					\mathscr{E}_{ (1 - \eta) }\!\left[ 
					\tfrac{ 
						T^{ 1 - \eta }
						\sqrt{ 2 }
						\,
						\left| F \right|_{
							\C^1(
								H_{ \gamma }, \left\| \cdot \right\|_{  H_{ \gamma - \eta }  } 
							)
						}
					}{ \sqrt{ 1 - \eta } }
						+
						\sqrt{ 
							T^{ 1 - \eta }
							p ( p - 1 ) 
						} 
					\left| B \right|_{
						\C^1(
							H_{ \gamma }, \left\| \cdot \right\|_{  HS( \U, \H_{ \gamma - \eta / 2 } )  } 
						)
					}
					\| \mathscr{P}_0 \|_{ L( U ) }
					\right] \\
					<{} &\infty
					.
			\end{split}
		\end{equation}
	\begin{proof}[Proof of Lemma \ref{l: Lp estimate for difference}]
		First of all, observe that
		Lemma 4.1 in Cox et al.\@ \cite{CoxHutzenthalerJentzenVanNeervenWelti2020}
		ensures that
		\begin{align}
		\label{eq:finiteness_XN}
			\sup_{ t \in [0,T] }
				\max \bigl\{
					\| X^{0, \hat{x}}_t \|_{ L^p( \P; \left\| \cdot \right\|_{ H_{ \gamma } } ) },
					\| X^{N,x}_t \|_{ L^p( \P; \left\| \cdot \right\|_{ H_{ \gamma } } ) }
				\bigr\}
				< \infty
				.
		\end{align}
		Hence we can apply 
		Proposition~7.1.6 in Jentzen \cite{Jentzen2015} 
		to obtain
		\begin{align}
		& 
			\sup_{ t \in [0,T] }
			\big\|
				X^{0,\hat{x}}_t - X^{N, x}_t
			\big\|_{ 
				L^p( 
					\P ; 
					\left\| \cdot \right\|_{  H_{ \gamma }  }
				) 
			}
		\\  \nonumber\leq{}&
			\mathscr{E}_{ (1 - \eta) } \bigg[ 
			\tfrac{ 
				T^{ 1 - \eta }
				\sqrt{ 2 }
				\,
				\big| \pi^{H_{\gamma-\alpha}}_{V_N} F( \cdot ) \big|_{
					\C^1(
						H_{ \gamma }, \left\| \cdot \right\|_{  H_{ \gamma - \eta }  } 
					)
				}
			}{ \sqrt{ 1 - \eta } }
				+
				\sqrt{ 
					T^{ 1 - \eta }
					p ( p - 1 ) 
				} 
			\big| \pi^{H_{\gamma-\beta}}_{V_N} B( \cdot ) \mathscr{P}_0 \big|_{
				\C^1(
					H_{ \gamma }, \left\| \cdot \right\|_{  HS( \U, \H_{ \gamma - \eta / 2 } )  } 
				)
			}
			\bigg]
		\\ \nonumber & \quad 
			\cdot
			\sqrt{2} 
			\,
			\sup_{ t \in [0,T] }
			\Big\| 
				X_t^{0, \hat{x}}
				-
					\int_0^t
					e^{ ( t - s ) A } \pi^{H_{\gamma-\alpha}}_{V_N} F( X^{0, \hat{x}}_s ) 
					\ds
					-
					\int_0^t
					e^{ ( t - s ) A } \pi^{H_{\gamma-\beta}}_{V_N} B( X^{0, \hat{x}}_s ) \mathscr{P}_0
					\dWs
		\\ \nonumber & \qquad\qquad
			+
					\int_0^t
					e^{ ( t - s ) A } \pi^{H_{\gamma-\alpha}}_{V_N} F( X^{N, x}_s ) 
					\ds
					- X^{N,x}_t
					+
					\int_0^t
					e^{ ( t - s ) A } \pi^{H_{\gamma-\beta}}_{V_N} B( X^{N,x}_s ) \mathscr{P}_0
					\dWs
			\Big\|_{ 
				L^p( 
					\P ; 
					\left\| \cdot \right\|_{  H_{ \gamma }  } 
				) 
			}
		\\ \nonumber & \leq{}
			\mathscr{E}_{ (1 - \eta) } \bigg[ 
			\tfrac{ 
				T^{ 1 - \eta }
				\sqrt{ 2 }
				\,
				\big| \pi^{H_{\gamma-\alpha}}_{V_N} F( \cdot ) \big|_{
					\C^1(
						H_{ \gamma }, \left\| \cdot \right\|_{  H_{ \gamma - \eta }  } 
					)
				}
			}{ \sqrt{ 1 - \eta } } 
				+
				\sqrt{ 
					T^{ 1 - \eta }
					p ( p - 1 ) 
				} 
			\big| \pi^{H_{\gamma-\beta}}_{V_N} B( \cdot ) \mathscr{P}_0 \big|_{
				\C^1(
					H_{ \gamma }, \left\| \cdot \right\|_{  HS( \U, \H_{ \gamma - \eta / 2 } )  } 
				)
			}
			\bigg]
		\\ \nonumber & \quad 
			\cdot
			\sqrt{2} 
			\,
			\sup_{ t \in [0,T] }
			\Big\|
				X_t^{0, \hat{x}}
				-
					\int_0^t
					e^{ ( t - s ) A } \pi^{H_{\gamma-\alpha}}_{V_N} F( X^{0, \hat{x}}_s ) 
					\ds
					-
					\int_0^t
					e^{ ( t - s ) A } \pi^{H_{\gamma-\beta}}_{V_N} B( X^{0, \hat{x}}_s ) \mathscr{P}_0
					\dWs
		\\ \nonumber & \qquad\qquad
			- e^{tA} \pi^{H_\gamma}_{V_N} \hat{x}
			+ e^{tA} \pi^{H_\gamma}_{V_N} \hat{x}
			- e^{tA} \pi^{H_\gamma}_{V_N} x
			- \int_0^t e^{(t-s) A}  \pi^{H_{\gamma-\beta}}_{V_N} B(X_s^{N, x}) \mathscr{P}_N \dWs \\ 
					& \nonumber \qquad\qquad
					+
					\int_0^t
					e^{ ( t - s ) A } \pi^{H_{\gamma-\beta}}_{V_N} B( X^{N,x}_s ) \mathscr{P}_0
					\dWs
			\Big\|_{ 
				L^p( 
					\P ; 
					\left\| \cdot \right\|_{  H_{ \gamma }  } 
				) 
			}	.
		\end{align}
		This shows that
		\begin{equation}
		\begin{split}
		&
			\sup_{ t \in [0,T] }
			\big\|
				X^{0, \hat{x}}_t - X^{N,x}_t
			\big\|_{ 
				L^p( 
					\P ; 
					\left\| \cdot \right\|_{  H_{ \gamma }  } 
				) 
			}
		\\ \leq{}
			&\mathscr{E}_{ (1 - \eta) } \bigg[ 
			\tfrac{ 
				T^{ 1 - \eta }
				\sqrt{ 2 }
				\,
				\left| F \right|_{
					\C^1(
						H_{ \gamma }, \left\| \cdot \right\|_{  H_{ \gamma - \eta }  } 
					)
				}
			}{ \sqrt{ 1 - \eta } }
				+
				\sqrt{ 
					T^{ 1 - \eta }
					p ( p - 1 ) 
				} 
			\left| B \right|_{
				\C^1(
					H_{ \gamma }, \left\| \cdot \right\|_{  HS( \U, \H_{ \gamma - \eta / 2 } )  } 
				)
			}
			\| \mathscr{P}_0 \|_{ L( U ) }
			\bigg]
		\\ & 
			\cdot
			\sqrt{2} 
			\,
			\sup_{ t \in [0,T] } \bigg [
				\Big \|
					( \id_{H_\gamma} - \pi^{H_\gamma}_{V_N} ) X^{0, \hat{x}}_t
			\\ & \qquad 
					+
						\int_0^t
							e^{ ( t - s ) A } \pi^{H_{\gamma-\beta}}_{V_N} B( X^{N, x}_s ) 
								( \mathscr{P}_0 - \mathscr{P}_N )
						\dWs
				\Big\|_{ 
					L^p( 
						\P ; 
						\left\| \cdot \right\|_{  H_{ \gamma }  } 
					) 
				}
				+ \|e^{tA} (\pi^{H_\gamma}_{V_N} \hat{x} - \pi^{H_\gamma}_{V_N} x) \|_{H_\gamma}
			\bigg].
		\end{split}
		\end{equation}
		The Burkholder-Davis-Gundy type inequality
		in Lemma~7.7 in Da Prato \& Zabczyk~\cite{DaPratoZabczyk1992}
		hence implies that
		\begin{align}
		& \nonumber
			\sup_{ t \in [0,T] }
			\big\|
				X^{0, \hat{x}}_t - X^{N, x}_t
			\big\|_{ 
				L^p( 
					\P ; 
					\left\| \cdot \right\|_{ 
						H_{ \gamma } 
					} 
				) 
			}
		\\ \nonumber \leq{}
			&\mathscr{E}_{ (1 - \eta) } \bigg[ 
			\tfrac{ 
				T^{ 1 - \eta }
				\sqrt{ 2 }
				\,
				\left| F \right|_{
					\C^1(
						H_{ \gamma }, \left\| \cdot \right\|_{  H_{ \gamma - \eta }  } 
					)
				}
			}{ \sqrt{ 1 - \eta } }
				+
				\sqrt{ 
					T^{ 1 - \eta }
					p ( p - 1 ) 
				} 
			\left| B \right|_{
				\C^1(
					H_{ \gamma }, \left\| \cdot \right\|_{  HS( \U, \H_{ \gamma - \eta / 2 } )  } 
				)
			}
			\| \mathscr{P}_0 \|_{ L( U ) }
			\bigg]
		\\ & \quad
			\cdot
			\sqrt{2} 
			\,
			\bigg[
				\sup_{ t \in [0,T] }
				\left\| ( \id_{H_\gamma} - \pi^{H_\gamma}_{V_N} ) X_t^{0, \hat{x}} \right\|_{
					L^p( \P ; \left\| \cdot \right\|_{  H_{ \gamma }  } )
				}
			+ \|\pi^{H_\gamma}_{V_N} (\hat{x} -  x) \|_{H_\gamma}
		\\ \nonumber & \quad
				+
					\sup_{ s \in [0,T] }
					\left\|
						B( X^{N, x}_s ) [ \mathscr{P}_0 - \mathscr{P}_N ]
					\right\|_{
						L^p( \P ; \left\| \cdot \right\|_{ HS( \U, \H_{ \gamma - \chi } ) } )
					}
				\sqrt{
					\tfrac{ p \, ( p - 1 ) }{ 2 }
					\sup_{ t \in [0,T] }
					\int_0^t
					\left( t - s \right)^{ - 2 \chi }
					\ds
				}
			\bigg]
			.
		\end{align}
		This and \eqref{eq:finiteness_XN} complete the proof
		of Lemma \ref{l: Lp estimate for difference}. 
	\end{proof}
	\end{lemma}
	The following Lemma \ref{l: Lp convergence for X_N} 
	(resp.\@ Corollary \ref{cor:convGlobLip}) complements the result
	on positive convergence rates established in Corollary 4.3 
	(resp.\@ Corollary 4.4)
	in Cox et al.\@
	\cite{CoxHutzenthalerJentzenVanNeervenWelti2020} with converging starting points
	(without positive rate).
	The proof follows closely 
	Corollary 4.3 in Cox et al.\@ \cite{CoxHutzenthalerJentzenVanNeervenWelti2020}
	and in some parts literally.
	\begin{lemma}
		[$L^p$-convergence 
			of Galerkin approximations for SPDEs with globally Lipschitz continuous coefficients]
	\label{l: Lp convergence for X_N}
		Assume the setting in 
		Section~\ref{ssec:setting},
		let 
		$\vartheta \in \big( 0, \min \{1- \alpha, \tfrac 12 - \beta \} \big)$,
		$ p \in [2,\infty) $,
		$
			(x_N)_{N \in \N_0} 
				\subseteq H_{\gamma+ \vartheta}
		$,
		and assume that 
		$
				\limsup_{N \to \infty} 
					\|
						x_N - x_0 
					\|_{{  H_{ \gamma}  } }
			= 0
		$,
		\begin{equation}
			\sup_{N \in \N_0} \left(
					\| x_N \|_{ H_{ \gamma + \vartheta } }
			\right)
			+
			| F |_{
				\C^1( H_{ \gamma }, \left\| \cdot \right\|_{  H_{ \gamma - \alpha }  } )
			}
			+
			| B |_{
				\C^1( H_{ \gamma }, \left\| \cdot \right\|_{  HS( \U, \H_{ \gamma - \beta } )  } )
			}
			< \infty,
		\end{equation}
		and that
		\begin{equation}
		\label{eq: B C1 convergence}
			\limsup_{N \to \infty} \left(
				\sup_{ v \in H_{ \gamma } }
					\frac
						{ 
							\|
								B( v ) ( \mathscr{P}_0 - \mathscr{P}_N ) 
							\|_{HS( \U, \H_{ \gamma - \chi } )}
						}
						{
							1 + \| v \|_{ H_{ \gamma } }
						}
						\right)
			=0.
		\end{equation}
		Then it holds that
		\begin{equation}
			\begin{split}
					\sup_{N \in \N_0} \sup_{t \in [0,T]}(
						\| F(X_t^N) \|_{ 
								L^p( 
									\P ; 
									\left\| \cdot \right\|_{  H_{ \gamma -\alpha }  }
								) 
							}
						+\| B(X_t^N) \mathscr{P}_N\|_{ 
								L^p( 
									\P ; 
									\left\| \cdot \right\|_{  HS(\U,\H_{\gamma -\chi} )  }
								) 
							}
					)
				< \infty
			\end{split}
		\end{equation}
		and that
		\begin{equation}
				\limsup_{N \to \infty}
					\sup_{ t \in [0,T] }
						\left\|
							X^{0,x_0}_t - X^{N,x_N}_t
						\right\|_{ 
							L^p( 
								\P ; 
								\left\| \cdot \right\|_{  H_{ \gamma }  }
							) 
						}
				= 0.
		\end{equation}
	\end{lemma}
	\begin{proof}
		Fix 
		$ 
			\eta \in \big[ \max\{ \alpha, 2 \beta \}, 1 \big)
		$
		for the rest of the proof
		and note that we get from Lemma \ref{l: Lp estimate for difference} 
		that for all $N \in \N_0$ it holds that
		\begin{align}
		\label{eq: lemma estimate}
			&
				\sup_{ t \in [0,T] }
				\big\|
					X^{0,x_0}_t - X^{N,x_N}_t
				\big\|_{ 
					L^p( 
						\P ; 
						\left\| \cdot \right\|_{  H_{ \gamma }  }
					) 
				} \\ \nonumber
				\leq{} 
				&\bigg[
					\sqrt{ 2 }
					\sup_{ t \in [0,T] }
					\left\| ( \id_{H_\gamma} - \pi^{H_\gamma}_{V_N} ) X_t^{0, x_0} \right\|_{
						L^p( \P ; \left\| \cdot \right\|_{  H_{ \gamma }  } )
					}
				+ \sqrt{ 2 } \,
					\|
						\pi^{H_\gamma}_{V_N} (x_0 - x_N)
					\|_{ H_{ \gamma } }
			\\ \nonumber & 
					+
					\tfrac{
						T^{ 1 / 2 - \chi }
						\sqrt{ p \, ( p - 1 ) } 
					}{
						\sqrt{ 1 - 2 \chi }
					}
						\big(
							1 +
							\sup_{ t \in [0,T] }
							\| X^{N, x_N}_t \|_{
								L^p( \P ; \left\| \cdot \right\|_{  H_{ \gamma }  } )
							}
						\big)
						\big(
							\sup_{ v \in H_{ \gamma } }
							\tfrac{
								\| B( v ) ( \mathscr{P}_0 - \mathscr{P}_N ) \|_{
									HS( \U, \H_{ \gamma - \chi } )
								}
							}{
								1 + \| v \|_{ H_{ \gamma } }
							}
						\big)
				\bigg]
			\\ \nonumber & 
				\cdot 
				\mathscr{E}_{ (1 - \eta) } \bigg[ 
				\tfrac{ 
					T^{ 1 - \eta }
					\sqrt{ 2 }
					\,
					\left| F \right|_{
						\C^1(
							H_{ \gamma }, \left\| \cdot \right\|_{  H_{ \gamma - \eta }  } 
						)
					}
				}{ \sqrt{ 1 - \eta } }
					+
					\sqrt{ 
						T^{ 1 - \eta }
						p ( p - 1 ) 
					} 
				\left| B \right|_{
					\C^1(
						H_{ \gamma }, \left\| \cdot \right\|_{  HS( \U, \H_{ \gamma - \eta / 2 } )  } 
					)
				}
				\| \mathscr{P}_0 \|_{ L( U ) }
				\bigg]
				< \infty
				.
		\end{align}
		In addition, observe that Lemma 4.1 in 
		Cox et al.\@ \cite{CoxHutzenthalerJentzenVanNeervenWelti2020}
		and $\sup_{N \in \N_0} \|x_N\|_{H_\gamma} < \infty$ 
		imply
		that
		\begin{equation}
		\label{eq: X_N bound in N}
				\sup_{N \in \N_0} \sup_{ t \in [0,T] }
					\| X^{N, x_N}_t \|_{ L^p( \P; \left\| \cdot \right\|_{ H_{ \gamma } } ) }
			< \infty.
		\end{equation}
		Moreover, we get from 
		$
			x_0
				\in H_{\gamma+ \vartheta}
		$
		that
		$
			\forall t \in [0,T] \colon 
				\P(X^{0, x_0}_t \in H_{\gamma + \vartheta})=1
		$
		and that
		\begin{equation}	
				\sup_{N \in \N_0} \sup_{t \in [0, T]} 
					\E[ 
						\| 
							\1_{X_t^{0, x_0} \in H_{\gamma+\vartheta}} X_t^{0, x_0}
						\|^p_{H_{\gamma + \vartheta}} 
					]
			< \infty
		\end{equation}
		(see eg. Proposition 7.1.19 in Jentzen \cite{Jentzen2015} 
		and Proposition 7.1.10 in Jentzen \cite{Jentzen2015}).
		Therefore
		it follows from 
		$
				\limsup_{N \to \infty} 
					\|
						x_N - x_0 
					\|_{H_{ \gamma } }
			= 0
		$,
		from the definition of $\mathfrak{A}_N$, $N \in \N$
		and from
		$
				\limsup_{N \to \infty}
					\lambda_N
			= - \infty
		$
		that
		\begin{align}
		\label{eq: projection and start value convergence}
				& 
					\limsup_{ N \to \infty }
					\sup_{ t \in [0,T] }
					\bigg[
						\left\|
							( \id_{H_\gamma} - \pi^{H_\gamma}_{V_N} ) X^{0, x_0}_t
						\right\|_{
							L^p( \P ; \left\| \cdot \right\|_{  H_{ \gamma }  } )
						}
					+\|
						\pi^{H_\gamma}_{V_N} (x_0 - x_N)
					\|_{ H_{ \gamma }}
					\bigg]
				\\ \nonumber 
				\leq{} 
					&\limsup_{ N \to \infty }
					\sup_{ t \in [0,T] }
					\bigg[
						\left\|
							( 
								\id_{H_{\gamma}}
								- \pi^{H_\gamma}_{V_N}
							)
							( - A )^{ - \vartheta }
						\right\|_{
							L( \mathbb{H}_{ \gamma}, \mathbb{H}_{ \gamma} )
						}
						\big\|
							\mathbbm{1}_{ \{ X^{0, x_0}_t \in H_{ \gamma + \vartheta } \} } 
								X^{0, x_0}_t
						\big\|_{
							L^p( \P ; \left\| \cdot \right\|_{ H_{ \gamma + \vartheta } } )
						}
			\\ \nonumber & \qquad\qquad
					+\|
						(x_0 - x_N)
					\|_{ H_{ \gamma }}
					\bigg]
				\\ \nonumber 
				\leq{} 
					&\bigg[
						\limsup_{ N \to \infty }
						\left \|
							( \operatorname{Id}_{ H_{ \gamma } } - \pi^{H_\gamma}_{V_N} )
							( - A )^{-1}
						\right \|^{\vartheta }_{L( \mathbb{H}_{ \gamma}, \mathbb{H}_{ \gamma} )}
					\bigg]
					\bigg[
						\sup_{N \in \N}
						\sup_{ t \in [0,T] }
						\big\|
							\mathbbm{1}_{ \{ X^{0, x_0}_t \in H_{ \gamma + \vartheta } \} } 
								X^{0, x_0}_t
						\big\|_{
							L^p( \P ; \left\| \cdot \right\|_{ H_{ \gamma + \vartheta } } )
						}
					\bigg]
				\\ \nonumber 
				={} 
					&\bigg[
						\limsup_{ N \to \infty }
						\big[
							\sup \{ \nicefrac {1}{|\lambda_N|} \colon e_N \in \mathfrak{A}_N \} 
						\big]^{ \vartheta }
					\bigg]
					\bigg[
						\sup_{ N \in \N }
						\sup_{ t \in [0,T] }
						\big\|
							\mathbbm{1}_{ \{ X^{0, x_0}_t \in H_{ \gamma + \vartheta } \} }
								X^{0, x_0}_t
						\big\|_{
							L^p( \P ; \left\| \cdot \right\|_{ H_{ \gamma + \vartheta } } )
						}
					\bigg]
				\\ \nonumber 
				={} &0 .
		\end{align}
		Furthermore, we get from \eqref{eq: B C1 convergence} that
		\begin{align}
		\nonumber
					&\sup_{N \in \N_0} \sup_{t \in [0,T]}
					\left(
						\|
							B( X_t^N ) \mathscr{P}_N 
						\|_{
							L^p( \P ; \left\| \cdot \right\|_{ HS( \U, \H_{ \gamma - \chi } ) } )
						}
					\right)\\ \nonumber
				\leq{}
					&\left(
						1
						+\sup_{N \in \N_0} \sup_{t \in [0,T]} 
							\|X_t^N\|_{
								L^p( \P ; \left\| \cdot \right\|_{ H_{ \gamma} } )
							}
					\right)
					\left(
						\sup_{N \in \N_0}
						\sup_{ v \in H_{ \gamma } }
							\frac
								{ 
									\|
										B( v ) \mathscr{P}_N 
									\|_{HS( \U, \H_{ \gamma - \chi } )}
								}
								{
									1 + \| v \|_{ H_{ \gamma } }
								}
					\right) \\ \label{eq: B C1 bound}
				\leq{}
					&\left(
						1
						+\sup_{N \in \N_0} \sup_{t \in [0,T]} 
							\|X_t^N\|_{
								L^p( \P ; \left\| \cdot \right\|_{ H_{ \gamma} } )
							}
					\right)
			\\ \nonumber &
					\cdot \left(
						\sup_{ v \in H_{ \gamma } }
							\frac
								{ 
									\|
										B( v ) \mathscr{P}_0 
									\|_{HS( \U, \H_{ \gamma - \chi } )}
								}
								{
									1 + \| v \|_{ H_{ \gamma } }
								}
						+\sup_{N \in \N_0}
						\sup_{ v \in H_{ \gamma } }
							\frac
								{ 
									\|
										B( v ) (\mathscr{P}_0 - \mathscr{P}_N)
									\|_{HS( \U, \H_{ \gamma - \chi } )}
								}
								{
									1 + \| v \|_{ H_{ \gamma } }
								}
					\right) \\ \nonumber
				<{} &\infty.
		\end{align}
		Combining then \eqref{eq: B C1 convergence},
		\eqref{eq: lemma estimate}, \eqref{eq: X_N bound in N}, and
		\eqref{eq: projection and start value convergence}
		shows that
		\begin{equation}
		\label{eq: Lp convergence for X}
				\limsup_{N \to \infty}
					\sup_{ t \in [0,T] }
						\left\|
							X^{0,x_0}_t - X^{N,x_N}_t
						\right\|_{ 
							L^p( 
								\P ; 
								\left\| \cdot \right\|_{  H_{ \gamma }  }
							) 
						}
				= 0.
		\end{equation}
		Finally note that \eqref{eq: X_N bound in N} verifies that
		$
				\sup_{N \in \N_0} \sup_{t \in [0,T]}(
					\| F(X_t^N) \|_{ 
						L^p( 
							\P ; 
							\left\| \cdot \right\|_{  H_{ \gamma -\alpha }  }
						) 
					}
				)
			< \infty
		$
		and this together with \eqref{eq: B C1 bound}
		and \eqref{eq: Lp convergence for X}
		completes the proof
		of Lemma \ref{l: Lp convergence for X_N}.
	\end{proof}
	The proof of the following corollary follows closely 
	Corollary 4.4 in Cox et al.\@ \cite{CoxHutzenthalerJentzenVanNeervenWelti2020}
	and in some parts literally.
	\begin{corollary}
			[$L^p$-H\"older-convergence 
			of Galerkin approximations for SPDEs with globally Lipschitz continuous coefficients]
	\label{cor:convGlobLip}
	Assume the setting in Section~\ref{ssec:setting},
	let 
	$
		\vartheta \in \big( 0, \min\{ 1 - \alpha, \nicefrac{ 1 }{ 2 } - \beta \} \big)
	$,
	$ p \in ( \nicefrac{ 1 }{ \vartheta }, \infty ) $,
	$
		(x_N)_{N \in \N_0}
			\subseteq H_{\gamma+ \vartheta}
	$,
	and assume that,
		$
				\limsup_{N \to \infty} 
					\|
						x_N - x_0 
					\|_{{  H_{ \gamma}  } }
			= 0
		$,
		$
			\sup_{ N \in \N } \sup_{ v \in H_{ \gamma } }
			\left[
				\frac{
				\|
					B( v ) \mathscr{P}_N 
				\|_{ HS( \U, \H_{ \gamma - \beta } ) }
				}{
					1 + \| v \|_{ H_{ \gamma } }
				}
			\right]
			< \infty
		$,
		\begin{equation}
			\sup_{N \in \N_0} \left(
					\| x_N \|_{ H_{ \gamma + \vartheta } }
			\right)
			+
			| F |_{
				\C^1( H_{ \gamma }, \left\| \cdot \right\|_{  H_{ \gamma - \alpha }  } )
			}
			+
			| B |_{
				\C^1( H_{ \gamma }, \left\| \cdot \right\|_{  HS( \U, \H_{ \gamma - \beta } )  } )
			}
			< \infty,
		\end{equation}
	and that
	\begin{equation}
	\label{eq:ass_B1.general_0b}
		\limsup_{N \to \infty} \left(
				\sup_{ v \in H_{ \gamma } }
					\frac
						{ 
							\|
								B( v ) ( \mathscr{P}_0 - \mathscr{P}_N ) 
							\|_{HS( \U, \H_{ \gamma - \chi } )}
						}
						{
							1 + \| v \|_{ H_{ \gamma } }
						}
						\right)
			=0.
	\end{equation}
	Then for all 
	$ \delta \in [ 0, \vartheta - \nicefrac{ 1 }{ p } ) $
	it holds that
	\begin{equation}
		\sup_{ N \in \N_0 } 
			\|
				X^{N,x_N}
			\|_{ 
				L^p( 
					\P ; 
					\left\| \cdot \right\|_{
					\C^{ \delta }( [0,T], \left\| \cdot \right\|_{  H_{ \gamma }  } )
					}
				)
			}
		< \infty
	\end{equation}
	and that
	\begin{equation}
		\limsup_{ N \to \infty } 
		\|
			X^{0,x_0} - X^{N,x_N}
		\|_{ 
			L^p( 
				\P ;
				\left\| \cdot \right\|_{
				\C^{ \delta }( [0,T], \left\| \cdot \right\|_{  H_{ \gamma }  } )
				}
			)
		}
		=0
		.
	\end{equation}
	\end{corollary}
	\begin{proof}[Proof of Corollary \ref{cor:convGlobLip}]
	Throughout this proof denote by $ \eta \in \R $ 
	the real number satisfying that
	$
		\eta = \max\{ \alpha, 2 \beta \}
	$.
	Note that 
	Lemma~\ref{l: Lp convergence for X_N}
	proves that
	\begin{equation}
	\label{eq:proof_Lip_to_apply_hoelder}
		\begin{split}
				\limsup_{ N \to \infty } \left[
					\sup_{ t \in [0,T] }
					\|
						X^{0, x_0}_t - X^{N, x_N}_t
					\|_{
						L^p( \P; \left\| \cdot \right\|_{  H_{ \gamma }  } ) 
					}
				\right] 
			= 0.
		\end{split}
	\end{equation}
	In the next step note that, e.g.,
	(Corollary~11.3.2 in Jentzen \cite{Jentzen2015a})
	yields that
	for all $ N \in \N_0 $,
	$ 
		\varepsilon \in 
		\big( 
			0, \min\{ 1 + \gamma - \eta , \nicefrac{ 1 }{ 2 } + \gamma - \beta \}
			- \gamma
		\big)
	$
	it holds that
	\begin{align}
	&
		\sup_{ 
			\substack{  
				t_1, t_2 \in [0,T] 
				,
				\\
				t_1 \neq t_2
			}
		}
		\left(
			\frac{
				\left|
					\min\{ t_1, t_2 \}
				\right|^{ \max\{ \gamma + \varepsilon - ( \gamma + \vartheta ) , 0 \} }
				\left\| X^{N,x_N}_{ t_1 } - X^{N,x_N}_{ t_2 } \right\|_{
					L^p( \P; \left\| \cdot \right\|_{  H_{ \gamma }  } )
				}
			}{
				\left| t_1 - t_2 \right|^{ \varepsilon }
			}
		\right)
	\\ \nonumber & 
		\leq{} 
			\left\| 
				x_N 
			\right\|_{  H_{ \min\{ \gamma + \vartheta, \gamma + \varepsilon \} }  }
		+
		\left[
			\sup_{ s \in [0,T] }
			\left\|
				\pi^{H_{\gamma-\alpha}}_{V_N} F( X_s^{N,x_N} ) 
			\right\|_{
				L^p( \P; \left\| \cdot \right\|_{  H_{ \gamma - \eta }  } )
			}
		\right]
			\frac{
				2
				\,
				T^{ 
					( 1 + \gamma - \eta - \min\{ \gamma + \vartheta, \gamma + \varepsilon \} ) 
				}
			}{
				\left( 1 - \eta - \varepsilon \right) 
			}
	\\ \nonumber & \quad
		+
		\left[
			\sup_{ s \in [0,T] }
			\left\|
				\pi^{H_{\gamma-\beta}}_{V_N} B( X_s^{N, x_N} ) \mathscr{P}_N
			\right\|_{
				L^p( \P; \left\| \cdot \right\|_{  HS( \U, \H_{ \gamma - \beta } )  } )
			}
		\right]
			\frac{ 
				\sqrt{ 2 }
				\,
				T^{ ( 1 / 2 + \gamma - \beta - \min\{ \gamma + \vartheta , \gamma + \varepsilon \} ) }
			}{
				\left( 
					1 - 2 \beta - 2 \varepsilon 
				\right)^{ 1 / 2 }
			}
		< \infty
		.
	\end{align}
	This and the fact that
	$
			\min\{ 1 + \gamma - \eta , \nicefrac{ 1 }{ 2 } + \gamma - \beta \}
			- \gamma
		=
			\min\{ 1 - \max\{ \alpha, 2 \beta \} , \nicefrac{ 1 }{ 2 } - \beta \}
		=
			\min\{ 1 - \alpha, \nicefrac{ 1 }{ 2 } - \beta \}
		> \vartheta > 0
	$
	show that
	\begin{align}
	\nonumber
	&
		\sup_{ N \in \N_0 }
		\sup_{ 
			\substack{  
				t_1, t_2 \in [0,T] 
				,
				\\
				t_1 \neq t_2
			}
		}
		\left(
			\frac{
				\left\| X^{N,x_N}_{ t_1 } - X^{N,x_N}_{ t_2 } \right\|_{
					L^p( \P; \left\| \cdot \right\|_{  H_{ \gamma }  } )
				}
			}{
				\left| t_1 - t_2 \right|^{ \vartheta }
			}
		\right)
	\\ & 
		\leq 
			\sup_{ N \in \N_0 }
			\left\| 
				x_N 
			\right\|_{  H_{ \gamma + \vartheta }  }
		+
		\left[
			\sup_{ N \in \N_0 }
			\sup_{ s \in [0,T] }
			\left\|
				F( X_s^{N,x_N} ) 
			\right\|_{
				L^p( \P; \left\| \cdot \right\|_{  H_{ \gamma - \eta }  } )
			}
		\right]
			\frac{
				2
				\,
				T^{ 
					( 1 - \eta - \vartheta ) 
				}
			}{
				\left( 1 - \eta - \vartheta \right) 
			}
	\\ \nonumber & \quad
		+
		\left[
			\sup_{ N \in \N_0 }
			\sup_{ s \in [0,T] }
			\left\|
				B( X_s^{N,x_N} ) \mathscr{P}_N
			\right\|_{
				L^p( \P; \left\| \cdot \right\|_{  HS( \U, \H_{ \gamma - \beta } )  } )
			}
		\right]
			\frac{ 
				\sqrt{ 2 }
				\,
				T^{ 
					( 
						1 / 2 - \beta - \vartheta
					) 
				}
			}{
				\left( 
					1 - 2 \beta - 2 \vartheta 
				\right)^{ 1 / 2 }
			}
		.
	\end{align}
	Therefore
	Lemma~\ref{l: Lp convergence for X_N}
	and estimate~\eqref{eq:ass_B1.general_0b}
	ensure that
	\begin{align}
	& \nonumber
		\sup_{ N \in \N_0 }
		\left|
			X^{N,x_N}
		\right|_{
			\C^{ \vartheta }( [0,T] ,
					\left\| \cdot \right\|_{ L^p( \P ; \left\| \cdot \right\|_{  H_{ \gamma }  } )  } )
		}
		\\ &
	\leq{}
			\sup_{ N \in \N_0 }
				\left\| 
					x_N
				\right\|_{  H_{ \gamma + \vartheta }  }
		+
		\left[
			\sup_{ N \in \N_0 }
			\sup_{ s \in [0,T] }
			\left\|
				F( X_s^{N,x_N} ) 
			\right\|_{
				L^p( \P; \left\| \cdot \right\|_{  H_{ \gamma - \eta }  } )
			}
		\right]
			\tfrac{
				2
				\,
				T^{ 
					( 1 - \eta - \vartheta ) 
				}
			}{
				\left( 1 - \eta - \vartheta \right) 
			}
	\\ \nonumber & \quad
		+
		\left[
			\sup_{ N \in \N_0 }
			\sup_{ s \in [0,T] }
			\left\|
				B( X_s^{N,x_N} ) \mathscr{P}_N
			\right\|_{
				L^p( \P; \left\| \cdot \right\|_{  HS( \U, \H_{ \gamma - \beta } )  } )
			}
		\right]
			\tfrac{ 
				\sqrt{ 2 }
				\,
				T^{ 
					( 
						1 / 2 - \beta - \vartheta
					) 
				}
			}{
				\left( 
					1 - 2 \beta - 2 \vartheta 
				\right)^{ 1 / 2 }
			}
		< \infty
		.
	\end{align}
	Combining this with \eqref{eq:proof_Lip_to_apply_hoelder},
	the fact that
	$
		\vartheta \in ( \nicefrac{ 1 }{ p } , 1 ]
	$,
	and
	Corollary~2.11 in Cox et al.\@ \cite{CoxHutzenthalerJentzenVanNeervenWelti2020}
	shows that for all 
	$ \delta \in [ 0, \vartheta - \nicefrac{ 1 }{ p } ) $,
	$ \varepsilon \in (0,\infty) $,
	and all $(d_N)_{N \in \N_0} \subseteq (0, \infty)$ satisfying that
	$
			\lim_{N \to \infty} \Big(
				(d_N)^{-\vartheta}
				\|
					X^{0,\hat{x}_N} - X^{N,x_N}
				\|_{ 
					L^p( 
						\P ;
						\left\| \cdot \right\|_{
							\C^{ \delta }( [0,T], \left\| \cdot \right\|_{  H_{ \gamma }  } )
						}
					)
				}
			\big)
		= 0
	$
	and that
	$\lim_{N \to \infty} d_N =0$
	it holds that
	\begin{equation}
	\begin{split}
		\sup_{ N \in \N } \,
		\Big[
		& \|
			X^{N,x_N}
		\|_{ 
			L^p( 
				\P ;
				\left\| \cdot \right\|_{
				\C^{ \delta }( [0,T],\left\| \cdot \right\|_{  H_{ \gamma }  } )
				}
			)
		}
		\\ & 
		+ (d_N)^{ - ( \vartheta - \delta - 1 / p - \eps ) }
		\|
			X^{0,x_0} - X^{N,x_N}
		\|_{ 
			L^p( 
				\P ;
				\left\| \cdot \right\|_{
				\C^{ \delta }( [0,T], \left\| \cdot \right\|_{  H_{ \gamma }  } )
				}\texttt{}
			)
		}
		\Big]
		< \infty.
	\end{split}
	\end{equation}
	Hence we have for all 
	$ \delta \in [ 0, \vartheta - \nicefrac{ 1 }{ p } ) $ that
	\begin{equation}
		\begin{split}
				\limsup_{N \to \infty} \Big(
					\|
						X^{0,x_0} - X^{N,x_N}
					\|_{ 
						L^p( 
							\P ;
							\left\| \cdot \right\|_{
							\C^{ \delta }( [0,T], \left\| \cdot \right\|_{  H_{ \gamma }  } )
							}\texttt{}
						)
					}
				\Big)
			= 0
		\end{split}
	\end{equation}
	and this completes the proof of
	Corollary~\ref{cor:convGlobLip}.\\
	\end{proof}
	The next theorem 
	complements
	Corollary 4.5 in Cox et al.\@
	\cite{CoxHutzenthalerJentzenVanNeervenWelti2020}
	with 
	uniform convergence in probability on bounded subsets.
	It will be needed in the proof of Theorem \ref{thm: existence} below.
	The proof follows closely 
	Corollary 4.5 in Cox et al.\@ \cite{CoxHutzenthalerJentzenVanNeervenWelti2020}.
	\begin{theorem}
		[Convergence 
			and boundedness of
			Galerkin projections 
			in probability]
	\label{t: uniform convergence}
		Assume the setting in Section~\ref{ssec:setting}, 
		let $R \in (0, \infty)$,
		let
		$ 
			\vartheta \in (0,\infty) 
		$,
		let $ \delta \in [ 0, \vartheta ) $,
		and assume for all $ \mathbb{H}_{ \gamma }$-bounded 
		sets $ E \subseteq H_{ \gamma } $ that
		\begin{equation}
			\sup_{ N \in \N } \sup_{ v \in E }
			\left[
				\frac{
				\|
					B( v ) \mathscr{P}_N 
				\|_{ HS( \U, \H_{ \gamma - \beta } ) }
				}{
					1 + \| v \|_{ H_{ \gamma } }
				}
			\right]
			< \infty
		\end{equation}
		and that
		\begin{equation}
		\label{eq:convLocLip_assumption}
			\limsup_{ N \to \infty }
			\sup_{ v \in E }
			\left[
			\frac{
			\|
				B( v ) ( \mathscr{P}_0 - \mathscr{P}_N )
			\|_{ HS( \U, \H_{ \gamma - \chi } ) }
			}{
				1 + \| v \|_{ H_{ \gamma } }
			}
			\right]
			=0
			.
		\end{equation}
		Then it holds that
		\begin{equation}
			\limsup_{N \to \infty}
				\sup_{\substack{x \in H_{\gamma + \vartheta} \\  \|x\|_{H_{\gamma + \vartheta}} \leq R}}
						\E \Big[
							\|X^{N, x} - X^{0,x}\|_{ 
								\C^{ \delta }( [0,T], \left\| \cdot \right\|_{  H_{ \gamma }  } )
							}
							\wedge 1
						\Big] 
			= 0
		\end{equation}
		and that
		\begin{equation}
			\limsup_{N \to \infty}
				\sup_{\substack{x \in H_{\gamma + \vartheta} \\  \|x\|_{H_{\gamma + \vartheta}} \leq R}}
					\P(
						\| X^{0,x} \|_{
							\C^{ 0 }( [0,T], \left\| \cdot \right\|_{  H_{ \gamma }  } )
						}
						\geq N
					) 
			= 0.
		\end{equation}
	\end{theorem}
	\begin{proof}
		First note that the fact that for all $r\in \R$ and all $s\in (-\infty, r]$
		there exists a $C \in (0,\infty)$ such that for all 
		$u \in H_{s}$ it holds that
		$\|u\|_{H_s} \leq C \|u\|_{H_r}$
		ensures that it is enough to prove the assertion for
		$\vartheta \in \big( 0, \min\{ 1 - \alpha , \nicefrac{ 1 }{ 2 } - \beta \} \big)$.
		Therefore assume for the rest of the proof that
		$\vartheta \in \big( 0, \min\{ 1 - \alpha , \nicefrac{ 1 }{ 2 } - \beta \} \big)$.
		Next
		denote for all
		$ r \in \R $ and all
		$ M \in (0,\infty) $
		by
		$
			\phi_{ r, M } \colon H_r \to H_r
		$
		the function satisfying
		for all $ r \in \R $,
		$ M \in (0,\infty) $,
		and all $ v \in H_r $ that
		\begin{equation}
				\phi_{ r, M }( v )
			=
				v \cdot
				\min\biggl\{ 
					1, \frac{M + 1}{1 + \| v \|_{ H_r }}
				\biggr\},
		\end{equation}
		by 
		$ F_M \colon H_{ \gamma } \to H_{ \gamma - \alpha } $
		and 
		$ B_M \colon H_{ \gamma } \to HS( \U, \H_{ \gamma - \beta } ) $
		the functions satisfying
		that
		$
			F_M = F \circ \phi_{ \gamma, M } 
		$
		and that
		$
			B_M = B \circ \phi_{ \gamma, M } 
		$
		and by 
		$
			S_M \subseteq H_{ \gamma }
		$
		the set satisfying 
		that
		$
			S_M = \{ 
				v \in H_{ \gamma } \colon 
				\| v \|_{ H_{ \gamma } } \leq M + 1
			\}
		$.
		Note that for all $v, w \in H_\gamma$ and all $M \in \N$ it holds that
		\begin{align}
					&\left\|
						\phi_{ \gamma, M }( v )
						-
						\phi_{ \gamma, M }( w )
					\right\|_{ H_{ \gamma } }
				\\ \nonumber  ={}
					&\left\|
						\frac{
							v 
							\,
							( 1 + \| w \|_{ H_{ \gamma } } )
							\min\{ 1 + \| v \|_{ H_{ \gamma } } , M + 1 \}
							-
							w 
							\,
							( 1 + \| v \|_{ H_{ \gamma } } )
							\min\{ 1 + \| w \|_{ H_{ \gamma } } , M + 1 \}
						}{
							( 1 + \| v \|_{ H_{ \gamma } } )
							\,
							( 1 + \| w \|_{ H_{ \gamma } } )
						}
					\right\|_{ H_{ \gamma } }
				\\ \nonumber \leq{}
					&\left\| v - w \right\|_{ H_{ \gamma } }
				\\ \nonumber &
					+
					\left\|
						\frac{
							w 
							\,
							\big[
								( 1 + \| w \|_{ H_{ \gamma } } )
								\min\{ 1 + \| v \|_{ H_{ \gamma } } , M + 1 \}
								-
								( 1 + \| v \|_{ H_{ \gamma } } )
								\min\{ 1 + \| w \|_{ H_{ \gamma } } , M + 1 \}
							\big]
						}{
							( 1 + \| v \|_{ H_{ \gamma } } )
							\,
							( 1 + \| w \|_{ H_{ \gamma } } )
						}
					\right\|_{ H_{ \gamma } }
				\\ \nonumber \leq{}
					&\left\| v - w \right\|_{ H_{ \gamma } }
				\\ \nonumber & 
					+
						\frac{
							\left|
								( 1 + \| w \|_{ H_{ \gamma } } )
								\min\{ 1 + \| v \|_{ H_{ \gamma } } , M + 1 \}
								-
								( 1 + \| v \|_{ H_{ \gamma } } )
								\min\{ 1 + \| w \|_{ H_{ \gamma } } , M + 1 \}
							\right|
						}{
							( 1 + \| v \|_{ H_{ \gamma } } )
						}
		\end{align}
		and this implies that for all $v, w \in H_\gamma$ and all $M \in \N$ it holds that
		\begin{align}
				&
				\left\|
					\phi_{ \gamma, M }( v )
					-
					\phi_{ \gamma, M }( w )
				\right\|_{ H_{ \gamma } }
			\\ \nonumber \leq{}
				&\left\| v - w \right\|_{ H_{ \gamma } }
				+
					\frac{
							\big| 
								\| w \|_{ H_{ \gamma } } 
								-
								\| v \|_{ H_{ \gamma } } 
							\big|
							\min\{ 1 + \| v \|_{ H_{ \gamma } } , M + 1 \}
					}{
						( 1 + \| v \|_{ H_{ \gamma } } )
					}
			\\ \nonumber & \quad
				+
					\frac{
							( 1 + \| v \|_{ H_{ \gamma } } )
						\left|
							\min\{ 1 + \| v \|_{ H_{ \gamma } } , M + 1 \}
							-
							\min\{ 1 + \| w \|_{ H_{ \gamma } } , M + 1 \}
						\right|
					}{
						( 1 + \| v \|_{ H_{ \gamma } } )
					}
			\\ \nonumber \leq{}
				&\left\| v - w \right\|_{ H_{ \gamma } }
				+
							\left| 
								\| w \|_{ H_{ \gamma } } 
								-
								\| v \|_{ H_{ \gamma } } 
							\right|
						+
						\left|
							\min\{ 1 + \| v \|_{ H_{ \gamma } } , M + 1 \}
							-
							\min\{ 1 + \| w \|_{ H_{ \gamma } } , M + 1 \}
						\right|
			\\ \nonumber \leq{}
				&3
				\left\| v - w \right\|_{ H_{ \gamma } }
				.
		\end{align}
		Thus it follows that for all $M \in \N$ it holds that
		$
			|
				\phi_{ \gamma, M }
			|_{
				\C^1( H_{ \gamma }, \left\| \cdot \right\|_{  H_{ \gamma }  } )
			}
			\leq 3
		$.
		Combining this with
		the fact that for all $ M \in \N$ it holds that
		$
			|
				(F|_{ S_M })
			|_{
				\C^1( S_M, \left\| \cdot \right\|_{  H_{ \gamma - \alpha }  } )
			}
			+
			|
				(B|_{ S_M })
			|_{
				\C^1( S_M, \left\| \cdot \right\|_{  HS( \U, \H_{ \gamma - \beta } )  } )
			}
			+
			|
				\phi_{ \gamma, M }
			|_{
				\C^1( H_{ \gamma }, \left\| \cdot \right\|_{  H_{ \gamma }  } )
			}
			< \infty
		$,
		and the fact that
		for all $M \in \N$ it holds that 
		$ 
			\phi_{ \gamma, M }( H_{ \gamma } )
			\subseteq
			S_M
		$
		show that
		for all $ M \in \N $, $ p \in [1,\infty) $
		it holds that
		\begin{equation}
		\label{eq: FM, BM bound}
			|
				F_M
			|_{
				\C^1( H_{ \gamma }, \left\| \cdot \right\|_{  H_{ \gamma - \alpha }  } )
			}
			+
			|
				B_M
			|_{
				\C^1( H_{ \gamma }, \left\| \cdot \right\|_{  HS( \U, \H_{ \gamma - \beta } )  } )
			}
			< \infty.
		\end{equation}
		Hence,
		for all $ N \in \N_0 $, $ M \in \N $,
		and all $x \in H_{\gamma}$
		there exists an
		$ ( \mathbb{F}_t )_{ t \in [0,T] } $-adapted 
		stochastic process
		$
			\mathscr{X}^{ N, M, x } \colon [0,T] \times \Omega \rightarrow H_{ \gamma }
		$
		with continuous sample paths
		such that for all 
		$t \in [0,T] $
		it holds $\P$-a.s.\@ that
		\begin{equation}
		\label{eq:GlobLipGalerkin_solution2}
		\begin{split}
				\mathscr{X}_t^{ N, M, x}
		 & = 
				e^{ t A } \pi^{H_\gamma}_{V_N} \phi_{\gamma + \vartheta,M}(x)
				+
				\int_0^t
					e^{ ( t - s ) A }
					\pi^{H_{\gamma-\alpha}}_{V_N} F_M( \mathscr{X}^{ N, M, x}_s )
				\ds
			\\ & \quad
				+
				\int_0^t
					e^{ ( t - s ) A }
					\pi^{H_{\gamma-\beta}}_{V_N} B_M( \mathscr{X}_s^{ N, M, x} ) \mathscr{P}_N
				\dWs
		\end{split}
		\end{equation}
		(cf.\@ e.g., Proposition~7.1.19 in Jentzen \cite{Jentzen2015}
		or Theorem 6.1 in van~Neerven, Veraar \& Weis \cite{VanNeervenVeraarWeis2008}).
		We now proof the assertion by contradiction. Therefore
		assume that
		\begin{equation}
				\limsup_{N \to \infty}
					\sup_{\substack{x \in H_{\gamma + \vartheta} \\  \|x\|_{H_{\gamma + \vartheta}} \leq R}}
						\left(
							\E [
								\|X^{N, x} - X^{0,x}\|_{ 
									\C^{ \delta }( [0,T], \left\| \cdot \right\|_{  H_{ \gamma }  } )
								}
								\wedge 1
							] 
							\vee \P [\| X^{0,x} \|_{\C^0([0,T], \| \cdot \|_{H_\gamma})} \geq N]
						\right)
			> 0.
		\end{equation}
		Then there exist a sequence $(x_N)_{N \in \N_0} \subseteq H_{\gamma + \vartheta}$
		such that $\sup_{N \in \N_0} \| x_N \|_{H_\gamma + \vartheta} \leq R$
		and that
		\begin{equation}
		\label{eq: contradiction assumption}
				\limsup_{N \to \infty} \left(
						\E [
							\|X^{N, x_N} - X^{0,x_N}\|_{ 
								\C^{ \delta }( [0,T], \left\| \cdot \right\|_{  H_{ \gamma }  } )
							}
							\wedge 1
						]  
						\vee \P [
							\| X^{0,x_N} \|_{
								\C^{ 0 }( [0,T], \left\| \cdot \right\|_{  H_{ \gamma }  } )
							} 
							\geq N
						]
					\right )
			> 0.
		\end{equation}
		Since $H_{\gamma + \vartheta}$ is compactly embedded in $H_{\gamma}$
		we can assume that 
		$\lim_{N \to \infty} \|x_0 - x_N \|_{H_{\gamma}} =0$
		(otherwise take a subsequence and choose $x_0$ accordingly).
		Next denote for all
		$N \in \N_0$ and all $M \in \N$ 
		by $\tau_{N,M} \colon \Omega \to [0,T]$ the 
		function
		satisfying  
		that
		\begin{equation}
		\label{eq: def of tau}
				\tau_{N,M} 
			= 
				\left \{ 
					T \, \1_{ R \leq M}
					\wedge 
					\inf \{t \in [0,T] \colon 
						\| \mathscr{X}_t^{N,M,x_N} \|_{H_\gamma} 
						\geq M 
					\} 
				\right \}.
		\end{equation}
		Moreover, note that \eqref{eq: def of tau}
		and the fact that $\sup_{N \in \N} \|x_N\|_{\gamma + \vartheta} \leq R$
		ensure that for all  
		$N \in \N_0$, $M \in \N \cap (R, \infty)$,
		and all $t \in [0,T]$
		it holds $\P$-a.s.\@ that
		\begin{align}
					&\biggl(
						\mathscr{X}^{ N, M, x_N }_t
						-
						e^{ t A } \pi^{H_\gamma}_{V_N} x_N
					\biggr)
					\mathbbm{1}_{
						\{ t \leq \tau_{ N, M } \}
					} 
				=
					\biggl(
						\mathscr{X}^{ N, M, x_N }_t
							-
							e^{ t A } \pi^{H_\gamma}_{V_N} \phi_{\gamma + \vartheta,M}(x_N)
					\biggr)
						\mathbbm{1}_{
							\{ t \leq \tau_{ N, M } \}
						}
					\\ 
				\nonumber
				={}
					&\biggl(
						\int_0^t 
							e^{ ( t - s ) A }
							\pi^{H_{\gamma-\alpha}}_{V_N}
							F_M( \mathscr{X}^{ N, M, x_N}_s )
						\ds
				\\ \nonumber & \quad
						+
						\int_0^t 
							e^{ ( t - s ) A }
							\pi^{H_{\gamma-\beta}}_{V_N}
							B_M( \mathscr{X}^{ N, M, x_N}_s )
						\mathscr{P}_N
						\dWs
					\biggr)
					\mathbbm{1}_{
						\{ t \leq \tau_{ N, M } \}
					} \\
				\nonumber
				={}
					&\biggl(
						\int_0^t 
							\mathbbm{1}_{
								\{ s < \tau_{ N, M } \}
							}
							\,
							e^{ ( t - s ) A }
							\pi^{H_{\gamma-\alpha}}_{V_N}
							F_M( \mathscr{X}^{ N, M, x_N}_s ) 
						\ds
				\\ \nonumber & \quad
						+
						\int_0^t 
							\mathbbm{1}_{
								\{ s < \tau_{ N, M} \}
							}
							e^{ ( t - s ) A }
							\pi^{H_{\gamma-\beta}}_{V_N} 
							B_M( \mathscr{X}^{ N, M, x_N}_s )
						\mathscr{P}_N
						\dWs
					\biggr)
					\mathbbm{1}_{
						\{ t \leq \tau_{ N, M } \}
					} \\ 
				\nonumber  
				={} 
					&\biggl(
						\int_0^t 
							\mathbbm{1}_{
								\{ s < \tau_{ N, M } \}
							}
							\,
							e^{ ( t - s ) A }
							\pi^{H_{\gamma-\alpha}}_{V_N} 
							F( \mathscr{X}^{ N, M, x_N}_s )
						\ds
				\\ \nonumber & \quad 
						+    
						\int_0^t 
							\mathbbm{1}_{
								\{ s < \tau_{ N, M } \}
							}
							\,
							e^{ ( t - s ) A }
							\pi^{H_{\gamma-\beta}}_{V_N} 
							B( \mathscr{X}^{ N, M ,x_N}_s )
							\mathscr{P}_N
						\dWs
					\biggr)
					\mathbbm{1}_{
						\{ t \leq \tau_{ N, M } \}
					}.
		\end{align}
		Thus, we get for all 
		$ N \in \N_0 $ 
		and all $ M \in \N$ that
		\begin{equation}
				\P\Bigl(
					\forall \, t \in [0,T] \colon 
					\mathbbm{1}_{
						\{ R \leq M \}
					}
					\mathscr{X}^{ N, M, x_N}_{ \min\{ t, \tau_{ N, M } \} }
				=
					\mathbbm{1}_{
						\{ R \leq M \}
					}
					X^{N,x_N}_{ \min\{ t, \tau_{ N, M } \} }
				\Bigr)
			= 1
		\end{equation}
		(cf.\@ e.g., Proposition~7.1.10 in Jentzen \cite{Jentzen2015}
		or Lemma 8.2 in van~Neerven, Veraar \& Weis~\cite{VanNeervenVeraarWeis2008}).
		This ensures that for all
		$M\in \N \cap (R, \infty)$ and all $N \in \N_0$
		it holds that
		\begin{equation}
			\P(
				\forall t \in [0, T] \colon
						\mathscr{X}^{ N, M, x_N}_{ \min\{ t, \tau_{ N, M } \} }
					=	
						X^{N,x_N}_{ \min\{ t, \tau_{ N, M } \} }) = 1.
		\end{equation}
		Combining this with \eqref{eq: def of tau} shows that
		for all $M \in \N \cap (R, \infty)$,
		$m \in (0,M]$,
		and all $N \in \N_0$ it holds $\P$-a.s.\@ that
		\begin{equation} 
				\inf \big(
					\{t \in [0,T] \colon \| X_t^{N,x_N} \|_{H_\gamma} \geq m \} \cup \{T\} 
				\big)
			=
				\inf \big(
					\{t \in [0,T] \colon \| \mathscr{X}_t^{N,M,x_N} \|_{H_\gamma} \geq m \} \cup \{T\} 
				\big).
		\end{equation}
		Therefore, we obtain for all
		$ N \in \N_0 $ and $M \in \N \cap (R, \infty)$ that
		\begin{align}
		\nonumber
					&\P \left( 
						\sup_{t \in [0,T]} \| \mathscr{X}_t^{N,M,x_N} \|_{H_\gamma} < M 
					\right) \\ \nonumber
				\geq{}
					&\P \left(
								\sup_{t \in [0,T]} (
									\| \mathscr{X}_t^{0,M,x_0} \|_{H_\gamma} 
								)
						+
							\sup_{t \in [0,T]} (
								\| \mathscr{X}_t^{N,M,x_N} - \mathscr{X}_t^{0,M,x_0}\|_{H_\gamma}
							)
						< M
					\right) 
		\\ \label{eq: H_N bounded in Probability}
				\geq{}
					&\P \left(
								\sup_{t \in [0,T]} (
									\| \mathscr{X}_t^{0,M,x_0} \|_{H_\gamma} 
								)
						+
							\| \mathscr{X}^{N,M,x_N} - \mathscr{X}^{0,M,x_0}\|_{
								\C^{\delta}([0,T], \| \cdot \|_{H_\gamma})}
						< M
					\right) \\ \nonumber
				\geq{}
					&\P \left(
						\sup_{t \in [0,T]} \| \mathscr{X}_t^{0,M,x_0} \|_{H_\gamma} 
						\leq M -1, 
								\| \mathscr{X}^{N,M,x_N} - \mathscr{X}^{0,M,x_0}\|_{
									\C^{\delta}([0,T], \| \cdot \|_{H_\gamma})} 
							< 1 
					\right) \\ \nonumber
				={}
					&\P \left(
						\sup_{t \in [0,T]} \| X_t^{0,x_0} \|_{H_\gamma} 
						\leq M -1, 
								\| \mathscr{X}^{N,M,x_N} - \mathscr{X}^{0,M,x_0}\|_{
									\C^{\delta}([0,T], \| \cdot \|_{H_\gamma})} 
							< 1
					\right).
		\end{align}
		Analogously it follows for all $N \in \N_0$ and $M \in \N \cap (R, \infty)$ that
		\begin{equation}
		\label{eq: H_0 bounded in Probability}
			\begin{split}
					&\P \left( 
						\sup_{t \in [0,T]} \| \mathscr{X}_t^{0,M,x_N} \|_{H_\gamma} < M 
					\right) \\
				\geq{}
					&\P \left(
						\sup_{t \in [0,T]} \| X_t^{0,x_0} \|_{H_\gamma} 
						\leq M -1, 
								\| \mathscr{X}^{0,M,x_N} - \mathscr{X}^{0,M,x_0}\|_{
									\C^{\delta}([0,T], \| \cdot \|_{H_\gamma})} 
							< 1
					\right).
			\end{split}
		\end{equation}
		Furthermore, the fact that for all 
		$x \in H_\gamma$ and all $M \in \N$ it holds that
		$\|\phi_{\gamma, M} (x)\|_{H_\gamma} \leq M+1$
		together with the assumptions
		ensures that for all $M \in \N$
		it holds that
		\begin{equation}
			\begin{split}
					&\limsup_{N \to \infty} \left(
						\sup_{ v \in H_{ \gamma } }
							\frac
								{ 
									\|
										B_M( v ) ( \mathscr{P}_0 - \mathscr{P}_N ) 
									\|_{HS( \U, \H_{ \gamma - \chi } )}
								}
								{
									1 + \| v \|_{ H_{ \gamma } }
								}
					\right) \\
				\leq{}
					&\limsup_{N \to \infty} \left(
						\sup_{ v \in S_{M} }
							\frac
								{ 
									\|
										B( v ) ( \mathscr{P}_0 - \mathscr{P}_N ) 
									\|_{HS( \U, \H_{ \gamma - \chi } )}
								}
								{
									1 + \| v \|_{ H_{ \gamma } }
								}
					\right)
				=0
			\end{split}
		\end{equation}
		and that 
		\begin{equation}
				\sup_{N \in \N} \left(
					\sup_{ v \in H_{ \gamma } }
						\frac
							{ 
								\|
									B_M( v ) \mathscr{P}_N 
								\|_{HS( \U, \H_{ \gamma - \beta } )}
							}
							{
								1 + \| v \|_{ H_{ \gamma } }
							}
				\right)
			\leq
				\sup_{N \in \N} \left(
					\sup_{ v \in S_{M} }
						\frac
							{ 
								\|
									B( v ) \mathscr{P}_N 
								\|_{HS( \U, \H_{ \gamma - \beta} )}
							}
							{
								1 + \| v \|_{ H_{ \gamma } }
							}
				\right)
			<\infty.
		\end{equation}
		This, \eqref{eq: FM, BM bound}, 
		and Corollary \ref{cor:convGlobLip} then imply
		for all $M \in \N$ 
		and all
		$p \in (\nicefrac 1 \vartheta, \infty)$ that
		\begin{equation}
		\label{eq: H NMN C delta convergence}
			\limsup_{ N \to \infty } 
			\|
				\mathscr{X}^{0,M,x_0} - \mathscr{X}^{N,M,x_N}
			\|_{ 
				L^p( 
					\P ;
					\left\| \cdot \right\|_{
					\C^{ \delta }( [0,T], \left\| \cdot \right\|_{  H_{ \gamma }  } )
					}
				)
			}
			=0
		\end{equation}
		and that
		\begin{equation}
		\label{eq: H 0MN C delta convergence}
			\limsup_{ N \to \infty } 
			\|
				\mathscr{X}^{0,M,x_0} - \mathscr{X}^{0,M,x_N}
			\|_{ 
				L^p( 
					\P ;
					\left\| \cdot \right\|_{
					\C^{ \delta }( [0,T], \left\| \cdot \right\|_{  H_{ \gamma }  } )
					}
				)
			}
			=0.
		\end{equation}
		Thus we get for all $M \in \N$ that
		\begin{equation}
		\label{eq: H_N H_0 difference small in probability}
				\lim_{N \to \infty} \left(
					\P \left(
							\| \mathscr{X}^{N,M,x_N} - \mathscr{X}^{0,M,x_0}\|_{
								\C^{\delta}([0,T], \| \cdot \|_{H_\gamma})
							} 
							< 1
						\right)
				\right)
			= 1
		\end{equation}
		and that
		\begin{equation}
		\label{eq: H_0 difference small in probability}
				\lim_{N \to \infty} \left(
					\P \left(
							\| \mathscr{X}^{0,M,x_N} - \mathscr{X}^{0,M,x_0}\|_{
								\C^{\delta}([0,T], \| \cdot \|_{H_\gamma})
							} 
							< 1
						\right)
				\right)
			= 1.
		\end{equation}
		Combining now \eqref{eq: def of tau},
		\eqref{eq: H_N bounded in Probability},
		\eqref{eq: H_N H_0 difference small in probability},
		and the fact that
		$
			\forall \, \omega \in \Omega
			\colon
			\sup_{ t \in [0,T] }
			\| X_t^{0, x_0}( \omega ) \|_{ H_{ \gamma } }
			< \infty
		$
		yields that
		\begin{align}
		\nonumber
					&\liminf_{M \to \infty} \liminf_{N \to \infty} \,
						\P(\{ \tau_{N,M} = T \})
				\geq
					\liminf_{M \to \infty} \liminf_{N \to \infty} \,
						\P \left( 
							\sup_{t \in [0, T]}\| \mathscr{X}_t^{N,M,x_N} \|_{H_\gamma} < M 
						\right) \\
		\label{eq: tau is T in Probability}
				\geq{}
					&\liminf_{M \to \infty} \liminf_{N \to \infty} \, 
						\P \left(
							\sup_{t \in [0, T]} \| X_t^{0,x_0} \|_{H_\gamma} \leq M -1, 
								\| \mathscr{X}^{N,M,x_N} - \mathscr{X}^{0,M,x_0}\|_{
									\C^{\delta}([0,T], \| \cdot \|_{H_\gamma})
								} 
								< 1
						\right) \\ \nonumber
				={}
					&\lim_{M \to \infty}
						\P \left(
							\sup_{t \in [0, T]} \| X_t^{0,x_0} \|_{H_\gamma} \leq M -1
						\right)
				= 1.
		\end{align}
		Moreover, we derive from
		\eqref{eq: H_0 bounded in Probability},
		\eqref{eq: H_0 difference small in probability},
		and from the fact that for all $\omega \in \Omega$
		it holds that
		$
			\sup_{ t \in [0,T] }
			\| X_t^{0, x_0}( \omega ) \|_{ H_{ \gamma } }
			< \infty
		$
		that
		\begin{align}
					&\liminf_{N \to \infty}  \,
						\P \left( 
							\| X^{0,x_N} \|_{\C^{0}([0,T], \| \cdot \|_{H_\gamma}} < N
						\right) 
				\geq
					\liminf_{M \to \infty} \liminf_{N \to \infty} \,
						\P \left( 
							\| X^{0,x_N} \|_{\C^{0}([0,T], \| \cdot \|_{H_\gamma}} \leq M 
						\right) \\ \nonumber
				={}
					&\liminf_{M \to \infty} \liminf_{N \to \infty} \,
						\P \left( 
							\sup_{t \in [0,T]}
								\| \mathscr{X}_t^{0,M,x_N} \|_{H_\gamma}
							\leq M 
						\right) \\ \nonumber
				\geq{}
					&\liminf_{M \to \infty} \liminf_{N \to \infty} \, 
						\P \left(
							\sup_{t \in [0, T]} \| X_t^{0,x_0} \|_{H_\gamma} \leq M -1, 
									\| \mathscr{X}_t^{0,M,x_N} - \mathscr{X}_t^{0,M,x_0}\|_{
								\C^{\delta}([0,T], \| \cdot \|_{H_\gamma})
							} 
								< 1
						\right) \\ \nonumber
				={}
					&\liminf_{M \to \infty} \,
						\P \left(
							\sup_{t \in [0, T]} \| X_t^{0,x_0} \|_{H_\gamma} \leq M -1
						\right)
				= 1
		\end{align}
		and therefore it holds that
		\begin{equation}
		\label{eq: X bounded on Probability}
			\begin{split}
					&\limsup_{N \to \infty} 
						\P \left( 
							\sup_{t \in [0, T]}\| X_t^{0,x_N} \|_{H_\gamma} \geq N
						\right)	
				={}
					1
					-\liminf_{N \to \infty} 
						\P \left( 
							\sup_{t \in [0, T]}\| X_t^{0,x_N} \|_{H_\gamma} < N
						\right) 
				= 0.
			\end{split}
		\end{equation}
		Finally we obtain from \eqref{eq: H NMN C delta convergence},
		\eqref{eq: H 0MN C delta convergence},
		\eqref{eq: tau is T in Probability}, 
		and from \eqref{eq: X bounded on Probability}
		that
		\begin{align}
		\nonumber
					&\limsup_{N \to \infty}  \,
						\E [
							\|X^{N, x_N} - X^{0,x_N}\|_{ 
								\C^{ \delta }( [0,T], \left\| \cdot \right\|_{  H_{ \gamma }  } )
							}
							\wedge 1
						]\\	\nonumber	 	
				\leq{}
					&\limsup_{M \to \infty} \limsup_{N \to \infty} \bigg(
						\E[ 
							\|X^{N, x_N} - \mathscr{X}^{N,M,x_N}\|_{ 
								\C^{ \delta }( [0,T], \left\| \cdot \right\|_{  H_{ \gamma }  } )
							}
							\wedge 1
						]
				\\ \nonumber & \qquad\qquad
						+\E[ 
							\|\mathscr{X}^{N, M, x_N} - \mathscr{X}^{0,M,x_0}\|_{ 
									\C^{ \delta }( [0,T], \left\| \cdot \right\|_{  H_{ \gamma }  } )
							} 
							\wedge 1
						]
					\\ & \qquad \qquad
						+\E[
							\|\mathscr{X}^{0, M, x_0} - \mathscr{X}^{0,M,x_N}\|_{ 
								\C^{ \delta }( [0,T], \left\| \cdot \right\|_{  H_{ \gamma }  } )
							} 
							\wedge 1
						]
					\\ \nonumber & \qquad \qquad
						+\E[
							\|\mathscr{X}^{0, M, x_N} - X^{0,x_N}\|_{ 
								\C^{ \delta }( [0,T], \left\| \cdot \right\|_{  H_{ \gamma }  } )
							}
							\wedge 1
						]
					\bigg) \\ \nonumber
				\leq{}
					&\limsup_{M \to \infty} \limsup_{N \to \infty} \left(
						\P(\{\tau_{N,M} < T\}) 
						+ \P(\sup_{t \in [0,T]} \|X_t^{0,x_N}\|_{H_\gamma} \geq M)
					\right) 
				= 0,
		\end{align}
	which together with \eqref{eq: X bounded on Probability} 
	is a contradiction to \eqref{eq: contradiction assumption}.
		This finishes the proof of Theorem \ref{t: uniform convergence}.
	\end{proof}
\section{Existence of viscosity solutions of Kolmogorov equations}
\label{sec: Existence of viscosity solutions}
In this section we prove the main theorem of this chapter
(Theorem \ref{thm: existence}).  This 
shows the existence of viscosity
solutions of Kolmogorov equations of SPDEs and 
establishes a
representation of the viscosity solutions
using the corresponding solutions of the SPDEs.
Theorem \ref{thm: existence} generalizes
Theorem 4.16 in Hairer, Hutzenthaler \& Jentzen
\cite{HairerHutzenthalerJentzen2015}
to infinite-dimensional Hilbert spaces.
\begin{theorem}
	[Existence of viscosity solutions of Kolmogorov equations]
\label{thm: existence}
	Assume the setting in Section~\ref{ssec:setting}, assume that
	$ 
		\gamma =0
	$, 
	and assume 
	for all $\eps \in (0,\infty)$ and all $x \in H$ that
	\begin{equation}
	\label{eq: cont assumption}
		\limsup_{y \to x} 	\sup_{r \in [0,T]}
			\P(\|X^{0,x}_r - X^{0,y}_r\|_{H} \geq \eps) =0,
	\end{equation}
	let $\alpha_1, \beta_1, \theta \in (0, \infty)$,
	$\varz \in [\nicefrac 12, \infty)$,
	let $\varphi \in \C_{\mathbb{H}}(H,\R)$ 
	be bounded on $\mathbb{H}$-bounded sets,
	let
	$u_0 \colon [0,T] \times H \to \R$ be the function
	satisfying for all $(t,x) \in [0,T] \times H$ that
	\begin{equation}
		u_{0}(t,x) = \E[\varphi(X^{0,x}_t)],
	\end{equation}
	let
	$ 
		V \in 
			\C_{\mathbb{H}}^2( \mathbb{H}, (1,\infty))
	$
	satisfy
	that $D_{\mathbb{H}}^2 \, V$ 
	is uniformly continuous 
	with respect to the 
	$\| \cdot \|_{H}$ and the
	$\| \cdot \|_{L(\mathbb{H}, \mathbb{H}')}$-norm
	on $\mathbb{H}$-bounded subsets of $H$,
	assume for all $R \in (0, \infty)$ that
	\begin{equation}
	\label{eq: V bound for XN}
			\sup_{ N \in \N }	
				\sup_{ 
					\substack{
						(t,x) \in (0,T) \times H \\
						\| x \|_{H} \leq R
					}
				}
					\E[V(X_t^{N,x})]
		< \infty,
	\end{equation}
	assume that 
	$
      \lim_{
				r \to \infty
      }
			\inf \{
				V(x) \colon x \in H,~ \|x\|_{H} \geq r	
			\}
    = \infty, 
$
	and
	$
      \lim_{
				r \to \infty
      }
			\sup \big \{
				\frac
					{| \varphi(x) |}
					{V(x)}
				\colon x \in H,~ \|x\|_H \geq r
			\big \}
    = 0, 
	$
	assume that
	\begin{equation}
	\label{eq: F and B continuity}
		\begin{split}
			&F|_{H_{\varz + 1/2 - \alpha_1}} \in
				\C_{
					\mathbb{H}_{\varz + 1/2 - \alpha_1},
					\mathbb{H}_{\varz - 1/2}
				}(
					H_{\varz + 1/2 - \alpha_1}, 
					H_{\varz - 1/2}
				)
		\quad \textrm{ and that } \\
			&B|_{H_{\varz +1/2 -\beta_1}} \in
				\C_{
					\mathbb{H}_{\varz + 1/2 - \beta_1},
					\mathbb{HS}(\mathbb{U}, \mathbb{H}_{\varz})
				}(
					H_{\varz + 1/2 - \beta_1},
					HS(\mathbb{U}, \mathbb{H}_{\varz})
				),
		\end{split}
	\end{equation}
	assume
	that for all $x \in H_{2\varz}$ there exist
	$r_x\in (0, \infty)$, $R_x \in (0, \infty)$
	such that for all $\xi \in H_{\varz}$
	with $\|x- \xi \|_{H_{\varz}} \leq r_x$ it holds that
	\begin{equation}
	\label{eq: F H 1/2 bound}
			\|F(\xi)\|_{H_{\varz - 1/2}} 
		\leq
			R_x (\|\xi\|_{H_{\varz + 1/2-\alpha_1}} + 1)
	\end{equation}
	and that
	\begin{equation}
	\label{eq: B H 1 bound}
			\|B(\xi)\|^2_{HS(\mathbb{U}, \mathbb{H}_{\varz})} 
		\leq
			R_x (\|\xi\|^2_{H_{\varz + 1/2 -\beta_1}} + 1), 
	\end{equation}
	assume that $\mathscr{P}_0 = \id_U$, that
	$\mathscr{P}_N$, $N \in \N$, are finite-dimensional projections,
	and assume that for all $\mathbb{H}$-bounded sets 
	$ E \subseteq H $ it holds that
	\begin{equation}
	\label{eq: B C1 bound thm}
			\sup_{ N \in \N } \sup_{ v \in E }
			\left[
				\frac{
				\|
					B( v ) \mathscr{P}_N 
				\|_{ HS( \mathbb{U}, \mathbb{H}_{- \beta } ) }
				}{
					1 + \| v \|_{ H }
				}
			\right]
			< \infty
	\end{equation}
	and that
	\begin{equation}
	\label{eq:convLocLip_assumption thm}
			\limsup_{ N \to \infty }
			\sup_{ v \in E }
			\left[
			\frac{
			\|
				B( v ) ( \mathscr{P}_0 - \mathscr{P}_N )
			\|_{ HS( \mathbb{U}, \mathbb{H}_{ - \chi } ) }
			}{
				1 + \| v \|_{ H }
			}
			\right]
		=0.
	\end{equation}
	Then it holds that 
	$u_0$ is continuous with respect to the $\| \cdot \|_{\R \times H}$-norm,
	that
	$u_0$ is bounded on $\R \times \mathbb{H}$-bounded 
	subsets of $(0,T) \times H$,
	that for all $\tilde{R} \in (0, \infty)$ it holds that
	\begin{align}
	\label{eq: uniformly bounded at 0 with V thm}
				\lim_{R \to \infty} \lim_{r \downarrow 0} \lim_{\eps \downarrow 0}	
				\sup \Bigg \{
					&u_0(t,x) -  \varphi(\hat{x})
					\colon
					~x, \hat{x} \in H_{2\varz},
					~t\in (0,\eps),
					~\tfrac{\|x\|_{H_{\varz}}^2}{V(x)} 
							\vee 
								\tfrac{
									\|\hat{x}\|_{H_{\varz}}^2}{V(\hat{x})}  
						\leq R, \\ \nonumber
					~&\|x \|_{H} \leq \tilde{R}, 
					~\| x - \hat{x}\|_{H} \leq r
				\Bigg \}
			= 
				0,
	\end{align}
	that
  \begin{equation}
	\label{eq: grow condition on u}
    \lim_{ r \to \infty }
      \sup_{ t \in [0,T]}
			\sup \left \{
				\frac
					{| u_0(t, x) |}
					{V(x)}
				\colon x \in H, ~\|x\|_H \geq r
			\right \}
    = 0,
  \end{equation}
  and that
  $ 
    u_0|_{ (0,T) \times H} 
  $
  is a viscosity solution of
  \begin{equation} 
    \tfrac{ \partial }{ \partial t }
    u(t,x) -
    \big\langle
      F(x) + A(x), (D_\mathbb{H} u)(t,x)
    \big\rangle_{H, H'}
    -\nicefrac 12
		\big \langle 
			B(x),
			I_{\mathbb{H}_{\varz}}^{-1}
				\big( (D^2_{\mathbb{H}} \, u)(t,x) \, B(x) \big)
    \big \rangle_{HS(\mathbb{U},\mathbb{H}_{\varz})} 
    = 0
  \end{equation}
	for $ (t,x) \in (0,T) \times H $ relative to 
	$
		(
			(0,T) \times H \ni (t,x) 
			\to \nicefrac 12 \|x \|^2_{H_{\varz}} \in [0, \infty],
			\R \times \mathbb{H}, 
			\R \times \mathbb{H}_{\varz}
		)
	$.
\end{theorem}
\begin{proof}
	We will show this theorem by an application of
	Lemma \ref{l:limits.of.viscosity.solutions}.
	In the first step we
	introduce some notation.
	Denote 
	by 
	$\phi \colon [0, \infty) \to [0,\infty)$
	the function satisfying for all $r \in [0,\infty)$ that
	\begin{equation}
	\label{eq: definition of phi}
			\phi(r) 
		= 
			\sqrt{r} 
			\cdot \inf 
				\big( \{ \sqrt{V(y)} \colon y \in H, ~ |\varphi(y)| =r\} \cup \{r\} \big),
	\end{equation}
	by
	$u_{N} \colon [0,T] \times H \to \R$, $N \in \N$,
	the functions satisfying for all
	$t \in [0, T]$,
	$x \in H$, and all $N \in \N$, that
	\begin{equation}
		u_{N}(t,x) = \E[\varphi(X^{N,x}_t)],
	\end{equation}
	by $W \subseteq (0,T) \times H$ the set satisfying that
	$W = (0,T) \times H_{2\varz}$,
	by  
	$
    F_N \colon 
			((0,T) \times H_{1}) 
			\times \R 
			\times (\R \times H)' 
			\times 
				\mathbb{S}_{
					\R \times \mathbb{H}_{\varz}, 
					(\R \times \mathbb{H}_{\varz})'
				}
    \to \R 
  $,
	$N \in \N_0$,
  the functions satisfying for all
	$
		((t, x), r, p, C, \tilde{C}) 
			\in 
				((0,T) \times H_{1}) 
				\times \R 
				\times (\R \times H)'
				\times 
					\mathbb{S}_{
						\R \times \mathbb{H}_{\varz}, 
						(\R \times \mathbb{H}_{\varz})'
					}
				\times 
					\mathbb{S}_{\mathbb{H}_{\varz}, \mathbb{H}_{\varz}'}
	$ 
	and all $N \in \N_0$
	with
	$
		\forall y \in H_{\varz} \colon
			\tilde{C} y = \pi^{(\R \times H_{\varz})'}_{2} C (0,y)
	$
	that
  \begin{equation}
	\label{eq: def of FN}
		\begin{split}
				F_N( (t, x), r, p, C)
			={}
				 &I^{-1}_{\R} \pi^{(\R \times H)'}_1 p 
					-\left\langle
						F(	\pi^{H}_{V_N} \, x) 
						+A(	\pi^{H}_{V_N} \, x), 
						\pi_{V_N'}^{H'} \pi^{(\R \times H)'}_2 p
					\right\rangle_{H, H' } \\
					&-\nicefrac 12 \,
					\Big \langle
						B(	\pi^{H}_{V_N} \, x) ,
						I_{\mathbb{H}_{\varz}}^{-1} \tilde{C} 
							\, B(	\pi^{H}_{V_N} \, x)  
					\Big \rangle_{HS(\mathbb{U}, \mathbb{H}_{\varz})},
		\end{split}
  \end{equation}
	by 
	$
    \tilde{F}_N \colon 
			((0,T) \times V_N) 
			\times \R 
			\times (\R \times V_N)' 
			\times \mathbb{S}_{\R \times \mathbb{V}_N, (\R \times \mathbb{V}_N)'}
    \to \R 
  $,
	$N \in \N$,
  the functions satisfying for all
	$
		((t, x), r, p, C, \tilde{C}) 
			\in 
				((0,T) \times V_N) 
				\times \R 
				\times (\R \times V_N)'
				\times \mathbb{S}_{\R \times \mathbb{V}_{N}, (\R \times \mathbb{V}_{N})'}
				\times \mathbb{S}_{\mathbb{V}_{N}, \mathbb{V}_{N}'}
	$ 
	and all $N \in \N$
	with 
	$
		\forall y \in V_{N} \colon
			\tilde{C} y = \pi^{(\R \times V_N)'}_{2} C (0,y)
	$
	that
  \begin{equation}
	\label{eq: def of tilde FN}
		\begin{split}
				\tilde{F}_N((t, x), r, p, C)
			={}
				 &I^{-1}_{\R} \pi^{(\R \times V_N)'}_1 p 
					-\left\langle
						\pi^{H}_{V_N} F(x) 
						+ \pi^{H}_{V_N} A(x), 
						\pi^{(\R \times V_N)'}_2 p
					\right\rangle_{V_N, V_N' } \\
					&- \nicefrac 12
					\Big \langle
						\pi^{H}_{V_N} B(x) ,
						I_{\mathbb{V}_{N}}^{-1} \tilde{C} 
							\, \pi^{H}_{V_N} B(x)  
					\Big \rangle_{HS(\mathbb{U}, \mathbb{V}_{N})},
		\end{split}
  \end{equation}
	by $h \colon (0,T) \times H \to [0, \infty]$
	the function satisfying for all
	$(t,x) \in (0,T) \times H$ that
	\begin{equation}
	\label{eq: def of h thm}
			h(t,x) 
		= 
			\tfrac 12 \|x\|^2_{H_{\varz}},
	\end{equation}
	let $d \in (0,\vartheta)$ satisfy that
	\begin{equation}
	\label{eq: d sufficient small}
		\alpha + d < 1, \qquad d + \chi < \nicefrac 12,
	\end{equation}
	let $\mathcal{I} \subseteq \N$,
	and let $(\tilde{e}_i)_{i \in \mathcal{I}}$ be an orthonormal basis of $\mathbb{U}$.
	Note that \eqref{eq: F and B continuity} assures that $\tilde{F}_N$, $N \in \N$
	is well defined.
	Moreover,
	it holds for all $j \in (0,\infty)$ that
	the set $\{x \in H \colon \|x\|_{H_d} \leq j\}$ is compact in $H$.
	Thus
	it follows from the continuity of $\varphi$ that for all $j \in (0,\infty)$
	it holds that
	$\varphi$ is uniformly continuous on the set
	$\{x \in H \colon \|x\|_{H_d} \leq j\}$ 
	with respect to the $\| \cdot \|_{H_d}$-norm
	and therefore
	it holds that for all $j, \delta \in (0, \infty)$ 
	there exists an $\eps_{j, \delta} \in (0, 1)$ such that
	for all $x, y \in H$ with 
	$\|x\|_{H_d} \vee \|y\|_{H_d} \leq j$ 
 	and with $\|x-y\|_{H} \leq \eps_{j, \delta}$ it holds that
	\begin{equation}
	\label{eq: phi Lipschitz}
		|\varphi(x) - \varphi(y)| \leq \delta.
	\end{equation}
	In addition, note that 
	\eqref{eq: d sufficient small}
	and
	Theorem \ref{t: uniform convergence} (with $\gamma \leftarrow d$,
		$\alpha \leftarrow \alpha +d$,
		$\beta \leftarrow \beta +d$, $\chi \leftarrow \chi +d$,
		$\vartheta \leftarrow \vartheta - d$,
		$F \leftarrow F|_{H_d}$, and with $B \leftarrow B|_{H_d}$)
	ensure
	for all $R \in (0, \infty)$ that
	\begin{equation}
		\limsup_{N \to \infty}
			\sup_{\substack{x \in H_{\vartheta} \\  \|x\|_{H_{\vartheta}} \leq R}}
					\E [
							\|X^{N, x} - X^{0,x}\|_{ 
								\C^{ 0 }( [0,T], \left\| \cdot \right\|_{  H_d  } )
							}
							\wedge 1
						] 
			= 0 
	\end{equation}
	and that
	\begin{equation}
		\limsup_{N \to \infty}
			\sup_{\substack{x \in H_{\vartheta} \\  \|x\|_{H_{\vartheta}} \leq R}}
				\P(
					\| X^{0,x} \|_{
						\C^{ 0 }( [0,T], \left\| \cdot \right\|_{  H_d  } )
					}
					\geq N
				) 
		= 0
	\end{equation}
	and therefore we get for all $R, \eps \in (0, \infty)$ that
	\begin{align}
	\label{eq: XN convergence in probability}	
		\limsup_{N \to \infty}
			\sup_{\substack{x \in H_\vartheta \\  \|x\|_{H_\vartheta} \leq R}}
				\sup_{t \in [0, T]}
	 				\P (\| X_t^{N,x} - X_t^{0,x} \|_{H_d} \geq \eps) 
		= 0
	\end{align}
	and that
	\begin{align}
	\label{eq: X bounded in probability}
		\limsup_{j \to \infty}
			\sup_{\substack{x \in H_\vartheta \\  \|x\|_{H_\vartheta} \leq R}}
				\sup_{t \in [0, T]}
						\P (\| X_t^{0,x} \|_{H_d} \geq j) 
		= 0.
	\end{align}
	Furthermore, note that 
	\eqref{eq: XN convergence in probability}, Fatou's lemma
	and \eqref{eq: V bound for XN}
	show
	that for all $R \in (0,\infty)$ it holds that 
	\begin{equation}
	\label{eq: V bound for X0}
		\begin{split}
			&\sup_{
				\substack{
					(t,x) \in (0,T) \times H \\
					\| x \|_{H} \leq R
				}
			}
				\E [ V(X_{t}^{0,x})] 
		=
			\sup_{
				\substack{
					(t,x) \in (0,T) \times H \\
					\| x \|_{H} \leq R
				}
			}
				\E [\, \underset{N \to \infty}{\P\textrm{-}\lim} \, \,  V(X_{t}^{N,x})] \\
		\leq{}
			&\sup_{
				\substack{
					(t,x) \in (0,T) \times H \\
					\| x \|_{H} \leq R
				}
			}
				\liminf_{N \to \infty} \E [V(X_{t}^{N,x})]
		< \infty.
		\end{split}
	\end{equation}
	Moreover,
	we get from
	$
			\lim_{r \to \infty} \sup
			\Big \{
				\frac
					{| \varphi(x) |}
					{ V(x)}
				\colon x \in H, \|x\|_H \geq r
			\Big \}
    = 0
	$
	and from the fact that $\varphi$ is bounded on $\mathbb{H}$-bounded subsets of $H$
	that for all $\delta \in (0, \infty)$ there exist a $C_\delta \in (0,\infty)$
	such that for all $x \in H$ it holds that
	\begin{align}
	\label{eq: V bound for phi}
		| \varphi(x) | \leq C_\delta + \delta V(x)
	\end{align}
	and this together with \eqref{eq: V bound for XN} 
	and with \eqref{eq: V bound for X0} ensures that $u_N$, $N \in \N_0$,
	are well-defined. 
	In the second step we will show \eqref{eq: uniform convergence of u_n}.
	We derive with 
	the fact that for all $j, \delta \in (0,\infty)$
	it holds that
	$\eps_{j,\delta} <1$,
	with \eqref{eq: phi Lipschitz}
	and with \eqref{eq: V bound for phi}
	that for all
	$N \in \N,$
	$j, \delta \in (0, \infty)$, 
	$t \in (0,T)$ and all
	$x \in H$ it holds that 
	\begin{align}
	\label{eq: u difference bound}
					&\left|
						u_{N} (t,x) - u_0(t,x) 
					\right|
			= 
						\left|
							\E[\varphi(X^{N,x}_t)] - \E[\varphi(X^{0,x}_t)] 
						\right| \\ \nonumber
			\leq{}	
						&
							\E \left [
								\big| \varphi(X^{N,x}_t) - \varphi(X^{0,x}_t) \big|
								\1_{
									\|X^{N,x}_t - X^{0,x}_t\|_{H_d} \leq \eps_{j, \delta}
								}
								\1_{\|X^{0,x}_t\|_{H_d} \leq j-1}
							\right ] \\ \nonumber
							&+\E \left[
								\big| \varphi(X^{N,x}_t) - \varphi(X^{0,x}_t) \big|
								\1_{\|X^{0,x}_t\|_{H_d} \geq j-1}
							\right ] 
							+\E \left [ \big| \varphi(X^{N,x}_t)- \varphi(X^{0,x}_t) \big|
								\1_{
									\|X^{N,x}_t - X^{0,x}_t\|_{H_d} \geq \eps_{j, \delta}
								}
							\right] \\ \nonumber
			\leq{}
						&\delta
							+\E \left[
								(2 C_\delta + \delta V(X^{N,x}_t) + \delta V(X^{0,x}_t))
								\1_{\|X^{0,x}_t\|_{H_d} \geq j-1}
							\right ]
				\\ \nonumber & 
							+\E \left [ 
								(2 C_\delta + \delta V(X^{N,x}_t) + \delta V(X^{0,x}_t))
								\1_{
									\|X^{N,x}_t - X^{0,x}_t\|_{H} \geq \eps_{j, \delta}
								}
							\right]\\ \nonumber
			\leq{}
						&\delta
							+2 C_\delta (
								\P(\|X^{0,x}_t\|_{H_d} \geq j-1)
								+ \P(\|X^{N,x}_t - X^{0,x}_t\|_{H} \geq \eps_{j, \delta})
							)
							+2\delta \, \E \left[
								V(X^{N,x}_t) + V(X^{0,x}_t)
							\right ].
	\end{align}
	Combining 
	\eqref{eq: def of h thm},
	the fact for all $r\in \R$ and all $s\in (-\infty, r]$
	there exists a $C \in (0,\infty)$ such that for all 
	$u \in H_{s}$ it holds that
	$\|u\|_{H_s} \leq C \|u\|_{H_r}$,
	\eqref{eq: V bound for XN}, 
	\eqref{eq: V bound for X0}, 
	\eqref{eq: u difference bound}, 
	\eqref{eq: X bounded in probability}, and
	\eqref{eq: XN convergence in probability} then shows for all $R, \tilde{R} \in (0,\infty)$
	that
	\begin{align}
	\nonumber
				&\limsup_{ N \to \infty }	
				\sup_{ 
					(t,x) \in 
						\{ 
							(s,y) \in (0,T) \times H \colon 
								\|(s,y)\|_{\R \times H} \leq \tilde{R},~
								h(s,y) \leq R
						\} 
				}
					\left|
						u_N (t,x) - u_0(t,x) 
					\right| \\ \nonumber
			\leq{}
				&\limsup_{ N \to \infty }	
				\sup_{ 
					(t,x) \in 
						\{ 
							(s,y) \in (0,T) \times H \colon 
								\|(s,y)\|_{\R \times H_\vartheta} \leq R 
						\} 
				}
					\left|
						u_N (t,x) - u_0(t,x) 
					\right| \\
	\label{eq: uniform convergence on bounded subsets}
			\leq{} 
				&\limsup_{\delta \to 0} \limsup_{j \to \infty}
				\limsup_{ N \to \infty }	
					\sup_{ 
						x \in \{ 
							y \in H \colon \|y\|_{H_\vartheta} \leq R
						\} 
					}
					\sup_{t \in [0,T]} 
					\Big(
						\delta
						+4 C_\delta (
							\P(\|X^{0,x}_t\|_{H_d} \geq j-1)
				\\ \nonumber & \qquad \qquad \qquad \qquad 
							+ \P(\|X^{N,x}_t - X^{0,x}_t\|_{H} \geq \eps_{j, \delta})
						) 
							+2\delta \, \E \left[
								V(X^{N,x}_t) + V(X^{0,x}_t)
							\right ] 
					\Big )\\ \nonumber
			={} &0.
	\end{align}
	and this proves \eqref{eq: uniform convergence of u_n}.
	Furthermore, it follows from 
	\eqref{eq: V bound for phi}
	and from \eqref{eq: V bound for XN}
	that
	for all $(t_0,x_0) \in (0,T) \times H$ it holds that
	\begin{align}
	\nonumber
				&\sup_{ N \in \N }	
				\sup_{ 
					\substack{
						(t,x) \in (0,T) \times H \\
						\| x \|_{H} \leq R
					}
				}
					\left|
						u_{N} (t,x) - u_N(t_0,x_0) 
					\right|\\ \nonumber
			= {}
				&\sup_{ N \in \N }	
				\sup_{ 
					\substack{
						(t,x) \in (0,T) \times H \\
						\| x \|_{H} \leq R
					}
				}
						\left|
							\E[\varphi(X^{N,x}_t)] - \E[\varphi(X^{N,x_0}_{t_0})] 
						\right| \\
			\leq{}
				&\sup_{ N \in \N }	
				\sup_{ 
					\substack{
						(t,x) \in (0,T) \times H \\
						\| x \|_{H} \leq R
					}
				}
							\E \left[
								(2C_1 + V(X^{N,x}_t) + V(X^{N,x_0}_{t_0}))
							\right ]\\ \nonumber
			\leq{}
				&\sup_{ N \in \N }	
				\sup_{ 
					\substack{
						(t,x) \in (0,T) \times H \\
						\| x \|_{H} \leq R
					}
				}
				\left(
						2 C_1 
							+\E \left[
								V(X^{N,x}_t) + V(X^{N,x_0}_{t_0})
							\right ] 
				\right)
				<
				\infty
	\end{align}
	and this verifies \eqref{eq: uniform bound of u_n}.
	Moreover it follows from \eqref{eq: cont assumption}
	and the continuity of $\varphi$ that
	for all $\eps \in (0,\infty)$ and all $(t,x) \in [0,T] \times H$ it holds that
	\begin{equation}
	\label{eq: convergence of varphi X}
		\limsup_{(s,y) \to (t,x)} \P(\| \varphi(X_t^{0,x}) -\varphi(X_s^{0,y}) \|_H > \eps) =0.
	\end{equation}
	In addition note that we get from \eqref{eq: definition of phi}
	the fact that $\varphi$ is bounded on $\mathbb{H}$-bounded sets
	and from the fact that
	$
      \lim_{
				r \to \infty
      }
			\sup \Big \{
				\frac
					{| \varphi(x) |}
					{V(x)}
				\colon x \in H,~ \|x\|_H \geq r
			\Big \}
    = 0 
	$
	that 
	\begin{equation}
	\label{eq: phi grow more than linear}
		\begin{split}
				\lim_{r \to \infty}
					\tfrac{\phi(r)}{r} 
			={} 
				&\lim_{r \to \infty} \big(
					\tfrac {1}{\sqrt{r}}
					\inf 
						\big( \{ \sqrt{V(y)} \colon y \in H, ~ |\varphi(y)| =r\} \cup \{r\} \big)
				\big) \\
			={}
				&\lim_{r \to \infty} \Big(
					\inf 
						\Big( 
							\big\{ 
								\sqrt{\tfrac{V(y)}{|\varphi(y)|}} 
								\colon y \in H, ~ |\varphi(y)| =r
							\big\} 
							\cup \{\sqrt{r}\} 
						\Big)
				\Big)
			=
				\infty
		\end{split}
	\end{equation}
	and that
	\begin{align}
			&\lim_{r \to \infty}
			\sup \Big \{
				\frac
					{\phi(|\varphi(x)|)}
					{V(x)}
				\colon x \in H,~ \|x\|_H \geq r
			\Big \} \\ \nonumber
			={} 
				&\lim_{r \to \infty}
					\sup \Big \{
						\tfrac{\sqrt{|\varphi(x)|}}{V(x)} \cdot
						\inf 
							\big(
							\{ \sqrt{V(y)} \colon y \in H, ~ |\varphi(y)| =|\varphi(x)|\} 
								\cup \{r\} 
							\big)
					\colon x \in H,~ \|x\|_H \geq r
				\Big \}\\ \nonumber
			\leq{} 
				&\lim_{r \to \infty}
				\sup \Big \{
					\tfrac{\sqrt{|\varphi(x)|}}{V(x)} \cdot \sqrt{V(x)}
					\colon x \in H,~ \|x\|_H \geq r
				\Big \}
			=
				0.
	\end{align}
	Thus the fact that $\varphi$ 
	is bounded on $\mathbb{H}$-bounded sets
	and \eqref{eq: definition of phi} show that there exists a $C \in (0,\infty)$
	such that for all $x \in H$ it holds that
	\begin{equation}
		\phi(|\varphi(x)|) \leq V(x) + C
	\end{equation}
	and therefore it follows from 
	\eqref{eq: V bound for X0}
	that there exist a $C \in (0,\infty)$
	such that for all $x \in H$ it
	holds that
	\begin{equation}
		\label{eq: X uniform E bound}
			\sup_{
				\substack{
					(s,y) \in (0,T) \times H \\
					\| y \|_{H} \leq \| x \|_{H} +1
				}
			} 
				\E[\phi(|\varphi(X_s^{0,y})|)]
		\leq
			\sup_{
				\substack{
					(s,y) \in (0,T) \times H \\
					\| y \|_{H} \leq \| x \|_{H} +1
				}
			} 
				\E[V(X_s^{0,y})+C]
			<
			 \infty.
	\end{equation}
	Combining \eqref{eq: X uniform E bound}
	and \eqref{eq: phi grow more than linear}
	ensures that for all $x \in H$ it holds that
	$
		\{ 
			\varphi(X_s^{0,y}) 
				\colon (s,y) \in [0,T] \times H, ~\|y\|_H \leq \|x\|_H+1
		\}
	$
	is uniformly integrable and combining this and
	\eqref{eq: convergence of varphi X} proves that
	\begin{equation}
			\lim_{(s,y) \to (t,x)}
				|u_0(t,x)-u_0(s,y)|
		\leq
			\lim_{(s,y) \to (t,x)}
				\E[|\varphi(X^{0,y}_s) - \varphi(X^{0,x}_t)| ]
		=0.
	\end{equation}
	This verifies that $u_0 \in \C_{\R \times \mathbb{H}}([0,T] \times H, \R)$.
	In the next step we show \eqref{eq:convergence F+}.
	Observe that
	similarly as in Theorem 4.16
	in Hairer, Hutzenthaler \& Jentzen
	\cite{HairerHutzenthalerJentzen2015} 
	where we can replace the assumption on V by \eqref{eq: V bound for XN}
	and the time set $[0,\infty)$ by $[0,T]$
	it follows that
	for all $N \in \N$ it holds that
	$u_N|_{(0,T) \times V_N}$ is a classical viscosity solution of
	\begin{equation}    
			\tilde{F}_N(
				(t,x), 
				u(t,x), 
				(D_{\R \times \mathbb{V}_N} \, u)(t,x)),
				(D^2_{\R \times \mathbb{V}_N} \, u)(t,x)
			)
		= 0,
	\end{equation}
	for $(t,x) \in (0,T) \times V_N$ and therefore is also continuous with
	respect to the $\|\cdot \|_{\R \times V_N}$-norm.
	Furthermore, it holds that
	$
		h|_{(0,T) \times H_{\varz}} 
			\in 
				\C^2_{\R \times \mathbb{H}_{\varz}}
					((0,T) \times H_{\varz},\R).
	$
	Moreover, we get for all $x,y\in H_{\varz}$ 
	that
	$
			\|x \|_{H_{\varz}}^2 
		= 
			\|y\|_{H_{\varz}}^2 
			+ 2 \langle y, x-y \rangle_{H_{\varz}} 
			+ \langle x-y, x-y \rangle_{H_{\varz}}
		\geq
			\|y|_{H_{\varz}}^2 
			+2\langle 
				I_{\mathbb{H}_{\varz}} (x_0), 
				x-y 
			\rangle_{H'_{\varz},H_{\varz}}
	$
	and thus we get 
	from Definition
	\ref{d:semijets}
	and from the fact that for all
	$x \in H_{\vartheta}$ and all $y\in H_{2\varz}$ it holds that
	$
			\langle I_{\mathbb{H}_{\varz}} (y) , x \rangle_{H_{\vartheta}',H_{\varz}}	
		=
			\langle y , x \rangle_{H_{\varz}}
		=
			\langle (-A)^{2\varz } y , x \rangle_H
		= 
			\langle I_{\mathbb{H}} ( (-A)^{2\varz } y), x \rangle_{H',H}
	$
	that 
	$
			\{ 
				(t, y) \in (0,T) \times H \colon 
					(J^{2}_{\mathbb{\R \times H}, -} h) (t,y) \neq \emptyset 
			\} 
		=
			\{ 
				(t,y) \in (0,T) \times H 
					\colon I_{\mathbb{H}_{\varz}} (y) 
						\in D(E_{\mathbb{H}_{\varz}', \mathbb{H}'})
				\}
		=
			(0,T) \times H_{2\varz}
		=
			W.
	$ 
	In addition, note that we have for all $(t,x) \in W$ that
	\begin{equation}
	\label{eq: h derivative existence}
		\begin{split}
				&E_{
					(\R \times \mathbb{H}_{\varz})',
					(\R \times  \mathbb{H})'
				} 
					\big((D_{\R \times \mathbb{H}_{\varz}} 
						(h|_{\R \times H_{\varz}}))(t,x) \big) 
			= 
				I_{\R \times \mathbb{H}} (0,(-A)^{2\varz}(x)), \\ 
				&(D^2_{\R \times \mathbb{H}_{\varz}} 
					(h|_{\R \times H_{\varz}}))(t,x)
			=	
				I_{\R \times \mathbb{H}_{\varz}} 
					\pi
					^{
						\R \times \mathbb{H}_{\varz}
					}
					_{
						\{0\} \times \mathbb{H}_{\varz}
					}.
		\end{split}
	\end{equation}
	In addition, we obtain from \eqref{eq: def of FN} that
	for all 
	$N \in \N_0$, $(t,x) \in W$, $r \in \R$, $p \in (\R \times H)'$,
	and all
	$
		(C_\lambda)_{\lambda \in (0, \infty)} 
			\subseteq 
				\mathbb{S}_{
					\R \times \mathbb{H}, 
					(\R \times \mathbb{H})'
				}
	$ 
	it holds that
	\begin{align}
				&\liminf_{\lambda \downarrow 0} \left(
					F_N(
						(t,x),
						r,
						p, 
						(C_\lambda|_{\R \times H_{\varz}})
							|_{\R \times H_{\varz}}
						-\lambda I_{\R \times \mathbb{H}_{\varz}} 
							\pi_{\{ 0\} \times V_N^\perp}^{\R \times H_\varz}
					)
				\right) \\ \nonumber
			={} 
				&\liminf_{\lambda \downarrow 0} \left(
					F_N(
						(t,x),
						r,
						p, 
						(C_\lambda|_{\R \times H_{\varz}})
							|_{\R \times H_{\varz}}
					)
					+\nicefrac 12
					\langle
						B(\pi^{H}_{V_N} x),
						I_{\mathbb{H}_{\varz}}^{-1} \lambda 
							\, I_{\mathbb{H}_{\varz}} \pi_{V_N^\perp}^{H_\varz}
								\, B(\pi^{H}_{V_N} x)  
					\rangle_{HS(\mathbb{U}, \mathbb{H}_{\varz})}
				\right) \\ \nonumber
			\leq{}  
				&\limsup_{\lambda \downarrow 0} 
				F_N(
					(t,x),
					r,
					p, 
					(C_\lambda|_{\R \times H_{\varz}})
						|_{\R \times H_{\varz}}
				).
	\end{align}
	Analogous it follows for all 
	$N \in \N_0$, $z \in W$, $r \in \R$, $p \in (\R \times H)'$,
	and all
	$
		(C_\lambda)_{\lambda \in (0, \infty)} 
			\subseteq 
			\mathbb{S}_{
				\R \times \mathbb{H}, 
				(\R \times \mathbb{H})'
			}
	$ 
	it hold that
	\begin{equation}
		\begin{split}
			&\limsup_{\lambda \downarrow 0} \left(
				F_N(
					z,
					r,
					p, 
					(C_\lambda|_{\R \times H_{\varz}})
						|_{\R \times H_{\varz}} 
					+ \lambda I_{\R \times \mathbb{H}_{\varz}} 
							\pi_{\{ 0\} \times V_N^\perp}^{\R \times H_\varz}
				)
			\right) \\
		\geq{}
			&\liminf_{\lambda \downarrow 0} 
				F_N(
					z,
					r,
					p, 
					(C_\lambda|_{\R \times H_{\varz}})
						|_{\R \times H_{\varz}} 
				).
		\end{split}
	\end{equation}
	Moreover, we derive that 
	for all $N \in \N$,
	$(t,x) \in W$, $r \in \R$, $\tilde{p} \in (\R \times H)'$, 
	$\tilde{C} \in \mathbb{S}_{\R \times \mathbb{H}_\varz, (\R \times \mathbb{H}_\varz)'}$,
	$p \in (\R \times V_N)'$, 
	and all
	$C \in \mathbb{S}_{\R \times \mathbb{V}_N, (\R \times \mathbb{V}_N)'}$
	with 
	\begin{equation}
					\forall (t_1,x_1) \in \R \times H_\varz \colon
					\langle p, (t_1, \pi^{H}_{V_N} x_1) \rangle_{(\R \times V_N)', \R \times V_N} 
				= 
					\langle 
						\tilde{p}|_{\R \times H_\varz}, (t_1, x_1) 
					\rangle_{(\R \times H_\varz)', \R \times H_\varz} 
		\end{equation} and with
		\begin{equation}
			\begin{split}
					&\forall (t_1, x_1),(t_2, x_2) \in \R \times H_\varz \colon \\
					&\big \langle 
						C ( t_1, \pi^{H}_{V_N}  x_1), 
						(t_2, \pi^{H}_{V_N} x_2)
					\big \rangle_{V_N', V_N} 
				={}
					\langle 
						(t_1, x_1), \tilde{C} (t_2, x_2) 
					\rangle_{\R \times H_\varz, (\R \times H_\varz)'}
			\end{split}
		\end{equation}
	it holds that
	\begin{align}
		\nonumber
				&F_N((t,x),r,\tilde{p},\tilde{C}) \\ \nonumber
			={}
				 &I^{-1}_{\R} \pi^{(\R \times H)'}_1 \tilde{p} 
						-\left\langle
							F(	\pi^{H}_{V_N} \, x) 
							+A(	\pi^{H}_{V_N} \, x), 
							\pi_{V_N'}^{H'} \pi^{(\R \times H)'}_2 \tilde{p}
						\right\rangle_{H, H' } \\ \nonumber
						&- \nicefrac 12
						\left \langle
							B(	\pi^{H}_{V_N} \, x) ,
							I_{\mathbb{H}_{\varz}}^{-1} 
								\pi_2^{(\R \times H_\varz)'} \tilde{C} 
								\, (0, B(	\pi^{H}_{V_N} \, x))  
						\right \rangle_{HS(\mathbb{U}, \mathbb{H}_{\varz})}\\ \nonumber
			={}
				&I^{-1}_{\R} \pi^{(\R \times V_N)'}_1 p
						-\left\langle
							\pi_{V_N}^{H} F(	\pi^{H}_{V_N} \, x) 
							+\pi_{V_N}^{H} A(	\pi^{H}_{V_N} \, x), 
							\pi^{(\R \times H)'}_2 \tilde{p}
						\right\rangle_{H, H' } \\
						&- \nicefrac 12
						\sum_{i \in \mathcal{I}}
						\left \langle
							B(	\pi^{H}_{V_N} \, x) \tilde{e}_i,
							\pi_2^{(\R \times H_\varz)'} \tilde{C} 
								\, (0, B(	\pi^{H}_{V_N} \, x) \tilde{e}_i ) 
						\right \rangle_{H_{\varz}, H'_{\varz}}\\ \nonumber
			={}
				&I^{-1}_{\R} \pi^{(\R \times V_N)'}_1 p
						-\left\langle
							\pi_{V_N}^{H} F(	\pi^{H}_{V_N} \, x) 
							+\pi_{V_N}^{H} A(	\pi^{H}_{V_N} \, x), 
							\pi^{(\R \times V_N)'}_2 p
						\right\rangle_{V_N, V_N' } \\ \nonumber 
						&-\nicefrac 12
						\sum_{i \in \mathcal{I}}
						\left \langle
							\pi_{V_N}^{H} B(	\pi^{H}_{V_N} \, x) \tilde{e}_i,
							\pi_2^{(\R \times H_\varz)'} C 
								\, (0, \pi_{V_N}^{H} B(	\pi^{H}_{V_N} \, x) \tilde{e}_i ) 
						\right \rangle_{V_N, V'_N}\\ \nonumber
			={}
				&I^{-1}_{\R} \pi^{(\R \times V_N)'}_1 p
						-\left\langle
							\pi_{V_N}^{H} F(	\pi^{H}_{V_N} \, x) 
							+\pi_{V_N}^{H} A(	\pi^{H}_{V_N} \, x), 
							\pi^{(\R \times V_N)'}_2 p
						\right\rangle_{V_N, V'_N } \\ \nonumber
						&- \nicefrac 12
						\left \langle
							\pi_{V_N}^{H} B(	\pi^{H}_{V_N} \, x),
							I^{-1}_{\mathbb{V}_N} \pi_2^{(\R \times H_\varz)'} C 
								\, (0, \pi_{V_N}^{H} B(	\pi^{H}_{V_N} \, x)) 
						\right \rangle_{HS(\mathbb{U}, \mathbb{V}_N)}\\ \nonumber
			={}
				&\tilde{F}(\pi_{V_N}^{H} x,r,p, C).
	\end{align}
	Therefore, Proposition \ref{prop: u proj is viscosity solution}
		(with $f \leftarrow \pi^{\R \times H}_{\R \times V_N}$,
		$\mathbb{H} \leftarrow \R \times \mathbb{H}$,
		$\mathbb{V} \leftarrow \R \times \mathbb{V}_N$,
		$\mathbb{X} \leftarrow \R \times \mathbb{H}_\varz$)
	yields that
	for all $N \in  \N$ it holds that
	$u_N$ is a viscosity solution of 
	$
			F_N|_{
				W 
				\times \R 
				\times H' 
				\times \mathbb{S}_{\mathbb{H}_{\varz}, \mathbb{H}_{\varz}'}
			} 
		= 0
	$ 
	relative to $(h, \mathbb{H}, \mathbb{H}_{\varz})$.
	Next observe that the continuity of $u_0$ 
	with respect to the $\| \cdot \|_{\R \times H}$-norm
	together with the definition of $h$
	shows that for all $(t, x) \in (0,T) \times H$ 
	and all $\delta \in (0, \infty)$ it holds that
	$
			(u_0)^{-, W}_{\mathbb{H}, \delta, h}(t,x) 
		=
			\overline{(u_n - \delta h)}_{\mathbb{H}}^W (t,x)
		= 
			u_0(t,x) - \tfrac \delta2 \|x\|^2_{H_{\varz}}.
	$
	Thus the continuity of $u_0$ 
	with respect to the $\| \cdot \|_{\R \times H}$-norm
	implies that
	for all
	$ 
		(t_n, x_n)_{n \in \N_0} 
			\subseteq 	
				(0,T) \times H
	$
	and all $\delta \in (0, \infty)$ with
	\begin{equation}
			\lim_{n \to \infty} \left(
				|
					(u_0)^{-, W}_{\mathbb{H}, \delta, h}(t_0, x_0)
					-(u_0)^{-, W}_{\mathbb{H}, \delta, h}(t_n, x_n)
				|
				\vee |t_0 - t_n| \vee \|x_n - x_0\|_{H}
			\right)
		= 0
	\end{equation}
	it holds that
	\begin{equation}
		\begin{split}
				0
			={}
				&\lim_{n \to \infty} \left(
					|
						(u_0)^{-, W}_{\mathbb{H}, \delta, h}(t_0, x_0)
						-(u_0)^{-, W}_{\mathbb{H}, \delta, h}(t_n, x_n)
					|
				\right) \\
			={}
				&\lim_{n \to \infty} \left(
					|
						u_0(t_0,x_0) - \tfrac \delta2 \|x_0\|^2_{H_{\varz}}
						-u_0(t_n,x_n) + \tfrac \delta2 \|x_n\|^2_{H_{\varz}}
					|
				\right)
			=
				\lim_{n \to \infty} \left(
					\tfrac \delta2
					\big |
						\|x_0\|^2_{H_{\varz}}
						- \|x_n\|^2_{H_{\varz}}
					\big |
				\right).
		\end{split}
	\end{equation}
	This the fact that for all
	$ 
		(x_n)_{n \in \N} 
			\subseteq 	
				H
	$
	and all $x_0 \in H_{\varz}$
	it holds that
	$
			\big(
				\lim_{n \to \infty} \|x_n -x_0 \|_{H_{\varz}} =0
			\big)
	\Leftrightarrow
		\big(
			\lim_{n \to \infty} \big(
				\|x_n -x_0 \|_{H} 
				\vee 
				\big |
					\|x_0\|^2_{H_{\varz}}- \|x_n\|^2_{H_{\varz}}
				\big |
			\big) =0 
		\big)
	$,
	the continuity of $u_0$ 
	with respect to the $\| \cdot \|_{\R \times \mathbb{H}}$-norm,
	and the continuity of
	$h|_{(0,T) \times H_{\varz}}$
	with respect to the $\| \cdot \|_{\R \times \mathbb{H}_\varz}$-norm
	show that
	for all
	$ 
		(t_n, x_n)_{n \in \N} 
			\subseteq 	
				(0,T) \times H
	$,
	$ 
		(t_0, x_0) \in (0,T) \times H_{\varz},
	$
	and all $\delta \in (0, \infty)$ it holds that
	\begin{equation}
		\begin{split}
			&\Big(
					\lim_{n \to \infty} \left(
						|
							(u_0)^{-, W}_{\mathbb{H}, \delta, h}(t_0, x_0)
							-(u_0)^{-, W}_{\mathbb{H}, \delta, h}(t_n, x_n)
						|
						\vee |t_0 - t_n| \vee \|x_n - x_0\|_{H}
					\right)
				= 0
			\Big) \\
		\Leftrightarrow
			&\Big (
				\lim_{n \to \infty}
					(\|x_n -x_0\|_{H_{\varz}} \vee |t_0 - t_n|) = 0
			\Big ).
		\end{split}
	\end{equation}
	Therefore we obtain with \eqref{eq:def of d delta - without t} that
	for all
	$ 
		((t_n, x_n),r_n,p_n,C_n)_{N \in \N_0} 
			\subseteq 	
				((0,T) \times H) 
				\times \R 
				\times (\R \times H)' 
				\times 
					\mathbb{S}_{
						\mathbb{\R \times H}, 
						(\R \times \mathbb{H})'
					},
	$
	and all $\delta \in (0, \infty)$ it holds that 
	\begin{equation}
	\label{eq: norm equivalence}
		\begin{split}
				&\Big(
						\lim_{n \to \infty} 
						d^{-, W}_{\mathbb{H}, \delta, h, u_0}(
							((t_0,x_0),r_0,p_0,C_0), ((t_n,x_n),r_n,p_n,C_n)
						)
					= 0
				\Big) \\
			\Leftrightarrow{}
				&\Big(
					\big(
							\lim_{n \to \infty} \big(
								\|x_n -x_0\|_{H_{\varz}} 
								\vee |t_n - t_0| 
								\vee \|  p_0 - p_n \|_{(\R \times H)'}
								\vee 
								\|C_0 - C_n\|_{
									L(\R \times \mathbb{H},(\R \times \mathbb{H})')
								}
							\big)
						=
							0
					\big) 
				\\ & \qquad
					\wedge (x_0 \in H_{\varz})
				\Big).
		\end{split}
	\end{equation}
	In addition, note that for all
	$\eps \in (0, \infty)$
	there exists an $M_\eps \in (0, \infty)$
	such that for all $x \in H_{\varz}$
	it holds that
	\begin{equation}
	\label{eq: norm interpolation}
		\begin{split}
				&\|x\|_{H_{\varz + 1/2 -\alpha_1}}
			\leq
				\eps \|x\|_{H_{\varz + 1/2}} + M_\eps \|x\|_{H_\varz},
			\textrm{ and that } \quad
				\|x\|^2_{H_{\varz +1/2 -\beta_1}}
			\leq
				\eps \|x\|^2_{H_{\varz + 1/2}} + M_\eps \|x\|^2_{H_\varz}	
		\end{split}
	\end{equation}
	and this together with
	\eqref{eq: F H 1/2 bound} 
	and with $\varz \geq \nicefrac 12$ implies that
	for all 
	$\eps, \delta \in (0, \infty)$, $\xi \in H_{2\varz}$, $x \in H_\varz$, 
	and all $\rho \in (\R \times H)'$
	with 
	$\|x- \xi\|_\varz \leq r_x$
	it holds that 
	\begin{align}
	\nonumber
				&-\left\langle
					F(\xi)+A(\xi), 
					\pi^{(\R \times H)'}_2 \rho 
					+\delta I_{\mathbb{H}} (-A)^{2\varz}(\xi)
				\right\rangle_{H, H' }\\ \nonumber
			\geq{}
					&-\delta \|F(\xi)\|_{H_{\varz - 1/2}} 
							\, \|\xi\|_{H_{1/2 + \varz}}
					-\|F(\xi)\|_{H} \, \|\rho\|_{(\R \times H)'}
					-\|\xi\|_{H_1} \, \|\rho\|_{(\R \times H)'}
					+\delta \|\xi\|^2_{H_{1/2 + \varz}} \\ \nonumber
				\geq{}
					&-\delta \|F(\xi)\|_{H_{\varz - 1/2}}
						\, \|\xi\|_{H_{1/2 + \varz}} 
					-(\sup_{n \in \N} |\lambda_n|^{1/2 -\varz}) 
						\|F(\xi)\|_{H_{\varz -1/2}} \, \|\rho\|_{(\R \times H)'} \\ \nonumber
					&-(\sup_{n \in \N} |\lambda_n|^{1/2 -\varz}) 
						\|\xi\|_{H_{1/2+ \varz}} \, \|\rho\|_{(\R \times H)'} 
					+\delta \|\xi\|^2_{H_{1/2 + \varz}} \\ \nonumber
				\geq{} 
					&-\|F(\xi)\|_{H_{\varz - 1/2}} \big(
						\delta \|\xi\|_{H_{1/2 + \varz}}  
						+ (\sup_{n \in \N} |\lambda_n|^{1/2 -\varz}) \|\rho\|_{(\R \times H)'}
					\big) \\ 
	\label{eq: F upper bound}
					&-(\sup_{n \in \N} |\lambda_n|^{1/2 -\varz}) 
						\|\xi\|_{H_{1/2+ \varz}} \, \|\rho\|_{(\R \times H)'} 
					+\delta \|\xi\|^2_{H_{1/2 + \varz}}\\ \nonumber
				\geq{}
					&-R_x (\|\xi \|_{H_{\varz + 1/2 - \alpha_1}} +1)\big(
						\delta \| \xi \|_{H_{1/2 + \varz}}  
						+ (\sup_{n \in \N} |\lambda_n|^{1/2 -\varz}) \| \rho \|_{(\R \times H)'}
					\big) \\ \nonumber
					&-(\sup_{n \in \N} |\lambda_n|^{1/2 -\varz}) 
						\| \xi \|_{H_{1/2+ \varz}} \, \| \rho \|_{(\R \times H)'} 
					+\delta \|\xi \|^2_{H_{1/2 + \varz}} \\ \nonumber
				\geq{}
					&-R_x (\eps \|\xi \|_{H_{\varz + 1/2}} 
					+ M_\eps\|\xi \|_{H_\varz} +1)
					\big(
						\delta \| \xi \|_{H_{1/2 + \varz}}  
						+ (\sup_{n \in \N} |\lambda_n|^{1/2 -\varz}) \| \rho \|_{(\R \times H)'}
					\big) \\ \nonumber
					&-(\sup_{n \in \N} |\lambda_n|^{1/2 -\varz}) 
						\|\xi \|_{H_{1/2+ \varz}} \, \| \rho \|_{(\R \times H)'} 
					+\delta \| \xi \|^2_{H_{1/2 + \varz}}. 
	\end{align}
	Moreover, it follows from
	\eqref{eq: B H 1 bound} and \eqref{eq: norm interpolation} that
	for all 
	$\eps, \delta \in (0, \infty)$, $\xi \in H_{2\varz}$, $x \in H_\varz$, 
	$\mathfrak{C} \in \mathbb{S}_{\R \times \mathbb{H}, (\R \times \mathbb{H})'}$
	and all $\mathfrak{D} \in \mathbb{S}_{\mathbb{H}, \mathbb{H}'}$ with
	$
		\forall y,z \in H \colon 
				\langle z, \mathfrak{D} y \rangle_{H,H'} 
			=
				\langle (0,z) , \mathfrak{C} (0,y) \rangle_{\R \times H,(\R \times H)'}
	$
	and with 
	$\|x- \xi\|_\varz \leq r_x$
	it holds that 
	\begin{align}
	\nonumber
					&-
					\left \langle
						B(\xi) ,
						I_{\mathbb{H}_{\varz}}^{-1} (\mathfrak{D} + \delta I_{\mathbb{H}_{\varz}})
							\, B(\xi)  
					\right \rangle_{HS(\mathbb{U}, \mathbb{H}_{\varz})} \\ \nonumber
				\geq{}
					&- \delta \|B( \xi)\|^2_{HS(\mathbb{U}, \mathbb{H}_{\varz})} 
					-\left \langle
						B(\xi) ,
						I_{\mathbb{H}}^{-1} \mathfrak{D}
							\, B(\xi)  
					\right \rangle_{HS(\mathbb{U}, \mathbb{H})} \\ \nonumber
				\geq{} 
					&- \delta \|B(\xi)\|^2_{HS(\mathbb{U}, \mathbb{H}_{\varz})} 
					-\|B(\xi)  \|^2_{HS(\mathbb{U}, \mathbb{H})}
					  \| \mathfrak{D}\|_{L(\mathbb{H}, \mathbb{H}')} \\ \label{eq: B upper bound}
				\geq{} 
					&-\|B(\xi)\|^2_{HS(\mathbb{U}, \mathbb{H}_{\varz})} \big(
						\delta 
						+(\sup_{n \in \N} |\lambda_n|^{-2\varz}) \,
							\| \mathfrak{C} \|_{L(\R \times \mathbb{H}, (\R \times \mathbb{H})')}
					\big) \\ \nonumber
				\geq{} 
					&-R_x (\|\xi \|^2_{H_{\varz + 1/2 -\beta_1}} +1) \big(
						\delta 
						+(\sup_{n \in \N} |\lambda_n|^{-2\varz}) \,
							\| \mathfrak{C} \|_{L(\R \times \mathbb{H}, (\R \times \mathbb{H})')}
					\big) \\ \nonumber
				\geq{} 
					&-R_x (\eps \|\xi \|^2_{H_{\varz + 1/2}} 
					+ M_\eps\|\xi \|^2_{H_\varz} +1) 
					\big(
						\delta 
						+(\sup_{n \in \N} |\lambda_n|^{-2\varz}) \,
							\| \mathfrak{C} \|_{L(\R \times \mathbb{H}, (\R \times \mathbb{H})')}
					\big). 
	\end{align}
	Therefore, we obtain 
	from \eqref{eq: F upper bound} and from \eqref{eq: B upper bound} that
	for all $\eps, \delta \in (0, \infty)$, $\xi \in H_{2\varz}$, 
	$x \in H_{\varz}$,
	$\nu \in \R$, $p,\rho \in (\R \times H)'$,
	$C, \mathfrak{C} \in \mathbb{S}_{(\R \times \mathbb{H}),  (\R \times \mathbb{H})'}$,
	$ \mathfrak{D} \in \mathbb{S}_{\mathbb{H},  \mathbb{H}'}$,
	and all $K_1, K_2, K_3 \in \R$
	satisfying that
	$\|x -\xi\|_{H_\varz} \leq 1 \wedge r_x$,
	that
	$\|p -\rho\|_{(\R \times H)'} \leq 1$,
	that
	$\|C - \mathfrak{C} \|_{L(\R \times \mathbb{H},(\R \times \mathbb{H})')} \leq 1$,
	that
	$
		\forall y,z \in H \colon 
				\langle z, \mathfrak{D} y \rangle_{H,H'} 
			=
				\langle (0,z) , \mathfrak{C} (0,y) \rangle_{\R \times H,(\R \times H)'}
	$,
	that
	\begin{align}
	\nonumber
				K_{1}
			={} 
				&\tfrac 32 R_x \, \delta 
				+\tfrac 12 R_x \, (\sup_{n \in \N} |\lambda_n|^{-2\varz}) \,
				(\| C\|_{L(\R \times \mathbb{H}, (\R \times \mathbb{H})')} +1),
		\\  \nonumber
				K_{2}
			={} 
				&R_x \, \eps (\sup_{n \in \N} |\lambda_n|^{1/2 -\varz}) (\|p\|_{(\R \times H)'}+1)
				+\delta R_x \, (M_\eps (\|x\|_{H_\varz} +1) +1) \\
				&+ (\sup_{n \in \N} |\lambda_n|^{1/2 -\varz}) (\|p\|_{(\R \times H)'} +1),
		\\  \nonumber
				K_{3}
			={} 
				&(\| p\|_{(\R \times H)'} +1) 
				+ R_x \, (M_\eps (\|x\|_{H_\varz} +1)+1) (\sup_{n \in \N} |\lambda_n|^{1/2 -\varz}) 
					(\|p\|_{(\R \times H)'} +1) \\ \nonumber
				&+\tfrac 12 R_x \, (M_\eps (\|x\|_{H_\varz} +1)^2+1) 
					(\delta + (\sup_{n \in \N} |\lambda_n|^{-2\varz})) 
					(\|C\|_{L(\R \times \mathbb{H}, (\R \times \mathbb{H})')} +1),
	\end{align}
	and satisfying that $\delta - 2 \eps K_1 > 0$
	it holds that
	\begin{align}
	\nonumber
	\label{eq: F and B upper bound}
				&I^{-1}_{\R} \pi^{(\R \times H)'}_1 \rho 
					-\left\langle
						F(\xi)+A(\xi), 
						\pi^{(\R \times H)'}_2 \rho 
						+\delta I_{\mathbb{H}} (-A)^{2\varz}(\xi)
					\right\rangle_{H, H' } 
			\\ \nonumber &
					-\nicefrac 12
					\left \langle
						B(\xi) ,
						I_{\mathbb{H}_{\varz}}^{-1} (\mathfrak{D} + \delta I_{\mathbb{H}_{\varz}})
							\, B(\xi)  
					\right \rangle_{HS(\mathbb{U}, \mathbb{H}_{\varz})} \\
			\geq{} 
				&\|\xi\|^2_{H_{1/2 + \varz}} (\delta -\eps K_{1})
				-\|\xi\|_{H_{1/2 + \varz}}  K_{2}
				-K_{3} \\ \nonumber
			={} 
				&\tfrac \delta2 \|\xi\|^2_{H_{1/2 + \varz}}
				+(\tfrac \delta2 -\eps K_{1}) 
					\left (
						\|\xi\|_{H_{1/2 + \varz}} 
						- \frac{K_{2}}{(\delta -2\eps K_{1})}
					\right )^2
				-\frac{K^2_{2}}{2(\delta -2\eps K_{1})} 
				-K_{3} \\ \nonumber
			\geq{} 
				&\tfrac \delta2 \|\xi\|^2_{H_{1/2 + \varz}}
				-\frac{K^2_{2}}{2(\delta -2\eps K_{1})} 
				-K_{3}. 
	\end{align}
	Furthermore, \eqref{eq: def F+ delta h} and
	\eqref{eq: h derivative existence} yield that
	for all 
	$\delta \in (0, \infty)$, $\xi \in H_{2\varz}$,
	$x \in H_\varz$, $\tau \in (0,T)$, 
	$\nu \in \R$, $\rho \in (\R \times H)'$,
	$\mathfrak{C} \in \mathbb{S}_{\R \times \mathbb{H}, (\R \times \mathbb{H})'}$
	and all $\mathfrak{D} \in \mathbb{S}_{\mathbb{H}, \mathbb{H}'}$ with
	$
		\forall y,z \in H \colon 
				\langle z, \mathfrak{D} y \rangle_{H,H'} 
			=
				\langle (0,z) , \mathfrak{C} (0,y) \rangle_{H,H'}
	$
	it holds that 
	\begin{align}
	\nonumber
				&(F_0)^+_{\R \times \mathbb{H},\R \times \mathbb{H}_\varz, \delta, h}
					((\tau,\xi),\nu,\rho,\mathfrak{C}) \\
	\label{eq: F0 + delta h}
			={}
				&F_0(
					(\tau,\xi),
					\nu+ \tfrac{\delta}{2} \|\xi\|^2_{H_\varz},
					\rho + \delta I_{\R \times \mathbb{H}}
					(0,(-A)^{2\varz}(\xi)),
					(\mathfrak{C} |_{\R \times H_{\varz}}) |_{\R \times H_{\varz}}
					+ \delta I_{\R \times \mathbb{H}_{\varz}} 
						\pi
						^{
							\R \times \mathbb{H}_{\varz}
						}
						_{
							\{0\} \times \mathbb{H}_{\varz}
						}	
				)\\ \nonumber
			={} 
				&I^{-1}_{\R} \pi^{(\R \times H)'}_1 \rho 
					-\left\langle
						F(\xi)+A(\xi), 
						\pi^{(\R \times H)'}_2 \rho 
						+\delta I_{\mathbb{H}} (-A)^{2\varz}(\xi)
					\right\rangle_{H, H' } \\ \nonumber
					&- \nicefrac 12
					\left \langle
						B(\xi) ,
						I_{\mathbb{H}_{\varz}}^{-1} (\mathfrak{D} + \delta I_{\mathbb{H}_{\varz}})
							\, B(\xi)  
					\right \rangle_{HS(\mathbb{U}, \mathbb{H}_{\varz})}.
	\end{align}
	Combining \eqref{eq: def F+ delta h u},
	\eqref{eq: norm equivalence},
	\eqref{eq: F and B upper bound},
	and \eqref{eq: F0 + delta h}
	implies then that for all $\delta \in (0, \infty)$, 
	$x \in H_{\varz}$,
	$r \in \R$, $p \in (\R \times H)'$,
	and all
	$C \in \mathbb{S}_{(\R \times \mathbb{H}),  (\R \times \mathbb{H})'}$
	there exist an $L_{\delta, x,p,C} \in (0, \infty)$ such that
	\begin{align}
	\label{eq: F+ delta h u with extra norm bound}
				&(F_0)^+_{\R \times \mathbb{H}, \R \times \mathbb{H}_\varz, \delta, h, u_0}
					((t,x),r ,p ,C) \\ \nonumber
			={} 
				&\lim_{\eps \downarrow 0} \inf 
					\big \{ 
						(F_0)_{\R \times \mathbb{H}, \R \times \mathbb{H}_\varz, \delta, h}^{+}
							((\tau,\xi),\nu ,\rho ,\mathfrak{C}) \colon 
						~\xi \in H_{2\varz},
						~\tau \in (0,T),
						~\nu \in \R,
						~\rho \in (\R \times H)', 
					\\ \nonumber &\qquad \qquad
						~\mathfrak{C} \in \mathbb{S}_{\R \times \mathbb{H}, (\R \times \mathbb{H})'}, 
						~d^{-,W}_{\R \times \mathbb{H}, \delta, h, u_0}
							(
								((t,x),r ,p ,C), 
								((\tau,\xi),\nu ,\rho ,\mathfrak{C})
							) 
						\leq \eps 
					\big \}		\\ \nonumber
			={} 
				&\lim_{\eps \downarrow 0} \inf 
					\big \{ 
						(F_0)_{\R \times \mathbb{H},\R \times \mathbb{H}_\varz, \delta, h}^{+}
							((\tau,\xi),\nu ,\rho ,\mathfrak{C}) \colon 
						~\xi \in H_{2\varz},
						~\tau \in (0,T),
						~\nu \in \R,
						~\rho \in (\R \times H)',
					\\ \nonumber &\qquad \qquad
						~\mathfrak{C} \in \mathbb{S}_{\R \times \mathbb{H}, (\R \times \mathbb{H})'}, 
						~\|\xi\|_{H_{1/2 + \varz}} \leq L_{\delta, x, p, C},
					\\ \nonumber &\qquad \qquad
						~d^{-,W}_{\R \times \mathbb{H}, \delta, h, u_0}
							(
								((t,x),r ,p ,C), 
								((\tau,\xi),\nu ,\rho ,\mathfrak{C})
							) 
						\leq \eps 
					\big \}.
	\end{align}
	Moreover, note that for all
	$(x_n)_{n \in \N_0} \subseteq H$ with 
	$\limsup_{n \to \infty} \|x_n\|_{H_{\varz + 1/2}} < \infty$
	and with
	$
		\lim_{n \to \infty} \|x_n 
	$
	$
		-x_0\|_{H} = 0
	$ 
	it holds that
	$x_0 \in H_{\varz +1/2}$ and that
	$(x_n)_{n \in \N}$ converges weakly to $x_0$ with respect to
	the $\| \cdot \|_{H_{\varz + 1/2}}$-norm.
	Hence, we derive from the fact that 
	$H_{\varz + 1/2 -\alpha_1}$ and 
	$H_{\varz + 1/2 -\beta_1}$ are compactly embedded in
	$H_{\varz + 1/2 }$ that for all 
	$y \in H_{\varz + 1/2}$ and all
	$(x_n)_{n \in \N_0} \subseteq H$ with 
	$\limsup_{n \to \infty} \|x_n\|_{H_{\varz + 1/2}} < \infty$
	and with
	$\lim_{n \to \infty} \|x_n -x_0\|_{H} = 0$
	it holds that
	\begin{equation}
	\label{eq: x convergence in H 3/2-}
			\lim_{n \to \infty} \left(
				\big \| 
					x_0- x_n
				\big \|^2_{H_{\varz + 1/2 - \alpha_1}}
				+\big \| 
					x_0- x_n
				\big \|^2_{H_{\varz + 1/2 - \beta_1}}
				+\big| \langle
					y, x-x_0
				\rangle_{H_{\varz + 1/2}} \big |
			\right)
		= 0.
	\end{equation}
	Thus we get from \eqref{eq: def F+ delta h u}, 
	\eqref{eq: F and B continuity},
	\eqref{eq: h derivative existence},
	\eqref{eq: x convergence in H 3/2-},
	the lower semicontinuity of 
	$ H \ni y \to \|y\|^2_{\vartheta + 1/2} \in [0,\infty]$
	with respect to the $\|\cdot\|_H$-norm,
	and from $1 \leq \varz + 1/2$,
	that 
	for all $\delta \in (0, \infty)$, 
	$
		((\tau_n, \xi_n), \nu_n, \rho_n, \mathfrak{C}_n, \mathfrak{D}_n)_{n \in \N} 
			\subseteq
				W \times \R \times (\R \times H)'
				\times \mathbb{S}_{(\R \times \mathbb{H}),  (\R \times \mathbb{H})'}
				\times  \mathbb{S}_{\mathbb{H}, \mathbb{H}'}
	$, 
	$x \in H_{\varz}$,
	$p\in (\R \times H)'$,
	$C \in \mathbb{S}_{(\R \times \mathbb{H}),  (\R \times \mathbb{H})'}$,
	and all
	$\tilde{C} \in \mathbb{S}_{\mathbb{H}, \mathbb{H}'}$,
	with 
	\begin{align}	
		&\forall y, z\in H \colon
				\langle (0,x), C (0,y) \rangle_{\R \times H, (\R \times H)'}
			=
				\langle y, \tilde{C} z \rangle_{H,H'}, \\
		&\forall y, z \in H, \forall n \in \N \colon
				\langle (0,y), \mathfrak{C}_n (0,z) \rangle_{\R \times H, (\R \times H)'}
			=
				\langle y, \mathfrak{D}_n z \rangle_{H,H'}, \\
	&\limsup_{n \to \infty} \|\xi_n\|_{\varz + 1/2} < \infty, \\
			&\lim_{n \to \infty}
				\|x-\xi_n\|_H \vee \|p -\rho_n\|_{(\R \times H)'} 
				\vee \|C - \mathfrak{C}\|_{L(\R \times \mathbb{H}, (\R \times \mathbb{H})')}
		= 0
	\end{align}
	it holds that
	\begin{align}
	\label{eq: F_0 + h lower bound}
				&\lim_{n \to \infty}
					(F_0)_{\R \times \mathbb{H}, \R \times \mathbb{H}_\varz, \delta, h}^{+}
						((\tau_n,\xi_n),\nu_n ,\rho_n ,\mathfrak{C}_n) \\ \nonumber
			={} 
				&\lim_{n \to \infty}
					F_0(
						(\tau_n,\xi_n),
						\nu_n+ \delta \|\xi_n\|^2_\varz,
						\rho_n + I_{\R \times \mathbb{H}},
						(0,(-A)^{2\varz}(\xi_n)),
						\mathfrak{C}_n 
						+ \delta I_{\R \times \mathbb{H}_{\varz}} 
							\pi
							^{
								\R \times \mathbb{H}_{\varz}
							}
							_{
								\{0\} \times \mathbb{H}_{\varz}
							}	
					) \\ \nonumber
			={} 
				&\lim_{n \to \infty} \Big(
					I^{-1}_{\R} \pi^{(\R \times H)'}_1 \rho_n
						-\left\langle
							F(\xi_n)+A(\xi_n), 
							\pi^{(\R \times H)'}_2 \rho_n +\delta I_{\mathbb{H}} (-A)^{2\varz}(\xi_n)
						\right\rangle_{H, H' } 
						\\ \nonumber & \qquad
						- \nicefrac 12
						\left \langle
							B(\xi_n) ,
							I_{\mathbb{H}_{\varz}}^{-1} (\mathfrak{D}_n + \delta I_{\mathbb{H}_{\varz}})
								\, B(\xi_n)  
						\right \rangle_{HS(\mathbb{U}, \mathbb{H}_{\varz})} 
				\Big) \\ \nonumber
			={} 
				&\lim_{n \to \infty} \Big(
					I^{-1}_{\R} \pi^{(\R \times H)'}_1 \rho_n
						-\left\langle
							F(\xi_n)+A(\xi_n), 
							\pi^{(\R \times H)'}_2 p 
						\right\rangle_{H, H' } 
					\\ \nonumber &\qquad 
						-\left\langle
							F(\xi_n)+A(\xi_n), 
							\pi^{(\R \times H)'}_2 (\rho_n - p) 
						\right\rangle_{H, H' } 
						-\delta 
						\left\langle
							(-A)^{\varz - 1/2} F(x) , 
							(-A)^{\varz + 1/2}(\xi_n)
						\right\rangle_{H } 
					\\ \nonumber &\qquad 
						-\delta 
						\left\langle
							(-A)^{\varz - 1/2} (F(x) -  F(\xi_n)) , 
							(-A)^{\varz + 1/2}(\xi_n)
						\right\rangle_{H } 
						+\delta 
						\left \|
							\xi_n
						\right \|^2_{H_{\varz+1/2} } 
					\\ \nonumber &\qquad 
						- \nicefrac 12
						\left \langle
							B(\xi_n) ,
							I_{\mathbb{H}_{\varz}}^{-1} (\mathfrak{D}_n + \delta I_{\mathbb{H}_{\varz}})
								\, B(\xi_n)  
						\right \rangle_{HS(\mathbb{U}, \mathbb{H}_{\varz})} 
				\Big) \\ \nonumber 
			\geq{} 
				&I^{-1}_{\R} \pi^{(\R \times H)'}_1 p
					-\left\langle
						F(x)+A(x), 
						\pi^{(\R \times H)'}_2 p 
					\right\rangle_{H, H' } 
					-\delta 
					\left\langle
						(-A)^{\varz - 1/2} F(x), 
						(-A)^{\varz + 1/2}(x)
					\right\rangle_{H } \\ \nonumber
					&+\delta \| x \|^2_{H_{\varz + 1/2}}
					- \nicefrac 12
					\left \langle
						B(x) ,
						I_{\mathbb{H}_{\varz}}^{-1} (\tilde{C} + \delta I_{\mathbb{H}_{\varz}})
							\, B(x)  
					\right \rangle_{HS(\mathbb{U}, \mathbb{H}_{\varz})}.
	\end{align}
	Therefore \eqref{eq: F+ delta h u with extra norm bound}
	and \eqref{eq: F_0 + h lower bound}
	ensure that
	for all $\delta \in (0, \infty)$, 
	$x \in H_{\varz + 1/2}$,
	$p\in (\R \times H)'$,
	$C \in \mathbb{S}_{(\R \times \mathbb{H}),  (\R \times \mathbb{H})'}$,
	and all
	$\tilde{C} \in \mathbb{S}_{\mathbb{H}, \mathbb{H}'}$,
	with 
	$	
		\forall y, z \in H \colon
				\langle (0,y), C (0,z) \rangle_{\R \times H, (\R \times H)'}
			=
				\langle y, \tilde{C} z \rangle_{H,H'},
	$
	it holds that
	\begin{align}
	\nonumber
				&(F_0)^+_{\R \times \mathbb{H}, \R \times \mathbb{H}_\varz, \delta, h, u_0}
					((t,x),r ,p ,C) \\ 
			\geq{} 
				&I^{-1}_{\R} \pi^{(\R \times H)'}_1 p
						-\left\langle
							F(x)+A(x), 
							\pi^{(\R \times H)'}_2 p 
						\right\rangle_{H, H' } 
						-\delta 
						\left\langle
							(-A)^{\varz - 1/2} F(x), 
							(-A)^{\varz + 1/2}(x)
						\right\rangle_{H } \\ \nonumber
						&+\delta \| x \|^2_{H_{\varz + 1/2}}
						- \nicefrac 12
						\left \langle
							B(x) ,
							I_{\mathbb{H}_{\varz}}^{-1} (\tilde{C} + \delta I_{\mathbb{H}_{\varz}})
								\, B(x)  
						\right \rangle_{HS(\mathbb{U}, \mathbb{H}_{\varz})}.
	\end{align}
	In addition, it holds for all $x \in H_{\varz + 1/2}$ that
	$ \forall N \in \N \colon \pi^{H}_{V_N} x \in H_{2 \varz}$ and that
	$\lim_{N \to \infty} \|x- \pi^{H}_{V_N} x\|_{H_{\varz + 1/2}} =0$.
	Thus \eqref{eq: def F+ delta h u}, \eqref{eq: def F+ delta h},
	\eqref{eq: F and B continuity},
	\eqref{eq: h derivative existence}, \eqref{eq: norm equivalence},
	and \eqref{eq: F0 + delta h} shows
	for all $\delta \in (0, \infty)$, 
	$x \in H_{\varz + 1/2}$,
	$p\in (\R \times H)'$,
	$C \in \mathbb{S}_{(\R \times \mathbb{H}),  (\R \times \mathbb{H})'}$,
	and all
	$\tilde{C} \in \mathbb{S}_{\mathbb{H}, \mathbb{H}'}$,
	with 
	$	
		\forall y, z \in H \colon
				\langle (0,y), C (0,z) \rangle_{\R \times H, (\R \times H)'}
			=
				\langle y, \tilde{C} z \rangle_{H,H'},
	$
	that
	\begin{align}
	\label{eq: F_0 + h upper bound}
			&(F_0)^+_{\R \times \mathbb{H}, \R \times \mathbb{H}_\varz, \delta, h, u_0}
					((t,x),r ,p ,C) \\ \nonumber
		\leq{} 
			&\lim_{N \to \infty}
					(F_0)_{\R \times \mathbb{H}, \R \times \mathbb{H}_\varz, \delta, h}^{+}
						((t,\pi^{H}_{V_N} x),r ,p ,C) \\ \nonumber
			={} 
				&\lim_{N \to \infty}
					F_0(
						(t,\pi^{H}_{V_N} x),
						r+ \delta \|\pi^{H}_{V_N} x\|^2_\varz,
						p + I_{\R \times \mathbb{H}},
						(0,(-A)^{2\varz}(\pi^{H}_{V_N} x)),
						C 
						+ \delta I_{\R \times \mathbb{H}_{\varz}} 
							\pi
							^{
								\R \times \mathbb{H}_{\varz}
							}
							_{
								\{0\} \times \mathbb{H}_{\varz}
							}	
					) \\ \nonumber
			={} 
				&\lim_{N \to \infty} \Big(
					I^{-1}_{\R} \pi^{(\R \times H)'}_1 p
						-\left\langle
							F(\pi^{H}_{V_N} x)+A(\pi^{H}_{V_N} x), 
							\pi^{(\R \times H)'}_2 p 
							+\delta I_{\mathbb{H}} (-A)^{2\varz}(\pi^{H}_{V_N} x)
						\right\rangle_{H, H' } \\ \nonumber
					&\qquad - \nicefrac 12
						\left \langle
							B(\pi^{H}_{V_N} x) ,
							I_{\mathbb{H}_{\varz}}^{-1} (\tilde{C} + \delta I_{\mathbb{H}_{\varz}})
								\, B(\pi^{H}_{V_N} x)  
						\right \rangle_{HS(\mathbb{U}, \mathbb{H}_{\varz})} 
				\Big) \\ \nonumber
			={} 
				&\lim_{N \to \infty} \Big(
					I^{-1}_{\R} \pi^{(\R \times H)'}_1 p
						-\left\langle
							F(\pi^{H}_{V_N} x)+A(\pi^{H}_{V_N} x), 
							\pi^{(\R \times H)'}_2 p 
						\right\rangle_{H, H' } 
						\\ \nonumber & \qquad
						-\delta 
						\left\langle
							(-A)^{\varz - 1/2} F(\pi^{H}_{V_N} x) , 
							(-A)^{\varz + 1/2}(\pi^{H}_{V_N} x)
						\right\rangle_{H } 
						+\delta 
						\big \|
							\pi^{H}_{V_N} x
						\big \|^2_{H_{\varz+1/2} }
					\\ \nonumber & \qquad
						- \nicefrac 12
						\left \langle
							B(\pi^{H}_{V_N} x) ,
							I_{\mathbb{H}_{\varz}}^{-1} (\tilde{C} + \delta I_{\mathbb{H}_{\varz}})
								\, B(\pi^{H}_{V_N} x)  
						\right \rangle_{HS(\mathbb{U}, \mathbb{H}_{\varz})} 
				\Big) \\ \nonumber
			={} 
				&I^{-1}_{\R} \pi^{(\R \times H)'}_1 p
					-\left\langle
						F(x)+A(x), 
						\pi^{(\R \times H)'}_2 p 
					\right\rangle_{H, H' } 
					-\delta 
					\left\langle
						(-A)^{\varz - 1/2} F(x), 
						(-A)^{\varz + 1/2}(x)
					\right\rangle_{H } 
				\\ \nonumber
				& \qquad
					+\delta \| x \|^2_{H_{\varz + 1/2}}
					- \nicefrac 12
					\left \langle
						B(x) ,
						I_{\mathbb{H}_{\varz}}^{-1} (\tilde{C} + \delta I_{\mathbb{H}_{\varz}})
							\, B(x)  
					\right \rangle_{HS(\mathbb{U}, \mathbb{H}_{\varz})}.
	\end{align}
	Combining \eqref{eq: F_0 + h lower bound} and \eqref{eq: F_0 + h upper bound}
	yields that for all $\delta \in (0, \infty)$, 
	$x \in H_{\varz + 1/2}$,
	$p\in (\R \times H)'$,
	$C \in \mathbb{S}_{(\R \times \mathbb{H}),  (\R \times \mathbb{H})'}$,
	and all
	$\tilde{C} \in \mathbb{S}_{\mathbb{H}, \mathbb{H}'}$,
	with 
	$	
		\forall y, z \in H \colon
				\langle (0,y), C (0,z) \rangle_{\R \times H, (\R \times H)'}
			=
				\langle y, \tilde{C} z \rangle_{H,H'},
	$
	it holds that
	\begin{align}
	\nonumber
				&(F_0)^+_{\R \times \mathbb{H}, \R \times \mathbb{H}_\varz, \delta, h, u_0}
					((x,t),r ,p ,C) \\
	\label{eq: F0 delta + u}
			={} 
				&I^{-1}_{\R} \pi^{(\R \times H)'}_1 p
					-\left\langle
						F(x)+A(x), 
						\pi^{(\R \times H)'}_2 p 
					\right\rangle_{H, H' } 
					-\delta 
					\left\langle
						(-A)^{\varz - 1/2} F(x), 
						(-A)^{\varz + 1/2}(x)
					\right\rangle_{H } \\ \nonumber
					&+\delta \| x \|^2_{H_{\varz + 1/2}}
					-\nicefrac 12
					\left \langle
						B(x) ,
						I_{\mathbb{H}_{\varz}}^{-1} (\tilde{C} + \delta I_{\mathbb{H}_{\varz}})
							\, B(x)  
					\right \rangle_{HS(\mathbb{U}, \mathbb{H}_{\varz})}.
	\end{align}	
	Furthermore, the fact that
	for all $N \in \N$ it holds that
	$V_N$ is finite-dimensional implies
	for all $N \in \N$, 
	$r \in \R$,
	$(x_n)_{n \in \N_0} \subseteq H$
	with 
	$
		\lim_{n \to \infty} 
				\| \pi^{H}_{V_N} (x_n - x_0)\|_{H} 
		=0
	$
	that
	\begin{equation}
	\label{eq: finite dim norm equivalence}
			\lim_{n \to \infty} 
				\| \pi^{H}_{V_N} (x_n - x_0)\|_{H_{r}} 
		=0.
	\end{equation}
	Therefore, we get from \eqref{eq: def F+ delta h},
	\eqref{eq: def F+ delta h u},
	\eqref{eq: F and B continuity},
	and from \eqref{eq: finite dim norm equivalence}
	that for all $N \in \N$,
	$\delta \in (0, \infty)$, 
	$x \in H$,
	$p\in (\R \times H)'$,
	$C \in \mathbb{S}_{(\R \times \mathbb{H}),  (\R \times \mathbb{H})'}$,
	and all
	$\tilde{C} \in \mathbb{S}_{\mathbb{H}, \mathbb{H}'}$
	with 
	$	
		\forall y, z \in H \colon
				\langle y, \tilde{C} z \rangle_{H,H'}
			=
				\langle (0,y), C (0,z) \rangle_{\R \times H, (\R \times H)'}
	$
	it holds that
	\begin{align}
	\label{eq: FN delta + u}
			&
			(F_N)^+_{\R \times \mathbb{H}, \R \times \mathbb{H}_\varz, \delta, h, u_0}
					((t,x),r ,p ,C) \\
			={} 
				&\nonumber
				\lim_{\eps \downarrow 0} \inf 
					\big \{ 
						(F_N)_{\R \times \mathbb{H}, \R \times \mathbb{H}_\varz, \delta, h}^{+}
							((\tau,\xi),\nu ,\rho ,\mathfrak{C}) \colon 
						\xi \in H_{2\varz},
						~\tau \in (0,T),
						~\nu \in \R,
						~\rho \in (\R \times H)', 
				\\& \nonumber \qquad
						~\mathfrak{C} \in \mathbb{S}_{\R \times \mathbb{H}, (\R \times \mathbb{H})'},
						~d^{-,W}_{\R \times \mathbb{H}, \delta, h, u_0}
							(
								((x,t),r ,p ,C), 
								((\xi,\tau),\nu ,\rho ,\mathfrak{C})
							) 
						\leq \eps 
					\big \}		\\
			={} 
				& \nonumber
				\lim_{\eps \downarrow 0} \inf 
					\big \{ 
						F_N(
						(\tau,\pi^{H}_{V_N} \xi),
						\nu+ \delta \|\pi^{H}_{V_N} \xi \|^2_\varz,
						\rho + I_{\R \times \mathbb{H}}
						(0,(-A)^{2\varz}(\xi)),
						\mathfrak{C} 
						+ \delta I_{\R \times \mathbb{H}_{\varz}} 
							\pi
							^{
								\R \times \mathbb{H}_{\varz}
							}
							_{
								\{0\} \times \mathbb{H}_{\varz}
							}	
					) \colon 
					\\& \nonumber \qquad
						~\xi \in H_{2\varz},
						~\tau \in (0,T),
						~\nu \in \R, 
						~\rho \in (\R \times H)',
						~\mathfrak{C} \in \mathbb{S}_{\R \times \mathbb{H}, (\R \times \mathbb{H})'}, 
					\\& \nonumber \qquad 
						~d^{-,W}_{\R \times \mathbb{H}, \delta, h, u_0}
							(
								((x,t),r ,p ,C), 
								((\xi,\tau),\nu ,\rho ,\mathfrak{C})
							) 
						\leq \eps 
					\big \}		\\
			={} 
				& \nonumber
				\lim_{\eps \downarrow 0} \inf 
					\big \{ 
						\rho_1
						-\left\langle
							F(\pi^{H}_{V_N} \xi)+A(\pi^{H}_{V_N} \xi), 
							\pi^{H'}_{V'_N} \rho_2 
							+\delta \pi^{H'}_{V'_N} I_{\mathbb{H}} (-A)^{2\varz}(\xi)
						\right\rangle_{H, H' } 
					\\& \nonumber \qquad
						-\nicefrac 12
						\left \langle
							B(\pi^{H}_{V_N} \xi) ,
							I_{\mathbb{H}_{\varz}}^{-1} (\mathfrak{C} + \delta I_{\mathbb{H}_{\varz}})
								\, B(\pi^{H}_{V_N} \xi)  
						\right \rangle_{HS(\mathbb{U}, \mathbb{H}_{\varz})} 
						\colon 
						\xi \in H_{2\varz},
						~\tau \in (0,T),
					\\& \nonumber \qquad
						~\nu, \rho_1 \in \R, 
						~\rho_2 \in H',
						~\mathfrak{C} \in \mathbb{S}_{\mathbb{H}, \mathbb{H}'}, 
						~|\rho_1 - I^{-1}_{\R} \pi^{(\R \times H)'}_1 p| \leq \eps, 
					\\&\nonumber \qquad 
						~d^{-,W}_{\R \times \mathbb{H}, \delta, h, u_0}
							(
								(t,x,r , \pi^{(\R \times H)'}_2 p , \tilde{C}), 
								(\tau, \xi, \nu ,\rho_2 ,\mathfrak{C})
							) 
						\leq \eps 
					\big \}		\\
			={} 
				& \nonumber
				\lim_{\eps \downarrow 0} \inf 
					\big \{ 
						\rho_1
						-\left\langle
							F(\pi^{H}_{V_N} \xi)+A(\pi^{H}_{V_N} \xi), 
							\pi^{H'}_{V'_N} \rho_2
						\right\rangle_{H, H' }
						+\delta 
						\big \|
							\pi^{H}_{V_N} \xi
						\big \|^2_{H_{\varz+1/2} }
					\\&  \nonumber \qquad 
						-\delta 
						\left\langle
							(-A)^{\varz - 1/2} F(\pi^{H}_{V_N} \xi) , 
							(-A)^{\varz + 1/2}(\pi^{H}_{V_N} \xi)
						\right\rangle_{H } 
					\\ \nonumber & \qquad 
						- \nicefrac 12
						\left \langle
							B(\pi^{H}_{V_N} \xi) ,
							I_{\mathbb{H}_{\varz}}^{-1} (\mathfrak{C} + \delta I_{\mathbb{H}_{\varz}})
								\, B(\pi^{H}_{V_N} \xi)  
						\right \rangle_{HS(\mathbb{U}, \mathbb{H}_{\varz})} 
						\colon 
						\xi \in H_{2\varz},
						~\tau \in (0,T),
					\\& \nonumber \qquad
						~\nu, \rho_1 \in \R, 
						~\rho_2 \in H',
						~\mathfrak{C} \in \mathbb{S}_{\mathbb{H}, \mathbb{H}'}, 
						~|\rho_1 - I^{-1}_{\R} \pi^{(\R \times H)'}_1 p| \leq \eps,
					\\ \nonumber&\qquad 
						~d^{-,W}_{\R \times \mathbb{H}, \delta, h, u_0}
							(
								(t,x,r , \pi^{(\R \times H)'}_2 p , \tilde{C}), 
								(\tau, \xi, \nu ,\rho_2 ,\mathfrak{C})
							) 
						\leq \eps 
					\big \}		 \\
			={} 
				&\nonumber 
				I^{-1}_{\R} \pi^{(\R \times H)'}_1 p
					-\left\langle
						F(\pi^{H}_{V_N} x)+A(\pi^{H}_{V_N} x), 
						\pi^{H'}_{V'_N}  \pi^{(\R \times H)'}_2 p 
					\right\rangle_{H, H' } 
					+\delta \| \pi^{H}_{V_N} x \|^2_{H_{\varz + 1/2}}
				\\& \nonumber 
					-\delta 
					\left\langle
						(-A)^{\varz - 1/2} F(\pi^{H}_{V_N} x), 
						(-A)^{\varz + 1/2}(\pi^{H}_{V_N} x)
					\right\rangle_{H } \\ \nonumber
					&- \nicefrac 12
					\left \langle
						B(\pi^{H}_{V_N} x) ,
						I_{\mathbb{H}_{\varz}}^{-1} (\tilde{C} + \delta I_{\mathbb{H}_{\varz}})
							\, B(\pi^{H}_{V_N} x)  
					\right \rangle_{HS(\mathbb{U}, \mathbb{H}_{\varz})}.
	\end{align}
	Combining
	\eqref{eq: norm equivalence},
	\eqref{eq: F and B upper bound},
	and \eqref{eq: FN delta + u}
	yields
	for all $\eps, \delta \in (0, \infty)$,
	$ 
		((t_N, x_N),r_N,p_N,C_N)_{N \in \N_0} 
	$
	$
			\subseteq 	
				(\R \times H_{\varz}) 
				\times \R 
				\times (\R \times H)' 
				\times 
					\mathbb{S}_{
						\mathbb{\R \times H}, 
						(\R \times \mathbb{H})'
					}
	$,
	and all $K_1, K_2, K_3 \in \R$
	with 
	\begin{align}
			&\limsup_{N \to \infty} \|\pi^{H}_{V_N} x_N\|_{H_{\varz + 1/2}}
		= \infty, \\
		&\lim_{N \to \infty} 
			d^{-, W}_{\mathbb{H}, \delta, h, u_0}(
				((t_0, x_0),r_0,p_0,C_0), ((t_N, x_N),r_N,p_N,C_N)
			)
		= 0, \\
				&K_{1}
			={} 
				2R_{x_0} \, \delta 
				+R_{x_0} \, (\sup_{n \in \N} |\lambda_n|^{-2\varz}) \,
				(\| C_0\|_{L(\R \times \mathbb{H}, (\R \times \mathbb{H})')} +1),
		\\ 
		\begin{split}
				&K_{2}
			={} 
				R_{x_0} \, \eps (\sup_{n \in \N} |\lambda_n|^{1/2 -\varz}) 
					(\|p_0\|_{(\R \times H)'}+1) 
				+\delta R_{x_0} \, (M_\eps (\|x_0\|_{H_\varz} +1) +1) 
				\\& \qquad+ (\sup_{n \in \N} |\lambda_n|^{1/2 -\varz}) (\|p_0\|_{(\R \times H)'} +1),
			\end{split}
		\end{align}
		\begin{align}
			\begin{split}
				&K_{3}
			={} 
				(\| p_0\|_{(\R \times H)'} +1) 
				+ R_{x_0} \, (M_\eps (\|x_0\|_{H_\varz} +1)+1) 
					(\sup_{n \in \N} |\lambda_n|^{1/2 -\varz}) 
					(\|p_0\|_{(\R \times H)'} +1) 
				\\ &\qquad
				+R_{x_0} \, (M_\eps (\|x_0\|_{H_\varz} +1)^2+1) 
					(\delta + (\sup_{n \in \N} |\lambda_n|^{-2\varz})) 
					(\|C_0\|_{L(\R \times \mathbb{H}, (\R \times \mathbb{H})')} +1),
		\end{split}
	\end{align}
	and with $\delta - 2 \eps K_1 > 0$
	that
	\begin{align}
	\nonumber
				&\limsup_{N \to \infty} \,
					(F_N)^+_{
						\R \times \mathbb{H}, 
						\R \times \mathbb{H}_{\varz}, 
						\delta,
						h, 
						u_0
					}((t_N, x_N),r_N,p_N,C_N) \\ \nonumber
			={} 
				&\limsup_{N \to \infty} \Big(
					I^{-1}_{\R} \pi^{(\R \times H)'}_1 p_N
					-\left\langle
						F(\pi^{H}_{V_N} x_N)+A(\pi^{H}_{V_N} x_N), 
						\pi^{H'}_{V'_N}  \pi^{(\R \times H)'}_2 p_N 
					\right\rangle_{H, H' } 
				\\& \qquad
					+\delta \| \pi^{H}_{V_N} x_N \|^2_{H_{\varz + 1/2}} 
					-\delta 
					\left\langle
						(-A)^{\varz - 1/2} F(\pi^{H}_{V_N} x_N), 
						(-A)^{\varz + 1/2}(\pi^{H}_{V_N} x_N)
					\right\rangle_{H } 
				\\ \nonumber & \qquad
					- \nicefrac 12
					\left \langle
						B(\pi^{H}_{V_N} x_N) ,
						I_{\mathbb{H}_{\varz}}^{-1} (\tilde{C} + \delta I_{\mathbb{H}_{\varz}})
							\, B(\pi^{H}_{V_N} x_N)  
					\right \rangle_{HS(\mathbb{U}, \mathbb{H}_{\varz})}
				\Big) \\ \nonumber
			\geq{} 
				&\limsup_{N \to \infty} \Big(
					\tfrac \delta2 \| \pi^{H}_{V_N} x_N\|^2_{H_{1/2 + \varz}}
					-\frac{K^2_{2}}{2(\delta -2\eps K_{1})} 
					-K_{3}
				\Big)
			= \infty. 
	\end{align}
	Thus we derive from
	the lower semicontinuity of 
	$ 
		H \ni y 
			\to \| y \|_{H_{\varz + 1/2}} \in [0, \infty]
	$
	with respect to the $\| \cdot \|_{H}$-norm
	that for all
	$\delta \in (0, \infty)$ and all
	$ 
		((t_N, x_N),r_N,p_N,C_N)_{N \in \N_0} 
			\subseteq 	
				((0,T) \times H) 
				\times \R 
				\times (\R \times H)' 
				\times 
					\mathbb{S}_{
						\mathbb{\R \times H}, 
						(\R \times \mathbb{H})'
					}
	$
	with
	$
		\lim_{N \to \infty} 
			d^{-, W}_{\R \times \mathbb{H}, \delta, h, u_0}(
				((t_0, x_0),r_0,p_0,C_0), ((t_N, x_N),r_N,p_N,C_N)
			)
		= 0
	$
	and with 
	$
		\limsup_{N \to \infty}
			(F_N)^+_{
				\R \times \mathbb{H}, 
				\R \times \mathbb{H}_{\varz},
				\delta, 
				h, 
				u_0
			}
		((t_N, x_N),r_N,p_N,C_N) \leq 0
	$
	it holds that
	\begin{equation}
	\label{eq: xN bound in H 3/2}
			\|x_0\|_{H_{\varz + 1/2}}
		\leq
			\limsup_{N \to \infty}
				\|\pi^{H}_{V_N} x_N\|_{H_{\varz + 1/2}}
		< \infty
	\end{equation}
	and therefore
	\begin{equation}
	\label{eq: x0 in H_3/2}
		x_0 \in H_{\varz + 1/2}.
	\end{equation} 
	Combining
	\eqref{eq: F and B continuity}
	\eqref{eq: norm equivalence},
	\eqref{eq: x convergence in H 3/2-},
	\eqref{eq: FN delta + u},
	\eqref{eq: F0 delta + u},
	\eqref{eq: xN bound in H 3/2},
	\eqref{eq: x0 in H_3/2},
	the fact that $\varz + 1/2 \geq 1$,
	and the lower semicontinuity of 
	$H \ni x \to \|x \|^2_{H_{\varz + 1/2}} \in [0, \infty]$
	with respect to the $\| \cdot \|_{H}$-norm
	shows that
	for all 
	$ 
		((t_N, x_N),r_N,p_N,C_N)_{N \in \N_0} 
			\subseteq 	
				((0,T) \times H) 
				\times \R 
				\times (\R \times H)' 
				\times 
					\mathbb{S}_{
						\mathbb{\R \times H}, 
						(\R \times \mathbb{H})'
					}
	$
	and all $\delta \in (0, \infty)$ with
	\begin{equation}
		\lim_{N \to \infty} 
			d^{-, W}_{\R \times \mathbb{H}, \delta, h, u_0}(
				((t_0, x_0),r_0,p_0,C_0), ((t_n, x_N),r_N,p_N,C_N)
			)
		= 0
	\end{equation}
	and with 
	$
		\limsup_{N \to \infty}
			(F_N)^+_{
				\R \times \mathbb{H}, 
				\R \times \mathbb{H}_{\varz}, 
				\delta, 
				h, 
				u_0
			}((t_N, x_N),r_N,p_N,C_N) \leq 0
	$
	it holds that 
	\begin{align}
	\nonumber
				&\liminf_{ N \to \infty } \big(
					(F_0)^+_{
						\R \times \mathbb{H}, 
						\R \times \mathbb{H}_{\varz}, 
						\delta, 
						h, 
						u_0
					}
						((t_0, x_0),r_0,p_0,C_0)
				\\ \nonumber & \qquad\qquad
					-
					(F_N)^+_{
						\R \times \mathbb{H}, 
						\R \times \mathbb{H}_{\varz}, 
						\delta, 
						h, 
						u_0
					}
						((t_N, x_N),r_N,p_N,C_N)
				\big)	\\ \nonumber
			={} 
				&\liminf_{ N \to \infty } \Big(
				I^{-1}_{\R} \pi^{(\R \times H)'}_1 p_0
					-\left\langle
						F(x_0)+A(x_0), 
						\pi^{(\R \times H)'}_2 p_0
					\right\rangle_{H, H' } 
				\\ \nonumber & \qquad
					-\delta 
					\left\langle
						(-A)^{\varz - 1/2} F(x_0), 
						(-A)^{\varz + 1/2}(x_0)
					\right\rangle_{H } 
					+\delta \| x_0 \|^2_{H_{\varz + 1/2}}
				\\ \nonumber & \qquad
					-\nicefrac 12
					\left \langle
						B(x_0) ,
						I_{\mathbb{H}_{\varz}}^{-1} (\tilde{C}_0 + \delta I_{\mathbb{H}_{\varz}})
							\, B(x_0)  
					\right \rangle_{HS(\mathbb{U}, \mathbb{H}_{\varz})} 
					-I^{-1}_{\R} \pi^{(\R \times H)'}_1 p_N
			\\& \nonumber \qquad
					+\left\langle
						F(\pi^{H}_{V_N} x_N)+A(\pi^{H}_{V_N} x_N), 
						\pi^{H'}_{V'_N}  \pi^{(\R \times H)'}_2 p_N
					\right\rangle_{H, H' } \\
				&	 \nonumber \qquad
					+\delta 
					\left\langle
						(-A)^{\varz - 1/2} F(\pi^{H}_{V_N} x_N), 
						(-A)^{\varz + 1/2}(\pi^{H}_{V_N} x_N)
					\right\rangle_{H } 
					-\delta \| \pi^{H}_{V_N} x_N \|^2_{H_{\varz + 1/2}}\\
	\label{eq: F0+ - FN+ limit on H0}
				& \qquad
					+\nicefrac 12
					\left \langle
						B(\pi^{H}_{V_N} x_N) ,
						I_{\mathbb{H}_{\varz}}^{-1} (\tilde{C}_N + \delta I_{\mathbb{H}_{\varz}})
							\, B(\pi^{H}_{V_N} x_N)  
					\right \rangle_{HS(\mathbb{U}, \mathbb{H}_{\varz})}
				\Big) \\ \nonumber
			={} 
				&\liminf_{ N \to \infty } \Big(
					-\left\langle
						F(x_0) - F(x_N) + A(x_0) - A(x_N), 
						\pi^{(\R \times H)'}_2 p_0
					\right\rangle_{H, H' }
				\\ \nonumber & \qquad
					-\left\langle
						F(x_N) + A(x_N), 
						\pi^{(\R \times H)'}_2 (p_0 - p_N) 
					\right\rangle_{H, H' } 
				\\ \nonumber & \qquad
					-\delta 
					\left\langle
						(-A)^{\varz - 1/2} F(x_0), 
						(-A)^{\varz + 1/2}(x_0 -x_N)
					\right\rangle_{H } 
				\\  \nonumber & \qquad
					-\delta 
					\left\langle
						(-A)^{\varz - 1/2} (F(x_0) -F(x_N)), 
						(-A)^{\varz + 1/2} (x_N)
					\right\rangle_{H } 
					+\delta \| x_0 \|^2_{H_{\varz + 1/2}} 
				\\ \nonumber & \qquad
					-\delta \| x_N \|^2_{H_{\varz + 1/2}}
					-\nicefrac 12
					\left \langle
						B(x_0) ,
						I_{\mathbb{H}_{\varz}}^{-1} (\tilde{C}_0 + \delta I_{\mathbb{H}_{\varz}})
							\, B(x_0)  
					\right \rangle_{HS(\mathbb{U}, \mathbb{H}_{\varz})}
			\\ \nonumber & \qquad
					+\nicefrac 12
					\left \langle
						B(\pi^{H}_{V_N} x_N) ,
						I_{\mathbb{H}_{\varz}}^{-1} (\tilde{C}_N + \delta I_{\mathbb{H}_{\varz}})
							\, B(\pi^{H}_{V_N} x_N)  
					\right \rangle_{HS(\mathbb{U}, \mathbb{H}_{\varz})}
				\Big) \\ \nonumber
			={} 
				&\liminf_{ N \to \infty } \Big(
					\delta \| x_0 \|^2_{H_{\varz + 1/2}}
					- \delta 
							\|
								\pi^{H}_{V_N} x_N 
							\|^2_{H_{\varz + 1/2}} 
				\Big)
			\leq 0
	\end{align}
	which shows \eqref{eq:convergence F+}.
	Analogously it follows that
	for all
	$ 
		((t_N, x_N),r_N,p_N,C_N)_{N \in \N_0} 
			\subseteq 	
				((0,T) \times H) 
				\times \R 
				\times (\R \times H)' 
				\times 
					\mathbb{S}_{
						\mathbb{\R \times H}, 
						(\R \times \mathbb{H})'
					}
	$
	and all $\delta \in (0, \infty)$ with
	\begin{equation}
		\lim_{N \to \infty} 
			d^{+, W}_{\R \times \mathbb{H}, \delta, h, u_0}(
				((t_0, x_0),r_0,p_0,C_0), ((t_N, x_N),r_N,p_N,C_N)
			)
		= 0
	\end{equation}
	and with
	$
			\limsup_{N \to \infty}
				(F_N)^-_{
					\R \times \mathbb{H}, 
					\R \times \mathbb{H}_{\varz}, 
					\delta, 
					h, 
					u_0
				}
			((t_N, x_N),r_N,p_N,C_N) 
		\geq 0
	$
	it holds that
	\begin{equation}
	\label{eq: F0- - FN- limit on H}
		\begin{split}
			\limsup_{ N \to \infty } \big(
				&(F_0)^-_{
					\R \times \mathbb{H}, 
					\R \times \mathbb{H}_{\varz}, 
					\delta, 
					h, 
					u_0
				}
					((t_0, x_0),r_0,p_0,C_0) \\
				&-
				(F_N)^-_{
					\R \times \mathbb{H}, 
					\R \times \mathbb{H}_{\varz}, 
					\delta, 
					h, 
					u_0
				}
					((t_N, x_N),r_N,p_N,C_N)
			\big)
		\geq
			0.
		\end{split}
	\end{equation}
	It remains to prove that $u_0$ is bounded 
	on $\mathbb{H}$-bounded subsets of 
	$H$ and 
	fulfills
	\eqref{eq: uniformly bounded at 0 with V thm}
	and \eqref{eq: grow condition on u}.
	Therefore note that it follows from \eqref{eq: V bound for X0}
	and from \eqref{eq: V bound for phi}  that
	for all $R \in (0, \infty)$ it holds that
	\begin{equation}
		\begin{split}
				&\sup_{ 
					x \in \{ y \in H \colon \|y\|_{H} \leq R \} 
				} 
				\sup_{t \in [0,T]}
					|u_0(t,x)|
			\leq
				\sup_{ 
					x \in \{ y \in H \colon \|y\|_{H} \leq R \} 
				} 
				\sup_{t \in [0,T]}
					\E[|\varphi(X^{0,x}_t)|] \\
			\leq{} 
				&\sup_{ 
					x \in \{ y \in H \colon \|y\|_{H} \leq R \} 
				} 
				\sup_{t \in [0,T]}
					\E[C_1 + V(X^{0,x}_t)] 
			\leq 
				C_1 + V(R) e^{\theta T}
			< \infty.
		\end{split}
	\end{equation}
	Moreover,
	combining \eqref{eq: V bound for X0}
	and \eqref{eq: V bound for phi} yields that for all $\delta \in (0, \infty)$
	it holds that 
	\begin{align}
	\nonumber
				&\lim_{ r \to \infty }
				\sup_{t \in [0,T]}
				\sup \Big \{
					\frac
						{| u_0(t, x) |}
						{V(x)}
					\colon x \in H, ~ \|x\|_H \geq r
				\Big \} \\ \nonumber
			\leq{}
				&\lim_{ r \to \infty }
				\sup_{t \in [0,T]}
				\sup \Big \{
					\frac
						{\E[|\varphi(X^{0,x}_t)|] }
						{V(x)}
					\colon x \in H, ~ \|x\|_H \geq r
				\Big \} \\
			\leq{} 
				&\lim_{ r \to \infty }
				\sup_{t \in [0,T]}
				\sup \Big \{
					\frac
						{\E[C_\delta + \delta V(X^{0,x}_t)]}
						{V(x)} 
					\colon x \in H, ~ \|x\|_H \geq r
				\Big \} \\ \nonumber
			\leq{}
				&\lim_{ r \to \infty }
				\sup_{t \in [0,T]}
				\sup \Big \{
					\frac
						{C_\delta + \delta V(x) e^{\theta t}}
						{V(x)} 
					\colon x \in H, ~ \|x\|_H \geq r
				\Big \} \\ \nonumber
			={} 
				&\delta e^{\theta T}
  \end{align}
	and thus it holds that
	\begin{equation}
			\lim_{ r \to \infty }
			\sup_{t \in [0,T]}
			\sup \Big \{
				\frac
					{| u_0(t, x) |}
					{V(x)}
				\colon x \in H, \|x\|_{H} \geq r
			\Big \}
		= 0.
	\end{equation}
	In addition, the
	continuity of $\varphi$
	with respect to the $\| \cdot \|_{H}$-norm
	on $\mathbb{H}$-bounded subsets of $H$,
	the
	continuity of $u_0$
	with respect to the $\| \cdot \|_{\R \times H}$-norm
	on $\R \times \mathbb{H}$-bounded subsets of $[0,T] \times H$,
	the fact that for all $R \in (0,\infty)$ it holds that
	the set $\{(t,x) \in [0,T] \times H \colon \|x\|_{H_\vartheta} \leq R\}$
	is compact
	and the fact that $V$ is bounded on $\mathbb{H}$-bounded subsets of $H$
	imply for all $R, \tilde{R} \in (0, \infty)$ that
	\begin{align}
				&\lim_{r \downarrow 0} \lim_{\eps \downarrow 0}	
				\sup \bigg \{
					u_0(t,x) -  \varphi(\hat{x})
					\colon
					x, \hat{x} \in H_{2\varz},
					~t, \hat{t} \in (0,T),
					~\tfrac{e^{-\theta t} \|x\|_{H_{\varz}}^2}{V(x)} 
						\vee \tfrac{e^{-\theta \hat{t}} 
							\|\hat{x}\|_{H_{\varz}}^2}{V(\hat{x})} \leq R, 
				\\ \nonumber & \qquad \qquad \qquad \qquad
					~\|x \|_{H} \leq \tilde{R}, 
					~\| x - \hat{x}\|_{H} \leq r,
					~t \vee \hat{t} \leq \eps 
				\bigg \} \\ \nonumber
			={} 
				&\lim_{r \downarrow 0}
				\sup \bigg \{
					u_0(0,x) -  \varphi(\hat{x})
					\colon
					x, \hat{x} \in H_{2\varz},
					~\tfrac {\|x\|_{H_{\varz}}^2} {V(x)} 
						\vee \tfrac {\|\hat{x}\|_{H_{\varz}}^2} {V(\hat{x})} \leq R,
					~\|x \|_{H} \leq \tilde{R},
					~\| x - \hat{x}\|_{H} \leq r
				\bigg \} \\ \nonumber
			={}  
				&\lim_{r \downarrow 0}
				\sup \bigg \{
					\varphi(x) -  \varphi(\hat{x})
					\colon
					x, \hat{x} \in H_{2\varz},
					~\tfrac {\|x\|_{H_{\varz}}^2} {V(x)} 
						\vee \tfrac {\|\hat{x}\|_{H_{\varz}}^2} {V(\hat{x})} \leq R,
					~\|x \|_{H} \leq \tilde{R},
					~\| x - \hat{x}\|_{H} \leq r
				\bigg \} \\ \nonumber
			={}
				&0
	\end{align}
	and this finishes the proof of Theorem \ref{thm: existence}.
\end{proof}
\chapter{Regularity estimates}
\label{sec: Auxillary results}
In this chapter we establish some auxiliary results, which we need for the next chapter.
In particular in Section \ref{ssec: holder sob inequalities}
we estimate the product of two functions in a Besov space.
More precisely, given $d \in \N$,
$\alpha, \beta \in \N_0^d$, $p \in (1,2]$,
$s \in [0,\infty)$,
a Lipschitz domain $\Omega \subseteq \R^d$, 
and a suitable family of
Hilbert spaces $H_r$, $r \in \R$, having the interpolation
properties  (e.g.\@ Sobolev spaces, see
also Setting \ref{ssec: application})
we establish sufficient conditions for
 $r_1, r_2, \tilde{r}_1, \tilde{r}_2 \in \R$
such that there exists a $C \in (0,\infty)$
such that for all $u \in H_{r_1}$
and all $v \in H_{\tilde{r}_1} \cap H_{\tilde{r}_2}$ it holds that
\begin{equation}
\label{eq: sobolev holder type inequality}
		\| 
			(D_{\R^d}^{\alpha}u)  (D_{\R^d}^{\beta}v)
		\|_{B^s_{p,2}(\Omega)} 
	\leq 
		C(
			\|u\|_{H_{r_1}} \|v\|_{H_{\tilde{r}_2}} 	
			+ \|u\|_{H_{r_2}} \|v\|_{H_{\tilde{r}_1}}
		),
\end{equation}
where
$B^s_{p,2}(\Omega)$ denotes the Besov space (see Notation \ref{sec: notation chap 5}).
In Section \ref{ssec: Estimates for Nemytskij operator}
we establish estimates for compositions of functions.
Lemma \ref{lem: b boundedness inequality}
and Lemma \ref{lem: b continuity inequality} 
give bounds on the 
Slobodeckij semi-norm,
Lemma \ref{lem: b L2 boundedness inequality} 
and Lemma \ref{lem: L2 b continuity inequality}
on the $L^2$-norm,
and Corollary \ref{cor: Bs b bound}
and Corollary \ref{cor: Bs b continuity inequality} 
combine these results.

Throughout this chapter the following notation is used.
\section{Notation}
\label{sec: notation chap 5}
	Let $| \cdot | \colon \bigcup^{\infty}_{d=1} \R^d \to \R$
	be the function satisfying for all 
	$d \in \N \backslash \{ 1\}$,
	$x \in \R$, 
	and all $(y_1, \cdots , y_d) \in \R^d$ that
	\begin{equation}
			| x| 
		= 
			\begin{cases}
				x & \textrm{ if } x \geq 0  \\
				-x & \textrm{ if } x < 0
			\end{cases}
	\qquad \textrm{ and that } \quad
			| y| 
		= 
			\sum_{i=1}^d
				|y_i|.
	\end{equation}
	For $d \in \N$ we say that $\Omega \subseteq \R^d$ is a Lipschitz set
	if $\Omega$ is either $\R^d$,
	the graph above a Lipschitz function, or a bounded Lipschitz domain.
	Moreover, for $d,n \in \N$
	the d-dimensional Lebesgue measure $\mu^d$, 
	and an Lipschitz set $\Omega \subseteq \R^d$,
	we denote by 
	$L^0( \Omega )$ the set of all Lebesgue measurable functions,
	by $\mathbb{B}^{l}_{p,q}(\Omega)$, $l \in \R$, $p \in [1, \infty]$, $q \in [1, \infty]$
	the Besov spaces (see, e.g.,  Definition 2 on page 8
	in Runst \& Sickel \cite{RunstSickel1996}
	and Definition 1.1 and Section 1.2 in Rychkov \cite{Rychkov1998})
	satisfying for all $l \in \R$
	and all
	$(p,q) \in [1, \infty]^2 \backslash \{2,2\}$
	that
	$
			\mathbb{B}^{l}_{p,q}(\Omega)
		= 
			(
				B^{l}_{p,q}(\Omega),
				\| \cdot \|_{B^{l}_{p,q}(\Omega)}
			),
	$
	and that
	$
			\mathbb{B}^{l}_{2,2}(\Omega)
		= 
			(
				B^{l}_{2,2}(\Omega),
				\langle \cdot, \cdot \rangle_{B^{l}_{2,2}(\Omega)},
				\| \cdot \|_{B^{l}_{2,2}(\Omega)}
			)
	$, 
	by
	$\mathbb{L}^{p}(\Omega)$, $p \in [1, \infty]$,
	the Banach spaces satisfying for all
	$p \in [1, \infty] \backslash \{2\}$
	that
	$
			\mathbb{L}^{p}(\Omega)
		= 
			(
				L^{p}(\Omega),
				\| \cdot \|_{L^{p}(\Omega)}
			)
		= 
			(
				L^{p}(\mu^d,\Omega),
				\| \cdot \|_{L^{p}(\mu^d,\Omega)}
			)
	$
	and that
	$
			\mathbb{L}^{2}(\Omega)
		= 
			(
				L^{2}(\Omega),
				\langle \cdot, \cdot \rangle_{L^{2}(\Omega)},
				\| \cdot \|_{L^{2}(\Omega)}
			)
		= 
			(
				L^{2}(\mu^d,\Omega),
				\langle \cdot, \cdot \rangle_{L^{2}(\mu^d,\Omega)},
				\| \cdot \|_{L^{2}(\mu^d,\Omega)}
			),
	$
	by  
	$\mathbb{B}^{l}_{p,q}(\Omega, \R^n)$, 
	$l \in \R$, $p \in [1, \infty]$, $q \in [1, \infty]$
	the Besov spaces satisfying for all
	$l \in \R$, $p \in [1, \infty]$, and all $q \in [1, \infty]$ that
	$
			\mathbb{B}^{l}_{p,q}(\Omega, \R^n) 
		= 
			\otimes_{i=1}^n \mathbb{B}^{l}_{p,q}(\Omega),
	$
	and by $\mathbb{L}^{p}(\Omega,\R^n)$, $p \in [1, \infty]$,
	the Banach spaces satisfying for all
	$p \in [1, \infty]$
	that
	$
			\mathbb{L}^{p}(\Omega, \R^n) 
		= 
			\otimes_{i=1}^n \mathbb{L}^{p}(\Omega).
	$	
\section{Hölder-Sobolev type inequalities}
\label{ssec: holder sob inequalities}
Throughout this subsection the following setting
	is frequently used.
\begin{sett}
	\label{ssec: application}
	Let $d \in \N$, 
	let $\Omega \subseteq \R^d$ be a Lipschitz domain,
	let $K \colon \R \to (0,\infty)$,
	let 
	$
			\mathbb{H}_r
		=
			(H_r, \langle \cdot, \cdot \rangle_{H_r}, \| \cdot \|_{H_r})
	$,
	$r \in \R$, be a family of Hilbert spaces satisfying for all 
	$r \in [0,\infty)$, $s \in \R$, and all $t\in [s,\infty)$ that
	$
		H_r \subseteq B^{2r}_{2,2}(\Omega)
	$,
	that
	the $\| \cdot \|_{H_r}$-norm is equivalent
	to the $\| \cdot \|_{B^{2r}_{2,2}(\Omega)}$-norm on $H_r$,
	and that
	$H_t \subseteq H_s$
	and satisfying that	for all $r, r_1, r_2 \in \R$, $s_1, s_2 \in [0,1]$ 
	and for all $x \in H_r$ with
	$s_1 + s_2 =1,$ and with $s_1 r_1 + s_2 r_2 =r$
	it holds that
	\begin{equation}
	\label{eq: interpolation inequality eq}
			\|x\|_{H_r}
		\leq
			(K_{r_1}  \|x\|_{H_{r_1}})^{s_1} \cdot (K_{r_2} \|x\|_{H_{r_2}})^{s_2}.
	\end{equation}
	Let $(H, \langle \cdot, \cdot \rangle_{H}, \| \cdot \|_{H})=\mathbb{H}_0$,
	and let $\mathbb{H}'_r=(H'_r, \| \cdot \|_{H'_r})$,
	$r \in \R$,
	be the corresponding dual spaces.
	By abuse of notation we will also denote by
		$\| \cdot \|_{H_r}$, $r \in \R$,
		the extended functions 
		$
			\| \cdot \|_{H_r} \colon \bigcup_{i=1}^{\infty} H_{-i} \to [0, \infty]
		$,
		$r \in \R$,
		satisfying for all
		$r \in \R$, $x \in H_{r}$,
		and all $y \in \bigcup_{i=1}^{\infty} H_{-i}$
		that
		\begin{equation}
					\|y\|_{H_r}
				=
					\begin{cases} 
						\|y\|_{H_{r}} 
							& \textrm{ if } y \in H_{r} \\
						\infty 
							& \textrm{ if } y \notin H_{r}.
					\end{cases}
		\end{equation}
\end{sett}
\begin{lemma}
	\label{l: norm equivalence with interpolation}
	Let $d \in \N$, 
	let $\Omega \subseteq \R^d$ be a Lipschitz domain,
	let 
	$
			\mathbb{H} 
		=
			(H, \langle \cdot, \cdot \rangle \rangle_H, \|\cdot\|_H)
	$
	be a separable Hilbert space,
	let $e_n$, $n\in \N$, be an orthonormal basis of $\H$,
	let $\lambda \colon \N \to (- \infty, 0)$ be a function with 
	$\lim_{n \to \infty} \lambda_n = - \infty$,
	let $A \colon D(A) \subseteq H \to H$ 
	be the linear operator such that
	$
			D(A) 
		= 
			\{ 
				x \in H \colon 
					\sum_{n =1}^\infty |\lambda_n \langle e_n, x \rangle_{H}|^2 < \infty 
			\}
	$
	and such that for all $x \in D(A)$ it holds that
	$
		Av = \sum_{n =1}^\infty \lambda_n \langle e_n, x \rangle_{H} e_n
	$,
	let 
	$ 
			\mathbb{H}_{t}
		=
			( 
				H_{t} , 
				\left< \cdot , \cdot \right>_{ H_{t} }, 
				\left\| \cdot \right\|_{ H_{t} } 
			) 
	$,
	$ t \in \R $,
	be a family of interpolation spaces associated with
	$ -A $ 
	(see, e.g., Definition~3.6.30 in Jentzen \cite{Jentzen2015}),
	let $r_1\in [0,\infty)$, $r_2 \in [r_1,\infty)$, 
	$r \in [r_1,r_2]$, $\alpha \in (0,\infty)$, satisfy that
	$H_{r_1} \subseteq B^{\alpha r_1}_{2,2}(\Omega)$
	and that $H_{r_2} \subseteq B^{\alpha r_2}_{2,2}(\Omega)$,
	and let $C_1, C_2 \in (0, \infty)$
	satisfy that for all $u \in H_{r_1}$
	it holds that
	\begin{equation}
	\label{eq: norm equivalence r1}
			C_1 \|u\|_{B^{\alpha r_1}_{2,2}(\Omega)} 
		\leq 
			\|u\|_{H_{r_1}}
		\leq
			C_2 \|u\|_{B^{\alpha r_1}_{2,2}(\Omega)} 
	\end{equation}
	and that for all $u \in H_{r_2}$
	it holds that
	\begin{equation}
	\label{eq: norm equivalence r2}
		C_1 \|u\|_{B^{\alpha r_2}_{2,2}(\Omega)} 
		\leq 
			\|u\|_{H_{r_2}}
		\leq
			C_2 \|u\|_{B^{\alpha r_2}_{2,2}(\Omega)}.
	\end{equation}
	Then 
	there exist $C_3, C_4 \in (0,\infty)$
	such that for all $u \in H_{r}$ 
	it holds that
	\begin{equation}
			C_3 \|u\|_{B^{\alpha r}_{2,2}(\Omega)} 
		\leq 
			\|u\|_{H_{r}}
		\leq
			C_4 \|u\|_{B^{\alpha r}_{2,2}(\Omega)}.
	\end{equation}
\end{lemma}
\begin{proof}
	Throughout this proof we will use the notation in Section 1.3.2
	on page 24 in Triebel \cite{Triebel1978}.
	For $r \in \{r_1,r_2\}$ the assertion is trivial.
	Therefore assume for the rest of the proof that $r \in (r_1, r_2)$.
	Moreover,
	denote
	by $\theta \in (0,1)$ the value satisfying that
	$r = (1-\theta) r_1 + \theta r_2$.
	Next note that the theorem on page 182 in Triebel \cite{Triebel1978} yields
	that
	\begin{equation}
			(B^{\alpha r_1}_{2,2}(\R^d), B^{\alpha r_2}_{2,2}(\R^d))_{\theta, 2}
		=
			B^{\alpha r}_{2,2}(\R^d).
	\end{equation}
	Combining this with the Theorem on page 25 in Triebel \cite{Triebel1978}
	and with the fact that 
	there exist a
	linear bounded operator 
	$\xi \in L(\mathbb{B}^{\alpha r_1}_{2,2}(\Omega), \mathbb{B}^{\alpha r_1}_{2,2}(\R^d))$
	such that
	$
			\xi|_{B^{\alpha r_2}_{2,2}(\Omega)} 
		\in L(\mathbb{B}^{\alpha r_2}_{2,2}(\Omega), \mathbb{B}^{\alpha r_2}_{2,2}(\R^d)),
	$
	that
	$
			\xi|_{B^{\alpha r}_{2,2}(\Omega)} 
		\in L(\mathbb{B}^{\alpha r}_{2,2}(\Omega), \mathbb{B}^{\alpha r}_{2,2}(\R^d))
	$
	and that for all $u \in \mathbb{B}^{\alpha r_1}_{2,2}(\Omega)$
	it holds that
	$(\xi (u))|_{\Omega} = u$
	(see, e.g., Theorem 4.1 in Rychkov \cite{Rychkov1998})
	shows that
	\begin{equation}
		\begin{split}
				&B^{\alpha r}_{2,2}(\Omega)
			=
				(
					B^{\alpha r}_{2,2}(\R^d)
				)|_{\Omega}
			=
				\big(
						(B^{\alpha r_1}_{2,2}(\R^d), B^{\alpha r_2}_{2,2}(\R^d))_{\theta, 2}
					\big) \big|_{\Omega} 
			\subseteq
				(B^{\alpha r_1}_{2,2}(\Omega), B^{\alpha r_2}_{2,2}(\Omega))_{\theta, 2} \\
			={} 
				& 
				\Big(
					\xi \big(
						(B^{\alpha r_1}_{2,2}(\Omega), B^{\alpha r_2}_{2,2}(\Omega))_{\theta, 2}
					\big)
				\Big) \Big|_{\Omega} 
			\subseteq
				\big(
					(B^{\alpha r_1}_{2,2}(\R^d), B^{\alpha r_2}_{2,2}(\R^d))_{\theta, 2}
				\big) \big|_{\Omega}
			=
				(
					B^{\alpha r}_{2,2}(\R^d)
				)|_{\Omega}
			=
				B^{\alpha r}_{2,2}(\Omega)
		\end{split}
	\end{equation}
	and that there exist $K_1, K_2, K_3, K_4 \in (0,\infty)$
	such that
	\begin{equation}
	\label{eq: inter equiv to B}
		\begin{split} 
				&\|
					u
				\|_{B^{\alpha r}_{2,2}(\Omega)}
			\leq
				\|
					\xi u
				\|_{B^{\alpha r}_{2,2}(\R^d)} 
			\leq
				K_1 \|
					\xi u
				\|_{(B^{\alpha r_1}_{2,2}(\R^d), B^{\alpha r_2}_{2,2}(\R^d))_{\theta, 2}} \\
			\leq{}
				&K_2
				\|u\|_{(B^{\alpha r_1}_{2,2}(\Omega), B^{\alpha r_2}_{2,2}(\Omega))_{\theta, 2}}
			={}  
				K_2
				\|
					(\xi u) |_{\Omega}
				\|_{(B^{\alpha r_1}_{2,2}(\Omega), B^{\alpha r_2}_{2,2}(\Omega))_{\theta, 2}} \\
			\leq{}
				&K_2
				\|
					\xi u
				\|_{(B^{\alpha r_1}_{2,2}(\R^d), B^{\alpha r_2}_{2,2}(\R^d))_{\theta, 2}} 
			\leq{}
				K_3 \|
					\xi u
				\|_{B^{\alpha r}_{2,2}(\R^d)} 
			\leq  
				K_4
				\|
					u
				\|_{B^{\alpha r}_{2,2}(\Omega)}.
		\end{split}
	\end{equation}
	Furthermore, the theorem on page 141 in Triebel \cite{Triebel1978}
	ensures that
	\begin{equation}
			(H_{r_1}, H_{r_2})_{\theta, 2}
		=
			H_r
	\end{equation}
	and that there exist $K_1, K_2 \in (0,\infty)$
	such that for all $u \in H_r$ it holds that
	\begin{equation}
	\label{eq: inter equiv to H}
			K_1 \| u\|_{H_r}
		\leq
			\|u\|_{(H_{r_1}, H_{r_2})_{\theta, 2}}
		\leq
			K_2 \|u\|_{H_r}.
	\end{equation}
	Moreover, it follows from \eqref{eq: norm equivalence r1}
	and from \eqref{eq: norm equivalence r2}
	that for all $i \in \{1,2\}$ it holds that
	$
			\id_{H_{r_i}} 
		\in 
			L(
				(
					H_{r_i}, 
					\langle \cdot, \cdot \rangle_{{B}^{\alpha r_i}_{2,2}}, 
					\|\cdot \|_{{B}^{\alpha r_i}_{2,2}}
				), 
				\mathbb{H}_{r_i}
			)
	$
	and that
	$\id_{H_{r_i}} \in L(\mathbb{H}_{r_i}, \mathbb{B}^{\alpha r_i}_{2,2})$
	and combining this with \eqref{eq: inter equiv to B},
	the theorem on page 25 in Triebel \cite{Triebel1978}
	and with \eqref{eq: inter equiv to H}
	verifies
	that there exist $K_1, K_2, K_3, K_4 \in (0,\infty)$
	such that for all $u \in H_r$ it holds that
	\begin{equation}
		\begin{split}
				&C_1 \|u\|_{B^{\alpha r}_{2,2}(\Omega)} 
			\leq 
				C_1 K_1 
				\|
					\id_{H_{r_1}} u
				\|_{(B^{\alpha r_1}_{2,2}(\Omega), B^{\alpha r_2}_{2,2}(\Omega))_{\theta, 2}}
			\leq
				K_1 \|u\|_{(H_{r_1}, H_{r_2})_{\theta, 2}}
			\leq
				K_2 \|u\|_{H_{r}} \\
			\leq{}  & 
				K_3 \| \id_{H_{r_1}} u\|_{(H_{r_1}, H_{r_2})_{\theta, 2}} 
			\leq 
				C_2 K_3
				\|
					u
				\|_{(B^{\alpha r_1}_{2,2}(\Omega), B^{\alpha r_2}_{2,2}(\Omega))_{\theta, 2}}
			\leq
				K_4 \|u\|_{B^{\alpha r}_{2,2}(\Omega)}.
		\end{split}
	\end{equation}
	This finishes the proof of Lemma \ref{l: norm equivalence with interpolation}.
\end{proof}
Before we come to the main results of this subsection
we give two examples of families of Hilbert spaces
which satisfy the
interpolation property.
\begin{lemma}[Interpolation results for H] 
\label{lem: lemma: interpolation theorem for H}
		Let 
	$
			\mathbb{H} 
		=
			(H, \langle \cdot, \cdot \rangle \rangle_H, \|\cdot\|_H)
	$
	be a separable Hilbert space,
	let $e_n$, $n\in \N$, be an orthonormal basis of $\mathbb{H}$,
	let $\lambda \colon \N \to (- \infty, 0)$ be a function with 
	$\lim_{n \to \infty} \lambda_n = - \infty$,
	let $A \colon D(A) \subseteq H \to H$ 
	be the linear operator such that
	$
			D(A) 
		= 
			\{ 
				x \in H \colon 
					\sum_{n =1}^\infty |\lambda_n \langle e_n, x \rangle_{H}|^2 < \infty 
			\}
	$
	and such that for all $x \in D(A)$ it holds that
	$
		Av = \sum_{n =1}^\infty \lambda_n \langle e_n, x \rangle_{H} e_n
	$,
	let 
	$ 
			\mathbb{H}_{t}
		=
			( 
				H_{t} , 
				\left< \cdot , \cdot \right>_{ H_{t} }, 
				\left\| \cdot \right\|_{ H_{t} } 
			) 
	$,
	$ t \in \R $,
	be a family of interpolation spaces associated with
	$ -A $ 
	(see, e.g., Definition~3.6.30 in Jentzen \cite{Jentzen2015}).
	By abuse of notation we will also denote by
		$A$ and  by $\| \cdot \|_{H_t}$, $t \in \R$,
		the extended operators 
		$
			A \colon \bigcup_{i=1}^{\infty} H_{-i} \to \bigcup_{i=1}^{\infty} H_{-i}
		$ 
		and 
		$
			\| \cdot \|_{H_t} \colon \bigcup_{i=1}^{\infty} H_{-i} \to [0, \infty]
		$,
		$t \in \R$,
		satisfying for all
		$t \in \R$, $x \in H_{t}$, 
		and all $y\in \bigcup_{i=1}^{\infty} H_{-i}$
		that
		\begin{equation}
					\|y\|_{H_t}
				=
					\begin{cases} 
						\|y\|_{H_{t}} 
							& \textrm{ if } y \in H_{t} \\
						\infty 
							& \textrm{ if } y \notin H_{t}
					\end{cases}
		\end{equation}
		and that
		\begin{equation}
				(A(x) = y) 
			\Leftrightarrow 
				(
						\lim_{\eps \downarrow 0} \sup \{
							\|A(\xi) - y\|_{H_{t-1}} \colon \xi \in H_1, ~\|x- \xi\|_{H_{t}} \leq \eps 
						\}
					=
						0
				).
		\end{equation}
		Let $r, r_1, r_2 \in \R$, $s_1, s_2 \in [0,1]$, satisfy that
		$s_1 + s_2 =1,$ and that $s_1 r_1 + s_2 r_2 =r$.
	Then it holds for all $x \in H_r$ that
	\begin{equation}
	\label{eq: H interpolation inequality lem}
			\|x\|_{H_r}
		\leq
			\|x\|^{s_1}_{H_{r_1}} \cdot \|x\|^{s_2}_{H_{r_2}}.
	\end{equation}
\end{lemma}
\begin{proof}
	Observe that \eqref{eq: H interpolation inequality lem} is trivial
	for $r_1 = r_2$.
	Moreover it follows from the fact that $H_{r_2} \cap H$ is dense in $H_r$
	that it is enough to show the equation for $x \in H_{r_2} \cap H$.
	Therefore for the rest of the proof fix $x \in H_{r_2} \cap H$
	and assume that $r_1 < r_2$.
 	Now denote by $S \subseteq [0,1]$ the set satisfying that
	\begin{equation}
		\begin{split}
				S
			= 
				\Big \{
					s_3 \in [0,1] \colon 
					\Big(
						&\exists \tilde{s}_1, \tilde{s}_2 \in [0,1],
						\exists r_3 \in [r_1, r_2] \colon 
						\big(
							(
									\|x\|_{H_r}
								\leq
									\|x\|^{\tilde{s}_1}_{H_{r_1}} \cdot \|x\|^{\tilde{s}_2}_{H_{r_2}} 
									\cdot \|x\|^{s_3}_{H_{r_3}}
							) \\
							& \qquad \qquad \qquad
							\wedge (\tilde{s}_1 + \tilde{s}_2 +s_3 =1) 
							\wedge (\tilde{s}_1 r_1 + \tilde{s}_2 r_2 + s_3 r_3 =r)	
						\big)
					\Big)
				\Big \}.
		\end{split}
	\end{equation}
	As $s_1$, $s_2$ are uniquely determined by the equation
	$s_1 +s_2 =1$ and the equation $s_1r_1 + s_2 r_2 =r$ 
	it is sufficient to show $0 \in S$. 
	Note that 
	$
			\|x\|_{H_r}
		\leq
			\|x\|^{0}_{H_{r_1}} \cdot \|x\|^{0}_{H_{r_2}} 
			\cdot \|x\|^{1}_{H_{r}}
	$
	and hence it holds that
	$1 \in S$.
	This ensures that
	$S \neq \emptyset$. Next denote by $s \in [0,1]$ the value satisfying that
	$s = \inf S$. 
	This implies that there exist a sequence $(s^{(n)})_{n \in \N} \subseteq S$
	such that $\lim_{n \to \infty} s^{(n)} = s$. Moreover,
	it follows from the definition of $S$ that there exist sequences
	$(s_1^{(n)})_{n \in \N},(s_2^{(n)})_{n \in \N} \subseteq [0, 1]$,
	$(r^{(n)})_{n \in \N} \subseteq [r_1, r_2]$ such that for all
	$n \in \N$ it holds that
	$
			\|x\|_{H_r}
		\leq
			\|x\|^{s^{(n)}_1}_{H_{r_1}} \cdot \|x\|^{s^{(n)}_2}_{H_{r_2}} 
			\cdot \|x\|^{s^{(n)}}_{H_{r^{(n)}}}
	$,
	$s^{(n)}_1 + s^{(n)}_2 + s^{(n)}=1 $,
	and that
	$ 
		s^{(n)}_1 r_1 + s^{(n)}_2 r_2 + s^{(n)} r^{(n)} =r 
	$.
	By going to subsequences we can assume that there exist
	$\tilde{s}_1, \tilde{s}_2 \in [0,1]$, $r_3 \in [r_1, r_2]$ such that
	$ 
			\lim_{n \to \infty}
				(|s^{(n)}_1- \tilde{s}_1 |+ | s^{(n)}_2 - \tilde{s}_2| + |r^{(n)} - r_3|)
		= 0
	$.
	Thus we get from $x \in H_{r_2}$ 
	that
	\begin{equation}
		\begin{split}
					&\|x\|_{H_r}
				\leq
					\lim_{n \to \infty} 
					\big(
						\|x\|^{s^{(n)}_1}_{H_{r_1}} \cdot \|x\|^{s^{(n)}_2}_{H_{r_2}} 
						\cdot \|x\|^{s^{(n)}}_{H_{r^{(n)}}}
					\big)
				=
					\|x\|^{\tilde{s}_1}_{H_{r_1}} \cdot \|x\|^{\tilde{s}_2}_{H_{r_2}} 
						\cdot \|x\|^{s}_{H_{r_3}}, \\
				&\tilde{s}_1 + \tilde{s}_2 + s
			= 
				\lim_{n \to \infty} (s^{(n)}_1 + s^{(n)}_2 + s^{(n)})
			= 1, \quad \textrm{ and that} \\
				&\tilde{s}_1 r_1 + \tilde{s}_2 r_2 + s r_3
			=
				\lim_{n \to \infty} (s^{(n)}_1 r_1 + s^{(n)}_2 r_2 + s^{(n)} r^{(n)})
			=r
		\end{split}
	\end{equation}
	and therefore we obtain that $s \in S$.
	We will now show that $\tfrac s2 \in S$ and therefore
	we will consider two cases. \\
	{\it 1. case: $r_3 \geq \tfrac{r_1+r_2}{2}$.}
	First observe that we have
	that
	\begin{equation}
		\begin{split}
				&\|x\|_{H_r}
			\leq 
				\|x\|^{\tilde{s}_1}_{H_{r_1}} \cdot \|x\|^{\tilde{s}_2}_{H_{r_2}} 
				\cdot \|x\|^{s}_{H_{r_3}}
			=
				\|x\|^{\tilde{s}_1}_{H_{r_1}} \cdot \|x\|^{\tilde{s}_2}_{H_{r_2}} 
				\cdot (\langle 
					(-A)^{2r_3-r_2} x, (-A)^{r_2} x 
				\rangle_{H})^{\nicefrac s2}\\
			\leq{} 
				&\|x\|^{\tilde{s}_1}_{H_{r_1}} \cdot \|x\|^{\tilde{s}_2}_{H_{r_2}} 
				\cdot (\|x\|_{H_{2r_3-r_2}}
					\|x\|_{H_{r_2}})^{\nicefrac s2}
			=
				\|x\|^{\tilde{s}_1}_{H_{r_1}} \cdot \|x\|^{\tilde{s}_2 +\nicefrac s2}_{H_{r_2}} 
				\cdot \|x\|^{\nicefrac s2}_{H_{2r_3-r_2}}.
		\end{split}
	\end{equation}
	In addition, it follows from $r_3 \in [\tfrac{r_1+r_2}{2}, r_2]$ that
	$2r_3-r_2 \in [r_1, r_2]$.
	Furthermore, we obtain that
	$
			r
		=
			\tilde{s}_1 r_1 + \tilde{s}_2 r_2 + s r_3
		=
			\tilde{s}_1 r_1 + (\tilde{s}_2+ \tfrac s2) r_2 + \tfrac s2 (2r_3-r_2)
	$
	and this shows that $\tfrac s2 \in S$ and thus verifies the first case.\\
	{\it 2. case: $r_3 < \tfrac{r_1+r_2}{2}$.}
	We then get 
	that
	\begin{equation}
		\begin{split}
				&\|x\|_{H_r}
			\leq 
				\|x\|^{\tilde{s}_1}_{H_{r_1}} \cdot \|x\|^{\tilde{s}_2}_{H_{r_2}} 
				\cdot \|x\|^{s}_{H_{r_3}}
			=
				\|x\|^{\tilde{s}_1}_{H_{r_1}} \cdot \|x\|^{\tilde{s}_2}_{H_{r_2}} 
				\cdot (\langle 
					(-A)^{r_1} x, (-A)^{2r_3-r_1}  x 
				\rangle_{H})^{\nicefrac s2}\\
			\leq{} 
				&\|x\|^{\tilde{s}_1}_{H_{r_1}} \cdot \|x\|^{\tilde{s}_2}_{H_{r_2}} 
				\cdot (\|x\|_{H_{r_1}}
					\|x\|_{H_{2r_3-r_1}})^{\nicefrac s2}
			=
				\|x\|^{\tilde{s}_1+\nicefrac s2}_{H_{r_1}} \cdot \|x\|^{\tilde{s}_2}_{H_{r_2}} 
				\cdot \|x\|^{\nicefrac s2}_{H_{2r_3-r_1}}.
		\end{split}
	\end{equation}
	Moreover, the fact that $r_3 \in [r_1, \tfrac{r_1+r_2}{2}]$ yields that
	$2r_3-r_1\in [r_1, r_2]$.
	In addition, we have that
	$
			r
		=
			\tilde{s}_1 r_1 + \tilde{s}_2 r_2 + s r_3
		=
			(\tilde{s}_1 + \tfrac s2) r_1 + \tilde{s}_2 r_2 + \tfrac s2 (2r_3-r_1)
	$
	and this shows that $\tfrac s2 \in S$ and finishes the second case.\\
	Combining the two cases then ensures that $\tfrac s2 \in S$.
	Hence we derive from the fact
	that $s = \inf S$ that $s \leq \tfrac s2$ and therefore $s=0$.
	This finishes the proof of Lemma \ref{lem: lemma: interpolation theorem for H}.
\end{proof}
\begin{lemma}[Interpolation results for Sobolev spaces] 
\label{lem: lemma: interpolation theorem for B}
		Let $d \in \N$,
		let $\Omega \subseteq \R^d$ be a Lipschitz domain.
		Then there exist a $K \colon \R \to [0,\infty)$ such that for all
		$r, r_1, r_2 \in \R$, $s_1, s_2 \in [0,1]$, 
		and all $u \in B^r_{2,2}(\Omega)$ 
		with
		$s_1 + s_2 =1,$ and with $s_1 r_1 + s_2 r_2 =r$
		it holds that
	\begin{equation}
	\label{eq: B interpolation inequality lem}
			\|u\|_{B^r_{2,2}(\Omega)}
		\leq
			(K_{r_1} \|u\|_{B^{r_1}_{2,2}(\Omega)})^{s_1}
			\cdot (K_{r_2} \|u\|_{B^{r_2}_{2,2}(\Omega)})^{s_2}.
	\end{equation}
\end{lemma}
\begin{proof}
	First denote by 
	$	
		\xi \colon \cup^\infty_{r=-\infty} \mathbb{B}^{r}_{2,2}(\Omega)
			\to \cup^\infty_{r=-\infty} \mathbb{B}^{r}_{2,2}(\R^d)
	$
	an operator satisfying that for all $r \in \R$ it holds that
	that
	$
			\xi|_{B^{r}_{2,2}(\Omega)} 
		\in L(\mathbb{B}^{r}_{2,2}(\Omega), \mathbb{B}^{r}_{2,2}(\R^d))
	$,
	and that for all $u \in \cup^\infty_{r=-\infty} \mathbb{B}^{r}_{2,2}(\Omega)$
	it holds that
	$(\xi (u))|_{\Omega} = u$
	(the existence follows from Theorem 4.1 in Rychkov \cite{Rychkov1998}) 
	and by $K \colon \R \to [0,\infty)$ the function satisfying for all $r\in \R$
	that
	\begin{equation}
			K_r
		=
			\|\xi \|_{
				L(\mathbb{B}^{r}_{2,2}(\Omega), \mathbb{B}^{r}_{2,2}(\R^d))
			}
			\vee 1.
	\end{equation}
	Next note that 
	Remark 2 on page 87 in Runst \& Sickel \cite{RunstSickel1996}
		(with $s \leftarrow r$, $s_0 \leftarrow r_1$, $s_1 \leftarrow r_2$,
		$\theta \leftarrow s_2$, 
		$p \leftarrow 2$, $p_0 \leftarrow 2$, $p_1 \leftarrow 2$,
		$q \leftarrow 2$, $q_0 \leftarrow 2$, $q_1 \leftarrow 2$)
	implies that for all
		$r, r_1, r_2 \in \R$, $s_1, s_2 \in (0,1)$, 
		and all $u \in B^r_{2,2}(\Omega)$ 
		with
		$s_1 + s_2 =1$, $r_1 \neq r_2$, and with $s_1 r_1 + s_2 r_2 =r$
	it holds that
	\begin{equation}
				\|u\|_{B^r_{2,2}(\Omega)}
			\leq
				\|\xi u\|_{B^r_{2,2}(\R^d)}
			\leq
				\|\xi u\|^{s_1}_{B^{r_1}_{2,2}(\R^d)} 
				\cdot \|\xi u\|^{s_2}_{B^{r_2}_{2,2}(\R^d)} 
			\leq{}
				( 
					K_{r_1}
					\| u\|_{B^{r_1}_{2,2}(\Omega)}
				)^{s_1}
				\cdot (
					K_{r_2}
					\|u\|_{B^{r_2}_{2,2}(\Omega)}
				)^{s_2}.
	\end{equation}
	Thus, the fact that for all $r \in \R$ it holds that $K_r \geq 1$
	ensures that for all
		$r, r_1, r_2 \in \R$, $s_1, s_2 \in [0,1]$, 
		and all $u \in B^r_{2,2}(\Omega)$ 
		with
		$s_1 + s_2 =1$ and with $s_1 r_1 + s_2 r_2 =r$
	it holds that
	\begin{equation}
		\begin{split}
				\|u\|_{B^r_{2,2}(\Omega)}
			\leq
				(
					K_{r_1}
					\| u\|_{B^{r_1}_{2,2}(\Omega)}
				)^{s_1}
				\cdot (
					K_{r_2}
					\|u\|_{B^{r_2}_{2,2}(\Omega)}
				)^{s_2}.
		\end{split}
	\end{equation}
	This finishes the proof of Lemma \ref{lem: lemma: interpolation theorem for B}.
\end{proof}
The next lemma establishes another interpolation result.
\begin{lemma}[Another interpolation result]
\label{lem: another interpolation result}
	Assume the setting in Section \ref{ssec: application}, 
	let $r_1$, $r_2$, $\tilde{r}_1$, $\tilde{r}_2$, $s_1$, $s_2 \in \R$ satisfy that
	$r_2 < s_1 \leq r_1$ and that 
	$\tilde{r}_2 = s_2 \frac{r_1-r_2}{s_1-r_2}-\tilde{r}_1 \frac{r_1-s_1}{s_1-r_2}$
	then it holds 
	for all $u \in H_{r_1}$ 
	and all $v \in H_{\tilde{r}_1\vee \tilde{r}_2}$ 
	that
	\begin{equation}
			\|u\|_{H_{s_1}} \|v\|_{H_{s_2}}
		\leq
			\tfrac{s_1-r_2}{r_1-r_2} 
				K_{r_1}\|u\|_{H_{r_1}} 
				K_{\tilde{r}_2}\|v\|_{H_{\tilde{r}_2}}
			+\tfrac{r_1-s_1}{r_1-r_2} 
				K_{r_2} \|u\|_{H_{r_2}} K_{\tilde{r}_1} \|v\|_{H_{\tilde{r}_1}}.
	\end{equation}
\end{lemma}
\begin{proof}
Note that 
it holds that
$r_1 \frac{s_1-r_2}{r_1-r_2} + r_2 \frac{r_1-s_1}{r_1-r_2} =s_1$
and thus
\eqref{eq: interpolation inequality eq}
	(with 
	$
		r \leftarrow s_1
	$,
	$
		r_1 \leftarrow r_1
	$, 
	$
		r_2 \leftarrow r_2
	$,
	$
		s_1 \leftarrow \frac{s_1-r_2}{r_1-r_2}
	$,
	and with
	$
		s_2 \leftarrow \frac{r_1-s_1}{r_1-r_2}
	$)
implies for all $u \in H_{r_1}$
that
\begin{equation}
\label{eq: u estimate}
		\|u\|_{H_{s_1}} 
	\leq 
		(K_{r_1} \|u\|_{H_{r_1}})^{\tfrac{s_1-r_2}{r_1-r_2}}
		(K_{r_2} \|u\|_{H_{r_2}})^{\tfrac{r_1-s_1}{r_1-r_2}}.
\end{equation}
Furthermore, we get from the definition of $\tilde{r}_2$ that
\begin{equation}
	\begin{split}
			&\tilde{r}_1 \frac{r_1-s_1}{r_1-r_2}+\tilde{r}_2\frac{s_1-r_2}{r_1-r_2}
		=
			\tilde{r}_1 \frac{r_1-s_1}{r_1-r_2}
			+ \left( s_2 \frac{r_1-r_2}{s_1-r_2}-\tilde{r}_1 \frac{r_1-s_1}{s_1-r_2} \right)
				\frac{s_1-r_2}{r_1-r_2} \\
		={}
			&\tilde{r}_1 \frac{r_1-s_1}{r_1-r_2}
			+s_2 - \tilde{r}_1 \frac{r_1-s_1}{r_1-r_2}
		=
			s_2.
	\end{split}
\end{equation}
Hence \eqref{eq: interpolation inequality eq}
	(with
	$
		r \leftarrow s_2
	$,
	$
		r_1 \leftarrow \tilde{r}_1
	$,
	$
		r_2 \leftarrow \tilde{r}_2
	$, 
	$
		s_1 \leftarrow \frac{r_1-s_1}{r_1-r_2}
	$,
	and with
	$
		s_2 \leftarrow \frac{s_1-r_2}{r_1-r_2} 
	$)
shows for all $v \in H_{\tilde{r}_1\vee \tilde{r}_2}$
that
\begin{equation}
\label{eq: v estimate}
		\|v\|_{H_{s_2}} 
	\leq 
		(K_{\tilde{r}_1} \|v\|_{H_{\tilde{r}_1}})^{\tfrac{r_1-s_1}{r_1-r_2}}
		(K_{\tilde{r}_2} \|v\|_{H_{\tilde{r}_2}})^{\tfrac{s_1-r_2}{r_1-r_2}}.
\end{equation}
Combining \eqref{eq: u estimate}, \eqref{eq: v estimate} together with
the fact that 
$\forall a,b \in (0,\infty), \forall t \in [0,1]\colon a^t b^{1-t} \leq ta +(1-t)b$ 
yields that
for all $u \in H_{r_1}$ and all $v \in H_{\tilde{r}_1\vee \tilde{r}_2}$ 
it holds that
\begin{equation}
	\begin{split}
			\|v\|_{H_{s_2}}\|u\|_{H_{s_1}}
		\leq{}
			&(
				K_{r_1} \|u\|_{H_{r_1}}  
				K_{\tilde{r}_2} \|v\|_{H_{\tilde{r}_2}}
			)^{\tfrac{s_1-r_2}{r_1-r_2}}
			(
				K_{r_2} \|u\|_{H_{r_2}} K_{\tilde{r}_1} \|v\|_{H_{\tilde{r}_1}}
			)^{\tfrac{r_1-s_1}{r_1-r_2}} \\
		\leq{}
			&\tfrac{s_1-r_2}{r_1-r_2} 
				K_{r_1}\|u\|_{H_{r_1}} 
				K_{\tilde{r}_2}\|v\|_{H_{\tilde{r}_2}}
			+\tfrac{r_1-s_1}{r_1-r_2} 
				K_{r_2} \|u\|_{H_{r_2}} K_{\tilde{r}_1} \|v\|_{H_{\tilde{r}_1}}
	\end{split}
\end{equation}
and this completes the proof of Lemma \ref{lem: another interpolation result}.
\end{proof}
The next lemma is a combination of the H\"older inequality
and the Sobolev inequality.
\begin{lemma}
	\label{lem: sobolev holder inequality without interpolation for Lp}
	Let $d \in \N$, 
	let $\Omega \subseteq \R^d$ be a Lipschitz set,
	let $p\in [1,\infty]$,
	let 
	$
		\eps \in 
			( 
				0, 
				\infty
			)
	$,
	let
	$r_1,r_2 \in \R$ satisfy 
	that
	$
			r_1 \wedge r_2 
		\geq 
			(0 \vee  (\tfrac {d}{2}-\tfrac {d}{p})) + \eps \1_{\{ \infty \}}(p)
	$
	and that 
	$
			r_1+r_2 
		\geq  
			d-\tfrac dp
			+\eps \1_{\{ 0 \}}((r_1\wedge r_2)-\tfrac{d}{2}+\tfrac{d}{p}).
	$
	Then there exists a real number $C \in (0, \infty)$ such that
	for all $u \in B^{r_1}_{2,2}(\Omega)$ and all $v \in B^{r_2}_{2,2}(\Omega)$ it holds that
	\begin{equation}
			\|u  v\|_{L^{p}(\Omega)}  
		\leq
			C \|u\|_{B^{r_1}_{2,2}(\Omega)} \|v\|_{B^{r_2}_{2,2}(\Omega)}.
	\end{equation}
\end{lemma}
\begin{proof}
	First note that we get from Hölder's inequality that
	for all $p_1,p_2 \in [1,\infty]$  with $\tfrac{1}{p_1}+\tfrac{1}{p_2}=\tfrac{1}{p}$,
	for all $u\in L^{p_1}(\Omega)$
	and for all $v \in L^{p_2}(\Omega)$ it holds that
	\begin{equation}
		\begin{split}
				\| u  v\|_{L^{p}(\Omega)}  
			\leq
				&\|u\|_{L^{p_1}(\Omega)} 
				\|v\|_{L^{p_2}(\Omega)}. 
		\end{split}		
	\end{equation}
	Combining this with the Sobolev inequality
	(see, e.g., the Theorem on page 31 in Runst \& Sickel \cite{RunstSickel1996}
		(applied with 
			$s \leftarrow 0$, 
			$s_0 \leftarrow \tfrac{d}{2}-\tfrac{d}{p_1}+\delta_1 \1_{\{\infty\}}(p_1)$, 
			$p \leftarrow p_1$, 
			$p_0 \leftarrow 2$, 
			$q \leftarrow 2$, 
			$u \leftarrow 2$,
			and applied with
			$s \leftarrow 0$,
			$s_0 \leftarrow \tfrac{d}{2}-\tfrac{d}{p_2}+\delta_2 \1_{\{\infty\}}(p_2)$, 
			$p \leftarrow p_2 $, 
			$p_0 \leftarrow 2$, 
			$q \leftarrow 2$,
			$u \leftarrow 2$)
		and Theorem 1 on page 32
		(applied with 
			$s \leftarrow \tfrac{d}{2}-\tfrac{d}{p_1}+\delta_1 \1_{\{\infty\}}(p_1)$, 
			$p \leftarrow 2$, 
			$q \leftarrow 2$
			and applied with 
			$s \leftarrow \tfrac{d}{2}-\tfrac{d}{p_2}+\delta_2 \1_{\{\infty\}}(p_2)$, 
			$p_0 \leftarrow 2$, $q \leftarrow 2$)) 
	shows that for all
	$\delta_1,\delta_2 \in (0,\infty)$ and all
	$p_1,p_2 \in [2,\infty]$  
	with $\tfrac{1}{p_1}+\tfrac{1}{p_2}=\tfrac{1}{p}$ 
	there exists a real number
	$C \in (0,\infty)$ such that for all
	$u\in B^{d/2-d/p_1+\delta_1 \1_{\{\infty\}}(p_1)}_{2,2}(\Omega)$
	and all $v \in B^{d/2-d/p_2+\delta_2 \1_{\{\infty\}}(p_2)}_{2,2}(\Omega)$ 
	it holds that
	\begin{equation}
	\label{eq: sobolev inequality p<p_1<infty}
		\begin{split}
				\|u v\|_{L^{p}(\Omega)}  
			\leq
				C\|u\|_{B^{d/2-d/p_1+\delta_1 \1_{\{\infty\}}(p_1)}_{2,2}(\Omega)} 
				\|v\|_{B^{d/2-d/p_2+\delta_2 \1_{\{\infty\}}(p_2)}_{2,2}(\Omega)}.
		\end{split}		
	\end{equation}
	We will now divide the proof into 4 cases.\\
	{\it 1. case: $r_1 = 0 \vee ( \tfrac{d}{2}-\tfrac dp)$.}
	Note that the assumption
	$
			r_1 
		\geq 
			(0 \vee  (\tfrac {d}{2}-\tfrac {d}{p})) + \eps \1_{\{ \infty \}}(p)
	$
	then implies that $p \neq \infty$. 
	Moreover we get from the assumption
	$
			r_1+r_2 
		\geq 
			d-\tfrac dp
			+\eps \1_{\{ 0 \}}((r_1\wedge r_2)-\tfrac{d}{2}+\tfrac{d}{p})
	$
	and from the assumption
	$
			r_1 \wedge r_2 
		\geq 
			(0 \vee  (\tfrac {d}{2}-\tfrac {d}{p})) + \eps \1_{\{ \infty \}}(p)
	$
	that
	\begin{equation}
		\begin{split}
				r_2
			\geq{} {} 
				&d-\tfrac dp
				+\eps \1_{\{ 0 \}}((r_1\wedge r_2)-\tfrac{d}{2}+\tfrac{d}{p})
				-r_1
			\geq
				d-\tfrac dp-(0 \vee ( \tfrac{d}{2}-\tfrac dp))
				+\eps \1_{\{ 0 \}}(r_1-\tfrac{d}{2}+\tfrac{d}{p})\\
			={} 
				&((d-\tfrac dp) \wedge \tfrac{d}{2})
				+\eps \1_{[2,\infty]}(p) 
			>
				((d-\tfrac dp) \wedge \tfrac{d}{2})
				- \1_{[1,2)}(p).
		\end{split}
	\end{equation}
	Therefore \eqref{eq: sobolev inequality p<p_1<infty}
	(with 
		$p_1 \leftarrow p \vee 2$,
		$\tfrac{1}{p_2} \leftarrow \tfrac{1}{p} - \tfrac {1}{p \vee 2}$,
		$\delta_1 \leftarrow 1$, 
		and with
		$\delta_2 \leftarrow r_2-((d-\tfrac dp) \wedge \tfrac{d}{2}) + \1_{[1,2)}(p)$)
	and the fact that 
	for all $s\in \R$ and all $t\in (-\infty, s]$
	there exists a real number $C \in (0,\infty)$ such that for all 
	$u \in B^{s}_{2,2}(\Omega)$ it holds that
	$\|u\|_{B^{t}_{2,2}(\Omega)} \leq C \|u\|_{B^{s}_{2,2}(\Omega)}$
	(see, e.g., the Proposition on page 29 in Runst \& Sickel \cite{RunstSickel1996} )
	shows that there exist $C_1,C_2 \in (0,\infty)$ such that for all
	$u\in B^{r_1}_{2,2}(\Omega)$
	and all $v \in B^{r_2}_{2,2}(\Omega)$ 
	it holds that
	\begin{equation}
		\begin{split}
				&\|u v\|_{L^{p}(\Omega)} \\
			\leq{} 
				&C_1\|u\|_{B^{d/2-d/(p \vee 2)+ \1_{\{ \infty \} }(p)}_{2,2}(\Omega)} 
				\|v\|_{
					B^{
						d/2-(d/p-d/(p \vee 2))
						+(
							r_2 -(d/2 \wedge (d-d/p))
							+ \1_{[1,2)}(p)
						) \1_{[2,\infty]}(p)
					}_{2,2}(\Omega)
				} \\
			={} 
				&C_1\|u\|_{B^{d/2-((d/p) \wedge (d/2))}_{2,2}(\Omega)} 
				\|v\|_{B^{
					d/2-d/p+((d/p) \wedge (d/2))
					+( r_2-((d/2) \wedge (d-d/p)) ) \1_{[2,\infty]}(p)
				}_{2,2}(\Omega)} \\
			={} 
				&C_1\|u\|_{B^{d/2-((d/p) \wedge (d/2))}_{2,2}(\Omega)} 
				\|v\|_{B^{
					((d/2) \wedge (d-d/p))
					+( r_2-((d/2) \wedge (d-d/p)) ) (1-\1_{[1,2)}(p))
				}_{2,2}(\Omega)} \\
			={} 
				&C_1\|u\|_{B^{(d/2-d/p) \vee 0}_{2,2}(\Omega)} 
				\|v\|_{B^{
					r_2-(r_2 -(d/2 \wedge (d-d/p))) \1_{[1,2)}(p)
				}_{2,2}(\Omega)}
			\leq
				C_2\|u\|_{B^{r_1}_{2,2}(\Omega)} 
				\|v\|_{B^{r_2}_{2,2}(\Omega)}.\\
		\end{split}		
	\end{equation}
	{\it 2. case: 
		$
				(0 \vee ( \tfrac{d}{2}-\tfrac dp))
			< 
				r_1
			<
				(\tfrac{d}{2} \wedge (d-\tfrac{d}{p}))
		$.}
	In this case we get that
	\begin{equation}
		\begin{split}
				\frac{d}{\tfrac{d}{2}-r_1}
			>
				\frac{d}{\tfrac{d}{2}-(0 \vee ( \tfrac{d}{2}-\tfrac dp))}
			=
				\frac{d}{\tfrac{d}{2} \wedge \tfrac dp}
			= 
				2 \vee p
		\end{split}
	\end{equation}
	and that
	\begin{equation}
		\begin{split}
				\frac 1p - \frac{\tfrac{d}{2}-r_1}{d}
			<
				\frac 1p 
				- \frac{
					\tfrac{d}{2}
					-(\tfrac{d}{2} \wedge (d-\tfrac{d}{p}))
					}
					{d}
			=
				\frac 1p 
				- \frac{
					(0 \vee (\tfrac{d}{p}-\tfrac{d}{2}))
					}
					{d}
			=
				(\tfrac 1p \wedge\tfrac{1}{2}).
		\end{split}
	\end{equation}
	Thus \eqref{eq: sobolev inequality p<p_1<infty}
	(with 
		$p_1 \leftarrow \tfrac{d}{\frac{d}{2}-r_1}$, 
		$\tfrac{1}{p_2} \leftarrow \tfrac{1}{p} - \tfrac{\frac{d}{2}-r_1}{d}$,
		$\delta_1 \leftarrow 1$, 
		and with
		$\delta_2 \leftarrow 1$),
	the fact that 
	for all $s\in \R$ and all $t\in (-\infty, s]$
	there exists a real number $C \in (0,\infty)$ such that for all 
	$u \in B^{s}_{2,2}(\Omega)$ it holds that
	$\|u\|_{B^{t}_{2,2}(\Omega)} \leq C \|u\|_{B^{s}_{2,2}(\Omega)}$
	(see, e.g., the Proposition on page 29 in Runst \& Sickel \cite{RunstSickel1996}),
	and the assumption 
	$
			r_1+r_2 
		\geq 
			d-\tfrac dp
			+\eps \1_{\{ 0 \}}((r_1\wedge r_2)-\tfrac{d}{2}+\tfrac{d}{p})
	$
	ensure that there exist $C_1,C_2 \in (0,\infty)$ such that for all
	$u\in B^{r_1}_{2,2}(\Omega)$
	and all $v \in B^{r_2}_{2,2}(\Omega)$ 
	it holds that
	\begin{equation}
		\begin{split}
				&\| u v\|_{L^{p}(\Omega)} 
			\leq
				C_1\|u\|_{B^{d/2-(d/2-r_1)}_{2,2}(\Omega)} 
				\|v\|_{
					B^{
						d/2-(d/p-(d/2-r_1))
					}_{2,2}(\Omega)
				} 
			=
				C_1\|u\|_{B^{r_1}_{2,2}(\Omega)} 
				\|v\|_{
					B^{
						d-d/p-r_1
					}_{2,2}(\Omega)
				} \\
			\leq{} 
				&C_2\|u\|_{B^{r_1}_{2,2}(\Omega)} 
				\|v\|_{B^{r_2}_{2,2}(\Omega)}.\\
		\end{split}		
	\end{equation}
	{\it 3. case: 
		$
				r_1
			=
				\tfrac{d}{2} \wedge (d-\tfrac{d}{p})
		$.} 
	Note that the assumption 
	$
			r_1 
		\geq 
			(0 \vee  (\tfrac {d}{2}-\tfrac {d}{p})) + \eps \1_{\{ \infty \}}(p)
	$
	shows that
	$p\neq \infty$
	and the assumption
	$
			r_1+r_2 
		\geq 
			d-\tfrac dp
			+\eps \1_{\{ 0 \}}((r_1\wedge r_2)-\tfrac{d}{2}+\tfrac{d}{p})
	$
	that
	\begin{equation}
		\begin{split}
				r_2 
			\geq{} 
				&d-\tfrac dp -r_1
				+\eps \1_{\{ 0 \}}((r_1\wedge r_2)-\tfrac{d}{2}+\tfrac{d}{p}) 
			=
				d-\tfrac dp 
				-(\tfrac{d}{2}\wedge (d-\tfrac{d}{p}))
				+\eps \1_{\{ 0 \}}((r_1\wedge r_2)-\tfrac{d}{2}+\tfrac{d}{p}) \\
			\geq{} 
				&d-\tfrac dp 
				-\tfrac{d}{2}
				+\eps \1_{\{ 0 \}}(r_2-\tfrac{d}{2}+\tfrac{d}{p})
			=
				\tfrac d2-\tfrac dp 
				+\eps \1_{\{ 0 \}}(r_2-\tfrac{d}{2}+\tfrac{d}{p}).
		\end{split}
	\end{equation}
	This implies that
	$
		r_2 > \tfrac d2-\tfrac dp
	$
	and hence it holds that 
	$
			\tfrac{\frac{d}{2}-r_2}{d}
		<
			\tfrac{\frac{d}{2}-(\frac d2-\frac dp)}{d}
		=
			\tfrac 1p.
	$
	Furthermore we get from
	$
			r_2 
		\geq 
			(0 \vee  (\tfrac {d}{2}-\tfrac {d}{p})) + \eps \1_{\{ \infty \}}(p)
	$
	that
	$
			\tfrac{\frac{d}{2}-r_2}{d}
		\leq
			\tfrac 12.
	$
	Therefore \eqref{eq: sobolev inequality p<p_1<infty}
	(with 
		$
				\tfrac{1}{p_1}
			\leftarrow 
				\tfrac{1}{p} - (\tfrac{\frac{d}{2}-r_2}{d} \vee \tfrac{1}{2p})
		$,
		$
				\tfrac{1}{p_2}
			\leftarrow 
				\tfrac{\frac{d}{2}-r_2}{d} \vee \tfrac{1}{2p}
		$,
		$\delta_1 \leftarrow 1$, 
		$\delta_2 \leftarrow 1$),
	the assumption that it holds that 
	$
		 r_2 
		\geq 
			(0 \vee  (\tfrac {d}{2}-\tfrac {d}{p})) + \eps \1_{\{ \infty \}}(p)
	$,
	the fact that $ p \neq \infty$,
	and the fact that 
	for all $s\in \R$ and all $t\in (-\infty, s]$
	there exists a real number $C \in (0,\infty)$ such that for all 
	$u \in B^{s}_{2,2}(\Omega)$ it holds that
	$\|u\|_{B^{t}_{2,2}(\Omega)} \leq C \|u\|_{B^{s}_{2,2}(\Omega)}$
	(see, e.g., the Proposition on page 29 in Runst \& Sickel \cite{RunstSickel1996})
	yield that there exist $C_1,C_2 \in (0,\infty)$ such that for all
	$u\in B^{r_1}_{2,2}(\Omega)$
	and all $v \in B^{r_2}_{2,2}(\Omega)$ 
	it holds that
	\begin{equation}
		\begin{split}
				&\|u v\|_{L^{p}(\Omega)}  
			\leq
				C_1\|u\|_{B^{
					d/2-(d/p-((d/2-r_2) \vee (d/(2p))))
				}_{2,2}(\Omega)} 
				\|v\|_{B^{d/2-((d/2-r_2) \vee (d/(2p)))}_{2,2}(\Omega)} \\
			\leq{} 
				&C_2\|u\|_{B^{
					d/2-d/p+((d/2-(0 \vee  (d/2-d/p))) \vee (d/(2p)))
				}_{2,2}(\Omega)} 
				\|v\|_{B^{d/2-(d/2-r_2)}_{2,2}(\Omega)} \\
			={} 
				&C_2\|u\|_{B^{
					d/2-d/p+(((d/p) \wedge (d/2)) \vee (d/(2p)))
				}_{2,2}(\Omega)} 
				\|v\|_{B^{d/2-(d/2-r_2)}_{2,2}(\Omega)} \\
			={}
				&C_2\|u\|_{B^{
					d/2-d/p+((d/p) \wedge (d/2) )
				}_{2,2}(\Omega)} 
				\|v\|_{B^{r_2}_{2,2}(\Omega)} 
			={} 
				C_2\|u\|_{B^{
					r_1
				}_{2,2}(\Omega)} 
				\|v\|_{B^{r_2}_{2,2}(\Omega)}.
		\end{split}		
	\end{equation}
	{\it 4. case: 
		$
				r_1
			>
				\tfrac{d}{2}\wedge (d-\tfrac{d}{p})
		$.}
	Then
	it follows from \eqref{eq: sobolev inequality p<p_1<infty}
	(with 
		$
				\tfrac{1}{p_1}
			\leftarrow 
				\tfrac{1}{p} - \tfrac{1}{2\vee p}
		$
		and with
		$
				p_2
			\leftarrow
				2 \vee p
		$),
	the assumption that
	$
			r_1 \wedge r_2 
		\geq 
			(0 \vee  (\tfrac {d}{2}-\tfrac {d}{p})) + \eps \1_{\{ \infty \}}(p)
	$,
	and from the fact that 
	for all $s\in \R$ and all $t\in (-\infty, s]$
	there exists a real number $C \in (0,\infty)$ such that for all 
	$u \in B^{s}_{2,2}(\Omega)$ it holds that
	$\|u\|_{B^{t}_{2,2}(\Omega)} \leq C \|u\|_{B^{s}_{2,2}(\Omega)}$
	(see, e.g., the Proposition on page 29 in Runst \& Sickel \cite{RunstSickel1996})
	that for all $\delta_1,\delta_2 \in (0,\infty)$
	satisfying that
	$\delta_1 = r_1-(\tfrac{d}{2}\wedge (d-\tfrac{d}{p}))$
	and that
	$
			\delta_2 
		= 
			r_2-	(0 \vee (\tfrac{d}{2}-\tfrac{d}{p})) +\1_{[1,\infty)}(p)
	$
	there exist $C_1,C_2 \in (0,\infty)$ such that for all
	$u\in B^{r_1}_{2,2}(\Omega)$
	and all $v \in B^{r_2}_{2,2}(\Omega)$ 
	it holds that
	\begin{equation}
		\begin{split}
				&\| u v\|_{L^{p}(\Omega)} 
			\leq
				C_1\|u\|_{B^{
					d/2-(d/p-(d/2 \wedge  d/p))
					+\delta_1 \1_{[2,\infty]}(p)
				}_{2,2}(\Omega)} 
				\|v\|_{B^{
					d/2-((d/2) \wedge  (d/p))
					+\delta_2
					\1_{\{\infty\}}(p)
				}_{2,2}(\Omega)} \\
			={} 
				&C_1\|u\|_{B^{
					((d-d/p) \wedge  (d/2))
					+\delta_1 \1_{[2,\infty]}(p)
				}_{2,2}(\Omega)} 
				\|v\|_{B^{
					(0 \vee (d/2 - d/p)) 
					+(r_2-	(0 \vee (d/2-d/p)))
					\1_{\{\infty\}}(p)
				}_{2,2}(\Omega)} \\
			={} 
				&C_1\|u\|_{B^{
					r_1
					-\delta_1 \1_{[1,2)}(p)
				}_{2,2}(\Omega)} 
				\|v\|_{B^{
					r_2 
					-(r_2-	(0 \vee (d/2-d/p)))
					\1_{[1,\infty)}(p)
				}_{2,2}(\Omega)}
			\leq
				C_2\|u\|_{B^{
					r_1
				}_{2,2}(\Omega)} 
				\|v\|_{B^{
					r_2
				}_{2,2}(\Omega)}.
		\end{split}		
	\end{equation}
	Combining the 4 cases then proves
	Lemma
	\ref{lem: sobolev holder inequality without interpolation for Lp}.
\end{proof}
The next corollary combines 
Lemma \ref{lem: another interpolation result} and 
Lemma \ref{lem: sobolev holder inequality without interpolation for Lp}.
\begin{corollary}[Sobolev-H\"older inequality for the $L^p$-norm]
	\label{cor: sobolev holder inequality with interpolation for Lp}
	Assume the setting in Section \ref{ssec: application},
	let $ p \in [1,\infty]$,
	let $\eps \in (0,\infty)$,
	let $\delta \in [0,\infty)$,
	let $\alpha, \beta \in (\N_0)^d$,
	let $r_1,r_2, \tilde{r}_1,\tilde{r}_2 \in \R$
	satisfy that
	\begin{equation}
	\label{eq: r2 r1tilde sum}
			2r_2 + 2\tilde{r}_1 
		= 
			|\alpha| +|\beta|+ d -\tfrac{d}{p} +2\delta+ 2\eps \1_{\{\infty \}}(p),
	\end{equation}
	that
	\begin{equation}
	\label{eq: r1 r2tilde sum}
				2r_1 + 2\tilde{r}_2 
			\geq 
				2r_2 + 2\tilde{r}_1 
				-\frac
					{2\delta \cdot (2r_1-2r_2)}
					{
						|\alpha| +(0 \vee  (\frac {d}{2}-\frac {d}{p}))-2r_2
					}
				+ \eps \1_{(0,\infty) \times [2,\infty)}((\delta,p)),
	\end{equation}
	that
	$
			2r_1 
		>
			|\alpha| +(0 \vee  (\tfrac {d}{2}-\tfrac {d}{p})),
	$
	and that
	$
			2r_2
		\leq
			|\alpha| +(0 \vee  (\tfrac {d}{2}-\tfrac {d}{p})) 
			- \eps \1_{[1,2) }(p).
	$
	Then there exists a real number $C \in (0, \infty)$ such that
	for all $u \in H_{r_1}$ and all $v \in H_{\tilde{r}_1 \vee \tilde{r}_2}$ 
	it holds that
	\begin{equation}
			\|(D^{\alpha}_{\R^d} u) (D^{\beta}_{\R^d} v)\|_{L^p(\Omega)} 
		\leq
			C(
				\|u\|_{H_{r_1}} \|v\|_{H_{\tilde{r}_2}}
				+ \|u\|_{H_{r_2}} \|v\|_{H_{\tilde{r}_1}}
			).
	\end{equation}
\end{corollary}
\begin{proof}
	First note that  
	for all $s_1,s_2, \tilde{\delta} \in \R$ with
	$
			r_2 +\tilde{r}_1-s_1 -s_2
		=
			\tilde{\delta}
	$
	and with
	$s_1 \neq r_2$
	it holds that
	\begin{equation}
	\label{eq: interpol help equation}
		\begin{split}
				&s_2 
					\frac
						{r_1-r_2}
						{
							s_1-r_2
						}
				-\tilde{r}_1
					\frac
						{
							r_1-s_1
						}
						{
							s_1-r_2
						}
			=
				\tilde{r}_1
				+
				(r_1-r_2)
				\Big (
					\frac
						{s_2}
						{
							s_1-r_2
						}
					-
					\frac
						{
							\tilde{r}_1
						}
						{
							s_1-r_2
						}
				\Big) \\
			={}
				&\tilde{r}_1
				+
				(r_1-r_2)
				\Big (
					\frac
						{s_1+s_2-\tilde{r}_1 -r_2}
						{
							s_1-r_2
						}
					-1
				\Big)
			={} 
				\tilde{r}_1 +r_2-r_1
				-\frac
					{\tilde{\delta} (r_1-r_2)}
						{
							s_1-r_2
						}. \\
		\end{split}
	\end{equation}
	We will now divide the proof into 3 cases. \\
	{\it 1. case: $\delta = 0$.}
	Observe that it holds that
	$
			\tfrac{d}{2}-\tfrac{d}{2p} + \eps \1_{\{\infty\}}(p)
		>
			\tfrac{d}{2}-\tfrac{d}{p}
	$, that
	$
			\tfrac{d}{2}-\tfrac{d}{2p}
		 \geq
			\tfrac{d}{2}-\tfrac{d}{p}
	$, that
	$
			\tfrac{d}{2}-\tfrac{d}{2p} 
		\geq
			0
	$,
	and that
	$
			\tfrac{d}{2}-\tfrac{d}{2p} 
		>
			-\eps \1_{[1,2)}(p)
	$.
	Hence we get that
	\begin{equation}	
		\begin{split}
				&|\alpha| +|\beta|+ d -\tfrac{d}{p}-
				\big(
					(|\alpha| +\tfrac{d}{2}-\tfrac{d}{2p} + \eps \1_{\{\infty\}}(p)) \wedge (2r_1)
				\big) 
				+2\eps \1_{\{\infty\}}(p)
				-|\beta| \\
			\geq{} 
				&|\alpha| + d -\tfrac{d}{p}-
				(|\alpha| +\tfrac{d}{2}-\tfrac{d}{2p} + \eps \1_{\{\infty\}}(p)) 
				+2\eps \1_{\{\infty\}}(p)
			=
				\tfrac{d}{2} -\tfrac{d}{2p}
				+\eps \1_{\{\infty\}}(p) \\
			\geq{} 
				&((\tfrac{d}{2}-\tfrac{d}{p}) \vee 0) + \eps \1_{\{\infty\}}(p),
		\end{split}
	\end{equation}
	that
	\begin{equation}	
		\begin{split}
				&|\alpha| +|\beta|+ d -\tfrac{d}{p}-
				\big(
					(|\alpha| +\tfrac{d}{2}-\tfrac{d}{2p} + \eps \1_{\{\infty\}}(p)) \wedge (2r_1)
				\big) 
				+2\eps \1_{\{\infty\}}(p)
				-|\beta| \\
			\geq{} 
				&|\alpha| + d -\tfrac{d}{p}-
				(|\alpha| +\tfrac{d}{2}-\tfrac{d}{2p} + \eps \1_{\{\infty\}}(p)) 
				+2\eps \1_{\{\infty\}}(p)
			=
				\tfrac{d}{2} -\tfrac{d}{2p}
				+\eps \1_{\{\infty\}}(p) \\
			>{} 
				&\tfrac{d}{2}-\tfrac{d}{p},
		\end{split}
	\end{equation}
	that
	\begin{equation}	
		\begin{split}
				&\big( 
					(|\alpha| +\tfrac{d}{2}-\tfrac{d}{2p} + \eps \1_{\{\infty\}}(p)) \wedge (2r_1) 
				\big)
				-|\alpha| 
			=
				(\tfrac{d}{2}-\tfrac{d}{2p} + \eps \1_{\{\infty\}}(p)) \wedge (2r_1-|\alpha|) \\
			\geq{} 
				&\big( ((\tfrac{d}{2}-\tfrac{d}{p}) \vee 0) + \eps \1_{\{\infty\}}(p) \big)
					\wedge (2r_1-|\alpha|) \\
			={}
				&((\tfrac{d}{2}-\tfrac{d}{p}) \vee 0)
				+\Big( 
					(\eps \1_{\{\infty\}}(p))
					\wedge \big(2r_1-|\alpha|-((\tfrac{d}{2}-\tfrac{d}{p}) \vee 0) \big)
				\Big),
		\end{split}
	\end{equation}
	and from the assumption 
	$
			2r_1 
		>
			|\alpha| +(0 \vee  (\tfrac {d}{2}-\tfrac {d}{p}))
	$ 
	that
	\begin{equation}	
		\begin{split}
				&\big( 
					(|\alpha| +\tfrac{d}{2}-\tfrac{d}{2p} + \eps \1_{\{\infty\}}(p)) \wedge (2r_1) 
				\big)
				-|\alpha| 
			=
				(\tfrac{d}{2}-\tfrac{d}{2p} + \eps \1_{\{\infty\}}(p)) \wedge (2r_1-|\alpha|) 
			>
				\tfrac{d}{2}-\tfrac{d}{p}.
		\end{split}
	\end{equation}
	Moreover, the assumption
	$
			r_2
		\leq
			|\alpha| +(0 \vee  (\tfrac {d}{2}-\tfrac {d}{p})) 
			- \eps \1_{[1,2) }(p)
	$
	implies
	that
	\begin{equation}
		\begin{split}
				r_2
			\leq
				|\alpha| +(0 \vee  (\tfrac {d}{2}-\tfrac {d}{p})) 
				- \eps \1_{[1,2) }(p)
			<
				|\alpha| +\tfrac{d}{2}-\tfrac{d}{2p} + \eps \1_{\{\infty\}}(p).
		\end{split}
	\end{equation}
	Therefore the assumption 
	$
			2r_1 
		>
			|\alpha| +(0 \vee  (\tfrac {d}{2}-\tfrac {d}{p})),
	$ 
	Lemma \ref{lem: sobolev holder inequality without interpolation for Lp}
		(with $r_1 \leftarrow 2s_1-|\alpha|$, 
		$r_2 \leftarrow 2s_2-|\beta|$, 
		and with 
		$
				\eps 
			\leftarrow 
				\eps \wedge \big(2r_1-|\alpha|-((\tfrac{d}{2}-\tfrac{d}{p}) \vee 0) \big)
		$),
	the fact that for all $r\in \R$ it holds that the
	$\| \cdot \|_{H_r}$-norm is equivalent
	to the $\| \cdot \|_{B^{2r}_{2,2}(\Omega)}$-norm on $H_r$,
	\eqref{eq: interpol help equation}
		(with $\tilde{\delta} \leftarrow 0$),
	the assumption $2r_1 > |\alpha| +(0 \vee (\tfrac d2 -\tfrac dp)) \geq 2r_2$,
	Lemma \ref{lem: another interpolation result} 
		(with $\tilde{r}_2 \leftarrow \tilde{r}_1 +r_2-r_1$),
	the fact that for all $r\in \R$ and all $t\in (-\infty, r]$
	there exists a $C \in (0,\infty)$ such that for all 
	$u \in H_{t}$ it holds that
	$\|u\|_{H_t} \leq C \|u\|_{H_r}$,
	and \eqref{eq: r1 r2tilde sum}
	show that for all $s_1,s_2 \in \R$ with 
	$
			2s_1
		=
			((|\alpha| +\tfrac{d}{2}-\tfrac{d}{2p} + \eps \1_{\{\infty\}}(p)) \wedge (2r_1)) $
	and with
	$
			2s_2
		=
			|\alpha| +|\beta|+ d -\tfrac{d}{p} -2s_1
			+2\eps \1_{\{\infty\}}(p)
	$
	there exist $C_1,C_2,C_3,C_4,C_5 \in (0,\infty)$
	such that for all
	$u \in H_{r_1}$ and all $v \in H_{\tilde{r}_1 \vee \tilde{r}_2}$ 
	it holds that
	\begin{equation}
		\begin{split}
				&\|(D^{\alpha}_{\R^d} u) (D^{\beta}_{\R^d} v)\|_{L^p(\Omega)} 
			\leq
				C_1\|D^{\alpha}_{\R^d} u\|_{B^{2s_1-|\alpha|}_{2,2}(\Omega)}
				\|D^{\beta}_{\R^d} v\|_{B^{2s_2-|\beta|}_{2,2}(\Omega)} \\
			\leq{}
				&C_2\|u\|_{B^{2s_1}_{2,2}(\Omega)} \|v\|_{B^{2s_2}_{2,2}(\Omega)} 
			\leq{} 
				C_3 \|u\|_{H_{s_1}} \|v\|_{H_{s_2}} 
			\leq
				C_4 (
					\|u\|_{H_{r_1}} \|v\|_{H_{\tilde{r}_1+r_2-r_1}}
					+ \|u\|_{H_{r_2}} \|v\|_{H_{\tilde{r}_1}}
				) \\
			\leq{}
				&C_5(
					\|u\|_{H_{r_1}} \|v\|_{H_{\tilde{r}_2}}
					+ \|u\|_{H_{r_2}} \|v\|_{H_{\tilde{r}_1}}
				).
		\end{split}
	\end{equation}
	{\it 2. case: $\delta > 0$ and
	$2r_2 <|\alpha| +(0 \vee  (\tfrac {d}{2}-\tfrac {d}{p}))$.}
	We hence get from the assumption
	$
			2r_1 
		>
			|\alpha| +(0 \vee  (\tfrac {d}{2}-\tfrac {d}{p}))
	$
	that there exists a real number 
	$
		\tilde{\eps} 
			\in \big(0, r_1-\tfrac{|\alpha|}{2}-(0 \vee  (\tfrac {d}{4}-\tfrac {d}{2p})) \big)
	$ 
	such that
	\begin{equation}
	\label{eq: def of tilde eps}
		\begin{split}
				&-\frac
					{
						(
							\delta +\eps\1_{\{\infty\}}(p) 
							-\tilde{\eps} \1_{[2,\infty]}(p)-\tilde{\eps} \1_{\{\infty\}}(p)
						)
							\cdot (2r_1-2r_2)
					}
					{
						|\alpha| +(0 \vee  (\frac {d}{2}-\frac {d}{p}))
						+2\tilde{\eps} \1_{\{\infty\}}(p)-2r_2
					}\\
			\leq{} 
				&-\frac
					{\delta \cdot (2r_1-2r_2)}
					{
						|\alpha| +(0 \vee  (\frac {d}{2}-\frac {d}{p}))-2r_2
					}
				+ \eps \1_{(0,\infty) \times [2,\infty)}((\delta,p)).
		\end{split}
	\end{equation}
	Thus
	the fact that it holds that
	$ 
		(0 \vee  (\tfrac {d}{2}-\tfrac {d}{p})) \leq (\tfrac {d}{2} \wedge (d-\tfrac {d}{p}))
	$,
	that
	$ 
			(0 \vee  (\tfrac {d}{2}-\tfrac {d}{p})) + (\tfrac {d}{2} \wedge (d-\tfrac {d}{p}))
		=
			d-\tfrac dp
	$,
	and that
	$
		\big( (0 \vee  (\tfrac {d}{2}-\tfrac {d}{p})) = \tfrac {d}{2}-\tfrac {d}{p} \big)
		\Leftrightarrow (p \in [2,\infty])
	$,
	Lemma \ref{lem: sobolev holder inequality without interpolation for Lp}
		(with $r_1 \leftarrow 2s_1-|\alpha|$, 
		$r_2 \leftarrow 2s_2-|\beta|$,
		$\eps \leftarrow \tilde{\eps}$),
	the fact that for all $r\in \R$ it holds that the
	$\| \cdot \|_{H_r}$-norm is equivalent
	to the $\| \cdot \|_{B^{2r}_{2,2}(\Omega)}$-norm on $H_r$,
	\eqref{eq: interpol help equation}
		(with 
		$
				\tilde{\delta} 
			\leftarrow 
				\delta +\eps\1_{\{\infty\}}(p) 
				-\tilde{\eps} \1_{[2,\infty]}(p)-\tilde{\eps} \1_{\{\infty\}}(p)
		$),
	the fact that
	$
		2\tilde{\eps} < 2r_1-|\alpha|-(0 \vee  (\tfrac {d}{2}-\tfrac {d}{p}))
	$,
	Lemma \ref{lem: another interpolation result} 
		(with 
		$
				\tilde{r}_2 
			\leftarrow
				r_2 + \tilde{r}_1 -r_1 
				-\frac
					{
						(
							\delta +\eps\1_{\{\infty\}}(p) 
							-\tilde{\eps} \1_{[2,\infty]}(p)-\tilde{\eps} \1_{\{\infty\}}(p)
						)
							\cdot (2r_1-2r_2)
					}
					{
						|\alpha| +(0 \vee  (\frac {d}{2}-\frac {d}{p}))
						+2\tilde{\eps} \1_{\{\infty\}}(p)-2r_2 
					}
		$),
	the fact that for all $r\in \R$ and all $t\in (-\infty, r]$
	there exists a $C \in (0,\infty)$ such that for all 
	$u \in H_{t}$ it holds that
	$\|u\|_{H_t} \leq C \|u\|_{H_r}$,
	\eqref{eq: def of tilde eps},
	and \eqref{eq: r1 r2tilde sum}
	ensure that for all $s_1,s_2 \in \R$ with 
	$
			2s_1
		=
			|\alpha| 
			+(0 \vee  (\tfrac {d}{2}-\tfrac {d}{p}))
			+2\tilde{\eps} \1_{\{\infty\}}(p) 
	$
	and with
	$	
			2s_2
		=
			|\beta| +(\tfrac {d}{2} \wedge (d-\tfrac {d}{p})) 
			+ 2\tilde{\eps} \1_{[2,\infty]}(p)
	$
	there exist $C_1,C_2,C_3,C_4, C_5 \in (0,\infty)$
	such that for all
	$u \in H_{r_1}$ and all $v \in H_{\tilde{r}_1 \vee \tilde{r}_2}$ 
	it holds that
	\begin{align}
				&\|(D^{\alpha}_{\R^d} u) (D^{\beta}_{\R^d} v)\|_{L^p(\Omega)} 
			\leq
				C_1 \| D^{\alpha}_{\R^d} u \|_{B^{2s_1-|\alpha|}_{2,2}(\Omega)} 
					\|D^{\beta}_{\R^d} v\|_{B^{2s_2-|\beta|}_{2,2}(\Omega)} \\
			\leq{} \nonumber
				&C_2\|u\|_{B^{2s_1}_{2,2}(\Omega)} \|v\|_{B^{2s_2}_{2,2}(\Omega)} 
			\leq
				C_3 \|u\|_{H_{s_1}} \|v\|_{H_{s_2}} \\
			\leq{} \nonumber
				&C_4 (
					\|u\|_{H_{r_1}} 
					\|v\|_{H_{
						r_2 + \tilde{r}_1 -r_1 
						-(
							\delta +\eps\1_{\{\infty\}}(p) 
							-\tilde{\eps} \1_{[2,\infty]}(p)-\tilde{\eps} \1_{\{\infty\}}(p)
						)
						\cdot (2r_1-2r_2)
						/(
							|\alpha| +(0 \vee  (d/2-d/p))
							+2\tilde{\eps} \1_{\{\infty\}}(p)-2r_2
						)	
					}} \\ \nonumber
					&+ \|u\|_{H_{r_2}} \|v\|_{H_{\tilde{r}_1}}
				) \\ \nonumber
			\leq{} 
				&C_5(
					\|u\|_{H_{r_1}} \|v\|_{H_{\tilde{r}_2}}
					+ \|u\|_{H_{r_2}} \|v\|_{H_{\tilde{r}_1}}
				).
	\end{align}
	{\it 3. case: $\delta > 0$ and
	$2r_2 =|\alpha_1| +(0 \vee  (\tfrac {d}{2}-\tfrac {d}{p}))$.}
	Note that it follows
	from the assumption
	$
			2r_1 
		>
			|\alpha| +(0 \vee  (\tfrac {d}{2}-\tfrac {d}{p}))
	$
	that there exists a real number 
	$
		\hat{\eps} 
			\in 
				\big(
					0, r_1-\tfrac{|\alpha|}{2}-(0 \vee  (\tfrac {d}{4}-\tfrac {d}{2p})) 
				\big)
	$
	such that
	\begin{equation}
	\label{eq: def of hat eps}
				r_2+\tilde{r}_1-r_1
				-\frac
					{
						(
							\delta +\eps\1_{\{\infty\}}(p) 
							-2\hat{\eps}
						)
							\cdot (r_1-r_2)
					}
					{
						\hat{\eps}
					}
			\leq
				\tilde{r}_2.
	\end{equation}
	Moreover,
	the fact that it holds that
	$ 
		(0 \vee  (\tfrac {d}{2}-\tfrac {d}{p})) \leq (\tfrac {d}{2} \wedge (d-\tfrac {d}{p}))
	$
	and that
	$ 
			(0 \vee  (\tfrac {d}{2}-\tfrac {d}{p})) + (\tfrac {d}{2} \wedge (d-\tfrac {d}{p}))
		=
			d-\tfrac dp
	$,
	Lemma \ref{lem: sobolev holder inequality without interpolation for Lp}
		(with $r_1 \leftarrow 2s_1-|\alpha|$, 
		$r_2 \leftarrow 2s_2-|\beta|$,
		$\eps \leftarrow \hat{\eps}$),
	the fact that for all $r\in \R$ it holds that the
	$\| \cdot \|_{H_r}$-norm is equivalent
	to the $\| \cdot \|_{B^{2r}_{2,2}(\Omega)}$-norm on $H_r$,
	\eqref{eq: interpol help equation}
		(with 
		$
				\tilde{\delta} 
			\leftarrow 
				\delta +\eps\1_{\{\infty\}}(p) 
				-2\hat{\eps}
		$),
	the fact that
	$
		2\hat{\eps} < 2r_1-|\alpha|-(0 \vee  (\tfrac {d}{2}-\tfrac {d}{p}))
	$,
	Lemma \ref{lem: another interpolation result} 
		(with 
		$
				\tilde{r}_2 
			\leftarrow 
				r_2 + \tilde{r}_1 -r_1 
				-\frac
					{
						(
							\delta +\eps\1_{\{\infty\}}(p) 
							-2\hat{\eps}
						)
							\cdot (r_1-r_2)
					}
					{
						\hat{\eps}
					}
		$),
	the fact that for all $r\in \R$ and all $t\in (-\infty, r]$
	there exists a $C \in (0,\infty)$ such that for all 
	$u \in H_{t}$ it holds that
	$\|u\|_{H_t} \leq C \|u\|_{H_r}$,
	and
	\eqref{eq: def of hat eps}
	show that for all $s_1,s_2 \in \R$ with 
	$
			2s_1
		=
			|\alpha| 
			+(0 \vee  (\tfrac {d}{2}-\tfrac {d}{p}))
			+2\hat{\eps} 
	$
	and with
	$	
			2s_2
		=
			|\beta| +(\tfrac {d}{2} \wedge (d-\tfrac {d}{p})) 
			+ 2\hat{\eps}
	$
	there exist $C_1,C_2,C_3,C_4, C_5 \in (0,\infty)$
	such that for all
	$u \in H_{r_1}$ and all $v \in H_{\tilde{r}_1 \vee \tilde{r}_2}$ 
	it holds that
	\begin{equation}
		\begin{split}
				&\|(D^{\alpha}_{\R^d} u) (D^{\beta}_{\R^d} v)\|_{L^p(\Omega)}
			\leq
				C_1\|D^{\alpha}_{\R^d} u\|_{B^{2s_1-|\alpha|}_{2,2}(\Omega)}
				\|D^{\beta}_{\R^d} v\|_{B^{2s_2-|\beta|}_{2,2}(\Omega)}
			\leq
				C_2\|u\|_{B^{2s_1}_{2,2}(\Omega)} \|v\|_{B^{2s_2}_{2,2}(\Omega)} \\
			\leq{} 
				&C_3 \|u\|_{H_{s_1}} \|v\|_{H_{s_2}}
			\leq
				C_4 (
					\|u\|_{H_{r_1}} 
					\|v\|_{H_{
						r_2 + \tilde{r}_1 -r_1 
						-(
							\delta +\eps\1_{\{\infty\}}(p) 
							-2\hat{\eps}
						)
						\cdot (r_1-r_2)
						/\hat{\eps} 
					}}
					+ \|u\|_{H_{r_2}} \|v\|_{H_{\tilde{r}_1}}
				) \\
			\leq{} 
				&C_5(
					\|u\|_{H_{r_1}} \|v\|_{H_{\tilde{r}_2}}
					+ \|u\|_{H_{r_2}} \|v\|_{H_{\tilde{r}_1}}
				).
		\end{split}
	\end{equation}
	Combining the 3 cases then proves 
	Corollary \ref{cor: sobolev holder inequality with interpolation for Lp}.
\end{proof}
The next lemma establishes a similar result as in 
Corollary \ref{cor: sobolev holder inequality with interpolation for Lp}
for positive Besov spaces. It
generalizes Lemma 1 on page 345 in Runst \& Sickel
 \cite{RunstSickel1996}.
\begin{lemma}[Sobolev-H\"older inequality for positive Besov spaces]
	\label{lem: sobolev holder inequality for positive s}
	Assume the setting in Section \ref{ssec: application},
	let $ p\in (1,\infty)$,
	let $s \in (0,\infty)$,
	let $\eps \in (0,\infty)$,
	let $\delta \in [0,\infty)$,
	let $\alpha, \beta \in (\N_0)^d$,
	let $r_1,r_2, \tilde{r}_1,\tilde{r}_2 \in \R$
	satisfy that
	\begin{equation}
	\label{eq: r2 r1tilde sum s>0}
			2r_2 + 2\tilde{r}_1 
		= 
			2s+|\alpha| +|\beta|+ d -\tfrac{d}{p} +2\delta,
	\end{equation}
	that
	\begin{equation}
	\label{eq: r1 r2tilde sum lem s>0}
			2r_1 + 2\tilde{r}_2 
		\geq 
			2s+|\alpha| +|\beta|+d -\tfrac{d}{p}
			-\frac
				{
					2\delta \cdot (2r_1-2s-|\alpha| -(0 \vee  (\frac {d}{2}-\frac {d}{p})))
				}
				{
					2s+|\alpha| +(0 \vee  (\frac {d}{2}-\frac {d}{p}))-2r_2
				} 
			+ \eps \1_{(0,\infty) \times [2,\infty)}((\delta,p)),
	\end{equation}
	that
	$
			2r_1 
		\geq
			2s+|\alpha| +(0 \vee  (\tfrac {d}{2}-\tfrac {d}{p})),
	$
	that
	$
			2\tilde{r}_1 
		\geq
			2s+|\beta| +(0 \vee  (\tfrac {d}{2}-\tfrac {d}{p})),
	$
	that
	$r_1 \geq r_2$,
	that either
	$
			2r_2
		\leq
			2s+|\alpha| +(0 \vee  (\tfrac {d}{2}-\tfrac {d}{p}))
	$
	or
	$\delta =0$,
	that either
	$r_1 > r_2$
	or $\delta =0$
	and that 
	$
			\delta 
			\vee (
				(2r_1-2s-|\alpha| -\tfrac {d}{2}+\tfrac {d}{p}) 
				\cdot (2\tilde{r}_1-2s-|\beta| -\tfrac {d}{2}+\tfrac {d}{p})
			)
		>0.
	$
	Then there exists a $C \in (0, \infty)$ such that
	for all $u \in H_{r_1}$ and all $v \in H_{\tilde{r}_1 \vee \tilde{r}_2}$ 
	it holds that
	\begin{align}
			\|(D^{\alpha}_{\R^d} u) (D^{\beta}_{\R^d} v)\|_{B^{2s}_{p,2}(\Omega)} 
		\leq
			C(
				\|u\|_{H_{r_1}} \|v\|_{H_{\tilde{r}_2}}
				+ \|u\|_{H_{r_2}} \|v\|_{H_{\tilde{r}_1}}
			).
	\end{align}
\end{lemma}
\begin{proof}
	First observe that  
	for all $s_1,s_2, \tilde{\delta} \in \R$ with
	$
			r_2 +\tilde{r}_1-s_1 -s_2
		=
			\tilde{\delta}
	$
	and with
	$s_1 \neq r_2$
	it holds that
	\begin{equation}
	\label{eq: interpol help equation lem s>0}
		\begin{split}
				&s_2 
					\frac
						{r_1-r_2}
						{
							s_1-r_2
						}
				-\tilde{r}_1
					\frac
						{
							r_1-s_1
						}
						{
							s_1-r_2
						}
			=
				\tilde{r}_1
				+
				(r_1-r_2)
				\Big (
					\frac
						{s_2}
						{
							s_1-r_2
						}
					-
					\frac
						{
							\tilde{r}_1
						}
						{
							s_1-r_2
						}
				\Big) \\
			={}
				&\tilde{r}_1
				+
				(r_1-r_2)
				\Big (
					\frac
						{s_1+s_2-\tilde{r}_1 -r_2}
						{
							s_1-r_2
						}
					-1
				\Big) 
			=
				\tilde{r}_1 +r_2-r_1
				-\frac
					{\tilde{\delta} (r_1-r_2)}
						{
							s_1-r_2
						}.
		\end{split}
	\end{equation}
	In addition it follows from Theorem 4.1 in Rychkov \cite{Rychkov1998} 
	that there exists a 
	\begin{equation}
		\mathscr{E} \colon 
			\cup_{\hat{s} \in \R} \cup_{\hat{p} \in [1,\infty]} 
				B^{\hat{s}}_{\hat{p},2}(\Omega) 
				\to \cup_{\hat{s} \in \R} \cup_{\hat{p} \in [1,\infty]} 
					B^{\hat{s}}_{\hat{p},2}(\R^d) 
	\end{equation}
	such that for all $\hat{s} \in \R$ and all $\hat{p} \in [1,\infty]$
	it holds that $\mathscr{E}|_{B^{\hat{s}}_{\hat{p},2}(\Omega)}$ 
	is a linear and bounded
	extension operator from $B^{\hat{s}}_{\hat{p},2}(\Omega)$ to 
	$B^{\hat{s}}_{\hat{p},2}(\R^d)$.
	Moreover, we get from Lemma 1 on page 345 in Runst \& Sickel \cite{RunstSickel1996} 
		(with $r_1 \leftarrow q_2$, $r_2 \leftarrow q_1$,
		$p_1 \leftarrow p_1$, $p_2 \leftarrow p_2$,
		and $q \leftarrow 2$)
	and from 
	the Sobolev inequality 
	(see, e.g., the Theorem on page 31 in Runst \& Sickel \cite{RunstSickel1996} 
	and Theorem 1 on page 32)
	that for all $q_i, p_i \in [2,\infty]$ and all $\delta_i \in (0,\infty)$,
	$i \in \{ 1, 2\}$ with 
	$
			\tfrac 1p
		=
			\tfrac {1}{p_1} + \tfrac {1}{q_2}
		=
			\tfrac {1}{p_2} + \tfrac {1}{q_1}
	$
	there exist $C_1,C_2,C_3 \in (0,\infty)$
	such that for all 
	$u \in B^{d/2+((2s-d/p_1) \vee (-d/q_1+\delta_1 \1_{\{\infty\}}(q_1)))}_{2,2}(\Omega)$
	and all $v \in B^{d/2+((2s-d/p_2) \vee (-d/q_2+\delta_2 \1_{\infty}(q_2)))}_{2,2}(\Omega)$
	it holds that
	\begin{align}
				&\|u v\|_{B^{2s}_{p,2}(\Omega)}  
			=
				\inf_{
					z \in B^{2s}_{p,2}(\R^d), 
					z|_{\Omega} = u v
				}
					\|z\|_{B^{2s}_{p,2}(\R^d)} \\
			\leq{}  \nonumber
				&\inf_{
					\tilde{u} \in 
						B^{d/2+((2s-d/p_1) \vee (-d/q_1+\delta_1 \1_{\{\infty\}}(q_1)))}_{2,2}(\R^d), 
					\tilde{u}|_\Omega = u}
				\inf_{
					\tilde{v} \in 
						B^{d/2+((2s-d/p_2) \vee (-d/q_2+\delta_2 \1_{\{\infty\}}(q_2)))}_{2,2}(\R^d), 
					\tilde{v}|_\Omega = u} \\
			& \nonumber \qquad
					\|\tilde{u} \cdot \tilde{v}\|_{B^{2s}_{p,2}(\R^d)} \\
			\leq{}  \nonumber
				&\inf_{
					\tilde{u} \in 
						B^{d/2+((2s-d/p_1) \vee (-d/q_1+\delta_1 \1_{\{\infty\}}(q_1)))}_{2,2}(\R^d), 
					\tilde{u}|_\Omega = u}
				\inf_{
					\tilde{v} \in 
						B^{d/2+((2s-d/p_2) \vee (-d/q_2+\delta_2 \1_{\{\infty\}}(q_2)))}_{2,2}(\R^d), 
					\tilde{v}|_\Omega = u} \\
			& \nonumber \qquad
				C_1
				\big(
					\|
						\tilde{u}
					\|_{B^{2s}_{p_1,2}(\R^d)}
					\|
						\tilde{v}
					\|_{L^{q_2}(\R^d)}
					+
					\|
						\tilde{u}
					\|_{L^{q_1}(\R^d)}
					\|
						\tilde{v}
					\|_{B^{2s}_{p_2,2}(\R^d)} 
				\big)\\
			\leq{} \nonumber
				&\inf_{
					\tilde{u} \in 
						B^{d/2+((2s-d/p_1) \vee (-d/q_1+\delta_1 \1_{\{\infty\}}(q_1)))}_{2,2}(\Omega), 
					\tilde{u}|_\Omega = u}
				\inf_{
					\tilde{v} \in 
						B^{d/2+((2s-d/p_2) \vee (-d/q_2+\delta_2 \1_{\{\infty\}}(q_2)))}_{2,2}(\Omega), 
					\tilde{v}|_\Omega = u} \\
			& \nonumber \qquad
				C_2
				\big(
					\|
						\tilde{u} 
					\|_{B^{2s+d/2-d/p_1}_{2,2}(\R^d)}
					\|
						\tilde{v}
					\|_{B^{d/2-d/q_2+\delta_2 \1_{\{\infty\}}(q_2)}_{2,2}(\R^d)} 
				\\& \nonumber \qquad\qquad
					+\|
						\tilde{u} 
					\|_{B^{d/2-d/q_1+\delta_1 \1_{\{\infty\}}(q_1)}_{2,2}(\R^d)}
					\|
						\tilde{v}
					\|_{B^{2s+d/2-d/p_2}_{2,2}(\R^d)} 
				\big) \\
			\leq{} \nonumber
				&C_2 \big(
					\|
						\mathscr{E} u 
					\|_{B^{2s+d/2-d/p_1}_{2,2}(\R^d)}
					\|
						\mathscr{E} v
					\|_{B^{d/2-d/q_2+\delta_2 \1_{\{\infty\}}(q_2)}_{2,2}(\R^d)} \\
					& \nonumber \qquad
					+\|
						\mathscr{E} u
					\|_{B^{d/2-d/q_1+\delta_1 \1_{\{\infty\}}(q_1)}_{2,2}(\R^d)}
					\|
						\mathscr{E} v
					\|_{B^{2s+d/2-d/p_2}_{2,2}(\R^d)}
				\big)\\
			\leq{} \nonumber
				&C_3\big(
					\|
						u 
					\|_{B^{2s+d/2-d/p_1}_{2,2}(\Omega)}
					\|
						v
					\|_{B^{d/2-d/q_2+\delta_2 \1_{\{\infty\}}(q_2)}_{2,2}(\Omega)}
				\\ \nonumber & \qquad
					+\|
						u
					\|_{B^{d/2-d/q_1+\delta_1 \1_{\{\infty\}}(q_1)}_{2,2}(\Omega)}
					\|
						v
					\|_{B^{2s+d/2-d/p_2}_{2,2}(\Omega)} 
				\big)
	\end{align}
	and this shows that
	for all $q_i, p_i \in [2,\infty]$ and all $\delta_i \in (0,\infty)$,
	$i \in \{ 1, 2\}$, with 
	$
			\tfrac 1p
		=
			\tfrac {1}{p_1} + \tfrac {1}{q_2}
		=
			\tfrac {1}{p_2} + \tfrac {1}{q_1}
	$
	there exist 
	$C_1, C_2, C_3 \in (0,\infty)$
	such that for all
	$
		u \in 
			H_{(s-d/(2p_1)) \vee (-d/(2q_1)+\delta_1 \1_{\{\infty\}}(q_1))}
	$
	and all 
	$
		v \in 
			H_{(s-d/(2p_2)) \vee (-d/(2q_2)+\delta_2 \1_{\{\infty\}}(q_2))}
	$
	it holds that
	\begin{align}
	\nonumber
				&\|(D^{\alpha}_{\R^d} u) (D^{\beta}_{\R^d} v)\|_{B^{2s}_{p,2}(\Omega)}  \\ \nonumber
			\leq{} 
				&C_1(
					\|
						D^{\alpha}_{\R^d} u 
					\|_{B^{2s+d/2-d/p_1}_{2,2}(\Omega)}
					\|
						D^{\beta}_{\R^d} v
					\|_{B^{d/2-d/q_2+2\delta_2 \1_{\{\infty\}}(q_2)}_{2,2}(\Omega)}
				\\ \nonumber & \quad
					+\|
						D^{\alpha}_{\R^d} u
					\|_{B^{d/2-d/q_1+2\delta_1 \1_{\{\infty\}}(q_1)}_{2,2}(\Omega)}
					\|
						D^{\beta}_{\R^d} v
					\|_{B^{2s+d/2-d/p_2}_{2,2}(\Omega)}
				)\\
	\label{eq: basic estimate s>0}
			\leq{} 
				&C_2(
					\|
						u 
					\|_{B^{2s+|\alpha|+d/2-d/p_1}_{2,2}(\Omega)}
					\|
						v
					\|_{B^{|\beta|+d/2-d/q_2+2\delta_2 \1_{\{\infty\}}(q_2)}_{2,2}(\Omega)}
			\\ \nonumber & \quad
					+\|
						u
					\|_{B^{|\alpha|+d/2-d/q_1+2\delta_1 \1_{\{\infty\}}(q_1)}_{2,2}(\Omega)}
					\|
						v
					\|_{B^{2s+|\beta|+d/2-d/p_2}_{2,2}(\Omega)} 
				)\\ \nonumber
			\leq{} 
				&C_3(
					\|
						u 
					\|_{H_{s+|\alpha|/2+d/4-d/(2p_1)}}
					\|
						v
					\|_{H_{|\beta|/2+d/4-d/(2q_2)+\delta_2 \1_{\{\infty\}}(q_2)}}
				\\ \nonumber & \quad
					+\|
						u
					\|_{H_{|\alpha|/2+d/4-d/(2q_1)+\delta_1 \1_{\{\infty\}}(q_1)}}
					\|
						v
					\|_{H_{s+|\beta|/2+d/4-d/(2p_2)}}
				).
	\end{align}
	We will now divide the proof into 5 cases. \\
	{\it 1. case: 
		$
				\delta=0, r_1=r_2.
		$} 
	First observe that the assumption
	$
			\delta 
			\vee (
				(2r_1-2s-|\alpha| -\tfrac {d}{2}+\tfrac {d}{p}) 
				\cdot (2\tilde{r}_1-2s-|\beta| -\tfrac {d}{2}+\tfrac {d}{p})
			)
		>0
	$
	ensures that
	\begin{equation}
			2r_2
		=
			2r_1 
		>
			2s+|\alpha| + \tfrac {d}{2}-\tfrac {d}{p}
	\end{equation}
	and that
	\begin{equation}
			2\tilde{r}_1 
		>
			2s+|\beta| + \tfrac {d}{2}-\tfrac {d}{p}.
	\end{equation}
	Moreover,  \eqref{eq: r2 r1tilde sum s>0} 
	the assumptions	
	$
			2\tilde{r}_1 
	 \geq
			2s+|\beta| +(0 \vee  (\tfrac {d}{2}-\tfrac {d}{p})),
	$
	and the assumptions
	$
			2r_2=2r_1 
	 \geq
			2s+|\alpha| +(0 \vee  (\tfrac {d}{2}-\tfrac {d}{p}))
	$
	show that
	$
			2r_2
		= 
			2s+|\alpha| +|\beta|+ d -\tfrac{d}{p} +2\delta -\tilde{r}_1
		\leq
			|\alpha|  + ((d -\tfrac{d}{p}) \wedge  \tfrac {d}{2})
	$
	and that
	$
			2\tilde{r}_1
		= 
			2s+|\alpha| +|\beta|+ d -\tfrac{d}{p} +2\delta -r_2
		\leq
			|\beta|  + ((d -\tfrac{d}{p}) \wedge  \tfrac {d}{2}).
	$
	Hence there exists 
	$\hat{p}_1,\hat{p}_2 \in (p ,\infty) \cap [2,\tfrac{2p}{(2-p) \vee 0}]$
	such that
	\begin{equation}
		s+\tfrac{|\alpha|}{2}+\tfrac{d}{4}-\tfrac{d}{2\hat{p}_1} = r_2,
	\end{equation}	
	and that
	\begin{equation}
		s+\tfrac{|\beta|}{2}+\tfrac{d}{4}-\tfrac{d}{2\hat{p}_2} = \tilde{r}_1.
	\end{equation}	
	Thus \eqref{eq: basic estimate s>0}
		(with 
		$p_1 \leftarrow \hat{p}_1$, $p_2 \leftarrow \hat{p}_2$,
		$\tfrac{1}{q_1} \leftarrow \tfrac{1}{p}-\tfrac{1}{\hat{p}_2}$,
		and with
		$\tfrac{1}{q_2} \leftarrow \tfrac{1}{p}-\tfrac{1}{\hat{p}_1}$)
	and \eqref{eq: r2 r1tilde sum s>0}
	yield that there exists a $C \in (0,\infty)$
	such that for all
	$
		u \in 
			H_{r_2}
	$
	and all 
	$
		v \in 
			H_{\tilde{r}_1}
	$
	it holds that
	\begin{equation}
		\begin{split}
				&\|(D^{\alpha}_{\R^d} u) (D^{\beta}_{\R^d} v)\|_{B^{2s}_{p,2}(\Omega)}  \\
			\leq{} 
				&C(
					\|
						u 
					\|_{H_{s+|\alpha|/2+d/4-d/(2\hat{p}_1)}}
					\|
						v
					\|_{H_{|\beta|/2+d/4-d/(2p)+d/(2\hat{p}_1)}}
			\\ & \quad
					+\|
						u
					\|_{H_{|\alpha|/2+d/4-d/(2p)+d/(2\hat{p}_2)}}
					\|
						v
					\|_{H_{s+|\beta|/2+d/4-d/(2\hat{p}_2)}}
				) \\
			={} 
				&C(
					\|
						u 
					\|_{H_{r_2}}
					\|
						v
					\|_{H_{s+|\alpha|/2+|\beta|/2+d/2-d/(2p)-r_2}}
					+\|
						u
					\|_{H_{s+|\alpha|/2+|\beta|/2+d/2-d/(2p)-\tilde{r}_1}}
					\|
						v
					\|_{H_{\tilde{r}_1}}
				) \\
			={}
				&2C \cdot
					\|
						u 
					\|_{H_{r_2}}
					\|
						v
					\|_{H_{\tilde{r}_1}}.
		\end{split}
	\end{equation}
	\\
	{\it 2. case: 
		$
				\delta=0, r_1>r_2.
		$} 
	Again the assumption
	$
			\delta 
			\vee (
				(2r_1-2s-|\alpha| -\tfrac {d}{2}+\tfrac {d}{p}) 
				\cdot (2\tilde{r}_1-2s-|\beta| -\tfrac {d}{2}+\tfrac {d}{p})
			)
		>0
	$
	verifies that
	\begin{equation}
	\label{eq: r1 big}
			2r_1 
		>
			2s+|\alpha| + \tfrac {d}{2}-\tfrac {d}{p}
	\end{equation}
	and that
	\begin{equation}
	\label{eq: r1 tilde big}
			2\tilde{r}_1 
		>
			2s+|\beta| + \tfrac {d}{2}-\tfrac {d}{p}.
	\end{equation}
	In addition, \eqref{eq: r2 r1tilde sum s>0} 
	and the assumption 
	$
			2\tilde{r}_1 
		\geq
			2s+|\beta| +(0 \vee  (\tfrac {d}{2}-\tfrac {d}{p}))
	$
	show that
	\begin{equation}
			2r_2  
		= 
			2s+|\alpha| +|\beta|+ d -\tfrac{d}{p} +2\delta -\tilde{r}_1
		\leq
			|\alpha|  + ((d -\tfrac{d}{p}) \wedge  \tfrac {d}{2}).
	\end{equation}
	Therefore
	the assumption $r_1>r_2$
	and the assumption
	$
			2r_1 
		\geq
			2s+|\alpha| +(0 \vee  (\tfrac {d}{2}-\tfrac {d}{p}))
	$
	imply that
	there exists a $\tilde{p}_1 \in (p ,\infty) \cap [2,\tfrac{2p}{(2-p) \vee 0}]$ 
	such that
	\begin{equation}
	\label{eq: def of hat p1}
			r_2
		< 
			s+\tfrac{|\alpha|}{2} +\tfrac d4 - \tfrac{d}{2\tilde{p}_1}
		\leq
			r_1.
	\end{equation}
	Moreover, we get from the assumption $r_1 > r_2$, 
	\eqref{eq: r2 r1tilde sum s>0},
	the assumption
	$
			2r_1
		\geq
			2s+|\alpha| +(0 \vee  (\tfrac {d}{2}-\tfrac {d}{p})),
	$
	and from 
	\eqref{eq: r1 big}
	that it holds that
	\begin{equation}
		\tilde{r}_1+r_2-r_1 < \tilde{r}_1,
	\end{equation}
	that
	\begin{align}
			2\tilde{r}_1+2r_2-2r_1 
		=
			2s+|\alpha|+|\beta| +d - \tfrac dp -2r_1
		\leq
			|\beta| + ((d - \tfrac dp) \wedge \tfrac d2)
	\end{align}
	and that
	\begin{align}
			2\tilde{r}_1+2r_2-2r_1 
		=
			2s+|\alpha|+|\beta| +d - \tfrac dp -2r_1
		<
			|\beta| +  \tfrac d2.
	\end{align}
	Combining this with the assumption that
	$
			2\tilde{r}_1
		\geq
			2s+|\beta| +(0 \vee  (\tfrac {d}{2}-\tfrac {d}{p}))
	$
	and 
	\eqref{eq: r1 tilde big}
	yields
	that
	there exists a $\tilde{p}_2 \in [2,\infty)$ such that
	$
		\tilde{p}_2 > p,
	$
	that
	$
		\tfrac{1}{\tilde{p}_2} \geq \tfrac{1}{p} -\tfrac{1}{2},
	$
	and that
	\begin{equation}
	\label{eq: def of hat p2}
			\tilde{r}_1+r_2-r_1
		< 
			s+\tfrac{|\beta|}{2} +\tfrac d4 - \tfrac{d}{2\tilde{p}_2}
		\leq
			\tilde{r}_1.
	\end{equation}
	Thus \eqref{eq: basic estimate s>0}
	(with 
		$p_1 \leftarrow \tilde{p}_1$, 
		$\tfrac{1}{q_2} \leftarrow \tfrac{1}{p}-\tfrac{1}{\tilde{p}_1}$, 
		$p_2 \leftarrow \tilde{p}_2$, 
		$\tfrac{1}{q_1} \leftarrow \tfrac{1}{p}-\tfrac{1}{\tilde{p}_2}$),
	Lemma \eqref{lem: another interpolation result}
		(applied with 
		$s_1 \leftarrow s+\tfrac{|\alpha|}{2}+\tfrac{d}{4}-\tfrac{d}{2\tilde{p}_1}$,
		$
				s_2 
			\leftarrow 
				\tfrac{|\beta|}{2}+\tfrac{d}{4}-\tfrac{d}{2p}+\tfrac{d}{2\tilde{p}_1}
		$,
		$\tilde{r}_2 \leftarrow \tilde{r}_1+r_2-r_1$
		and applied
		with $s_1 \leftarrow s+\tfrac{|\beta|}{2}+\tfrac{d}{4}-\tfrac{d}{2\tilde{p}_2}$,
		$
				s_2 
			\leftarrow
				\tfrac{|\alpha|}{2}+\tfrac{d}{4}-\tfrac{d}{2p}+\tfrac{d}{2\tilde{p}_2}
		$,
		$r_1 \leftarrow \tilde{r}_1$,
		$\tilde{r}_1 \leftarrow r_1$,
		$r_2 \leftarrow \tilde{r}_1+r_2-r_1$,
		$\tilde{r}_2 \leftarrow r_2$),
	\eqref{eq: interpol help equation lem s>0}
		(with $\tilde{\delta} \leftarrow 0$),
	\eqref{eq: def of hat p1},
	\eqref{eq: def of hat p2},
	\eqref{eq: r1 r2tilde sum lem s>0},
	and the fact that 
	for all $r\in \R$ and all $t\in (-\infty, r]$
	there exists a $C \in (0,\infty)$ such that for all 
	$u \in H_{t}$ it holds that
	$\|u\|_{H_t} \leq C \|u\|_{H_r}$
	yields that there exist $C_1,C_2,C_3 \in (0,\infty)$
	such that for all 
	$u \in H_{r_1}$ and all $v \in H_{\tilde{r}_1 \vee \tilde{r}_2}$
	it holds that
	\begin{equation}
		\begin{split}
				&\|(D^{\alpha}_{\R^d} u) (D^{\beta}_{\R^d} v)\|_{B^{2s}_{p,2}(\Omega)}\\
			\leq{} 
				&C_1\big(
					\|
						u 
					\|_{H_{s+|\alpha|/2+d/4-d/(2\tilde{p}_1)}}
					\|
						v
					\|_{H_{|\beta|/2+d/4-d/(2p)+d/(2\tilde{p}_1)}} \\
					& \quad
					+\|
						u
					\|_{H_{|\alpha|/2+d/4-d/(2p)+d/(2\tilde{p}_2)}}
					\|
						v
					\|_{H_{s+|\beta|/2+d/4-d/(2\tilde{p}_2)}} 
				\big)\\
			\leq{} 
				&C_2 \big(
					\|u \|_{H_{r_1}}
					\|v \|_{H_{\tilde{r}_1+r_2-r_1}}
					+\|u \|_{H_{r_2}}
					\|v \|_{H_{\tilde{r}_1}} 
				\big)
			\leq
				C_3 \big(
					\|u \|_{H_{r_1}}
					\|v \|_{H_{\tilde{r}_2}}
					+\|u \|_{H_{r_2}}
					\|v \|_{H_{\tilde{r}_1}}
				\big).
		\end{split}
	\end{equation}
	{\it 3. case 
	$\delta >0,$
	$
			2r_2 
		=
			2s+|\alpha| +(0 \vee  (\tfrac {d}{2}-\tfrac {d}{p}))
	$.}
	Denote by 
	$
		\tilde{q} \in [2 \vee p, \tfrac{2p}{(2-p) \vee 0}] 
	$
	the real number satisfying that
	\begin{equation}
			\tilde{q}
		= 
			\sup \{
				q \in [2 \vee p, \tfrac{2p}{(2-p) \vee 0}] \colon  
					\tfrac{|\alpha|}{2}+\tfrac{d}{4}-\tfrac{d}{2q}+(\tfrac s2 \wedge \delta)
					\leq r_2
			\}.
	\end{equation}
	Then it follows from
	$
		2r_2=2s+|\alpha| +(0 \vee  (\tfrac {d}{2}-\tfrac {d}{p}))
	$
	that
	$\tilde{q}$ is well defined
	and that
	$\tilde{q}> 2 \vee p$
	and from the definition that
	\begin{equation}
	\label{eq: inequ r2 tilde q s>0}
			\tfrac{|\alpha|}{2}+\tfrac{d}{4}-\tfrac{d}{2\tilde{q}}+(\tfrac{s}{2} \wedge \delta )
		\leq 
			r_2.
	\end{equation}
	Moreover, we obtain
	from
	$
			\tfrac{|\alpha|}{2}+\tfrac{d}{4}-\tfrac{d}{2\tilde{q}}
			+s+\tfrac{|\beta|}{2}+\tfrac{d}{4}-\tfrac{d}{2p}+\tfrac{d}{2\tilde{q}}
		=
			\tilde{r}_1+r_2 -\delta
	$
	and from
	$
			s+\tfrac{|\beta|}{2}+\tfrac{d}{4}-\tfrac{d}{2p}+\tfrac{d ((2-p)\vee 0)}{4p}
		=
			s+\tfrac{|\beta|}{2}+(0\vee (\tfrac{d}{4}-\tfrac{d}{2p}))
		\leq
			\tilde{r}_1
	$
	that
	\begin{equation}
	\label{eq: inequ r1tilde tilde q s>0}
			s+\tfrac{|\beta|}{2}+\tfrac{d}{4}-\tfrac{d}{2p}+\tfrac{d}{2\tilde{q}}
		\leq
			\tilde{r}_1.
	\end{equation}
	Therefore \eqref{eq: basic estimate s>0}
	(with 
		$p_1 \leftarrow 2 \vee p$, 
		$
				\tfrac{1}{q_2} 
			\leftarrow 
				\tfrac{1}{p}-\tfrac{1}{2 \vee p}
		$, 
		$\tfrac{1}{p_2} \leftarrow \tfrac{1}{p}-\tfrac{1}{\tilde{q}}$, 
		$q_1 \leftarrow \tilde{q}$,
		$\delta_1 \leftarrow (\tfrac s2 \wedge \delta)$,
		and $\delta_2 \leftarrow \delta$),
	\eqref{eq: inequ r2 tilde q s>0},
	\eqref{eq: inequ r1tilde tilde q s>0},
	the fact that 
	for all $r\in \R$ and all $t\in (-\infty, r]$
	there exists a $C \in (0,\infty)$ such that for all 
	$u \in H_{t}$ it holds that
	$\|u\|_{H_t} \leq C \|u\|_{H_r}$,
	the fact that
	$
			2r_2 
		=
			2s+|\alpha| +(0 \vee  (\tfrac {d}{2}-\tfrac {d}{p})),
	$
	and \eqref{eq: r2 r1tilde sum s>0}
	verify that there exist $C_1,C_2,C_3 \in (0,\infty)$
	such that for all 
	$u \in H_{r_2}$ and all $v \in H_{\tilde{r}_1}$
	it holds that
	\begin{equation}
		\begin{split}
				&\|(D^{\alpha}_{\R^d} u) (D^{\beta}_{\R^d} v)\|_{B^{2s}_{p,2}(\Omega)}\\
			\leq{} 
				&C_1\big(
					\|
						u 
					\|_{H_{s+|\alpha|/2+d/4-d/(4 \vee (2p))}}
					\|
						v
					\|_{H_{|\beta|/2+d/4-d/(2p)+d/(4 \vee (2p)) + \delta \1_{[2,\infty)}(p)}}
				\\ &
					+\|
						u
					\|_{H_{|\alpha|/2+d/4-d/(2\tilde{q}) + ((s/2) \wedge \delta) \1_{\{\infty\}}(\tilde{q})}}
					\|
						v
					\|_{H_{s+|\beta|/2+d/4-d/(2p)+d/(2\tilde{q})}} 
				\big)\\
			\leq{} 
				&C_2 \big(
					\|u \|_{H_{r_2}}
					\|v \|_{
						H_{
							s+|\alpha|/2+|\beta|/2+d/2-d/(2p)
							-r_2 + \delta \1_{[2,\infty)}(p)
						}
					}
					+\|u \|_{H_{r_2}}
					\|v \|_{H_{\tilde{r}_1}} 
				\big) \\
			\leq{} 
				&C_3 \big(
					\|u \|_{H_{r_2}}
					\|v \|_{
						H_{
							\tilde{r}_1
						}
					}
					+\|u \|_{H_{r_2}}
					\|v \|_{H_{\tilde{r}_1}} 
				\big).
		\end{split}
	\end{equation}
	{\it 4. case 
	$\delta >0$ and
	$
			|\alpha| +(0 \vee  (\tfrac {d}{2}-\tfrac {d}{p}))
		<
			2r_2 
		<
			2s+|\alpha| +(0 \vee  (\tfrac {d}{2}-\tfrac {d}{p}))
	$.}
	First observe that there exists an $\tilde{\eps} \in (0,\infty)$ such that
	\begin{equation}
	\label{eq: def of tilde eps case 3}
		\begin{split}
				-\frac
					{(\delta-\tilde{\eps} \1_{[2,\infty)}(p)) (r_1-r_2)}
						{
							s+\tfrac{|\alpha|}{2}+(0 \vee (\tfrac{d}{4}-\tfrac{d}{2p}))-r_2
						}
			\leq
				-\frac
					{\delta (r_1-r_2)}
						{
							s+\tfrac{|\alpha|}{2}+(0 \vee (\tfrac{d}{4}-\tfrac{d}{2p}))-r_2
						}
				+ \eps \1_{[2,\infty)}(p).
		\end{split}
	\end{equation}
	Next
	denote by 
	$
		\hat{q} \in [2 \vee p, \tfrac{2p}{(2-p) \vee 0}] 
	$
	the real number satisfying that
	\begin{equation}
			\hat{q}
		= 
			\sup \{
				q \in [2 \vee p, \tfrac{2p}{(2-p) \vee 0}] \colon  
					\tfrac{|\alpha|}{2}+\tfrac{d}{4}-\tfrac{d}{2q}
					+((r_2-\tfrac{|\alpha|}{2} -(0 \vee  (\tfrac {d}{4}-\tfrac {d}{2p})))  \wedge \delta)
					\leq r_2
			\}.
	\end{equation}
	Then we get from
	$
		2r_2>|\alpha| +(0 \vee  (\tfrac {d}{2}-\tfrac {d}{p}))
	$
	that
	$\hat{q}$ is well defined
	and from the definition that
	\begin{equation}
	\label{eq: inequ r2 hat q s>0}
			\tfrac{|\alpha|}{2}+\tfrac{d}{4}-\tfrac{d}{2\hat{q}}
			+((r_2-\tfrac{|\alpha|}{2} -(0 \vee  (\tfrac {d}{4}-\tfrac {d}{2p})))  \wedge \delta)
		\leq 
			r_2.
	\end{equation}
	Furthermore, the fact that
	$
			\tfrac{|\alpha|}{2}+\tfrac{d}{4}-\tfrac{d}{2\hat{q}}
			+s+\tfrac{|\beta|}{2}+\tfrac{d}{4}-\tfrac{d}{2p}+\tfrac{d}{2\hat{q}}
		=
			\tilde{r}_1+r_2 -\delta
	$
	and the fact that
	$
			s+\tfrac{|\beta|}{2}+\tfrac{d}{4}-\tfrac{d}{2p}+\tfrac{d ((2-p)\vee 0)}{4p}
		=
			s+\tfrac{|\beta|}{2}+(0\vee (\tfrac{d}{4}-\tfrac{d}{2p}))
		\leq
			\tilde{r}_1
	$
	ensure that
	\begin{equation}
	\label{eq: inequ r1tilde hat q s>0}
			s+\tfrac{|\beta|}{2}+\tfrac{d}{4}-\tfrac{d}{2p}+\tfrac{d}{2\hat{q}}
		\leq
			\tilde{r}_1.
	\end{equation}
	Hence \eqref{eq: basic estimate s>0}
	(with 
		$p_1 \leftarrow 2 \vee p$, 
		$
				\tfrac{1}{q_2} 
			\leftarrow 
				\tfrac{1}{p}-\tfrac{1}{2 \vee p}
		$, 
		$\tfrac{1}{p_2} \leftarrow \tfrac{1}{p}-\tfrac{1}{\hat{q}}$, 
		$q_1 \leftarrow \hat{q}$,
		and $\delta_2 \leftarrow \tilde{\eps}$),
	Lemma \eqref{lem: another interpolation result}
		(with 
		$s_1 \leftarrow s+\tfrac{|\alpha|}{2}+(0 \vee (\tfrac{d}{4}-\tfrac{d}{2p}))$,
		$
				s_2 
			\leftarrow 
				\tfrac{|\beta|}{2}+(\tfrac{d}{4} \wedge (\tfrac{d}{2}-\tfrac{d}{2p})) 
				+ \tilde{\eps} \1_{[2,\infty)}(p)
		$,
		and 
		$
				\tilde{r}_2 
			\leftarrow 
				r
		$),
	\eqref{eq: interpol help equation lem s>0} 
		(with $\tilde{\delta} \leftarrow \delta -\tilde{\eps} \1_{[2,\infty)}(p)$)
	\eqref{eq: inequ r2 hat q s>0},
	\eqref{eq: inequ r1tilde hat q s>0},
	the fact that 
	for all $r\in \R$ and all $t\in (-\infty, r]$
	there exists a $C \in (0,\infty)$ such that for all 
	$u \in H_{t}$ it holds that
	$\|u\|_{H_t} \leq C \|u\|_{H_r}$,
	\eqref{eq: def of tilde eps case 3}
	and \eqref{eq: r1 r2tilde sum lem s>0}
	shows that for all 
	$
			\delta_1
				\in 
					\big (
						0,
						(r_2-\tfrac{|\alpha|}{2} -(0 \vee  (\tfrac {d}{4}-\tfrac {d}{2p})))  
						\wedge \delta
					\big )
	$
	and all
	$r \in \R$
	with
	$
			r
		=
			\tilde{r}_1+r_2-r_1
				-\frac
					{(\delta -\tilde{\eps} \1_{[2,\infty)}(p)) (r_1-r_2)}
					{s+\frac{|\alpha|}{2}+(0 \vee (\frac{d}{4}-\frac{d}{2p}))-r_2}
	$
	there exist $C_1,C_2,C_3 \in (0,\infty)$
	such that for all 
	$u \in H_{r_1}$ and all $v \in H_{\tilde{r}_1 \vee \tilde{r}_2}$
	it holds that
	\begin{equation}
		\begin{split}
				&\|(D^{\alpha}_{\R^d} u) (D^{\beta}_{\R^d} v)\|_{B^{2s}_{p,2}(\Omega)}\\
			\leq
				&C_1\big(
					\|
						u 
					\|_{H_{s+|\alpha|/2+d/4-d/(4 \vee (2p))}}
					\|
						v
					\|_{H_{|\beta|/2+d/4-d/(2p)+d/(4 \vee (2p)) + \tilde{\eps} \1_{[2,\infty)}(p)}}
				\\ &
					+\|
						u
					\|_{H_{|\alpha|/2+d/4-d/(2\hat{q}) + \delta_1 \1_{\{\infty\}}(\hat{q})}}
					\|
						v
					\|_{H_{s+|\beta|/2+d/4-d/(2p)+d/(2\hat{q})}} 
				\big)\\
			\leq{} 
				&C_2 \big(
					\|u \|_{H_{r_1}}
					\|v \|_{
						H_{
							r
						}
					}
					+
					\|u \|_{H_{r_2}}
					\|v \|_{\tilde{r}_1}
					+\|u \|_{H_{r_2}}
					\|v \|_{H_{\tilde{r}_1}} 
				\big) \\
			\leq{} 
				&C_3 \big(
					\|u \|_{H_{r_1}}
					\|v \|_{
						H_{
							\tilde{r}_2
						}
					}
					+\|u \|_{H_{r_2}}
					\|v \|_{H_{\tilde{r}_1}} 
				\big).
		\end{split}
	\end{equation}
		{\it 5. case 
		$\delta >0$ and
		$
				2r_2 
			\leq
				|\alpha| +(0 \vee  (\tfrac {d}{2}-\tfrac {d}{p}))
		$.}
		The Assumption $s>0$ and $p \in (1,\infty)$
		then ensure that there exists a
		$q^* \in  (2 \vee p, \tfrac{2p}{(2-p) \vee 0})$
		such that
	\begin{equation}
	\label{eq: inequ r2 star q s>0}
			\tfrac{|\alpha|}{2}+(0 \vee  (\tfrac {d}{4}-\tfrac {d}{2p}))
		<
			\tfrac{|\alpha|}{2}+\tfrac{d}{4}-\tfrac{d}{2q^*}
		<
			s+\tfrac{|\alpha|}{2}+(0 \vee  (\tfrac {d}{4}-\tfrac {d}{2p})).
	\end{equation}
	Moreover, there exists an $\tilde{\eps} \in (0,\infty)$ such that
	\begin{equation}
	\label{eq: def of tilde eps case 4}
		\begin{split}
				-\frac
					{(\delta-\tilde{\eps} \1_{[2,\infty)}(p)) (r_1-r_2)}
						{
							s+\tfrac{|\alpha|}{2}+(0 \vee (\tfrac{d}{4}-\tfrac{d}{2p}))-r_2
						}
			\leq
				-\frac
					{\delta (r_1-r_2)}
						{
							s+\tfrac{|\alpha|}{2}+(0 \vee (\tfrac{d}{4}-\tfrac{d}{2p}))-r_2
						}
				+ \eps \1_{[2,\infty)}(p).
		\end{split}
	\end{equation}
	Combining \eqref{eq: basic estimate s>0}
	(with 
		$p_1 \leftarrow 2 \vee p$,  
		$
				\tfrac{1}{q_2} 
			\leftarrow 
				\tfrac{1}{p}-\tfrac{1}{2 \vee p}
		$, 
		$\tfrac{1}{p_2} \leftarrow \tfrac{1}{p}-\tfrac{1}{q^*}$, 
		$q_1 \leftarrow q^*$,
		and $\delta_2 \leftarrow \tilde{\eps}$),
	Lemma \eqref{lem: another interpolation result}
		(with 
		$s_1 \leftarrow s+\tfrac{|\alpha|}{2}+(0 \vee (\tfrac{d}{4}-\tfrac{d}{2p}))$, 
		$
				s_2 
			\leftarrow 
				\tfrac{|\beta|}{2}+(\tfrac{d}{4} \wedge (\tfrac{d}{2}-\tfrac{d}{2p})) 
				+ \tilde{\eps} \1_{[2,\infty)}(p)
		$,
		$
				\tilde{r}_2 
			\leftarrow 
				\tilde{r}
		$
		and with
		$s_1 \leftarrow s+\tfrac{|\alpha|}{2}+\tfrac{d}{4}-\tfrac{d}{2q^*}$,
		$
				s_2 
			\leftarrow 
				\tfrac{|\beta|}{2}+\tfrac{d}{4}-\tfrac{d}{2p}+\tfrac{d}{2q^*}
		$,
		$
				\tilde{r}_2 
			\leftarrow 
				\hat{r}
		$),
	\eqref{eq: interpol help equation lem s>0} 
		(with 
		$\tilde{\delta} \leftarrow \delta -\tilde{\eps} \1_{[2,\infty)}(p)$
		and with $\tilde{\delta} \leftarrow \delta$),
	\eqref{eq: inequ r2 star q s>0},
	the fact that 
	for all $r\in \R$ and all $t\in (-\infty, r]$
	there exists a $C \in (0,\infty)$ such that for all 
	$u \in H_{t}$ it holds that
	$\|u\|_{H_t} \leq C \|u\|_{H_r}$,
	\eqref{eq: def of tilde eps case 4},
	and \eqref{eq: r1 r2tilde sum lem s>0}
	show that for all
	$\tilde{r}, \hat{r} \in \R$
	with
	$
			\tilde{r}
		=
			\tilde{r}_1+r_2-r_1
				-\frac
					{(\delta -\tilde{\eps} \1_{[2,\infty)}(p)) (r_1-r_2)}
					{s+\frac{|\alpha|}{2}+(0 \vee (\frac{d}{4}-\frac{d}{2p}))-r_2}
	$
	and with
	$
			\hat{r}
		=
			\tilde{r}_1+r_2-r_1
				-\frac
					{\delta (r_1-r_2)}
					{s+\frac{|\alpha|}{2}+\frac{d}{4}-\frac{d}{2q^*}-r_2}
	$
	there exist $C_1,C_2,C_3 \in (0,\infty)$
	such that for all 
	$u \in H_{r_1}$ and all $v \in H_{\tilde{r}_1}$
	it holds that
	\begin{equation}
		\begin{split}
				&\|(D^{\alpha}_{\R^d} u) (D^{\beta}_{\R^d} v)\|_{B^{2s}_{p,2}(\Omega)}\\
			\leq{} 
				&C_1\big(
					\|
						u 
					\|_{H_{s+|\alpha|/2+d/4-d/(4 \vee (2p))}}
					\|
						v
					\|_{H_{|\beta|/2+d/4-d/(2p)+d/(4 \vee (2p)) + \tilde{\eps} \1_{[2,\infty)}(p)}}
				\\ &
					+\|
						u
					\|_{H_{|\alpha|/2+d/4-d/(2q^*)}}
					\|
						v
					\|_{H_{s+|\beta|/2+d/4-d/(2p)+d/(2q^*)}} 
				\big)\\
			\leq{} 
				&C_2 \big(
					\|u \|_{H_{r_1}}
					\|v \|_{H_{\tilde{r}}}
					+
					\|u \|_{H_{r_2}}
					\|v \|_{H_{\tilde{r}_1}} 
					+\|u \|_{H_{r_1}}
					\|v \|_{H_{\hat{r}}}
					+
					\|u \|_{H_{r_2}}
					\|v \|_{\tilde{r}_1} 
				\big) \\
			\leq{} 
				&C_3 \big(
					\|u \|_{H_{r_1}}
					\|v \|_{
						H_{
							\tilde{r}_2
						}
					}
					+\|u \|_{H_{r_2}}
					\|v \|_{H_{\tilde{r}_1}} 
				\big).
		\end{split}
	\end{equation}
	Combining the 4 cases then finishes the proof of
	Lemma \ref{lem: sobolev holder inequality for positive s}.
\end{proof}
The next corollary combines 
Corollary \ref{cor: sobolev holder inequality with interpolation for Lp}
and Lemma \ref{lem: sobolev holder inequality for positive s}.
\begin{corollary}[Sobolev-H\"older inequality for non-negative Besov spaces]
	\label{cor: sobolev holder inequality for positive s}
	Assume the setting in Section \ref{ssec: application},
	let $ p\in (1,2]$,
	let $s \in [0,\infty)$,
	let $\eps \in (0,\infty)$,
	let $\delta \in [0,\infty)$,
	let $\alpha, \beta \in (\N_0)^d$,
	let $r_1,r_2, \tilde{r}_1,\tilde{r}_2 \in \R$
	satisfy that 
	\begin{equation}
	\label{eq: r2 r1tilde sum s>=0 cor}
			2r_2 + 2\tilde{r}_1 
		= 
			2s+|\alpha| +|\beta|+ d -\tfrac{d}{p} +2\delta,
	\end{equation}
	that
	\begin{equation}
	\label{eq: r1 r2tilde sum lem s>0 cor}
			2r_1 + 2\tilde{r}_2 
		\geq 
			2s+|\alpha| +|\beta|+d -\tfrac{d}{p}
			-\frac
				{
					2\delta \cdot (2r_1-2s-|\alpha|)
				}
				{
					2s+|\alpha|-2r_2
				} 
			+ \eps \1_{  (0,\infty) \times \{ 2 \}}((\delta,p)),
	\end{equation}
	that
	$
			2r_1 
		\geq
			2s+|\alpha|,
	$
	that
	$
			2\tilde{r}_1 
		\geq
			2s+|\beta|,
	$
	that
	$r_1 \geq r_2$,
	that either
	$
			2r_2
		\leq
			2s+|\alpha|
	$
	or
	$\delta =0$,
	that either
	$r_1 > r_2$
	or $\delta =0$
	and that 
	$
			\delta 
			\vee (
				(2r_1-2s-|\alpha| -\tfrac {d}{2}+\tfrac {d}{p}) 
				\cdot (2\tilde{r}_1-2s-|\beta| -\tfrac {d}{2}+\tfrac {d}{p})
			)
		>0.
	$
	Then there exists a $C \in (0, \infty)$ such that
	for all $u \in H_{r_1}$ and all $v \in H_{\tilde{r}_1 \vee \tilde{r}_2}$ 
	it holds that
	\begin{equation}
			\|(D^{\alpha}_{\R^d} u) (D^{\beta}_{\R^d} v)\|_{B^{2s}_{p,2}(\Omega)} 
		\leq
			C(
				\|u\|_{H_{r_1}} \|v\|_{H_{\tilde{r}_2}}
				+ \|u\|_{H_{r_2}} \|v\|_{H_{\tilde{r}_1}}
			).
	\end{equation}
\end{corollary}
\begin{proof}
First note that it follows from  
the Theorem on page 30 in Runst \& Sickel \cite{RunstSickel1996}, from
Proposition (iii) on page 14 in Runst \& Sickel \cite{RunstSickel1996} 
and from $p\leq 2$ that there exists a $C \in (0,\infty)$
such that for all $u \in L^p(\Omega)$ it holds that
\begin{equation}
\label{eq: Lp versus besov}
	\|u\|_{B^0_{p,2}(\Omega)} \leq C \|u\|_{L^p(\Omega)}.
\end{equation}
We will now divide the proof into 4 cases. \\
	{\it 1. case: 
		$
				s>0.
		$}
	This follows directly from Lemma \ref{lem: sobolev holder inequality for positive s}.\\
	{\it 2. case: 
		$s=0$, 
		$
				2r_1 
			>
				|\alpha|
		$,
		$
			\exists \tilde{\eps} \in (0,\eps] \colon
					2r_2
				\leq
					|\alpha| -\tilde{\eps}\1_{(1,2)}(p).
		$}
	This follows from Corollary \ref{cor: sobolev holder inequality with interpolation for Lp}
		(with $\eps \leftarrow \tilde{\eps}$)
	and from \eqref{eq: Lp versus besov}. \\
	{\it 3. case: 
		$s=0$, 
		$
				2r_1 
			\leq
				|\alpha|.
		$}
	Then the assumption 
	$
			2r_1 
		\geq
			2s+|\alpha|
	$
	ensures that
	$
			2r_1 
		=
			|\alpha|.
	$
	In addition 
	$
		\delta 
			\vee (
				(2r_1-2s-|\alpha| -\tfrac {d}{2}+\tfrac {d}{p}) 
				\cdot (2\tilde{r}_1-2s-|\beta| -\tfrac {d}{2}+\tfrac {d}{p})
			)
		>0
	$
	yields that 
	$( p=2 ) \Rightarrow (\delta >0)$
	and thus \eqref{eq: r1 r2tilde sum lem s>0 cor} shows that
	\begin{equation}
		\begin{split}
			2\tilde{r}_2 -|\beta|
		\geq{}
			&2s+|\alpha| +|\beta|+d -\tfrac{d}{p}
			-\frac
				{
					2\delta \cdot (2r_1-2s-|\alpha|)
				}
				{
					2s+|\alpha|-2r_2
				} 
			+ \eps \1_{  (0,\infty) \times \{ 2 \}}((\delta,p))
			-2r_1 -|\beta| \\
		={}
			&d -\tfrac{d}{p}
			+ \eps \1_{\{ 2 \}}(p)
		>
			0.
		\end{split}
	\end{equation}
	Therefore, \eqref{eq: Lp versus besov},
	Lemma \ref{lem: sobolev holder inequality without interpolation for Lp}
	(with 
	$r_1 \leftarrow 0$, 
	$r_2 \leftarrow 2\tilde{r}_2-|\beta|$,
	$u \leftarrow D^{\alpha}_{\R^d} u$,
	and with $v \leftarrow D^{\beta}_{\R^d} v$),
	and the fact that for all $r\in \R$ it holds that the
	$\| \cdot \|_{H_r}$-norm is equivalent
	to the $\| \cdot \|_{B^{2r}_{2,2}(\Omega)}$-norm on $H_r$
	imply that there exist $C_1,C_2,C_3,C_4 \in (0,\infty)$
	such that for all $u \in H_{r_1}$ 
	and all $v \in H_{\tilde{r}_2}$  it holds that
	\begin{equation}
		\begin{split}
				&\|(D^{\alpha}_{\R^d} u) (D^{\beta}_{\R^d} v)\|_{B^{0}_{p,2}(\Omega)} 
			\leq{} 
				C_1\|(D^{\alpha}_{\R^d} u) (D^{\beta}_{\R^d} v)\|_{L^{p}(\Omega)}  \\
			\leq
				&C_2(
					\|D^{\alpha}_{\R^d} u\|_{B^{2r_1-|\alpha|}_{2,2}(\Omega)}
					\|D^{\beta}_{\R^d} v\|_{B^{2\tilde{r}_2-|\beta|}_{2,2}(\Omega)}
				) 
			\leq
				C_3(
					\|u\|_{B^{2r_1}_{2,2}(\Omega)}
					\|v\|_{B^{2\tilde{r}_2}_{2,2}(\Omega)}
				)
			\leq
				C_4(
					\|u\|_{H_{r_1}}
					\|v\|_{H_{\tilde{r}_2}}
				).
		\end{split}
	\end{equation}	
	{\it 4. case: 
		$s=0$, 
		$
			\forall \tilde{\eps} \in (0,\eps] \colon
					2r_2
				>
					|\alpha| -\tilde{\eps}\1_{(1,2)}(p).
		$}
	The assumption then verifies that there exists a $\hat{\eps} \in (0,\eps]$
	such that
	$
			2r_2
		\geq
			|\alpha| +\hat{\eps}\1_{\{2\}}(p)
		>
			|\alpha| + \tfrac{d}{2} - \tfrac dp.
	$
	Furthermore the assumption
	$
			\delta 
			\vee (
				(2r_1-2s-|\alpha| -\tfrac {d}{2}+\tfrac {d}{p}) 
				\cdot (2\tilde{r}_1-2s-|\beta| -\tfrac {d}{2}+\tfrac {d}{p})
			)
		>0
	$
	assures that $\delta>0$ or $2\tilde{r}_1>|\beta| +\tfrac {d}{2}-\tfrac {d}{p}.$
	Thus \eqref{eq: Lp versus besov},
	Lemma \ref{lem: sobolev holder inequality without interpolation for Lp}
		(with $r_1 \leftarrow 2r_2-|\alpha|$,
		$r_2 \leftarrow 2\tilde{r}_1-|\beta|$, 
		$\eps \leftarrow \min \big( (\{\delta\} \cap (0,\infty)) \cup\{1\} \big)$,
		$u \leftarrow D^{\alpha}_{\R^d} u$,
		and with $v \leftarrow D^{\beta}_{\R^d} v$),
	\eqref{eq: r2 r1tilde sum s>=0 cor},
	and the fact that for all $r\in \R$ it holds that the
	$\| \cdot \|_{H_r}$-norm is equivalent
	to the $\| \cdot \|_{B^{2r}_{2,2}(\Omega)}$-norm on $H_r$
	show that there exist $C_1,C_2,C_3,C_4 \in (0,\infty)$
	such that for all $u \in H_{r_1}$ 
	and all $v \in H_{\tilde{r}_2}$  it holds that
	\begin{equation}
		\begin{split}
				&\|(D^{\alpha}_{\R^d} u) (D^{\beta}_{\R^d} v)\|_{B^{0}_{p,2}(\Omega)} 
			\leq{} 
				C_1\|(D^{\alpha}_{\R^d} u) (D^{\beta}_{\R^d} v)\|_{L^{p}(\Omega)} \\
			\leq{}
				&C_2(
					\|D^{\alpha}_{\R^d} u\|_{B^{2r_2-|\alpha|}_{2,2}(\Omega)}
					\|D^{\beta}_{\R^d} v\|_{B^{2\tilde{r}_1-|\beta|}_{2,2}(\Omega)}
				) 
			\leq{} 
				C_3(
					\|u\|_{B^{2r_2}_{2,2}(\Omega)}
					\|v\|_{B^{2\tilde{r}_1}_{2,2}(\Omega)}
				) \\
			\leq{}
				&C_4(
					\|u\|_{H_{r_2}}
					\|v\|_{H_{\tilde{r}_1}}
				).
		\end{split}
	\end{equation}
	Combining the 4 cases then finishes the proof of
	Corollary \ref{cor: sobolev holder inequality for positive s}.
\end{proof}
The next corollary is an implication of 
Corollary \ref{cor: sobolev holder inequality for positive s}.
It will be needed in the proof of Lemma \ref{l: norm equivalence Burger}
and Lemma \ref{l: norm equivalence Navier}.
\begin{corollary}
\label{cor: F assumption}
	Assume the setting in Section \ref{ssec: application},
	let $s,\gamma \in [0,\infty)$,
	$\eps \in (0,\infty)$,
	$\alpha \in \N_0^d$,
	$\alpha_1,\alpha_2 \in (0,\nicefrac 12)$, 
	and let
	$\vartheta \in [\nicefrac 12, \infty)$
	satisfy that
	$|\alpha|\leq 1$,
	that
	$\gamma > \tfrac d4 -\tfrac 12$,
	and that
	$\alpha_1 \vee \alpha_2 < \gamma-\tfrac d4 +\tfrac 12$.
	Then there exist a $C \in (0,\infty)$
	such that
	for all $u,v \in H_{\gamma}$
	it holds that
	\begin{align}
	\label{eq: h bound for F Cor}
			&\| u^2\|_{B^{2\gamma+2\vartheta-1}_{2,2}(\Omega)}
		\leq 
			C \|u\|^2_{H_{\gamma+\vartheta}}\\
	\label{eq: h+1/2 bound for F Cor}
			&\| (D_{\R^d}^\alpha u)  u\|_{B^{2\gamma+2\vartheta-1}}
		\leq 
			C	\|u\|_{H_{\gamma}} 
			\|u\|_{H_{\gamma+\vartheta+1/2-\alpha_1}} \\
	\label{eq: Lipschitz continuity of F Cor}
			&\|u v \|_{B^{2\gamma}_{2,2}(\Omega)}
		\leq
				C \|u\|_{H_{\gamma}}
					\|v\|_{H_{\gamma+1/2-\alpha_2}}, \\ 
		\label{eq: h+1/2 continuity of F Cor2}
			&\|  (D_{\R^d}^\alpha u) v\|_{B^{2s}_{2,2}(\Omega)} 
		\leq 
			\begin{cases}
				C 
					\|u\|_{H_{s/2+|\alpha|/4+d/8}} 
					\|v\|_{H_{s/2+|\alpha|/4+d/8}}
					& \textrm{if } 4s < d-2|\alpha|\\
					C\|u\|_{H_{s+|\alpha|/2 +\eps \1_{\{4s\}}(d-2|\alpha|)}} 
					\|v\|_{H_{s+|\alpha|/2 +\eps \1_{\{4s\}}(d-2|\alpha|)}}
					& \textrm{if } 4s \geq d-2|\alpha|
			\end{cases}.
	\end{align}
	Moreover if $\gamma+\vartheta-1\geq0$ then it holds 
	for all $u\in H_{\gamma}$ that
	\begin{equation}
	\label{eq: h bound for F Cor2}
			\| (D_{\R^d}^\alpha u) u\|_{B^{2\gamma+2\vartheta-2}_{2,2}(\Omega)}
		\leq 
			C \|u\|^2_{H_{\gamma+\vartheta}}
	\end{equation}
	and if $2\gamma-1 \geq 0$ and $|\alpha|=1$ then it hold 
	for all $u,v \in H_{\gamma}$ that
	\begin{equation}
		\label{eq: Lipschitz continuity of F Cor2}
				\|(D_{\R^d}^\alpha u) u - (D_{\R^d}^\alpha v) v\|_{B^{2\gamma-1}_{2,2}(\Omega)}
			\leq
				C (\|u\|_{H_\gamma}+\|v\|_{H_\gamma}) 
					\|u-v\|_{H_{\gamma+1/2-\alpha_2}}.
	\end{equation}
\end{corollary}
\begin{proof}
	We will first proof \eqref{eq: h bound for F Cor}. Therefore note that
	$\gamma \in [0,\infty)$ and $\vartheta \in [\nicefrac 12, \infty)$
	imply that $2\gamma +2\vartheta -1 \geq 0$.
	Moreover, the assumption $\gamma>\tfrac d4 -\tfrac 12$ yields that
	$\tfrac d4-\tfrac 12 < \gamma <  \gamma +\vartheta$.
	Therefore 
	Corollary \ref{cor: sobolev holder inequality for positive s}
	(with
		$p \leftarrow 2$, $\alpha \leftarrow 0$, $\beta \leftarrow0$, 
		$\delta \leftarrow 0$, $s \leftarrow \gamma + \vartheta -\tfrac 12$,
		$r_1 \leftarrow \gamma + \vartheta$, $\tilde{r}_1 \leftarrow \gamma + \vartheta$,
		$r_2 \leftarrow \tfrac d4 -\tfrac 12$, 
		$\tilde{r}_2 \leftarrow \tfrac d4 -\tfrac 12$,
		$u \leftarrow u$, 
		$v \leftarrow u$) 
	and the fact that $\tfrac d4-\tfrac 12 <  \gamma +\vartheta$
	yields that there exist $C_1,C_2 \in (0,\infty)$ such that
	for all $u \in H_{\gamma}$ it holds that
	\begin{equation}
			\|u^2\|_{B^{2\gamma+2\vartheta-1}_{2,2}(\Omega)}
		\leq
			C_1
			\|u\|_{H_{\gamma+\vartheta}}
			\|u\|_{H_{d/4-1/2}}
		\leq
			C_2
			\|u\|^2_{H_{\gamma+\vartheta}}.
	\end{equation}
	Next we will show \eqref{eq: h+1/2 bound for F Cor}.
Note that we get from the assumption $\alpha_1 < \nicefrac 12$ 
and and from the assumption $|\alpha| \leq 1$
that 
$
		\gamma+\vartheta+\tfrac 12-\alpha_1
	>
		(\gamma+\vartheta-\tfrac 12)+\tfrac {|\alpha|}{2}
$.
In addition it follows from $\alpha_1 < \gamma-\tfrac d4+\tfrac 12$,
from $|\alpha| \leq 1$, 
and from $\alpha_1 < \tfrac 12$ that
$
		\tfrac d4 -1 +\tfrac {|\alpha|}{2} +\alpha_1
	< 
		\gamma 
	< 
		\gamma + \vartheta + \tfrac 12 - \alpha_1
$.
Hence, Corollary \ref{cor: sobolev holder inequality for positive s}
	(with
		$p \leftarrow 2$, $\alpha \leftarrow \alpha$, $\beta \leftarrow 0$, 
		$\delta \leftarrow 0$, $s \leftarrow \gamma + \vartheta -\tfrac 12$,
		$r_1 \leftarrow \gamma + \vartheta+\tfrac 12 -\alpha_1$, 
		$\tilde{r}_1 \leftarrow \gamma + \vartheta+\tfrac 12 -\alpha_1$,
		$r_2 \leftarrow \tfrac d4 -1+\tfrac{|\alpha|}{2}+\alpha_1$, 
		$\tilde{r}_2 \leftarrow \tfrac d4 -1+\tfrac{|\alpha|}{2}+\alpha_1$,
		$u \leftarrow u$, 
		$v \leftarrow u$) 
	and the assumption $\alpha_1 < \gamma-\tfrac d4+\tfrac 12$
	establishes that there exist $C_1,C_2 \in (0,\infty)$ such that
	for all $u\in H_\gamma$ it holds that
	\begin{equation}
		\begin{split}
				&\|(D_{\R^d}^\alpha u) u \|_{B^{2\gamma +2\vartheta -1}_{2,2}(\Omega)}
			\leq{} 
				C_1
				\|u\|_{H_{\gamma + \vartheta+1/2 -\alpha_1}} 
				\|u\|_{H_{d/4 -1+|\alpha|/2+\alpha_1}} \\
			\leq{}
				&C_2
				\|u\|_{H_{\gamma + \vartheta+1/2 -\alpha_1}} 
				\|u\|_{H_{d/4 -1/2+\gamma-d/4+1/2}} 
			={} 
				C_2
				\|u\|_{H_{\gamma + \vartheta+1/2 -\alpha_1}} 
				\|u\|_{H_{\gamma}}.
		\end{split}
	\end{equation}
	Next we proof \eqref{eq: Lipschitz continuity of F Cor}. 
	Note that it holds that
	$
			(2\gamma +1 -2\alpha_2)+(\gamma +\tfrac d4 -\tfrac 12+\alpha_2)
		=
			2\gamma +\tfrac d2+(\gamma-\tfrac d4+\tfrac 12-\alpha_2).
	$
	Moreover, we have that
	$
			2\gamma + (2\gamma+1-2\alpha_2)
		=
			2\gamma +\tfrac d2 + (2\gamma -\tfrac d2+1-2\alpha_2).
	$
	In addition, the assumption $\alpha_2 < \tfrac 12$
	assures that
	$2\gamma+1-2\alpha_2 >2\gamma$.
	Furthermore, the assumption $\gamma-\tfrac d4 +\tfrac 12 >\alpha_2 $
	verifies that
	$
		\gamma > \tfrac{\gamma}{2} +\tfrac d8 -\tfrac 14+\tfrac{\alpha_2}{2}.
	$
	Thus
	Corollary \ref{cor: sobolev holder inequality for positive s}
	(with
		$p \leftarrow 2$, $\alpha \leftarrow 0$, $\beta \leftarrow 0$, 
		$\delta \leftarrow \tfrac{\gamma}{2} -\tfrac d8+\tfrac 14-\tfrac{\alpha_2}{2}$, 
		$\eps \leftarrow 2\gamma -\tfrac d2+1-2\alpha_2$,
		$s \leftarrow \gamma$,
		$r_1 \leftarrow \gamma $, 
		$\tilde{r}_1 \leftarrow \gamma +\tfrac 12 -\alpha_2$,
		$r_2 \leftarrow \tfrac{\gamma}{2} +\tfrac d8 -\tfrac 14+\tfrac{\alpha_2}{2}$,
		$\tilde{r}_2 \leftarrow \gamma+\tfrac 12-\alpha_2$)
	and the fact that
	$
		\gamma > \tfrac{\gamma}{2} +\tfrac d8 -\tfrac 14+\tfrac{\alpha_2}{2}
	$
	show that there exist $C_1,C_2 \in (0,\infty)$
	such that for all $u,v \in H_\gamma$ it holds that
	\begin{equation}
		\begin{split}
				&\|uv\|_{B^{2\gamma}_{2,2}(\Omega)}
			\leq
				C_1(
					\|u\|_{H_\gamma} \|v\|_{H_{\gamma+1/2-\alpha_2}}
					+\|u\|_{H_{\gamma/2 + d/8 - 1/4+\alpha_2/2}} 
					\|v\|_{H_{\gamma +1/2 -\alpha_2}}
				) \\
			\leq{}
				&C_2
				\|u\|_{H_\gamma} \|v\|_{H_{\gamma+1/2-\alpha_2}}.
		\end{split}
	\end{equation}
	We will divide the proof of \eqref{eq: h+1/2 continuity of F Cor2} into 2 cases.\\ 
	{\it 1. case $4s \geq d - 2|\alpha|$.}
	We then get that 
	$
			2(s + \tfrac {|\alpha|}{2} + \eps \1_{\{4s\}}(d-2|\alpha|))
			+2(\tfrac s2+ \tfrac {|\alpha|}{4} +\tfrac d8)
		=
			2s +|\alpha|+ \tfrac d2 
			+(s-\tfrac d4+ \tfrac {|\alpha|}{2} + 2\eps \1_{\{4s\}}(d-2|\alpha|))
	$,
	that
	$
			s+ \tfrac {|\alpha|}{2} +\tfrac d4
		\leq
			2s + |\alpha|	
	$,
	and that
	$
			s + \tfrac {|\alpha|}{2} + \eps \1_{\{4s\}}(d-2|\alpha|)
		>
			\tfrac s2+ \tfrac {|\alpha|}{4} +\tfrac d8.
	$
	Therefore,
	Corollary \ref{cor: sobolev holder inequality for positive s}
	(with
		$p \leftarrow 2$, $\alpha \leftarrow \alpha$, $\beta \leftarrow 0$, 
		$
				\delta 
			\leftarrow 
				\tfrac s2 -\tfrac d8+\tfrac{|\alpha|}{4}+\eps \1_{\{4s\}}(d-2|\alpha|)
		$, 
		$
				\eps 
			\leftarrow 
				s -\tfrac d4+\tfrac{|\alpha|}{2} + 2 \eps \1_{\{4s\}}(d-2|\alpha|)
		$,
		$s \leftarrow s$,
		$r_1 \leftarrow s + \tfrac {|\alpha|}{2}+ \eps \1_{\{4s\}}(d-2|\alpha|)$, 
		$\tilde{r}_1 \leftarrow s + \tfrac {|\alpha|}{2}+ \eps \1_{\{4s\}}(d-2|\alpha|)$,
		$r_2 \leftarrow \tfrac s2 + \tfrac {|\alpha|}{4}+ \tfrac d8$,
		$\tilde{r}_2 \leftarrow \tfrac s2 + \tfrac {|\alpha|}{4}+ \tfrac d8$)
	and the fact that
	$
			s + \tfrac {|\alpha|}{2} + \eps \1_{\{4s\}}(d-2|\alpha|)
		>
			\tfrac s2+ \tfrac {|\alpha|}{4} +\tfrac d8
	$
	yield that there exist $C_1,C_2 \in (0,\infty)$
	such that for all 
	$
		u,v \in 
			H_{
				\gamma
			}
	$
	it holds that 
	\begin{equation}
		\begin{split}
				&\| (D_{\R^d}^\alpha u) v\|_{B^{2s}_{2,2}(\Omega)} \\
			\leq{} 
				&C_1
					(
						\|u\|_{H_{s + |\alpha|/2+ \eps \1_{\{4s\}}(d-2|\alpha|)}} 
						\|v\|_{H_{s/2 + |\alpha|/4+ d/8}} 
					+
						\|v\|_{H_{s + |\alpha|/2+ \eps \1_{\{4s\}}(d-2|\alpha|)}} 
						\|u\|_{H_{s/2 +|\alpha|/4+ d/8}} 
					) \\
			\leq{} 
				&C_2
					\|u\|_{H_{s + |\alpha|/2+ \eps \1_{\{4s\}}(d-2|\alpha|)}} 
					\|v\|_{H_{s + |\alpha|/2+ \eps \1_{\{4s\}}(d-2|\alpha|)}}.
		\end{split}
	\end{equation}
	{\it 2. case $4s < d - 2|\alpha|$.}
	This shows that
	$
			\tfrac s2 + \tfrac {|\alpha|}{4} + \tfrac d8
		>
			s +\tfrac {|\alpha|}{2}
	$
	and thus 
	Corollary \ref{cor: sobolev holder inequality for positive s}
	(with
		$p \leftarrow 2$, $\alpha \leftarrow \alpha$, $\beta \leftarrow 0$, 
		$\delta \leftarrow 0$, 
		$s \leftarrow s$,
		$r_1 \leftarrow \tfrac s2 + \tfrac {|\alpha|}{4} + \tfrac d8$, 
		$\tilde{r}_1 \leftarrow \tfrac s2 + \tfrac {|\alpha|}{4} + \tfrac d8$,
		$r_2 \leftarrow \tfrac s2 + \tfrac {|\alpha|}{4} + \tfrac d8$,
		$\tilde{r}_2 \leftarrow \tfrac s2 + \tfrac {|\alpha|}{4} + \tfrac d8$)
	verifies that there exist $C_1,C_2 \in (0,\infty)$
	such that for all 
	$
		u,v \in 
			H_{
				\gamma
			}
	$
	it holds that
	\begin{equation}
				\| (D_{\R^d}^\alpha u) v\|_{B^{2s}_{2,2}(\Omega)}
			\leq 
				C_1
					\|u\|_{H_{s/2 + |\alpha|/4 + d/8}} 
					\|v\|_{H_{s/2 + |\alpha|/4 + d/8}}.
	\end{equation}
	For the proof of \eqref{eq: h bound for F Cor2} we assume that 
	$\gamma+\vartheta-1 \geq 0$.
	Observe that
	the assumption $\gamma < \tfrac d4 - \tfrac 12$ 
	and the assumption that
	$ |\alpha| \leq 1$
	ensure that
	$\tfrac d2 -2+|\alpha| < 2\gamma <2\gamma + 2\vartheta$. Thus
	Corollary \ref{cor: sobolev holder inequality for positive s}
	(with
		$p \leftarrow 2$, $\alpha \leftarrow \alpha$, $\beta \leftarrow 0$, 
		$\delta \leftarrow 0$, 
		$s \leftarrow \gamma +\vartheta-1$,
		$r_1 \leftarrow \gamma +\vartheta$, 
		$\tilde{r}_1 \leftarrow \gamma +\vartheta$,
		$r_2 \leftarrow \tfrac d4 -1 +\tfrac {|\alpha|}{2}$,
		$\tilde{r}_2 \leftarrow \tfrac d4 -1 +\tfrac {|\alpha|}{2}$,
		$u \leftarrow u$, 
		$v \leftarrow u$)
	and the fact that $\tfrac d2 -2 +|\alpha| <2\gamma + 2\vartheta$
	imply that there exist $C_1,C_2 \in (0,\infty)$ such that for all
	$u \in H_\gamma$ it holds that
	\begin{equation}
				\| (D_{\R^d}^\alpha u) u\|_{B^{2\gamma+2\vartheta-2}_{2,2}(\Omega)}
			\leq 
				C_1
					\|u\|_{H_{\gamma +\vartheta}} 
					\|u\|_{H_{d/4-1 +|\alpha|/2}} 
			\leq
				C_2
					\|u\|_{H_{\gamma +\vartheta}} 
					\|u\|_{H_{\gamma +\vartheta}}.
	\end{equation}
	Next assume that $2\gamma-1 \geq 0$
	and that $|\alpha|=1$.
	Then it holds that
	$
			(2\gamma +1 -2\alpha_2)+(\gamma +\tfrac d4 -\tfrac 12+\alpha_2)
		=
			(2\gamma -1)+1+\tfrac d2+(\gamma-\tfrac d4+\tfrac 12-\alpha_2).
	$
	Furthermore, we have that
	$
			2\gamma + (2\gamma+1-2\alpha_2)
		=
			(2\gamma -1)+1+\tfrac d2 + (2\gamma -\tfrac d2+1-2\alpha_2).
	$
	In addition, the assumption $\alpha_2 < \tfrac 12$
	verifies that
	$2\gamma+1-2\alpha_2 >2\gamma-1$.
	Moreover, it follows from the assumption $\gamma-\tfrac d4 +\tfrac 12 >\alpha_2 $
	that
	$
		\gamma > \tfrac{\gamma}{2} +\tfrac d8 -\tfrac 14+\tfrac{\alpha_2}{2}.
	$
	Thus
	Corollary \ref{cor: sobolev holder inequality for positive s}
	(with
		$p \leftarrow 2$, $\alpha \leftarrow \alpha$, $\beta \leftarrow 0$, 
		$\delta \leftarrow \tfrac{\gamma}{2} -\tfrac d8+\tfrac 14-\tfrac{\alpha_2}{2}$, 
		$\eps \leftarrow 2\gamma -\tfrac d2+1-2\alpha_2$,
		$s \leftarrow \gamma-\tfrac 12$,
		$r_1 \leftarrow \gamma $, 
		$\tilde{r}_1 \leftarrow \gamma +\tfrac 12 -\alpha_2$,
		$r_2 \leftarrow \tfrac{\gamma}{2} +\tfrac d8 -\tfrac 14+\tfrac{\alpha_2}{2}$,
		$\tilde{r}_2 \leftarrow \gamma+\tfrac 12-\alpha_2$,
		$u \leftarrow u$, 
		$v \leftarrow u-v$)
	and the fact that
	$
		\gamma > \tfrac{\gamma}{2} +\tfrac d8 -\tfrac 14+\tfrac{\alpha_2}{2}
	$
	show that there exist $C_1,C_2 \in (0,\infty)$
	such that for all $u,v \in H_\gamma$ it holds that
	\begin{equation}
	\label{eq: (vi) first part}
		\begin{split}
				&\|(D^\alpha_{\R^d} u) (u-v)\|_{B^{2\gamma-1}_{2,2}(\Omega)} \\
			\leq{} 
				&C_1(
					\|u\|_{H_\gamma} \|u-v\|_{H_{\gamma+1/2-\alpha_2}}
					+\|u\|_{H_{\gamma/2 + d/8 - 1/4+\alpha_2/2}} 
					\|u-v\|_{H_{\gamma +1/2 -\alpha_2}}
				) \\
			\leq{} 
				&C_2
				\|u\|_{H_\gamma} \|u-v\|_{H_{\gamma+1/2-\alpha_2}}.
		\end{split}
	\end{equation}
	Furthermore, the assumption $\alpha_2 < \tfrac 12$
	ensures that $2\gamma +1 - 2\alpha_2 > (2\gamma -1) +1$ 
	and the assumption $\alpha_2<\gamma +\tfrac 12 -\tfrac d4$
	ensures that $\gamma > \tfrac d4 -\tfrac 12 +\alpha_2$.
	Therefore,
	Corollary \ref{cor: sobolev holder inequality for positive s}
	(with
		$p \leftarrow 2$, $\alpha \leftarrow 0$, $\beta \leftarrow \alpha$, 
		$\delta \leftarrow 0$, 
		$s \leftarrow \gamma-\tfrac 12$,
		$r_1 \leftarrow \gamma $, 
		$\tilde{r}_1 \leftarrow \gamma +\tfrac 12 -\alpha_2$,
		$r_2 \leftarrow \tfrac d4 -\tfrac 12+\alpha_2$,
		$\tilde{r}_2 \leftarrow \tfrac d4$,
		$u \leftarrow v$, 
		$v \leftarrow u-v$)
	and the fact that
	$
		\gamma > \tfrac d4 -\tfrac 12 +\alpha_2
	$
	yield that there exist $C_1,C_2 \in (0,\infty)$
	such that for all $u,v \in H_\gamma$ it holds that
	\begin{equation}
	\label{eq: (vi) second part}
		\begin{split}
				\|v \, D^\alpha_{\R^d}(u-v)\|_{B^{2\gamma-1}_{2,2}(\Omega)}
			\leq{} 
				&C_1(
					\|v\|_{H_\gamma} \|u-v\|_{H_{d/4}}
					+\|v\|_{H_{d/4 -1/2+\alpha_2}} 
					\|u-v\|_{H_{\gamma +1/2 -\alpha_2}}
				) \\
			\leq{} 
				&C_2
				\|v\|_{H_\gamma} \|u-v\|_{H_{\gamma+1/2-\alpha_2}}.
		\end{split}
	\end{equation}
	Combining \eqref{eq: (vi) first part} and \eqref{eq: (vi) second part} proves
	\eqref{eq: Lipschitz continuity of F Cor2}
	and thus finishes the proof of 
	Corollary \ref{cor: F assumption}.
\end{proof}
\section{Estimates for Nemytskij operators on Besov spaces}
\label{ssec: Estimates for Nemytskij operator}
In this subsection we establish bounds for Nemytskij operators.
For 
a more detailed overview on this subject
we refer to
Chapter 5 in Runst \& Sickel \cite{RunstSickel1996}.
\begin{lemma}[Upper bound for the Slobodeckij semi-norm of Nemytskij operators]
\label{lem: b boundedness inequality}
	Let $d,m \in \N$, 
	let $\Omega \subseteq \R^d$ be a bounded Lipschitz domain,
	denote by $B(r) \subseteq \R^m$, $r \in [0,\infty]$,
	the ball satisfying that for all $r \in [0,\infty]$ it holds
	that
	$
			B(r)
		= 
			\{ 
				y \in \R^m \colon 
					\| y\|_{\R^m} \leq r
			\}
	$,
	let $K \in (0,\infty)$,
	$p \in [0,\infty]$,
	$q \in [2,\infty]$,
	$s \in (0,1)$,
	$\eps \in (0, 1-s)$,
	satisfy that
	$pq \geq 2 \cdot \1_{ (0,\infty)}(p)$,
	let
	$b \in \C_{\R^d \times \R^m}^1( \Omega \times \R^m, \R)$,
	$K_1\colon [0,\infty] \times L^0(\Omega,\R^m) \to [0, \infty]$,
	satisfy that for all $x \in \Omega$, $y \in \R^m$, $r \in [0,\infty]$,
	and all $u \in L^0(\Omega,\R^m)$
	it holds that
	$
			\|D_{\R^d \times \R^m} b(x,y)\|_{(\R^d \times \R^m)'} 
		\leq 
			K (\|y\|_{\R^m} +1)^{p}
	$
	and that
	\begin{equation}
			K_1(r,u)
		=
		\begin{cases}
			1 & \textrm{ if } r= 0 \\
			\|u\|^{2r}_{B^{d/2-d/(r q) + \eps \1_{\{\infty\}}(q)}_{2,2}(\Omega,\R^m)} +1	
					& \textrm{ if } r \in (0,\infty) \\
					\|
						D_{\R^d \times \R^m} b
					\|^2_{
						\C(
							\Omega
							\times 
							B(\|u \|_{L^\infty(\Omega,\R^m)}),
							\| \cdot \|_{(\R^d \times \R^m)'}
						)
					}
					& \textrm{ if } r =\infty. 
		\end{cases}
	\end{equation}
	Then there exists a
	$C \in (0,\infty)$ such that 
	for all $u \in L^0(\Omega,\R^m)$,
	it holds that
	\begin{equation}
				\int_{\Omega} \int_{\Omega}
				\frac
					{
						|
							b( x ,u(x)) -b( y ,u(y)) 
						|^2
					}
					{\|x-y\|_{\R^d}^{2s+d}}
			\ud x \ud y
		\leq
			C \, K_1(p,u) \cdot 
				(\|u\|^2_{B^{s+\eps \1_{[2,\infty)}(q)+d/q}_{2,2}(\Omega,\R^m)} +1).
	\end{equation}
 \end{lemma}
\begin{proof}
We will divide the proof into 3 cases. \\
{\it 1. case $p=0$.}
First note that we get from
the assumption that for all
$x \in \Omega$ and all $y\in \R^m$ it holds that
$\|D_{\R^d \times \R^m} b(x,y)\|_{(\R^d \times \R^m)'} \leq K (\|y\|_{\R^m} +1)^p$,
from the fact that $\Omega$ is bounded,
from Proposition 2.28 in Mitrea \cite{Mitrea2013},
and
from the fact that 
for all $r\in \R$ and all $t\in (-\infty, r]$
there exists a $C \in (0,\infty)$ such that for all 
$u \in B^t_{2,2}(\Omega,\R^m)$ it holds that
$\|u\|_{B^t_{2,2}(\Omega,\R^m)} \leq C \|u\|_{B^r_{2,2}(\Omega,\R^m)}$,
that there exist $C_1, C_2 \in (0,\infty)$
such that for all 
$u \in L^0(\Omega,\R^m)$
it holds that
\begin{align}
		\nonumber
			&\int_{\Omega} \int_{\Omega}
				\frac
					{
						|
							b( x ,u(x)) -b( y ,u(y)) 
						|^2
					}
					{\|x-y\|_{\R^d}^{2s+d}}
			\ud x \ud y \\
		\leq{} \nonumber
			&\int_{\Omega} \int_{\Omega} \int_0^1 
				\tfrac
					{1}
					{\|x-y\|_{\R^d}^{2s+d}}
						\Big |
							\Big \langle 
								(D_{\R^d \times \R^m} b) \big(tx + (1-t)y,t u(x) + (1-t) u(y) \big),
					\\ & \nonumber \qquad
								(x-y, u(x)-u(y) ),
							\Big \rangle_{(\R^d \times \R^m)',\R^d \times \R^m}
						\Big |^2
			\ud t \ud x \ud y \\
	\label{eq: bound assertion for p=0}
		\leq{} 
			&\int_{\Omega} \int_{\Omega} \int_0^1
				\frac
					{
							\| 
								(D_{\R^d \times \R^m} b) \big(tx + (1-t)y,t u(x) + (1-t) u(y) \big)
							\|^2_{(\R^d \times \R^m)'}
					}
					{\|x-y\|_{\R^d}^{2s+d}}  
		\\ & \nonumber \qquad 
				\cdot (\|x-y\|^2_{\R^d}+\|u(x)-u(y)\|_{\R^m}^2)
			\ud t \ud x \ud y \\
		\leq{} 
			&K^2 \nonumber
			\int_{\Omega} \int_{\Omega} \int_0^1
				\frac
					{
						1
					}
					{\|x-y\|_{\R^d}^{2(s-1)+d}} 
			+
				\frac
					{
						\|u(x)-u(y)\|_{\R^m}^2
					}
					{\|x-y\|_{\R^d}^{2s+d}} 
			\ud t \ud x \ud y \\
		\leq{} \nonumber
			&C_1 (
				1+\|u\|^2_{B^{s}_{2,2}(\Omega,\R^m)} 
			)
		\leq
			C_2 (
				1+\|u\|^2_{B^{s+\eps \1_{[2,\infty)}(q)+d/q}_{2,2}(\Omega,\R^m)}
			) \\
		={} \nonumber
			&C_2 \, K_1(0,u) \cdot 
				(\|u\|^2_{B^{s+\eps \1_{[2,\infty)}(q)+d/q}_{2,2}(\Omega,\R^m)} +1).
\end{align}
{\it 2. case $p\in (0,\infty)$.}
Note that we get from Hölder's inequality,
the fact that $\Omega$ is bounded,
from Proposition 2.28 in Mitrea \cite{Mitrea2013},
and from the Sobolev inequality
(see, e.g., the Theorem on page 31 in Runst \& Sickel \cite{RunstSickel1996})
that for all $\delta \in (0, 1-s)$
and all $q_1,q_2 \in (2,\infty)$ with $1/q_1 + 1/q_2 =1/2$
there exist $ C_1, C_2, C_3 \in (0,\infty)$ such that
for all 
$u \in L^0(\Omega,\R^m)$
and all $v \in L^0(\Omega,\R^m)$ 
it holds that
	\begin{equation}
	\label{eq: holder estimate with p,q frac s}
		\begin{split}
				&\int_{\Omega} \int_{\Omega}
					\frac
						{\|u(x)\|^2_{\R^m} \|v(x) -v(y)\|^2_{\R^m}|}
						{\|x-y\|_{\R^d}^{2s+d}}
				\ud x \ud y 
				\\
			\leq{} 
				&\bigg(
					\int_{\Omega} \int_{\Omega}
						\frac
							{
								\|
									u(x)
								\|_{\R^m}^{q_1}
							}
							{\|x-y\|_{\R^d}^{d-\delta q_1}}
					\ud x \ud y 
				\bigg)^{\nicefrac{2}{q_1}}
				\bigg(
					\int_{\Omega} \int_{\Omega}
						\frac
							{
								\|
									v(x) -v(y)
								\|_{\R^m}^{q_2}
							}
							{\|x-y\|_{\R^d}^{q_2/2 (2s+d-(d-\delta q_1) 2/q_1)}}
					\ud x \ud y 
				\bigg)^{\nicefrac {2}{q_2}}\\
			\leq{} 
				&\bigg(
					\int_{\Omega} 
						\int_{\Omega}
							\frac
								{
									1
								}
								{\|x-y\|_{\R^d}^{d-\delta q_1}}
						\ud y \,
						\|
							u(x)
						\|_{\R^m}^{q_1}
					\ud x 
				\bigg)^{\nicefrac{2}{q_1}}
			\\ & \qquad \cdot 
				\bigg(
					\int_{\Omega} \int_{\Omega}
						\frac
							{
								\|
									v(x) -v(y)
								\|_{\R^m}^{q_2}
							}
							{\|x-y\|_{\R^d}^{q_2s+d q_2/2-d q_2(1/2-1/q_2)+\delta q_2}}
					\ud x \ud y 
				\bigg)^{\nicefrac {2}{q_2}}\\
			\leq{} 
				&C_1
				\Big(
					\int_{\Omega} 
						\|
							u(x)
						\|_{\R^m}^{q_1}
					\ud x 
				\Big)^{\nicefrac{2}{q_1}}
				\bigg(
					\int_{\Omega} \int_{\Omega}
						\frac
							{
								\|
									v(x) -v(y)
								\|_{\R^m}^{q_2}
							}
							{\|x-y\|_{\R^d}^{q_2 (s+\delta) +d}}
					\ud x \ud y 
				\bigg)^{\nicefrac {2}{q_2}}\\
			\leq{} 
				&C_2
					\|u\|^2_{L^{q_1}(\Omega,\R^m)}
					\|v\|^2_{B^{s+\delta}_{q_2,q_2}(\Omega,\R^m)}
			\leq
				C_3 \|u\|^2_{L^{q_1}(\Omega,\R^m)}
					\|v\|^2_{B^{s+\delta+d/q_1}_{2,2}(\Omega,\R^m)}.
		\end{split}
	\end{equation}
	In addition, it follows from Hölder's inequality,
	the fact that $\Omega$ is bounded, 
	and from the
	Sobolev inequality 
	(see, e.g., item (e) in the Theorem on page 327 in
	Triebel \cite{Triebel1978}) 
	that 
	for all $\delta \in (0, 1-s)$
	there exists a $ C \in (0,\infty)$ such that
	for all $u, v \in L^0(\Omega,\R^m)$
	it holds that
	\begin{equation}
	\label{eq: holder estimate with p,infty frac s}
		\begin{split}
				&\int_{\Omega} \int_{\Omega}
					\frac
						{\|u(x)\|^2_{\R^m} \| v(x) -v(y)\|^2_{\R^m}}
						{\|x-y\|_{\R^d}^{2s+d}}
				\ud x \ud y 
			\leq
				\int_{\Omega} \int_{\Omega}
						\frac
							{
								\|
									u(x)
								\|_{\R^m}^{2}
							}
							{\|x-y\|_{\R^d}^{d-2\delta}}
					\ud x \ud y 
				\cdot \sup_{x,y \in \Omega}
					\frac
							{
								\|
									v(x) -v(y)
								\|_{\R^m}^{2}
							}
							{\|x-y\|_{\R^d}^{2s + 2\delta}}\\
			\leq{} 
				&\int_{\Omega} 
						\int_{\Omega}
							\frac
								{
									1
								}
								{\|x-y\|_{\R^d}^{d-2\delta}}
						\ud y \,
						\|
							u(x)
						\|_{\R^m}^{2}
					\ud x 
				\cdot \|v\|^2_{C^{s+\delta}(\Omega, \| \cdot \|_{\R^m})} 
			\leq
				C \|u\|^2_{L^{2}(\Omega,\R^m)}
					\|v\|^2_{B^{s+\delta+d/2}_{2,2}(\Omega,\R^m)}.
		\end{split}
	\end{equation}
	Similar we obtain that
	from Proposition 2.28 in Mitrea \cite{Mitrea2013} that
	for all $\delta \in (0, 1-s)$
	there exists a $ C \in (0,\infty)$ such that
	for all $u,v \in L^0(\Omega,\R^m)$ 
	it holds that
	\begin{equation}
	\label{eq: holder estimate with infty,p frac s}
		\begin{split}
				&\int_{\Omega} \int_{\Omega}
					\frac
						{\|u(x)\|^2_{\R^m} \|v(x) -v(y)\|_{\R^m}^2}
						{\|x-y\|_{\R^d}^{2s+d}}
				\ud x \ud y \\
			\leq{}
				&\|u\|^2_{L^{\infty}(\Omega,\R^m)}
				\int_{\Omega} \int_{\Omega}
					\frac
						{\| v(x) -v(y)\|_{\R^m}^2}
						{\|x-y\|_{\R^d}^{2s+d}}
				\ud x \ud y 
			\leq{}
				C
				\|u\|^2_{L^{\infty}(\Omega,\R^m)}
				\|v\|^2_{B^{s}_{2,2}(\Omega,\R^m)}.
		\end{split}
	\end{equation}
Combining \eqref{eq: holder estimate with p,q frac s}
(with $u \leftarrow u \cdot \|u\|^{p-1}_{\R^m} $),
\eqref{eq: holder estimate with p,infty frac s}
(with $u \leftarrow u \cdot \|u\|^{p-1}_{\R^m} $),
\eqref{eq: holder estimate with infty,p frac s}
(with $u \leftarrow u \cdot \|u\|^{p-1}_{\R^m} $),
the fact that $pq \geq 2 \cdot \1_{ (0,\infty)}(p)$,
and the Sobolev inequality 
(see, e.g., the Theorem on page 31 and Theorem 1 on page 32 
in Runst \& Sickel \cite{RunstSickel1996})
verifies that
there exist $C_1, C_2, C_3, C_4 \in (0,\infty)$ such that for all
$u, v \in L^0(\Omega,\R^m)$
it holds that
\begin{align}
\nonumber
				&\int_{\Omega} \int_{\Omega}
					\frac
						{\|u(x)\|_{\R^m}^{2p} \|v(x) -v(y)\|_{\R^m}^2}
						{\|x-y\|_{\R^d}^{2s+d}}
				\ud x \ud y 
				\\ \nonumber
			\leq{} 
				&C_1
					\| u \cdot \|u\|_{\R^m}^{p-1}\|^2_{L^{q}(\Omega,\R^m)}
					\|v\|^2_{B^{s+\eps \1_{[2,\infty)}(q)+d/q}_{2,2}(\Omega,\R^m)} \\
	\label{eq: combined holder estimate}
			\leq{}
				&C_2
				\| \|u\|^p_{\R^m}\|^2_{L^{q}(\Omega,\R)}
				\|v\|^2_{B^{s+\eps \1_{[2,\infty)}(q)+d/q}_{2,2}(\Omega,\R^m)} \\ \nonumber
			\leq{}
				&C_3
				\|u\|^{2p}_{L^{p q}(\Omega,\R^m)}
				\|v\|^2_{B^{s+\eps \1_{[2,\infty)}(q)+d/q}_{2,2}(\Omega,\R^m)} \\ \nonumber
			\leq{} 
				&C_4 \|u\|^{2p}_{B^{d/2-d/(p q) + \eps \1_{\{\infty \}}(q)}_{2,2}(\Omega,\R^m)}
					\|v\|^2_{B^{s+\eps \1_{[2,\infty)}(q)+d/q}_{2,2}(\Omega,\R^m)}.
	\end{align}
	Therefore we get from the assumption that for all
	$x \in \Omega$ and all $y\in \R^m$ it holds that
	$\|D_{\R^d \times \R^m} b(x,y)\|_{(\R^d \times \R^m)'} \leq K (\|y\|_{\R^m} +1)^p$,
	from \eqref{eq: combined holder estimate},
	from Proposition 2.28 in Mitrea \cite{Mitrea2013},
	from the fact that $q \geq 2$,
	from the fact that $\Omega$ is bounded,
	from the fact that 
	for all $r\in \R$ and all $t\in (-\infty, r]$
	there exists a $C \in (0,\infty)$ such that for all 
	$u \in B^t_{2,2}(\Omega,\R^m)$ it holds that
	$\|u\|_{B^t_{2,2}(\Omega,\R^m)} \leq C \|u\|_{B^r_{2,2}(\Omega,\R^m)}$,
	the fact that $pq \geq 2 \cdot \1_{ (0,\infty)}(p)$,
	and from the Sobolev inequality 
		(see, e.g., the Theorem on page 31 in Runst \& Sickel \cite{RunstSickel1996})
	that
	there exist $C_1,C_2, C_3, C_4 \in (0,\infty)$
	such that for all 
	$u \in L^0(\Omega,\R^m)$
	it holds that
\begin{align}
	\label{eq: bound assertion for p in (0,infty)}
			&\int_{\Omega} \int_{\Omega}
				\frac
					{
						|
							b( x ,u(x)) -b( y ,u(y)) 
						|^2
					}
					{\|x-y\|_{\R^d}^{2s+d}}
			\ud x \ud y \\
		\leq{} \nonumber
			&\int_{\Omega} \int_{\Omega} \int_0^1
				\frac
					{1}
					{\|x-y\|_{\R^d}^{2s+d}}
						\Big |
							\Big \langle
								(D_{\R^d \times \R^m} b) \big(tx + (1-t)y,t u(x) + (1-t) u(y) \big),
				\\ & \nonumber  \qquad
								(x-y, u(x)-u(y) )
							\Big \rangle_{(\R^d \times \R^m)',\R^d \times \R^m}
						\Big |^2
			\ud t \ud x \ud y \\
		\leq{} \nonumber
			&\int_{\Omega} \int_{\Omega} \int_0^1
				\frac
					{
							\| 
								(D_{\R^d \times \R^m} b) \big(tx + (1-t)y,t u(x) + (1-t) u(y) \big)
							\|^2_{(\R^d \times \R^m)'}
					}
					{\|x-y\|_{\R^d}^{2s+d}} 
		\\ & \nonumber \qquad
				\cdot (\|x-y\|^2_{\R^d}+\|u(x)-u(y)\|_{\R^m}^2)
			\ud t \ud x \ud y \\
		\leq{} \nonumber
			&K^2\int_{\Omega} \int_{\Omega} \int_0^1
				\frac
					{
						( 
							\|t u(x) + (1-t) u(y)\|_{\R^m}+1
						)^{2p}
					}
					{\|x-y\|_{\R^d}^{2(s-1)+d}} 
		\\& \nonumber \qquad
			+
				\frac
					{
						\|u(x)-u(y)\|_{\R^m}^2 \cdot
						(\|t u(x) + (1-t) u(y) \|_{\R^m}+1)^{2p}
					}
					{\|x-y\|_{\R^d}^{2s+d}} 
			\ud t \ud x \ud y \\
		\leq{}  \nonumber
			&C_1 \int_{\Omega} \int_{\Omega}
				\frac
					{
						\| 
							u(x)
						\|_{\R^m}^{2p}
						+1
					}
					{\|x-y\|_{\R^d}^{2(s-1)+d}} 
			+
				\frac
					{
						\|u(x)-u(y)\|_{\R^m}^2
						(\|u(x)\|_{\R^m}^{2p} +1)
					}
					{\|x-y\|_{\R^d}^{2s+d}} 
			\ud x \ud y \\
		\leq{} \nonumber
			&C_2 (
				1+ \|u\|^{2p}_{L^{2p}(\Omega,\R^m)} 
				+\|u\|^{2p}_{B^{d/2-d/(p q) + \eps \1_{\{ \infty \}}(q)}_{2,2}(\Omega,\R^m)}
					\|u\|^2_{B^{s+\eps \1_{ [2, \infty )}(q)+d/q}_{2,2}(\Omega,\R^m)} 
			\\ & \nonumber \qquad
				+\|u\|^2_{B^{s}_{2,2}(\Omega,\R^m)}
			) \\
		\leq{} \nonumber
			&C_3 (
				1+ \|u\|^{2p}_{L^{pq}(\Omega,\R^m)} 
				+(
					\|u\|^{2p}_{B^{d/2-d/(p q) 
					+ \eps \1_{\{ \infty \}}(q)}_{2,2}(\Omega,\R^m)} +1
				)
					\|u\|^2_{B^{s+\eps \1_{ [2, \infty )}(q)+d/q}_{2,2}(\Omega,\R^m)} 
			) \\
		\leq{} \nonumber
			&C_4 
				(\|u\|^{2p}_{B^{d/2-d/(p q) + \eps \1_{\{ \infty \}}(q)}_{2,2}(\Omega,\R^m)} +1)
				(\|u\|^2_{B^{s+\eps \1_{ [2, \infty )}(q)+d/q}_{2,2}(\Omega,\R^m)} +1) \\
	\nonumber
		={}
			&C_4 \, 
				K_1(p,u) \cdot 
				(\|u\|^2_{B^{s+\eps \1_{[2,\infty)}(q)+d/q}_{2,2}(\Omega,\R^m)} +1).
\end{align}
{\it 3. case $p = \infty$.}
Then, it follows from
the fact that $\Omega$ is bounded,
from Proposition 2.28 in Mitrea \cite{Mitrea2013},
and from the fact that 
for all $r\in \R$ and all $t\in (-\infty, r]$
there exists a $C \in (0,\infty)$ such that for all 
$u \in B^t_{2,2}(\Omega,\R^m)$ it holds that
$\|u\|_{B^t_{2,2}(\Omega,\R^m)} \leq C \|u\|_{B^r_{2,2}(\Omega,\R^m)}$
that
there exist $C_1, C_2 \in (0, \infty)$ such that for all 
$
	u \in L^0(\Omega,\R^m)
$
it holds that
\begin{align}
\nonumber
			&\int_{\Omega} \int_{\Omega}
				\frac
					{
						|
							b( x ,u(x)) -b( y ,u(y)) 
						|^2
					}
					{\|x-y\|_{\R^d}^{2s+d}}
			\ud x \ud y \\ \nonumber
		\leq{} 
			&\int_{\Omega} \int_{\Omega} \int_0^1
				\frac
					{1}
					{\|x-y\|_{\R^d}^{2s+d}}
						\Big |
							\Big \langle 
								(D_{\R^d \times \R^m} b) \big(tx + (1-t)y,t u(x) + (1-t) u(y) \big),
				\\ \nonumber & \qquad
								(x-y, u(x)-u(y) )
							\Big \rangle_{(\R^d \times \R^m)',\R^d \times \R^m}
						\Big |^2
			\ud t \ud x \ud y \\ \nonumber
		\leq{} 
			&\int_{\Omega} \int_{\Omega} \int_0^1
				\frac
					{
							\| 
								(D_{\R^d \times \R^m} b) \big(tx + (1-t)y,t u(x) + (1-t) u(y) \big)
							\|^2_{(\R^d \times \R^m)'}
					}
					{\|x-y\|_{\R^d}^{2s+d}}  
		\\ & \qquad 
	\label{eq: bound assertion for p=infty}
				\cdot (\|x-y\|^2_{\R^d}+\|u(x)-u(y)\|_{\R^m}^2)
			\ud t \ud x \ud y \\ \nonumber
		\leq{} 
			&\|
				D_{\R^d \times \R^m} b
			\|^2_{
				\C(
					\Omega \times 
					B(\|u \|_{L^\infty(\Omega,\R^m)}),
					\| \cdot \|_{(\R^d \times \R^m)'} 
				)
			}
		\\  \nonumber & \qquad 
			\int_{\Omega} \int_{\Omega} \int_0^1
				\frac
					{
						1
					}
					{\|x-y\|_{\R^d}^{2(s-1)+d}} 
			+
				\frac
					{
						\|u(x)-u(y)\|_{\R^m}^2
					}
					{\|x-y\|_{\R^d}^{2s+d}} 
			\ud t \ud x \ud y \\ \nonumber
		\leq{} 
			&C_1 \|
				D_{\R^d \times \R^m} b
			\|^2_{
				\C(
					\Omega \times 
					B(\|u \|_{L^\infty(\Omega,\R^m)}),
					\| \cdot \|_{(\R^d \times \R^m)'}  
				)
			}
			(
				1+\|u\|^2_{B^{s}_{2,2}(\Omega,\R^m)} 
			) \\ \nonumber
		\leq{} 
			&C_2 \|
				D_{\R^d \times \R^m} b
			\|^2_{
				\C(
					\Omega \times 
					B(\|u \|_{L^\infty(\Omega,\R^m)}),
					\| \cdot \|_{(\R^d \times \R^m)'}  
				)
			}
			(
				1+\|u\|^2_{B^{s+\eps \1_{[2,\infty)}(q)+d/q}_{2,2}(\Omega,\R^m)} 
			) \\ \nonumber
		={} 
			&C_2 \, K_1(\infty,u) 
				\cdot (\|u\|^2_{B^{s+\eps \1_{[2,\infty)}(q)+d/q}_{2,2}(\Omega,\R^m)} +1)
\end{align}
and this together with \eqref{eq: bound assertion for p=0}
and with \eqref{eq: bound assertion for p in (0,infty)}
completes the poof of Lemma \ref{lem: b boundedness inequality}.
\end{proof}
The next lemma establish a bound on the $L^2$-norm.
\begin{lemma}[Upper bound for the $L^2$-norm of Nemytskij operators]
\label{lem: b L2 boundedness inequality}
	Let $d,m \in \N$, 
	let $\Omega \subseteq \R^d$ be a bounded Lipschitz domain,
	denote by $B(r) \subseteq \R^m$, $r \in [0,\infty]$,
	the ball satisfying that for all $r \in [0,\infty]$ it holds
	that
	$
			B(r)
		= 
			\{ 
				y \in \R^m \colon 
					\| y\|_{\R^m} \leq r
			\}
	$,
	let $K \in (0,\infty)$,
	$p \in [0,\infty]$,
	$b \in \C_{\R^d \times \R^m}^1(\Omega \times \R^m, \R)$,
	$K_2\colon [0,\infty] \times L^0(\Omega,\R^m) \to [0, \infty]$,
	satisfy that for all $x \in \Omega$, $y \in \R^m$, $r \in [0, \infty]$
	and all $u \in L^0(\Omega,\R^m)$
	it holds that
	$
			\|D_{\R^d \times \R^m} b(x,y)\|_{(\R^d \times \R^m)'} 
		\leq 
			K (\|y\|_{\R^m} +1)^{p}
	$
	and that
	\begin{equation}
			K_2(r,u)
		=
			\begin{cases}
					\big (
						\|u \|^2_{L^2(\Omega,\R^m)} 
						+ \|u \|^{2r+2}_{B_{2,2}^{d/2 - d/(2r+2)}(\Omega,\R^m)}
					\big)
					& \textrm{ if } r \in [0,\infty)\\
					\|
						D_{\R^d \times \R^m} b
					\|^2_{
						\C(
							\Omega
							\times B(\|u\|_{L^\infty(\Omega,\R^m)}),
							\| \cdot \|_{(\R^d \times \R^m)'} 
						)
					}
				\|u\|_{L^2(\Omega,\R^m)}^2
					& \textrm{ if } r =\infty. \\
				\end{cases}
	\end{equation}
	Then 
	there exists a
	$C \in (0,\infty)$ such that 
	for all $u \in L^0(\Omega,\R^m)$,
	it holds that
	\begin{equation}
		\begin{split}
				\int_{\Omega} 
					|
						b( x ,u(x)) 
					|^2
			\ud x 
		\leq
			C \cdot K_2(p,u)
			+2\int_\Omega
				|b(x,0)|^2
			\ud x.
		\end{split}
	\end{equation}
 \end{lemma}
\begin{proof}
	We will divide the proof into 2 cases \\
	{\it 1. case $p \in [0,\infty)$.}
	Then it follows from 
	$
			\|D_{\R^d \times \R^m} b(x,y)\|_{(\R^d \times \R^m)'} 
		\leq 
			K (\|y\|_{\R^m} +1)^{p}
	$
	that for all $u \in L^0(\Omega,\R^m)$ it holds that
	\begin{align}
  \nonumber
				&\int_\Omega 
					|b(x,u(x))|^2
				\ud x \\
			={} 
				&\int_\Omega
					\Big |
						\int_0^1 
							\big \langle
								(D_{\R^d \times \R^m} b) (x, t u(x)),
								(0, u(x))
							\big \rangle_{(\R^d \times \R^m)', (\R^d \times \R^m)}
						\ud t
						+b(x,0)
					\Big |^2
				\ud x \\ \nonumber
			\leq{} 
				&2\int_\Omega
					\Big |
						\int_0^1 
							\big \|
								(D_{\R^d \times \R^m} b) (x, t u(x))
							\big \|_{(\R^d \times \R^m)'}
							\|u(x)\|_{\R^m}
						\ud t
					\Big |^2
					+|b(x,0)|^2
				\ud x \\ \nonumber
			\leq{} 
				&2K^2 \int_\Omega
					\Big |
						\int_0^1
								(1+t |u(x)|)^p
						\ud t
							\cdot
							\|u(x) \|_{\R^m}
					\Big |^2
				\ud x
				+2\int_\Omega
					|b(x,0)|^2
				\ud x. 
	\end{align}
	Thus we derive from the Sobolev inequality 
	(see, e.g., the Theorem on page 31 in Runst \& Sickel \cite{RunstSickel1996})
	that there exist $C_1, C_2 \in (0,\infty)$
	such that for all $u \in L^0(\Omega,\R^m)$
	it holds that
	\begin{align}
				\int_\Omega 
					|b(x,u(x))|^2
				\ud x
			\leq{} 
				&2^{p+1} K^2
				\int_\Omega
					\Big |
						\int_0^1
							1+t^p \|u(x)\|_{R^m}^p
						\ud t
							\cdot
							\| u(x) \|_{\R^m}
					\Big |^2
				\ud x
				+2\int_\Omega
					|b(x,0)|^2
				\ud x \\ \nonumber
			\leq{} 
				&2^{p+1} K^2
				\int_\Omega
					\Big |
						\|u(x)\|_{\R^m}+\tfrac{\|u(x)\|_{\R^m}^{p+1}}{p+1}	
					\Big |^2
				\ud x
				+2\int_\Omega
					|b(x,0)|^2
				\ud x \\ \nonumber
			\leq{} 
				&2^{p+2} K^2
				\int_\Omega
					\|u(x)\|_{\R^m}^2+\tfrac{\|u(x)\|_{\R^m}^{2p+2}}{(p+1)^2}	
				\ud x
				+2\int_\Omega
					|b(x,0)|^2 
				\ud x \\ \nonumber
			\leq{} 
				&C_1 \big (
					\|u \|^2_{L^2(\Omega,\R^m)} +\|u \|^{2p+2}_{L^{2p+2}(\Omega,\R^m)}
				\big)
				+2\int_\Omega
					|b(x,0)|^2
				\ud x \\ \nonumber
			\leq{} 
				&C_2 \Big (
					\|u \|^2_{L^2(\Omega,\R^m)} 
					+ \|u \|^{2p+2}_{B_{2,2}^{d/2 - d/(2p+2)}(\Omega,\R^m)}
				\Big)
				+2\int_\Omega
					|b(x,0)|^2
				\ud x \\ \nonumber
			=
				&C_2 \cdot K_2(p,u)
				+2\int_\Omega
					|b(x,0)|^2
				\ud x.
	\end{align}
	{\it 2. case $p = \infty$.}
	Then it holds for all $u \in L^0(\Omega,\R^m)$ that
	\begin{align}
	\nonumber
				&\int_\Omega 
					|b(x,u(x))|^2
				\ud x \\
			={} \nonumber 
				&\int_\Omega
					\Big |
						\int_0^1 
							\big \langle
								(D_{\R^d \times \R^m} b) (x, t u(x)),
								(0, u(x))
							\big \rangle_{(\R^d \times \R^m)', (\R^d \times \R^m)}
						\ud t
						+b(x,0)
					\Big |^2
				\ud x \\
			\leq{} \nonumber 
				&2\int_\Omega
					\Big |
						\int_0^1 
							\big \|
								(D_{\R^d \times \R^m} b) (x, t u(x))
							\big \|_{(\R^d \times \R^m)'}
							\|u(x)\|_{\R^m}
						\ud t
					\Big |^2
					+|b(x,0)|^2
				\ud x \\
			\leq{} 
				&2\int_\Omega
					\Big |
						\int_0^1 
							\|
								D_{\R^d \times \R^m} b
							\|_{
								\C(
									\Omega
									\times 
										B(\|u\|_{L^\infty(\Omega,\R^m)}),
									\| \cdot \|_{(\R^d \times \R^m)'} 
								)
							}
							\|u(x)\|_{\R^m}
						\ud t
					\Big |^2
					+|b(x,0)|^2
				\ud x \\
			={} \nonumber 
				&2\int_\Omega
					\|
						D_{\R^d \times \R^m} b
					\|^2_{
						\C(
							\Omega
							\times B(\|u\|_{L^\infty(\Omega,\R^m)}),
							\| \cdot \|_{(\R^d \times \R^m)'} 
						)
					}
					\|u(x)\|_{\R^m}^2
				\ud x
				+2\int_\Omega
					|b(x,0)|^2
				\ud x\\
			={} \nonumber 
				&2\|
					D_{\R^d \times \R^m} b
				\|^2_{
					\C(
						\Omega
						\times B(\|u\|_{L^\infty(\Omega,\R^m)}),
						\| \cdot \|_{(\R^d \times \R^m)'} 
					)
				}
				\|u\|_{L^2(\Omega,\R^m)}^2
				+2\int_\Omega
					|b(x,0)|^2
				\ud x \\
			={} \nonumber 
				&2 \cdot K_2(\infty,u)
				+2\int_\Omega
					|b(x,0)|^2
				\ud x.
	\end{align}
	Combining the two cases finishes the proof of 
	Lemma \ref{lem: b L2 boundedness inequality}.
\end{proof}
As a direct consequence of Lemma \ref{lem: b boundedness inequality}
and 
Lemma \ref{lem: b L2 boundedness inequality} we get
the next corollary.
\begin{corollary}[Upper bound for the Besov-norm of Nemytskij operators]
\label{cor: Bs b bound}
	Let $d,m \in \N$, 
	let $\Omega \subseteq \R^d$ be a bounded Lipschitz domain,
	denote by $B(r) \subseteq \R^m$, $r \in [0,\infty]$,
	the ball satisfying that for all $r \in [0,\infty]$ it holds
	that
	$
			B(r)
		= 
			\{ 
				y \in \R^m \colon 
					\| y\|_{\R^m} \leq r
			\}
	$,
	let $K \in (0,\infty)$,
	$s \in (0,1)$,
	$\eps \in (0, 1-s)$,
	$p \in [0,\infty]$,
	$q \in [2,\infty]$,
	with $pq \geq 2 \cdot \1_{ (0,\infty)}(p)$,
	let
	$b \in \C_{\R^d \times \R^m}^1( \Omega \times \R^m, \R)$,
	$K_1, K_2 \colon [0,\infty] \times L^0(\Omega,\R^m) \to [0, \infty]$,
	satisfy that for all $x \in \Omega$, $y \in \R^m$, $r \in [0,\infty]$,
	and all $u \in L^0(\Omega,\R^m)$
	it holds that
	$
			\|D_{\R^d \times \R^m} b(x,y)\|_{(\R^d \times \R^m)'} 
		\leq 
			K (\|y\|_{\R^m} +1)^{p},
	$
	that
	\begin{equation}
			K_1(r,u)
		=
		\begin{cases}
			1 & \textrm{ if } r= 0 \\
			\|u\|^{2r}_{B^{d/2-d/(r q) + \eps \1_{\{\infty\}}(q)}_{2,2}(\Omega,\R^m)} +1	
					& \textrm{ if } r \in (0,\infty) \\
					\|
						D_{\R^d \times \R^m} b
					\|_{
						\C(
							\Omega
							\times B(\|u\|_{L^\infty(\Omega,\R^m)}),
							\| \cdot \|_{(\R^d \times \R^m)'}
						)
					}
					& \textrm{ if } r =\infty,
		\end{cases}
	\end{equation}
	and that
	\begin{equation}
			K_2(r,u)
		=
			\begin{cases}
					\big (
						\|u \|^2_{L^2(\Omega,\R^m)} 
						+ \|u \|^{2r+2}_{B_{2,2}^{d/2 - d/(2r+2)}(\Omega,\R^m)}
					\big)
					& \textrm{ if } r \in [0,\infty)\\
					\|
						D_{\R^d \times \R^m} b
					\|^2_{
						\C(
							\Omega \times
							B(\|u\|_{L^\infty(\Omega,\R^m)}),
							\| \cdot \|_{(\R^d \times \R^m)'} 
						)
					}
				\|u\|_{L^2(\Omega,\R^m)}^2
					& \textrm{ if } r =\infty.
				\end{cases}
	\end{equation}
	Then there exists a
	$C \in (0,\infty)$ such that 
	for all $u \in L^0(\Omega,\R^m)$,
	it holds that
	\begin{equation}
				\|
					b( \cdot ,u( \cdot ))
				\|_{B^s_{2,2}(\Omega)}^2
		\leq
			C \, K_1(p,u) \cdot 
				(\|u\|^2_{B^{s+\eps \1_{[2,\infty)}(q)+d/q}_{2,2}(\Omega,\R^m)} +1)
			+C \, K_2(p,u) 
			+C\int_\Omega
				|b(x,0)|^2
			\ud x.
	\end{equation}
\end{corollary}
\begin{proof}
	This Corollary follows directly from Lemma \ref{lem: b boundedness inequality}, 
	Lemma \ref{lem: b L2 boundedness inequality},
	and from Proposition 2.28 in Mitrea \cite{Mitrea2013}.
\end{proof}
The next lemma shows the Lipschitz continuity with respect to 
Slobodeckij semi-norms
of Nemytskij operators
on bounded sets.
\begin{lemma}[Lipschitz continuity wrt. Slobodeckij semi-norms
 of Nemytskij operators on bounded sets]
\label{lem: b continuity inequality}
	Let $d, m \in \N$, 
	let $\Omega \subseteq \R^d$ be a bounded Lipschitz domain,
	denote by $B(r) \subseteq \R^m$, $r \in [0,\infty]$,
	the ball satisfying that for all $r \in [0,\infty]$ it holds
	that
	$
			B(r)
		= 
			\{ 
				y \in \R^m \colon 
					\| y\|_{\R^m} \leq r
			\}
	$,
	let 
	$s \in (0,1)$,
	$\eps \in (0,1-s)$,
	$p_1,p_2 \in [0,\infty]$,
	$q_1,q_2,q_3,q_4,q_5 \in [2,\infty]$,
	satisfy that
	$
			\nicefrac {1}{q_1} + \nicefrac{1}{q_2} + \nicefrac{1}{q_3}
		=
			\nicefrac {1}{q_4} + \nicefrac{1}{q_5}
		=
			\nicefrac {1}{2}
	$, that
	$p_1 \, q_3 \geq 2 \cdot \1_{(0,\infty)}(p_1)$,
	and that $p_2 \, q_5 \geq 2 \cdot \1_{(0,\infty)}(p_2)$,
	let
	$K \in (0,\infty)$,
	$b \in \C_{\R^d \times \R^m}^{1,2}(\Omega \times \R^m, \R)$,
	$b''_{2} \in \C(\Omega \times \R^m, L(\R^m,(\R^d \times \R^m)'))$,
	$K_1 \colon [0,\infty] \times L^0(\Omega,\R^m) \times L^0(\Omega,\R^m) \to [0, \infty]$,
	$K_2 \colon [0,\infty] \times L^0(\Omega,\R^m) \to [0, \infty]$,
	satisfy that for all $p \in [0,\infty]$,
	$x \in \Omega$, $y \in \R^m$, and all $u_1, u_2 \in L^0(\Omega,\R^m)$ 
	it holds that
	\begin{equation}
		\begin{split}
				b''_{2}(x,y)
			=
				\Big(
					D_{\R^m}\big( 
						\R^m \ni z \to D_{\R^d \times \R^m} b(x,z) \in (\R^d \times \R^m)'
					\big)
				\Big)(y),
		\end{split}
	\end{equation}
	that
	\begin{equation}
			K_1(p,u_1, u_2)
		=
			\begin{cases}
				1 & \textrm{if } p =0 \\
				1
				+\sum^2_{i=1}
					\|u_i\|^{2p}_{B^{d/2-d/(p q_3)+\eps \1_{\{\infty\}}(q_3)}_{2,2}(\Omega,\R^m)}
				& \textrm{if } p \in (0,\infty) \\
				\|
					b''_2
				\|^2_{
					\C(
						\Omega 
						\times 
						B(\|u_1\|_{L^\infty(\Omega,\R^m)} \vee \|u_2\|_{L^\infty(\Omega,\R^m)}),
						\| \cdot \|_{L(\R^m,(\R^d \times \R^m)')}
					)
				} & \textrm{if } p =\infty
			\end{cases},
	\end{equation}
	that
	\begin{equation}
			K_2(p,u_1)
		=
			\begin{cases}
				1 & \textrm{ if } p =0 \\
				\|u_1\|^{2p}_{B^{d/2-d/(p q_5)+\eps \1_{\{\infty\}}(q_5)}_{2,2}(\Omega,\R^m)}
				+1 & \textrm{ if } p \in (0,\infty) \\
				\|
					D_{\R^d \times \R^m} b
				\|^2_{
					\C(
						\Omega
						\times 
						B(\|u_1\|_{L^\infty(\Omega,\R^m)}),
						\| \cdot \|_{(\R^d \times \R^m)'}
					)
				} & \textrm{ if } p =\infty
			\end{cases},
	\end{equation}
	that
	\begin{equation}
	\label{eq: b estimate 1st derivative} 
			\|
				D_{\R^d \times \R^m} b(x,y)
			\|_{(\R^d \times \R^m)'} 
		\leq 
			K \cdot (1+\|y\|_{\R^m})^{p_2},
	\end{equation}
	and that
	\begin{equation}
	\label{eq: b estimate 2nd derivative} 
			\|
				b''_2(x,y)
			\|_{L(\R^m, (\R^d \times \R^m)')}
		\leq
			K \cdot (1+\|y\|_{\R^m})^{p_1}.
	\end{equation}
	Then there exists a
	$C \in (0,\infty)$ such that 
	for all $u, \tilde{u} \in L^0(\Omega,\R^m)$
	it holds that
	\begin{equation}
		\begin{split}
				&\int_{\Omega} \int_{\Omega}
					\frac
						{
							|
								b( x ,u(x))-b (x,\tilde{u}(x)) -b( y ,u(y)) + b(y, \tilde{u}(y)) 
							|^2
						}
						{\|x-y\|_{\R^d}^{2s+d}} 
				\ud x \ud y \\
			\leq 
				&C \, \|
					u-\tilde{u}
				\|^{2}_{B^{d/2-d/q_1+\eps \1_{\{\infty\}}(q_1)}_{2,2}(\Omega,\R^m)}
				\big(
					\|
						u
					\|^2_{B^{s+d/2-d/q_2+\eps \1_{(2,\infty]}(q_2)}_{2,2}(\Omega,\R^m)}
					+1
				\big) \cdot K_1(p_1,u,\tilde{u}) \\
			& \qquad
				+C \, 
					\|
						u - \tilde{u}
					\|^2_{B^{s+d/2-d/q_4+\eps \1_{(2,\infty]}(q_4)}_{2,2}(\Omega,\R^m)}
				\cdot K_2(p_2, \tilde{u}).
		\end{split}
	\end{equation}
 \end{lemma}
\begin{proof}
First denote by $ {\bf 1} \in \C(\Omega,\R^m)$ the function satisfying 
for all $x \in \Omega$ that
\begin{align}
	{\bf 1}(x) = (1,0, \cdots, 0).
\end{align}
Next note that Hölder's inequality and the fact that
$\Omega$ is bounded verify that
for all $\tilde{q}_1, \tilde{q}_2 \in (2, \infty)$ and all $t \in (s,\infty)$ with 
$ \nicefrac {1}{\tilde{q}_1} + \nicefrac{1}{\tilde{q}_2} = \nicefrac 12$
there exists a $C \in (0, \infty)$ such that
all $ v_1, v_2 \in L^0(\Omega,\R^m)$
it holds that
\begin{align}
\label{eq: simple estm xy}
			&\int_\Omega \int_\Omega 
				\frac
					{
						\|v_1(x) \|^2_{\R^m} \| v_2(y)\|^2_{\R^m} \|x-y\|_{\R^d}^{2t}
					}
					{\|x-y\|_{\R^d}^{2s+d}}
			\ud x \ud y \\
		={} \nonumber 
			&\int_\Omega \int_\Omega
				\frac{\|v_1(x)\|_{\R^m}^2}{\|x-y\|_{\R^d}^{2(d+2s-2t)/\tilde{q}_1}}
				\frac{\|v_2(y)\|_{\R^m}^2}{\|x-y\|_{\R^d}^{2(d+2s-2t)/\tilde{q}_2}}
			\ud x \ud y \\
		\leq{} \nonumber 
			&\Big(	\int_\Omega \int_\Omega
				\frac{\|v_1(x) \|_{\R^m}^{\tilde{q}_1}}{\|x-y\|_{\R^d}^{d+2s-2t}}
			\ud x \ud y \Big)^{\nicefrac {2}{\tilde{q}_1}}
			\Big (\int_\Omega \int_\Omega
				\frac{\|v_2(y)\|_{\R^m}^{\tilde{q}_2}}{\|x-y\|_{\R^d}^{d+2s-2t}}
			\ud x \ud y \Big)^{\nicefrac {2}{\tilde{q}_2}} \\
		={} \nonumber 
			&\Big(	\int_\Omega \int_\Omega
				\frac{1}{\|x-y\|_{\R^d}^{d+2(s-t)}}
			\ud y \, \|v_1(x)\|_{\R^m}^{\tilde{q}_1} \ud x \Big)^{\nicefrac {2}{\tilde{q}_1}}
			\Big (\int_\Omega \int_\Omega
				\frac{1}{\|x-y\|_{\R^d}^{d+2(s-t)}}
			\ud x \, \|v_2(y)\|_{\R^m}^{\tilde{q}_2} \ud y \Big)^{\nicefrac {2}{\tilde{q}_2}} \\
	\nonumber
		\leq{} 
			&C \|v_1\|^2_{L^{\tilde{q}_1}(\Omega,\R^m)} \|v_2\|^2_{L^{\tilde{q}_2}(\Omega,\R^m)}.
\end{align}
In addition, it holds that for all $t \in (s,\infty)$
there exists a $C \in (0, \infty)$ such that
all $ v_1, v_2 \in L^0(\Omega,\R^m)$
it holds that
\begin{equation}
\label{eq: simple estm. infty}
	\begin{split}
			&\int_\Omega \int_\Omega 
				\frac
					{\|v_1(x)\|^2_{\R^m} \|v_2(y)\|_{\R^m}^2 \|x-y\|_{\R^d}^{2t}}
					{\|x-y\|_{\R^d}^{2s+d}} \\
			\ud x \ud y 
		\leq{}
			&\|v_1\|^2_{L^{\infty}(\Omega,\R^m)}
			\int_\Omega \int_\Omega
				\frac{1}{\|x-y\|_{\R^d}^{d+2(s-t)}}
			\ud x \, \|v_2(y)\|_{\R^m}^{2} \ud y 
		\leq{} 
			C \|v_1\|^2_{L^{\infty}(\Omega,\R^m)} \|v_2\|^2_{L^{2}(\Omega,\R^m)}.
	\end{split}
\end{equation}
Analogously it follows that
for all $\tilde{q}_1, \tilde{q}_2 \in [2, \infty]$ and all $t \in (s,\infty)$ with 
$ \nicefrac {1}{\tilde{q}_1} + \nicefrac{1}{\tilde{q}_2} = \nicefrac 12$
there exists a $C \in (0, \infty)$ such that
all $ v_1, v_2 \in L^0(\Omega,\R^m)$
it holds that
\begin{align}
\label{eq: simple estm xx}
			\int_\Omega \int_\Omega 
				\frac
					{\|v_1(x)\|^2_{\R^m} \| v_2(x)\|_{\R^m}^2 \|x-y\|_{\R^d}^{2t}}
					{\|x-y\|_{\R^d}^{2s+d}}
			\ud x \ud y 
		\leq
			C \|v_1\|^2_{L^{\tilde{q}_1}(\Omega,\R^m)} \|v_2\|^2_{L^{\tilde{q}_2}(\Omega,\R^m)}.
\end{align}
Moreover, Hölder's inequality,
the fact that $\Omega$ is bounded,
Proposition 2.28 in Mitrea \cite{Mitrea2013},
and the Sobolev inequality
(see, e.g., the Theorem on page 31 in Runst \& Sickel \cite{RunstSickel1996})
imply that for all
$\tilde{q}_1, \tilde{q}_2, \tilde{q}_3 \in [2,\infty)$
with 
$ 
		\nicefrac {1}{\tilde{q}_1} + \nicefrac{1}{\tilde{q}_2} + \nicefrac {1}{\tilde{q}_3} 
	= 
		\nicefrac 12
$
there exist $C_1, C_2, C_3 \in (0,\infty)$ such that for 
all $ v_1, v_2, v_3 \in L^0(\Omega,\R^m)$
it holds that
\begin{align}
\label{eq: 3 func xy estimate q1q2 < infty}
			&\int_\Omega \int_\Omega
				\frac
					{
						\|v_1(x)\|^2_{\R^m} \|v_2(y)\|^2_{\R^m} \|v_3(x)- v_3(y)\|_{\R^m}^2
					}
					{\|x-y\|_{\R^d}^{2s+d}}
			\ud x \ud y \\ \nonumber
		={} 
			&\int_\Omega \int_\Omega
				\frac
					{\|v_1(x)\|_{\R^m}^2}
					{	
						\|x-y\|_{\R^d}^{(d-\eps \tilde{q}_1 /2)2/\tilde{q}_1}
					}
				\frac
					{\|v_2(y)\|_{\R^m}^2}
					{	
						\|x-y\|_{\R^d}^{(d-\eps \tilde{q}_2 /2)2/\tilde{q}_2}
					}
			\\ & \nonumber
				\cdot \frac
					{ \|v_3(x)- v_3(y)\|_{\R^m}^2}
					{	
						\|x-y\|_{\R^d}^{
							2s+d-(d-\eps \tilde{q}_1 /2)2/\tilde{q}_1
							-(d-\eps \tilde{q}_2 /2)2/\tilde{q}_2
						}
					}
			\ud x \ud y \\
		\leq{} \nonumber
			&\Big(	\int_\Omega \int_\Omega
				\frac{\|v_1(x)\|_{\R^m}^{\tilde{q}_1}}{\|x-y\|_{\R^d}^{d- \eps\tilde{q}_1/2}}
			\ud x \ud y \Big)^{\nicefrac {2}{\tilde{q}_1}}
			\Big (\int_\Omega \int_\Omega
				\frac{\|v_2(y)\|_{\R^m}^{\tilde{q}_2}}{\|x-y\|_{\R^d}^{d- \eps \tilde{q}_2/2}}
			\ud x \ud y \Big)^{\nicefrac {2}{\tilde{q}_2}} \\ \nonumber
			&\Big( \int_\Omega \int_\Omega
				\frac
					{\|v_3(x)- v_3(y)\|_{\R^m}^{\tilde{q}_3}}
					{	
						\|x-y\|_{\R^d}^{
							(
								2s+d
								-(d-\eps \tilde{q}_1 /2)2/\tilde{q}_1
								-(d-\eps \tilde{q}_2 /2)2/\tilde{q}_2
							)
							\tilde{q}_3/2
						}
					}
			\ud x \ud y \Big)^{\nicefrac {2}{\tilde{q}_3}} \\
		\leq{} \nonumber 
			&\Big(	\int_\Omega \int_\Omega
				\frac{1}{\|x-y\|_{\R^d}^{d- \eps \tilde{q}_1/2}}
				\ud y \,
				\|v_1(x)\|_{\R^m}^{\tilde{q}_1}
			\ud x \Big)^{\nicefrac {2}{\tilde{q}_1}}
			\Big (\int_\Omega \int_\Omega
				\frac{1}{\|x-y\|_{\R^d}^{d- \eps \tilde{q}_2/2}}
			\ud x \, 
			\|v_2(y)\|_{\R^m}^{\tilde{q}_2} \ud y \Big)^{\nicefrac {2}{\tilde{q}_2}} 
		\\ \nonumber
			&\Big( \int_\Omega \int_\Omega
				\frac
					{\|v_3(x)- v_3(y)\|_{\R^m}^{\tilde{q}_3}}
					{	
						\|x-y\|_{\R^d}^{
							(
								2s+d (1- 2/\tilde{q}_1 -2/\tilde{q}_2)
								+2\eps
							)
							\tilde{q}_3/2
						}
					}
			\ud x \ud y \Big)^{\nicefrac {2}{\tilde{q}_3}} \\
		\leq{} \nonumber
			&C_1 \|v_1\|^2_{L^{\tilde{q}_1}(\Omega,\R^m)} \|v_2\|^2_{L^{\tilde{q}_2}(\Omega,\R^m)} 
			\Big( \int_\Omega \int_\Omega
				\frac
					{\|v_3(x)- v_3(y)\|_{\R^m}^{\tilde{q}_3}}
					{	
						\|x-y\|_{\R^d}^{
							\tilde{q}_3 (s+ \eps)+d
						}
					}
			\ud x \ud y \Big)^{\nicefrac {2}{\tilde{q}_3}} \\
		\leq{} \nonumber
			&C_2 \|v_1\|^2_{L^{\tilde{q}_1}(\Omega,\R^m)} \|v_2\|^2_{L^{\tilde{q}_2}(\Omega,\R^m)} 
			\|v_3\|^2_{B^{s+\eps}_{\tilde{q}_3,\tilde{q}_3}(\Omega,\R^m)} \\
		\leq{} 
	\nonumber
			&C_3 \|v_1\|^2_{L^{\tilde{q}_1}(\Omega,\R^m)} \|v_2\|^2_{L^{\tilde{q}_2}(\Omega,\R^m)} 
			\|v_3\|^2_{B^{s+d/2-d/\tilde{q}_3+ \eps}_{2,2}(\Omega,\R^m)}.
\end{align}
Analogously to \eqref{eq: 3 func xy estimate q1q2 < infty} 
and to \eqref{eq: holder estimate with p,q frac s}
it follows that for all 
$\tilde{q}_2, \tilde{q}_3 \in (2,\infty)$
with $\nicefrac{1}{\tilde{q}_2} + \nicefrac {1}{\tilde{q}_3} = \nicefrac 12$
there exists a $C \in (0,\infty)$ such that for 
all $ v_1, v_2, v_3 \in L^0(\Omega,\R^m)$
it holds that
\begin{equation}
\label{eq: 3 func xy estimate q1 = infty q2 < infty}
	\begin{split}
			&\int_\Omega \int_\Omega
				\frac
					{
						\|v_1(x)\|^2_{\R^m} 
						\|v_2(y)\|^2_{\R^m} 
						\|v_3(x)- v_3(y)\|_{\R^m}^2
					}	
					{\|x-y\|_{\R^d}^{2s+d}}
			\ud x \ud y \\
		\leq{}
			&\|v_1\|^2_{L^\infty(\Omega,\R^m)}
			\int_\Omega \int_\Omega
				\frac{\|v_2(y)\|^2_{\R^m} \|v_3(x)- v_3(y)\|_{\R^m}^2}{\|x-y\|_{\R^d}^{2s+d}}
			\ud x \ud y \\ 
		\leq{}  
			&C \|v_1\|^2_{L^{\infty}(\Omega,\R^m)} \|v_2\|^2_{L^{\tilde{q}_2}(\Omega,\R^m)}  
				\|v_3\|^2_{B^{s+d/2 - d/\tilde{q}_3+\eps}_{2,2}(\Omega,\R^m)}.
	\end{split}
\end{equation}
In addition, we get
from Proposition 2.28 in Mitrea \cite{Mitrea2013}
that there exists a $C \in (0,\infty)$ such that for 
all $ v_1, v_2, v_3 \in L^0(\Omega,\R^m)$
it holds that
\begin{equation}
\label{eq: 3 func xy estimate q1q2 = infty}
	\begin{split}
			&\int_\Omega \int_\Omega
				\frac
					{\|v_1(x)\|^2_{\R^m} \| v_2(y)\|^2_{\R^m} \|v_3(x)- v_3(y)\|_{\R^m}^2}
					{\|x-y\|_{\R^d}^{2s+d}}
			\ud x \ud y \\
		\leq{}
			&\|v_1\|^2_{L^{\infty}(\Omega,\R^m)} \|v_2\|^2_{L^{\infty}(\Omega,\R^m)}
			\int_\Omega \int_\Omega
				\frac{\|v_3(x)- v_3(y)\|_{\R^m}^2}{\|x-y\|_{\R^d}^{2s+d}}
			\ud x \ud y \\
		\leq{} 
		&C \|v_1\|^2_{L^{\infty}(\Omega,\R^m)} \|v_2\|^2_{L^{\infty}(\Omega,\R^m)}  
				\|v_3\|^2_{B^{s}_{2,2}(\Omega,\R^m)}.
	\end{split}
\end{equation}
Finally observe, that \eqref{eq: simple estm xy}, 
\eqref{eq: simple estm. infty} (with $t \leftarrow s +\eps$)
and the Sobolev inequality (see, e.g., item (e)
in the Theorem on page 327 in Triebel \cite{Triebel1978})
show that for all $\tilde{q}_1, \tilde{q}_2 \in [2,\infty]$
with $\nicefrac {1}{\tilde{q}_1} + \nicefrac {1}{\tilde{q}_2} = \nicefrac 12$
there exists a $C \in (0,\infty)$ such that for 
all $ v_1, v_2, v_3 \in L^0(\Omega,\R^m)$ it holds that
\begin{align}
\nonumber
			&\int_{\Omega} \int_{\Omega}
				\frac
					{
						\|v_1(x)\|^2_{\R^m} \|v_2(y)\|^2_{\R^m} \|v_3(x) -v_3(y)\|_{\R^m}^2
					}
					{\|x-y\|_{\R^d}^{2s+d}}
			\ud x \ud y \\
	\label{eq: 3 func xy estimate q3 = infty}
		\leq{}
			&\int_{\Omega} \int_{\Omega}
					\frac
						{
							\|
								v_1(x)
							\|_{\R^m}^2
							\|
								v_2(y)
							\|_{\R^m}^{2}
						}
						{\|x-y\|_{\R^d}^{d-2\eps}}
				\ud x \ud y 
			\cdot \sup_{x,y \in \Omega}
				\frac
						{
							\|
								v_3(x) -v_3(y)
							\|_{\R^m}^{2}
						}
						{\|x-y\|_{\R^d}^{2s + 2\eps}}\\ \nonumber
		\leq{} 
			&C\|v_1\|^2_{L^{\tilde{q}_1}(\Omega,\R^m)} \|v_2\|^2_{L^{\tilde{q}_2}(\Omega,\R^m)}
			\cdot \|v_3\|^2_{C^{s+\eps}(\Omega,\| \cdot \|_{\R^m})} \\ \nonumber
		\leq{}
			&C\|v_1\|^2_{L^{\tilde{q}_1}(\Omega,\R^m)} \|v_2\|^2_{L^{\tilde{q}_2}(\Omega,\R^m)}
				\|v_3\|^2_{B^{s+\eps+d/2}_{2,2}(\Omega,\R^m)}.
\end{align}
Combining \eqref{eq: 3 func xy estimate q1q2 < infty}
	(with $v_1 \leftarrow v_1 \|v_1\|_{\R^m}^{p-1}$),
\eqref{eq: 3 func xy estimate q1 = infty q2 < infty}
	(with $v_1 \leftarrow v_1 \|v_1\|_{\R^m}^{p-1}$),
\eqref{eq: 3 func xy estimate q1q2 = infty}
	(with $v_1 \leftarrow v_1 \|v_1\|_{\R^m}^{p-1}$),
\eqref{eq: 3 func xy estimate q3 = infty}
	(with $v_1 \leftarrow v_1 \|v_1\|_{\R^m}^{p-1}$),
and the Sobolev inequality 
(see, e.g., the Theorem on page 31
and Theorem 1 on page 32 in Runst \& Sickel \cite{RunstSickel1996})
yields that for all 
$p \in (0,\infty)$
and all $\tilde{q}_1, \tilde{q}_2, \tilde{q}_3 \in [2,\infty]$
with 
$ 
		\nicefrac {1}{\tilde{q}_1} + \nicefrac{1}{\tilde{q}_2} + \nicefrac {1}{\tilde{q}_3} 
	= 
		\nicefrac 12
$
and with $p \cdot \tilde{q}_1 \geq 2$
there exist $C_1, C_2, C_3, C_4 \in (0,\infty)$ such that for 
all $ v_1, v_2, v_3 \in L^0(\Omega,\R^m)$
it holds that
\begin{equation}
\label{eq: 3 func xy estimate}
	\begin{split}
			&\int_\Omega \int_\Omega
				\frac
					{
						\|v_1(x)\|_{\R^m}^{2p} \|v_2(y)\|^2_{\R^m} \|v_3(x)- v_3(y)\|_{\R^m}^2
					}
					{\|x-y\|_{\R^d}^{2s+d}}
			\ud x \ud y \\
		\leq{}
			&C_1 \big \|v_1 \cdot  \|v_1\|_{\R^m}^{p-1} \big \|^2_{L^{\tilde{q}_1}(\Omega,\R^m)} 
			\|v_2\|^2_{L^{\tilde{q}_2}(\Omega,\R^m)}  
			\|v_3\|^2_{
				B^{s+d/2-d/\tilde{q}_3+\eps \1_{(2,\infty]}(\tilde{q}_3)}_{2,2}(\Omega,\R^m)
			} \\
		\leq{}
			&C_2 \big \| \|v_1\|_{\R^m}^{p} \big \|^2_{L^{\tilde{q}_1}(\Omega,\R)} 
			\|v_2\|^2_{L^{\tilde{q}_2}(\Omega,\R^m)}  
			\|v_3\|^2_{
				B^{s+d/2-d/\tilde{q}_3+\eps \1_{(2,\infty]}(\tilde{q}_3)}_{2,2}(\Omega,\R^m)
			} \\
		\leq{} 
			&C_3 \|v_1\|^{2p}_{L^{p\tilde{q}_1}(\Omega,\R^m)} 
			\|v_2\|^2_{L^{\tilde{q}_2}(\Omega,\R^m)}  
			\|v_3\|^2_{
				B^{s+d/2-d/\tilde{q}_3+\eps \1_{(2,\infty]}(\tilde{q}_3)}_{2,2}(\Omega,\R^m)
			} \\
		\leq{} 
			&C_4
				\|
					v_1
				\|^{2p}_{
					B^{d/2-d/(p \tilde{q}_1)+\eps \1_{\{\infty\}}(\tilde{q}_1)}_{2,2}(\Omega,\R^m)
				}
				\|v_2\|^2_{
					B^{d/2-d/\tilde{q}_2+\eps \1_{\{\infty\}}(\tilde{q}_2)}_{2,2}(\Omega,\R^m)
				}
			\\ &
				\cdot \|v_3\|^2_{
					B^{s+d/2-d/\tilde{q}_3+\eps \1_{(2,\infty]}(\tilde{q}_3)}_{2,2}(\Omega,\R^m)
				}.
	\end{split}
\end{equation}
Analogously it follows that for all 
$p \in (0,\infty)$
and all $\tilde{q}_1, \tilde{q}_2, \tilde{q}_3 \in [2,\infty]$
with 
$ 
		\nicefrac {1}{\tilde{q}_1} + \nicefrac{1}{\tilde{q}_2} + \nicefrac {1}{\tilde{q}_3} 
	= 
		\nicefrac 12
$
and with $p \cdot \tilde{q}_1 \geq 2$
there exists a $C \in (0,\infty)$ such that for 
all $ v_1, v_2, v_3 \in L^0(\Omega,\R^m)$
it holds that
\begin{align}
\nonumber
			&\int_\Omega \int_\Omega
				\frac
					{
						\|v_1(x)\|_{\R^m}^{2p}
						\|v_2(x)\|^2_{\R^m} 
						\|v_3(x)- v_3(y)\|_{\R^m}^2
					}
					{\|x-y\|_{\R^d}^{2s+d}}
			\ud x \ud y \\
	\label{eq: 3 func xx estimate}
		\leq{} 
			&C 
				\|
					v_1
				\|^{2p}_{
					B^{d/2-d/(p \tilde{q}_1)+\eps \1_{\{\infty\}}(\tilde{q}_1)}_{2,2}(\Omega,\R^m)
				}
				\|v_2\|^2_{
					B^{d/2-d/\tilde{q}_2+\eps \1_{\{\infty\}}(\tilde{q}_2)}_{2,2}(\Omega,\R^m)
				}
			\\ \nonumber & \cdot
				\|v_3\|^2_{
					B^{s+d/2-d/\tilde{q}_3+\eps \1_{(2,\infty]}(\tilde{q}_3)}_{2,2}(\Omega,\R^m)
				}.
\end{align}
Next note that for all $x,y \in \Omega$ and all $u, \tilde{u} \in L^0(\Omega,\R^m)$
it holds that
\begin{align}
\label{eq: b representation}
			&b(x,u(x)) - b(y,u(y)) - b(x,\tilde{u}(x)) + b(y, \tilde{u}(y)) \\ \nonumber
		={} 
			&\int_0^1 
				\big \langle 
					(D_{\R^d \times \R^m} b )\big( tx + (1-t) y, t u(x) + (1-t) u(y) \big),
			\\& \nonumber \qquad
					\big( x-y, u(x)-u(y) \big)
				\big \rangle_{(\R^d \times \R^m)', \R^d \times \R^m}
			\ud t \\ \nonumber
			&-\int_0^1
				\big \langle 
					(D_{\R^d \times \R^m}b) 
						\big( tx + (1-t) y, t \tilde{u}(x) + (1-t) \tilde{u}(y) \big),
				\\& \nonumber \qquad
					\big( x-y, \tilde{u}(x)-\tilde{u}(y) \big)
				\big \rangle_{(\R^d \times \R^m)', \R^d \times \R^m}
			\ud t \\ \nonumber
		={} 
			&\int_0^1 \int_0^1
				\Big \langle 
					b''_2 \big (
						t x + (1-t) y, 
						\lambda t u(x) + \lambda (1-t) u(y)
						+ (1-\lambda) t \tilde{u}(x) + (1-\lambda) (1-t) \tilde{u}(y)
					\big) 
				\\ & \nonumber \qquad
						\big (t (u(x) - \tilde{u}(x))+ (1-t) (u(y) - \tilde{u}(y)) \big),
					\big( x-y, u(x)-u(y) \big)
				\Big \rangle_{(\R^d \times \R^m)', \R^d \times \R^m}
			\ud \lambda \ud t \\
	\nonumber
			&+\int_0^1
				\big \langle 
					(D_{\R^d \times \R^m} b) 
						\big( tx + (1-t) y, t \tilde{u}(x) + (1-t) \tilde{u}(y) \big),
				\\& \qquad
			\nonumber
					\big( 0, u(x)-\tilde{u}(x)-u(y) +\tilde{u}(y) \big)
				\big \rangle_{(\R^d \times \R^m)', \R^d \times \R^m}
			\ud t.
\end{align}
We will now estimate the first term. Therefore we divide the proof into 3 cases.
\\
{\it 1. case $p_1=0$.}
The assumption \eqref{eq: b estimate 2nd derivative} 
implies that for all $x,y\in \Omega$ 
and all $u, \tilde{u} \in L^0(\Omega,\R^m)$ it holds that
\begin{align}
\label{eq: b 2nd derivative estimate p=0}
			&\int_0^1 \int_0^1
				\Big \langle 
					b''_2 \big (
						t x + (1-t) y, 
						\lambda t u(x) + \lambda (1-t) u(y)
						+ (1-\lambda) t \tilde{u}(x) + (1-\lambda) (1-t) \tilde{u}(y)
					\big) 
				\\ \nonumber & \qquad 
					 \big (t (u(x) - \tilde{u}(x))+ (1-t) (u(y) - \tilde{u}(y)) \big),
					\big( x-y, u(x)-u(y) \big)
				\Big \rangle_{(\R^d \times \R^m)', \R^d \times \R^m}
			\ud \lambda \ud t \\ \nonumber
		\leq{} 
			&\int_0^1 \int_0^1
				\Big \|
					b''_2 \big (
						t x + (1-t) y, 
						\lambda t u(x) + \lambda (1-t) u(y)
						+ (1-\lambda) t \tilde{u}(x) 
					\\ \nonumber & \qquad \qquad	
						+ (1-\lambda) (1-t) \tilde{u}(y)
					\big)
				\Big \|_{L(\R^m,(\R^d \times \R^m)'}
				\\ \nonumber & \qquad  
						\cdot \big \| t (u(x) - \tilde{u}(x))+ (1-t) (u(y) - \tilde{u}(y)) \big\|_{\R^m}
					\cdot \big \| (x-y, u(x)-u(y)) \big \|_{\R^d \times \R^m}
			\ud \lambda \ud t\\ \nonumber
		\leq{} 
			&\int_0^1 \int_0^1
				K
			\ud \lambda
				\cdot \big( t\|u(x) - \tilde{u}(x)\|_{\R^m} 
				+ (1-t)  \| u(y) - \tilde{u}(y) \|_{\R^m} \big)
			\\ \nonumber & \qquad
				\cdot \big \| (x-y, u(x)-u(y)) \big \|_{\R^d \times \R^m}
			\ud t\\ \nonumber
		={} 
			&K \cdot\big \| (x-y, u(x)-u(y)) \big \|_{\R^d \times \R^m}
			\int_0^1
				\big( 
					t\|u(x) - \tilde{u}(x)\|_{\R^m} + (1-t) \|  u(y) - \tilde{u}(y) \|_{\R^m} 
				\big)
			\ud t\\ \nonumber
		={} 
			&\tfrac {K}{2} \cdot\big \| (x-y, u(x)-u(y)) \big \|_{\R^d \times \R^m}
			\big(\|u(x) - \tilde{u}(x)\|_{\R^m} + \|  u(y) - \tilde{u}(y) \|_{\R^m} \big).
\end{align}
We thus get from \eqref{eq: b 2nd derivative estimate p=0},
\eqref{eq: simple estm xx}
	(with $v_1 \leftarrow u-\tilde{u}$, $v_2 \leftarrow {\bf 1}$, 
	$\tilde{q}_1 \leftarrow 2$, 
	$\tilde{q}_2 \leftarrow \infty$, and with $t \leftarrow 1$),
from
\eqref{eq: 3 func xy estimate} 
	(with $p \leftarrow 1$, $v_1 \leftarrow u - \tilde{u}$,
	$v_2 \leftarrow {\bf 1}$, $v_3 \leftarrow u$, 
	$
			\tilde{q}_1 
		\leftarrow 
			q_1
	$,
	and with
	$
			\tilde{q}_2 
		\leftarrow 
			\infty 
	$),
from the fact that 
for all $r\in \R$ and all $t\in (-\infty, r]$
there exists a $C \in (0,\infty)$ such that for all 
$u \in B^t_{2,2}(\Omega,\R^m)$ it holds that
$\|u\|_{B^t_{2,2}(\Omega,\R^m)} \leq C \|u\|_{B^r_{2,2}(\Omega,\R^m)}$,
and	from the fact that $\tilde{q}_3 \leq q_2$
that  for all  
$\tilde{q}_3 \in [2,\infty]$
with $\nicefrac {1}{q_1}+ \nicefrac{1} {\tilde{q}_3} = \nicefrac 12$
there exist $C_1, C_2 \in (0,\infty)$
such that for all $u, \tilde{u} \in L^0(\Omega,\R^m)$ it holds that
\begin{align}
\label{eq: first term p= 0}
			&\int_\Omega \int_\Omega \bigg | \int_0^1 \int_0^1
				\Big \langle 
					b''_2 \big (
						t x + (1-t) y, 
						\lambda t u(x) + \lambda (1-t) u(y)
						+ (1-\lambda) t \tilde{u}(x) + (1-\lambda) (1-t) \tilde{u}(y)
					\big) 
				\\ & \nonumber \qquad
						\big (t (u(x) - \tilde{u}(x))+ (1-t) (u(y) - \tilde{u}(y)) \big),
					( x-y, u(x)-u(y))
				\Big \rangle_{(\R^d \times \R^m)', \R^d \times \R^m}
			\ud \lambda \ud t \bigg | ^2
		\\& \nonumber \qquad \qquad
			 \tfrac{1}{\|x-y \|^{2s+d}_{\R^d}}\ud x \ud y\\
		\leq{} \nonumber
			&\tfrac {K^2}{4}
			\int_\Omega \int_\Omega
				\frac
					{
						\| (x-y, u(x)-u(y)) \|^2_{\R^d \times \R^m}
						\big(\|u(x) - \tilde{u}(x)\|_{\R^m} + \|u(y) - \tilde{u}(y)\|_{\R^m}\big)^2
					}
					{
						\|x-y\|^{2s+d}_{\R^d}
					}
			\ud x \ud y \\
		\leq{} \nonumber
			&\tfrac {K^2}{2}
			\int_\Omega \int_\Omega
				\frac
					{
						\| x-y\|^2_{\R^d }
						\big(\|u(x) - \tilde{u}(x)\|_{\R^m}^2 
						+ \|  u(y) - \tilde{u}(y) \|_{\R^m}^2 \big)
					}
					{
						\|x-y\|^{2s+d}_{\R^d}
					}
			\ud x \ud y \\ \nonumber
			&+\tfrac {K^2}{2}
			\int_\Omega \int_\Omega
				\frac
					{
						\| u(x)-u(y)\|_{\R^m}^2
						\big( \|u(x) - \tilde{u}(x)\|_{\R^m}^2 
						+ \|  u(y) - \tilde{u}(y) \|_{\R^m}^2 \big)
					}
					{
						\|x-y\|^{2s+d}_{\R^d}
					}
			\ud x \ud y \\ 
		={} \nonumber
			&K^2
			\int_\Omega \int_\Omega
				\frac
					{
						\| x-y\|^2_{\R^d }
						\|u(x) - \tilde{u}(x)\|_{\R^m}^2
					}
					{
						\|x-y\|^{2s+d}_{\R^d}
					}
			\ud x \ud y 
		\\ \nonumber &
			+K^2
			\int_\Omega \int_\Omega
				\frac
					{
						\| u(x)-u(y)\|_{\R^m}^2
						\|u(x) - \tilde{u}(x)\|_{\R^m}^2
					}
					{
						\|x-y\|^{2s+d}_{\R^d}
					}
			\ud x \ud y \\ 
		\leq{} \nonumber
			&C_1 \|u - \tilde{u}\|^2_{L^2(\Omega,\R^m)}
			+C_1
				\|
					u-\tilde{u}
				\|^{2}_{B^{d/2-d/q_1+\eps \1_{\{\infty\}}(q_1)}_{2,2}(\Omega,\R^m)}
				\|
					u
				\|^2_{
					B^{s+d/2-d/\tilde{q}_3+\eps \1_{(2,\infty]}(\tilde{q}_3)}_{2,2}(\Omega,\R^m)
				} \\
		\leq{} \nonumber
			&C_2 
				\|
					u-\tilde{u}
				\|^{2}_{B^{d/2-d/q_1+\eps \1_{\{\infty\}}(q_1)}_{2,2}(\Omega,\R^m)}
				\big(
					\|u\|^2_{B^{s+d/2-d/q_2+\eps \1_{(2,\infty]}(q_2)}_{2,2}(\Omega,\R^m)}
					+1
				\big) \\ \nonumber
		={} 
			&C_2
			\|
				u-\tilde{u}
			\|^{2}_{B^{d/2-d/q_1+\eps\1_{\{\infty\}}(q_1)}_{2,2}(\Omega,\R^m)}
				\big(
					\|
						u
					\|^2_{B^{s+d/2-d/q_2+\eps \1_{(2,\infty]}(q_2)}_{2,2}(\Omega,\R^m)}
					+1
				\big)
			\cdot K_1(0, u , \tilde{u}).
\end{align}
\\
{\it 2. case $p_1 \in (0, \infty)$.}
We then get from the assumption \eqref{eq: b estimate 2nd derivative} 
that for all $x,y \in \Omega$, $t \in [0,1]$ and all $u, \tilde{u} \in L^0(\Omega,\R^m)$ 
it holds that
\begin{align}
\label{eq: int of lambda 0<p<infty} 
			&\int_0^1
				\Big \langle 
					b''_2 \big (
						t x + (1-t) y, 
						\lambda t u(x) + \lambda (1-t) u(y)
						+ (1-\lambda) t \tilde{u}(x) + (1-\lambda) (1-t) \tilde{u}(y)
					\big) 
				\\ \nonumber & \qquad 
						\big (t (u(x) - \tilde{u}(x))+ (1-t) (u(y) - \tilde{u}(y)) \big),
					\big( x-y, u(x)-u(y) \big)
				\Big \rangle_{(\R^d \times \R^m)', \R^d \times \R^m}
			\ud \lambda \\ \nonumber
		\leq{} 
			&\int_0^1
				\Big \|
					b''_2 \big (
						t x + (1-t) y, 
						\lambda t u(x) + \lambda (1-t) u(y)
				\\ \nonumber& \qquad \qquad
						+ (1-\lambda) t \tilde{u}(x) + (1-\lambda) (1-t) \tilde{u}(y)
					\big)
				\Big \|_{L(\R^m, (\R^d \times \R^m)')}
				\\ \nonumber & \qquad
						\cdot \big \| t (u(x) - \tilde{u}(x))+ (1-t) (u(y) - \tilde{u}(y)) \big\|_{\R^m}
					\cdot \big \| (x-y, u(x)-u(y)) \big \|_{\R^d \times \R^m}
			\ud \lambda \\ \nonumber
		\leq{} 
			&\int_0^1
				K\big( 1 + 
					\big \|
						\lambda t u(x) + \lambda (1-t) u(y)
						+ (1-\lambda) t \tilde{u}(x) + (1-\lambda) (1-t) \tilde{u}(y)
					\big\|_{\R^m}
				\big)^{p_1}
			\ud \lambda
			\\ \nonumber & \qquad
					\cdot \big \| t (u(x) - \tilde{u}(x))+ (1-t) (u(y) - \tilde{u}(y)) \big\|_{\R^m}
				\cdot \big \| (x-y, u(x)-u(y)) \big \|_{\R^d \times \R^m} \\ \nonumber
		\leq{} 
			&\int_0^1
				5^{p_1} \cdot K \big( 
						1 + 
						\lambda^{p_1} t^{p_1} \|u(x)\|_{\R^m}^{p_1} 
						+ \lambda^{p_1} (1-t)^{p_1} \|u(y)\|_{\R^m}^{p_1}
				\\ \nonumber &  \qquad
						+ (1-\lambda)^{p_1} t^{p_1} \| \tilde{u}(x) \|_{\R^m}^{p_1} 
						+ (1-\lambda)^{p_1} (1-t)^{p_1} \|\tilde{u}(y) \|_{\R^m}^{p_1}
				\big)
			\ud \lambda
			\\ \nonumber &
					\cdot \big \| t (u(x) - \tilde{u}(x))+ (1-t) (u(y) - \tilde{u}(y)) \big\|_{\R^m}
				\cdot \big \| (x-y, u(x)-u(y)) \big \|_{\R^d \times \R^m} \\ \nonumber
		={} 
			&\tfrac {5^{p_1}}{p_1+1} \cdot K \big( 
						p_1+1 
						+t^{p_1} \|u(x)\|_{\R^m}^{p_1} + (1-t)^{p_1} \|u(y)\|_{\R^m}^{p_1}
						+ t^{p_1} \| \tilde{u}(x) \|_{\R^m}^{p_1} 
						+ (1-t)^{p_1} \|\tilde{u}(y) \|_{\R^m}^{p_1}
				\big)
			\\ \nonumber & \qquad
					\cdot \big \| t (u(x) - \tilde{u}(x))+ (1-t) (u(y) - \tilde{u}(y)) \big\|_{\R^m}
				\cdot \big \| (x-y, u(x)-u(y)) \big \|_{\R^d \times \R^m}.
\end{align}
In addition, it holds for all $x, y \in \Omega$ and all $u, \tilde{u} \in L^0(\Omega,\R^m)$ 
that 
\begin{align}
\nonumber
			&\int_0^1
				\big( 
						p_1+1 
						+t^{p_1} \|u(x)\|_{\R^m}^{p_1} + (1-t)^{p_1} \|u(y)\|_{\R^m}^{p_1}
						+ t^{p_1} \| \tilde{u}(x) \|_{\R^m}^{p_1} 
						+ (1-t)^{p_1} \|\tilde{u}(y) \|_{\R^m}^{p_1}
					\big)
			\\ & \nonumber \qquad
					\cdot \big \| t (u(x) - \tilde{u}(x))+ (1-t) (u(y) - \tilde{u}(y)) \big\|_{\R^m}
			\ud t \\
		\leq{}  \nonumber
		 &\int_0^1
				t (p_1+1) \cdot \| u(x) - \tilde{u}(x)\|_{\R^m} 
				+ (1-t) (p_1+1) \cdot\| u(y) - \tilde{u}(y) \|_{\R^m}
			\ud t \\ \nonumber
			&+\int_0^1
				t^{p_1+1} ( \|u(x) \|_{\R^m}^{p_1}+ \| \tilde{u}(x) \|_{\R^m}^{p_1} )
					\cdot \| u(x) - \tilde{u}(x) \|_{\R^m}
			\ud t \\ \nonumber
			&+\int_0^1
				(1-t) t^{p_1} (\| u(x) \|_{\R^m}^{p_1}+ \| \tilde{u}(x) \|_{\R^m}^{p_1} )
					\cdot \| u(y) - \tilde{u}(y) \|_{\R^m}
			\ud t \\ 
		\label{eq: int of t 0<p<infty}
		\begin{split}
			&+\int_0^1
				(1-t)^{p_1}t (\|u(y)\|_{\R^m}^{p_1}+ \| \tilde{u}(y) \|_{\R^m}^{p_1} )
					\cdot \| u(x) - \tilde{u}(x) \|_{\R^m}
			\ud t  \\ 
			&+\int_0^1
				(1-t)^{p_1+1} (\|u(y)\|_{\R^m}^{p_1}+ \| \tilde{u}(y) \|_{\R^m}^{p_1} )
					\cdot \| u(y) - \tilde{u}(y) \|_{\R^m}
			\ud t 
		\end{split} \\ \nonumber
		={} 
		 &\tfrac{(p_1+1)}{2} \cdot
			(
				\| u(x) - \tilde{u}(x)\|_{\R^m}
				+\|u(y) - \tilde{u}(y) \|_{\R^m}
			)  \\ \nonumber
			&+\tfrac{1}{p_1+2} (\|u(x)\|_{\R^m}^{p_1}+ \| \tilde{u}(x) \|_{\R^m}^{p_1} )
					\cdot \| u(x) - \tilde{u}(x) \|_{\R^m} \\ \nonumber
			&+\tfrac{1}{(p_1+1)(p_1+2)}
				(\|u(x)\|_{\R^m}^{p_1}+ \| \tilde{u}(x) \|_{\R^m}^{p_1} )
					\cdot \| u(y) - \tilde{u}(y) \|_{\R^m} \\ \nonumber
			&+\tfrac{1}{(p_1+1)(p_1+2)}
				(\|u(y)\|_{\R^m}^{p_1}+ \| \tilde{u}(y) \|_{\R^m}^{p_1} )
					\cdot \| u(x) - \tilde{u}(x) \|_{\R^m} \\ \nonumber
			&+\tfrac{1}{p_1+2}
				(\|u(y)\|_{\R^m}^{p_1}+ \| \tilde{u}(y) \|_{\R^m}^{p_1} )
					\cdot \| u(y) - \tilde{u}(y) \|_{\R^m}.
\end{align}
Furthermore it follows from
\eqref{eq: simple estm xx}
	(with $v_1 \leftarrow u-\tilde{u}$, $v_2 \leftarrow {\bf 1}$, $\tilde{q}_1 \leftarrow 2$, 
	$\tilde{q}_2 \leftarrow \infty$, and with $t \leftarrow 1$),
\eqref{eq: 3 func xy estimate} 
	(with $p \leftarrow 1$, $v_1 \leftarrow u - \tilde{u}$,
	$v_2 \leftarrow {\bf 1}$, $v_3 \leftarrow u$, 
	$\tilde{q}_1 \leftarrow q_1$,
	and with 
	$
		\tilde{q}_2 \leftarrow \infty
	$),
from the fact that 
for all $r\in \R$ and all $t\in (-\infty, r]$
there exists a $C \in (0,\infty)$ such that for all 
	$u \in B^t_{2,2}(\Omega,\R^m)$ it holds that
$\|u\|_{B^t_{2,2}(\Omega,\R^m)} \leq C \|u\|_{B^r_{2,2}(\Omega,\R^m)}$,
and from the fact that $\tilde{q}_3 \leq q_2$
that  for all 
$\tilde{q}_3 \in [2,\infty]$
with $\nicefrac {1}{q_1}+ \nicefrac{1} {\tilde{q}_3} = \nicefrac 12$
there exist $C_1, C_2 \in (0,\infty)$
such that for all $u, \tilde{u} \in L^0(\Omega,\R^m)$ it holds that
\begin{align}
	\label{eq: 1st term of b estimate 0<p<infty} 
			&\int_\Omega \int_\Omega
				\frac
					{
						\| (x-y, u(x)-u(y)) \|^2_{\R^d \times \R^m}
						\big(
							\|u(x) - \tilde{u}(x)\|_{\R^m} + \|  u(y) - \tilde{u}(y)\|_{\R^m}
						\big)^2
					}
					{
						\|x-y\|^{2s+d}_{\R^d}
					}
			\ud x \ud y \\ \nonumber
		\leq{} 
			&2
			\int_\Omega \int_\Omega
				\frac
					{
						\| x-y\|^2_{\R^d }
						\big( \|u(x) - \tilde{u}(x)\|_{\R^m}^2
						+ \|  u(y) - \tilde{u}(y) \|_{\R^m}^2 \big)
					}
					{
						\|x-y\|^{2s+d}_{\R^d}
					}
			\ud x \ud y \\ \nonumber
			&+2
			\int_\Omega \int_\Omega
				\frac
					{
						\| u(x)-u(y)\|_{\R^m}^2
						\big(\|u(x) - \tilde{u}(x)\|_{\R^m}^2 
						+ \|  u(y) - \tilde{u}(y) \|_{\R^m}^2 \big)
					}
					{
						\|x-y\|^{2s+d}_{\R^d}
					}
			\ud x \ud y \\ \nonumber
		={} 
			&4
			\int_\Omega \int_\Omega
				\frac
					{
						\| x-y\|^2_{\R^d }
						\|u(x) - \tilde{u}(x)\|_{\R^m}^2
					}
					{
						\|x-y\|^{2s+d}_{\R^d}
					}
			\ud x \ud y \\ \nonumber
			&+4
			\int_\Omega \int_\Omega
				\frac
					{
						\| u(x)-u(y)\|_{\R^m}^2
						\|u(x) - \tilde{u}(x)\|_{\R^m}^2
					}
					{
						\|x-y\|^{2s+d}_{\R^d}
					}
			\ud x \ud y \\ \nonumber
		\leq{} 
			&C_1 \|u - \tilde{u}\|^2_{L^2(\Omega,\R^m)}
			+C_1
				\|
					u-\tilde{u}
				\|^{2}_{B^{d/2-d/q_1+\eps \1_{\{\infty\}}(q_1)}_{2,2}(\Omega,\R^m)}
				\|u\|^2_{
					B^{s+d/2-d/\tilde{q}_3+\eps \1_{(2,\infty]}(\tilde{q}_3)
				}_{2,2}(\Omega,\R^m)} \\
	\nonumber
		\leq{} 
			&C_2
				\|
					u-\tilde{u}
				\|^{2}_{B^{d/2-d/q_1+\eps \1_{\{\infty\}}(q_1)}_{2,2}(\Omega,\R^m)}
				\big(
					\|u\|^2_{B^{s+d/2-d/q_2+\eps \1_{(2,\infty]}(q_2)}_{2,2}(\Omega,\R^m)} 
					+1
				\big).
\end{align}
Moreover, we obtain from \eqref{eq: simple estm xx}
	(with $v_1 \leftarrow u \cdot\|u\|_{\R^m}^{p_1-1}$ 
		resp.\@ $ \tilde{u} \cdot \|\tilde{u}\|_{\R^m}^{p_1-1}$,
		$v_2 \leftarrow u-\tilde{u}$, $\tilde{q}_1 \leftarrow q_3$,
		and with
		$t \leftarrow 1$) and from
\eqref{eq: 3 func xx estimate} 
	(with $p \leftarrow p_1$, $v_1 \leftarrow u$ resp.\@ $\tilde{u}$,
	$v_2 \leftarrow u-\tilde{u}$, $v_3 \leftarrow u$,
	$\tilde{q}_1 \leftarrow q_3$,
	$\tilde{q}_2 \leftarrow q_1$,
	and with $\tilde{q}_3 \leftarrow q_2$),
that for all 
$\tilde{q}_2 \in [2, \infty]$
with 
$
		\nicefrac {1}{\tilde{q}_2} + \nicefrac{1}{q_3} 
	= 
		\nicefrac 12
$
there exist $C_1,C_2 \in (0, \infty)$
such that for all $u, \tilde{u} \in L^0(\Omega,\R^m)$ it holds that
\begin{align}
			&\int_\Omega \int_\Omega
			\frac
				{
					(\|u(x)\|_{\R^m}^{p_1}+ \| \tilde{u}(x) \|_{\R^m}^{p_1} )^2
					\cdot \| u(x) - \tilde{u}(x) \|_{\R^m}^2 
					\cdot \big \| (x-y, u(x)-u(y)) \big \|^2_{\R^d \times \R^m}
				}
				{
					\|x-y\|^{2s+d}_{\R^d}
				}
			\ud x \ud y \\ \nonumber
		\leq{} 
			&2\int_\Omega \int_\Omega
			\frac
				{
					(\|u(x)\|_{\R^m}^{2p_1}+ \| \tilde{u}(x) \|_{\R^m}^{2p_1} )
					\cdot \| u(x) - \tilde{u}(x) \|_{\R^m}^2 \cdot \|x-y\|^{2}_{\R^d}
				}
				{
					\|x-y\|^{2s+d}_{\R^d}
				}
			\ud x \ud y \\ \nonumber
			&+2\int_\Omega \int_\Omega
			\frac
				{
					(\|u(x)\|_{\R^m}^{2p_1}+ \| \tilde{u}(x) \|_{\R^m}^{2p_1} )
					\cdot \| u(x) - \tilde{u}(x) \|_{\R^m}^2 
					\cdot \big \| u(x)-u(y) \big \|_{\R^m}^2
				}
				{
					\|x-y\|^{2s+d}_{\R^d}
				}
			\ud x \ud y \\ \nonumber
		\leq{} 
			&C_1(
				\| u \cdot  \|u\|_{\R^m}^{p_1-1}\|^2_{L^{q_3}(\Omega,\R^m)}
				+ \| \tilde{u} \cdot \|\tilde{u}\|_{\R^m}^{p_1-1}\|^2_{L^{q_3}(\Omega,\R^m)}
			) 
				\cdot \| u - \tilde{u} \|^2_{L^{\tilde{q}_2}(\Omega,\R^m)} \\ \nonumber
			&+C_1\big(
				\|
					u
				\|^{2p_1}_{
					B^{d/2-d/(p_1 q_3)+\eps \1_{\{\infty\}}(q_3)}_{2,2}(\Omega,\R^m)
				}
				+\|
					\tilde{u}
				\|^{2p_1}_{
					B^{d/2-d/(p_1 q_3)+\eps \1_{\{\infty\}}(q_3)}_{2,2}(\Omega,\R^m)
				}
			\big ) 
		\\ \nonumber &
				\cdot 
				\|
					u-\tilde{u}
				\|^2_{B^{d/2-d/q_1+\eps \1_{\{\infty\}}(q_1)}_{2,2}(\Omega,\R^m)} 
				\|
					u
				\|^2_{B^{s+d/2-d/q_2+\eps \1_{(2,\infty]}(q_2)}_{2,2}(\Omega,\R^m)} \\ \nonumber
		\leq{} 
			&C_2 (
				\| u\|^{2p_1}_{L^{p_1 q_3}(\Omega,\R^m)}
				+ \| \tilde{u}\|^{2p_1}_{L^{p_1 q_3}(\Omega,\R^m)}
			) 
				\cdot \| u - \tilde{u}  \|^2_{L^{\tilde{q}_2}(\Omega,\R^m)} \\ \nonumber
			&+C_1\big(
				\|
					u
				\|^{2p_1}_{
					B^{d/2-d/(p_1 q_3)+\eps \1_{\{\infty\}}(q_3)}_{2,2}(\Omega,\R^m)
				}
				+\|
					\tilde{u}
				\|^{2p_1}_{
					B^{d/2-d/(p_1 q_3)+\eps \1_{\{\infty\}}(q_3)}_{2,2}(\Omega,\R^m)
				}
			\big ) \\ \nonumber
				&\cdot 
				\|
					u-\tilde{u}
				\|^2_{B^{d/2-d/q_1+\eps \1_{\{\infty\}}(q_1)}_{2,2}(\Omega,\R^m)} 
				\|
					u
				\|^2_{B^{s+d/2-d/q_2+\eps \1_{(2,\infty]}(q_2)}_{2,2}(\Omega,\R^m)}. 
\end{align}
Thus the Sobolev inequality 
(see, e.g., the Theorem on page 31 and Theorem 1 on page 32 
in Runst \& Sickel \cite{RunstSickel1996}), the fact that $\tilde{q}_2 \leq q_1$,
and the fact that 
for all $r\in \R$ and all $t\in (-\infty, r]$
there exists a $C \in (0,\infty)$ such that for all 
	$u \in B^t_{2,2}(\Omega,\R^m)$ it holds that
$\|u\|_{B^t_{2,2}(\Omega,\R^m)} \leq C \|u\|_{B^r_{2,2}(\Omega,\R^m)}$
shows that for all
$\tilde{q}_2 \in [2, \infty]$
with 
$
		\nicefrac {1}{q_3} + \nicefrac{1}{\tilde{q}_2} 
	=
		\nicefrac 12
$
there exist $C_1, C_2, C_3 \in (0, \infty)$
such that for all $u, \tilde{u} \in L^0(\Omega,\R^m)$ it holds that
\begin{align}
\nonumber
			&\int_\Omega \int_\Omega
			\frac
				{
					(\|u(x)\|_{\R^m}^{p_1}+ \| \tilde{u}(x) \|_{\R^m}^{p_1} )^2
					\cdot \| u(x) - \tilde{u}(x) \|_{\R^m}^2 
				}
				{
					\|x-y\|^{2s+d}_{\R^d}
				}
		\\ & \nonumber
				\cdot \big \| (x-y, u(x)-u(y)) \big \|^2_{\R^d \times \R^m}
			\ud x \ud y \\
		\leq{} \nonumber
			&C_1(
				\| u\|^{2p_1}_{L^{p_1 q_3}(\Omega,\R^m)}
				+ \| \tilde{u}\|^{2p_1}_{L^{p_1 q_3}(\Omega,\R^m)}
			) 
				\cdot \| u - \tilde{u}  \|^2_{L^{\tilde{q}_2}(\Omega,\R^m)} \\ \nonumber
			&+C_1 \big(
				\|
					u
				\|^{2p_1}_{
					B^{d/2-d/(p_1 q_3)+\eps \1_{\{\infty\}}(q_3)}_{2,2}(\Omega,\R^m)
				}
				+\|
					\tilde{u}
				\|^{2p_1}_{
					B^{d/2-d/(p_1 q_3)+\eps \1_{\{\infty\}}(q_3)}_{2,2}(\Omega,\R^m)
				}
			\big) \\ \nonumber
			&	\cdot 
				\|
					u-\tilde{u}
				\|^2_{B^{d/2-d/q_1+\eps \1_{\{\infty\}}(q_1)}_{2,2}(\Omega, \R^m)} 
				\|
					u
				\|^2_{B^{s+d/2-d/q_2+\eps \1_{(2,\infty]}(q_2)}_{2,2}(\Omega, \R^m)} \\
		\leq{} 
	\label{eq: 2nd term of b estimate 0<p<infty} 
			&C_2(
				\|
					u
				\|^{2p_1}_{
					B^{d/2-d/(p_1 q_3)+\eps \1_{\{\infty\}}(q_3)}_{2,2}(\Omega,\R^m)
				}
				+\|
					\tilde{u}
				\|^{2p_1}_{
					B^{d/2-d/(p_1 q_3)+\eps \1_{\{\infty\}}(q_3)}_{2,2}(\Omega,\R^m)
				}
			) \\ \nonumber
			& \cdot
				\|
					u-\tilde{u}
				\|^2_{
					B^{d/2-d/\tilde{q}_2+\eps \1_{\{\infty\}}(\tilde{q}_2)}_{2,2}(\Omega,\R^m)
				} 
			+C_2(
				\|
					u
				\|^{2p_1}_{
					B^{d/2-d/(p_1 q_3)+\eps \1_{\{\infty\}}(q_3)}_{2,2}(\Omega,\R^m)
				}
			\\\nonumber &  
				+\|
					\tilde{u}
				\|^{2p_1}_{
					B^{d/2-d/(p_1 q_3)+\eps \1_{\{\infty\}}(q_3)}_{2,2}(\Omega,\R^m)
				}
			) 
				\cdot 
				\|
					u-\tilde{u}
				\|^2_{B^{d/2-d/q_1+\eps \1_{\{\infty\}}(q_1)}_{2,2}(\Omega,\R^m)} \\ \nonumber
				&\cdot 
				\|
					u
				\|^2_{B^{s+d/2-d/q_2+\eps \1_{(2,\infty]}(q_2)}_{2,2}(\Omega,\R^m)} \\ \nonumber
		\leq{} 
			&C_3(
				\|
					u
				\|^{2p_1}_{B^{
					d/2-d/(p_1 q_3)+\eps \1_{\{\infty\}}(q_3)}_{2,2}(\Omega,\R^m)
				}
				+\|
					\tilde{u}
				\|^{2p_1}_{
					B^{d/2-d/(p_1 q_3)+\eps \1_{\{\infty\}}(q_3)}_{2,2}(\Omega,\R^m)
				}
			) \\ \nonumber
			&	\cdot 
				\|
					u-\tilde{u}
				\|^2_{
					B^{d/2-d/q_1+\eps \1_{\{\infty\}}(q_1)}_{2,2}(\Omega,\R^m)
				} 
				\cdot 
				\Big(
					\|
						u
					\|^2_{B^{s+d/2-d/q_2+\eps \1_{(2,\infty]}(q_2)}_{2,2}(\Omega,\R^m)}
					+1
				\Big).
\end{align}
Analogously it follows that 
there exists a $C \in (0, \infty)$
such that for all $u, \tilde{u} \in L^0(\Omega,\R^m)$ it holds that
\begin{align}
\label{eq: 3rd term of b estimate 0<p<infty} 
			&\int_\Omega \int_\Omega
			\frac
				{
					(\|u(x)\|_{\R^m}^{p_1}+ \| \tilde{u}(x)\|_{\R^m}^{p_1} )^2
					\cdot \| u(y) - \tilde{u}(y) \|_{\R^m}^2
					\cdot \big \| (x-y, u(x)-u(y)) \big \|^2_{\R^d \times \R^m}
				}
				{
					\|x-y\|^{2s+d}_{\R^d}
				}
			\ud x \ud y \\ \nonumber
		\leq{} 
			&C(
				\|
					u
				\|^{2p_1}_{B^{
					d/2-d/(p_1 q_3)+\eps \1_{\{\infty\}}(q_3)}_{2,2}(\Omega,\R^m)
				}
				+\|
					\tilde{u}
				\|^{2p_1}_{
					B^{d/2-d/(p_1 q_3)+\eps \1_{\{\infty\}}(q_3)}_{2,2}(\Omega,\R^m)
				}
			) \\ \nonumber
				&\cdot 
				\|
					u-\tilde{u}
				\|^2_{
					B^{d/2-d/q_1+\eps \1_{\{\infty\}}(q_1)}_{2,2}(\Omega,\R^m)
				} 
				\cdot 
				\Big(
					\|
						u
					\|^2_{B^{s+d/2-d/q_2+\eps \1_{(2,\infty]}(q_2)}_{2,2}(\Omega,\R^m)}
					+1
				\Big),
\end{align}
that
\begin{align}
\label{eq: 4th term of b estimate 0<p<infty} 
			&\int_\Omega \int_\Omega
			\frac
				{
					(\|u(y)\|_{\R^m}^{p_1}+ \| \tilde{u}(y) \|_{\R^m}^{p_1} )^2
					\cdot \| u(x) - \tilde{u}(x) \|_{\R^m}^2
					\cdot \big \| (x-y, u(x)-u(y)) \big \|^2_{\R^d \times \R^m}
				}
				{
					\|x-y\|^{2s+d}_{\R^d}
				}
			\ud x \ud y \\ \nonumber
		\leq{} 
			&C(
				\|
					u
				\|^{2p_1}_{B^{
					d/2-d/(p_1 q_3)+\eps \1_{\{\infty\}}(q_3)}_{2,2}(\Omega,\R^m)
				}
				+\|
					\tilde{u}
				\|^{2p_1}_{
					B^{d/2-d/(p_1 q_3)+\eps \1_{\{\infty\}}(q_3)}_{2,2}(\Omega,\R^m)
				}
			) \\ \nonumber
				&\cdot 
				\|
					u-\tilde{u}
				\|^2_{
					B^{d/2-d/q_1+\eps \1_{\{\infty\}}(q_1)}_{2,2}(\Omega,\R^m)
				} 
				\cdot 
				\Big(
					\|
						u
					\|^2_{B^{s+d/2-d/q_2+\eps \1_{(2,\infty]}(q_2)}_{2,2}(\Omega,\R^m)}
					+1
				\Big),
\end{align}
and that
\begin{align}
\label{eq: 5th term of b estimate 0<p<infty} 
			&\int_\Omega \int_\Omega
			\frac
				{
					(\|u(y)\|_{\R^m}^{p_1}+ \| \tilde{u}(y) \|_{\R^m}^{p_1} )
					\cdot \| u(y) - \tilde{u}(y) \|_{\R^m}^2
					\cdot \big \| (x-y, u(x)-u(y)) \big \|^2_{\R^d \times \R^m}
				}
				{
					\|x-y\|^{2s+d}_{\R^d}
				}
			\ud x \ud y \\ \nonumber
		\leq{} 
			&C(
				\|
					u
				\|^{2p_1}_{B^{
					d/2-d/(p_1 q_3)+\eps \1_{\{\infty\}}(q_3)}_{2,2}(\Omega,\R^m)
				}
				+\|
					\tilde{u}
				\|^{2p_1}_{
					B^{d/2-d/(p_1 q_3)+\eps \1_{\{\infty\}}(q_3)}_{2,2}(\Omega,\R^m)
				}
			) \\ \nonumber
				&\cdot 
				\|
					u-\tilde{u}
				\|^2_{
					B^{d/2-d/q_1+\eps \1_{\{\infty\}}(q_1)}_{2,2}(\Omega,\R^m)
				} 
				\cdot 
				\Big(
					\|
						u
					\|^2_{B^{s+d/2-d/q_2+\eps \1_{(2,\infty]}(q_2)}_{2,2}(\Omega,\R^m)}
					+1
				\Big).
\end{align}
Combining 
\eqref{eq: int of lambda 0<p<infty},
\eqref{eq: int of t 0<p<infty}, 
\eqref{eq: 1st term of b estimate 0<p<infty},
\eqref{eq: 2nd term of b estimate 0<p<infty},
\eqref{eq: 3rd term of b estimate 0<p<infty},
\eqref{eq: 4th term of b estimate 0<p<infty},
and
\eqref{eq: 5th term of b estimate 0<p<infty}
establishes that
there exist a $C \in (0, \infty)$
such that for all $u, \tilde{u} \in L^0(\Omega,\R^m)$ it holds that
\begin{align}
\label{eq: first term 0 < p < infty}
			&\int_\Omega \int_\Omega \bigg | \int_0^1 \int_0^1
				\Big \langle 
					b''_2 \big (
						t x + (1-t) y, 
						\lambda t u(x) + \lambda (1-t) u(y)
						+ (1-\lambda) t \tilde{u}(x) + (1-\lambda) (1-t) \tilde{u}(y)
					\big) 
				\\ \nonumber & \qquad 
						\big (t (u(x) - \tilde{u}(x))+ (1-t) (u(y) - \tilde{u}(y)) \big),
					( x-y, u(x)-u(y))
				\Big \rangle_{(\R^d \times \R^m)', \R^d \times \R^m}
			\ud \lambda \ud t \bigg | ^2
		\\ \nonumber & \qquad
			\cdot \tfrac{1}{\|x-y \|^{2s+d}_{\R^d}}\ud x \ud y\\ \nonumber
		\leq{} 
			&C
				\|
					u-\tilde{u}
				\|^{2}_{B^{d/2-d/q_1+\eps \1_{\{\infty\}}(q_1)}_{2,2}(\Omega,\R^m)}
				\big(
					\|u\|^2_{B^{s+d/2-d/q_2+\eps \1_{(2,\infty]}(q_2)}_{2,2}(\Omega,\R^m)} 
					+1
				\big) \\ \nonumber
			&+C(
				\|
					u
				\|^{2p_1}_{B^{
					d/2-d/(p_1 q_3)+\eps \1_{\{\infty\}}(q_3)}_{2,2}(\Omega,\R^m)
				}
				+\|
					\tilde{u}
				\|^{2p_1}_{
					B^{d/2-d/(p_1 q_3)+\eps \1_{\{\infty\}}(q_3)}_{2,2}(\Omega,\R^m)
				} \\ \nonumber
			)&
				\cdot 
				\|
					u-\tilde{u}
				\|^2_{
					B^{d/2-d/q_1+\eps \1_{\{\infty\}}(q_1)}_{2,2}(\Omega,\R^m)
				} 
			\cdot 
				\Big(
					\|
						u
					\|^2_{B^{s+d/2-d/q_2+\eps \1_{(2,\infty]}(q_2)}_{2,2}(\Omega,\R^m)}
					+1
				\Big) \\ \nonumber
		={} 
			&C(
				\|
					u
				\|^{2p_1}_{B^{
					d/2-d/(p_1 q_3)+\eps \1_{\{\infty\}}(q_3)}_{2,2}(\Omega,\R^m)
				}
				+\|
					\tilde{u}
				\|^{2p_1}_{
					B^{d/2-d/(p_1 q_3)+\eps \1_{\{\infty\}}(q_3)}_{2,2}(\Omega,\R^m)
				}
				+1
			) \\ \nonumber
				&\cdot 
				\|
					u-\tilde{u}
				\|^2_{
					B^{d/2-d/q_1+\eps \1_{\{\infty\}}(q_1)}_{2,2}(\Omega,\R^m)
				} 
			\cdot 
				\Big(
					\|
						u
					\|^2_{B^{s+d/2-d/q_2+\eps \1_{(2,\infty]}(q_2)}_{2,2}(\Omega,\R^m)}
					+1
				\Big) \\ \nonumber
		={} 
			&C
			\|
				u-\tilde{u}
			\|^{2}_{B^{d/2-d/q_1+\eps\1_{\{\infty\}}(q_1)}_{2,2}(\Omega,\R^m)}
				\big(
					\|
						u
					\|^2_{B^{s+d/2-d/q_2+\eps \1_{(2,\infty]}(q_2)}_{2,2}(\Omega,\R^m)}
					+1
				\big)
			\cdot K_1(p_1, u , \tilde{u}).
\end{align}
\\
{\it 3. case $p_1= \infty$.}
We then get that for all $x,y\in \Omega$, $R \in [0, \infty]$, 
and all $u, \tilde{u} \in L^0(\Omega,\R^m)$ 
with 
$\|u\|_{L^\infty(\Omega,\R^m)} \vee \|\tilde{u}\|_{L^\infty(\Omega,\R^m)} \leq R $ 
it holds that
\begin{align}
\label{eq: b 2nd derivative estimate p= infty}
			&\int_0^1 \int_0^1
				\Big \langle 
					b''_2 \big (
						t x + (1-t) y, 
						\lambda t u(x) + \lambda (1-t) u(y)
						+ (1-\lambda) t \tilde{u}(x) 
						+ (1-\lambda) (1-t) \tilde{u}(y)
					\big) 
				\\ \nonumber & \qquad 
						\big (t (u(x) - \tilde{u}(x))+ (1-t) (u(y) - \tilde{u}(y)) \big),
					\big( x-y, u(x)-u(y) \big)
				\Big \rangle_{(\R^d \times \R^m)', \R^d \times \R^m}
			\ud \lambda \ud t\\ \nonumber
		\leq{} 
			&\int_0^1 \int_0^1
				\Big \|
					b''_2 \big (
						t x + (1-t) y, 
						\lambda t u(x) + \lambda (1-t) u(y)
						+ (1-\lambda) t \tilde{u}(x)
				\\ \nonumber & \qquad
						+ (1-\lambda) (1-t) \tilde{u}(y)
					\big)
				\Big \|_{L(\R^m,(\R^d \times \R^m)')}
				\\ \nonumber & \qquad
						\cdot 
							\big \| 
								t (u(x) - \tilde{u}(x))+ (1-t) (u(y) - \tilde{u}(y)) 
							\big\|_{\R^m}
					\cdot \big \| (x-y, u(x)-u(y)) \big \|_{\R^d \times \R^m}
			\ud \lambda \ud t \\ \nonumber
		\leq{} 
			&\int_0^1 \int_0^1
				\|
					 b''_2
				\|_{
						\C(
							\Omega 
							\times 
							B(R),
							\| \cdot \|_{L(\R^m,(\R^d \times \R^m)')}
						)
					}
			\ud \lambda
				\cdot 
					\big(
						t\|u(x) - \tilde{u}(x)\|_{\R^m} + (1-t) \|  u(y) - \tilde{u}(y) \|_{\R^m}
					\big)
		\\ \nonumber & \qquad
				\cdot \big \| (x-y, u(x)-u(y)) \big \|_{\R^d \times \R^m}
			\ud t\\ \nonumber
		={} 
			&\|
					 b''_2
				\|_{
						\C(
							\Omega 
							\times B(R),
							\| \cdot \|_{L(\R^m,(\R^d \times \R^m)')}
						)
					}
			\cdot\big \| (x-y, u(x)-u(y)) \big \|_{\R^d \times \R^m}
		\\ \nonumber &
			\int_0^1
				\big( 
					t \|u(x) - \tilde{u}(x) \|_{\R^m} + (1-t) \| u(y) - \tilde{u}(y) \|_{\R^m}
				\big)
			\ud t\\ \nonumber
		={} 
			&\nicefrac {1}{2}
			\|
					 b''_2
				\|_{
						\C(
							\Omega 
							\times (-R,R),
							\| \cdot \|_{L(\R^m,(\R^d \times \R)')}
						)
					} 
			\big \| (x-y, u(x)-u(y)) \big \|_{\R^d \times \R^m}
		\\  \nonumber & \cdot
			\big(\|u(x) - \tilde{u}(x)\|_{\R^m} + \|  u(y) - \tilde{u}(y) \|_{\R^m}\big).
\end{align}
Hence, 
\eqref{eq: b 2nd derivative estimate p= infty},
\eqref{eq: simple estm xx}
	(with $v_1 \leftarrow u-\tilde{u}$, $v_2 \leftarrow { \bf 1}$, 
	$\tilde{q}_1 \leftarrow 2$,
	$\tilde{q}_2 \leftarrow \infty$, 
	$t \leftarrow 1$),
\eqref{eq: 3 func xy estimate} 
	(with $p \leftarrow 1$, $v_1 \leftarrow u - \tilde{u}$,
	$v_2 \leftarrow {\bf 1}$, $v_3 \leftarrow u$, 
	$
			\tilde{q}_1 
		\leftarrow 
			q_1
	$,
	$
			\tilde{q}_2 
		\leftarrow 
			\infty 
	$),
the fact that 
for all $r\in \R$ and all $t\in (-\infty, r]$
there exists a $C \in (0,\infty)$ such that for all 
	$u \in B^t_{2,2}(\Omega,\R^m)$ it holds that
$\|u\|_{B^t_{2,2}(\Omega,\R^m)} \leq C \|u\|_{B^r_{2,2}(\Omega,\R^m)}$
and the fact that $\tilde{q}_3 \leq q_2$
assures that  for all
$\tilde{q}_3 \in [2,\infty]$
with $\nicefrac {1}{q_1}+ \nicefrac{1} {\tilde{q}_3} = \nicefrac 12$
there exist $C_1, C_2 \in (0,\infty)$
such that for all $R \in [0, \infty]$
and all $u, \tilde{u} \in L^0(\Omega,\R^m)$ with 
$\|u\|_{L^\infty(\Omega,\R^m)} \vee \|\tilde{u}\|_{L^\infty(\Omega,\R^m)} = R $ 
it holds that 
\begin{align}
\label{eq: first term p= infty}
			\int_\Omega &\int_\Omega \bigg | \int_0^1 \int_0^1
				\Big \langle 
					b''_2 \big (
						t x + (1-t) y, 
						\lambda t u(x) + \lambda (1-t) u(y)
						+ (1-\lambda) t \tilde{u}(x) + (1-\lambda) (1-t) \tilde{u}(y)
					\big) 
				\\ \nonumber & \qquad 
						\big (t (u(x) - \tilde{u}(x))+ (1-t) (u(y) - \tilde{u}(y)) \big),
					( x-y, u(x)-u(y))
				\Big \rangle_{(\R^d \times \R^m)', \R^d \times \R^m}
			\ud \lambda \ud t \bigg | ^2
		\\ \nonumber & \qquad \qquad
			 \cdot \tfrac{1}{\|x-y \|^{2s+d}_{\R^d}}\ud x \ud y\\ \nonumber
		\leq{} 
			&\nicefrac {1}{4}
			\|
				 b''_2
			\|^2_{
					\C(
						\Omega 
						\times B(R),
						\| \cdot \|_{L(\R^m,(\R^d \times \R)')}
					)
				}
		\\ \nonumber &
			\int_\Omega \int_\Omega
				\frac
					{
						\| (x-y, u(x)-u(y)) \|^2_{\R^d \times \R^m}
						\big( 
							\|u(x) - \tilde{u}(x)\|_{\R^m} 
							+ \|  u(y) - \tilde{u}(y) \|_{\R^m}  \big)^2
					}
					{
						\|x-y\|^{2s+d}_{\R^d}
					}
			\ud x \ud y \\ \nonumber
		\leq{} 
			&\nicefrac {1}{2}
			\|
				 b''_2
			\|^2_{
					\C(
						\Omega 
						\times B(R),
						\| \cdot \|_{L(\R^m,(\R^d \times \R^m)')} 
					)
				} \\ \nonumber
			&\bigg(
			\int_\Omega \int_\Omega
				\frac
					{
						\| x-y\|^2_{\R^d }
						\big( 
							\|u(x) - \tilde{u}(x)\|_{\R^m} ^2 
							+ \| u(y) - \tilde{u}(y) \|_{\R^m} ^2
						\big)
					}
					{
						\|x-y\|^{2s+d}_{\R^d}
					}
			\ud x \ud y \\ \nonumber \qquad
			&+
			\int_\Omega \int_\Omega
				\frac
					{
						\| u(x)-u(y)\|_{\R^m} ^2
						\big( 
							\|u(x) - \tilde{u}(x)\|_{\R^m} ^2 
							+ \|  u(y) - \tilde{u}(y) \|_{\R^m} ^2 
						\big)
					}
					{
						\|x-y\|^{2s+d}_{\R^d}
					}
			\ud x \ud y 
			\bigg)\\ \nonumber
		={} 
			&\|
				 b''_2
			\|^2_{
					\C(
						\Omega
						\times B(R),
						\| \cdot \|_{L(\R^m, (\R^d \times \R^m)')}
					)
				}
			\bigg(
			\int_\Omega \int_\Omega
				\frac
					{
						\| x-y\|^2_{\R^d }
						\|u(x) - \tilde{u}(x)\|_{\R^m} ^2
					}
					{
						\|x-y\|^{2s+d}_{\R^d}
					}
			\ud x \ud y  
		\\ \nonumber & \qquad
			+
			\int_\Omega \int_\Omega
				\frac
					{
						\| u(x)-u(y)\|_{\R^m} ^2
						\|u(x) - \tilde{u}(x)\|_{\R^m} ^2
					}
					{
						\|x-y\|^{2s+d}_{\R^d}
					}
			\ud x \ud y 
			\bigg)\\ \nonumber
		\leq{} 
			&C_1
			\|
				 b''_2
			\|^2_{
					\C(
						\Omega 
						\times B(R),
						\| \cdot \|_{L(\R^m, (\R^d \times \R^m)')} 
					)
				}
			\bigg(
			\|u - \tilde{u}\|^2_{L^2(\Omega,\R^m)}
		\\ \nonumber & \qquad
			+\|
				u-\tilde{u}
			\|^{2}_{B^{d/2-d/q_1+\eps \1_{\{\infty\}}(q_1)}_{2,2}(\Omega,\R^m)}
				\|u\|^2_{
					B^{s+d/2-d/\tilde{q}_3+\eps \1_{(2,\infty]}(\tilde{q}_3)}_{2,2}(\Omega,\R^m)
				}
			\bigg) \\ \nonumber
		\leq{} 
			&C_2 
			\|
				 b''_2
			\|^2_{
					\C(
						\Omega 
						\times B(R),
						\| \cdot \|_{L(\R^m,(\R^d \times \R^m)')} 
					)
				}
			\|
				u-\tilde{u}
			\|^{2}_{B^{d/2-d/q_1+\eps \1_{\{\infty\}}(q_1)}_{2,2}(\Omega,\R^m)} \\ \nonumber
			& \cdot
				\big(
					\|
						u
					\|^2_{B^{s+d/2-d/q_2+\eps \1_{(2,\infty]}(q_2)}_{2,2}(\Omega,\R^m)}
					+1
				\big) \\ \nonumber
		={} 
			&C_2
			\|
				u-\tilde{u}
			\|^{2}_{B^{d/2-d/q_1+\eps\1_{\{\infty\}}(q_1)}_{2,2}(\Omega,\R^m)}
				\big(
					\|
						u
					\|^2_{B^{s+d/2-d/q_2+\eps \1_{(2,\infty]}(q_2)}_{2,2}(\Omega,\R^m)}
					+1
				\big)
			\cdot K_1(\infty, u , \tilde{u}).
\end{align}
Combining the 3 cases 
\eqref{eq: first term p= 0}, 
\eqref{eq: first term 0 < p < infty}, 
and \eqref{eq: first term p= infty}
proves that there exist a $C \in (0, \infty)$
such that for all $u, \tilde{u} \in L^0(\Omega,\R^m)$
it holds that 
\begin{align}
\label{eq: first term estimate}
			\int_\Omega &\int_\Omega \bigg | \int_0^1 \int_0^1
				\Big \langle 
					b''_2 \big (
						t x + (1-t) y, 
						\lambda t u(x) + \lambda (1-t) u(y)
						+ (1-\lambda) t \tilde{u}(x) + (1-\lambda) (1-t) \tilde{u}(y)
					\big) 
				\\ \nonumber & \qquad
						\cdot \big (t (u(x) - \tilde{u}(x))+ (1-t) (u(y) - \tilde{u}(y)) \big),
					( x-y, u(x)-u(y))
				\Big \rangle_{(\R^d \times \R^m)', \R^d \times \R^m}
			\ud \lambda \ud t \bigg | ^2 
		\\  \nonumber & \qquad \qquad
			\cdot \tfrac{1}{\|x-y \|^{2s+d}_{\R^d}}\ud x \ud y\\ \nonumber
		\leq{} 
			&C
			\|
				u-\tilde{u}
			\|^{2}_{B^{d/2-d/q_1+\eps\1_{\{\infty\}}(q_1)}_{2,2}(\Omega,\R^m)}
				\big(
					\|
						u
					\|^2_{B^{s+d/2-d/q_2+\eps \1_{(2,\infty]}(q_2)}_{2,2}(\Omega,\R^m)}
					+1
				\big)
			\cdot K_1(p_1, u , \tilde{u})
\end{align}
and this finishes the estimation of the first term in \eqref{eq: b representation}.
Next we will estimate the 2nd term in \eqref{eq: b representation}.
Therefore we will again divide the proof into 3 cases.
\\
{\it 1. case $p_2= 0$.}
Note that assumption \eqref{eq: b estimate 1st derivative} 
verifies that for all $x, y \in \Omega$ and all
$u, \tilde{u} \in L^0(\Omega,\R^m)$ it holds that
\begin{align}
\nonumber
			&\int_0^1
				\big \langle 
					(D_{\R^d \times \R^m} b)
						\big( tx + (1-t) y, t \tilde{u}(x) + (1-t) \tilde{u}(y) \big),
				\\ \nonumber & \qquad
					\big( 0, u(x)-\tilde{u}(x)-u(y) +\tilde{u}(y) \big)
				\big \rangle_{(\R^d \times \R^m)', \R^d \times \R^m}
			\ud t \\
		\leq{} 
			&\int_0^1
				\|
					(D_{\R^d \times \R^m} b)
						\big( tx + (1-t) y, t \tilde{u}(x) + (1-t) \tilde{u}(y) \big)
				\|_{(\R^d \times \R^m)'}
			\\ \nonumber & \qquad \cdot 
				\|
					( 0, u(x)-\tilde{u}(x)-u(y) +\tilde{u}(y) )
				\|_{\R^d \times \R^m}
			\ud t \\ \nonumber 
		\leq{} 
			&\int_0^1
					K 
				\cdot 
				\|
					u(x)-\tilde{u}(x)-u(y) +\tilde{u}(y)
				\|_{\R^m} 
			\ud t 
		=
			K \cdot
			\|
				u(x)-\tilde{u}(x)- (u(y) -\tilde{u}(y))
			\|_{\R^m}.
\end{align}
Hence we derive 
from Proposition 2.28 in Mitrea \cite{Mitrea2013}
and from the fact that 
for all $r\in \R$ and all $t\in (-\infty, r]$
there exists a $C \in (0,\infty)$ such that for all 
	$u \in B^t_{2,2}(\Omega,\R^m)$ it holds that
$\|u\|_{B^t_{2,2}(\Omega,\R^m)} \leq C \|u\|_{B^r_{2,2}(\Omega,\R^m)}$
that
there exist $C_1, C_2 \in (0,\infty)$
such that for all
$u, \tilde{u} \in L^0(\Omega,\R^m)$ 
it holds that
\begin{align}
\nonumber
			&\int_\Omega \int_\Omega 
				\frac
					{1}
					{\|x-y \|^{2s+d}_{\R^d}}
						\Big | \int_0^1
							\Big \langle 
								(D_{\R^d \times \R^m} b)
									\big( tx + (1-t) y, t \tilde{u}(x) + (1-t) \tilde{u}(y) \big),
						\\ \nonumber & \qquad 
								\big( 0, u(x)-\tilde{u}(x)-u(y) +\tilde{u}(y) \big)
							\Big \rangle_{(\R^d \times \R^m)', \R^d \times \R^m}
						\ud t \Big | ^2
			\ud x \ud y\\ 
	\label{eq: second term p= 0}
		\leq{} 
			&K^2
			\int_\Omega \int_\Omega
				\frac
					{
						\|u(x)-\tilde{u}(x)- (u(y) -\tilde{u}(y)) \|_{\R^m} ^2
					}
					{
						\|x-y\|^{2s+d}_{\R^d}
					}
			\ud x \ud y \\ \nonumber
		\leq{}
			&C_1 \|u - \tilde{u} \|^2_{B^s_{2,2}(\Omega,\R^m)}
		\leq 
			C_2 
				\|
					u - \tilde{u}
				\|^2_{B^{s+d/2-d/q_4+\eps \1_{(2,\infty]}(q_4)}_{2,2}(\Omega,\R^m)} \\ \nonumber
		={} 
			&C_2
			\|
				u - \tilde{u} 
			\|^2_{B^{s+d/2-d/q_4+\eps \1_{(2,\infty]}(q_4)}_{2,2}(\Omega,\R^m)}
			\cdot K_2(0, \tilde{u}).
\end{align}
{\it 2. case $p_2 \in (0,\infty)$.}
Assumption \eqref{eq: b estimate 1st derivative} 
yields that for all $x, y \in \Omega$ and all
$u, \tilde{u} \in L^0(\Omega,\R^m)$ it holds that
	\begin{align}
	\nonumber
			&\int_0^1
				\big \langle 
					(D_{\R^d \times \R^m} b)
						\big( tx + (1-t) y, t \tilde{u}(x) + (1-t) \tilde{u}(y) \big),
			\\ \nonumber & \qquad
					\big( 0, u(x)-\tilde{u}(x)-u(y) +\tilde{u}(y) \big)
				\big \rangle_{(\R^d \times \R^m)', \R^d \times \R^m}
			\ud t \\ \nonumber
		\leq{} 
			&\int_0^1
				\|
					(D_{\R^d \times \R^m} b)
						\big( tx + (1-t) y, t \tilde{u}(x) + (1-t) \tilde{u}(y) \big)
				\|_{(\R^d \times \R^m)'}
		\\ & \label{eq: b' int of t} \qquad \cdot
				\|
					( 0, u(x)-\tilde{u}(x)-u(y) +\tilde{u}(y) )
				\|_{\R^d \times \R^m}
			\ud t \\ \nonumber
		\leq{} 
			&\int_0^1
					K (1+ \| t \tilde{u}(x) + (1-t) \tilde{u}(y) \|_{\R^m} )^{p_2}
				\cdot 
				\|
					u(x)-\tilde{u}(x)-u(y) +\tilde{u}(y)
				\|_{\R^m} 
			\ud t \\ \nonumber
		\leq{} 
			&\int_0^1
					3 ^{p_2} K 
						(
							1+ \| t \tilde{u}(x)\|_{\R^m} ^{p_2} 
							+ \|(1-t) \tilde{u}(y) \|_{\R^m} ^{p_2}
						)
			\ud t \,
			\|
				u(x)-\tilde{u}(x)-u(y) +\tilde{u}(y)
			\|_{\R^m}  \\ \nonumber
		={} 
			&\tfrac{3^{p_2}}{p_2+1} 
				K (p_2+1+ \|\tilde{u}(x)\|_{\R^m} ^{p_2} + \|\tilde{u}(y)\|_{\R^m} ^{p_2})
			\|
				u(x)-\tilde{u}(x)- (u(y) -\tilde{u}(y))
			\|_{\R^m}.
\end{align}
Moreover, we obtain 
from Proposition 2.28 in Mitrea \cite{Mitrea2013}
and from the fact that 
for all $r\in \R$ and all $t\in (-\infty, r]$
there exists a $C \in (0,\infty)$ such that for all 
	$u \in B^t_{2,2}(\Omega,\R^m)$ it holds that
$\|u\|_{B^t_{2,2}(\Omega,\R^m)} \leq C \|u\|_{B^r_{2,2}(\Omega,\R^m)}$
that there exist $C_1, C_2 \in (0,\infty)$ such that for all
$u, \tilde{u} \in L^0(\Omega,\R^m)$ it holds that 
\begin{equation}
\label{eq: 1st term of b' estimate 0<p<infty} 
	\begin{split}
			&\int_\Omega \int_\Omega
				\frac
					{
						\|
							u(x)-\tilde{u}(x)- (u(y) -\tilde{u}(y))
						\|_{\R^m} ^2
					}
					{
						\|x-y\|^{2s+d}_{\R^d}
					}
			\ud x \ud y 
		\leq
			C_1 \|u -\tilde{u}\|^2_{B^s_{2,2}(\Omega,\R^m)} \\
		\leq{}
			&C_2 
			\|
				u-\tilde{u}
			\|^2_{B^{s+d/2-d/q_4+\eps \1_{(2, \infty]}(q_4)}_{2,2}(\Omega,\R^m)}. 
	\end{split}
\end{equation}
In addition, \eqref{eq: 3 func xy estimate}
	(with $v_1 \leftarrow \tilde{u}$, $v_2 \leftarrow {\bf 1}$, 
	$v_3 \leftarrow u - \tilde{u}$,
	$p \leftarrow p_2$, $\tilde{q}_1 \leftarrow q_5$,
	$\tilde{q}_2 \leftarrow \infty$,
	$\tilde{q}_3 \leftarrow q_4$)
implies that 
there exists a $C \in (0,\infty)$ such that for all 
$u, \tilde{u} \in L^0(\Omega,\R^m)$
it holds that
\begin{equation}
\label{eq: 2nd term of b' estimate 0<p<infty} 
	\begin{split}
			&\int_\Omega \int_\Omega
				\frac
					{
						(\|\tilde{u}(x)\|_{\R^m} ^{2 p_2} + \|\tilde{u}(y) \|_{\R^m} ^{2 p_2})
						\|
							u(x)-\tilde{u}(x)- (u(y) -\tilde{u}(y))
						\|_{\R^m} ^2
					}
					{
						\|x-y\|^{2s+d}_{\R^d}
					}
			\ud x \ud y \\
		={} 
			&2\int_\Omega \int_\Omega
				\frac
					{
						\| \tilde{u}(x) \|_{\R^m} ^{2 p_2}
						\|
							u(x)-\tilde{u}(x)- (u(y) -\tilde{u}(y))
						\|_{\R^m} ^2
					}
					{
						\|x-y\|^{2s+d}_{\R^d}
					}
			\ud x \ud y \\
		\leq{} 
			&C\|
				\tilde{u}
			\|^{2p_2}_{
				B^{d/2-d/(p_2 q_5)+\eps \1_{\{\infty\}}(q_5)}_{2,2}(\Omega,\R^m)
			}
				\|
					u-\tilde{u}
				\|^2_{B^{s+d/2-d/q_4+\eps \1_{(2, \infty]}(q_4)}_{2,2}(\Omega,\R^m)}.
	\end{split}
\end{equation}
Combining \eqref{eq: b' int of t}, \eqref{eq: 1st term of b' estimate 0<p<infty}, and
\eqref{eq: 2nd term of b' estimate 0<p<infty}
thus ensures that
there exists a $C \in (0,\infty)$ such that 
for all $u, \tilde{u} \in L^0(\Omega,\R^m)$
it holds that
\begin{align}
\nonumber
			&\int_\Omega \int_\Omega 
				\frac
					{1}
					{\|x-y \|^{2s+d}_{\R^d}}
						\big | \int_0^1
							\big \langle 
								(D_{\R^d \times \R^m} b)
									\big( tx + (1-t) y, t \tilde{u}(x) + (1-t) \tilde{u}(y) \big),
						\\ 
	\label{eq: second term 0 < p < infty}
		\begin{split}
			& \qquad
								\big( 0, u(x)-\tilde{u}(x)-u(y) +\tilde{u}(y) \big)
							\big \rangle_{(\R^d \times \R^m)', \R^d \times \R^m}
						\ud t \big | ^2
			\ud x \ud y\\
		\leq{}  
			&C\|
				u-\tilde{u}
			\|^2_{B^{s+d/2-d/q_4+\eps \1_{(2,\infty]}(q_4)}_{2,2}(\Omega,\R^m)}
			\big(
				\|
					\tilde{u}
				\|^{2p_2}_{
					B^{d/2-d/(p_2 q_5)+\eps \1_{\{\infty\}}(q_5)}_{2,2}(\Omega,\R^m)
				}
				+1
			\big) 
	\end{split}\\ \nonumber
		={} 
			&C
			\|
				u - \tilde{u} 
			\|^2_{B^{s+d/2-d/q_4+\eps \1_{(2,\infty]}(q_4)}_{2,2}(\Omega,\R^m)}
			\cdot K_2(p_2, \tilde{u}).
\end{align}
{\it 3. case $p_2= \infty$.}
Note that we have for all $x, y \in \Omega$ 
and all
$u, \tilde{u} \in L^0(\Omega,\R^m)$ 
that
\begin{equation}
		\begin{split}
			&\int_0^1
				\big \langle 
					(D_{\R^d \times \R^m} b)
						\big( tx + (1-t) y, t \tilde{u}(x) + (1-t) \tilde{u}(y) \big),
		\\ & \qquad
					\big( 0, u(x)-\tilde{u}(x)-u(y) +\tilde{u}(y) \big)
				\big \rangle_{(\R^d \times \R^m)', \R^d \times \R^m}
			\ud t \\
		\leq{} 
			&\int_0^1
				\|
					(D_{\R^d \times \R^m} b)
						\big( tx + (1-t) y, t \tilde{u}(x) + (1-t) \tilde{u}(y) \big)
				\|_{(\R^d \times \R^m)'}
		\\ & \qquad \cdot
				\|
					( 0, u(x)-\tilde{u}(x)-u(y) +\tilde{u}(y) )
				\|_{\R^d \times \R^m}
			\ud t \\
		\leq{} 
			&\|
				D_{\R^d \times \R^m} b
			\|_{
					\C(
						\Omega
						\times B(\|\tilde{u} \|_{L^{\infty}(\Omega,\R^m)}),
						\| \cdot \|_{(\R^d \times \R^m)'}
					)
				}
			\int_0^1 
				\|
					u(x)-\tilde{u}(x)-u(y) +\tilde{u}(y)
				\|_{\R^m} 
			\ud t \\
		={} 
			&\|
				D_{\R^d \times \R^m} b
			\|_{
					\C(
						\Omega
						\times B(\|\tilde{u} \|_{L^{\infty}(\Omega,\R^m)}),
						\| \cdot \|_{(\R^d \times \R^m)'}
					)
				} \cdot
			\|
				u(x)-\tilde{u}(x)- (u(y) -\tilde{u}(y))
			\|_{\R^m}.
	\end{split}
\end{equation}
Therefore
Proposition 2.28 in Mitrea \cite{Mitrea2013}
and the fact that 
for all $r\in \R$ and all $t\in (-\infty, r]$
there exists a $C \in (0,\infty)$ such that for all 
	$u \in B^t_{2,2}(\Omega,\R^m)$ it holds that
$\|u\|_{B^t_{2,2}(\Omega,\R^m)} \leq C \|u\|_{B^r_{2,2}(\Omega,\R^m)}$
establish that
there exist $C_1, C_2 \in (0,\infty)$ such that
for all 
$u, \tilde{u} \in L^0(\Omega,\R^m)$ 
it holds that
\begin{equation}
\label{eq: second term p= infty}
	\begin{split}
			&\int_\Omega \int_\Omega 
				\frac
					{1}
					{\|x-y \|^{2s+d}_{\R^d}}
						\Big | \int_0^1
							\Big \langle 
								(D_{\R^d \times \R^m} b)
									\big( tx + (1-t) y, t \tilde{u}(x) + (1-t) \tilde{u}(y) \big),
					\\& \qquad
								\big( 0, u(x)-\tilde{u}(x)-u(y) +\tilde{u}(y) \big)
							\Big \rangle_{(\R^d \times \R^m)', \R^d \times \R^m}
						\ud t \Big | ^2
			\ud x \ud y\\
		\leq{} 
			&\|
				D_{\R^d \times \R^m} b
			\|^2_{
					\C(
						\Omega 
						\times B(\|\tilde{u} \|_{L^{\infty}(\Omega,\R^m)}),
						\| \cdot \|_{(\R^d \times \R^m)'}
					)
				}
			\int_\Omega \int_\Omega
				\frac
					{
						\|u(x)-\tilde{u}(x)- (u(y) -\tilde{u}(y))\|_{\R^m}^2
					}
					{
						\|x-y\|^{2s+d}_{\R^d}
					}
			\ud x \ud y \\
		\leq{} 
			&C_1\|
				D_{\R^d \times \R^m} b
			\|^2_{
					\C(
						\Omega 
						\times B(\|\tilde{u} \|_{L^{\infty}(\Omega,\R^m)}),
						\| \cdot \|_{(\R^d \times \R^m)'} 
					)
				}
				\|u - \tilde{u} \|^2_{B^s_{2,2}(\Omega,\R^m)} \\
		\leq{} 
			&C_2\|
				D_{\R^d \times \R^m} b
			\|^2_{
					\C(
						\Omega
						\times B(\|\tilde{u} \|_{L^{\infty}(\Omega,\R^m)}),
						\| \cdot \|_{(\R^d \times \R^m)'} 
					)
				}
				\|
					u - \tilde{u} 
				\|^2_{B^{s+d/2-d/q_4+\eps \1_{(2,\infty]}(q_4)}_{2,2}(\Omega,\R^m)} \\
		={} 
			&C_2
			\|
				u - \tilde{u} 
			\|^2_{B^{s+d/2-d/q_4+\eps \1_{(2,\infty]}(q_4)}_{2,2}(\Omega,\R^m)}
			\cdot K_2(\infty, \tilde{u}).
	\end{split}
\end{equation}
Combining the 3 cases \eqref{eq: second term p= 0},
\eqref{eq: second term 0 < p < infty}, and
\eqref{eq: second term p= infty}
shows that there exists a $C \in (0,\infty)$ such that
for all $u, \tilde{u} \in L^0(\Omega,\R^m)$
it holds that
\begin{equation}
\label{eq: second term estimate}
	\begin{split}
			&\int_\Omega \int_\Omega 
				\frac
					{1}
					{\|x-y \|^{2s+d}_{\R^d}}
						\Big | \int_0^1
							\Big \langle 
								(D_{\R^d \times \R^m} b)
									\big( tx + (1-t) y, t \tilde{u}(x) + (1-t) \tilde{u}(y) \big),
					\\ & \qquad
								\big( 0, u(x)-\tilde{u}(x)-u(y) +\tilde{u}(y) \big)
							\Big \rangle_{(\R^d \times \R^m)', \R^d \times \R^m}
						\ud t \Big | ^2
			\ud x \ud y\\
		\leq{} 
			&C
			\|
				u - \tilde{u} 
			\|^2_{B^{s+d/2-d/q_4+\eps \1_{(2,\infty]}(q_4)}_{2,2}(\Omega,\R^m)}
			\cdot K_2(p_2, \tilde{u}).
	\end{split}
\end{equation}
Thus \eqref{eq: b representation}, \eqref{eq: first term estimate}, 
and \eqref{eq: second term estimate} 
finishes the proof of Lemma \ref{lem: b continuity inequality}.
\end{proof}
The next lemma shows the Lipschitz continuity of Nemytskij operators
on bounded sets
in the $L^2$-norm.
\begin{lemma}[Lipschitz continuity wrt. the $L^2$-norm
 of Nemytskij opertors on bounded sets]
\label{lem: L2 b continuity inequality}
	Let $d, m \in \N$, 
	let $\Omega \subseteq \R^d$ be a bounded Lipschitz domain,
	denote by $B(r) \subseteq \R^m$, $r \in [0,\infty]$,
	the ball satisfying that for all $r \in [0,\infty]$ it holds
	that
	$
			B(r)
		= 
			\{ 
				y \in \R^m \colon 
					\| y\|_{\R^m} \leq r
			\}
	$,
	let 
	$\eps \in (0,\infty)$,
	$p \in [0,\infty]$,
	$q \in [2,\infty]$,
	satisfy that
	$p \, q \geq 2 \cdot \1_{(0,\infty)}(p)$,
	let
	$K \in (0,\infty)$,
	$b \in \C_{\R^d \times \R^m}^{1}(\Omega \times \R^m, \R)$,
	$
		K_1 \colon 
			[0,\infty] \times L^0(\Omega,\R^m) \times L^0(\Omega,\R^m) \to [0, \infty]
	$,
	satisfy that for all $r \in [0,\infty]$,
	$x \in \Omega$, $y \in \R^m$, and all $u_1, u_2 \in L^0(\Omega,\R^m)$ 
	it holds that
	\begin{equation}
			K_1(r,u_1, u_2)
		=
			\begin{cases}
				1 & \textrm{if } r =0 \\
				1+ \sum^2_{i=1}\|u_i\|^{2r}_{B^{d/2-d/(r q)+\eps \1_{\{\infty\}}(q)}_{2,2}(\Omega,\R^m)}
				& \textrm{if } r \in (0,\infty) \\
				\|
					D_{\R^d \times \R^m} b
				\|^2_{
					\C(
						\Omega
						\times 
						B(
							\|u\|_{L^\infty(\Omega,\R^m)} 
							\vee \|\tilde{u}\|_{L^\infty(\Omega,\R^m)}
						),
						\| \cdot \|_{(\R^d \times \R^m)'}
					)
				} 
				& \textrm{if } r =\infty
			\end{cases}
	\end{equation}
	and that
	\begin{equation}
	\label{eq: b estimate 1st derivative L2 lem} 
			\|(D_{\R^d \times \R^m} b)(x,y)\|_{(\R^d \times \R^m)'} 
		\leq 
			K \cdot (1+\|y\|_{\R^m})^{p}.
	\end{equation}
	Then there exists a
	$C \in (0,\infty)$ such that 
	for all $u, \tilde{u} \in L^0(\Omega,\R^m)$
	it holds that
	\begin{equation}
		\begin{split}
				&\int_{\Omega}
					|
						b( x ,u(x))-b (x,\tilde{u}(x)) 
					|^2
				\ud x
			\leq 
				C \, \|u - \tilde{u} \|^2_{B^{d/q+\eps \1_{\{2 \}}(q)}_{2,2}(\Omega,\R^m)}
				K_1(p,u, \tilde{u}).
		\end{split}
	\end{equation}
 \end{lemma}
\begin{proof}
	We will divide the proof into 3 cases. \\
	{\it 1. case $p= 0$.} We thus get from \eqref{eq: b estimate 1st derivative L2 lem} 
	and from the fact that
	for all $r\in \R$ and all $t\in (-\infty, r]$
	there exists a $C \in (0,\infty)$ such that for all 
	$u \in B^t_{2,2}(\Omega,\R^m)$ it holds that
	$\|u\|_{B^t_{2,2}(\Omega,\R^m)} \leq C \|u\|_{B^r_{2,2}(\Omega,\R^m)}$
	that there exists a $C \in (0, \infty)$ such that for all 
	$u, \tilde{u} \in L^0(\Omega,\R^m)$ it holds that
	\begin{align}
	\label{eq: b L2 diff estimate p=0}
				&\int_\Omega 
					|b(x,u(x))-b(x, \tilde{u}(x))|^2
				\ud x \\ \nonumber
			={} 
				&\int_\Omega 
					\Big |
						\int_0^1 
							\langle
								(D_{\R^d \times \R^m} b)
									(x,tu(x) + (1-t) \tilde{u}(x)),
								(0, u(x) - \tilde{u}(x))
							\rangle_{(\R^d \times \R^m)', \R^d \times \R^m}
						\ud t
					\Big |^2
				\ud x \\ \nonumber
			\leq{} 
				&\int_\Omega 
					\Big |
						\int_0^1 
							\|
								(D_{\R^d \times \R^m} b)
									(x,tu(x) + (1-t) \tilde{u}(x))
							\|_{(\R^d \times \R^m)'}
							\cdot \|u(x) - \tilde{u}(x)\|_{\R^m}
						\ud t
					\Big |^2
				\ud x \\ \nonumber
			\leq{} 
				&\int_\Omega 
						\Big |
							\int_0^1 
								K
								\cdot \|u(x) - \tilde{u}(x) \|_{\R^m}
							\ud t
						\Big |^2
					\ud x 
			=
				K^2 \|u - \tilde{u}\|_{L^2(\Omega,\R^m)}^2
			\leq
				C \, \|u - \tilde{u} \|^2_{B^{d/q+\eps \1_{\{2 \}}(q)}_{2,2}(\Omega,\R^m)} \\ \nonumber
			={} 
				&C \, \|u - \tilde{u} \|^2_{B^{d/q+\eps \1_{\{2 \}}(q)}_{2,2}(\Omega,\R^m)} 
				K_1(0,u,\tilde{u}).
	\end{align}
	{\it 2. case $p \in (0, \infty)$.}
	Hence, it follows from \eqref{eq: b estimate 1st derivative L2 lem}
	that for all $u, \tilde{u} \in L^0(\Omega,\R^m)$ and all 
	$x \in \Omega$ it holds that
	\begin{align}
	\nonumber
					&|b(x,u(x))-b(x, \tilde{u}(x))|^2 \\
			={} \nonumber
					&\Big |
						\int_0^1 
							\langle
								(D_{\R^d \times \R^m} b)
									(x,tu(x) + (1-t) \tilde{u}(x)),
								(0, u(x) - \tilde{u}(x))
							\rangle_{(\R^d \times \R^m)', \R^d \times \R^m}
						\ud t
					\Big |^2 \\ 
			\leq{}  
					&\Big |
						\int_0^1 
							\|
								(D_{\R^d \times \R^m} b)
									(x,tu(x) + (1-t) \tilde{u}(x))
							\|_{(\R^d \times \R^m)'}
							\cdot \|u(x) - \tilde{u}(x)\|_{\R^m}
						\ud t
					\Big |^2 \\  \nonumber
			\leq{} 
					&\Big |
						\int_0^1 
							(1+\|tu(x) + (1-t) \tilde{u}(x)\|_{\R^m})^{p}
							\cdot \|u(x) - \tilde{u}(x)\|_{\R^m}
						\ud t
					\Big |^2 
		\\ \nonumber
			\leq{} 
					&
					\Big |
						3^p
						\int_0^1 
							1+t^p \|u(x)\|_{\R^m}^p + (1-t)^p \|\tilde{u}(x)\|_{\R^m}^p
						\ud t
						\cdot \|u(x) - \tilde{u}(x)\|_{\R^m}
					\Big |^2 \\ 
			={} \nonumber
					&9^p
						(
							1
							+ \tfrac{\|u(x)\|_{\R^m}^p}{p+1} 
							+ \tfrac{\|\tilde{u}(x)\|_{\R^m}^p}{p+1}
						)^2
						\cdot \|u(x) - \tilde{u}(x)\|^2_{\R^m}.
	\end{align}
	This together with
	the Sobolev inequality 
		(see, e.g., the Theorem on page 31
		and Theorem 1 on page 32 in Runst \& Sickel \cite{RunstSickel1996})
	and with the fact that
	for all $r\in \R$ and all $t\in (-\infty, r]$
	there exists a $C \in (0,\infty)$ such that for all 
	$u \in B^t_{2,2}(\Omega,\R^m)$ it holds that
	$\|u\|_{B^t_{2,2}(\Omega,\R^m)} \leq C \|u\|_{B^r_{2,2}(\Omega,\R^m)}$
	implies that there exist $C_1, C_2, C_3, C_4 \in (0, \infty)$ 
	such that for all $u, \tilde{u} \in L^0(\Omega,\R^m)$
	and all $q' \in [2, \infty]$ with $\tfrac {1}{q'} + \tfrac 1q = \tfrac 12$
	it holds that
	\begin{align}
	\nonumber
				&\int_\Omega 
					|b(x,u(x))-b(x, \tilde{u}(x))|^2 
				\ud x \\ \nonumber
			\leq{} 
				&\int_\Omega 
					9^p 
						(
							1
							+ \tfrac{\|u(x)\|_{\R^m}^p}{p+1} 
							+ \tfrac{\|\tilde{u}(x)\|_{\R^m}^p}{p+1}
						)^2
						\cdot \|u(x) - \tilde{u}(x)\|^2_{\R^m}
				\ud x \\ \nonumber
			\leq{} 
				&\int_\Omega 
					9^{p+1}
					(
						1
						+\tfrac{\|u(x)\|_{\R^m}^{2p}}{(p+1)^2} 
						+ \tfrac{\|\tilde{u}(x)\|_{\R^m}^{2p}}{(p+1)^2})
					\cdot \|u(x) - \tilde{u}(x)\|_{\R^m}^2
				\ud x \\ \nonumber
			\leq{} 
				&C_1 \big(
					\|u - \tilde{u}\|_{L^2(\Omega,\R^m)}
					+
					\big(
						\big \| \|u\|_{\R^m}^{p} \big \|^2_{L^{q}(\Omega)} 
						+ \big \| \|\tilde{u}\|_{\R^m}^{p} \big \|^2_{L^{q}(\Omega)}
					\big)
						\cdot \big \| \|u - \tilde{u}\|_{\R^m} \big\|^2_{L^{q'}(\Omega)}
				\big) \\ \nonumber
			\leq{} 
				&C_2\big(
					\|u - \tilde{u}\|_{L^2(\Omega,\R^m)}
					+
					(
						\|u \|^{2p}_{L^{p q}(\Omega,\R^m)}
						+ \|\tilde{u} \|^{2p}_{L^{p q}(\Omega,\R^m)}
					)
						\cdot \| u - \tilde{u}\|^2_{L^{q'}(\Omega,\R^m)} 
				\big) \\
	\label{eq: b L2 diff estimate p in (0, infty)}
			\leq{} 
				&C_3 \big(
					\|u - \tilde{u}\|_{L^2(\Omega,\R^m)}
					+(
						\|u \|^{2p}_{B_{2,2}^{d/2 - d/(p q)+\eps \1_{ \{\infty \}}(q)}(\Omega,\R^m)}
						+\|\tilde{u} \|^{2p}_{
							B_{2,2}^{d/2 - d/(p q)+\eps \1_{ \{\infty \}}(q)}(\Omega,\R^m)
						}
					)
				\\ \nonumber & \qquad
						\cdot 
							\| 
								u - \tilde{u}
							\|^2_{B_{2,2}^{d/2 - d/q'+\eps \1_{ \{\infty \}}(q')}(\Omega,\R^m)} 
				\big)\\ \nonumber
			\leq{} 
				&C_4 
				(
					\|u \|^{2p}_{B_{2,2}^{d/2 - d/(p q)+\eps \1_{ \{\infty \}}(q)}(\Omega,\R^m)} 
					+\|
						\tilde{u} 
					\|^{2p}_{B_{2,2}^{d/2 - d/(p q)+\eps \1_{ \{\infty \}}(q)}(\Omega,\R^m)}
					+1
				)
			\\ \nonumber & \qquad
					\cdot 
						\| 
							u - \tilde{u}
						\|^2_{B_{2,2}^{d/q+\eps \1_{ \{2 \}}(q)}(\Omega,\R^m)} \\ \nonumber
			={} 
				&C_4 \,
					\| 
						u - \tilde{u}
					\|^2_{B_{2,2}^{d/q+\eps \1_{ \{2 \}}(q)}(\Omega,\R^m)}
					K_1(p,u,\tilde{u}).
	\end{align}
	{\it 3. case $p = \infty$.}
	Note that the fact that
	for all $r\in \R$ and all $t\in (-\infty, r]$
	there exists a $C \in (0,\infty)$ such that for all 
	$u \in B^t_{2,2}(\Omega,\R^m)$ it holds that
	$\|u\|_{B^t_{2,2}(\Omega,\R^m)} \leq C \|u\|_{B^r_{2,2}(\Omega,\R^m)}$
	establishes
	that there exists a $C \in (0, \infty)$ such that for all 
	$u, \tilde{u} \in L^0(\Omega,\R^m)$ it holds that
	\begin{equation}
	\label{eq: b L2 diff estimate p=infty}
		\begin{split}
				&\int_\Omega 
					|b(x,u(x))-b(x, \tilde{u}(x))|^2
				\ud x \\
			={} 
				&\int_\Omega 
					\Big |
						\int_0^1 
							\langle
								(D_{\R^d \times \R^m} b)
									(x,tu(x) + (1-t) \tilde{u}(x)),
								(0, u(x) - \tilde{u}(x))
							\rangle_{(\R^d \times \R^m)', \R^d \times \R^m}
						\ud t
					\Big |^2
				\ud x \\
			\leq{} 
				&\int_\Omega 
					\Big |
						\int_0^1 
							\|
								(D_{\R^d \times \R^m} b)
									(x,tu(x) + (1-t) \tilde{u}(x))
							\|_{(\R^d \times \R^m)'}
							\cdot \|u(x) - \tilde{u}(x)\|_{\R^m}
						\ud t
					\Big |^2
				\ud x \\
			\leq{} 
				&\int_\Omega 
						\Big |
							\int_0^1 
								\|
									D_{\R^d \times \R^m} b
								\|_{
									\C(
										\Omega \times 
											B(
												\|u\|_{L^\infty(\Omega,\R^m)} 
												\vee \|\tilde{u}\|_{L^\infty(\Omega,\R^m)}
											),
										\| \cdot \|_{(\R^d \times \R^m)'}
									)
								}
								\cdot \|u(x) - \tilde{u}(x)\|_{\R^m}
							\ud t
						\Big |^2
					\ud x \\
			={} 
				&\|u - \tilde{u}\|_{L^2(\Omega,\R^m)}^2 K_1(\infty, u, \tilde{u})
			\leq
				C \, 
				\|u - \tilde{u} \|^2_{B^{d/q+\eps \1_{\{2 \}}(q)}_{2,2}(\Omega,\R^m)}
				K_1(\infty, u, \tilde{u}). 
		\end{split}
	\end{equation}
	Combining the 3 cases \eqref{eq: b L2 diff estimate p=0},
	\eqref{eq: b L2 diff estimate p in (0, infty)},
	and
	\eqref{eq: b L2 diff estimate p=infty}
	completes the proof of
	Lemma \ref{lem: L2 b continuity inequality}.
\end{proof}
The next corollary is a direct consequence of Lemma \ref{lem: b continuity inequality}
and Lemma \ref{lem: L2 b continuity inequality}.
\begin{corollary}[Lipschitz continuity wrt. Besov-norms
 of Nemytskij operators on bounded sets]
	\label{cor: Bs b continuity inequality}
	Let $d, m \in \N$, 
	let $\Omega \subseteq \R^d$ be a bounded Lipschitz domain,
	denote by $B(r) \subseteq \R^m$, $r \in [0,\infty]$,
	the ball satisfying that for all $r \in [0,\infty]$ it holds
	that
	$
			B(r)
		= 
			\{ 
				y \in \R^m \colon 
					\| y\|_{\R^m} \leq r
			\}
	$,
	let 
	$s \in (0,1)$,
	$\eps \in (0,1-s)$,
	$p_1,p_2 \in [0,\infty]$,
	$q, q_1,q_2,q_3,q_4,q_5 \in [2,\infty]$,
	satisfy that
	$
			\nicefrac {1}{q_1} + \nicefrac{1}{q_2} + \nicefrac{1}{q_3}
		=
			\nicefrac {1}{q_4} + \nicefrac{1}{q_5}
		=
			\nicefrac {1}{2}
	$, that
	$p_1 \, q_3 \geq 2 \cdot \1_{(0,\infty)}(p_1)$,
	that $p_2 \, q_5 \geq 2 \cdot \1_{(0,\infty)}(p_2)$,
	and that $p_2 \, q \geq 2 \cdot \1_{(0,\infty)}(p_2)$,
	let
	$K \in (0,\infty)$,
	$b \in \C_{\R^d \times \R^m}^{1,2}(\Omega \times \R^m, \R)$,
	$b''_{2} \in \C(\Omega \times \R^m, L(\R^m,(\R^d \times \R^m)'))$,
	$
		K_1, K_3 \colon 
			[0,\infty] \times L^0(\Omega,\R^m) \times L^0(\Omega,\R^m) \to [0, \infty]
	$,
	$K_2 \colon [0,\infty] \times L^0(\Omega,\R^m) \to [0, \infty]$,
	satisfy that for all $p \in [0,\infty]$,
	$x \in \Omega$, $y \in \R^m$, $u_1 \in L^0(\Omega,\R^m)$ 
	and all $u_2 \in L^0(\Omega,\R^m)$ it holds that
	\begin{equation}
		\begin{split}
				b''_{2}(x,y)
			=
				\Big(
					D_{\R^m}\big( 
						\R^m \ni z \to D_{\R^d \times \R^m} b(x,z) \in (\R^d \times \R^m)'
					\big)
				\Big)(y),
		\end{split}
	\end{equation}
	that
	\begin{equation}
			K_1(p,u_1,u_2)
		=
			\begin{cases}
				1 & \textrm{ if } p =0 \\
				1+ \sum_{i=1}^2
					\|u_i\|^{2p}_{B^{d/2-d/(p q_3)+\eps \1_{\{\infty\}}(q_3)}_{2,2}(\Omega,\R^m)}
				& \textrm{if } p \in (0,\infty) \\
				\|
					b''_2
				\|^2_{
					\C(
						\Omega 
						\times 
							B(
								\|u_1\|_{L^\infty(\Omega,\R^m)} 
								\vee \|u_2\|_{L^\infty(\Omega,\R^m)}),
						\| \cdot \|_{L(\R^m, (\R^d \times \R^m)')}
					)
				} & \textrm{if } p =\infty
			\end{cases},
	\end{equation}
	that
	\begin{equation}
			K_2(p,u_1)
		=
			\begin{cases}
				1 & \textrm{ if } p =0 \\
				\|u_1\|^{2p}_{B^{d/2-d/(p q_5)+\eps \1_{\{\infty\}}(q_5)}_{2,2}(\Omega,\R^m)}
				+1 & \textrm{ if } p \in (0,\infty) \\
				\|
					D_{\R^d \times \R^m} b
				\|^2_{
					\C(
						\Omega
						\times
							B(\|u_1\|_{L^\infty(\Omega,\R^m)}),
						\| \cdot \|_{(\R^d \times \R^m)'} 
					)
				} & \textrm{ if } p =\infty
			\end{cases},
	\end{equation}
	that
	\begin{equation}
			K_3(p,u_1, u_2)
		=
			\begin{cases}
				1 & \textrm{if } p =0 \\
				1+ \sum_{i=1}^2
					\|u_i\|^{2p}_{B^{d/2-d/(p q)+\eps \1_{\{\infty\}}(q)}_{2,2}(\Omega,\R^m)}
				& \textrm{if } p \in (0,\infty) \\
				\|
					D_{\R^d \times \R^m} b
				\|^2_{
					\C(
						\Omega
						\times
							B(\|u_1\|_{L^\infty(\Omega,\R^m)} \vee \|u_2\|_{L^\infty(\Omega,\R^m)}),
						\| \cdot \|_{(\R^d \times \R^m)'} 
					)
				} 
				& \textrm{if } p =\infty
			\end{cases},
	\end{equation}
	that
	\begin{equation}
			\| (D_{\R^d \times \R^m} b) (x,y)\|_{(\R^d \times \R^m)'}
		\leq 
			K \cdot (1+\|y\|_{\R^m})^{p_2},
	\end{equation}
	and that
	\begin{equation}
			\|
				b''_2(x,y)
			\|_{L(\R^m,(\R^d \times \R^m)')}
		\leq
			K \cdot (1+\|y\|_{\R^m})^{p_1}.
	\end{equation}
	Then there exists a
	$C \in (0,\infty)$ such that 
	for all $u, \tilde{u} \in L^0(\Omega,\R^m)$
	it holds that
	\begin{equation}
		\begin{split}
				&\| 
					b(\cdot , u(\cdot)) - b(\cdot, \tilde{u}(\cdot))
				\|^2_{B^s_{2,2}(\Omega,\R)}
			\\ & 
			\leq 
				C \, \|
					u-\tilde{u}
				\|^{2}_{B^{d/2-d/q_1+\eps \1_{\{\infty\}}(q_1)}_{2,2}(\Omega,\R^m)}
				\big(
					\|u\|^2_{B^{s+d/2-d/q_2+\eps \1_{(2,\infty]}(q_2)}_{2,2}(\Omega,\R^m)}
					+1
				\big) \cdot K_1(p_1,u,\tilde{u}) \\
			& 
				+C \, \|
					u - \tilde{u} 
				\|^2_{B^{s+d/2-d/q_4+\eps \1_{(2,\infty]}(q_4)}_{2,2}(\Omega,\R^m)}
				\cdot K_2(p_2, \tilde{u}) 
			 \\&
				+C \, \|u - \tilde{u} \|^2_{B^{d/q+\eps \1_{\{2 \}}(q)}_{2,2}(\Omega,\R^m)}
				\cdot K_3(p_2,u, \tilde{u}).
		\end{split}
	\end{equation}
\end{corollary}
\begin{proof}
	The proof follows immediately from Lemma \ref{lem: b continuity inequality},
	Lemma \ref{lem: L2 b continuity inequality},
	and from Proposition 2.28 in Mitrea \cite{Mitrea2013}.
\end{proof}
\chapter{Application}
\label{sec: Application}
In this chapter we apply the results of Chapter 3 and Chapter 4 
to stochastic Burgers equations and stochastic 2-D Navier-Stokes equations.
In Section \ref{sec: Burgers' equation} we apply 
Corollary \ref{cor:uniqueness2}
and Theorem \ref{thm: existence}
to Kolmogorov backward equations of stochastic Burgers equations
and in Section \ref{sec: 2-D Navier-Stokes equations}
to Kolmogorov backward equations of stochastic 2-D Navier-Stokes equations.
\section{Stochastic Burgers equations}
\label{sec: Burgers' equation}
In this section we prove that Kolmogorov equations
of stochastic Burgers equations
have a unique viscosity solution having at most polynomial growth.
\begin{sett}
\label{sett: Burger's}
 Let 
	$
			\mathbb{H}
		= (H, \langle \cdot, \cdot \rangle_{H},  \|\cdot\|_{H} ) 
		= \mathbb{L}^2(0,1)
	$
	be the $L^2$ space,
	let $\Delta \colon L^2(0,1) \supseteq D(\Delta) \to L^2(0,1)$
	be the Laplace operator
	with Dirichlet boundary condition
	i.e.\@ the operator satisfying that
	$
			D(\Delta) 
		= 
			\{ 
				u \in B^{2}_{2,2}(0,1) \colon 
					u(0)=u(1)=0
			\}
	$,
	and that for all $u \in D(\Delta)$
	it holds that
	$\Delta(u) = u''$,
	let 
		$ 
				\mathbb{H}_{t}
			=
				( 
					H_{t} , 
					\left< \cdot , \cdot \right>_{ H_{t} }, 
					\left\| \cdot \right\|_{ H_{t} } 
				) 
		$,
		$ t \in \R $,
		be a family of interpolation spaces associated with
		$ -\Delta$ 
		(see, e.g., Definition~3.6.30 in Jentzen \cite{Jentzen2015}).
		By abuse of notation we will also denote by
		$\Delta$ and  by $\| \cdot \|_{H_t}$, $t \in \R$,
		the extended operators 
		$
			\Delta \colon \bigcup_{i=1}^{\infty} H_{-i} \to \bigcup_{i=1}^{\infty} H_{-i}
		$ 
		and 
		$
			\| \cdot \|_{H_t} \colon \bigcup_{i=1}^{\infty} H_{-i} \to [0, \infty]
		$,
		$t \in \R$,
		satisfying for all
		$t \in \R$, $x \in H_{t}$, 
		and all $y\in \bigcup_{i=1}^{\infty} H_{-i}$
		that
		\begin{equation}
					\|y\|_{H_t}
				=
					\begin{cases} 
						\|y\|_{H_{t}} 
							& \textrm{ if } y \in H_{t} \\
						\infty 
							& \textrm{ if } y \notin H_{t}
					\end{cases}
		\end{equation}
		and that
		\begin{equation}
				(\Delta(x) = y) 
			\Leftrightarrow 
				(
						\lim_{\eps \downarrow 0} \sup \{
							\|\Delta(\xi) - y\|_{H_{t-1}} \colon \xi \in H_1, ~\|x- \xi\|_{H_{t}} \leq \eps 
						\}
					=
						0
				).
		\end{equation}
	Let
	$e_k \in L^2(0,1)$, $k\in \N$, be the eigenfunction of $\Delta$ satisfying
	for all $x \in (0,1)$ and all $k\in \N$ that
	$e_k(x) = \sqrt{2} \sin(k \pi x)$,
	let $\lambda_k \in (-\infty, 0)$, $k \in \N$, be the eigenvalues of $e_k$
	i.e.\@ the real number satisfying for all $k \in \N$ that
	$\lambda_k = -\pi^2 k^2$,
	let
	$\alpha \in (\nicefrac 34, 1)$,
	$\beta \in [0,\nicefrac 12)$,
	$
		c \in \R \backslash \{ 0 \}
	$
	let $\mathbb{U}=(U, \langle \cdot, \cdot \rangle_{U}, \| \cdot \|_U)$
	be a real separable Hilbert space,
	let
	$
		B \in 
		\C_{
			\mathbb{H}, 
			\mathbb{HS}(\mathbb{U}, \mathbb{H}_{-\beta})
		}( 
			H, 
			HS( \mathbb{U}, \mathbb{H}_{-\beta} ) 
		),
		$
	and denote by 
	$
		F \colon 
			H \to
			H_{-\alpha}
	$
	the function satisfying for all $u \in H$
	that
	$F(u) = -\tfrac{c}{2} (u^2)'$.
\end{sett}
\begin{lemma}
\label{l: norm equivalence Burger}
	Assume the 
	Setting \ref{sett: Burger's} and let $r \in [0,\infty)$.
	Then there exist $C_1, C_2 \in (0,\infty)$
	such that for all $u \in H_r$ it holds that
	\begin{equation}
		C_1 \|u\|_{B^{2r}_{2,2}(0,1)} \leq \|u\|_{H_r} \leq C_2 \|u\|_{B^{2r}_{2,2}(0,1)}.
	\end{equation}
\end{lemma}
\begin{proof}
	Note that we obtain from (2.200) on page 54 in Mitrea \cite{Mitrea2013}
	and from
	Theorem 1.8 in Nečas \cite{Necas2012} that for all $n \in \N_0$
	there exist $K_1,K_2,K_3,K_4 \in (0,\infty)$ such that
	for all $u \in H_{n}$
	it holds that
	\begin{equation}
		\begin{split}
				&\|u\|_{B^{2n}_{2,2}(0,1)}
			\leq
				K_1(
					\|u\|_{L^2(0,1)} 
					+ \sum^{2n}_{i=1} \|D_{\R}^{i} u\|_{L^2(0,1)}
				)
			\leq
				K_2 (
					\|u\|_{L^2(0,1)} + \|D_{\R}^{2n} u \|_{L^2(0,1)}
				) \\
			={} &
				K_2 (\|u\|_{H} + \|\Delta^n u \|_{H})
			\leq
				K_2 ( |\lambda_1|^{-n} \|u\|_{H_n} + \|u \|_{H_n})
			\leq
				K_3 \|u\|_{H_n} 
			= 
				K_3 \| D_{\R}^{2n} u \|_{L^2(0,1)} \\
			\leq{} &
				K_3(
					\|u\|_{L^2(0,1)} 
					+ \sum^{2n}_{i=1} \|D_{\R}^{i} u\|_{L^2(0,1)}
				)
			\leq
				K_4 \|u\|_{B^{2n}_{2,2}(0,1)}.
		\end{split}
	\end{equation}
	Combining this and Lemma \ref{l: norm equivalence with interpolation}
	finishes the proof of Lemma \ref{l: norm equivalence Burger}.
\end{proof}
The next corollary establishes our main result of this section. It follow from
Corollary \ref{cor:uniqueness2}
and Theorem \ref{thm: existence}.
\begin{corollary}
\label{cor: Burgers equation}
	Assume 
	the Setting \ref{sett: Burger's},
	let
	$\vartheta \in [\nicefrac 12, \infty)$,
	$T , \theta, \beta_1, \beta_2 \in (0,\infty)$,
	$\beta_3 \in (0, 1)$,
	$ \chi \in [ \beta, \nicefrac{ 1 }{ 2 } ) $,
	let $\mathscr{P}_0 \in L(\mathbb{U},\mathbb{U})$ be the function satisfying that
	$\mathscr{P}_0=\id_U$,
	let $\mathscr{P}_N \in L(\mathbb{U},\mathbb{U})$, $N \in \N$, 
	be finite-dimensional projection,
	let $\varphi \in \C_{\mathbb{H}}(H, \R)$ have at most polynomial $\H$-growth, 
	let $K \colon \R \to (0, \infty)$ be increasing,
	and assume that
	\begin{equation}
	\label{eq: B continuity Burger}
		B|_{H_{1/2 +\vartheta - \beta_3}} \in
				\C_{
					\mathbb{H}_{1/2 +\vartheta - \beta_3},
					\mathbb{HS}(\mathbb{U}, \mathbb{H}_{\vartheta})
				}(
					H_{1/2 +\vartheta - \beta_3},
					HS(\mathbb{U}, \mathbb{H}_{\vartheta})
				), 
	\end{equation}
	that for all $u, v \in H_{2\vartheta}$ 
	it holds that
	\begin{align}
		\label{eq: H bound B Burger}
			&\| B(u) \|^2_{HS(\mathbb{U},\mathbb{H})}
		\leq
			\theta \cdot ( 1 +\|u\|^2_{H}), \\
	\label{eq: B bounded wrt H_vartheta Burger}
			&\| B(u) \|^2_{HS(\mathbb{U},\mathbb{H}_{\vartheta})} 
		\leq
			K(\|u\|_{H}) 
				(\|u\|^2_{H_{1/2+\vartheta- \beta_1}} + 1), \\
	\label{eq: B Lipschitz wrt H_1/2 Burger}
			&\|B(u) - B(v)\|_{HS(\mathbb{U}, \mathbb{H})}
		\leq
			K(\|v\|_{H} \vee \|u \|_H) 
				\|u-v\|_{H_{1/2 -\beta_2}}, \\
	\label{eq: Lipschitz bound B Burger}
		&| B |_{
			\C^1( 
				H, 
				\left\| \cdot \right\|_{HS( \mathbb{U}, \mathbb{H}_{- \beta} )} 
			)
		}
		< \infty,
	\end{align}
	and that for all $\mathbb{H}$-bounded sets
	$ E \subseteq H$
	it holds that
	\begin{equation}
	\label{eq: B C1 bound Burger} 
			\sup_{ N \in \N } \sup_{ v \in E }
			\left[
				\frac{
				\|
					B( v ) \mathscr{P}_N 
				\|_{ HS( \mathbb{U}, \mathbb{H}_{- \beta } ) }
				}{
					1 + \| v \|_{ H }
				}
			\right]
			< \infty
	\end{equation}
	and that
	\begin{equation}
	\label{eq:convLocLip_assumption Burger}
			\limsup_{ N \to \infty }
			\sup_{ v \in E }
			\left[
			\frac{
			\|
				B( v ) ( \id_U - \mathscr{P}_N )
			\|_{ HS( \mathbb{U}, \mathbb{H}_{- \chi} ) }
			}{
				1 + \| v \|_{ H }
			}
			\right]
		=0,
	\end{equation}
		let $V_N \subseteq L^2(0,1)$, $N \in \N_0$, be the linear subspaces 
		satisfying that
		$V_0= H$ and that for all $N \in \N$ it holds that
		$V_N = \Span_{\mathbb{H}}(\{e_k \colon k \in \{1, \ldots, N \}\})$, let
		$ ( \Omega, \mathbb{F}, \P ) $
		be a probability space with a normal filtration 
		$ ( \mathbb{F}_t )_{ t \in [0,T] } $, 
		let 
		$
			( W_t )_{ t \in [0,T] } 
		$ 
		be an $ \operatorname{Id}_U $-cylindrical $ ( \mathbb{F}_t )_{ t \in [0,T] } $-Wiener
		process and let 
	$
		X^{N,u} \colon [0,T] \times \Omega \rightarrow H
	$,
	$ N \in \N_0 $,
	$u \in H$,
	be
	$
		(\mathbb{F}_t )_{ t \in [0,T] }
	$-adapted stochastic processes 
	with continuous sample paths satisfying
	for all $ N \in \N_0 $, $ t \in [0,T] $, and all $u \in H$
	that
	\begin{equation}
			X_t^{N,u}
		= 
			e^{ t A } \pi^{H}_{V_N} u
			+
			\int_0^t
				e^{ ( t - s ) A }
				\pi^{H_{-\alpha}}_{V_N} F( X_t^{N,u}  )
			\ds
			+
			\int_0^t
				e^{ ( t - s ) A }
				\pi^{H_{-\beta}}_{V_N} B( X_t^{N,u}  ) \mathscr{P}_N
			\dWs,
	\end{equation}
	and let $f \colon [0,T] \times H \to \R$
	satisfy that for all $(t,u) \in [0,T] \times H$
	it holds that
	$
		f(t,u) = \E[\varphi(X^{0,u}_t)].
	$
	Then
	$ 
    f|_{ (0,T) \times H} 
  $
  is the unique 
	viscosity solution of
  \begin{equation} 
	\label{eq: Kol equation Burger}
		\begin{split}
			&\tfrac{ \partial }{ \partial t }
			f(t,u) -
			\big\langle
				u''- \tfrac c2 (u^2)', 
				I^{-1}_{\mathbb{H}} (D_{\mathbb{H}} f)(t,u)
			\big\rangle_{H} \\
			&-
			\big \langle 
				B(u),
				I_{\mathbb{H}_{\vartheta}}^{-1}
					\big( (D^2_{\mathbb{H}} \, f)(t,u) \, B(u) \big)
			\big \rangle_{HS(\mathbb{U},\mathbb{H}_{\vartheta})} 
			= 0
		\end{split}
  \end{equation}
	for $ (t,u) \in (0,T) \times H $ relative to 
	$
		(
			(0,T) \times H \ni (t,u) 
			\to \nicefrac 12 \|u \|^2_{\mathbb{H}_{\vartheta}} \in [0, \infty],
			\R \times \mathbb{H}, 
			\R \times \mathbb{H}_{\vartheta}
		)
	$
	which satisfy that 
	$ 
    f \in \C_{\R \times \mathbb{H}}([0,T] \times H, \R)
  $,
	that $f$ is
	bounded on $\R \times \mathbb{H}$-bounded 
	subsets of $[0,T] \times H$, 
	that $f$ have at most polynomial $\H$-growth,
  and that for all $ u \in H $ it holds that
  $
    f(0, u) = \varphi(u)
  $.
	\end{corollary}
	\begin{proof}
	First note that the existence of $X^{N,u}$, $N \in \N_0$, $u \in H$, 
		follows from Theorem 1.1 and Remark 3.1. in Liu \& R\"ockner \cite{LiuRockner2010}.
		Moreover, 
		the fact that $\varphi$ have at most polynomial $\H$-growth verifies that
		there exists a $p \in [2,\infty)$, which we fix for the rest of the proof,
		such that 
		$
			\lim_{r \to \infty}
				\sup \big\{ \tfrac{\varphi(u)}{\|u\|_H^p} \colon u \in H, \|u\|_H \geq r \big\}
			=0.
		$
		Denote by 
		$ V \colon H \to (1,\infty)$ the function satisfying for all
		$u \in H$ that
		$V(u) = \|u\|^p_{L^2(0,1)} + 2$.
		We will now apply Theorem \ref{thm: existence} to show that
		$f$ is a viscosity solution.
	Therefore observe that it follows from the Sobolev inequality
	(see, e.g., Theorem 1 on page 32 in Runst \& Sickel \cite{RunstSickel1996}),
	from the fact that $2\alpha-1> \tfrac 12$, and from
	Lemma \ref{l: norm equivalence Burger}
	that there exist $C_1,C_2,C_3 \in (0,\infty)$ such that for all
	$w\in H_\alpha$ it holds that
	\begin{equation}
			\|w'\|_{L^\infty(0,1)}
		\leq
			C_1 \|w'\|_{B^{2\alpha-1}_{2,2}(0,1)}
		\leq
			C_2 \|w\|_{B^{2\alpha}_{2,2}(0,1)}
		\leq
			C_3 \|w\|_{H_\alpha}.
	\end{equation}
	Thus there exists a $C \in (0,\infty)$ such that for all
	$u,v \in H$ it holds that
	\begin{equation}
	\label{eq: F loc Lipschitz Burger}
		\begin{split}
				&\|F(u)-F(v)\|_{H_{-\alpha}}
			={} 
				\sup_{w \in H_\alpha} \bigg(
					\int_0^1 
						-\tfrac c2 (u^2 -v^2)' (x) \tfrac {w(x)}{\|w\|_{H_\alpha}} 
					\ud x
				\bigg) \\
			=
				&\sup_{w \in H_\alpha} \bigg(
					\tfrac c2 \int_0^1	(u^2-v^2)(x) \tfrac {w'(x)}{\|w\|_{H_\alpha}} \ud x 
				\bigg) \\
			\leq{} 
				&\sup_{w \in H_\alpha} \bigg(
					\tfrac c2 
					\|(u-v)(u+v)\|_{L^1(0,1)} 
					\cdot \tfrac {\|w'\|_{L^\infty(0,1)}}{\|w\|_{H_\alpha}}
				\bigg) 
			\leq
				C \|u-v\|_{H} \cdot \|u+v\|_{H}
		\end{split}
	\end{equation}
	and this 
	together with \eqref{eq: Lipschitz bound B Burger}
	shows \eqref{eq: Lipschitz bound for F and B}
		(with $\gamma\leftarrow 0$).
	In addition, we obtain that for all $N \in \N_0$ and 
	all $u \in H_{2 \vartheta} \cap V_N$ it holds that 
	\begin{equation}
		\begin{split}
				&\langle 
					\pi^{H_{-\alpha}}_{V_N}F(u), 
					I^{-1}_{\mathbb{H}} (D_{\mathbb{H}} V) (u) 
				\rangle_{H}
			=
				-\tfrac{c p}{2} \|u\|^{p-2}_{H} \,
				\langle 
					\pi^{H_{-\alpha}}_{V_N} (u^2)', 
					u 
				\rangle_{L^2(0,1)} \\
			={}
				&\tfrac{c p}{2} \|u\|^{p-2}_{H}
				\int_0^1 	
					u u'
					\cdot u
				\ud x 
			=
				\tfrac{c p}{6} \|u\|^{p-2}_{H}
				\int_0^1 	
					(u^3)' 
				\ud x
			=0.
		\end{split}
	\end{equation}
	This together with \eqref{eq: H bound B Burger}
	yields that for all $N \in \N_0$, $u \in H_{2\vartheta} \cap V_N$ it holds that
	\begin{equation}
	\label{eq: V is lypanov Burgers}
		\begin{split}
				&\big \langle
					\pi^{H_{-\alpha}}_{V_N} F(u) + \Delta u,
					I^{-1}_{\mathbb{H}} (D_{\mathbb{H}} V)(u)
				\big \rangle_{H} 
				+
				\nicefrac 12
				\big \langle 
					\pi^{H_{-\beta}}_{V_N} B(u),
					I_{\mathbb{H}_{\vartheta}}^{-1} 
						\big( 
							(D_{\mathbb{H}}^2 \, V)(u) \, 
							\pi^{H_{-\beta}}_{V_N} B(u) 
						\big)
				\big \rangle_{HS(\mathbb{U},\mathbb{H}_{\vartheta})} \\
			={} 
				&p \|u\|_H^{p-2}
				\big \langle
					\Delta u, u 
				\big \rangle_{H}
				+ p (p-2) \|u\|_H^{p-4}
				\big \langle 
					\pi^{H_{-\beta}}_{V_N} B(u),
					(-\Delta)^{-2 \vartheta} u \, \langle \pi^{H_{-\beta}}_{V_N} B(u), u \rangle_H
				\big \rangle_{HS(\mathbb{U},\mathbb{H}_{\vartheta})} \\
				&+p \|u\|_H^{p-2}
				\big \langle 
					\pi^{H_{-\beta}}_{V_N} B(u),
					(-\Delta)^{-2 \vartheta} \, \pi^{H_{-\beta}}_{V_N} B(u)
				\big \rangle_{HS(\mathbb{U},\mathbb{H}_{\vartheta})} \\
			\leq{} 
				&-p \|u\|_H^{p-2}
				\|
					u
				\|^2_{H_{1/2}}
				+ p (p-2) \|u\|_H^{p-4}
				\big \langle 
					\pi^{H_{-\beta}}_{V_N} B(u),
					 u \, \langle \pi^{H_{-\beta}}_{V_N} B(u), u \rangle_H 
				\big \rangle_{HS(\mathbb{U},\mathbb{H})} \\
				&+ p \|u\|_H^{p-2}
				\| 
					\pi^{H_{-\beta}}_{V_N} B(u)
				\|^2_{HS(\mathbb{U},\mathbb{H})} \\
			\leq{} 
				&p (p-2) \|u\|_H^{p-2}
				\|
					\pi^{H_{-\beta}}_{V_N} B(u)
				\|^2_{HS(\mathbb{U},\mathbb{H})} 
				+ p \|u\|_H^{p-2}
				\| 
					\pi^{H_{-\beta}}_{V_N} B(u)
				\|^2_{HS(\mathbb{U},\mathbb{H})} \\
			\leq{}
				&\theta \, p \, (p-1) \|u\|_H^{p-2} (1+\|u\|^2_{H}) \\
			\leq{} 
				&2 \, \theta \, p \, (p-1) (1+\|u\|^p_{H}) 
			\leq
				2 \theta \, p \, (p-1) V(u).
		\end{split}
	\end{equation}
	Furthermore, we get
	from the product-to-sum identity that
	for all $N \in \N$, $x \in (0,1)$ and all
	$c_i \in \R$, $i \in [1,N] \cap \N$, it holds that
	\begin{equation}
		\begin{split}
				&(F(\sum_{k=1}^N c_k e_k))(x)
			=  
				-c \sum_{k,l=1}^N 2\pi \, l \, c_k \, c_l \sin(k \pi x) \, \cos(l \pi x) \\
			={} 
				&-c \sum_{k,l=1}^N l \, \pi \, c_k \, c_l  \cdot 
					( \sin((k-l) \pi x) + \sin((k+l) \pi x))
		\end{split}
	\end{equation}
	and this ensures that for $N \in \N$ it holds that
	$F(V_N) \subseteq \bigcap_{r=1}^{\infty} H_r$.
	Combining this and
	Lemma \ref{l: norm equivalence Burger}
	shows that for all $r, s \in [0,\infty)$
	it holds that
	\begin{equation}
	\label{eq: same continuity}
			F|_{H_r} \in 
				C_{\mathbb{B}^{2r}_{2,2}(0,1) ,\mathbb{B}^{2s}_{2,2}(0,1)}
					(H_r, B^{2s}_{2,2}(0,1)) 
		\Leftrightarrow
			F|_{H_r} \in 
				C_{\mathbb{H}_r ,\mathbb{H}_s}
					(H_r, H_s).
	\end{equation}
	Next note that it follows from Lemma \ref{lem: lemma: interpolation theorem for H},
	Lemma \ref{l: norm equivalence Burger},
	and from Corollary \ref{cor: F assumption} \eqref{eq: h+1/2 continuity of F Cor2}
	(with $d \leftarrow 1$, $s \leftarrow \tfrac 12 +r$, 
	and with $\alpha \leftarrow 0$)
	that for all /
	$r \in [-\nicefrac 12, \infty) \backslash \{-\nicefrac 14\}$
	there exists a $C \in (0,\infty)$ such that
	for all $u,v \in H_{(1/2+r)\vee(r/2+3/8)}$ it holds that
	\begin{equation}
	\label{eq: product estimate}
				\|(u-v)(u+v)\|_{B^{1+2r}_{2,2}(0,1)}
			\leq
				C\|u+v\|_{H_{(1/2+r)\vee(r/2+3/8)}}
					\|u-v\|_{H_{(1/2+r)\vee(r/2+3/8)}}.
	\end{equation}
	Moreover, Lemma 4.9 in
	Jentzen, Lindner \& Pušnik \cite{JentzenLindnerPusnik2019}
	yields that for all $r \in [-\nicefrac12, 0]$
	there exists a $C \in (0,\infty)$ such that
	for all $u,v \in H_{(1/2+r)\vee(r/2+3/8)}$ it holds that
	\begin{equation}
	\label{eq: F estimate Burger r negative}
		\begin{split}
				&\|F(u)-F(v)\|_{H_{r}}
			=
				\|-c (uu'-vv')\|_{H_{r}}
			=
				\tfrac{c}{2} \|(u^2-v^2)'\|_{H_{r}}
			\leq
				C \|(u-v)(u+v)\|_{B^{1+2r}_{2,2}(0,1)}.
		\end{split}
	\end{equation}
	In addition, 
	for all $r \in [0,\infty)$
	there exists a $C \in (0,\infty)$ such that
	for all $u,v \in H_{1/2+r}$ it holds that
	\begin{equation}
		\begin{split}
	\label{eq: F estimate Burger r positive}
				&\|F(u)-F(v)\|_{B^{2r}_{2,2}(0,1)}
			=
			 \|-c (uu'-vv')\|_{B^{2r}_{2,2}(0,1)}
			=
				\nicefrac c2 \|((u-v)(u+v))'\|_{B^{2r}_{2,2}(0,1)} \\
			\leq{}
				&C  \|(u-v)(u+v)\|_{B^{2r+1}_{2,2}(0,1)}.
		\end{split}
	\end{equation}
	Combining \eqref{eq: F estimate Burger r negative},
	\eqref{eq: F estimate Burger r positive},
	\eqref{eq: product estimate},
	and \eqref{eq: same continuity} ensures that
	for all 
	$r \in [-\nicefrac 12, \infty) \backslash \{-\nicefrac 14\}$
	it holds that
	\begin{equation}
	\label{eq: F continuity r Burger}
		F|_{H_{(r+1/2)\vee (r/2+3/8)}} \in 
				C_{\mathbb{H}_{(r+1/2)\vee (r/2+3/8)},\mathbb{H}_r}
					(H_{(r+1/2)\vee (r/2+3/8)}, H_r)
	\end{equation}
	and combining \eqref{eq: F estimate Burger r negative},
	\eqref{eq: F estimate Burger r positive},
	\eqref{eq: F continuity r Burger}, and
	Lemma \ref{l: norm equivalence Burger}
	establishes that
	for all 
	$r \in [-\nicefrac 12, \infty)$
	and all $\eps \in (0,\infty)$
	there exists a $C \in (0,\infty)$ such that for all
	$u,v \in H_{(1/2+r+\eps 1_{\{r\}}(-\nicefrac 14))\vee(r/2+3/8)}$
	it holds that
	\begin{equation}
	\label{eq: F estimate Burger r}
			\|F(u)-F(v)\|_{H_{r}}
		\leq
			C \|(u-v)(u+v)\|_{B^{1+2r}_{2,2}(0,1)}.
	\end{equation}
	In particular we get from \eqref{eq: F continuity r Burger}
	(with $r \leftarrow \vartheta -1/2$)
	that
	\begin{equation}
	\label{eq: F continuity vartheta Burger}
		F|_{H_{\vartheta}} \in 
				C_{\mathbb{H}_{\vartheta},\mathbb{H}_{\vartheta -1/2}}
					(H_{\vartheta}, H_{\vartheta -1/2}).
	\end{equation}
	Combining 
	\eqref{eq: F continuity vartheta Burger},
	and
	\eqref{eq: B continuity Burger}
	verifies \eqref{eq: F and B continuity}
	(with $\alpha_1 \leftarrow \nicefrac 18$).
	Furthermore,
	it follows from \eqref{eq: F continuity vartheta Burger},
	Lemma \ref{l: norm equivalence Burger},
	Lemma \ref{lem: lemma: interpolation theorem for H},
	and from 
	\eqref{eq: h+1/2 bound for F Cor} in
	Corollary \ref{cor: F assumption}
		(with $d \leftarrow 1$, $\gamma \leftarrow 0$, $\alpha_1 \leftarrow \nicefrac 18$
		and $\alpha \leftarrow 1$) 
	that there exist $C_1, C_2 \in (0,\infty)$
	such that for all $u \in H_{\vartheta +3/8}$ it holds that
	\begin{equation}
		\label{eq: F vartheta -1/2 bound Burger}
			\|F(u)\|_{H_{\vartheta-1/2}}
		\leq
			C_1\|u'u\|_{B^{2\vartheta-1}_{2,2}(0,1)}
		\leq
			C_2 \|u\|_{H} \|u\|_{H_{\vartheta +1/2 -1/8}}
	\end{equation}
	and this 
	ensures \eqref{eq: F H 1/2 bound} 
		(with $\alpha_1 \leftarrow \nicefrac 18$,
		$r_x \leftarrow 1$ and $R_x \leftarrow C_2 (\|x\|_{H} +1)$).
  In addition, 
	\eqref{eq: V is lypanov Burgers} assures
	that for all $N \in \N$ and all
	all $(t,u) \in [0,T] \times V_N$ 
	it holds that 
	\begin{equation}
		\E [ V(X_{t}^{N,u})] \leq V(u) e^{2\theta \, p \, (p-1) \, t}
	\end{equation}
	which shows \eqref{eq: V bound for XN}.
Moreover, it follows from (5.14) 
in Cox, Hutzenthaler \& Jentzen \cite{CoxHutzenthalerJentzen2013}
(with $r \leftarrow 3$)
that there exist $C_1, C_2 \in (0,\infty)$ such that 
for all $\eps \in (0,\infty)$ it holds that
\begin{equation}
	\begin{split}
			&\limsup_{y \to x} 	\sup_{r \in [0,T]}
				\P(\|X^{0,x}_r - X^{0,y}_r\|_{H} \geq \eps) 
		\leq
			\limsup_{y \to x} 	\sup_{r \in [0,T]}
				\frac
					{
						\|X^{0,x}_r - X^{0,y}_r\|^3_{
							L^3( 
								\P ; 
								\left\| \cdot \right\|_{  H  }
						) 
						}
					}
					{\eps^3} \\
		\leq{} 
			&\limsup_{y \to x}
			 \big[ 
				\tfrac{C_1}{\eps^3} \|x-y\|^3_H \cdot 
					\exp\big( C_2 (\|x\|^2_H +\|y\|_H^2 +1\big) 
				\big]
		=0
	\end{split}
\end{equation}
and this proves \eqref{eq: cont assumption}.
Furthermore, note that 
\eqref{eq: B bounded wrt H_vartheta Burger} verifies 
\eqref{eq: B H 1 bound}
	(with $r_x \leftarrow 1$
	and $R_x \leftarrow K(\|x\|_H+1)$), 
that \eqref{eq: B C1 bound Burger} verifies
\eqref{eq: B C1 bound thm},
and that
\eqref{eq:convLocLip_assumption Burger} verifies 
\eqref{eq:convLocLip_assumption thm}.
Therefore we derive from Theorem \ref{thm: existence} that
  $ 
    f \in \C_{\R \times \mathbb{H}}(\R \times H, \R)
  $,
	that  
	$ 
    f 
  $
	is bounded on $\R \times \mathbb{H}$-bounded 
	subsets of $(0,T) \times H$,
	that for all $\tilde{R} \in (0, \infty)$ it holds that
	\begin{equation}
		\begin{split}
				\lim_{r \downarrow 0} \lim_{R \to \infty} \lim_{\eps \downarrow 0}	
				\sup \Bigg \{
					&f(t,u) -  \varphi(v)
					\colon
					~u, v \in H_{2 \vartheta},
					~t \in (0,T),
					~\|u\|_{H_{\vartheta}}^2
						\vee 
						\|v\|_{H_{\vartheta}}^2  
						\leq R, \\
					~&\|u \|_{H} \leq \tilde{R}, 
					~\| u - v\|_{H} \leq r,
					~t \leq \eps 
				\Bigg \}
			= 
				0,
		\end{split}
	\end{equation}
	that for all $q \in [p,\infty)$
  \begin{equation}
    \lim_{ r \to \infty }
      \sup_{ t \in [0,T]}
			\sup \left \{
				\frac
					{| f(t,u) |}
					{\|u\|^q_{L^2(0,1)}+2}
				\colon u \in H, ~\|u\|_{H} \geq r
			\right \}
    = 0,
  \end{equation}
  and that
  $ 
    f|_{ (0,T) \times H} 
  $
  is a viscosity solution of
  \begin{equation} 
		\begin{split}
			&\tfrac{ \partial }{ \partial t }
			f(t,u) -
			\big\langle
				\Delta u- \tfrac{c}{2} (u^2)', 
				I^{-1}_{\mathbb{H}} (D_{\mathbb{H}} f)(t,u)
			\big\rangle_{H} 
			-
			\big \langle 
				B(u),
				I_{\mathbb{H}_{\vartheta}}^{-1}
					\big( (D^2_{\mathbb{H}} \, f)(t,u) \, B(u) \big)
			\big \rangle_{HS(\mathbb{U},\mathbb{H}_{\vartheta})} 
			= 0
		\end{split}
  \end{equation}
	for $ (t,u) \in (0,T) \times H $ relative to 
	$
		(
			(0,T) \times H \ni (t,u) 
			\to \nicefrac 12 \|u \|^2_{\mathbb{H}_{\vartheta}} \in [0, \infty],
			\R \times \mathbb{H}, 
			\R \times \mathbb{H}_{\vartheta}
		)
	$. 
	For the uniqueness we apply Corollary \ref{cor:uniqueness2}
	(with $O \leftarrow H$, $v \leftarrow 0$, 
	and $\kappa \leftarrow 0$).
	Therefore assume that $g$ is another viscosity solution of 
	\eqref{eq: Kol equation Burger}, 
	which satisfy that 
	$ 
    g \in \C_{\R \times \mathbb{H}}([0,T] \times H, \R)
  $,
	that $g$ is
	bounded on $\R \times \mathbb{H}$-bounded 
	subsets of $[0,T] \times H$, 
	that $g$ have at most polynomial $\H$-growth,
  and that for all $ u \in H $ it holds that
  $
    g(0, u) = \varphi(u)
  $.
	Thus there exists a $q \in [p,\infty)$, which we fix for the rest of the proof,
	such that
	 \begin{equation}
    \lim_{ r \to \infty }
      \sup_{ t \in [0,T]}
			\sup \left \{
				\frac
					{| g(t,u) |}
					{\|u\|^q_{L^2(0,1)}+2}
				\colon u \in H, ~\|u\|_{H} \geq r
			\right \}
    = 0.
  \end{equation}
	Furthermore, the
	continuity of $\varphi$
	with respect to the $\| \cdot \|_{H}$-norm,
	the
	continuity of $g$
	with respect to the $\| \cdot \|_{\R \times H}$-norm,
	and the fact that for all $R \in (0,\infty)$ it holds that
	the set $\{(t,x) \in [0,T] \times H \colon \|x\|_{H_\vartheta} \leq R\}$
	is compact in $\H$
	ensure for all $R, \tilde{R} \in (0, \infty)$ that
	\begin{equation}
		\begin{split}
				&\lim_{r \downarrow 0} \lim_{\eps \downarrow 0}	
				\sup \bigg \{
					g(t,u) -  \varphi(\hat{u})
					\colon
					u, \hat{u} \in H_{2\varz},
					~t, \hat{t} \in (0,T),
					~\tfrac{e^{-\theta t} \|u\|_{H_{\varz}}^2}{\|u\|_H^q+2} 
						\vee \tfrac{e^{-\theta \hat{t}} 
							\|\hat{u}\|_{H_{\varz}}^2}{\|\hat{u}\|_H^q+2} \leq R, 
				\\& \qquad \qquad \qquad \qquad
					~\|u \|_{H} \leq \tilde{R}, 
					~\| u - \hat{u}\|_{H} \leq r,
					~t \vee \hat{t} \leq \eps 
				\bigg \} \\
			={} 
				&\lim_{r \downarrow 0}
				\sup \bigg \{
					g(0,u) -  \varphi(\hat{u})
					\colon
					u, \hat{u} \in H_{2\varz},
					~\tfrac {\|u\|_{H_{\varz}}^2} {\|u\|_H^q+2} 
						\vee \tfrac {\|\hat{u}\|_{H_{\varz}}^2} {\|\hat{u}\|_H^q+2} \leq R,
					~\|u \|_{H} \leq \tilde{R},
					~\| u - \hat{u}\|_{H} \leq r
				\bigg \} \\
			={}  
				&\lim_{r \downarrow 0}
				\sup \bigg \{
					\varphi(u) -  \varphi(\hat{u})
					\colon
					u, \hat{u} \in H_{2\varz},
					~\tfrac {\|u\|_{H_{\varz}}^2} {\|u\|_H^q+2} 
						\vee \tfrac {\|\hat{u}\|_{H_{\varz}}^2} {\|\hat{u}\|_H^q+2} \leq R,
					~\|u \|_{H} \leq \tilde{R},
					~\| u - \hat{u}\|_{H} \leq r
				\bigg \} \\ 
			={}
				&0.
		\end{split}
	\end{equation}
	Next note that we derive from \eqref{eq: F estimate Burger r}
		(with $r \leftarrow \vartheta -1$,
		$v \leftarrow 0$, and $\eps \leftarrow \nicefrac 12$),
	Lemma \ref{lem: lemma: interpolation theorem for H},
	and from \eqref{eq: h bound for F Cor} in	Corollary \ref{cor: F assumption} 
		(with $d \leftarrow 1$ and $\gamma \leftarrow 0$)
	that
	there exist $C_1,C_2 \in (0,\infty)$ such that for all $u \in H_{2\vartheta}$
	it holds that
	\begin{equation}
		\begin{split}
				&\| F (u)\|_{H_{\vartheta-1}}
			\leq
				C_1\|u^2\|_{B^{2\vartheta-1}_{2,2}(0,1)}
			\leq
				C_2 \|u\|^2_{H_\vartheta}.
		\end{split}
	\end{equation}
	Combining this and 
	\eqref{eq: H bound B Burger} verifies
	\eqref{eq: H bound for F, B and v}. 
	In addition, \eqref{eq: F vartheta -1/2 bound Burger}
	and
	\eqref{eq: B bounded wrt H_vartheta Burger}
	imply \eqref{eq: H_3/2 bound for F} 
		(with $\alpha_1 \leftarrow \nicefrac 18$)
	and \eqref{eq: H_3/2 bound for B}.
	In addition,
	we get from \eqref{eq: F estimate Burger r} 
		(with $r \leftarrow - \nicefrac 12$)
	and from \eqref{eq: Lipschitz continuity of F Cor} in
	Corollary \ref{cor: F assumption}
		(with $d \leftarrow 1$, $\gamma \leftarrow 0$,
		$\alpha_2 \leftarrow \nicefrac 18$,
		$u \leftarrow u+v$, and $v \leftarrow u-v$) 
	that
	there exists $C_1,C_2 \in (0,\infty)$ such that for all $u,v \in H_{2\vartheta}$
	it holds that
	\begin{equation}
		\begin{split}
				&\|F(u) - F(v)\|_{H_{-1/2}}
			\leq
				C_1 \|(u+v) (u-v)\|_{B^{0}_{2,2}}
			\leq
				C_2 \|u+v\|_H \|u-v\|_{H_{1/2-1/8}}
		\end{split}
	\end{equation}
	and this shows
	\eqref{eq: Lipschitz continuity of F}
		(with $\alpha_2 \leftarrow \nicefrac 18$).
	Furthermore, \eqref{eq: B Lipschitz wrt H_1/2 Burger}
	verifies \eqref{eq: Lipschitz continuity of B},
	\eqref{eq: V is lypanov Burgers} verifies \eqref{eq:assumption_V},
	and Remark \ref{rem: assumption on V} 
		(with $\kappa \leftarrow 0$, $l \leftarrow 0$, and with
		$f \leftarrow (\R \ni x \to  (2x)^{\nicefrac p2} + 2 \in (1,\infty))$)
	verifies
	\eqref{eq: V continuous with respect to H1},
	\eqref{eq: V Lipschitz with respect to H -gamma},
	\eqref{eq: V bounded with respect to H},
	and \eqref{eq: V bounded with respect to H 1/2}.
	Thus Corollary \ref{cor:uniqueness2} 
	(with $V \leftarrow \| \cdot\|^q_H+2$) 
	shows that $f=g$
	which finishes
	the proof of Corollary \ref{cor: Burgers equation}.
\end{proof}
In the next remark we give examples of functions
$
	B
 \in 
		\C_{
			\mathbb{H}, 
			\mathbb{HS}(\mathbb{U}, \mathbb{H}_{-\beta})
		}( 
			H, 
			HS( \mathbb{U}, \mathbb{H}_{-\beta} ) 
		)
$ satisfying the assumptions of Corollary \ref{cor: Burgers equation}.
\begin{remark}[Sufficient condition for B]
\label{rem: sufficient condition for B}
	Assume 
	the Setting \ref{sett: Burger's},
	let
	$C \in (0,\infty)$,
	$s \in (\nicefrac 34, \infty)$,
	let $\mathscr{P}_N \in L(\mathbb{H}_s,\mathbb{H}_s)$, $N \in \N$,
	be finite-dimensional projections
	satisfying for all $i \in \N$ that
	\begin{equation}
	\label{eq: PN converges in L infty Burgers}
		\lim_{N \to \infty} \|(\id_{H_s} - \mathscr{P}_N) e_i\|_{L^\infty(0,1)} = 0
	\end{equation}
	and that
	\begin{equation}
	\label{eq: PN is bounded in L infty Burgers}
			\sup_{N \in \N} \|\mathscr{P}_N e_i\|^2_{L^{\infty}(0,1)} 
		\leq 
			C |\lambda_i|,
	\end{equation}
	let $b \in \C^2([0,1] \times \R; \R)$ satisfy that
	for all $x \in [0,1]$ and all $y \in \R$ it holds that
	\begin{equation}
	\label{eq: b derivative bound}
		\|(D_{\R \times \R} b) (x,y) \|_{(\R^2)'} \leq C,
	\end{equation}
	let $B \colon H \to HS(\mathbb{H}_s, \mathbb{H})$ be the function
	satisfying for all $u \in H$, $v \in H_s$
	and all $x \in (0,1)$ that
	$(B(u) v)(x) = b(x, u(x)) v(x)$.
	Then B satisfies the condition of 
	Corollary \ref{cor: Burgers equation} 
	with $\vartheta \leftarrow \nicefrac 12$, $\beta \leftarrow 0$,
	$\beta_1 \leftarrow \nicefrac 12$, $\beta_2 \leftarrow \nicefrac 12$, 
	$\beta_3 \leftarrow \nicefrac 12$,
	$\chi \leftarrow 0$, 
	and with $\mathbb{U} \leftarrow \mathbb{H}_s$
	i.e.\@ there exists 
	a $\theta \in (0, \infty)$
	and an increasing
	$K \colon \R \to (0, \infty)$,
	such that
	\begin{equation}
	\label{eq: B continuity Burger rem}
		B|_{H_{1/2}} \in
				\C_{
					\mathbb{H}_{1/2},
					\mathbb{HS}(\mathbb{H}_s, \mathbb{H}_{1/2})
				}(
					H_{1/2},
					HS(\mathbb{H}_s, \mathbb{H}_{1/2})
				), 
	\end{equation}
	that for all $u, v \in H_{1}$ 
	it holds that
	\begin{align}
		\label{eq: H bound B Burger rem}
			&\| B(u) \|^2_{HS(\mathbb{H}_s,\mathbb{H})}
		\leq
			\theta \cdot ( 1 +\|u\|^2_{H}), \\
	\label{eq: B bounded wrt H_vartheta Burger rem}
			&\| B(u) \|^2_{HS(\mathbb{H}_s,\mathbb{H}_{1/2})} 
		\leq
			K(\|u\|_{H}) 
				(\|u\|^2_{H_{1/2}} + 1), \\
	\label{eq: B Lipschitz wrt H_1/2 Burger rem}
			&\|B(u) - B(v)\|_{HS(\mathbb{H}_s, \mathbb{H})}
		\leq
			K(\|v\|_{H} \vee \|u\|_{H}) 
				\|u-v\|_{H}, \\
	\label{eq: Lipschitz bound B Burger rem}
		&|B |_{
			\C^1( 
				H, 
				\left\| \cdot \right\|_{HS( \mathbb{H}_s, \mathbb{H} )} 
			)
		}
		< \infty,
	\end{align}
	and that for all $\mathbb{H}$-bounded sets
	$ E \subseteq H$
	it holds that
	\begin{equation}
	\label{eq: B C1 bound Burger rem} 
			\sup_{ N \in \N } \sup_{ v \in E }
			\left[
				\frac{
				\|
					B( v ) \mathscr{P}_N 
				\|_{ HS( \mathbb{H}_s, \mathbb{H} ) }
				}{
					1 + \| v \|_{ H }
				}
			\right]
			< \infty
	\end{equation}
	and that
	\begin{equation}
	\label{eq:convLocLip_assumption Burger rem}
			\limsup_{ N \to \infty }
			\sup_{ v \in E }
			\left[
			\frac{
			\|
				B( v ) ( \id_{H_s} - \mathscr{P}_N )
			\|_{ HS( \mathbb{H}_s, \mathbb{H} ) }
			}{
				1 + \| v \|_{ H }
			}
			\right]
		=0.
	\end{equation}
\end{remark}
\begin{proof}
	Denote by $W^1_{0}(0,1) \subseteq B^1_{2,2}(0,1)$
	the set satisfying that
	$
			W^1_{0}(0,1) 
		= 
			\overline{C_0^\infty(0,1)}_{\mathbb{B}^1_{2,2}(0,1)}
	$.
	Then Example 4.34 in Lunardi \cite{Lunardi2018}
	shows that $H_{1/2} = W^1_{0}(0,1)$.
	Hence the fact that $\forall k \in \N \colon e_k \in H_{1/2}$
	proves that for all $k \in \N$ there exists a sequence 
	$(u_{n,k})_{n \in \N} \subseteq \C_0^\infty(0,1)$
	such that for all $k\in \N$
	it holds that
	$\lim_{n \to \infty} \|e_k -u_{n,k}\|_{B^1_{2,2}(0,1)} =0$.
	Furthermore,
	Lemma \ref{lem: b L2 boundedness inequality}
		(with $d \leftarrow 1$, $m \leftarrow 1$,
		$p \leftarrow 0$),
	\eqref{eq: b derivative bound},
	and item (iv) in Lemma 4.3 in Jentzen, Lindner \& Pušnik \cite{JentzenLindnerPusnik2019}
	establish
	that there exist $C_1, C_2 \in (0,\infty)$ such
	that for all $u \in H_{1/2}$ it holds that
	\begin{equation}
	\label{eq: B estimate in H 1/2 Burger}
		\begin{split}
				&\|b(\cdot,u(\cdot)) \|^2_{B^{1}_{2,2}(0,1)} \\
			={}
				&\Big(
					\big \| 
						\big \langle
							(D_{\R \times \R} b)(\cdot,u(\cdot)),
							(1,u'(\cdot))
						\big \rangle_{(\R^2)', \R^2}
					\big \|^2_{H} 
					+\|b(\cdot,u(\cdot)) \|^2_{L^2(0,1)}
				\Big)\\
		\leq{} &
				C_1\Big(
					\big \|  
						\| (D_{\R \times \R} b)(\cdot,u(\cdot)) \|_{(\R^2)'} 
					\big \|^2_{L^\infty(0,1)} \, 
					\big\| 
						\|(1, u'(\cdot)) \|_{\R^2} 
					\big\|^2_{H}
					+\|u\|^2_{H}
					+\int_0^1 |b(x,0)|^2\ud x
				\Big) \\
			\leq{}
				&
				C_2\Big(
					1+
					\| 
						u' 
					\|^2_{H}
					+\|u\|^2_{H}
				\Big)
			\leq
				2C_2
				(\| u\|^2_{H_{1/2}}+1)
		\end{split}
	\end{equation}
	and this verifies that for all $u \in H_{1/2}$ it holds that
	$b( \cdot, u(\cdot)) \in B_{2,2}^1(0,1)$.
	Thus we obtain that for all $u \in H_{1/2}$ there exists a sequence
	$(v_n)_{n\in \N} \subseteq C^\infty(0,1)$ such that
	$\lim_{n \to \infty} \|v_n(\cdot)-b( \cdot, u(\cdot))  \|_{B^1_{2,2}(0,1)} =0$.
	Therefore we get
	from Theorem 2 on page 191 in Runst \& Sickel \cite{RunstSickel1996}
	that for all $k \in \N$
	and all $u \in H_{1/2}$ there exist sequences 
	$(u_{n,k})_{n \in \N} \subseteq C_0^\infty(0,1)$
	and $(v_{n})_{n \in \N} \subseteq C^\infty(0,1)$ such that
	\begin{equation}
		\begin{split}
				&\lim_{n \to \infty}
					\| b(\cdot,u(\cdot)) e_k(\cdot) - u_{n,k}(\cdot) v_n(\cdot)\|_{B_{2,2}^1(0,1)}\\
			\leq{}
				&\lim_{n \to \infty}
					\Big(
						\| \big( b(\cdot,u(\cdot)) -v_n(\cdot) \big) \, e_k(\cdot)\|_{B_{2,2}^1(0,1)}
						+\|  (e_k - u_{n,k}) v_n\|_{B_{2,2}^1(0,1)} 
					\Big) \\
			\leq{}
				&\lim_{n \to \infty}
					\Big(
						\| b(\cdot,u(\cdot)) -v_n(\cdot)\|_{B_{2,2}^1(0,1)}
						\| e_k\|_{B_{2,2}^1(0,1)}
						+\|  e_k - u_{n,k}\|_{B_{2,2}^1(0,1)}
						\|v_n\|_{B_{2,2}^1(0,1)}
					\Big)
			=0
		\end{split}
	\end{equation}
	and this implies that for all $k \in \N$ and all $u \in H_{1/2}$
	it holds that
	\begin{equation}
	\label{eq: b e in right space Burger}
		 b(\cdot,u(\cdot)) e_k(\cdot) \in W^1_0(0,1) =H_{1/2}.
	\end{equation}
	Next note that it follows from
	the fact that $\forall k \in \N \colon \|e_k\|_{L^\infty(0,1)} = \sqrt{2}$ 
	and from the fact that $s>\nicefrac 34$ that there exists a $C_1 \in (0,\infty)$
	such 
	that for all $u \in H$ it holds that
	\begin{equation}
	\label{eq: basic estimate H}
		\begin{split}
				&\|B(u)\|^2_{HS(\H_s,\H)}
			=
				\sum^\infty_{k=1}
					\|b(\cdot,u(\cdot)) \, |\lambda_k|^{-s} e_k (\cdot)\|^2_{H}
			\leq
				\sum^\infty_{k=1}
					\|b(\cdot,u(\cdot)) \|^2_H \, 
						\| |\lambda_k|^{-s} e_k \|^2_{L^{\infty}(0,1)} \\
			={} &  
				2\|b(\cdot,u(\cdot)) \|^2_H
					\sum^\infty_{k=1}
						|\lambda_k|^{-2s}
			=
				2\|b(\cdot,u(\cdot)) \|^2_H
					\sum^\infty_{k=1}
						(\pi k)^{-4s} 
			\leq
				C_1\|b(\cdot,u(\cdot)) \|^2_H.
		\end{split}
	\end{equation}
	Similar \eqref{eq: b e in right space Burger},
	Lemma \ref{l: norm equivalence Burger},
	Theorem 2 on page 191 in Runst \& Sickel \cite{RunstSickel1996},
	and the fact that $s> \nicefrac 34$
	imply that there exist $C_1, C_2, C_3, C_4 \in (0,\infty)$ such that 
	for all $u \in H_{1/2}$ it holds that
	\begin{equation}
	\label{eq: basic estimate H 1/2}
		\begin{split}
				&\|B(u)\|^2_{HS(\H_s,\H_{1/2})}
			=
				\sum^\infty_{k=1}
					\|b(\cdot,u(\cdot)) \, |\lambda_k|^{-s} e_k(\cdot)\|^2_{H_{1/2}}\\
			\leq{}
				&C_1 \sum^\infty_{k=1}
					\|b(\cdot,u(\cdot)) \, |\lambda_k|^{-s} e_k(\cdot)\|^2_{B^1_{2,2}(0,1)} 
			\leq{} 
				C_2
				\sum^\infty_{k=1}
					\|b(\cdot,u(\cdot)) \|^2_{B^{1}_{2,2}(0,1)} \, 
					\||\lambda_k|^{-s} e_k\|^2_{B^{1}_{2,2}(0,1)} \\
			\leq{}
				&C_3
				\sum^\infty_{k=1}
					\|b(\cdot,u(\cdot)) \|^2_{B^{1}_{2,2}(0,1)} \, 
					\||\lambda_k|^{-s} e_k\|^2_{H_{1/2}} 
			={} 
				C_3\|b(\cdot,u(\cdot)) \|^2_{B^{1}_{2,2}(0,1)}
					\sum^\infty_{k=1}
						|\lambda_k|^{-2s+1} \\
			={}  
				&C_3\|b(\cdot,u(\cdot)) \|^2_{B^{1}_{2,2}(0,1)}
					\sum^\infty_{k=1}
						(\pi k)^{2(-2s+1)} 
			\leq{} 
				C_4\|b(\cdot,u(\cdot)) \|^2_{B^{1}_{2,2}(0,1)}.
		\end{split}
	\end{equation}
	Analogously we obtain that there exist a 
	$C_1 \in (0,\infty)$ such that for all $i \in \{0,\nicefrac 12\}$ and all
	$u, v \in H_{i}$ it holds that
	\begin{equation}
	\label{eq: basic estimate diff}
				\|B(u)-B(v)\|_{HS(\H_s,\H_{i})} 
			\leq
				C_1 \|b(\cdot,u(\cdot)) -b(\cdot,v(\cdot))\|_{B^{2i}_{2,2}(0,1)}.
	\end{equation}
	In addition we deduce from \eqref{eq: basic estimate H} and 
	from Lemma \ref{lem: b L2 boundedness inequality} 
		(with $d \leftarrow 1$, $m \leftarrow 1$, and $p \leftarrow 0$)
	that there exist $C_1,C_2 \in (0,\infty)$
	such that for all $u \in H$ it holds that
	\begin{equation}
			\|B(u)\|^2_{HS(\H_s, \H)}
		\leq
			C_1 \|b(\cdot,u(\cdot)) \|^2_{H}
		\leq
			C_2 (\|u\|^2_H + \int_0^1 |b(x,0)|^2 \ud x)
	\end{equation}
	and this proves \eqref{eq: H bound B Burger rem}.
	Moreover, \eqref{eq: basic estimate diff} and
	Lemma \ref{lem: L2 b continuity inequality}
		(with $d \leftarrow 1$, $m \leftarrow 1$, $p \leftarrow 0$, and $q \leftarrow \infty$)
	verify that there exist $C_1,C_2 \in (0,\infty)$
	such that for all $u,v \in H$ it holds that
	\begin{equation}
	\label{eq: B diff H estimate Burger rem}
			\|B(u) -B(v)\|_{HS(\H_s, \H)}
		\leq
			C_1\|b(\cdot,u(\cdot)) -b(\cdot,v(\cdot))\|_{H}
		\leq
			C_2 \| u-v \|_{H}
	\end{equation}
	and this ensures \eqref{eq: B Lipschitz wrt H_1/2 Burger rem}
	and \eqref{eq: Lipschitz bound B Burger rem}.
	Combining
	\eqref{eq: basic estimate H 1/2}
	and
	\eqref{eq: B estimate in H 1/2 Burger}
	shows \eqref{eq: B bounded wrt H_vartheta Burger rem}.
	Next note that it holds for all $x \in (0,1)$
	and all $u,v \in H$ that
	\begin{equation}
	\label{eq: derivative b diff estimate}
		\begin{split}
				&(D_{\R \times \R} b)(x,u(x))
				- (D_{\R \times \R} b)(x,v(x)) \\
			={}
				&\int_0^1 
					(D^2_{\R \times \R} b)(x,t u(x) +(1-t) v(x) ) 
						(0,u(x)-v(x))
				\ud t \\
			\leq{} & 
				\| 
					D^2_{\R \times \R} b
				\|_{
					\C(
						[0,1]
						\times (
							-(\|u\|_{L^\infty(\Omega)} \vee \|v\|_{L^\infty(\Omega)}),
							\|u\|_{L^\infty(\Omega)} \vee \|v\|_{L^\infty(\Omega)}
						),
						\| \cdot \|_{L(\R^2, (\R^2)')}
					)
				} 
				|u(x)-v(x)|.
		\end{split}
	\end{equation}
	Moreover, the Sobolev inequality
	(see, e.g., Theorem 1 on page 32 in Runst \& Sickel \cite{RunstSickel1996}) and
	Lemma \ref{l: norm equivalence Burger}
	assure that there exist $C_1, C_2 \in (0,\infty)$
	such that for all $u \in H_{1/2}$ it holds that
	\begin{align}
	\label{eq: L infty estimat}
			\|u\|_{L^\infty(0,1)}
		\leq
			C_1\|u\|_{B^1_{2,2}(0,1)}
		\leq
			C_2\|u\|_{H_{1/2}}.
	\end{align}
	Therefore,
	\eqref{eq: b derivative bound},
	item (iv) in Lemma 4.3 in Jentzen, Lindner \& Pušnik \cite{JentzenLindnerPusnik2019},
	\eqref{eq: derivative b diff estimate},
	and \eqref{eq: L infty estimat}
	yield that there exist $C_1,C_2,C_3 \in (0,\infty)$
	such that for all $u, v \in H_{1/2}$
	it holds that
	\begin{align}
	\label{eq: B diff H 1/2 1st estimate Burger rem}
				&\big(
				 \big \| 
						\langle
							(D_{\R \times \R} b)(\cdot,u(\cdot)),
							(1,u'(\cdot))
						\rangle_{(\R^2)', \R^2}
					-\langle
						(D_{\R \times \R} b)(\cdot,v(\cdot)),
						(1,v'(\cdot))
					\rangle_{(\R^2)', \R^2}
				 \big \|_{H} 
				\big) \\ \nonumber
			\leq{} &
					\big(
					 \big \| 
						\langle
							(D_{\R \times \R} b)(\cdot,u(\cdot)),
							(0,  u'(\cdot) -v'(\cdot))
						\rangle_{(\R^2)', \R^2}
					\big\|_{H} 
		\\ \nonumber &
					+\big\|
						\langle
							(D_{\R \times \R} b)(\cdot,u(\cdot))
							- (D_{\R \times \R} b)(\cdot,v(\cdot)), 
							(1,v')
						\rangle_{(\R^2)', \R^2}
					\big\|_{H}
				\big) \\ \nonumber
			\leq{} &
				\big(
					C\| 
						u' -v'
					\|_{H} 
					+\big\|
						\|
							(D_{\R \times \R} b)(\cdot,u(\cdot))
							- (D_{\R \times \R} b)(\cdot,v(\cdot))
						\|_{(\R^2)'}
						(1+|v'(\cdot)|)
					\big\|_{H}
				\big) \\ \nonumber
			\leq{} &
				C_1 \big(
					\| 
						u -v
					\|_{H_{1/2}} 
			+
					\big\|
						\| 
							D^2_{\R \times \R} b
						\|_{
							\C(
								[0,1]
								\times (
									-(\|u\|_{L^\infty(\Omega)} \vee \|v\|_{L^\infty(\Omega)}),
									\|u\|_{L^\infty(\Omega)} \vee \|v\|_{L^\infty(\Omega)}
								),
								\| \cdot \|_{L(\R^2, (\R^2)')} 
							)
						} 
					\\ \nonumber & \quad \cdot
						|u-v| \,
						(1+|v'|)
					\big\|_{H}
				\big) \\ \nonumber
			\leq{} &
				C_1 \big(
					\| 
						u -v
					\|_{H_{1/2}} 
				+
					\| 
						D^2_{\R \times \R} b
					\|_{
						\C(
							[0,1]
							\times (
								-(\|u\|_{L^\infty(0,1)} \vee \|v\|_{L^\infty(0,1)}),
								\|u\|_{L^\infty(0,1)} \vee \|v\|_{L^\infty(0,1)}
							),
							\| \cdot \|_{L(\R^2, (\R^2)')} 
						)
					}
			\\ \nonumber & \quad \cdot
					\|
						u-v
					\|_{L^\infty(0,1)}
					(
						1+
						\|
							v'
						\|_{H}
					)
				\big) \\ \nonumber
			\leq{} &
				C_2 \big(
					\| 
						u -v
					\|_{H_{1/2}} 
					+
					\| 
						D^2_{\R \times \R} b
					\|_{
						\C(
							[0,1]
							\times (
								-C_3(\|u\|_{H_{1/2}} \vee \|v\|_{H_{1/2}}),
								C_3(\|u\|_{H_{1/2}} \vee \|v\|_{H_{1/2}})
							),
							\| \cdot \|_{L(\R^2, (\R^2)')} 
						)
					}
			\\ \nonumber & \quad
					\|
						u-v
					\|_{H_{1/2}}
					(
						1+
						\|
							v
						\|_{H_{1/2}}
					)
				\big).
	\end{align}
	In addition, \eqref{eq: basic estimate diff}
	assures that there exists a $C_1 \in (0,\infty)$
	such that for all $u, v \in H_{1/2}$
	it holds that
	\begin{equation}
	\label{eq: B diff H 1/2 basic estimate Burger rem}
		\begin{split}
				&\|B(u) -B(v) \|_{HS(\H_s, \H_{1/2})}
			\leq
				C_1 \|b(\cdot,u(\cdot)) -b(\cdot,v(\cdot))\|_{B^{1}_{2,2}(0,1)} \\
			={} 
				&C_1\big(
				 \big \| 
						\langle
							(D_{\R \times \R} b)(\cdot,u(\cdot)),
							(1,u'(\cdot))
						\rangle_{(\R^2)', \R^2}
					-\langle
						(D_{\R \times \R} b)(\cdot,v(\cdot)),
						(1,v'(\cdot))
					\rangle_{(\R^2)', \R^2}
				 \big \|_{H} 
		\\ & \quad
				+\|b(\cdot,u(\cdot)) -b(\cdot,v(\cdot))\|_{L^2(0,1)}
				\big).
		\end{split}
	\end{equation}
	Combining \eqref{eq: B diff H 1/2 basic estimate Burger rem},
	\eqref{eq: B diff H 1/2 1st estimate Burger rem},
	and \eqref{eq: B diff H estimate Burger rem}
	shows \eqref{eq: B continuity Burger rem}.
	Furthermore, it holds for all $v \in H$, and all $M,N \in \N$
	that
	\begin{equation}
	\label{eq: B id - proj estimate Burger}
		\begin{split}
				&\|
					B( v ) ( \id_{H_s} - \mathscr{P}_N )
				\|^2_{ HS( \mathbb{H}_s, \mathbb{H} ) }
			=	
					\sum^\infty_{k=1}
						\|
							b(\cdot, v(\cdot))  ( \id_{H_s} - \mathscr{P}_N ) 
								e_k(\cdot) |\lambda_k|^{-s} 
						\|^2_H  \\
			\leq{} & 
					\sum^\infty_{k=1}
						\|
							b(\cdot, v(\cdot)) \|^2_H \|(\id_{H_s} -\mathscr{P}_N) e_k |\lambda_k|^{-s}
						\|^2_{L^\infty(0,1)} \\
			\leq{} &
					\|b(\cdot, v(\cdot)) \|^2_H \big(
					\sum^\infty_{k=M+1}
						 \|(\id_{H_s} -\mathscr{P}_N) e_k |\lambda_k|^{-s} \|^2_{L^\infty(0,1)}
					+
					\sum^M_{k=1}
						 \|(\id_{H_s} -\mathscr{P}_N) e_k |\lambda_k|^{-s} \|^2_{L^\infty(0,1)}
					\big).
			\end{split}
		\end{equation}
		Furthermore it follows from \eqref{eq: B id - proj estimate Burger}, 
		Lemma \ref{lem: b L2 boundedness inequality} 
			(with $d \leftarrow 1$, $m \leftarrow 1$, and $p \leftarrow 0$),
		\eqref{eq: PN converges in L infty Burgers},
		the fact that 
		$\forall k \in \N \colon \|e_k\|_{L^\infty(0,1)}= \sqrt{2}$,
		\eqref{eq: PN is bounded in L infty Burgers},
		and from the fact that $s > \nicefrac 34$ that
		there exists a $C_1 \in (0,\infty)$
		such that for all $\mathbb{H}$-bounded sets $E$ it holds that
		\begin{equation}
			\begin{split}
				&\limsup_{ N \to \infty }
				\sup_{ v \in E }
				\left[
				\|
					B( v ) ( \id_{H_s} - \mathscr{P}_N )
				\|^2_{ HS( \mathbb{H}_s, \mathbb{H} ) }
				\right] \\
			\leq{} & 
				\sup_{ v \in E } (
					\|b(\cdot, v(\cdot)) \|^2_H
				)
				\Big(
					\lim_{M \to \infty}
					\limsup_{ N \to \infty }
					\big(
						\sum^\infty_{k=M+1}
							 \|(\id_{H_s} -\mathscr{P}_N) e_k |\lambda_k|^{-s} \|^2_{L^\infty(0,1)}
					\big) \\
			& \qquad 
						+
					\lim_{M \to \infty}
					\limsup_{ N \to \infty } \big(
						\sum^M_{k=1}
							 \|(\id_{H_s} -\mathscr{P}_N) e_k |\lambda_k|^{-s} \|^2_{L^\infty(0,1)}
					\big)
				\Big) \\
			\leq{} & 
				C_1
				\sup_{ v \in E } (
					(\|v\|^2_H + \int_0^1 |b(x,0)|^2 \ud x)
				)
		\\ & \qquad \cdot
				\lim_{M \to \infty}
					\limsup_{ N \to \infty }
					\big(
						\sum^\infty_{k=M+1}
							|\lambda_k|^{-2s} 
							(
								\| e_k  \|^2_{L^\infty}
								+ \sup_{n \in \N} 
									\|\mathscr{P}_n e_k  \|^2_{L^\infty(0,1)}
							)
					\big) \\
			\leq{} &
				C_1
				\sup_{ v \in E } (
					(\|v\|^2_H + \int_0^1 |b(x,0)|^2 \ud x)
				)
				\lim_{M \to \infty}
					\big(
						\sum^\infty_{k=M+1}
							(\pi k)^{-4s} 
							(
								2
								+ 
								C (\pi k)^2
							)
					\big) \\
			\leq{} &
				C_1
				\sup_{ v \in E } (
					(\|v\|^2_H + \int_0^1 |b(x,0)|^2 \ud x)
				)
				\lim_{M \to \infty}
					\big(
						(C+2)\sum^\infty_{k=M+1}
							(\pi k)^{-4s+2} 
					\big)
			= 
				0
		\end{split}
	\end{equation}
	and this verifies \eqref{eq:convLocLip_assumption Burger rem}.
	Finally note that
	Lemma \ref{lem: b L2 boundedness inequality} 
		(with $d \leftarrow 1$, $m \leftarrow 1$, and $p \leftarrow 0$),
	the assumption
	$
		\forall k \in \N \colon 
			\sup_{N \in \N}
				\|\mathscr{P}_N e_k\|^2_{L^\infty(0,1)} \leq C (\pi k)^2
	$,
	and the fact that $s > \nicefrac 34$
	implies that there exists a $C_1 \in (0,\infty)$
	such that for all $\mathbb{H}$-bounded sets $E$ it holds that
	\begin{align}
	\nonumber
				&\sup_{ N \in \N } \sup_{ v \in E }
					\|
						B( v ) \mathscr{P}_N 
					\|_{ HS( \mathbb{H}_s, \mathbb{H} ) }
			=	
				\sup_{ N \in \N }
				\sup_{ v \in E } \Big[
					\sum^\infty_{k=1}
						\|
							b(\cdot, v(\cdot)) \mathscr{P}_N e_k (\cdot) |\lambda_k|^{-s} 
						\|^2_H  
				\Big]\\
		\begin{split}
			\leq{} & 
				\sup_{ N \in \N }
				\sup_{ v \in E } \Big[
					\sum^\infty_{k=1}
						\|
							b(\cdot, v(\cdot)) \|^2_H \|\mathscr{P}_N e_k |\lambda_k|^{-s}
						\|^2_{L^\infty(0,1)} 
				\Big]	\\
			\leq{} & 
				\sup_{ v \in E } (
					\|
						b(\cdot, v(\cdot)) 
					\|^2_{H}
				)
				\cdot \sup_{ N \in \N }
				 \Big(
					\sum^\infty_{k=1}
						 \|\mathscr{P}_N e_k \|^2_{L^\infty(0,1)} (\pi k)^{-4s}
				\Big)	
		\end{split} \\  \nonumber
			\leq{} & 
				C_1 \sup_{ v \in E } (
					\|v\|^2_H + \int_0^1 |b(x,0)|^2 \ud x
				)
				\cdot 
					\sum^\infty_{k=1}
						 C (\pi k)^{-4s+2}
			< \infty.
	\end{align}
	This ensures \eqref{eq: B C1 bound Burger rem} 
	and thus finishes the proof of Remark \ref{rem: sufficient condition for B}.
\end{proof}
\section{Stochastic 2-D Navier-Stokes equations}
\label{sec: 2-D Navier-Stokes equations}
In this section we show that Kolmogorov backward equations
of stochastic 2-D Navier-Stokes equations
have a unique viscosity solution having at most polynomial growth.
\begin{sett}
\label{sett: Navier}
 Assume the setting in Section \ref{ssec: application},
 let $d =2$, $\eta \in (0,\infty)$,
 let $\C^{n, per} _{\R^2}([0,1]^2, \R^2)$, $n \in \N_0$, be the set satisfying 
 for all $n \in \N_0$ that
 \begin{equation}
	\begin{split}
			&\C^{n, per}_{\R^2}([0,1]^2, \R^2) 
		= 
			\big \{
				u \in \C^{n}_{\R^2}([0,1]^2, \R^2) \colon
				\forall i \in \{j \in \N^2_0 \colon |j| < n\}, \, \forall x \in (0,1) \colon \\
		& \qquad
					\big(
						(D^i_{\R^2} u)(x,0)=(D^i_{\R^2} u)(x,1)
					\big) 
					\wedge 
					\big(
						(D^i_{\R^2} u)(0,x)=(D^i_{\R^2} u)(1,x)
					\big) 
			\big\},
	\end{split}
 \end{equation}
 let 
	$H \subseteq \mathbb{L}^2((0,1)^2, \R^2)$ 
	be the set satisfying that
	$
			H 
		= 
			\{ 
				u \in L^2((0,1)^2, \R^2) \colon \Div(u) =0
			\}
	$
	and $\mathbb{H}= (H, \langle \cdot, \cdot \rangle_{H},  \|\cdot\|_{H} ) $ 
	the Hilbert space satisfying that
	$
			\mathbb{H}
		=
			(
				H, 
				\langle \cdot, \cdot \rangle_{L^2((0,1)^2, \R^2)}|_{H^2}, 
				\|\cdot\|_{L^2((0,1)^2, \R^2)}|_H
			) 
	$,
	let
	$\xi_k \in L^2(0,1)$, $k\in \Z$, be the function satisfying
	for all $x \in (0,1)$ and all $k \in \N$ that
	$\xi_k(x) = \sqrt{2} \sin(2k \pi x)$,
	that
	$\xi_0(x) = 1$,
	and that
	$\xi_{-k}(x) = \sqrt{2} \cos(-2k \pi x)$,
	let $\xi_{k,l} \in L^2((0,1)^2)$, $(k,l) \in \Z^2$, be the function satisfying
	for all $(k,l) \in \Z^2$ and all $(x,y) \in (0,1)^2$ that
	$\xi_{k,l} = \xi_k(x) \cdot \xi_{l}(y)$,
	let $e_{k,l} \in H$, $(k,l) \in \Z^2$, be the function satisfying 
	that
	\begin{equation}	
			e_{0,0} 
		\equiv 
			\left (
				\begin{array}{l} 
				1 \\ 0
				\end{array}
			\right)
	\end{equation}
	and that for all $(k,l) \in \Z^2\backslash \{(0,0)\}$ it holds
	that
	\begin{equation}
			e_{k,l}
		=
			\tfrac{1}{\sqrt{k^2+l^2}} 
				\left(
					\begin{array}{c}
						l \xi_{k,l} \\
						-k \xi_{-k,-l}
					\end{array}
				\right),
	\end{equation}
	let $\overline{e}_{0,0} \in H$ be the function satisfying 
	that
	\begin{equation}	
			\overline{e}_{0,0} 
		\equiv 
			\left (
				\begin{array}{l} 
				0 \\ 1
				\end{array}
			\right),
	\end{equation}
	let $\lambda_{k,l} \in (-\infty, 0]$, $(k,l) \in \Z^2$, be 
	the real numbers satisfying for all $(k,l) \in \Z^2$ that
	$\lambda_{k,l} = -(2\pi)^2 (k^2+ l^2)$,
	let
	$ \| \cdot  \|_{H_1} \colon B^2_{2,2}((0,1)^2,\R^2) \to [0,\infty)$
	be the norm satisfying that for all $(u_1,u_2) \in  B^2_{2,2}((0,1)^2,\R^2)$
	it holds that
	$
			\| (u_1,u_2) \|^2_{H_1}
		=
			\| \Delta u_1 \|^2_{L^2((0,1)^2)}+ \| \Delta u_2 \|^2_{L^2((0,1)^2)}
			+ \eta \| u_1 \|^2_{L^2((0,1)^2)}+ 
			\eta \| u_2 \|^2_{L^2((0,1)^2)}
	$,
	let $H_1 \subseteq  B^2_{2,2}((0,1)^2,\R^2)$ be the set satisfying that
	\begin{equation}
			H_1 
		= 
			\overline{
				\Span_{\mathbb{H}} 
					\big(
						\cup_{(k,l) \in \Z^2} \{ e_{k,l}\} \cup \{ \overline{e}_{0,0} \}
					\big)
			}_{\| \cdot \|_{H_1}},
	\end{equation}
	let $\Delta \colon H \supseteq D(\Delta) \to H$
	be the Laplace operator
	with periodic boundary condition,
	i.e.\@ the operator satisfying that
	$
			D(\Delta) 
		= 
			H_1
	$
	and that for all $(u_1,u_2) \in D(\Delta)$
	it holds that
	$\Delta(u_1,u_2) = (\Delta u_1,\Delta u_2)$,
	let 
		$ 
				\mathbb{H}_{t}
			=
				( 
					H_{t} , 
					\left< \cdot , \cdot \right>_{ H_{t} }, 
					\left\| \cdot \right\|_{ H_{t} } 
				) 
		$,
		$ t \in \R $,
		be a family of interpolation spaces associated with
		$ \eta-\Delta$ 
		(see, e.g., Definition~3.6.30 in \cite{Jentzen2015}).
		By abuse of notation we will also denote by
		$\Delta$ and  by $\| \cdot \|_{H_t}$, $t \in \R$,
		the extended operators 
		$
			\Delta \colon \bigcup_{i=1}^{\infty} H_{-i} \to \bigcup_{i=1}^{\infty} H_{-i}
		$ 
		and 
		$
			\| \cdot \|_{H_t} \colon \bigcup_{i=1}^{\infty} H_{-i} \to [0, \infty]
		$,
		$t \in \R$,
		satisfying for all
		$t \in \R$, $x \in H_{t}$, 
		and all $y\in \bigcup_{i=1}^{\infty} H_{-i}$
		that
		\begin{equation}
					\|y\|_{H_t}
				=
					\begin{cases} 
						\|y\|_{H_{t}} 
							& \textrm{ if } y \in H_{t} \\
						\infty 
							& \textrm{ if }y \notin H_{t}
					\end{cases}
		\end{equation}
		and that
		\begin{equation}
				(\Delta(x) = y) 
			\Leftrightarrow 
				(
						\lim_{\eps \downarrow 0} \sup \{
							\|\Delta(\xi) - y\|_{H_{t-1}} \colon \xi \in H_1, ~\|x- \xi\|_{H_{t}} \leq \eps 
						\}
					=
						0
				).
		\end{equation}
	Let
	$\alpha \in (\nicefrac 34, 1)$,
	$\beta \in [0,\nicefrac 12)$,
	$
		c \in \R \backslash \{ 0 \},
	$
	$ 
		\gamma \in 
		[ 1-\alpha, \nicefrac 14 )
	$,
	let $\mathbb{U}=\mathbb{L}^2((0,1)^2,\R^2)$
	be the $L^2$-space,
	let $\mathcal{I} \subseteq \N$,
	let $\tilde{e}_n$, $n \in \mathcal{I}$, be an orthonormal basis of $\U$,
	let
	\begin{equation}
		B \in 
		\C_{
			\mathbb{H}_\gamma, 
			\mathbb{HS}(\mathbb{U}, \mathbb{H}_{\gamma-\beta})
		}( 
			H_\gamma, 
			HS( \mathbb{U}, \mathbb{H}_{\gamma-\beta} ) 
		),
	\end{equation}
	denote by 
	\begin{equation}
		F \in \C_{\mathbb{H}_\gamma, \mathbb{H}_{\gamma - \alpha}}(H_\gamma, H_{\gamma-\alpha})
	\end{equation}
	the function satisfying for all $u \in H_{1}$
	that
	\begin{equation}
			F(u) (x)
		=
			c \pi^{L^2((0,1)^2,\R^2)}_{H}
				\big( (D_{\R^2 } u) u \big)
			+ \eta u.
	\end{equation}
\end{sett}
\begin{lemma}
\label{l: norm equivalence Navier}
	Assume
	the Setting \ref{sett: Navier} and let $r \in [0,\infty)$.
	Then there exist $C_1, C_2 \in (0,\infty)$
	such that for all $u \in H_r$ it holds that
	\begin{equation}
			C_1 \|u\|_{B^{2r}_{2,2}((0,1)^2,\R^2)} 
		\leq 
			\|u\|_{H_r} 
		\leq 
			C_2 \|u\|_{B^{2r}_{2,2}((0,1)^2,\R^2)}.
	\end{equation}
\end{lemma}
\begin{proof}
	First note, that it holds for all 
	$(u_1,u_2),(v_1,v_2) \in \C^{1,per}_{\R^2}([0,1]^2, \R^2)$
	that
	\begin{align}
	\nonumber
				&\langle (u_1,u_2) , D_{\R^2}^{(1,0)} (v_1,v_2) \rangle_{L^2((0,1)^2,\R^2)} \\ \nonumber
			={}
				&\int_0^1 \int_0^1 u_1(x,y) \, (D_{\R^2}^{(1,0)} v_1) (x,y) \ud x \ud y
				+\int_0^1 \int_0^1 u_2(x,y) \, (D_{\R^2}^{(1,0)} v_2) (x,y) \ud x \ud y \\
			={}
				&\int_0^1 u_1(0,y) v_1(0,y) -u_1(1,y) v_1(1,y) \ud y
				+\int_0^1 u_2(0,y) v_2(0,y) -u_2(1,y) v_2(1,y) \ud y \\ \nonumber
				&-\int_0^1 \int_0^1 (D_{\R^2}^{(1,0)} u_1) (x,y) \, v_1(x,y) \ud x \ud y
				-\int_0^1 \int_0^1 (D_{\R^2}^{(1,0)} u_2) (x,y) \, v_2(x,y) \ud x \ud y \\ \nonumber 
			={}
				&-\langle D_{\R^2}^{(1,0)} (u_1,u_2) , (v_1,v_2) \rangle_{L^2((0,1)^2,\R^2)}.
	\end{align}
	Analogously it follows that for all 
	$(u_1,u_2),(v_1,v_2) \in \C^{1,per}_{\R^2}([0,1]^2, \R^2)$
	it holds
	that
	\begin{equation}
				\langle (u_1,u_2) , D_{\R^2}^{(0,1)} (v_1,v_2) \rangle_{L^2((0,1)^2, \R^2)}
			=
				-\langle D_{\R^2}^{(0,1)} (u_1,u_2) , (v_1,v_2) \rangle_{L^2((0,1)^2, \R^2)}.
	\end{equation}
	Therefore we get for all $u \in \C^{3,per}_{\R^2}([0,1]^2, \R^2)$
	that
	\begin{equation}
		\begin{split}
				&\|\Delta u \|_{L^2((0,1)^2, \R^2)}^2
			=
				\langle 
					(D_{\R^2}^{(2,0)} + D_{\R^2}^{(0,2)}) u, 
					(D_{\R^2}^{(2,0)} + D_{\R^2}^{(0,2)}) u 
				\rangle_{L^2((0,1)^2, \R^2)} \\
			={}
				&\| D_{\R^2}^{(2,0)} u\|_{L^2((0,1)^2, \R^2)}^2
				+ \| D_{\R^2}^{(0,2)} u\|_{L^2((0,1)^2, \R^2)}^2 
				+2\langle 
					D_{\R^2}^{(1,0)} D_{\R^2}^{(1,0)}u, D_{\R^2}^{(0,1)} D_{\R^2}^{(0,1)} u 
				\rangle_{L^2((0,1)^2, \R^2)}\\
			={}
				&\| D_{\R^2}^{(2,0)} u\|_{L^2((0,1)^2, \R^2)}^2
				+ \| D_{\R^2}^{(0,2)} u\|_{L^2((0,1)^2, \R^2)}^2 
				+2\langle 
					D_{\R^2}^{(0,1)} D_{\R^2}^{(1,0)}u, D_{\R^2}^{(0,1)} D_{\R^2}^{(1,0)} u 
				\rangle_{L^2((0,1)^2, \R^2)} \\
			={}
				&\| D_{\R^2}^{(2,0)} u\|_{L^2((0,1)^2, \R^2)}^2
				+ \| D_{\R^2}^{(0,2)} u\|_{L^2((0,1)^2, \R^2)}^2 
				+2\|D_{\R^2}^{(1,1)} u \|^2_{L^2((0,1)^2, \R^2)} \\
			\geq{}
				&\sum_{i \in \N_0^2, |i|=2}
					\| D_{\R^2}^{i} u\|_{L^2((0,1)^2, \R^2)}^2.
		\end{split}
	\end{equation}
	Hence, we obtain that for all 
	$n \in \N_0$ there exists a $C \in (0,\infty)$ such that for all
	$u \in \C^{2n+1,per}_{\R^2}([0,1]^2, \R^2)$ it holds
	that
	\begin{equation}
				C \|\Delta^n u \|_H
			\geq
				\sum_{i \in \N_0^{2}, |i|=2n}
					\| D_{\R^2}^{i} u\|_{L^2((0,1)^2, \R^2)}^2.
	\end{equation}
	Combining this and the fact that for all $n \in \N$ it holds that
	$H_n \cap \C^{2n+1,per}_{\R^2}([0,1]^2, \R^2)$ is dense in $H_n$ shows that 
	for all $n \in \N_0$ there exists a $C \in (0,\infty)$ such that for all
	$u \in H_n$ it holds
	that
	\begin{equation}
				C \|\Delta^n u \|_H
			\geq
				\sum_{i \in \N_0^{2}, |i|=2n}
					\| D_{\R^2}^{i} u\|_{L^2((0,1)^2, \R^2)}^2.
	\end{equation}
	Next note that we obtain from (2.200) on page 54 in Mitrea \cite{Mitrea2013}
	and from
	Theorem 1.8 in Nečas \cite{Necas2012} that for all $n \in \N_0$
	there exist $C_1,C_2,C_3,C_4, C_5 \in (0,\infty)$ such that
	for all $u \in H_{n}$
	it holds that
	\begin{align}
		\nonumber
				&\|u\|_{B^{2n}_{2,2}((0,1)^2,\R^2)}
			\leq
				C_1
					\sum_{i\in \N_0^2, |i| \leq 2n} \|D_{\R^2}^{i} u\|_{L^2((0,1)^2,\R^2)} \\ 
			\leq{}
				&C_2
					(
						\|u\|_{L^2((0,1)^2,\R^2)}
						+ \sum_{i\in \N_0^2, |i| = 2n} \|D_{\R^2}^{i} u\|_{L^2((0,1)^2,\R^2)}
					) \\ \nonumber
			\leq{}
				&C_3 (\|u\|_{H} + \|\Delta^n u \|_{H})
			\leq
				C_3 ( \eta^{-n} \|u\|_{H_n} + \|u \|_{H_n}) 
			\leq
				C_3 (\eta^{-n} +1 ) \|u\|_{H_n} \\
		\nonumber
			\leq{}
				&C_4
						\sum_{i \in \N^2_0, |i| \leq 2n} \|D_{\R^2}^{i} u\|_{L^2((0,1)^2,\R^2)}
			\leq
				C_5 \|u\|_{B^{2n}_{2,2}((0,1)^2,\R^2)}.
	\end{align}
	Combining this and Lemma \ref{l: norm equivalence with interpolation}
	finishes the proof of Lemma \ref{l: norm equivalence Navier}.
\end{proof}
The next corollary applies Corollary \ref{cor:uniqueness2}
and Theorem \ref{thm: existence} to Kolmogorov backward equations
of stochastic 2-D Navier-Stokes equations.
\begin{corollary}
\label{cor: Navier equation}
	Assume
	the Setting \ref{sett: Navier},
	let
	$\vartheta \in [\nicefrac 12, \nicefrac 34 -\gamma)$,
	$\alpha_1 \in (0,\gamma)$,
	$T , \theta, \beta_1, \beta_2 \in (0,\infty)$,
	$\beta_3 \in (0, 1)$,
	$ \chi \in [ \beta, \nicefrac{ 1 }{ 2 } ) $,
	let $\mathscr{P}_0 \in L(\mathbb{U},\mathbb{U})$ be the function satisfying that
	$\mathscr{P}_0=\id_U$,
	let $\mathscr{P}_N \in L(\mathbb{U},\mathbb{U})$, $N \in \N$,
	be finite-dimensional projection,
	let $\varphi \in \C_{\mathbb{H}_\gamma}(H_\gamma, \R)$ 
	have at most polynomial $\H$-growth,
	let $K \colon \R \to (0, \infty)$ be increasing,
	and assume that
	\begin{equation}
	\label{eq: B continuity Navier}
		B|_{H_{1/2 +\gamma +\vartheta - \beta_3}} \in
				\C_{
					\mathbb{H}_{1/2 +\gamma +\vartheta - \beta_3},
					\mathbb{HS}(\mathbb{U}, \mathbb{H}_{\gamma+\vartheta})
				}(
					H_{1/2+\gamma +\vartheta - \beta_3},
					HS(\mathbb{U}, \mathbb{H}_{\gamma+\vartheta})
				), 
	\end{equation}
	that for all $u, v \in H_{1/2}$ 
	it holds that
	\begin{align}
		\label{eq: H bound B Navier}
			&\| B(u)-B(v) \|^2_{HS(\mathbb{U},\mathbb{H})}
		\leq
			K(1) \cdot (1+\|u\|^4_{L^4((0,1)^2,\R^2)}) \|u-v\|^2_{H},
	\end{align}
	that for all $u, v \in H_{\gamma+2\vartheta}$ 
	it holds that
	\begin{align}
	\label{eq: B bounded wrt H_vartheta Navier}
			&\| B(u) \|^2_{HS(\mathbb{U},\mathbb{H}_{\gamma+\vartheta})} 
		\leq
			K(\|u\|_{H_\gamma}) 
				(\|u\|^2_{H_{1/2+\gamma+\vartheta- \beta_1}} + 1), \\
	\label{eq: B Lipschitz wrt H_1/2 Navier}
			&\|B(u) - B(v)\|_{HS(\mathbb{U}, \mathbb{H}_\gamma)}
		\leq
			K(\|v\|_{H_\gamma} \vee \|u \|_{H_\gamma}) 
				\|u-v\|_{H_{1/2+\gamma -\beta_2}}
	\end{align}
	and that for all $\mathbb{H}_\gamma$-bounded sets
	$ E \subseteq H_\gamma$
	it holds that
	\begin{equation}
	\label{eq: Lipschitz bound B Navier}
		| (B|_E) |_{
			\C^1( 
				E, 
				\left\| \cdot \right\|_{HS( \mathbb{U}, \mathbb{H}_{\gamma- \beta} )} 
			)
		}
		< \infty,
	\end{equation}
	\begin{equation}
	\label{eq: B C1 bound Navier} 
			\sup_{ N \in \N } \sup_{ v \in E }
			\left[
				\frac{
				\|
					B( v ) \mathscr{P}_N 
				\|_{ HS( \mathbb{U}, \mathbb{H}_{\gamma- \beta } ) }
				}{
					1 + \| v \|_{ H_\gamma }
				}
			\right]
			< \infty,
	\end{equation}
	and that
	\begin{equation}
	\label{eq:convLocLip_assumption Navier}
			\limsup_{ N \to \infty }
			\sup_{ v \in E }
			\left[
			\frac{
			\|
				B( v ) ( \id_H - \mathscr{P}_N )
			\|_{ HS( \mathbb{U}, \mathbb{H}_{\gamma- \chi} ) }
			}{
				1 + \| v \|_{ H_\gamma }
			}
			\right]
		=0,
	\end{equation}
	let $\Z^2_N$, $N \in \N$, be the sets satisfying for all $N \in \N$ that
		$\Z^2_N =  \{ (l,k) \in \Z^2 \colon |l|+|k|\leq N \}$,
		let $V_N \subseteq H_\gamma$, $N \in \N_0$, be the linear subspaces 
		satisfying that
		$V_0= H_\gamma$ and that for all $N \in \N$ it holds that
		$
				V_N 
			= 
				\Span_{\mathbb{H}}(
					\{
						e_{k,l} 
						\colon (l,k) \in \Z_N^2
					\} 
					\cup 
					\{
						\overline{e}_{0,0}
					\}
				)
		$, 
		let $A \colon H_1 \to H$ be the operator satisfying for all 
		$u \in H_1$ that
		$A u = \Delta u -\eta u$,
		let
		$ ( \Omega, \mathcal{F}, \P ) $
		be a probability space with a normal filtration 
		$ ( \mathbb{F}_t )_{ t \in [0,T] } $, 
		let 
		$
			( W_t )_{ t \in [0,T] } 
		$ 
		be an $ \operatorname{Id}_U $-cylindrical $ ( \mathbb{F}_t )_{ t \in [0,T] } $-Wiener
		process and let 
	$
		X^{N,u} \colon [0,T] \times \Omega \rightarrow H_\gamma
	$,
	$ N \in \N_0 $,
	$u \in H_\gamma$,
	be
	$
		(\mathbb{F}_t )_{ t \in [0,T] }
	$-adapted stochastic processes 
	with continuous sample paths satisfying
	for all $ N \in \N_0 $, $ t \in [0,T] $, $u \in H_\gamma$,
	and all $ \eps \in (0,\infty)$
	that
	\begin{equation}
			X_t^{N,u}
		= 
			e^{ t A } \pi^{H_\gamma}_{V_N} u
			+
			\int_0^t
				e^{ ( t - s ) A }
				\pi^{H_{\gamma-\alpha}}_{V_N} F( X_t^{N,u}  )
			\ds
			+
			\int_0^t
				e^{ ( t - s ) A }
				\pi^{H_{\gamma-\beta}}_{V_N} B( X_t^{N,u}  ) \mathscr{P}_N
			\dWs,
	\end{equation}
	and that
	\begin{equation}
	\label{eq: X continuous in start point Navier}
		\begin{split}
			&\limsup_{\delta \to 0}  \sup_{v \in H_\gamma, \| v- u \|_{H_\gamma}\leq \delta}  
				\sup_{r \in [0,T]}
					\P( \| X^{0,u}_r - X^{0,v}_r \|_{H_\gamma} \geq \eps)
			=0,
		\end{split}
	\end{equation}
	and let $f \colon [0,T] \times H_\gamma \to \R$
	satisfy that for all $(t,x) \in [0,T] \times H_\gamma$
	it holds that
	$
		f(t,x) = \E[\varphi(X^{0,x}_t)]
	$. 
	Then
	$ 
    f|_{ (0,T) \times H_\gamma} 
  $
  is the unique viscosity solution of
  \begin{equation} 
			\tfrac{ \partial }{ \partial t }
			f(t,u) -
			\big\langle
				A u+ F(u), 
				I^{-1}_{\mathbb{H}_\gamma} (D_{\mathbb{H}_\gamma} f)(t,u)
			\big\rangle_{H_\gamma} 
			-
			\big \langle 
				B(u),
				I_{\mathbb{H}_{\gamma+\vartheta}}^{-1}
					\big( (D^2_{\mathbb{H}_\gamma} \, f)(t,u) \, B(u) \big)
			\big \rangle_{HS(\mathbb{U},\mathbb{H}_{\gamma+\vartheta})} 
			= 0
  \end{equation}
	for $ (t,u) \in (0,T) \times H_\gamma $ relative to 
	$
		(
			(0,T) \times H_\gamma \ni (t,u) 
			\to \nicefrac 12 \|u \|^2_{\mathbb{H}_{\gamma+\vartheta}} \in [0, \infty],
			\R \times \mathbb{H}_\gamma, 
			\R \times \mathbb{H}_{\gamma+\vartheta}
		)
	$
	which satisfy that 
	$ 
    f \in \C_{\R \times \mathbb{H}_\gamma}([0,T] \times H_\gamma, \R)
  $
	that $f$ is
	bounded on $\R \times \mathbb{H}_\gamma$-bounded 
	subsets of $[0,T] \times H_\gamma$,
	that $f$ have at most polynomial $\H$-growth,
	and that for all $ u \in H_\gamma $ it holds that
  $
    f(0, u) = \varphi(u)
  $.
	\end{corollary}
	\begin{proof}
		First note that the existence of
		$X^{N,x}$, $N \in \N_0$, $x \in H_\gamma$, 
		follows from Example 3.3 in Liu \& R\"ockner \cite{LiuRockner2010}.
		In addition, 
		the fact that $\varphi$ have at most polynomial $\H$-growth demonstrates that
		there exists a $p \in [2,\infty)$, which we fix for the rest of the proof,
		such that 
		$
			\lim_{r \to \infty}
				\sup \big\{ 
					\tfrac{\varphi(u)}{\|u\|_{H^p}} \colon u \in H_\gamma, \|u\|_{H} \geq r 
				\big\}
			=0.
		$
		Denote by 
		$ V \colon H_\gamma \to (1,\infty)$ the function satisfying for all
		$u \in H_\gamma$ that
		$V(u) = \|u\|^p_{H} + 2$.
		Next we will apply Theorem \ref{thm: existence} 
			(with $\mathbb{H} \leftarrow \mathbb{H}_\gamma$)
		to verify that $f$ is a viscosity solution.
	Therefore note that we obtain from the Sobolev inequality
	(see, e.g., the Theorem on page 31 in Runst \& Sickel \cite{RunstSickel1996}),
	from the fact that $\gamma \geq 1 -\alpha$, and from
	Lemma \ref{l: norm equivalence Navier}
	that there exist $C_1,C_2,C_3, C_4, C_5 \in (0,\infty)$ such that for all
	$w\in H_{\alpha-\gamma}$ it holds that
	\begin{equation}
		\begin{split}
				&\big\|
					\| D_{\R^2} w \|_{L(\R^2,\R^2)}
				\big \|_{L^{1/(2\gamma)}((0,1)^2)}
			\leq
				\big \| 
					\|D^{(1,0)}_{\R^2} w\|_{\R^2} + \|  D^{(0,1)}_{\R^2} w\|_{\R^2}
				\big \|_{L^{1/(2\gamma)}((0,1)^2)} \\
			\leq{}
				&C_1\big(
					\|D^{(1,0)}_{\R^2} w\|_{L^{1/(2\gamma)}((0,1)^2,\R^2)} 
					+ \|  D^{(0,1)}_{\R^2} w\|_{L^{1/(2\gamma)}((0,1)^2,\R^2)} 
				\big)\\
			\leq{}
				&C_2\big(
					\|
						D^{(1,0)}_{\R^2} w
					\|_{B^{1-4\gamma}_{2,2}((0,1)^2,\R^2)}
					+
					\|
						D^{(0,1)}_{\R^2} w
					\|_{B^{1-4\gamma}_{2,2}((0,1)^2,\R^2)}
				\big)
			\leq
				C_3 \|w\|_{B^{2-4\gamma}_{2,2}((0,1)^2,\R^2)}\\
			\leq{}
				&C_4 \|w\|_{B^{2\alpha -2\gamma}_{2,2}((0,1)^2,\R^2)} 
			\leq
				C_5 \|w\|_{H_{\alpha - \gamma}}.
		\end{split}
	\end{equation}
	Moreover, it holds for all
	$u=(u_{(1,0)},u_{(0,1)})\in H_{1/2}$,
	$v=(v_{(1,0)},v_{(0,1)})\in H_{1/2}$,
	and all
	$w=(w_{(1,0)},w_{(0,1)})\in H_{1/2}$
	that
	\begin{align}
	\nonumber
					&\int_0^1 \int_0^1
						\Big \langle
						 (D_{\R^2}u)(x,y) (v(x,y)),
							w(x,y)
						\Big \rangle_{\R^2}
					\ud x \ud y \\ \nonumber
				={}
					&\int_0^1 \int_0^1
						\sum_{\alpha, \beta \in \N_0^2, |\alpha|=|\beta|=1}
							(D^{\beta}_{\R^2}u_{\alpha})(x,y) \, v_{\beta}(x,y) \, w_{\alpha }(x,y)
					\ud x \ud y \\
	\label{eq: interchanging diff.}
		\begin{split}
				={}
					&-\int_0^1 \int_0^1
						\bigg(
							\sum_{\alpha \in \N_0^2, |\alpha|=1}
								u_{\alpha}(x,y)  \, w_{\alpha}(x,y)
							\sum_{\beta \in \N_0^2, |\beta|=1}
								(D^{\beta}_{\R^2} v_{\beta}) (x,y)
						\bigg)
						\\ & \qquad
							+\sum_{\alpha, \beta \in \N_0^2, |\alpha|=|\beta|=1}
								u_{\alpha}(x,y) \, v_{\beta} (x,y) \, (D^{\beta}_{\R^2} w_{\alpha})(x,y)
					\ud x \ud y 
		\end{split} \\ \nonumber
				={}
					&-\int_0^1 \int_0^1
						\sum_{\alpha, \beta \in \N_0^2, |\alpha|=|\beta|=1}
								u_{\alpha}(x,y) \, v_{\beta} (x,y) \, (D^{\beta}_{\R^2} w_{\alpha})(x,y)
					\ud x \ud y \\ \nonumber
			={}
				&-\int_0^1 \int_0^1
						\Big \langle
						 (D_{\R^2}w)(x,y) (v(x,y)),
							u(x,y)
						\Big \rangle_{\R^2}
					\ud x \ud y.
	\end{align}
	Hence the Sobolev inequality
	(see, e.g., the Theorem on page 31 in Runst \& Sickel \cite{RunstSickel1996}),
	and Lemma \ref{l: norm equivalence Navier} show  that
	there exist $C_1, C_2, C_3 \in (0,\infty)$ such that for all
	$u,v \in H_1$ it holds that 
	\begin{align}
	\label{eq: F loc Lipschitz Navier}
				&\|F(u)-F(v)\|_{H_{\gamma-\alpha}} \\
		\nonumber
			\leq{}
				&\sup_{w \in H_{\alpha-\gamma}} \bigg(
					\int_0^1 \int_0^1
						c\Big \langle
								\big( (D_{\R^2}u)(x,y) \big) (u(x,y))
								-\big( (D_{\R^2}v)(x,y) \big) (v(x,y)),
							\tfrac {w(x,y)}{\|w\|_{H_{\alpha-\gamma}}}
						\Big \rangle_{\R^2}
					\ud x \ud y
				\bigg) 
			\\ & \nonumber \qquad
				+ \eta \|u-v\|_{H_{\gamma-\alpha}}\\
			={} \nonumber
				&\sup_{w \in H_{\alpha-\gamma}} \bigg(
					-\int_0^1 \int_0^1
						c\Big \langle
							\big( (D_{\R^2}w)(x,y) \big) (u(x,y)), 
							\tfrac {u(x,y)}{\|w\|_{H_{\alpha-\gamma}}}
						\Big \rangle_{\R^2}
				\\ & \nonumber \qquad
						-c\Big \langle
							\big( (D_{\R^2}w)(x,y) \big) (v(x,y)), 
							\tfrac {v(x,y)}{\|w\|_{H_{\alpha-\gamma}}}
						\Big \rangle_{\R^2}
					\ud x \ud y
				\bigg) 
				+ \eta \|u-v\|_{H_{\gamma-\alpha}}\\
			={} \nonumber
				&\sup_{w \in H_{\alpha-\gamma}} \bigg(
					-\int_0^1 \int_0^1
						c\Big \langle
							\big( (D_{\R^2}w)(x,y) \big) (u(x,y)-v(x,y)), 
							\tfrac {u(x,y)}{\|w\|_{H_{\alpha-\gamma}}}
						\Big \rangle_{\R^2}
				\\ & \nonumber \qquad
						+c\Big \langle
							\big( (D_{\R^2}w)(x,y) \big) (v(x,y)), 
							\tfrac {u(x,y)-v(x,y)}{\|w\|_{H_{\alpha-\gamma}}}
						\Big \rangle_{\R^2}
					\ud x \ud y
				\bigg) 
				+ \eta \|u-v\|_{H_{\gamma-\alpha}}\\
			\leq{} \nonumber
				&\sup_{w \in H_{\alpha-\gamma}} \bigg(
					\int_0^1 \int_0^1
						|c| \cdot
							\| u(x,y)-v(x,y)\|_{\R^2}
							(\| u(x,y)\|_{\R^2} +\| v(x,y)\|_{\R^2})
				\\ & \nonumber \qquad\qquad \cdot
							\tfrac {\| D_{\R^2}w (x,y) \|_{L(\R^2, \R^2)}}{\|w\|_{H_{\alpha-\gamma}}}
					\ud x \ud y
				\bigg) 
				+ \eta \|u-v\|_{H_{\gamma-\alpha}}\\
			\leq{} \nonumber
				&\sup_{w \in H_{\alpha-\gamma}} \bigg(
					|c| 
					\big \| 
						\|u-v\|_{\R^2} (\| u \|_{\R^2} + \|v \|_{\R^2})
					\big \|_{L^{1/(1-2\gamma)}((0,1)^2)} 
					\cdot 
					\tfrac 
						{
							\|
								\| D_{\R^2}w\|_{L(\R^2, \R^2)} 
							\|_{L^{1/(2\gamma)} ((0,1)^2)}
						}
						{\|w\|_{H_{\alpha-\gamma}}}
				\bigg) 
		\\ & \nonumber \qquad
				+ \eta \|u-v\|_{H_{\gamma-\alpha}} \\
			\leq{} \nonumber
				&C_1 \|u-v\|_{L^{2/(1-2\gamma)}((0,1)^2,\R^2)} 
				\cdot 
					(\|u\|_{L^{2/(1-2\gamma)}((0,1)^2,\R^2)} 
					+\|v\|_{L^{2/(1-2\gamma)}((0,1)^2,\R^2)}) 
				+ \eta \|u-v\|_{H_{\gamma-\alpha}} \\
			\leq{} \nonumber
				&C_2 \|u-v\|_{B_{2,2}^{2\gamma}((0,1)^2,\R^2)} 
				\cdot (
					\|u\|_{B_{2,2}^{2\gamma}((0,1)^2,\R^2)}
					+\|v\|_{B_{2,2}^{2\gamma}((0,1)^2,\R^2)}
				)
				+ \eta \|u-v\|_{H_{\gamma-\alpha}} \\
			\leq{} \nonumber
				&C_3 \|u-v\|_{H_\gamma} 
				\cdot (
					\|u\|_{H_{\gamma}}
					+\|v\|_{H_{\gamma}}
				)
				+ \eta^{1-\alpha} \|u-v\|_{H_\gamma} \\
	\nonumber
			={}
				&C_3 \|u-v\|_{H_\gamma} 
				\cdot( \|u\|_{H_{\gamma}}+\|v\|_{H_{\gamma}} +\eta^{1-\alpha})
	\end{align}
	and this 
	together with the fact that
	$H_1$ is dense in $H_\gamma$ and
	\eqref{eq: Lipschitz bound B Navier}
	verifies \eqref{eq: Lipschitz bound for F and B}
		(with $H \leftarrow H_\gamma$ and $\gamma \leftarrow 0$). 
	Furthermore, we have for all $N \in \N_0$ and 
	all $u \in H_{\gamma+ 2 \vartheta} \cap V_N$ that
	\begin{equation}
		\begin{split}
				&\langle 
					\pi^{H_{\gamma-\alpha}}_{V_N} ( F(u) -\eta u ), 
					 (D_{\mathbb{H}_\gamma} V) (u) 
				\rangle_{H_\gamma,H'_\gamma}
			=
				c p \|u\|^{p-2}_{H} \,
				\langle 
					\pi^{L^2((0,1)^2,\R^2)}_{V_N} \big( ( D_{\R^2} u) u \big), 
					u 
				\rangle_{L^2((0,1)^2,\R^2)} \\
			={}
				&c p \|u\|^{p-2}_{H} \,
				\int_0^1 \int_0^1
						\sum_{\alpha, \beta \in \N_0^2, |\alpha|=|\beta|=1}
							(D^{\alpha}_{\R^2}u_{\beta})(x,y) \, u_{\alpha}(x,y) \, u_{\beta}(x,y)
					\ud x \ud y \\
			={}
				&\tfrac{c p}{2} \|u\|^{p-2}_{H}
				\int_0^1 \int_0^1
						\sum_{\alpha, \beta \in \N_0^2, |\alpha|=|\beta|=1}
							(D^{\alpha}_{\R^2} (u_{\beta}^2))(x,y) \, u_{\alpha}(x,y)
					\ud x \ud y \\
			={}
				&-\tfrac{c p}{2} \|u\|^{p-2}_{H}
				\int_0^1 \int_0^1
						\sum_{\beta \in \N_0^2, |\beta|=1}
							(u_{\beta}(x,y))^2 
							\sum_{\alpha\in \N_0^2, |\alpha|=1}
							(D^{\alpha}_{\R^2} u_{\alpha})(x,y)
					\ud x \ud y
			=0.
		\end{split}
	\end{equation}
	Thus it follows from \eqref{eq: H bound B Navier} 
		(with $u \leftarrow 0$ and $v \leftarrow u$)
	that 
	for all $N \in \N_0$, $u \in H_{\gamma +2\vartheta} \cap V_N$ it holds that 
	\begin{align}
	\label{eq: V is lypanov Navier}
				&\big \langle
					\pi^{H_{\gamma-\alpha}}_{V_N} F(u) + A u,
					(D_{\mathbb{H}_\gamma} V)(u)
				\big \rangle_{H_\gamma,H_\gamma'} 
		\\ & \nonumber
				+
				\nicefrac 12
				\big \langle 
					\pi^{H_{\gamma-\beta}}_{V_N} B(u),
					I_{\mathbb{H}_{\gamma+\vartheta}}^{-1} 
						\big( 
							(D_{\mathbb{H}_\gamma}^2 \, V)(u) \, 
							\pi^{H_{\gamma-\beta}}_{V_N} B(u) 
						\big)
				\big \rangle_{HS(\mathbb{U},\mathbb{H}_{\gamma+\vartheta})} \\
			={}  \nonumber
				&p \|u\|_H^{p-2}
				\big \langle
					(A+ \eta) u, u 
				\big \rangle_{H}
		\\ & \nonumber
				+ p (p-2) \|u\|_H^{p-4}
				\big \langle 
					\pi^{H_{\gamma-\beta}}_{V_N} B(u),
					(-A)^{-2 \vartheta-2\gamma} u \, 
						\langle \pi^{H_{\gamma-\beta}}_{V_N} B(u), u \rangle_H
				\big \rangle_{HS(\mathbb{U},\mathbb{H}_{\gamma+\vartheta})} \\ \nonumber
				&+p \|u\|_H^{p-2}
				\big \langle 
					\pi^{H_{\gamma-\beta}}_{V_N} B(u),
					(-A)^{-2 \vartheta-2\gamma} \, \pi^{H_{\gamma-\beta}}_{V_N} B(u)
				\big \rangle_{HS(\mathbb{U},\mathbb{H}_{\gamma+\vartheta})} \\
			\leq{}  \nonumber
				&-p \|u\|_H^{p-2}
				\|
					u
				\|^2_{H_{1/2}}
				+p \eta \|u\|_H^{p}
				+ p (p-2) \|u\|_H^{p-4}
				\big \langle 
					\pi^{H_{\gamma-\beta}}_{V_N} B(u),
					 u \, \langle \pi^{H_{\gamma-\beta}}_{V_N} B(u), u \rangle_H 
				\big \rangle_{HS(\mathbb{U},\mathbb{H})}
			\\ \nonumber & \qquad
				+ p \|u\|_H^{p-2}
				\| 
					\pi^{H_{\gamma-\beta}}_{V_N} B(u)
				\|^2_{HS(\mathbb{U},\mathbb{H})} \\
			\leq{}  \nonumber
				&p \eta \|u\|_H^{p}
				+p (p-2) \|u\|_H^{p-2}
				\|
					\pi^{H_{\gamma-\beta}}_{V_N} B(u)
				\|^2_{HS(\mathbb{U},\mathbb{H})}
				+ p \|u\|_H^{p-2}
				\| 
					\pi^{H_{\gamma-\beta}}_{V_N} B(u)
				\|^2_{HS(\mathbb{U},\mathbb{H})} \\
			\leq{}  \nonumber
				&p \eta \|u\|_H^{p}
				+2 p (p-1) \|u\|_H^{p-2}
				\big(
					\|
						B(u)-B(0)
					\|^2_{HS(\mathbb{U},\mathbb{H})} 
					+\|
						B(0)
					\|^2_{HS(\mathbb{U},\mathbb{H})} 
				\big)\\
			\leq{}  \nonumber
				&p \eta \|u\|_H^{p}
				+2 \, p (p-1) \|u\|_H^{p-2}
				(
					K(1) \, \|u\|^2_{H} 
					+\|
						B(0)
					\|^2_{HS(\mathbb{U},\mathbb{H})}
				) \\
			\leq{} \nonumber
				&p \eta (1+\|u\|_H^{p})
				+4 \, p (p-1) 
				\big(K(1) \vee \| B(0) \|^2_{HS(\mathbb{U},\mathbb{H})} \big)
				(1+\|u\|_H^{p}) \\
			\leq{} 
	\nonumber
				&\Big(
					p \eta
					+4 \, p (p-1) 
					\big(K(1) \vee \| B(0) \|^2_{HS(\mathbb{U},\mathbb{H})} \big)
				\Big)
				V(u).
	\end{align}
	Next note that for all $(k,l), (m,n) \in \Z^2 \backslash \{ (0,0) \}$
	it holds that 
	\begin{equation}
	\label{eq: F ohne 00}
		\begin{split}
				(D_{\R^2} e_{k,l}) e_{m,n}
			={}
				&\tfrac{2 \pi}{\sqrt{k^2+l^2}}
				\left(
					\begin{array}{cc}
						k \, l \, \xi_{-k,l} & l^2 \xi_{k,-l} \\
						k^2 \xi_{k,-l} & k \, l \, \xi_{-k,l}
					\end{array}
				\right)
				e_{m,n} \\
			={}
				&\tfrac{2 \pi}{\sqrt{k^2+l^2}}\tfrac{1}{\sqrt{m^2+n^2}}
				\left(
					\begin{array}{c}
						l \, k \, n \, \xi_{-k,l} \xi_{m,n} - l^2 \, m \xi_{k,-l} \xi_{-m,-n} \\
						k^2  \, n  \, \xi_{k,-l} \xi_{m,n} - l \, k \, m \, \xi_{-k,l} \xi_{-m,-n}
					\end{array}
				\right),
		\end{split}
	\end{equation}
	that
	\begin{equation}
	\label{eq: F 0}
		\begin{split}
				(D_{\R^2} e_{k,l}) e_{0,0}
			=
				\tfrac{2\pi}{\sqrt{k^2+l^2}}
				\left(
					\begin{array}{c}
						k \, l \, \xi_{-k,l} \xi_{0,0}\\
						k^2 \xi_{k,-l} \xi_{0,0}
					\end{array}
				\right),
		\end{split}
	\end{equation}
	and that
	\begin{equation}
	\label{eq: F 00}
		\begin{split}
				(D_{\R^2} e_{k,l}) \overline{e}_{0,0}
			=
				\tfrac{2 \pi}{\sqrt{k^2+l^2}}
				\left(
					\begin{array}{c}
						l^2 \xi_{k,-l} \xi_{0,0} \\
						k \, l \, \xi_{-k,l} \xi_{0,0}
					\end{array}
				\right).
		\end{split}
	\end{equation}
	Moreover, the product to sum identity implies that for all
	$k,l,m \in \Z$ it holds that
	\begin{equation}
	\label{eq: equation cos cos cos}
		\begin{split}
				&4\int_0^1 
					\cos ( 2 \pi \, k \, x) \cos ( 2 \pi \, l \, x) \cos ( 2 \pi \, m \, x) 
				\ud x \\
			={}
				&2
				\int_0^1 
					\big(
						\cos ( 2 \pi \, x (k-l)) + \cos ( 2 \pi \, x (k+l)) 
					\big)
					\cos ( 2 \pi \, m \, x) \ud x \\
			={}
				&\int_0^1 
					\cos ( 2 \pi \, x (k-l-m)) + \cos ( 2 \pi \, x (k-l+m))
			\\ & \quad
					+\cos ( 2 \pi \, x (k+l-m)) + \cos ( 2 \pi \, x (k+l+m))
				\ud x \\
			={}
				&\1_{\{0\}}(k-l-m) +\1_{\{0\}}(k-l+m) +\1_{\{0\}}(k+l-m) +\1_{\{0\}}(k+l+m).
		\end{split}
	\end{equation}
	Analogously it follows that for all $k,l,m \in \Z$ it holds that
	\begin{align}
	\label{eq: equation sin cos cos}
				&4\int_0^1 
					\sin ( 2 \pi \, k \, x) \cos ( 2 \pi \, l \, x) \cos ( 2 \pi \, m \, x) 
				\ud x
			=
				0, \\
	\label{eq: equation sin sin cos}
		\begin{split}
				&4\int_0^1 
					\sin ( 2 \pi \, k \, x) \sin ( 2 \pi \, l \, x) \cos ( 2 \pi \, m \, x) 
				\ud x
	\\& \qquad
			=
				\1_{\{0\}}(k-l-m) +\1_{\{0\}}(k-l+m) -\1_{\{0\}}(k+l-m)-\1_{\{0\}}(k+l+m),
		\end{split} \\
		\label{eq: equation sin sin sin}
				&4\int_0^1 
					\sin ( 2 \pi \, k \, x) \sin ( 2 \pi \, l \, x) \sin ( 2 \pi \, m \, x) 
				\ud x
			=
				0.
	\end{align}
	Combining \eqref{eq: equation cos cos cos},
	\eqref{eq: equation sin cos cos},
	\eqref{eq: equation sin sin cos},
	and
	\eqref{eq: equation sin sin sin} ensures that for all $k,l,n \in \Z$
	with $|n| > |k|+|l|$ it holds that
	\begin{equation}
		\int _0^1 \xi_{k}(x) \, \xi_{l}(x) \, \xi_{n}(x) \ud x = 0.
	\end{equation}
	Thus we get for all $(k,l), (m,n), (r,s) \in \Z^2$
	with $(|r|-|k|-|m|) \vee (|s|-|l|-|n|) > 0$ that
	\begin{equation}
	\label{eq: int 3 factor}
		\int_0^1 \int _0^1 \xi_{k,l}(x,y) \, \xi_{m,n}(x,y) \, \xi_{r,s}(x,y) \ud x \ud y = 0.
	\end{equation}
	In addition, we have for all $\overline{a}_{0,0} \in \R$, $N \in \N$
	and all $a_{k,l} \in \R$, $(k,l) \in  \Z_N^2$  that
	\begin{equation}
	\label{eq: F in H N}
		\begin{split}
				&\Big(
					D_\R \Big(
						\overline{a}_{0,0} \overline{e}_{0,0} + \sum_{(k,l)\in \Z_N^2} a_{k,l} e_{k,l}
					\Big)
				\Big)
					\Big
						(\overline{a}_{0,0} \overline{e}_{0,0} + \sum_{(k,l)\in \Z_N^2} a_{k,l} e_{k,l}
					\Big) \\
			={}
				&\Big(
					\sum_{(k,l)\in \Z_N^2 \backslash \{(0,0) \}} a_{k,l} (D_\R e_{k,l})
				\Big)
					\Big
						(\overline{a}_{0,0} \overline{e}_{0,0} + \sum_{(k,l)\in \Z_N^2} a_{k,l} e_{k,l}
					\Big) \\
			={}
				&\sum_{(k,l)\in \Z_N^2 \backslash \{(0,0) \}}
					a_{k,l} \overline{a}_{0,0} \, (D_\R e_{k,l}) \overline{e}_{0,0}
				\sum_{(k,l)\in \Z_N^2 \backslash \{(0,0) \}}
				\sum_{(m,n)\in \Z_N^2}
					a_{k,l} a_{m,n} \, (D_\R e_{k,l}) e_{m,n}.
		\end{split}
	\end{equation}
	Combining \eqref{eq: F ohne 00},
	\eqref{eq: F 0},
	\eqref{eq: F 00},
	\eqref{eq: int 3 factor}, \eqref{eq: F in H N} shows that for
	all $N \in \N$, $u \in H_N$, and all $(k,l) \in \Z^2$ with 
	$|k|+|l| > 2N$ it holds that
	\begin{equation}
		\langle (D_\R u) u, e_{k,l} \rangle_{L^2((0,1)^2,\R^2)} =0
	\end{equation}
	and this verifies that for all $u \in H_N$ it holds that
	$F(u) \in H_{2N}$.
	Therefore Lemma \eqref{l: norm equivalence Navier}
	implies that for all $r, s \in [0,\infty)$
	it holds that
	\begin{align}
	\label{eq: same continuity Navier}
			F|_{H_r} \in 
				C_{\mathbb{B}^{2r}_{2,2}((0,1)^2,\R^2) ,\mathbb{B}^{2s}_{2,2}((0,1)^2,\R^2)}
					(H_r, B^{2s}_{2,2}((0,1)^2,\R^2) 
		\Leftrightarrow
			F|_{H_r} \in 
				C_{\mathbb{H}_r ,\mathbb{H}_s}
					(H_r, H_s).
	\end{align}
	Moreover, Lemma \ref{lem: lemma: interpolation theorem for B} and
	\eqref{eq: h+1/2 continuity of F Cor2}
	in Corollary \ref{cor: F assumption} 
		(with $s \leftarrow r$, $d \leftarrow 2$, $\alpha \leftarrow i$
		and $H_t \leftarrow B_{2,2}^{2t}((0,1)^2)$, $t \in \R$)
	assures that for all
	$r \in [0,\infty)$ and all $\eps \in (0,\infty)$
	there exists a $C \in (0,\infty)$ such that
	for all $u=(u_{(1,0)}, u_{(0,1)}) \in B^{2r+1 + \eps \1_{\{0\}}(r)}_{2,2}((0,1)^2,\R^2)$
	and all $v=(v_{(1,0)}, v_{(0,1)}) \in B^{2r+1 + \eps \1_{\{0\}}(r)}_{2,2}((0,1)^2,\R^2)$
	it holds that
		\begin{align}
		\label{eq: F cont r>=0}
				&\|F(u) -F(v)\|_{B^{2r}_{2,2}((0,1)^2,\R^2)} \\
			\leq{} \nonumber
				&c \|
					(D_{\R^2} u) \, u -(D_{\R^2} v) \, v
				\|_{B^{2r}_{2,2}((0,1)^2,\R^2)} 
				+ \eta \|u-v \|_{B^{2r}_{2,2}((0,1)^2,\R^2)} \\
			\leq{} \nonumber
				&(
					\|
						(D_{\R^2} (u-v)) u
					\|_{B^{2r}_{2,2}((0,1)^2, \R^2)} 
					+\|
						(D_{\R^2} v) (u-v)
					\|_{B^{2r}_{2,2}((0,1)^2, \R^2)} 
					+  \|u-v\|_{B^{2r}_{2,2}((0,1)^2,\R^2)} 
				)\\
			\leq{} \nonumber
				&\bigg(
					\sum_{i,j \in \N^2_0, |i|=|j|=1} \Big(
							\|
								(D^{i}_{\R^2} (u_j-v_j)) u_i
							\|_{B^{2r}_{2,2}((0,1)^2)} 
							+\|
								(D^{i}_{\R^2} v_j) (u_i-v_i)
							\|_{B^{2r}_{2,2}((0,1)^2)} 
					\Big) 
			\\ & \nonumber \qquad
					+  \|u-v\|_{B^{2r}_{2,2}((0,1)^2,\R^2)} 
				\bigg)\\
			\leq{} \nonumber
				&C\bigg(
					\sum_{i,j \in \N^2_0, |i|=|j|=1} \Big(
						\|
							u_j-v_j
						\|_{B^{2r+1 + 2\eps \1_{\{0\}}(r)}_{2,2}((0,1)^2)} 
						\|
							u_i
						\|_{B^{2r+1 + 2\eps \1_{\{0\}}(r)}_{2,2}((0,1)^2)} 
				\\ & \nonumber \qquad
						+\|
							 v_j
						\|_{B^{2r+1 + 2\eps \1_{\{0\}}(r)}_{2,2}((0,1)^2)}
						\|
							u_i-v_i
						\|_{B^{2r+1 + 2\eps \1_{\{0\}}(r)}_{2,2}((0,1)^2)}
					\Big)
			\\ & \nonumber \qquad
					+  \|u-v\|_{B^{2r+1 + 2\eps \1_{\{0\}}(r)}_{2,2}((0,1)^2,\R^2)} 
				\bigg)\\
			\leq{} \nonumber
				&C\Big(
						4\|
							u-v
						\|_{B^{2r+1 + 2\eps \1_{\{0\}}(r)}_{2,2}((0,1)^2,\R^2)} 
						\big(
							\|
								u
							\|_{B^{2r+1 + 2\eps \1_{\{0\}}(r)}_{2,2}((0,1)^2,\R^2)} 
							+\|
								 v
							\|_{B^{2r+1 + 2\eps \1_{\{0\}}(r)}_{2,2}((0,1)^2,\R^2)}
						\big) \\
		\nonumber
				& \qquad
					+  \|u-v\|_{B^{2r+1 + 2\eps \1_{\{0\}}(r)}_{2,2}((0,1)^2,\R^2)} 
				\Big).
	\end{align}
	Hence \eqref{eq: F cont r>=0} and \eqref{eq: same continuity Navier}
	prove  that for all $r \in [0,\infty)$ and all $\eps \in (0,\infty)$ it holds that
	\begin{equation}
	\label{eq: F cont for all r}
			F|_{H_{r+1/2 +\eps \1_{\{0\}}(r)}} \in 
				C_{\mathbb{H}_{r+1/2 +\eps \1_{\{0\}}(r)} ,\mathbb{H}_{r}}
					(H_{r+1/2 +\eps \1_{\{0\}}(r)}, H_{r}).
	\end{equation}
	Thus \eqref{eq: F cont for all r} 
		(with $r \leftarrow \gamma + \vartheta -\nicefrac 12$)
	the fact that $\alpha_1 < \nicefrac 12$
	and 
	\eqref{eq: B continuity Navier}
	verifies \eqref{eq: F and B continuity} (with $H \leftarrow H_\gamma$).
	In addition, it holds for all $(u_1,u_2) \in H_{\gamma +\vartheta}$ that
	\begin{equation}
	\label{eq: repr. 1st bad term}
		\begin{split}
				(D^{(1,0)}_{\R^2} u_2) \, u_1
			={}
				&(D^{(1,0)}_{\R^2} (u_1+u_2)) \, (u_1+u_2)
				-(D^{(1,0)}_{\R^2} u_1) \, u_1
				-(D^{(1,0)}_{\R^2} u_1) \, u_2
				-(D^{(1,0)}_{\R^2} u_2) \, u_2 \\
			={}
				&(D^{(1,0)}_{\R^2} (u_1+u_2)) \, (u_1+u_2)
				-(D^{(1,0)}_{\R^2} u_1) \, u_1
				+(D^{(0,1)}_{\R^2} u_2) \, u_2
				-(D^{(1,0)}_{\R^2} u_2) \, u_2.
		\end{split}
	\end{equation}
	Analogously it follows that for all $(u_1,u_2) \in H_{\gamma +\vartheta}$
	it holds that
	\begin{equation}
	\label{eq: repr. 2nd bad term}
				(D^{(0,1)}_{\R^2} u_1) \, u_2
			=
				(D^{(0,1)}_{\R^2} (u_1+u_2)) \, (u_1+u_2)
				-(D^{(0,1)}_{\R^2} u_1) \, u_1
				+(D^{(1,0)}_{\R^2} u_1) \, u_1
				-(D^{(0,1)}_{\R^2} u_2) \, u_2.
	\end{equation}
	Thus \eqref{eq: F cont for all r} (with $r \leftarrow \gamma +\vartheta-\nicefrac 12)$,
	Lemma \ref{l: norm equivalence Navier},
	\eqref{eq: repr. 1st bad term}, 
	\eqref{eq: repr. 2nd bad term},  
	Lemma \ref{lem: lemma: interpolation theorem for B},
	\eqref{eq: h+1/2 bound for F Cor} in Corollary \ref{cor: F assumption}
	(with 
		$H_r \leftarrow B^{2r}_{2,2}((0,1)^2)$, $r \in \R$,
		$d \leftarrow 2$),
	and 
	the fact that
	$\alpha_1 < \gamma$
	establishes that  there exist
	$C_1, C_2, C_3 \in (0,\infty)$ such that for all 
	$u=(u_1, u_2) \in H_{\gamma+\vartheta}$
	it holds that
	\begin{align}
	\label{eq: F vartheta -1/2 bound Navier}
				&\|F(u) \|_{H_{\gamma + \vartheta-1/2}}
			\leq
				C_1(
					\|(D_{\R^2} u) \, u\|_{B^{2\gamma + 2\vartheta-1}_{2,2}((0,1)^2, \R^2)} 
					+  \|u\|_{H_{\gamma + \vartheta-1/2}} 
				)\\
			\leq{} \nonumber
				&C_1\big(
					\|
						(D^{(1,0)}_{\R^2} u_1) \, u_1 
					\|_{B^{2\gamma + 2\vartheta-1}_{2,2}((0,1)^2)}  
					+
					\|
						(D^{(1,0)}_{\R^2} u_2) \, u_1 
					\|_{B^{2\gamma + 2\vartheta-1}_{2,2}((0,1)^2)}  
			\\ & \nonumber \quad
					+
					\|
						(D^{(0,1)}_{\R^2} u_1) \, u_2 
					\|_{B^{2\gamma + 2\vartheta-1}_{2,2}((0,1)^2)}  
					+
					\|
						(D^{(0,1)}_{\R^2} u_2) \, u_2 
					\|_{B^{2\gamma + 2\vartheta-1}_{2,2}((0,1)^2)} 
					+\|u\|_{H_{\gamma + \vartheta-1/2}}
				\big) \\
			\leq{} \nonumber
				&C_1\big(
					3\|
						(D^{(1,0)}_{\R^2} u_1) \, u_1 
					\|_{B^{2\gamma + 2\vartheta-1}_{2,2}((0,1)^2)}  
					+
					3\|
						(D^{(0,1)}_{\R^2} u_2) \, u_2 
					\|_{B^{2\gamma + 2\vartheta-1}_{2,2}((0,1)^2)}
			\\ & \nonumber \quad
					+\|
						(D^{(0,1)}_{\R^2} u_1) \, u_1 
					\|_{B^{2\gamma + 2\vartheta-1}_{2,2}((0,1)^2)}  
					+
					\|
						(D^{(1,0)}_{\R^2} u_2) \, u_2 
					\|_{B^{2\gamma + 2\vartheta-1}_{2,2}((0,1)^2)} 
				\\ & \nonumber \quad
					+
					\|
						(D^{(1,0)}_{\R^2} (u_1+u_2) \, (u_1+u_2)
					\|_{B^{2\gamma + 2\vartheta-1}_{2,2}((0,1)^2)} 
				\\ & \nonumber \quad
					+
					\|
						(D^{(0,1)}_{\R^2} (u_1 +u_2)) \, (u_1+u_2) 
					\|_{B^{2\gamma + 2\vartheta-1}_{2,2}((0,1)^2)} 
					+\|u\|_{H_{\gamma + \vartheta-1/2}}
				\big) \\
			\leq{} \nonumber
				&C_2\big(
					3\|
						u_1 
					\|_{B^{2\gamma + 2\vartheta+1-2\alpha_1}_{2,2}((0,1)^2)}  
					\|
						u_1 
					\|_{B^{2\gamma}_{2,2}((0,1)^2)}
					+
					3\|
						u_2 
					\|_{B^{2\gamma + 2\vartheta+1-2\alpha_1}_{2,2}((0,1)^2)}  
					\|
						u_2 
					\|_{B^{2\gamma}_{2,2}((0,1)^2)}
				\\ & \nonumber \quad
					+\|
						u_1 
					\|_{B^{2\gamma + 2\vartheta+1-2\alpha_1}_{2,2}((0,1)^2)}  
					\|
						u_1 
					\|_{B^{2\gamma}_{2,2}((0,1)^2)}  
					+
					\|
						u_2 
					\|_{B^{2\gamma + 2\vartheta+1-2\alpha_1}_{2,2}((0,1)^2)}  
					\|
						u_2 
					\|_{B^{2\gamma}_{2,2}((0,1)^2)} 
				\\ & \nonumber \quad
					+
					\|
						u_1 +u_2
					\|_{B^{2\gamma + 2\vartheta+1-2\alpha_1}_{2,2}((0,1)^2)}  
					\|
						u_1 +u_2
					\|_{B^{2\gamma}_{2,2}((0,1)^2)} 
				\\ & \nonumber \quad
					+
					\|
						u_1 +u_2
					\|_{B^{2\gamma + 2\vartheta+1-2\alpha_1}_{2,2}((0,1)^2)}  
					\|
						u_1 +u_2
					\|_{B^{2\gamma}_{2,2}((0,1)^2)} 
					+\|u\|_{H_{\gamma + \vartheta+1/2 -\alpha_1}}
				\big) \\
			\leq{} \nonumber
				&C_2\big(
					16\|
						u
					\|_{B^{2\gamma + 2\vartheta+1-2\alpha_1}_{2,2}((0,1)^2,\R^2)}  
					\|
						u
					\|_{B^{2\gamma}_{2,2}((0,1)^2,\R^2)}
					+\|u\|_{H_{\gamma + \vartheta+1/2 -\alpha_1}}
				\big) \\
			\leq{} 
	\nonumber
				&C_3\big(
					\|
						u
					\|_{H_{\gamma + \vartheta+1/2-\alpha_1}}  
					\|
						u
					\|_{H_{\gamma}}
					+\|u\|_{H_{\gamma + \vartheta+1/2 -\alpha_1}}
				\big)
	\end{align}
	and this proves \eqref{eq: F H 1/2 bound}
		(with $\mathbb{H} \leftarrow \mathbb{H}_\gamma$
		$r_x \leftarrow 1$, 
		$R_x \leftarrow C_3 (\eta^{-\vartheta } +\|x\|_{H_\gamma} +1)$).
  Moreover, 
	\eqref{eq: V is lypanov Navier} implies
	that for all $N \in \N$ and all
	all $(t,u) \in [0,T] \times V_N$ 
	it holds that 
	\begin{equation}
		\E [ V(X_{t}^{N,u})] \leq V(u) 
			e^{
				(
					p \eta
					+4 \, p (p-1) 
					( K(1) \vee \| B(0) \|^2_{HS(\mathbb{U},\mathbb{H})} )
				) t
			}
	\end{equation}
	which verifies \eqref{eq: V bound for XN}.
In addition, \eqref{eq: X continuous in start point Navier}
 proves \eqref{eq: cont assumption},
\eqref{eq: B bounded wrt H_vartheta Navier} 
proves
\eqref{eq: B H 1 bound}
	(with $\mathbb{H} \leftarrow \mathbb{H}_\gamma$, 
	$r_x \leftarrow 1$, $R_x \leftarrow K(\|x\|_{H_\gamma}+\eta^{-\vartheta})$), 
\eqref{eq: B C1 bound Navier} 
proves
\eqref{eq: B C1 bound thm},
and 
\eqref{eq:convLocLip_assumption Navier}
proves
\eqref{eq:convLocLip_assumption thm}.
Thus Theorem \ref{thm: existence} yields that
  $ 
    f \in \C_{\R \times \mathbb{H}_\gamma}(\R \times H_\gamma, \R)
  $,
	that  
	$ 
    f 
  $
	is bounded on $\R \times \mathbb{H}_\gamma$-bounded 
	subsets of $(0,T) \times H_\gamma$,
	that for all $\tilde{R} \in (0, \infty)$ it holds that
	\begin{equation}
		\begin{split}
				\lim_{R \to \infty} \lim_{r \downarrow 0} \lim_{\eps \downarrow 0}	
				\sup \Bigg \{
					&f(t,u) -  \varphi(v)
					\colon
					~u, v \in H_{\gamma+2 \vartheta},
					~t \in (0,T),
					~\|u\|_{H_{\gamma+\vartheta}}^2
						\vee 
						\|v\|_{H_{\gamma+\vartheta}}^2  
						\leq R, \\
					~&\|u \|_{H_\gamma} \leq \tilde{R}, 
					~\| u - v\|_{H_\gamma} \leq r,
					~t \leq \eps 
				\Bigg \}
			= 
				0,
		\end{split}
	\end{equation}
	that for all $q \in [p,\infty)$ it holds that
  \begin{equation}
    \lim_{ r \to \infty }
      \sup_{ t \in [0,T]}
			\sup \left \{
				\frac
					{| f(t,u) |}
					{\|u\|^q_{L^2(0,1)}+2}
				\colon u \in H, ~\|u\|_{H_\gamma} \geq r
			\right \}
    = 0,
  \end{equation}
  and that
  $ 
    f|_{ (0,T) \times H_\gamma} 
  $
  is a viscosity solution of
  \begin{equation} 
		\begin{split}
			&\tfrac{ \partial }{ \partial t }
			f(t,u) -
			\big\langle
				A u + F(u), 
				I^{-1}_{\mathbb{H}_\gamma} (D_{\mathbb{H_\gamma}} f)(t,u)
			\big\rangle_{H_\gamma} 
		\\ & \qquad
			-\big \langle 
				B(u),
				I_{\mathbb{H}_{\gamma + \vartheta}}^{-1}
					\big( (D^2_{\mathbb{H_\gamma}} \, f)(t,u) \, B(u) \big)
			\big \rangle_{HS(\mathbb{U},\mathbb{H}_{\gamma+\vartheta})} 
			= 0
		\end{split}
  \end{equation}
	for $ (t,u) \in (0,T) \times H_\gamma $ relative to 
	$
		(
			(0,T) \times H_\gamma \ni (t,u) 
			\to \nicefrac 12 \|u \|^2_{\mathbb{H}_{\gamma+\vartheta}} \in [0, \infty],
			\R \times \mathbb{H}_\gamma, 
			\R \times \mathbb{H}_{\gamma+\vartheta}
		)
	$. 
	Next we apply Corollary \ref{cor:uniqueness2} 
		(with $\mathbb{H} \leftarrow \mathbb{H}_\gamma$, 
		$O \leftarrow H_\gamma$, $v \leftarrow 0$, 
		$\kappa \leftarrow -\gamma$,
		$\alpha_1 \leftarrow \alpha_1$,
		and $\alpha_2 \leftarrow \alpha_1$)
	to show that $f$ is the unique viscosity solution.
	Therefore assume that $g$ is another viscosity solution of 
	\eqref{eq: Kol equation Burger}, 
	which satisfy that 
	$ 
    g \in \C_{\R \times \mathbb{H}_\gamma}([0,T] \times H_\gamma, \R)
  $,
	that $g$ is
	bounded on $\R \times \mathbb{H}_\gamma$-bounded 
	subsets of $[0,T] \times H_\gamma$, 
	that $g$ have at most polynomial $\H$-growth,
  and that for all $ u \in H_\gamma $ it holds that
  $
    g(0, u) = \varphi(u)
  $.
	Hence, there exists a $q \in [p,\infty)$, which we fix for the rest of the proof,
	such that
	 \begin{equation}
    \lim_{ r \to \infty }
      \sup_{ t \in [0,T]}
			\sup \left \{
				\frac
					{| g(t,u) |}
					{\|u\|^q_{H}+2}
				\colon u \in H_\gamma, ~\|u\|_{H} \geq r
			\right \}
    = 0.
  \end{equation}
	Moreover, the
	continuity of $\varphi$
	with respect to the $\| \cdot \|_{H_\gamma}$-norm,
	the
	continuity of $g$
	with respect to the $\| \cdot \|_{\R \times H_\gamma}$-norm,
	and the fact that for all $R \in (0,\infty)$ it holds that
	the set $\{(t,x) \in [0,T] \times H_\gamma \colon \|x\|_{H_{\gamma+\vartheta}} \leq R\}$
	is compact in $\H_\gamma$
	yields for all $R, \tilde{R} \in (0, \infty)$ that
	\begin{equation}
		\begin{split}
				&\lim_{r \downarrow 0} \lim_{\eps \downarrow 0}	
				\sup \bigg \{
					g(t,u) -  \varphi(\hat{u})
					\colon
					u, \hat{u} \in H_{2\varz+\gamma},
					~t, \hat{t} \in (0,T),
					~\tfrac{e^{-\theta t} \|u\|_{H_{\varz+\gamma}}^2}{\|u\|_H^q+2} 
						\vee \tfrac{e^{-\theta \hat{t}} 
							\|\hat{u}\|_{H_{\varz+\gamma}}^2}{\|\hat{u}\|_H^q+2} \leq R, 
				\\& \qquad \qquad \qquad \qquad
					~\|u \|_{H_\gamma} \leq \tilde{R}, 
					~\| u - \hat{u}\|_{H_\gamma} \leq r,
					~t \vee \hat{t} \leq \eps 
				\bigg \} \\
			={}  
				&\lim_{r \downarrow 0}
				\sup \Big \{
					\varphi(u) -  \varphi(\hat{u})
					\colon
					u, \hat{u} \in H_{2\varz+\gamma},
					~\tfrac {\|u\|_{H_{\varz+\gamma}}^2} {\|u\|_H^q+2} 
						\vee \tfrac {\|\hat{u}\|_{H_{\varz+\gamma}}^2} {\|\hat{u}\|_H^q+2} \leq R,
					~\|u \|_{H_\gamma} \leq \tilde{R},
					~\| u - \hat{u}\|_{H_\gamma} \leq r
				\Big \} \\
			={}
				&0.
		\end{split}
	\end{equation}
	Next note that 
	the Sobolev inequality (see, e.g., the Theorem on page 31
	in Runst \& Sickel \cite{RunstSickel1996}),
	Lemma \ref{l: norm equivalence Navier},
	and Lemma \ref{lem: sobolev holder inequality without interpolation for Lp} 
		(with 
		$d \leftarrow 2$,
		$p \leftarrow \tfrac{2}{2r+1}$,
		$u \leftarrow D^{i}_{\R^2} u_j$,
		$v \leftarrow u_i$,
		$r_1 \leftarrow 2\gamma \1_{\{0\}}(r)$,
		$r_2 \leftarrow 1+2\gamma \1_{\{0\}}(r)$,
		$\eps \leftarrow 2\gamma$)
	proves that 
	for all $r \in [0,\nicefrac 12)$
	there exist $C_1, C_2 \in (0,\infty)$
	such that
	for all $u=(u_{(1,0)}, u_{(0,1)}) \in H_{\gamma +\vartheta}$
	and all $v=(v_{(1,0)}, v_{(0,1)}) \in H_{\gamma +\vartheta}$
	it holds that
	\begin{align}
	\label{eq: F cont r<=0}
		\begin{split}
				&\|F(u) \|_{H_{-r}} \\
			\leq{}
				&\sup_{w \in H_{r}} \bigg(
					\int_0^1 \int_0^1
						c\Big \langle	
								\big( (D_{\R^2}u)(x,y) \big) (u(x,y)),
							\tfrac {w(x,y)}{\|w\|_{H_{r}}}
						\Big \rangle_{\R^2}
					\ud x \ud y
				\bigg) 
				+ \eta \|u\|_{H_{-r}} \\ 
			={}
				&\sup_{w \in H_{r}} \bigg(
					\int_0^1 \int_0^1
						c\sum_{i,j \in \N^2_0, |i|=|j|=1} \Big(
								(D^{i}_{\R^2} u_j)(x,y) \, u_i(x,y)
							\tfrac {w_j(x,y)}{\|w\|_{H_{r}}}
						\Big)
					\ud x \ud y
				\bigg) 
				+ \eta \|u\|_{H_{-r}} \\ 
			\leq{}  
				&\sum_{i,j \in \N^2_0, |i|=|j|=1}
				\sup_{w \in H_{r}} \bigg(
					|c| \, 
					\|
							(D^{i}_{\R^2} u_j) \, u_i
					\|_{L^{2/(2r+1)}((0,1)^2)} 
					\cdot 
					\tfrac 
						{
							\|
								w_j
							\|_{L^{2/(1-2r)} ((0,1)^2)}
						}
						{\|w\|_{H_{r}}}
				\bigg) 
				+ \eta \|u\|_{H_{-r}} \\
			\leq{}
				&C_1 \sum_{i,j \in \N^2_0, |i|=|j|=1} \Big(
				\|
					D^{i}_{\R^2}u_j 
				\|_{B_{2,2}^{2\gamma \1_{\{0\}}(r)}((0,1)^2)}
				\|
					u_i
				\|_{B_{2,2}^{1 + 2\gamma \1_{\{0\}}(r) }((0,1)^2)} 
				\Big)
				\sup_{w \in H_{r}} \Big(
					\tfrac 
						{
							\|
								w
							\|_{B^{2r} ((0,1)^2,\R^2)}
						}
						{\|w\|_{H_{r}}}
				\Big) 
			\\ 
			& \qquad
				+\eta^{-1/2-r-\gamma \1_{\{0\}}(r)} \|u\|_{H_{1/2 +\gamma \1_{\{0\}}(r)}}\\
			\leq{}
				&C_2 \Big(
				\sum_{i,j \in \N^2_0, |i|=|j|=1} \big(
					\|u_i\|_{B_{2,2}^{1 + 2\gamma \1_{\{0\}}(r)}((0,1)^2)}
					\|u_j\|_{B_{2,2}^{1 + 2\gamma \1_{\{0\}}(r)}((0,1)^2)}
				\big) 
				+\|u\|_{B_{2,2}^{1 + 2\gamma \1_{\{0\}}(r)}((0,1)^2,\R^2)}
				\Big) \\
			\leq{}
				&4C_2 \big(
					\|u\|^2_{B_{2,2}^{1 + 2\gamma \1_{\{0\}}(r)}((0,1)^2,\R^2)}
					+\|u\|_{B_{2,2}^{1 + 2\gamma \1_{\{0\}}(r)}((0,1)^2,\R^2)}
				\big).
		\end{split}
	\end{align}
	Thus
	\eqref{eq: F cont r<=0} (with $r \leftarrow \gamma +\vartheta -1$),
	Lemma \ref{l: norm equivalence Navier},
	\eqref{eq: F cont r>=0} (with $r \leftarrow \gamma +\vartheta -1$, 
	$\eps \leftarrow \gamma$, and $v\leftarrow 0$),
	the fact that $2\gamma -1 <0$
	and the fact that $\vartheta \geq \nicefrac 12$
	yield that 
	there exist $C_1, C_2,C_3 \in (0,\infty)$ such that for all $u \in H_{\gamma +\vartheta}$
	it holds that
	\begin{align}
	\nonumber
			&\|F(u)\|_{H_{\gamma+\vartheta - 1}} \\
			\leq{}
				&C_1 \big(
					\|u\|^2_{B_{2,2}^{
						(1 \vee (2\gamma+2\vartheta-1))+ 2\gamma \1_{\{0\}}(\gamma+\vartheta-1)
					}((0,1)^2,\R^2)}
					+\|u\|_{B_{2,2}^{
						(1 \vee (2\gamma+2\vartheta-1))+ 2\gamma \1_{\{0\}}(\gamma+\vartheta-1)
					}((0,1)^2,\R^2)}
				\big) \\ \nonumber
			\leq{}
				&C_2 \big(
					\|u\|^2_{B_{2,2}^{
						2\gamma+ 2\vartheta
					}((0,1)^2,\R^2)}
					+\|u\|_{B_{2,2}^{
						2\gamma+ 2\vartheta
					}((0,1)^2,\R^2)}
				\big) \\ \nonumber
			\leq{}
				&C_3 \|u\|_{H_{\gamma+\vartheta}} (\|u\|_{H_{\gamma+\vartheta}} +1). 
	\end{align}
	Combining this and 
	\eqref{eq: H bound B Navier} (with $u \leftarrow 0$)
	proves
	\eqref{eq: H bound for F, B and v}
		(with $H \leftarrow H_\gamma$, $\kappa \leftarrow -\gamma$). 
	Moreover,
	it follows from 
	Theorem 4 on page 84 in Runst \& Sickel \cite{RunstSickel1996},
	the fact that $2 \gamma < \nicefrac 12$,
	Theorem 2 on page 191 in Runst \& Sickel \cite{RunstSickel1996}
		(with $s_1 \leftarrow 2\gamma$,
		$s_2 \leftarrow 1 +2\gamma -2\alpha_1$,
		$p \leftarrow 2$,
		$q_1 \leftarrow 2$,
		and $q_2 \leftarrow 2$)
	the Proposition on page 29 in Runst \& Sickel \cite{RunstSickel1996},
	the Sobolev inequality 
		(see, e.g., Theorem 1 on page 32 in Runst \& Sickel \cite{RunstSickel1996}),
	and from Lemma \ref{l: norm equivalence Navier}
	that  there exist $C_1, C_2 \in (0,\infty)$
	such that for all 
	$u=(u_{(1,0)},u_{(0,1)})\in H_{1/2}$,
	$v=(v_{(1,0)},v_{(0,1)})\in H_{1/2}$,
	and all
	$w=(w_{(1,0)},w_{(0,1)})\in H_{1/2}$
	it holds that
	\begin{align}
	\nonumber
				&\int_0^1 \int_0^1
					\Big \langle
						\big( (D_{\R^2}w)(x,y) \big) (u(x,y)), 
						v(x,y)
					\Big \rangle_{\R^2} 
				\ud x \ud y\\
			={} \nonumber
				&\int_0^1 \int_0^1
					\sum_{i,j \in \N^2_0, |i|=|j|=1}
						(D^{i}_{\R^2} w_{j})(x,y) \, u_{i}(x,y) \, v_{j}(x,y)
					\ud x \ud y\\
			\leq{} 
				&\sum_{i,j \in \N^2_0, |i|=|j|=1}
					\| D^{i}_{\R^2} w_{j}\|_{B^{-2\gamma}_{2,2}((0,1)^2)} 
					\| u_{i}\, v_{j}\|_{B^{2\gamma}_{2,2}((0,1)^2)}\\
			\leq{}  \nonumber
				&C_1
				\sum_{i,j \in \N^2_0, |i|=|j|=1}
					\| w_{j}\|_{B^{1-2\gamma}_{2,2}((0,1)^2)} 
					\| u_i \|_{B^{1+2\gamma-2\alpha_1}_{2,2}((0,1)^2)}
					\|v_{j}\|_{B^{2\gamma}_{2,2}((0,1)^2)}\\
			\leq{} \nonumber
				&4C_1
					\| w \|_{B^{1-2\gamma}_{2,2}((0,1)^2,\R^2)} 
					\| u \|_{B^{1+2\gamma-2\alpha_1}_{2,2}((0,1)^2,\R^2)}
					\| v \|_{B^{2\gamma}_{2,2}((0,1)^2,\R^2)}\\
			\leq{} \nonumber
				&C_2
					\| w \|_{H_{1/2-\gamma}} 
					\| u \|_{H_{1/2+\gamma-\alpha_1}}
					\| v \|_{H_{\gamma}}. 
	\end{align}
	Analogously it follows that
	there exists a $C \in (0,\infty)$
	such that for all 
	$u,v, w\in H_{1/2}$
	it holds that
	\begin{align}
				\int_0^1 \int_0^1
					\Big \langle
						\big( (D_{\R^2}w)(x,y) \big) (u(x,y)), 
						v(x,y)
					\Big \rangle_{\R^2} 
				\ud x \ud y
			\leq{} 
				C
					\| w \|_{H_{1/2-\gamma}} 
					\| v \|_{H_{1/2+\gamma-\alpha_1}}
					\| u \|_{H_{\gamma}}. 
	\end{align}
	Thus the fact that 
	$H_{1/2}$ is dense in $H_{\nicefrac 12-\gamma}$
	and
	\eqref{eq: interchanging diff.}
	establishes that 
	there exists a $C \in (0,\infty)$ such that 
	for all $u,v \in H_{1/2}$ it holds that
	\begin{align}
				&\|F(u) -F(v)\|_{H_{\gamma-1/2}} \\ \nonumber
			\leq{}
				&\sup_{w \in H_{1/2}} \bigg(
					\int_0^1 \int_0^1
						c\Big \langle
								\big( (D_{\R^2}u)(x,y) \big) (u(x,y))
								-\big( (D_{\R^2}v)(x,y) \big) (v(x,y)),
							\tfrac {w(x,y)}{\|w\|_{H_{1/2-\gamma}}}
						\Big \rangle_{\R^2}
					\ud x \ud y	
				\bigg) 
			\\ \nonumber & \quad
				+ \eta \|u-v\|_{H_{\gamma-1/2}}\\ \nonumber
			={}
				&\sup_{w \in H_{1/2}} \bigg(
					-\int_0^1 \int_0^1
						c\Big \langle
							\big( (D_{\R^2}w)(x,y) \big) (u(x,y)), 
							\tfrac {u(x,y)}{\|w\|_{H_{1/2-\gamma}}}
						\Big \rangle_{\R^2}
				\\ \nonumber & \quad
						-c\Big \langle
							\big( (D_{\R^2}w)(x,y) \big) (v(x,y)), 
							\tfrac {v(x,y)}{\|w\|_{H_{1/2-\gamma}}}
						\Big \rangle_{\R^2}
					\ud x \ud y
				\bigg) 
				+ \eta \|u-v\|_{H_{\gamma-1/2}}\\ \nonumber
			={}
				&\sup_{w \in H_{1/2}} \bigg(
					-\int_0^1 \int_0^1
						c\Big \langle
							\big( (D_{\R^2}w)(x,y) \big) (u(x,y)-v(x,y)), 
							\tfrac {u(x,y)}{\|w\|_{H_{1/2-\gamma}}}
						\Big \rangle_{\R^2}
				\\  \nonumber & \qquad
						+c\Big \langle
							\big( (D_{\R^2}w)(x,y) \big) (v(x,y)), 
							\tfrac {u(x,y)-v(x,y)}{\|w\|_{H_{1/2-\gamma}}}
						\Big \rangle_{\R^2}
					\ud x \ud y
				\bigg) 
				+ \eta \|u-v\|_{H_{\gamma-1/2}}\\ \nonumber
			\leq{}
				&C \big(
					\| u-v \|_{H_{1/2+\gamma-\alpha_1}}
					\| u \|_{H_{\gamma}}
					+
					\| u-v \|_{H_{1/2+\gamma-\alpha_1}}
					\| v \|_{H_{\gamma}}
					+\eta^{\alpha_1} \| u-v \|_{H_{1/2+\gamma-\alpha_1}}
				\big)
	\end{align}
	and this shows
	\eqref{eq: Lipschitz continuity of F}.
	In addition, \eqref{eq: F vartheta -1/2 bound Navier}
	ensure \eqref{eq: H_3/2 bound for F},
	\eqref{eq: B bounded wrt H_vartheta Navier}
	ensure \eqref{eq: H_3/2 bound for B},
	\eqref{eq: B Lipschitz wrt H_1/2 Navier}
	enures
	\eqref{eq: Lipschitz continuity of B},
	\eqref{eq: V is lypanov Navier} ensures \eqref{eq:assumption_V},
	and Remark \ref{rem: assumption on V} 
		(with $H \leftarrow H_\gamma$,
		$\kappa \leftarrow -\gamma$, $l \leftarrow \gamma$, and
		$f \leftarrow (\R \ni x \to  (2x)^{\nicefrac p2} + 2 \in (1,\infty))$)
	ensures
	\eqref{eq: V continuous with respect to H1},
	\eqref{eq: V Lipschitz with respect to H -gamma},
	\eqref{eq: V bounded with respect to H},
	and \eqref{eq: V bounded with respect to H 1/2}.
	Thus Corollary \ref{cor:uniqueness2}
		(with $\mathbb{H} \leftarrow \mathbb{H}_\gamma$, 
		$O \leftarrow H_\gamma$, $v \leftarrow 0$, 
		$\kappa \leftarrow -\gamma$,
		$\alpha_1 \leftarrow \alpha_1$,
		$\alpha_2 \leftarrow \alpha_1$,
		and $V \leftarrow \|\cdot \|^q_H+2$)	
		verifies that $g=f$
	which finishes
	the proof of Corollary \ref{cor: Navier equation}.
\end{proof}
The next remark establishes sufficient conditions for $B$
to satisfy the conditions of Corollary \ref{cor: Navier equation}.
\begin{remark}[Sufficient condition for B]
\label{rem: sufficient condition for B Navier}
	Assume 
	the Setting \ref{sett: Navier},
	let
	$\vartheta \in [\nicefrac 12, 1 -\gamma)$,
	$\beta_1 \in (0,\gamma)$,
	let
	$C \in (0,\infty)$,
	let $\mathscr{P}_N \in L(\mathbb{U},\mathbb{H})$, $N \in \N$, 
	be finite-dimensional projection
	satisfying for all $n \in \mathcal{I}$ and all $N \in \N$ that
	\begin{equation}
	\label{eq: def of PN 1}
		\mathscr{P}_N \tilde{e}_{n} = \1_{[1,N] }(n) \, \tilde{e}_{n},
	\end{equation}
	let $b_{ij} \in \C^{3}([0,1]^2 \times \R^2; \R)$,
	$i,j \in \{1,2\}$,
	let $b \in \C ([0,1]\times \R^2; \R^{2 \times 2})$,
	satisfy that for all $x \in [0,1]^2$, $y \in \R^2$,
	and all $i,j \in \{1,2\}$ it holds that
	$ \| (D_{\R^2 \times \R^2} b_{ij})(x,y) \|_{(\R^2 \times \R^2)'} \leq C$,
	$ 
			\| (D^2_{\R^2 \times \R^2} b_{ij})(x,y) \|_{
				L(\R^2 \times \R^2,(\R^2 \times \R^2)')
			} 
		\leq 
			C
	$,
	and that
	\begin{equation}
		b(x,y) =
		\left (
			\begin{array}{cc}
				b_{11}(x,y) & b_{12}(x,y) \\
				b_{21}(x,y) & b_{22}(x,y) 
			\end{array}
		\right),
	\end{equation}
	let $Q \in HS(\mathbb{U},\mathbb{H}_{\gamma+\vartheta})$,
	assume that
	for all $u \in H$, $v \in U$
	and all $x \in (0,1)^2$ it holds that
	$(B(u) v)(x) = \pi^{L^2((0,1)^2,\R^2)}_H \big( b(\cdot, u(\cdot))  (Qv)(\cdot) \big) (x)$
	and assume that for all $n \in \mathcal{I}$, $i \in \{\gamma,\vartheta+\gamma\}$,
	and all $u \in H_{i}$ 
	it holds that 
	\begin{equation}
	\label{eq: right codomain}
			B(u) \tilde{e}_{n} \in H_i.
	\end{equation}
	Then B satisfies the condition of 
	Corollary \ref{cor: Navier equation} 
	with 
	$\beta \leftarrow \gamma$,
	$\beta_1 \leftarrow \beta_1$, $\beta_2 \leftarrow \beta_1$, 
	$\beta_3 \leftarrow \beta_1$,
	and
	$\chi \leftarrow \gamma$ 
	i.e.\@ there exists 
	a $\theta \in (0, \infty)$
	and an increasing
	$K \colon \R \to (0, \infty)$,
	such that
	\begin{equation}
	\label{eq: B continuity Navier rem}
		B|_{H_{1/2 +\gamma +\vartheta - \beta_1}} \in
				\C_{
					\mathbb{H}_{1/2 +\gamma +\vartheta - \beta_1},
					\mathbb{HS}(\mathbb{U}, \mathbb{H}_{\gamma+\vartheta})
				}(
					H_{1/2+\gamma +\vartheta - \beta_1},
					HS(\mathbb{U}, \mathbb{H}_{\gamma+\vartheta})
				), 
	\end{equation}
	that for all $u, v \in H_{1/2}$ 
	it holds that
	\begin{align}
		\label{eq: H bound B Navier rem}
			&\| B(u) - B(v)\|^2_{HS(\mathbb{U},\mathbb{H})}
		\leq
			C \, \|u-v\|^2_H
	\end{align}
	and that for all $u, v \in H_{\gamma + 2\vartheta}$ 
	it holds that
	\begin{align}
	\label{eq: B bounded wrt H_vartheta Navier rem}
			&\| B(u) \|^2_{HS(\mathbb{U},\mathbb{H}_{\gamma+\vartheta})} 
		\leq
			K(\|u\|_{H_\gamma}) 
				(\|u\|^2_{H_{1/2+\gamma+\vartheta- \beta_1}} + 1), \\
	\label{eq: B Lipschitz wrt H_1/2 Navier rem}
			&\|B(u) - B(v)\|_{HS(\mathbb{U}, \mathbb{H}_\gamma)}
		\leq
			K(\|v\|_{H_\gamma} \vee \|u \|_{H_\gamma}) 
				\|u-v\|_{H_{1/2+\gamma -\beta_1}}, \\
	\label{eq: Lipschitz bound B Navier rem}
		&| B |_{
			\C^1( 
				H_\gamma, 
				\left\| \cdot \right\|_{HS( \mathbb{U}, \mathbb{H} )} 
			)
		}
		< \infty,
	\end{align}
	and that for all $\mathbb{H}_\gamma$-bounded sets
	$ E \subseteq H_\gamma$
	it holds that
	\begin{equation}
	\label{eq: B C1 bound Navier rem} 
			\sup_{ N \in \N } \sup_{ v \in E }
			\left[
				\frac{
				\|
					B( v ) \mathscr{P}_N 
				\|_{ HS( \mathbb{U}, \mathbb{H}) }
				}{
					1 + \| v \|_{ H_\gamma }
				}
			\right]
			< \infty
	\end{equation}
	and that
	\begin{equation}
	\label{eq:convLocLip_assumption Navier rem}
			\limsup_{ N \to \infty }
			\sup_{ v \in E }
			\left[
			\frac{
			\|
				B( v ) ( \id_U - \mathscr{P}_N )
			\|_{ HS( \mathbb{U}, \mathbb{H} ) }
			}{
				1 + \| v \|_{ H_\gamma }
			}
			\right]
		=0.
	\end{equation}
\end{remark}
\begin{proof}
	Let $\mathbf{B}(r) \subseteq \R^2$, $r \in [0,\infty]$ be the set satisfying for all 
	$r \in [0,\infty]$ that $\mathbf{B}(r) = \{x \in \R^2 \colon \|x \|_{\R^2} \leq r\}$.
	Moreover, note that \eqref{eq: right codomain} and
	Lemma \ref{l: norm equivalence Navier}
	imply that 
	for all $i \in \{ 0, \gamma, \gamma+\vartheta\}$ 
	there exists a $C_1 \in (0,\infty)$
	such 
	that for all $u \in H_i$ it holds that
	\begin{align}
	\nonumber
				&\|B(u)\|^2_{HS(\U,\H_i)}
			=
				\sum_{n\in \mathcal{I}}
					\|B(u) \tilde{e}_{n}\|^2_{H_i}  
			\leq{}
				C_1 \big(	
					\sum_{n\in \mathcal{I}}
						\|B(u) \tilde{e}_{n}\|^2_{B^{2i}_{2,2}((0,1)^2,\R^2)}
				\big) \\
		\begin{split}
			\leq{}
				&C_1 \big(
				\sum_{n\in \mathcal{I}}
					\|b(\cdot,u(\cdot)) \, (Q \tilde{e}_{n}) (\cdot)\|^2_{B^{2i}_{2,2}((0,1)^2,\R^2)} 
				\big) \\
			={}
				&C_1
					\sum_{n\in \mathcal{I}}\Big(
						\Big \|
							\big \langle
								\big( b_{11}(\cdot,u(\cdot)), b_{12}(\cdot,u(\cdot)) \big),
								(Q \tilde{e}_{n}) (\cdot)
							\big \rangle_{\R^2}
						\Big \|^2_{B^{2i}_{2,2}((0,1)^2)}
		\end{split}
			\\ \nonumber & \quad
						+\Big \|
						\big \langle
							\big( b_{21}(\cdot,u(\cdot)), b_{22}(\cdot,u(\cdot)) \big),
							(Q \tilde{e}_{n}) (\cdot)
						\big \rangle_{\R^2}
					\|^2_{B^{2i}_{2,2}((0,1)^2)}
					\Big).
		\end{align}
		Therefore
		Theorem 2 on page 191 in Runst \& Sickel \cite{RunstSickel1996}
			(with $s_1 \leftarrow i$, $s_2 \leftarrow \gamma +\vartheta$,
			$d \leftarrow 2$, $p \leftarrow 2$,
			$q_1 \leftarrow 2$, $q_2 \leftarrow 2$, $q \leftarrow 2$)
		and Lemma \ref{l: norm equivalence Navier}
		verify that 
		for all $i \in \{ 0, \gamma, \gamma+\vartheta\}$ 
		there exist $C_1, C_2 \in (0,\infty)$
		such 
		that for all $u \in H_i$ it holds that
		\begin{equation}
		\label{eq: basic estimate H Navier}
			\begin{split}
				&\|B(u)\|^2_{HS(\U,\H_i)} \\
			\leq{}
				&C_1	
					\sum_{n\in \mathcal{I}} \Big(
						\big(
						\| b_{11}(\cdot,u(\cdot)) \|^2_{B^{2i}_{2,2}((0,1)^2)}
						+\| b_{12}(\cdot,u(\cdot)) \|^2_{B^{2i}_{2,2}((0,1)^2)}
					\big)
					\| Q \tilde{e}_{n} \|_{B^{2\gamma+2\vartheta}_{2,2}((0,1)^2,\R^2)}
				\\ & \qquad
					+\big(
						\| b_{21}(\cdot,u(\cdot)) \|^2_{B^{2i}_{2,2}((0,1)^2)}
						+\| b_{22}(\cdot,u(\cdot)) \|^2_{B^{2i}_{2,2}((0,1)^2)}
					\big)
					\| Q \tilde{e}_{n} \|_{B^{2\gamma+2\vartheta}_{2,2}((0,1)^2,\R^2)}
					\Big)\\
			\leq{}
				&C_1
					\sum_{n\in \mathcal{I}}
						\|
							b(\cdot,u(\cdot))
						\|^2_{(B^{2i}_{2,2}((0,1)^2))^{2 \times 2}} \, 
							\|
								Q \tilde{e}_{n}
							\|^2_{B^{2\gamma+2\vartheta}_{2,2}((0,1)^2,\R^2)}\\
			\leq{} &  
				C_2\|b(\cdot,u(\cdot)) \|^2_{(B^{2i}_{2,2}((0,1)^2))^{2 \times 2}} 
					\|Q \|_{HS(\mathbb{U},\mathbb{H}_{\gamma+\vartheta})}.
		\end{split}
	\end{equation}
	Analogously we derive that there exists a 
	$C_1 \in (0,\infty)$ such that for all $i \in \{0,\gamma, \gamma + \vartheta\}$
	and all
	$u, v \in H_i$ it holds that
	\begin{equation}
	\label{eq: basic estimate diff Navier}
				\|B(u)-B(v)\|_{HS(\U,\H_{i})} 
			\leq
				C_1 \|Q \|_{HS(\mathbb{U},\mathbb{H}_{\gamma+\vartheta})}
				\|
					b(\cdot,u(\cdot)) -b(\cdot,v(\cdot))
				\|_{(B^{2i}_{2,2}((0,1)^2))^{2 \times 2}}.
	\end{equation}
	Combining \eqref{eq: basic estimate diff Navier} and 
	Lemma \ref{lem: L2 b continuity inequality} 
		(with $d \leftarrow 2$, $m \leftarrow 2$, $p \leftarrow 0$, $q\leftarrow \infty$)
	shows that there exist $C_1,C_2 \in (0,\infty)$
	such that for all $u, v \in H$ it holds that
	\begin{equation}
	\label{eq: B Lip cont Navier rem}
		\begin{split}
			\|B(u) -B(v)\|^2_{HS(\U, \H)}
		\leq{}
			&C_1 \|b(\cdot,u(\cdot)) - b(\cdot,v(\cdot))  \|^2_{(L^2((0,1)^2))^{2 \times 2}} \\
		\leq{}
			&C_2 \|u-\tilde{u}\|^2_{L^2((0,1)^2,\R^2)}
		=
			C_2 \|u-\tilde{u}\|^2_{H}
		\end{split}
	\end{equation}
	and this verifies \eqref{eq: H bound B Navier rem}
	and \eqref{eq: Lipschitz bound B Navier rem}.
	Furthermore,
	\eqref{eq: basic estimate diff Navier},
	Corollary \ref{cor: Bs b continuity inequality}
		(with $m \leftarrow 2$, $d \leftarrow 2$, 
		$s \leftarrow 2\gamma$,
		$p_2 \leftarrow 0$,
		$p_1 \leftarrow 0$,
		$q \leftarrow \infty$, $q_1 \leftarrow \infty$,
		$q_2 \leftarrow 2$, $q_3 \leftarrow \infty$, $q_4 \leftarrow 2$,
		$q_5 \leftarrow \infty$,
		$\eps \leftarrow 2\gamma-2\beta_1$),
	and Lemma \ref{l: norm equivalence Navier}
	establish that there exist $C_1,C_2, C_3 \in (0,\infty)$
	such that for all $u,v \in H_\gamma$ it holds that
	\begin{align}
	\nonumber
			&\|B(u) -B(v)\|_{HS(\U, \H_\gamma)}
		\leq
			C_1\|
				b(\cdot,u(\cdot)) -b(\cdot,v(\cdot))
			\|_{(B^{2\gamma}_{2,2}((0,1)^2))^{2\times 2}} \\
	\nonumber
		\leq{}
			&C_2 \, \|
					u-v
				\|^{2}_{B^{1+2\gamma-2\beta_1 }_{2,2}((0,1)^2,\R^2)}
				\big(
					\|u\|^2_{B^{2\gamma}_{2,2}((0,1)^2,\R^2)}
					+1
				\big) 
			\\& \quad
				+C_2 \, \|
					u - v 
				\|^2_{B^{2\gamma}_{2,2}((0,1)^2,\R^2)}
				+C_2 \, \|u - v \|^2_{L^2((0,1)^2,\R^2)} \\
		\leq{} \nonumber
			&C_3 \big(
				\|
					u-v
				\|^{2}_{H_{1/2+\gamma-\beta_1 }}
				\big(
					\|u\|^2_{H_\gamma}
					+1
				\big) 
				+\|
					u - v 
				\|^2_{H_\gamma}
				+\|u - v \|^2_{H}
			\big) \\
	\nonumber
		\leq{}
			&C_3 
				\|
					u-v
				\|^{2}_{H_{1/2+\gamma-\beta_1 }}
			\big(
				\|u\|^2_{H_\gamma}
				+1 
				+\eta^{\beta_1-\nicefrac 12}
				+\eta^{\beta_1-\gamma-\nicefrac 12}
			\big)
	\end{align}
	and this implies \eqref{eq: B Lipschitz wrt H_1/2 Navier rem}.
	In addition, note that it holds that
	\begin{equation}
		\begin{split}
				&2\gamma \tfrac{1+\vartheta +\gamma -2\beta_1}{2\vartheta +1 -2 \beta_1}
				+ 
					\tfrac
						{(2\vartheta +2\gamma +1 -2 \beta_1)(\vartheta -\gamma)}
						{2\vartheta +1 -2 \beta_1} 
			={}
				2\gamma +2\gamma \tfrac{\gamma-\vartheta}{2\vartheta +1 -2 \beta_1}
				+ \vartheta -\gamma +\tfrac{2\gamma(\vartheta -\gamma)}{2\vartheta +1 -2 \beta_1}
			={}
				\vartheta +\gamma
		\end{split}
	\end{equation}
	and that
	\begin{equation}
		\begin{split}
				&2\gamma \tfrac{\vartheta +\gamma -2\beta_1}{2\vartheta +1 -2 \beta_1}
				+ 
					\tfrac
						{(2\vartheta +2\gamma +1 -2 \beta_1)(1+\vartheta -\gamma)}
						{2\vartheta +1 -2 \beta_1} 
			={}
				2\gamma +2\gamma \tfrac{\gamma-\vartheta-1}{2\vartheta +1 -2 \beta_1}
				+1+ \vartheta -\gamma +\tfrac{2\gamma(1+\vartheta -\gamma)}{2\vartheta +1 -2 \beta_1}
			={}
				1+\vartheta +\gamma.
		\end{split}
	\end{equation}
	Thus
	Lemma \ref{lem: lemma: interpolation theorem for B}
		(with $r \leftarrow \vartheta +\gamma$,
		$r_1 \leftarrow 2\gamma$,
		$r_2 \leftarrow 2\vartheta +2\gamma +1 -2 \beta_1$,
		$s_1 \leftarrow \frac{1+\vartheta +\gamma -2\beta_1}{2\vartheta +1-2\beta_1}$,
		$s_2 \leftarrow \frac{\vartheta -\gamma}{2\vartheta +1-2\beta_1}$,
		and with
		$r \leftarrow 1+\vartheta +\gamma$,
		$r_1 \leftarrow 2\gamma$,
		$r_2 \leftarrow 2\vartheta +2\gamma +1 -2 \beta_1$,
		$s_1 \leftarrow \frac{\vartheta +\gamma -2\beta_1}{2\vartheta +1-2\beta_1}$,
		$s_2 \leftarrow \frac{1+\vartheta -\gamma}{2\vartheta +1-2\beta_1}$)
	demonstrates that
	there exists a
	$C_1 \in (0,\infty)$
	such that for all
	$u \in H_{\gamma +2\vartheta}$ it holds that
	\begin{align}
	\nonumber
				&\Big(
					\|u\|^2_{B^{\vartheta +\gamma}_{2,2}((0,1)^2,\R^2)} +1
				\Big)
					\|u \|^2_{B^{1+\vartheta +\gamma}_{2,2}((0,1)^2,\R^2)}  \\ \nonumber
			\leq{} &
				C_1\Big( 
					\Big(
							\|u\|^{
								2(1+\vartheta +\gamma-2\beta_1)/(2\vartheta+1-2\beta_1)
							}_{B^{2\gamma}_{2,2}((0,1)^2,\R^2)}
							\|u\|^{
								2(\vartheta-\gamma)/(2\vartheta+1-2\beta_1)
							}_{
								B^{2\vartheta +2\gamma+1-2\beta_1}_{2,2}((0,1)^2,\R^2)
							} 
							+1
					\Big)
				\\ \nonumber & \qquad \cdot
					\big( 
						\|u \|^{
								2(\vartheta +\gamma-2\beta_1)/(2\vartheta+1-2\beta_1)
							}_{B^{2\gamma}_{2,2}((0,1)^2,\R^2)}
						\|u \|^{
								2(1+\vartheta-\gamma)/(2\vartheta+1-\beta_1)
							}_{B^{2\vartheta +2\gamma+1-2\beta_1}_{2,2}((0,1)^2,\R^2)}
					\big ) 
				\Big)  \\ \nonumber
			={} &
				C_1\Big( 
							\|u\|^{
								2(2\vartheta +2\gamma+1-4\beta_1)/(2\vartheta+1-2\beta_1)
							}_{B^{2\gamma}_{2,2}((0,1)^2,\R^2)}
							\|u\|^{
								2(2\vartheta-2\gamma+1)/(2\vartheta+1-2\beta_1)
							}_{
								B^{2\vartheta +2\gamma+1-2\beta_1}_{2,2}((0,1)^2,\R^2)
							} 
				\\ & \qquad 
					+\big( 
						\|u \|^{
								2(\vartheta +\gamma-2\beta_1)/(2\vartheta+1-2\beta_1)
							}_{B^{2\gamma}_{2,2}((0,1)^2,\R^2)}
						\|u \|^{
								2(1+\vartheta-\gamma)/(2\vartheta+1-\beta_1)
							}_{B^{2\vartheta +2\gamma+1-2\beta_1}_{2,2}((0,1)^2,\R^2)}
					\big )  
				\Big)  \\ \nonumber
			\leq{} &
				C_1 \Big( 
						(1+\|u\|^{3}_{B^{2\gamma}_{2,2}((0,1)^2,\R^2)})
						(1+\|u\|^{2}_{B^{2\vartheta +2\gamma+1-2\beta_1}_{2,2}((0,1)^2,\R^2)})
				\\ \nonumber & \qquad
					+(1+\|u\|^{3}_{B^{2\gamma}_{2,2}((0,1)^2,\R^2)})
						(1+\|u\|^{2}_{B^{2\vartheta +2\gamma+1-2\beta_1}_{2,2}((0,1)^2,\R^2)})
				\Big)  \\ \nonumber
			={} &
				2C_1 \Big( 
						(1+\|u\|^{3}_{B^{2\gamma}_{2,2}((0,1)^2,\R^2)})
						(1+\|u\|^{2}_{B^{2\vartheta +2\gamma+1-2\beta_1}_{2,2}((0,1)^2,\R^2)})
				\Big).
	\end{align}
	Hence
	Lemma \ref{lem: lemma: interpolation theorem for B},
	Corollary \ref{cor: sobolev holder inequality for positive s}
		(with $H_r \leftarrow B^{2r}_{2,2}((0,1)^2)$, $r \in \R$,
		$d\leftarrow 2$, $p\leftarrow 2$,
		$s \leftarrow \vartheta + \gamma -\nicefrac 12$,
		$\beta \leftarrow \beta$,
		$\alpha \leftarrow 0$,
		$r_1 \leftarrow \tfrac{\vartheta+\gamma}{2}$,
		$r_2 \leftarrow \tfrac{\vartheta+\gamma}{2}$,
		$\tilde{r}_1 \leftarrow \tfrac{1+\vartheta+\gamma}{2}$,
		$\tilde{r}_2 \leftarrow \tfrac{1+\vartheta+\gamma}{2}$,
		$\delta \leftarrow 0$,
		$u \leftarrow (D_{\R^2 \times \R^2}^{\alpha} b_{kl})(\cdot,u(\cdot))$,
		and
		$v \leftarrow u_i$),
	the fact that $\vartheta+\gamma < 1$,
	and Corollary \ref{cor: Bs b bound}
		(with $d\leftarrow 2$, $m \leftarrow 2$,
		$b \leftarrow 	D_{\R^2 \times \R^2}^{\alpha} b_{kl}$,
		$s \leftarrow \gamma +\vartheta$,
		$q\leftarrow \infty$,
		$p\leftarrow 0$)
	assure
	that
	for all $\alpha \in \N^4_0$,
	$\beta \in \N_0^2$
	with $|\alpha| =|\beta| = 1$ there exist
	$C_1,C_2,C_3,C_4 \in (0,\infty)$
	such that for all
	$u = (u_{1},u_{2}) \in H_{\gamma +2\vartheta}$
	and all $i \in \{1, 2\}$ it holds that
	\begin{align}
	\nonumber
				&\big \| 
						(D_{\R^2 \times \R^2}^{\alpha} b)(\cdot,u(\cdot))
						(D_{\R^2}^{\beta} u_i)(\cdot)
					\big \|^2_{(B^{2\vartheta +2\gamma-1}_{2,2}((0,1)^2))^{2\times 2}} \\ \nonumber
			\leq{} &
				C_1 \sum_{k,l \in \{1,2\} }
				\Big(
					\big \|  
						(D^\alpha_{\R^2 \times \R^2} b_{kl})(\cdot,u(\cdot)) 
					\big \|^2_{(B^{\vartheta +\gamma}_{2,2}((0,1)^2))} \, 
					\big\| 
						 u_i
					\big\|^2_{B^{1+\vartheta +\gamma}_{2,2}((0,1)^2)} 
				\Big) \\
	\label{eq: b Du bound}
			\leq{} &
				C_2 \sum_{k,l \in \{1,2\} }\Big( 
					\Big(
						\|u\|^2_{B^{\vartheta +\gamma}_{2,2}((0,1)^2,\R^2)} +1
						+\|u\|^2_{B^{0}_{2,2}((0,1)^2,\R^2)}
				\\ \nonumber & \qquad
						+\int_{(0,1)^2}|	(D^\alpha_{\R^2 \times \R^2} b_{kl})(x,0)|^2 \ud x
					\Big)
					\|u_i \|^2_{B^{1+\vartheta +\gamma}_{2,2}((0,1)^2)} 
				\Big)  \\ \nonumber
			\leq{} &
				C_3 
					\Big(
							\|u\|^2_{B^{\vartheta +\gamma}_{2,2}((0,1)^2,\R^2)} +1
					\Big)
					\|u \|^2_{B^{1+\vartheta +\gamma}_{2,2}((0,1)^2,\R^2)}  \\ \nonumber
			\leq{} &
				C_4 \Big( 
						(1+\|u\|^{3}_{B^{2\gamma}_{2,2}((0,1)^2,\R^2)})
						(1+\|u\|^{2}_{B^{2\vartheta +2\gamma+1-2\beta_1}_{2,2}((0,1)^2,\R^2)})
				\Big).  
	\end{align}
	Combining this,
	\eqref{eq: basic estimate H Navier},
	Proposition 2 on page 19 in Runst \& Sickel \cite{RunstSickel1996},
	and Corollary \ref{cor: Bs b bound}
		(with $d\leftarrow 2$, $m\leftarrow 2$,
		$s \leftarrow 2\gamma +2\vartheta-1$,
		$q\leftarrow \infty$,
		$p\leftarrow 0$) establishes
	that there exist $C_1, C_2,C_3, C_4 \in (0,\infty)$ such
	that for all $u=(u_1,u_2) \in H_{\gamma +2\vartheta}$ it holds that
	\begin{align}
				&\|B(u)\|^2_{HS(\U, \H_{\vartheta +\gamma})}
			\leq
				C_1\|
					b(\cdot,u(\cdot)) 
				\|^2_{(B^{2\gamma+2\vartheta}_{2,2}((0,1)^2,\R))^{2\times 2}} \\ \nonumber
			\leq{}
				&C_2\big \|
					\sum_{\alpha \in \N^2_0, |\alpha| \leq 1}
						D_{\R^2}^{\alpha} \big( b (\cdot,u(\cdot)) \big)
				\big \|^2_{(B^{2\gamma+2\vartheta-1}_{2,2}((0,1)^2,\R))^{2\times 2}} \\ \nonumber
			={}
				&C_2 \Big(
					\big \| 
						(D_{\R^2 \times \R^2}^{(1,0,0,0)} b)(\cdot,u(\cdot))
						+(D_{\R^2 \times \R^2}^{(0,0,1,0)} b)(\cdot,u(\cdot))
						(D_{\R^2}^{(1,0)} u_1)(\cdot)
				\\ \nonumber & \quad
						+(D_{\R^2 \times \R^2}^{(0,0,0,1)} b)(\cdot,u(\cdot))
						(D_{\R^2}^{(1,0)} u_2)(\cdot)
					\big \|^2_{(B^{2\vartheta +2\gamma-1}_{2,2}((0,1)^2))^{2\times 2}}
					+\big \| 
						(D_{\R^2 \times \R^2}^{(0,1,0,0)} b)(\cdot,u(\cdot))
				\\ \nonumber & \quad
						+(D_{\R^2 \times \R^2}^{(0,0,1,0)} b)(\cdot,u(\cdot))
						(D_{\R^2}^{(0,1)} u_1)(\cdot)
						+(D_{\R^2 \times \R^2}^{(0,0,0,1)} b)(\cdot,u(\cdot))
						(D_{\R^2}^{(0,1)} u_2)(\cdot)
					\big \|^2_{(B^{2\vartheta +2\gamma-1}_{2,2}((0,1)^2))^{2\times 2}} 
				\\ \nonumber & \quad
					+\|
						b(\cdot,u(\cdot))
					\|^2_{(B^{2\vartheta +2\gamma-1}_{2,2}((0,1)^2))^{2\times 2}}
				\Big)\\ \nonumber
			\leq{} &
				C_3 \Big( 
						(1+\|u\|^{3}_{B^{2\gamma}_{2,2}((0,1)^2,\R^2)})
						(1+\|u\|^{2}_{B^{2\vartheta +2\gamma+1-2\beta_1}_{2,2}((0,1)^2,\R^2)})
						+\|u\|^2_{B^{2\vartheta +2\gamma-1}_{2,2}((0,1)^2,\R^2)} 
				\\ \nonumber & \quad
						+	\|u \|^2_{B^{0}_{2,2}((0,1)^2,\R^2)}  	
						+1
				\Big)			\\ \nonumber
			\leq{} &
				C_4 \Big( 
						(1+\|u\|^{3}_{H_\gamma})
						(1+\|u\|^{2}_{H_{\vartheta +\gamma+1/2-\beta_1}})
						+\|u\|^2_{H_{\vartheta +\gamma-1/2}} 
						+	\|u \|^2_{H}  	
						+1
				\Big)	\\ \nonumber
			\leq{} &
				C_4 (1+\|u\|^{2}_{H_{\vartheta +\gamma+1/2-\beta_1}})
					\Big( 
						1+\|u\|^{3}_{H_\gamma}
						+\eta^{2(\beta_1-1)}
						+\eta^{-2\gamma}	\|u \|^2_{H_\gamma}  	
						+1
				\Big)
	\end{align}
	and this proves \eqref{eq: B bounded wrt H_vartheta Navier rem}.
	Moreover, the Sobolev inequality
	(see, e.g., Theorem 1 on page 32 in Runst \& Sickel \cite{RunstSickel1996}),
	the fact that $\vartheta \geq \nicefrac 12$ and
	Lemma \ref{l: norm equivalence Navier}
	yield that for all $\eps >0$ there exist $C_1, C_2 \in (0,\infty)$
	such that for all $u \in H_{\eps+\vartheta}$ it holds that
	\begin{align}
	\label{eq: L infty estimat Navier}
			\|u\|_{L^\infty((0,1)^2,\R^2)}
		\leq
			C_1\|u\|_{B^{2\vartheta+2\eps}_{2,2}((0,1)^2;\R^2)}
		\leq
			C_2\|u\|_{H_{\vartheta+\eps}}.
	\end{align}
	Furthermore, it follows from
	Theorem 2 on page 191 in Runst \& Sickel \cite{RunstSickel1996},
	Corollary \ref{cor: Bs b bound}
		(with $d\leftarrow 2$, $m \leftarrow 2$,
		$b \leftarrow 	D_{\R^2 \times \R^2}^{\alpha} b_{kl}$,
		$s \leftarrow 2\gamma +2\vartheta-1$,
		$q\leftarrow \infty$,
		$p\leftarrow 0$),
	and from Lemma \ref{l: norm equivalence Navier} 
	that there exist $C_1, C_2, C_3 \in (0,\infty)$ such that for all
	$v=(v_1,v_2) \in H_{\gamma +2\vartheta}$,
	$u=(u_1,u_2) \in H_{\gamma +2\vartheta}$,
	$i \in \{1,2\}$,
	$\alpha \in \N^4_0$,
	and all $\beta \in \N_0^2$
	with $|\alpha| =|\beta| = 1$
	it holds that
	\begin{align}
				&\Big \|
					(D_{\R^2 \times \R^2}^{\alpha} b)(\cdot,v(\cdot))
					\big(
						(D_{\R^2}^{\beta} v_i)(\cdot)-(D_{\R^2}^{\beta} u_i)(\cdot)
					\big)
				\Big \|^2_{(B^{2\vartheta +2\gamma-1}_{2,2}((0,1)^2))^{2\times 2}} \\ \nonumber
			\leq{}
				&C_1
				\sum_{k,l \in \{1,2\}} \Big(
					\|
						(D_{\R^2 \times \R^2}^{\alpha} b_{kl})(\cdot,v(\cdot))
					\|^2_{B^{2\vartheta +2\gamma-1}_{2,2}((0,1)^2)}
					\|
						D_{\R^2}^{\beta} v_i-D_{\R^2}^{\beta} u_i
					\|^2_{B^{2\vartheta +2\gamma-2\beta_1}_{2,2}((0,1)^2)} 
				\Big)\\ \nonumber
			\leq{}
				&C_2
				\sum_{k,l \in \{1,2\}} \bigg(
					\big(
						(\|v\|^2_{B^{2\vartheta +2\gamma-1}_{2,2}((0,1)^2,\R^2)} +1)
						+\|v\|^2_{B^{0}_{2,2}((0,1)^2,\R^2)} 
						+\int_{(0,1)^2}
						|	(D_{\R^2 \times \R^2}^{\alpha} b_{kl})(x,0)|^2
						\ud x
					\big)
		\\ \nonumber & \quad \cdot
					\|
						v- u
					\|^2_{B^{2\vartheta +2\gamma+1-2\beta_1}_{2,2}((0,1)^2,\R^2)} 
				\bigg)\\ \nonumber
			\leq{}
				&C_3
				\sum_{k,l \in \{1,2\}} \bigg(
					\big(
						\|v\|^2_{H_{\vartheta +\gamma-1/2}} +1
						+\|v\|^2_{H} 
						+\int_{(0,1)^2}
						|	(D_{\R^2 \times \R^2}^{\alpha} b_{kl})(x,0)|^2
						\ud x
					\big)
					\|
						v- u
					\|^2_{H_{\vartheta +\gamma+1/2-\beta_1}} 
				\bigg).
	\end{align}
	In addition, we obtain from Theorem 2 on page 191 in
	Runst \& Sickel \cite{RunstSickel1996},
	Corollary \ref{cor: Bs b continuity inequality}
		(with 
		$d \leftarrow 2$,
		$m \leftarrow 2$,
		$b \leftarrow 	D_{\R^2 \times \R^2}^{\alpha} b_{kl}$,
		$s \leftarrow 2\vartheta +2\gamma -1$,
		$q_1 \leftarrow \infty$, $q_2 \leftarrow 2$, $q_3 \leftarrow \infty$,
		$q_4 \leftarrow 2$, $q_5 \leftarrow \infty$,
		$q \leftarrow \infty$,
		$p_1 \leftarrow \infty$, $p_2 \leftarrow 0$,
		$\eps \leftarrow 2\gamma -2\beta_1$),
	Lemma \ref{l: norm equivalence Navier},
	and \eqref{eq: L infty estimat Navier}
	(with $\eps \leftarrow \gamma -\beta_1$)
	that there exist $C_1, C_2, C_3, C_4\in (0,\infty)$ such that for all
	$u =(u_1,u_2)\in H_{\gamma +2\vartheta}$,
	$v =(v_1,v_2) \in H_{\gamma +2\vartheta}$,
	$i \in \{1,2\}$,
	$\alpha \in \N^4_0$,
	and all $\beta \in \N_0^2$
	with $|\alpha| =|\beta| = 1$
	it holds that
	\begin{align}
				&\Big \| 
					\big(
						(D_{\R^2 \times \R^2}^{\alpha} b)(\cdot,u(\cdot))
						-(D_{\R^2 \times \R^2}^{\alpha} b)(\cdot,v(\cdot))
					\big)
						(D_{\R^2}^{\beta} u_i)(\cdot)
				\Big \|^2_{(B^{2\vartheta +2\gamma-1}_{2,2}((0,1)^2))^{2\times 2}} \\ \nonumber
			\leq
				&C_1
					\sum_{k,l \in \{1,2\}} \bigg(
						\big \| 
							(D_{\R^2 \times \R^2}^{\alpha} b_{kl})(\cdot,u(\cdot))
							-(D_{\R^2 \times \R^2}^{\alpha} b_{kl})(\cdot,v(\cdot))
						\big \|^2_{B^{2\vartheta +2\gamma-1}_{2,2}((0,1)^2)}
				\\ \nonumber & \quad
						\big \| 
							D_{\R^2}^{\beta} u_i
						\big \|^2_{B^{2\vartheta +2\gamma-2 \beta_1}_{2,2}((0,1)^2)} 
					\bigg)\\ \nonumber
			\leq
				&C_2
				\sum_{k,l \in \{1,2\}} \bigg(
					\Big(
						\|
							u-v
						\|^{2}_{B^{1+2\gamma-2\beta_1 }_{2,2}((0,1)^2,\R^2)}
						\big(
							\|u\|^2_{B^{2\vartheta +2\gamma-1}_{2,2}((0,1)^2,\R^2)}
							+1
						\big) 
				\\ \nonumber & \quad \cdot 
						\| 
							D^2_{\R^2 \times \R^2} D^\alpha_{\R^2 \times \R^2} b_{kl}
						\|_{
							\C(
								[0,1]
								\times 
								\mathbf{B}(
									\|u\|_{L^\infty((0,1)^2,\R^2)} \vee \|v\|_{L^\infty((0,1)^2,\R^2)}
								),
								\| \cdot \|_{L(\R^2 \times \R^2, (\R^2 \times \R^2)')}
							)
						}  \\ \nonumber
				& \quad
						+\|
							u - v 
						\|^2_{B^{2\vartheta +2\gamma-1}_{2,2}((0,1)^2,\R^2)}
						+\|u - v \|^2_{L^2((0,1)^2,\R^2)}
					\Big)
					\| 
						u
					\|^2_{B^{2\vartheta +2\gamma+1-2 \beta_1}_{2,2}((0,1)^2,\R^2)}
				\bigg)\\ \nonumber
			\leq
				&C_3
				\sum_{k,l \in \{1,2\}} \bigg(
					\Big(
						\|
							u-v
						\|^{2}_{H_{1/2+\gamma-\beta_1 }}
						\big(
							\|u\|^2_{H_{\vartheta +\gamma-1/2}}
							+1
						\big) 
			\\ \nonumber & \quad \cdot 
						\| 
							D^2_{\R^2 \times \R^2} D^\alpha_{\R^2 \times \R^2} b_{kl}
						\|_{
							\C(
								[0,1]
								\times 
								\mathbf{B}(
									(C_4 \|u\|_{H_{\gamma + \vartheta + 1/2 -\beta_1}})
									\vee (C_4 \|v\|_{H_{\gamma + \vartheta + 1/2 -\beta_1}})
								),
								\| \cdot \|_{L(\R^2 \times \R^2, (\R^2 \times \R^2)')}
							)
						}  \\ \nonumber
			& \quad
						+\|
							u - v 
						\|^2_{H_{\vartheta +\gamma-1/2}}
						+\|u - v \|^2_{H}
					\Big)
					\| 
						u
					\|^2_{H_{\vartheta +\gamma+1/2-\beta_1}}  
				\bigg).
	\end{align}
	Thus we get that there exist $C_1, C_2 \in (0,\infty)$ such that for all
	$u =(u_1,u_2) \in H_{\gamma +2\vartheta}$,
	$v=(v_1,v_2) \in H_{\gamma +2\vartheta}$,
	$i \in \{1,2\}$,
	$\alpha \in \N^4_0$,
	and all $\beta \in \N_0^2$
	with $|\alpha| =|\beta| = 1$
	it holds that
	\begin{align}
	\label{eq: b diff Dv estimate}
				&\big \| 
						(D_{\R^2 \times \R^2}^{\alpha} b)(\cdot,u(\cdot))
						(D_{\R^2}^{\beta} u_i)(\cdot)
						-(D_{\R^2 \times \R^2}^{\alpha} b)(\cdot,v(\cdot))
						(D_{\R^2}^{\beta} v_i)(\cdot)
					\big \|^2_{(B^{2\vartheta +2\gamma-1}_{2,2}((0,1)^2))^{2\times 2}} \\ \nonumber
			\leq
				&\Big \| 
					\big(
						(D_{\R^2 \times \R^2}^{\alpha} b)(\cdot,u(\cdot))
						-(D_{\R^2 \times \R^2}^{\alpha} b)(\cdot,v(\cdot))
					\big)
						(D_{\R^2}^{\beta} u_i)(\cdot)
				\Big \|^2_{(B^{2\vartheta +2\gamma-1}_{2,2}((0,1)^2))^{2\times 2}}
			\\ \nonumber & \quad
				+\Big \|
					(D_{\R^2 \times \R^2}^{\alpha} b)(\cdot,v(\cdot))
					\big(
						(D_{\R^2}^{\beta} v_i)(\cdot)-(D_{\R^2}^{\beta} u_i)(\cdot)
					\big)
				\Big \|^2_{(B^{2\vartheta +2\gamma-1}_{2,2}((0,1)^2))^{2\times 2}} \\ \nonumber
			\leq
				&C_1
				\sum_{k,l \in \{1,2\}} \bigg(
					\Big(
						\|
							u-v
						\|^{2}_{H_{1/2+\gamma-\beta_3 }}
						\big(
							\|u\|^2_{H_{\vartheta +\gamma-1/2}}
							+1
						\big) 
			\\ \nonumber & \quad \cdot 
						\| 
							D^2_{\R^2 \times \R^2} D^\alpha_{\R^2 \times \R^2} b_{kl}
						\|_{
							\C(
								[0,1]
								\times 
								\mathbf{B}(
									(C_2 \|u\|_{H_{\gamma + \vartheta + 1/2 -\beta_1}})
									\vee (C_2 \|v\|_{H_{\gamma + \vartheta + 1/2 -\beta_1}})
								),
								\| \cdot \|_{L(\R^2 \times \R^2, (\R^2 \times \R^2)')}
							)
						}  \\ \nonumber
			& \quad
						+\|
							u - v 
						\|^2_{H_{\vartheta +\gamma-1/2}}
						+\|u - v \|^2_{H}
					\Big)
					\| 
						u
					\|^2_{H_{\vartheta +\gamma+1/2-\beta_1}}
			\\ 
	\nonumber
			& \quad
					+
					\big(
						\|u\|^2_{H_{\vartheta +\gamma-1/2}} +1
						+\|u\|^2_{H} 
						+\int_{(0,1)^2}
						|	(D_{\R^2 \times \R^2}^{\alpha} b_{kl})(x,0)|^2
						\ud x
					\big)
					\|
						v- u
					\|^2_{H_{\vartheta +\gamma+1/2-\beta_1}}  
				\bigg).
	\end{align}
	Furthermore, Corollary \ref{cor: Bs b continuity inequality}
		(with 
		$d \leftarrow 2$,
		$m \leftarrow 2$,
		$b \leftarrow 	D_{\R^2 \times \R^2}^{\alpha} b_{kl}$,
		$s \leftarrow 2\vartheta +2\gamma -1$,
		$q_1 \leftarrow \infty$, $q_2 \leftarrow 2$, $q_3 \leftarrow \infty$,
		$q_4 \leftarrow 2$, $q_5 \leftarrow \infty$,
		$q \leftarrow \infty$,
		$p_1 \leftarrow \infty$, $p_2 \leftarrow 0$,
		$\eps \leftarrow 2\gamma -2\beta_1$),
	Lemma \ref{l: norm equivalence Navier},
	and \eqref{eq: L infty estimat Navier} (with $\eps \leftarrow \gamma-\beta_1$)
	imply
	that there exist $C_1, C_2, C_3 \in (0,\infty)$ such that for all
	$u,v \in H_{\gamma +2\vartheta}$,
	$\alpha \in \N^4_0$,
	and all $\beta \in \N_0^2$
	with $|\alpha| =|\beta| = 1$
	it holds that
	\begin{align}
	\label{eq: b diff bound Navier}
				&\big \| 
						(D_{\R^2 \times \R^2}^{\alpha} b)(\cdot,u(\cdot))
						-(D_{\R^2 \times \R^2}^{\alpha} b)(\cdot,v(\cdot))
					\big \|^2_{(B^{2\vartheta +2\gamma-1}_{2,2}((0,1)^2))^{2\times 2}} \\ \nonumber
			\leq{} &
				C_1
				\sum_{k,l \in \{1,2\}} \Big(
					\|
						u-v
					\|^{2}_{B^{1+2\gamma-2\beta_3}_{2,2}((0,1)^2,\R^2)}
					\big(
						\|u\|^2_{B^{2\vartheta +2\gamma-1}_{2,2}((0,1)^2,\R^2)}
						+1
					\big)
			\\ \nonumber & \quad \cdot 
						\| 
							D^2_{\R^2 \times \R^2} D^\alpha_{\R^2 \times \R^2} b_{kl}
						\|_{
							\C(
								[0,1]
								\times 
								\mathbf{B}(
									\|u\|_{L^\infty((0,1)^2,\R^2)} \vee \|v\|_{L^\infty((0,1)^2,\R^2)}
								),
								\| \cdot \|_{L(\R^2 \times \R^2, (\R^2 \times \R^2)')}
							)
						} \\ \nonumber
			& \qquad\qquad \quad
					+\|
						u - v 
					\|^2_{B^{2\vartheta +2\gamma-1}_{2,2}((0,1)^2,\R^2)}
					+\|u - v \|^2_{L^{2}((0,1)^2,\R^2)}
				\Big) \\ \nonumber
			\leq{} &
				C_2 
				\sum_{k,l \in \{1,2\}} \Big(
					\|
						u-v
					\|^{2}_{H_{\gamma+1/2-\beta_3}}
					\big(
						\|u\|^2_{H_{\vartheta +\gamma-1/2}}
						+1
					\big)
			\\ \nonumber & \quad \cdot 
						\| 
							D^2_{\R^2 \times \R^2} D^\alpha_{\R^2 \times \R^2} b_{kl}
						\|_{
							\C(
								[0,1]
								\times 
								\mathbf{B}(
									(C_3 \|u\|_{H_{\gamma +\vartheta +1/2 -\beta_3}})
									\vee ( C_3 \|v\|_{H_{\gamma +\vartheta +1/2 -\beta_3}})
								),
								\| \cdot \|_{L(\R^2 \times \R^2, (\R^2 \times \R^2)')}
							)
						} \\ \nonumber
			& \qquad\qquad \quad
					+\|
						u - v 
					\|^2_{H_{\vartheta +\gamma-1/2}}
					+\|u - v \|^2_{H}
				\Big).
	\end{align}
	Therefore 
	\eqref{eq: basic estimate diff Navier}
	and
	Proposition 2 on page 19 in Runst \& Sickel \cite{RunstSickel1996}
	yields that there exist $C_1,C_2 \in (0,\infty)$
	such that for all $u =(u_1,u_2) \in H_{\gamma +2\vartheta}$,
	$v =(v_1,v_2) \in H_{\gamma +2\vartheta}$
	it holds that 
	\begin{align}
				&\|B(u)-B(v)\|^2_{HS(\U, \H_{\vartheta +\gamma})}
			\leq
				C_1\|
					b(\cdot,u(\cdot)) -b(\cdot,v(\cdot))
				\|^2_{(B^{2\gamma+2\vartheta}_{2,2}((0,1)^2,\R))^{2\times 2}} \\ \nonumber
			\leq{}
				&C_2
					\sum_{\alpha \in \N^2_0, |\alpha| \leq 1}
						\big \|
							D_{\R^2}^{\alpha} \big( b (\cdot,u(\cdot)) -b(\cdot,v(\cdot)) \big)
						\big \|^2_{(B^{2\gamma+2\vartheta-1}_{2,2}((0,1)^2,\R))^{2\times 2}} \\ \nonumber
			={}
				&C_2 \Big(
					\big \| 
						(D_{\R^2 \times \R^2}^{(1,0,0,0)} b)(\cdot,u(\cdot))
						-(D_{\R^2 \times \R^2}^{(1,0,0,0)} b)(\cdot,v(\cdot))
						+(D_{\R^2 \times \R^2}^{(0,0,1,0)} b)(\cdot,u(\cdot))
						(D_{\R^2}^{(1,0)} u_1)(\cdot) 
				\\ \nonumber & \quad
					-(D_{\R^2 \times \R^2}^{(0,0,1,0)} b)(\cdot,v(\cdot))
						(D_{\R^2}^{(1,0)} v_1)(\cdot)
						+(D_{\R^2 \times \R^2}^{(0,0,0,1)} b)(\cdot,u(\cdot))
						(D_{\R^2}^{(1,0)} u_2)(\cdot)
				\\ \nonumber & \quad
						-(D_{\R^2 \times \R^2}^{(0,0,0,1)} b)(\cdot,v(\cdot))
						(D_{\R^2}^{(1,0)} v_2)(\cdot)
					\big \|^2_{(B^{2\vartheta +2\gamma-1}_{2,2}((0,1)^2))^{2\times 2}}
					+\big \| 
						(D_{\R^2 \times \R^2}^{(0,1,0,0)} b)(\cdot,u(\cdot))
				\\ \nonumber  & \quad
						-(D_{\R^2 \times \R^2}^{(0,1,0,0)} b)(\cdot,v(\cdot))
						+(D_{\R^2 \times \R^2}^{(0,0,1,0)} b)(\cdot,u(\cdot))
						(D_{\R^2}^{(0,1)} u_1)(\cdot)
						-(D_{\R^2 \times \R^2}^{(0,0,1,0)} b)(\cdot,v(\cdot))
						(D_{\R^2}^{(0,1)} v_1)(\cdot)
					\\ \nonumber & \quad
						+(D_{\R^2 \times \R^2}^{(0,0,0,1)} b)(\cdot,u(\cdot))
						(D_{\R^2}^{(0,1)} u_2)(\cdot)
						-(D_{\R^2 \times \R^2}^{(0,0,0,1)} b)(\cdot,v(\cdot))
						(D_{\R^2}^{(0,1)} v_2)(\cdot)
					\big \|^2_{(B^{2\vartheta +2\gamma-1}_{2,2}((0,1)^2))^{2\times 2}} 
				\\ \nonumber & \quad
					+\|
						b(\cdot,u(\cdot)) - b(\cdot,v(\cdot))
					\|^2_{(B^{2\vartheta +2\gamma-1}_{2,2}((0,1)^2))^{2\times 2}}
				\Big).
	\end{align}
	Combining this, \eqref{eq: b diff Dv estimate},
	\eqref{eq: b diff bound Navier},
	the assumption that for all $k,l \in \{1,2\}$
	it holds that  $b_{kl} \in \C^{3}([0,1]^2 \times \R^2; \R)$,
	and the fact that for all $r \in (-\infty,\vartheta +\gamma -\beta_1]$
	and all $u \in H_r$
	it holds that
	$
			\|u\|_{H_r} 
		\leq 
			\eta^{r-\vartheta-\gamma+\beta_1} \|u\|_{H_{\vartheta +\gamma -\beta_1}}
	$
	proves \eqref{eq: B continuity Navier rem}.
	Moreover, 
	\eqref{eq: def of PN 1}
	and
	\eqref{eq: L infty estimat Navier} (with $\eps \leftarrow \gamma$)
	show that
	there exist $C_1,C_2 \in (0,\infty)$ such 
	that for all $v \in H$, and all $N \in \N$ it holds that
	\begin{align}
			\nonumber
				&\|
					B( v ) ( \id_{H} - \mathscr{P}_N )
				\|^2_{ HS( \mathbb{U}, \mathbb{H} ) } 
			\leq{}
					\sum_{n\in \mathcal{I}}
						\|
							b(\cdot, v(\cdot))  
								(Q ( \id_{H} - \mathscr{P}_N )  \tilde{e}_{n})(\cdot)
						\|^2_{L^2{((0,1)^2,\R^2)}} \\
			={} \nonumber
				&\sum_{n\in \mathcal{I} \backslash [1,N]}
						\|
							b(\cdot, v(\cdot))  
								(Q \tilde{e}_{n})(\cdot)
						\|^2_{L^2{((0,1)^2,\R^2)}} \\
			\begin{split}
			\leq{}
				&\sum_{n\in \mathcal{I} \backslash [1,N]}
						\big \|
							\| b(\cdot, v(\cdot)) \|_{HS(\R^2,\R^2)}
								\| (Q \tilde{e}_{n})(\cdot) \|_{\R^2}
						\big \|^2_{L^2{((0,1)^2)}} \\
			\leq{} 
					&C_1\sum_{n\in \mathcal{I} \backslash [1,N]}
						\big \|
							\| b(\cdot, v(\cdot)) \|_{HS(\R^2,\R^2)}  
						\big \|^2_{L^2{((0,1)^2)}}
						\big \|
								\| (Q \tilde{e}_{n})(\cdot) \|_{\R^2}
						\big \|^2_{L^\infty((0,1)^2)} 
		\end{split} \\    
		\nonumber
			={} 
					&C_1\sum_{n\in \mathcal{I} \backslash [1,N]}
						\|
							b(\cdot, v(\cdot))   
						\|^2_{(L^2{((0,1)^2))^{2 \times 2}}}
						\|
							(Q \tilde{e}_{n})(\cdot) 
						\|^2_{L^\infty((0,1)^2,\R^2)} \\   
			\leq{} \nonumber
					&C_2
					\|
						b(\cdot, v(\cdot))
					\|^2_{(L^2{((0,1)^2))^{2\times 2}}}
					\sum_{n\in \mathcal{I} \backslash [1,N]}
						\|
								Q  \tilde{e}_{n}
						\|^2_{H_{\gamma+\vartheta}}.
		\end{align}
		Combining this and 
		Lemma \ref{lem: b L2 boundedness inequality} 
			(with $d \leftarrow 2$, $m \leftarrow 2$, $p \leftarrow 0$),
		establishes that
		there exists a $C_1 \in (0,\infty)$
		such that for all $\mathbb{H}_\gamma$-bounded sets $E$ it holds that
		\begin{align}
		\nonumber
				&\limsup_{ N \to \infty }
				\sup_{ v \in E }
				\left[
				\|
					B( v ) ( \id_{H} - \mathscr{P}_N )
				\|^2_{ HS( \mathbb{U}, \mathbb{H} ) }
				\right] \\
			\leq{} \nonumber & 
				C_1\sup_{ v \in E } (
					\|
						b(\cdot, v(\cdot))
					\|^2_{(L^2{((0,1)^2))^{2\times 2}}}
				)
				\limsup_{ N \to \infty }
					\Big(
						\sum_{n\in \mathcal{I} \backslash [1,N]}
						\|
								Q  \tilde{e}_{n}
						\|^2_{H_{\gamma+\vartheta}}
					\Big) \\
			\leq{} \nonumber &
				C_1
				\sup_{ v \in E } \Big(
					\sum_{k,l \in \{1,2\}}
					\big(
						\| v\|^2_{L^2((0,1),\R^2)}
						+2\int_{(0,1)^2}
							|b_{kl}(x,0)|^2
						\ud x
					\big)
				\Big)
				\limsup_{ N \to \infty }
					\Big(
						\sum_{n\in \mathcal{I} \backslash [1,N]}
						\|
								Q  \tilde{e}_{n}
						\|^2_{H_{\gamma+\vartheta}}
					\Big) \\
			={} 
				&0
	\end{align}
	and this assures \eqref{eq:convLocLip_assumption Navier rem}.
	Finally observe that we derive from
	\eqref{eq: B Lip cont Navier rem}
	that there exists a $C_1 \in (0,\infty)$
	such that for all $\mathbb{H}_\gamma$-bounded sets $E$ it holds that
	\begin{equation}
		\begin{split}
				&\sup_{ N \in \N } \sup_{ v \in E }
					\|
						B( v ) \mathscr{P}_N 
					\|_{ HS( \mathbb{U}, \mathbb{H} ) }
			\leq	
				\sup_{ v \in E}
					\|
						B( v )
					\|_{ HS( \mathbb{U}, \mathbb{H} ) }\\
			\leq{}	
				&\sup_{ v \in E}
					\big(
						\|
							B( v ) -B(0)
						\|_{ HS( \mathbb{U}, \mathbb{H} ) }
						+\|
							B(0)
						\|_{ HS( \mathbb{U}, \mathbb{H} ) }
					\big)	\\
			\leq{}	
				&\sup_{ v \in E}
					C_1\big(
						\|
							v
						\|^2_{H}
						+\|
							B(0)
						\|_{ HS( \mathbb{U}, \mathbb{H} ) }
					\big)	
			< \infty.
		\end{split}
	\end{equation}
	This verifies \eqref{eq: B C1 bound Navier rem} 
	and thus finishes the proof of Remark \ref{rem: sufficient condition for B Navier}.
\end{proof}

\subsubsection*{Acknowledgements}
This work has been partially supported by the Deutsche Forschungsgemeinschaft (DFG) via RTG 2131
"High-dimensional phenomena in probability - Fluctuations and discontinuity".

\printbibliography[title={Bibliography}]
\end{document}